\documentclass[11pt]{scrartcl}

\newcounter{variable}
\usepackage{ifthen}
\newcommand{\twooptions}[2]{\ifthenelse{\value{variable}=0}{#1}{#2}}

\title{%
Optimal Bounds for the $k$-Disjoint Paths Problem%
\twooptions{}{\thanks{Maximilian Gorsky was supported by the Institute for Basic Science (IBS-R029-C1). Dimitrios M.\@ Thilikos was supported by French National Research Agency (ANR) under project GODASse ANR-24-CE48-4377 and under the France 2030 grant reference number ANR-24-RRII-0002 operated by the Inria Quadrant Program. Emails of the authors:
\href{mailto:d.cavallaro@tu-berlin.de}{d.cavallaro@tu-berlin.de},
\href{mailto:m.gorsky@pm.me}{m.gorsky@pm.me},
\href{mailto:stephan.kreutzer@tu-berlin.de}{stephan.kreutzer@tu-berlin.de},
\href{mailto:sedthilk@thilikos.info}{sedthilk@thilikos.info},
\href{mailto:wiederrecht@kaist.ac.kr}{wiederrecht@kaist.ac.kr}.}
}\vspace{-7mm}}

\author{%
\twooptions{Anonymous Author(s)}{%
 \parbox{\textwidth}{\centering
 Dario Cavallaro\thanks{Technische Universit\"at Berlin, Germany.}\qquad
 Maximilian Gorsky\thanks{Discrete Mathematics Group, Institute for Basic Science (IBS), Daejeon, South Korea.}\qquad
 Stephan Kreutzer\footnotemark[2]\\[1ex]
 \hspace{0.5em}Dimitrios M. Thilikos\thanks{LIRMM, Univ Montpellier, CNRS, Montpellier, France.}\qquad
 Sebastian Wiederrecht\thanks{School of Computing, KAIST, Daejeon, South Korea.}
 }%
}%
}
\date{}

\usepackage[utf8]{inputenc}
\usepackage[T1]{fontenc}
\usepackage{lmodern}
\usepackage{alphabeta}

\usepackage{titling}
\usepackage{xspace}
\usepackage{soul} %

\usepackage{mathrsfs}

\usepackage[a4paper, top=25.4mm, bottom=25.4mm, left=25.4mm, right=25.4mm, includefoot]{geometry}

\DeclareSectionCommand[%
level=4,
indent=0pt,
beforeskip=1ex plus 1ex minus .2ex,
afterskip=-1em,
font={},
tocindent=7em,
tocnumwidth=4.1em,
counterwithin=subsubsection
]{paragraph}

\usepackage[inline,shortlabels]{enumitem}

\usepackage{dsfont}
\usepackage{setspace}

\usepackage{bm}
\usepackage{microtype}
\usepackage{amsmath}
\usepackage{amssymb}
\usepackage{amsfonts}
\usepackage{mathtools}
\usepackage[mathscr]{euscript}
\usepackage{wasysym}

\usepackage[hyphens]{url}
\usepackage{amsthm}

\usepackage{tikz}
\usepackage{xcolor}
\usepackage{subcaption}
\usepackage{graphicx}
\usepackage{wrapfig}

\usepackage{cite}

\usepackage{hyperref}
\usepackage{zref-clever}

\usepackage{nicefrac}
\usepackage[textwidth=1.5in]{todonotes}

\newcommand{\cref}[1]{\zcref{#1}}
\renewcommand{\autoref}[1]{\zcref{#1}}

\definecolor{mediumslateblue}{rgb}{0.48, 0.41, 0.93}
\definecolor{lightsalmonpink}{rgb}{1.0, 0.6, 0.6}
\definecolor{amaranth}{rgb}{0.9, 0.17, 0.31}
\definecolor{jonquil}{rgb}{0.98, 0.85, 0.37}
\definecolor{orangepeel}{rgb}{1.0, 0.62, 0.0}
\definecolor{darkspringgreen}{rgb}{0.09, 0.45, 0.27}
\colorlet{myGreen}{green!50!black}
\colorlet{myLightgreen}{green}
\colorlet{myRed}{red!90!black}
\definecolor{myBlue}{rgb}{0.25, 0.0, 1.0}
\definecolor{myLightBlue}{rgb}{0.39, 0.58, 0.93}
\colorlet{myViolet}{myBlue!55!myRed}
\definecolor{myOrange}{rgb}{1.0, 0.66, 0.07}
\definecolor{ceruleanblue}{rgb}{0.16, 0.32, 0.75}
\definecolor{CornflowerBlue}{rgb}{0.39, 0.58, 0.93}
\definecolor{DarkGoldenrod}{rgb}{0.72, 0.53, 0.04}
\definecolor{BritishRacingGreen}{rgb}{0.0, 0.26, 0.15}
\definecolor{DarkMagenta}{rgb}{0.55, 0.0, 0.55}
\definecolor{AO}{rgb}{0.0, 0.5, 0.0}
\definecolor{BostonUniversityRed}{rgb}{0.8, 0.0, 0.0}
\definecolor{myRed}{rgb}{0.8, 0.0, 0.0}
\definecolor{DarkMidnightBlue}{rgb}{0.0, 0.2, 0.4}
\definecolor{DarkTangerine}{rgb}{1.0, 0.66, 0.07}
\definecolor{AppleGreen}{rgb}{0.55, 0.71, 0.0}
\definecolor{BrightUbe}{rgb}{0.82, 0.62, 0.91}
\definecolor{Amethyst}{rgb}{0.6, 0.4, 0.8}
\definecolor{DarkGray}{rgb}{0.52, 0.52, 0.51}
\definecolor{Gray}{rgb}{0.66, 0.66, 0.66}
\definecolor{BananaYellow}{rgb}{1.0, 0.88, 0.21}
\definecolor{Amber}{rgb}{1.0, 0.75, 0.0}
\definecolor{LightGray}{rgb}{0.83, 0.83, 0.83}
\definecolor{PrincetonOrange}{rgb}{1.0, 0.56, 0.0}
\definecolor{DeepCarrotOrange}{rgb}{0.91, 0.41, 0.17}
\definecolor{CarrotOrange}{rgb}{0.93, 0.57, 0.13}
\definecolor{MidnightBlue}{rgb}{0.1, 0.1, 0.44}
\definecolor{Magenta}{rgb}{0.50, 0.0, 0.50}
\definecolor{BrightPink}{rgb}{1.0, 0.0, 0.5}
\definecolor{BrilliantRose}{rgb}{1.0, 0.33, 0.64}
\definecolor{ChromeYellow}{rgb}{1.0, 0.65, 0.0}
\definecolor{HotMagenta}{rgb}{1.0, 0.11, 0.81}
\definecolor{Auburn}{rgb}{0.43, 0.21, 0.1}
\definecolor{BrightTurquoise}{rgb}{0.03, 0.91, 0.87}
\definecolor{DarkCyan}{rgb}{0.0, 0.55, 0.55}

\vfuzz \hfuzz
\setlist[itemize]{topsep=0pt,partopsep=0pt,itemsep=0pt,parsep=0pt}
\setlist[itemize,1]{label={\small\textbullet}}
\setlist[itemize,2]{label={\tiny\textbullet}}
\setlist[itemize,3]{label=$\cdot$}
\setlist[enumerate]{topsep=0pt,partopsep=0pt,itemsep=0pt,parsep=0pt}
\setlist[enumerate,1]{label=\roman*)}
\setlist[enumerate,2]{label=\alph*)}
\setlist[enumerate,3]{label=\arabic*)}

\hypersetup{
colorlinks=true,
linkcolor=AO!65!black,
citecolor=AO!65!black,
urlcolor=AppleGreen!65!black,
bookmarksopen=true,
bookmarksnumbered,
bookmarksopenlevel=2,
bookmarksdepth=3
}

\newenvironment{claimproof}[1][Proof of claim]
  {\begin{proof}[#1]}
  {\end{proof}}

\theoremstyle{definition}

\newtheorem{environment}{Environment}[section]

\newtheorem{lemma}[environment]{Lemma}
\AddToHook{env/lemma/begin}{\zcsetup{countertype={environment=lemma}}}
\zcRefTypeSetup{lemma}{
Name-sg = Lemma ,
name-sg = Lemma ,
Name-pl = Lemmas ,
name-pl = Lemmas ,
}

\newtheorem*{lemma*}{Lemma}
\AddToHook{env/lemma*/begin}{\zcsetup{countertype={environment=lemma*}}}
\zcRefTypeSetup{lemma*}{
Name-sg = Lemma ,
name-sg = Lemma ,
Name-pl = Lemmas ,
name-pl = Lemmas ,
}

\newtheorem{proposition}[environment]{Proposition}
\AddToHook{env/proposition/begin}{\zcsetup{countertype={environment=proposition}}}
\zcRefTypeSetup{proposition}{
Name-sg = Proposition ,
name-sg = Proposition ,
Name-pl = Propositions ,
name-pl = Propositions ,
}

\newtheorem{corollary}[environment]{Corollary}
\AddToHook{env/corollary/begin}{\zcsetup{countertype={environment=corollary}}}
\zcRefTypeSetup{corollary}{
Name-sg = Corollary ,
name-sg = Corollary ,
Name-pl = Corollaries ,
name-pl = Corollaries ,
}

\newtheorem{theorem}[environment]{Theorem}
\AddToHook{env/theorem/begin}{\zcsetup{countertype={environment=theorem}}}
\zcRefTypeSetup{theorem}{
Name-sg = Theorem ,
name-sg = Theorem ,
Name-pl = Theorems ,
name-pl = Theorems ,
}

\newtheorem*{theorem*}{Theorem}
\AddToHook{env/theorem*/begin}{\zcsetup{countertype={environment=theorem*}}}
\zcRefTypeSetup{theorem*}{
Name-sg = Theorem ,
name-sg = Theorem ,
Name-pl = Theorems ,
name-pl = Theorems ,
}

\AddToHook{env/conjecture/begin}{\zcsetup{countertype={environment=conjecture}}}
\zcRefTypeSetup{conjecture}{
Name-sg = Conjecture ,
name-sg = Conjecture ,
Name-pl = Conjectures ,
name-pl = Conjectures ,
}

\newtheorem*{hypothesis*}{Hypothesis}
\AddToHook{env/hypothesis*/begin}{\zcsetup{countertype={environment=hypothesis*}}}
\zcRefTypeSetup{hypothesis*}{
Name-sg = Hypothesis ,
name-sg = Hypothesis ,
Name-pl = Hypotheses ,
name-pl = Hypotheses ,
}

\newtheorem{observation}[environment]{Observation}
\AddToHook{env/observation/begin}{\zcsetup{countertype={environment=observation}}}
\zcRefTypeSetup{observation}{
Name-sg = Observation ,
name-sg = Observation ,
Name-pl = Observations ,
name-pl = Observations ,
}

\AddToHook{env/example/begin}{\zcsetup{countertype={environment=example}}}
\zcRefTypeSetup{example}{
Name-sg = Example ,
name-sg = Example ,
Name-pl = Examples ,
name-pl = Examples ,
}

\AddToHook{env/remark/begin}{\zcsetup{countertype={environment=remark}}}
\zcRefTypeSetup{remark}{
Name-sg = Remark ,
name-sg = Remark ,
Name-pl = Remarks ,
name-pl = Remarks ,
}

\zcRefTypeSetup{equation}{
Name-sg = Equation ,
name-sg = Equation ,
Name-pl = Equations ,
name-pl = Equations ,
}

\zcRefTypeSetup{chapter}{
Name-sg = Chapter ,
name-sg = Chapter ,
Name-pl = Chapters ,
name-pl = Chapters ,
}

\zcRefTypeSetup{section}{
Name-sg = Section ,
name-sg = Section ,
Name-pl = Sections ,
name-pl = Sections ,
}

\zcRefTypeSetup{algorithm}{
Name-sg = Algorithm ,
name-sg = Algorithm ,
Name-pl = Algorithms ,
name-pl = Algorithms ,
}

\AddToHook{env/notation/begin}{\zcsetup{countertype={environment=notation}}}
\zcRefTypeSetup{notation}{
Name-sg = Notation ,
name-sg = Notation ,
Name-pl = Notations ,
name-pl = Notations ,
}

\AddToHook{env/question/begin}{\zcsetup{countertype={environment=question}}}
\zcRefTypeSetup{question}{
Name-sg = Question ,
name-sg = Question ,
Name-pl = Questions ,
name-pl = Questions ,
}

\AddToHook{env/problem/begin}{\zcsetup{countertype={environment=problem}}}
\zcRefTypeSetup{problem}{
Name-sg = Problem ,
name-sg = Problem ,
Name-pl = Problems ,
name-pl = Problems ,
}

\newtheorem{definition}[environment]{Definition}
\AddToHook{env/definition/begin}{\zcsetup{countertype={environment=definition}}}
\zcRefTypeSetup{definition}{
Name-sg = Definition ,
name-sg = Definition ,
Name-pl = Definitions ,
name-pl = Definitions ,
}

\zcRefTypeSetup{figure}{
Name-sg = Figure ,
name-sg = Figure ,
Name-pl = Figures ,
name-pl = Figures ,
}

\newtheoremstyle{beautypalour}%
{3pt}%
{3pt}%
{}%
{}%
{\itshape}%
{.}%
{.5em}%
{}%
\theoremstyle{beautypalour}

\AddToHook{env/claim/begin}{\zcsetup{countertype={environment=claim}}}
\zcRefTypeSetup{claim}{
Name-sg = Claim ,
name-sg = Claim ,
Name-pl = Claims ,
name-pl = Claims ,
}

\newtheorem{beautifulclaim}{Claim}[environment]
\AddToHook{env/beautifulclaim/begin}{\zcsetup{countertype={environment=beautifulclaim}}}
\zcRefTypeSetup{beautifulclaim}{
Name-sg = Claim ,
name-sg = Claim ,
Name-pl = Claims ,
name-pl = Claims ,
}

\usetikzlibrary{calc}
\usetikzlibrary{fit}
\usetikzlibrary{decorations}
\usetikzlibrary{decorations.pathmorphing}
\usetikzlibrary{decorations.text}
\usetikzlibrary{shapes,hobby}

\tikzset{
	position/.style args={#1:#2 from #3}{
		at=($(#3)+(#1:#2)$)
	}
}

\tikzset{
  v:main/.style = {draw, circle, scale=0.8, thick,fill=black,inner sep=0.7mm},
  v:ghost/.style = {inner sep=0pt,scale=1},
  >={latex},
  e:marker/.style = {line width=8.5pt,line cap=round,opacity=0.35,color=DarkGoldenrod},
  e:main/.style = {line width=1pt},
}

\newcommand{\Coordinates}{%
    \def\R{6}

    \draw[step=.25, very thin, gray!40] (-\R,-\R) grid (\R,\R);

    \draw[step=1, thin, gray!70] (-\R,-\R) grid (\R,\R);

    \draw[->, thick] (-\R,0) -- (\R+0.35,0);
    \draw[->, thick] (0,-\R) -- (0,\R+0.35);

    \foreach \x in {-6,-5,...,6}{
        \draw (\x,0.08) -- (\x,-0.08);
        \ifnum\x=0\else
            \node[below,font=\scriptsize] at (\x,0) {\x};
        \fi
    }

    \foreach \y in {-6,-5,...,6}{
        \draw (0.08,\y) -- (-0.08,\y);
        \ifnum\y=0\else
            \node[left,font=\scriptsize] at (0,\y) {\y};
        \fi
    }

}

\newcommand{\N}{\mathbb{N}}
\newcommand{\BBB}{\mathcal{B}}

\newcommand{\LLL}{\mathcal{L}}
\newcommand{\CCC}{\mathcal{C}}
\newcommand{\SSS}{\mathcal{S}}
\newcommand{\PPP}{\mathcal{P}}
\newcommand{\TTT}{\mathcal{T}}
\newcommand{\VVV}{\mathcal{V}}

\newcommand{\Abs}[1]{|#1|}

\definecolor{BrilliantRose}{rgb}{1.0, 0.33, 0.64}
\definecolor{MidnightBlack}{rgb}{0.1,0.1,.34}
\definecolor{MidnightBlue}{rgb}{0.1,0.1,0.43}
\definecolor{Black}{rgb}{0,0, 0}
\definecolor{Blue}{rgb}{0, 0 ,1}
\definecolor{Red}{rgb}{1, 0 ,0}
\definecolor{White}{rgb}{1, 1, 1}
\definecolor{grey}{rgb}{.6, .6, .6}
\definecolor{Mygreen}{rgb}{.0, .7, .0}
\definecolor{Yellow}{rgb}{.55,.55,0}
\definecolor{Mustard}{rgb}{1.0, 0.86, 0.35}
\definecolor{applegreen}{rgb}{0.55, 0.71, 0.0}
\definecolor{darkturquoise}{rgb}{0.0, 0.81, 0.82}
\definecolor{celestialblue}{rgb}{0.29, 0.59, 0.82}
\definecolor{green_yellow}{rgb}{0.68, 1.0, 0.18}
\definecolor{crimsonglory}{rgb}{0.75, 0.0, 0.2}
\definecolor{darkmagenta}{rgb}{0.30, 0.0, 0.30}
\definecolor{magenta}{rgb}{0.50, 0.0, 0.50}
\definecolor{internationalorange}{rgb}{1.0, 0.31, 0.0}
\definecolor{darkorange}{rgb}{1.0, 0.55, 0.0}
\definecolor{ao}{rgb}{0.0, 0.5, 0.0}
\definecolor{awesome}{rgb}{1.0, 0.13, 0.32}
\definecolor{darkcyan}{rgb}{0.0, 0.50, 0.50}
\definecolor{violet}{rgb}{0.93, 0.51, 0.93}
\definecolor{brown}{rgb}{0.65, 0.16, 0.16}
\definecolor{orange}{rgb}{1.0, 0.65, 0.0}
\definecolor{cornflowerblue}{rgb}{0.39, 0.58, 0.93}

 \newcommand{\blue}[1]{{\color{Blue}#1}}
 \newcommand{\red}[1]{{\color{Red}#1}}

\newcommand{\Acal}{\mathcal{A}}

\newcommand{\Gcal}{\mathcal{G}}

\newcommand{\Lcal}{\mathcal{L}}

\newcommand{\Rcal}{\mathcal{R}}

\newcommand{\Nbbb}{\mathbb{N}}

\newcommand{\remove}[1]{}

\newcommand{\tw}{\mathsf{tw}\xspace}%
\newcommand{\poly}{\mathbf{poly}\xspace}%

\newcommand{\p}{\mathsf{p}\xspace}%

\newcommand{\bidim}{\mathsf{bidim}\xspace}%

\usepackage{xparse}

\predate{}
\date{}
\postdate{}

 \ifdefined\sth
 \else 
 \newcommand{\sth}{\mathrel : }
 \fi
\ifdefined\bd
\else
\newcommand{\bd}{\mathsf{bd}}
\fi

\newcommand{\R}{\mathbb{R}}

\newcommand{\JJJ}{\mathcal{J}}

\newcommand{\MMM}{\mathcal{M}}

\newcommand{\RRR}{\mathcal{R}}

\newcommand{\WWW}{\mathcal{W}}

 \newcommand{\QQQ}{\mathcal{Q}}
\begin{document}

\maketitle

\vspace{-.7cm}
\begin{center}
 \textbf{Abstract}
\end{center}

\begin{abstract}
\noindent The Graph Minors Series of Robertson and Seymour forms the foundation of algorithmic structural graph theory, yielding fixed-parameter algorithms for problems such as \textsc{Disjoint Paths}, \textsc{Rooted Minor Checking}, and \textsc{Folio}.
A key ingredient behind the fixed-parameter tractability of the \textsc{$k$-Disjoint Paths} problem is the irrelevant-vertex technique, which is a key mechanism in a wide range of algorithms.
The strength of this machinery is governed by the \textsl{Vital Linkage Theorem} and the so-called \textsl{Linkage Function}~$\ell$.
However, despite its foundational role, the best known bounds on the Linkage Function have remained enormous and are only implicitly understood.
The quantitative bounds behind these results have traditionally been so large that the resulting algorithms are often regarded as ``galactic.''

We show that this view is overly pessimistic and can be overturned in an essentially optimal way.
Our main result is a general irrelevant-vertex theorem for a common generalisation of \textsc{$k$-Disjoint Paths} and \textsc{Rooted Minor Checking} for graphs of size at most $d,$ commonly called the $(k,d)$-\textsc{Folio} problem.
Specifically, we show that for any graph $G$ in which the $k$ terminals are chosen from some set $R,$ if the treewidth of $G$ exceeds $\beta(k,b,d) \in 2^{\mathbf{poly}(b + d)} \cdot \mathbf{poly}(k)$ then we can locate an irrelevant vertex for the $(k,d)$-\textsc{Folio} problem.
Here, the quantity $b$ is the \textsl{bidimensionality} of $R,$ that is, the largest $b$ for which a $(b\times b)$-grid minor in $G$ can be rooted on $R$.
Thus, the exponential component of the irrelevant-vertex threshold is driven by the bound on the bidimensionality, rather than by the number of terminals, and we argue that this dependence is essentially optimal up to polynomial factors.
As a consequence, the Linkage Function satisfies $\ell(k) \in 2^{\mathbf{poly}(k)}$.
Beyond its structural significance, our result yields improved parameter dependencies for algorithms for \textsc{Disjoint Paths} and \textsc{Rooted Minor Checking}, and provides a radical quantitative improvement for a broad range of graph-minor-based algorithmic frameworks.
\end{abstract}

\medskip
\medskip
\medskip

\let\sc\itshape
\thispagestyle{empty}

\newpage

\newpage
\thispagestyle{empty}

\phantom{.} 
\tableofcontents
\thispagestyle{empty}
\newpage

\setcounter{page}{1}

\section{Introduction}\label{sec:carpetselling}
The $k$-\textsc{Disjoint Paths} and \textsc{$H$-Minor Checking} problems are among the central algorithmic problems in graph theory. In Graph Minors XIII \cite{RobertsonS1995Graph}, Robertson and Seymour showed that both \textsc{$k$-Disjoint Paths}, for fixed $k,$ and \textsc{$H$-Minor Checking}, for a fixed graph $H,$ can be solved in cubic time. The fact that these algorithmic results originate from the Graph Minors Series is not a coincidence. Rather, in this series Robertson and Seymour introduced a structural–algorithmic paradigm~--~built on tree-decompositions, excluded-grid dualities, and the irrelevant-vertex technique that has profoundly influenced modern graph algorithms, particularly in parameterized complexity \cite{FellowsL1988Nonconstructive,DemaineHK2005Algorithmic,Thilikos2012Graph,LokshtanovSZ2020Efficient}. In this sense, the significance of the Graph Minors Series extends far beyond a single routing problem. It established a general methodology for designing parameterized and structural algorithms (see \cite{LokshtanovSZ2020Efficient}).

\subsection{The algorithmics of Graph Minors}
From the outset, however, there was a striking gap between qualitative solvability and quantitative efficiency. Already in 1987, David Johnson, in ``\textsl{The {NP}-Completeness Column: An Ongoing Guide}'' \cite{Johnson1987Guide}, emphasised that while Robertson and Seymour's results were a major theoretical breakthrough, they were far from being practically applicable. He writes: \textsl{``For any instance $G=(V,E)$ that one could fit into the known universe, one would easily prefer $|V|^{70}$ to even constant time, if that constant had to be one of Robertson and Seymour's''}. 
Moreover, Johnson’s back-of-the-envelope reconstruction of just one structural constant in Graph Minors~V \cite{RobertsonS86GraphminorsV} already illustrates the scale of the problem. Unwinding his estimate leads to a quantity roughly of the form
\vspace{-1.5mm}

$$2^{\uparrow 2^{2^{2^{2^{\uparrow 2^{\uparrow \Theta(r)}}}}}}$$
\vspace{-5mm}

where $2^{\uparrow r}$ denotes a tower of $r$ twos. Thus, the difficulty is not merely a matter of implementation or low-order polynomial factors: the hidden constants arise from the structural theorems themselves. In a later retrospective, Lokshtanov, Saurabh, and Zehavi expressed as follows:
\smallskip

\begin{quote}
``\textsl{While Parameterized Complexity does provide an extremely rich toolkit for designing efficient parameterized algorithms, one of its foundations and still most powerful tools yields algorithms that are wildly impractical.}’’ \cite{LokshtanovSZ2020Efficient}.
\end{quote}
\smallskip

Indeed, the original Graph Minors dependencies are enormous, and the actual capacity for these theorems to come with reasonable bounds has largely remained opaque until recently \cite{Chuzhoy15impr,ChekuriC2016Polynomial,KawarabayashiTW2018New,KawarabayashiTW2021Quickly,ChuzhoyT2021Tighter,GorskySW2025Polynomial,GorskySW2026Price}. In particular, the constants governing the irrelevant-vertex machinery are not merely large, but also only implicitly constructible and deeply intertwined with the more arcane structural results deeper within the theory. Therefore making the Graph Minors paradigm efficient requires more than just faster implementations of existing arguments.
Instead it requires precise control and improvement of the structural thresholds that drive the entire framework. This supports Mike Fellows’ 1989 prescient appraisal: ``\textsl{it is likely to be many years before the practical significance of Robertson-Seymour theorems is fully understood}'' \cite{Fellows1989Survey}.

\paragraph{The Linkage Function.}
Given a graph $G,$ a \emph{linkage} in $G$ is a set $\Lcal$ of pairwise vertex-disjoint paths, and it is \emph{vital} if its paths span all vertices of $G$ and there is no other linkage in $G$ joining the same pairs of endpoints (see \cite{RobertsonS2009Graph}).

A central structural threshold in the Graph Minors Series is the \textsl{Linkage Function}. It stems from a theorem of Robertson and Seymour in Graph Minors~XXI \cite{RobertsonS2009Graph}, which states that every graph containing a vital linkage has treewidth bounded by a function of the number of endpoints of that linkage. This result lies at the heart of the celebrated \textsl{irrelevant-vertex technique}: Given an instance of the $k$-\textsc{Disjoint Paths} problem, when the treewidth of the graph is sufficiently large, one can delete a non-terminal vertex without affecting the outcome of the instance. By iterating this reduction, one eventually obtains an equivalent instance of bounded treewidth, that is amenable to dynamic programming. 

This paradigm lies at the heart of Robertson and Seymour's algorithm for \textsc{Disjoint Paths}, the more general folio machinery for rooted pattern detection, and \textsc{Minor Checking} as we discuss below, as well as numerous subsequent developments in structural algorithmics, parameterized complexity, and meta-algorithmics \cite{RobertsonS1995Graph,DemaineHK2005Algorithmic,RobertsonS2009Graph,RobertsonS2012Graph,Thilikos2012Graph,LokshtanovSZ2020Efficient,GolovachST2023ModelChecking,SchirrmacherSSYV24Modelchecking,FominFSST2025Compound,SauST2025Parameterizing}.

\medskip
 
The main contribution of this paper is that this core graph minors mechanism can in fact be made quantitatively reasonable. We prove an optimal estimate of the Linkage Function, up to polynomial dependencies in the relevant parameters. As a consequence, we obtain a drastically improved parameter dependence for \textsc{Disjoint Paths}, \textsc{Rooted Minor Checking}, \textsc{Folio}, and a broad range of algorithms based on the irrelevant-vertex technique. In this sense, we turn the algorithmic core of the Graph Minors Algorithm into a controlled tool, by optimally taming its ``galactic-scale’’ behaviour.

\paragraph{\textsc{Disjoint Paths} problem.}
An instance of the $k$-\textsc{Disjoint Paths} problem consists of a graph $G$ together with a \emph{pattern} $\tau=\{(s_i,t_i)\mid i\in\{1,\ldots,k\}\}$ of terminal pairs. The problem is to decide whether $G$ contains a linkage consisting of $k$ paths whose pairs of endpoints agree with~$\tau$. 

In Graph Minors XIII \cite{RobertsonS1995Graph}, Robertson and Seymour gave an algorithm for the $k$-\textsc{Disjoint Paths} problem running in $f(k)\cdot |G|^3$ time\footnote{We write $|G|$ as a shorthand for $|V(G)|$ and $|\!|G|\!|$ as a shorthand for $|E(G)|$.} for some function $f\colon \Nbbb\to\Nbbb$. At a high level, the irrelevant-vertex technique is applied as follows. The algorithm proceeds by repeatedly identifying and deleting irrelevant vertices, that is, vertices whose removal preserves equivalence of the instance. More precisely, if the input graph contains a large clique-minor, then one can find an irrelevant vertex inside a model of this clique. Otherwise, assuming that large clique-minors are excluded and that $G$ has sufficiently large treewidth, the graph contains a sufficiently large flat wall. Robertson and Seymour further proved in \cite{RobertsonS1995Graph} that the central vertices of such a flat wall are \textsl{irrelevant}. Thus removing these vertices again yields an equivalent instance. Iterating this reduction eventually produces an equivalent instance of bounded treewidth, which can then be solved by standard dynamic programming techniques.

Although the algorithm of \cite{RobertsonS1995Graph} outlined above is conceptually simple, its correctness hinges on a highly non-trivial structural claim: the central vertices of a sufficiently large flat wall are indeed irrelevant for the instance. Equivalently, every linkage that passes through such a vertex can be rerouted and replaced by one with the same pattern that avoids it. A proof of this statement only appeared later, in Graph Minors XXI \cite{RobertsonS2009Graph}, where Robertson and Seymour established it through a long and intricate sequence of arguments building on many of the structural results developed in the intervening volumes and combining their large parametric dependencies. 

We say that a linkage $\Lcal$ is a \emph{$T$-linkage} if $T$ is the set of endpoints of its paths.

\begin{proposition}[Robertson and Seymour \cite{RobertsonS2009Graph}]
\label{prop_main_RS_XXI}
There is a function $\ell\colon\Nbbb\to\Nbbb$ such that 
if a graph $G$ has a vital $T$-linkage $\Lcal,$ then the treewidth of $G$ is bounded 
by $\ell(|T|)$.
\end{proposition}

This function $\ell$ is the aforementioned \emph{Linkage Function}.
In view of the above, the algorithm of \cite{RobertsonS1995Graph} is correct and runs in $\mathbf{O}(f(k)\cdot |G|^3)$ time for some function $f(k)$ intrinsically depending on $\ell(k)$. The polynomial part of this running time was improved to quadratic in \cite{KawarabayashiKR2012Disjoint}, and more recently an almost-linear-time algorithm was given in \cite{KorhonenPS2024Minor}. As for the parametric dependence $f,$ any improvement in the dependencies of the algorithm of \cite{RobertsonS1995Graph} relies mainly on upper bounds for the Linkage Function $\ell$ in \zcref{prop_main_RS_XXI}.

\paragraph{Improving the Linkage Function}
As we already mentioned, the estimation of $\ell$ emerging from an analysis of the Graph Minors Series is colossal.
An important improvement appeared in \cite{KawarabayashiW2010Shorter},
where Kawarabayashi and Wollan gave a new proof of \zcref{prop_main_RS_XXI}
with better dependence on $k$. While they do not state any concrete 
estimate for their bound on $\ell,$ it appears to be 
no better than triple exponential in $k$ \cite{AdlerKKLST2017Irrelevant}, which remains huge, though still much better than colossal.

Another line of research has aimed to improve the dependence in \zcref{prop_main_RS_XXI} for particular graph classes. In the planar case, Adler, Kolliopoulos, Krause, Lokshtanov, Saurabh, and Thilikos proved in \cite{AdlerKKLST2017Irrelevant} that $\ell(k)\in 2^{\mathbf{O}(k)}$. Conversely, they showed in \cite{AdlerKKLST2011Tight} (see also \cite{AdlerK2019Lower}) that exponential growth of $\ell$ is \textsl{unavoidable} even for planar graphs, by giving a suitable counterexample (see \cref{fig_g_vit_lin}). Subsequently, Mazoit proved in \cite{Mazoit2013Single} that an exponential bound likewise holds for graphs embeddable in any fixed surface. His approach differs from that of \cite{AdlerKKLST2017Irrelevant} and builds on results of Geelen, Huynh, and Richter in \cite{GeelenHR2018Explicit}, who had earlier obtained bounds on $\ell$ for graphs embeddable in surfaces.

\subsection{Our results}
In this paper we prove that $\ell(k) \in 2^{\mathbf{poly}(k)}$.
In fact, we give a much more refined bound on the Linkage Function.
To state it, we need the notion of bidimensionality.

Given a graph $G$ and a set of vertices $R \subseteq V(G),$ we refer to the pair $(G,R)$ as an \emph{annotated graph}.
For a given graph $H,$ we say that $H$ is an \emph{$R$-minor of $G$} if there is a collection $\mathcal{S}=\{S_{v}\mid v\in V(H)\}$ of 
pairwise vertex-disjoint connected\footnote{A set $X\subseteq V(G)$ is \emph{connected} in~$G$ if the induced subgraph $G[X]$ is connected.} subsets of $V(G),$ each containing at least one vertex of $R,$ and such that, for every edge $xy \in E(H),$ the set $S_{x} \cup S_{y}$ is connected in~$G$.
We also work with a stronger version of annotated graphs: A \emph{$k$-rooted graph} is a pair $(G,\mathcal{R})$ where $\RRR$ is an ordered multiset of $k$ not necessarily distinct vertices in $V(G)$. Rooted graphs come with a respective notion of rooted minors; we defer the details to \cref{subsec:folio}.

The bidimensionality of $R$ in $G,$ denoted by $\textsf{bidim}(G,R),$ 
is the maximum $k$ for which $G$ contains a $(k\times k)$-grid as an $R$-minor.
Note that $\textsf{bidim}(G,R) \leq \sqrt{\Abs{R}}$.
This parameter was introduced in \cite{ThilikosW2023Graph,SauST2025Parameterizing} and intuitively captures the limited structure of a vertex set in a plane graph that lies on a bounded number of faces and translates this into a non-topological definition applicable in general annotated graphs.
Our main result is the following.
\begin{theorem}
\label{th_our_result}
There exists a function $\beta \colon \mathbb{N}^2 \to \mathbb{N}$ such that, for every annotated graph $(G,R),$ if there exists a vital $T$-linkage for some $T \subseteq R$ in $G,$ then $G$ has treewidth at most $\beta(k,b),$ where $k\coloneqq |T|$ and $b\coloneqq\textsf{bidim}(G,R)$.
Moreover, $\beta(k,b)\in \mathbf{poly}(k)\cdot 2^{\mathbf{poly}(b)}$.
\end{theorem}

The above shows that the source of the exponential behaviour of the Linkage Function 
is not the number $k$ of terminals, but rather the extent to which they are ``bidimensionally'' distributed in~$G,$ where this measure is itself bounded by $\sqrt{k}$. 
As a simple illustration of \zcref{th_our_result}, consider the following consequence:
\textsl{If $G$ is a graph embedded in the sphere and all $k$ terminals lie in $b$ faces of $G,$ then $\beta(k,b)\in \mathbf{poly}(k) \cdot 2^{\mathbf{poly}(b)}$}.
The same conclusion holds for graphs embedded in any fixed surface, even when the terminals are at a bounded distance from $b$ faces.
Since general graphs do not come equipped with faces in which terminals can be placed, \zcref{th_our_result} suggests that
bidimensionality is the right notion for expressing this phenomenon in full generality.

\paragraph{From linkages to folios.}
While the discussion above focused on the $k$-\textsc{Disjoint Paths} problem, our result holds at a more general level. Already in Graph Minors XIII \cite{RobertsonS1995Graph}, Robertson and Seymour worked in the broader framework of folios, of which $k$-\textsc{Disjoint Paths} is only a special case. For our purposes, the appropriate abstraction is the \textsc{$(k,d)$-Folio} problem.

Given a rooted graph $(G,\mathcal{R}),$ the $d$-folio of $(G,\mathcal{R}),$ denoted by $d\text{-}\mathsf{folio}(G,\mathcal{R}),$ is the set of all rooted graphs of detail at most $d$ which are rooted minors of $(G,\mathcal{R}),$ where the \emph{detail} $d$ sets an upper bound on the number of edges and non-root vertices that may be used in the rooted minors.
Notice that every rooted graph $(G,\mathcal{R})$ naturally gives rise to an annotated graph $(G,R_{\mathcal{R}})$ by simply forgetting the labels of the multiset $\mathcal{R}$ and removing duplicates.
This allows us to extend the definition of bidimensionality from annotated to rooted graphs by setting $\mathsf{bidim}(G,\mathcal{R}) \coloneqq \mathsf{bidim}(G,R_{\mathcal{R}})$.
Given an annotated graph $(G,R),$ the $(k,d)$-folio of $(G,R),$ denoted by $(k,d)\text{-}\mathsf{folio}(G,R)$ is the union of all $d$-folios taken over the $k$-rooted graphs~$(G,\mathcal{R}),$ created from any ordered multisubset $\mathcal{R}$ of $R$ such that $|\mathcal{R}| = k$.
Thus, the folio is a ``bounded-complexity summary'' of all rooted structures visible from the annotated set.

The $(k,d)$-\textsc{Folio} problem asks us to compute the $(k,d)$-folio of a given annotated graph $(G,R)$ while the \textsc{Rooted} $d$-\textsc{Folio} problem asks for the $d$-folio of a given rooted graph $(G,\mathcal{R})$.
We stress that computing the $(k,d)$-folio of an annotated graph $(G,R)$ corresponds to computing \textsl{all} $d$-folios of the rooted graphs with $k$ roots generated by selecting roots from the set $R$.

This perspective is natural for the Graph Minor Algorithm. Irrelevant-vertex arguments in fact do not merely preserve the existence of a specific routing pattern, but rather preserve the entire collection of bounded-detail rooted patterns that remain realisable after each reduction step. 

At a conceptual level, the \textsc{$(k,d)$-Folio} problem and its rooted sibling interpolate between several basic tasks. For example, \textsc{Rooted} $0$-\textsc{Folio} captures the family of feasible linkage patterns on a terminal set and, in particular, encodes the $k$-\textsc{Disjoint Paths} problem.
When the set of annotated vertices is empty, \textsc{$(0,d)$-Folio} expresses the computation of all minors of bounded size, therefore it expresses the \textsc{Minor Checking} problem, asking whether some graph $G$ contains some fixed graph $H$ as a minor.
More generally, it provides a unified language for routing, rooted minor checking, and bounded-pattern detection in graphs.
For this reason, an irrelevant-vertex theorem for folios should be viewed as a genuine strengthening of the classical Vital Linkage Theorem, and it is precisely the natural general form already suggested in Graph Minors XIII \cite{RobertsonS1995Graph}.

That the irrelevant-vertex argument applies to folios is shown in Graph Minors~XXII~\cite{RobertsonS2012Graph}.
The parameter dependence in \cite{RobertsonS2012Graph} 
remains enormous, since the result relies heavily 
on \zcref{prop_main_RS_XXI} (proved in Graph Minors~XXI \cite{RobertsonS2009Graph}) and on the implied Linkage Function.

Our main technical theorem shows that folios admit irrelevant vertices once the treewidth is large relative to an exponential function on the detail parameter and the bidimensionality of the annotated vertices and a polynomial on the number of roots that are chosen from the annotated set. Thus it lifts \zcref{th_our_result} from linkages to the full collection of bounded-detail rooted patterns.
Let $(G,R)$ be an annotated graph.
We say that a vertex $v \in V(G) \setminus R$ is \emph{strongly irrelevant} for $(k,d)\text{-}\mathsf{folio}(G,R)$ if for every ordered multiset $\mathcal{R}$ with $|\mathcal{R}| \leq k$ and $R_{\mathcal{R}} \subseteq R$ it holds that $d\text{-}\mathsf{folio}(G - v,\mathcal{R}) = d\text{-}\mathsf{folio}(G,\mathcal{R})$.

\begin{theorem} 
\label{thm:intro_folio_irrelevant_informal}
There exists a function $\beta\colon\mathbb{N}^3\to\mathbb{N}$ such that for every annotated graph $(G,R)$ with $\textsf{bidim}(G,R)\leq b,$ if $\mathsf{tw}(G)>\beta(k,b,d),$ then there exists a vertex $v\in V(G)\setminus R$ such that $v$ is strongly irrelevant for $(k,d)\text{-}\mathsf{folio}(G,R)$.
Moreover, $\beta(k,b,d)\in \mathbf{poly}(k) \cdot 2^{\mathbf{poly}(b + d)}.$
\end{theorem}

At this point, we would like to highlight the following feature of both \cref{th_our_result} and \cref{thm:intro_folio_irrelevant_informal}.
In both statements, we impose no upper bound on the size of the annotated set $R;$ we require only that its bidimensionality $b$ is bounded.
The parameter $k$ only bounds the number of the terminals in the linkage or the size of the $\Rcal$'s in the folio under consideration. Thus, $R$ should be viewed as a ``repository'' from which the $k$ terminals of $T$ or the multisets $\Rcal$ may be selected.
Our combinatorial bound depends \textsl{polynomially} on $k,$ whereas its exponential part depends only on the bidimensionality of this repository $R,$ while the size of $R$ itself may be arbitrarily large.
Another key point of \zcref{thm:intro_folio_irrelevant_informal} arises from our definition of strong irrelevancy: a strongly irrelevant vertex is naturally also irrelevant for the $d$-folio of any rooted graph $(G,\mathcal{R})$ where $\mathcal{R}$ consists of $k$ roots taken from $R$.

\cref{thm:intro_folio_irrelevant_informal} yields an algorithm for computing folios: If the graph has treewidth upper bounded by $\beta(k,b,d),$ one computes $(k,d)\text{-}\mathsf{folio}(G,R)$ via dynamic programming, in 
$2^{\mathbf{poly}(\beta(k,b,d))} \cdot |G|$ time; otherwise, one deletes an irrelevant vertex and recurses. This gives the following algorithmic application of \zcref{thm:intro_folio_irrelevant_informal}.

\begin{theorem}
\label{thm:intro_folio_algorithm_informal}
There exists an algorithm that, given a rooted graph $(G,\mathcal{R})$ with $\textsf{bidim}(G,\mathcal{R})\leq~b$ and $|\mathcal{R}| \leq k,$ computes $(d,k)\text{-}\mathsf{folio}(G,\mathcal{R})$ in 
$2^{\mathbf{poly}(k) \cdot 2^{\mathbf{poly}(b + d)}}\cdot |G|^2$ time.
\end{theorem}

As $\bidim(G,R)\leq \sqrt{|R|},$ we readily obtain the following corollary of \cref{thm:intro_folio_algorithm_informal}.

\begin{corollary}
\label{cor:other_intro_folio_algorithm_informal}
There exists an algorithm that, given a rooted graph $(G,\mathcal{R})$ where $|\mathcal{R}|\leq k,$ computes $(k,d)\text{-}\mathsf{folio}(G,\mathcal{R})$ in 
$2^{2^{\mathbf{poly}(k + d)}}\cdot |G|^2$ time. 
\end{corollary}

A nice feature of the above is that \zcref{thm:intro_folio_algorithm_informal,cor:other_intro_folio_algorithm_informal} also give rise to algorithms that compute the $(k,d)$\textsf{-folio} of $(G,R),$ when $|R| \leq k,$ in $2^{\mathbf{poly}(k) \cdot 2^{\mathbf{poly}(b + d)}}\cdot |G|^2$ time or in $2^{2^{\mathbf{poly}(k + d)}}\cdot |G|^2$ time (based on the fact that $\bidim(G,R)\leq\sqrt{|R|}$), see \cref{lem_with_r_in_it} and \cref{cor_wsithr_mor_fin_it}.
For a discussion on possible extensions of this without demanding $|R|\leq k,$ see \zcref{subsec_comp_folios_un_an}.

We emphasise an additional feature of how \cref{thm:intro_folio_irrelevant_informal} 
can be applied to prove \cref{thm:intro_folio_algorithm_informal}.
What the irrelevant vertex argument provides, under the light of \cref{thm:intro_folio_irrelevant_informal}, is a $2^{\mathbf{poly}(k)\cdot 2^{\mathbf{poly}(b + d)}} \cdot |G|^2$ time \textsl{preprocessing} algorithm that transforms an annotated graph $(G,R),$ where $\textsf{bidim}(G,R)\leq b,$ into an annotated graph $(G',R)$ where 
for every ordered multiset $\mathcal{R}$ with $|\mathcal{R}| \leq k$ and $R_{\mathcal{R}} \subseteq R$ it holds that $d\text{-}\mathsf{folio}(G',\mathcal{R}) = d\text{-}\mathsf{folio}(G,\mathcal{R})$.
Notice that this is precisely the definition of $v$ being \textsl{strongly irrelevant} for $(k,d)\text{-}\mathsf{folio}(G,R)$.
Here the annotated set $R$ is preserved and $G'$ is an induced subgraph of $G$ with treewidth in $\mathbf{poly}(k) \cdot 2^{\mathbf{poly}(b + d)} $.
Interestingly, this data reduction applies to every ordered multiset $\mathcal{R}$ that may be specified \textsl{after} the preprocessing.
Consequently, one can subsequently compute $d\text{-}\mathsf{folio}(G,\mathcal{R}),$ via dynamic programming, in $2^{\mathbf{poly}(k) \cdot 2^{\mathbf{poly}(b + d)}}\cdot |G'|$ time.
For an application of this feature, see the paragraph on ``meta-folios'' of the conclusion section, that is \zcref{subsec_ta_more_oro}.

\paragraph{Some consequences and applications.}
Our results have both structural and algorithmic consequences.
The key point is that the complexity is governed not merely by the number of terminals, but more precisely by the bidimensionality of the terminal set.
Some applications of this follow below.
\medskip

As already mentioned, the case $d=0$ of \cref{thm:intro_folio_algorithm_informal} specialises to the $k$-\textsc{Disjoint Paths} problem and yields a
$$2^{\mathbf{poly}(k)\cdot2^{\mathbf{poly}(b)}}\cdot |G|^2$$
time algorithm for $k$-\textsc{Disjoint Paths}, where $b = \mathsf{bidim}(G,T)\leq \sqrt{k}$ for the terminal set $T$.
Even in this special case, no better parametric dependence is known beyond the galactic bounds stemming from the Graph Minors Series and the still enormous~--~and not explicitly specified~--~bounds arising from the Linkage Function estimates in \cite{KawarabayashiW2010Shorter}.
For planar graphs, single-exponential bounds were obtained recently in \cite{LokshtanovMPSZ2020Exponential,ChoOO2023Parameterized} and kernelization algorithms have been given in \cite{WlodarczykZ2023Planar}.
As we discuss below, our bounds for the Linkage Function in \cref{th_our_result} are asymptotically optimal.
This implies that, in order to improve the presented running time for the $k$-\textsc{Disjoint Paths} problem, one would need to expand on the techniques from \cite{LokshtanovMPSZ2020Exponential,ChoOO2023Parameterized} which avoid explicit dynamic programming on tree-decompositions.

Another natural specialisation is \textsc{Rooted Minor Checking}.
This problem asks, given a rooted graph $(G,\mathcal{R})$ and a rooted graph $(H,\mathcal{T}),$ whether $(G,\mathcal{R})$ contains $(H,\mathcal{T})$ as a rooted minor.
Note that this is just a membership query to the folio of $(G,\mathcal{R})$.
Therefore, \zcref{thm:intro_folio_algorithm_informal} yields an algorithm for \textsc{Rooted Minor Checking} on general graphs with running time
$$2^{\mathbf{poly}(k) \cdot 2^{\mathbf{poly}(b + d)}}\cdot |G|^2,$$
where $b \coloneqq \mathsf{bidim}(G,\mathcal{R})$ and $d$ is the detail of $(H,\TTT)$.
This yields, for the first time, an algorithm for rooted minor checking with reasonable parametric dependencies. The special case $k=0$ gives a
$$2^{2^{\mathbf{poly}(|H|)}}\cdot |G|^2$$
time algorithm for deciding whether a graph $G$ contains a graph $H$ as a minor. We again emphasise that, to the best of our knowledge, this marks the first time an explicit running time for this problem can be stated.
In our view, this goes a long way toward dispelling the reputation of graph minors algorithms as being ``galactic-scale''.
\medskip

We now consider a refinement of the main meta-algorithmic consequence of the Graph Minors Series.
A \emph{graph parameter} is a function $\p$ mapping graphs to non-negative integers and it is \emph{minor-monotone} if its value on a graph is no smaller than its value on each of its minors.
Given such a parameter $\p,$ we define $\mathsf{Obs}(\p,k)$ as the set of all minor-minimal graphs $G$ with~$\p(G)>k$.

The seminal result of the Graph Minors Series is that graphs are well-quasi-ordered under the minor relation \cite{RobertsonS86GraphminorsXX}.
This implies that, for every minor-monotone parameter $\p,$ there exists a function $f_{\p}\colon \Nbbb\to\Nbbb$ such that for every $k\in\Nbbb,$ every graph in $\mathsf{Obs}(\p,k)$ has size at most $f_{\p}(k)$.
Since the proof that graphs are well-quasi-ordered by minors is non-constructive (see \cite{FriedmanRS87them,KrombholzR2020Upper}), there is no a priori general mechanism for constructing the function $f_{\p}$ from a description of the parameter~$\p$. 

In \cite{FellowsL1988Nonconstructive}, Fellows and Langston revealed the meta-algorithmic potential of the Graph Minors Series (along with its constructibility limitations).
We refine this line of reasoning as follows.

Suppose that, for some particular $\p,$ we are in position to \textsl{know}
some fixed upper bound on $f_{\p}(k)$
and, moreover, $\p$ is exponential-time computable, i.\@e.\@, there is an algorithm computing $\p(G)$ in $2^{\poly(|G|)}$ time.
We directly have an algorithm to compute $\mathsf{Obs}(\p,k)$ in $2^{\textbf{poly}(f_{\p}(k))},$ 
since there are $2^{O\big((f_{\p}(k))^2\big)}$ graphs in $\mathsf{Obs}(\p,k)$ each of size at most $f_{\p}(k)$.
By the definition of $\mathsf{Obs}(\p,k),$ for every $k$ and every graph $G,$ one can decide whether $\p(G)\leq k$ by checking whether $G$ contains some graph in $\mathsf{Obs}(\p,k)$ as a minor. 
Therefore, the check whether $G$ contains some graph in $\mathsf{Obs}(\p,k)$ can be done in time $2^{O\big((f_{\p}(k))^2\big)}\cdot 2^{2^{\mathbf{poly}(f_{\p}(k))}}\cdot |G|^2,$ by our results. We obtain the following.

\begin{theorem}\label{prop_rs_meta}
For every exponential-time computable minor-monotone parameter $\p,$ there is an algorithm that, given a graph $G,$ decides whether $\p(G)\leq k$ in time $ 2^{2^{\mathbf{poly}(f_{\p}(k))}}\cdot |G|^2.$
\end{theorem}

We view \cref{prop_rs_meta} as a refinement of the main algorithmic consequence of the fact that graphs are well-quasi-ordered by the minor relation, proved in \cite{RobertsonS86GraphminorsXX}.
\medskip

The significance of the above statements extends beyond the two endpoint problems discussed above, i.\@e.\@, \textsc{Disjoint Paths} and \textsc{Minor Checking}.
Folios are reusable summaries of bounded-detail rooted structure.
Consequently, any algorithm that repeatedly asks whether one of a bounded number of rooted patterns occurs in a graph can benefit from faster folio computation.
This includes rooted routing tasks and algorithms that use bounded rooted minors as certificates or obstructions.
More generally, every algorithmic framework whose reduction phase proceeds by finding an irrelevant vertex in a big clique-minor or in a big flat wall inherits its quantitative behaviour from the same threshold.
Prominent examples of such results whose {(meta-)parametric} dependencies can be improved using our bounds include finding embeddings in $\mathbb{R}^3$ \cite{KawarabayashiKM10link}, modification problems \cite{FominLPSZ2019Hitting,SauST2022kApices2,BasteST2023Hitting,SauST2023kApices,MorelleSST2024Faster}, as well as combinatorial upper bounds to minor obstructions \cite{SauST2022kApices2,PaulPTW2024Obstructionsa}.

\paragraph{Optimality of our bounds.}
We also show that our upper bounds are essentially optimal.
The lower bounds of \cite{AdlerKKLST2011Tight,AdlerK2019Lower} already show that the exponential dependence in \zcref{th_our_result} is unavoidable, even in planar graphs.
Moreover, as we show in \cref{sec:lowerbounds}, this phenomenon is not merely a consequence of having many terminals.
The same construction also yields an exponential lower bound when the parameter is expressed in terms of bidimensionality.
More precisely, if \zcref{th_our_result} holds for some threshold function $\beta,$ then necessarily $\beta(k,b)\in2^{\Omega(\sqrt{b})}$.

We prove a similar optimality statement for \zcref{thm:intro_folio_irrelevant_informal}, indicating that the contribution of the detail $d$ is essentially optimal. In this case, the argument relies on enriching the counterexample of \cite{AdlerKKLST2011Tight,AdlerK2019Lower} with many pairwise non-isomorphic bounded-degree gadgets of very small size, obtained via asymptotic enumeration of regular graphs. This implies that any threshold function $\beta$ as the one in \zcref{thm:intro_folio_irrelevant_informal} must satisfy
$\beta(k,b,d)\in 2^{\Omega\left(\frac{d\log\log d}{\log d}\right)}$.
 
\subsection{Proof overview}
There are two major steps towards our proof of \cref{thm:intro_folio_algorithm_informal}.
At its heart lies the general bound on the Linkage Function and thus the proof of \zcref{th_our_result}, which on its own spans several sections.
Our ultimate goal for this part can be intuitively understood as follows:
if there exists an annotated graph $(G,R)$ admitting a vital $R$-linkage $\LLL$ with $b \coloneqq \mathsf{bidim}(G,R),$ then there exists an annotated graph $(G',R')$ admitting a vital $R'$-linkage $\LLL',$ such that $\mathsf{bidim}(G',R') \in \mathbf{poly}(\Abs{\mathcal{L}})\cdot2^{\mathbf{poly}(b)}$ and $G'$ can be embedded in some surface $\Sigma$ of genus $g(b)$ admitting $h(b)$ holes with $g(b),c(b) \in \mathbf{poly}(b)$ such that the terminals $T'$ are embedded on the boundary components of $\Sigma$.
Given such an instance $(G',R')$ together with $\LLL',$ we show that the treewidth of $(G',R')$ must be bounded in $\mathbf{poly}(\Abs{R'} + b'),$ where $b' \coloneqq \mathsf{bidim}(G',R')$.
As a consequence of the above we have $b' \in\mathbf{poly}(\Abs{R})\cdot2^{\mathbf{poly}(b)},$ concluding the proof of \cref{th_our_result}.

In a second step we leverage \cref{th_our_result} to derive an algorithmic version of \zcref{thm:intro_folio_irrelevant_informal}.
This involves another two-step process.
First, we make use of \cref{th_our_result} directly to prove an analogue of the Annulus Combing Lemma of Golovach, Stamoulis, and Thilikos \cite{GolovachST2023Combing}, with better bounds, and respecting the bidimensionality of the terminals.
Later, we leverage this lemma to show that minor models of partially disc-embedded graphs may be redrawn inside a large flat wall~--~given certain homogeneity assumptions.
This hues close to a lemma by Baste, Sau, and Thilikos \cite{BasteST2023Hitting}.
However, their result focussed solely on the detail and does not analyse the role of the terminals, requiring more work from our side again.
Finally, our result may then be combined with a recent homogenisation technique by Gorsky, Seweryn, and Wiederrecht \cite{GorskySW2026Price} in order to complete the proof of \zcref{thm:intro_folio_irrelevant_informal}.
\smallskip

We start with a sketch of the proof of \cref{th_our_result}.
Let $(G,R),$ $\LLL,$ and $b$ be as above.
Assume, towards a contradiction, that $G$ has larger treewidth than $\beta_{\ref{th_our_result}}(k,b)$ where~$k \coloneqq |\mathcal{L}|$.

As a first step, we use standard machinery of Robertson and Seymour \cite{RobertsonS1995Graph} together with recent results of Protopapas, Thilikos, and Wiederrecht \cite{ProtopapasTW2025Colorful} on $R$-minors to prove that $G$ cannot admit a clique of order $\mathbf{poly}(b)$ as a minor, as otherwise $\LLL$ was not vital. 

\paragraph{Hunting down the apex.}
Applying the main result of \cite{GorskyPW2026Quickly} (see \cref{thm:localstructure}), we use the absence of a clique to derive the existence of a surface $\Sigma$ of genus at most $\mathbf{poly}(b),$ an \emph{apex set} $A$ of size $\Abs{A} \in \mathbf{poly}(b)$ such that the graph $(G-A,R\setminus A)$ admits a so-called ``almost embedding'' with $\leq \mathbf{poly}(b)$ regions~--~called the ``vortices''~--~that exhibit non-planar behaviour~--~called \emph{depth}~--~controlled in $\mathsf{poly}(b)$ and a ``flat'' (think of it as having ``planar behaviour'') $\beta(k,b)$-wall $W,$ such that all the vertices in $R\setminus A$ are part of vortices.
Note that since $\LLL$ was vital for $G,$ there is a natural vital linkage $\LLL'$ in $G' \coloneqq G - A$ obtained from $\LLL$ by deleting $A,$ where each vertex in $A$ may result in two new terminals in $\LLL'$.
Thus, $\LLL'$ is a vital $R'$-linkage in $G-A$ where $\Abs{R'} \leq \Abs{R} + 2\Abs{A}$.
Further $b' \coloneqq \mathsf{bidim}(G - A, R') \in \mathbf{O}(b)$.
However, note that the newly introduced terminals may not be part of vortices and it is a priori unclear where they are located with respect to the wall $W$.

Thus we start by finding a large subwall $W'$ of $W$ disjoint from these terminals and declaring each spot in which such a new terminal sits to be a new vortex, introducing at most $\mathbf{poly}(b)$ new vortices.
We then need to find disjoint sets of $\mathbf{poly}(b)$ concentric cycles around our vortices while maintaining a sizeable subwall of $W'$ disjoint from these cycles.
Here we build on ideas of Diestel, Kawarabayashi, M{\"u}ller, and Wollan \cite{DiestelKMW2012Excluded}.
However, their techniques may introduce new apices outside of the existing vortices, further creating even more vortices, and thus threatening a loop.
We therefore need to significantly alter their approach:
we accomplish this by directly manipulating the surface and the almost embedding instead of deleting parts of the graph, while maintaining our bounds and the vitality of the linkage.

\paragraph{Taming vortices.}
Given the set of $\mathbf{poly}(b)$ vortices we need to ``tame'' the behaviour of $G$ in these regions.
Here we build on ideas pioneered by Huynh \cite{Huynh2009Linkage} and refined by Diestel et al.\ \cite{DiestelKMW2012Excluded}, who take a vortex of depth at most $d$ and $(d+1)$ concentric cycles tightly around it and ``push'' the cycles into the vortex as deeply as possible, where the bounded depth guarantees that at least one cycle will remain completely outside of that region.
One then extends the vortex to the disc bounded by said cycle.

Unfortunately, we cannot use these results directly, since they again introduce apices that we cannot afford and require that vortices come with $(d+1)$ ``private'' concentric cycles, such that no vortex is part of the disc bounded by a cycle of another vortex.
We cannot guarantee this easily without introducing an exponential dependency on our parameters.
Fortunately we do not need too.
Expanding on the techniques in their proof, we show that if every vortex~--~call them $c_1,\ldots,c_m$~--~comes with a family $\CCC_i$ of $\mathbf{poly}(d)$ many concentric cycles such that either all of the cycles $\CCC_i$ lie in the disc bounded by some cycle of $\CCC_j$ or are completely disjoint from it, then an ``apex-less'' variation of the above procedure still works.

The main idea is that, if an untamed vortex $c_1$ together with its concentric cycles $\CCC_1$ is completely contained in the disc of some cycle of $\CCC_2$ say, then when we tame $c_2$ by pushing the $(d+1)$ closest cycles in $\CCC_2$ into $c_2,$ the resulting set of tight cycles will never interfere with most of the cycles in $\CCC_1,$ as otherwise we could deflect them resulting in tighter cycles and in particular, the cycles will never penetrate $c_1$.
Thus the vortex $c_1$ will be either completely absorbed by the vortex $c_2$ or, after taming $c_2,$ it will still be outside of the vortex and have most of the cycles of $\CCC_1$ left in order to tame $c_1$ in the future.
That way, we can tame our vortices by guaranteeing only $\mathbf{poly}(d) \in \mathbf{poly}(b)$ cycles around each of them, which the previous step of the proof provides.

Note that after linking the vortices we are still embedded with up to $\mathbf{poly}(b)$ many vortices, although some vortices might have been swallowed by others along the way.

\paragraph{No apex, no cry.}
We are left with an instance $(G,T),$ $\LLL$ embedded on a surface $\Sigma$ of genus $g(b)$ with $h(b)$ holes, at most $\mathbf{poly}(b)$ vortices, all of them tame, and a set of $\beta(k,b)$ many concentric cycles $\CCC$ coming from the wall $W$.
Building on ideas of \cite{RobertsonS2009Graph,CavallaroKK2024EdgeDisjoint} we study the structure of vortices in the presence of vital linkages, proving that the vital linkage interacts with vortices in a very restrictive way.
Robertson and Seymour \cite{RobertsonS2009Graph} essentially prove that either a tamed vortex has a short boundary or most of the vital linkage must be contained in the vortex.
However, the main technical lemma of their proof already produces a doubly exponential bound on the required length of the vortex that is involved.
Thus we need to avoid this.

Denote by $\SSS(\LLL)$ the set of all the maximal subpaths of $\LLL$ with endpoints in vortices that are otherwise disjoint from vortices.
We then prove that given a disc $\Delta(S_1,S_2)\subseteq \Sigma,$ where $S_1,S_2$ are disjoint regions on the boundary of some vortices (possibly the same), such that $\Delta(S_1,S_2)$ is otherwise disjoint from all vortices, $\Delta(S_1,S_2)$ cannot contain more than $\mathbf{poly}(k)\cdot2^{\mathbf{poly}(d)}$ paths in $\SSS(\LLL)$ with one endpoint in $S_1$ and the other in $S_2,$ where $d\in \mathbf{poly}(b)$ is the maximal depth of the vortices.
Hence the number of homotopic paths of $\SSS(\LLL)$ that are not null-homotopic~--~they do not bound a disc~--~is bounded in $\mathbf{poly}(k)\cdot 2^{\mathbf{poly}(b)}$.
Unfortunately, we cannot bound the number of null-homotopic paths~--~it may be more than single-exponential in $b$~--~and the respective techniques in \cite{RobertsonS2009Graph} are of no direct use to us.

Fortunately, it turns out that we do not need to bound the number of these paths.
It suffices to prove that every null-homotopic path does not interact with most of the cycles in $\CCC$.
We are able to guarantee this by combining cheap linkage techniques in the spirit of \cite{AdlerKKLST2011Tight,CavallaroKK2024EdgeDisjoint} in order to gather further structural insight on the interaction of $\mathcal{S}(\LLL)$ with $\CCC$.
Assuming one of these paths reaches too deeply into $\CCC,$ we prove the existence of a set of null-homotopic paths in $\mathcal{S}(\mathcal{L})$ forming a large linkage with one endpoint in $S_1$ and the other in $S_2$ for two disjoint regions on some common vortex.
Then we can again apply the result above, proving that $\Delta(S_1,S_2)$ cannot contain more than $\mathbf{poly}(k)\cdot2^{\mathbf{poly}(b)}$ paths in $\SSS(\LLL)$ with one endpoint in $S_1$ and the other in $S_2$.
In particular, at most $\mathbf{poly}(k)\cdot2^{\mathbf{poly}(b)}$ many circles of $\CCC$ were intersected.
This in turn implies that a large subset $\mathcal{C}'$ of the $\beta(k,b)$ cycles~--~for an appropriate choice as will be discussed in the next paragraph~--~have not been intersected.

We then delete all of the vertices that are part of vortices, all simple loops, and the edges of all cycles outside $\CCC'$.
This produces an instance $(G',T')$ embedded in $\Sigma$ without vortices, admitting $\CCC'$ concentric cycles and a vital $T'$-linkage $\LLL',$ such that for each homotopy type of paths in $\LLL$ there exist at most $\mathbf{poly}(k)\cdot 2^{\mathbf{poly}(b)}$ many paths, none of which are simple loops.

Finally, it is an easy combinatorial exercise to prove that a surface $\Sigma$ with genus $g(b)\in \mathbf{poly(b)}$ and $h(b)\in \mathbf{poly}(b)$ holes admits at most $\mathbf{poly}(b)$ paths that are pairwise not homotopic.
Thus,~$\LLL'$ admits at most $\mathbf{poly}(k)\cdot 2^{\mathbf{poly}(b)}$ paths and $\Abs{T'} \leq \mathbf{poly}(k)\cdot 2^{\mathbf{poly}(b)}$.
Note that this potentially results in $2^{\mathbf{poly}(b)}$ new terminals.
However, the bidimensionality of $(G',T')$ is still bounded in $\mathbf{poly}(g(b)+h(b))$ and hence in $\mathbf{poly}(b)$.

\paragraph{The bounded genus case.}
We now have an annotated graph $(G,R),$ with $G$ embedded on a surface $\Sigma$ of genus at most $g = g(b) \in \mathbf{poly}(b)$ and there exist at most $h = h(b) \in \mathbf{poly(b)}$ holes in $\Sigma$ such that all vertices of $R$ are drawn on the boundary of those holes.

We aim to show that, if $G$ contains a family $\mathcal{C}$ of $\lambda(g,h,k) \in \mathbf{poly}(k)\cdot 2^{\mathbf{poly}(g + h)}$ concentric cycles embedded in a disc $\Delta \subseteq \Sigma$ and $\mathcal{L}$ is an $R$-linkage of order at most $ \mathbf{poly}(k)\cdot2^{\mathbf{poly}(b)} ,$ then $\mathcal{L}$ avoids the innermost cycle from $\mathcal{C}$.
This implies that such a linkage $\mathcal{L}$ cannot be vital.
Notice that here the total number of terminals, while polynomial in $k,$ depends exponentially on $b$.

In order to handle this situation and to ensure that all remaining exponential dependencies remain polynomial in this already large number of terminals, we revisit an unpublished proof by Mazoit \cite{Mazoit2013Single} from 2013.
In his preprint, Mazoit gave a sketch of a proof establishing single-exponential bounds for the Linkage Function for graphs of Euler-genus at most $g$.
While the core ideas presented in \cite{Mazoit2013Single} turn out to be correct and crucial for our goal, the paper contains several flaws and has~--~at the time of writing~--~only appeared on \texttt{arxiv.\@org}.
For the sake of ensuring correctness and completeness, we include a full, detailed, and reworked proof of Mazoit's core result.
In particular, we extract and highlight a crucial dichotomy in the dependencies on $h,$ $g,$ and $k$ from Mazoit's ideas as indicated above:
\begin{center}
 \textsl{The exponential dependency of the function $\lambda$ is entirely confined to $h$ and $g$ while the dependency on the number $k$ of terminals remains polynomial in $k$.}
\end{center}
The main idea of the proof is to make use of conditions on when a linkage that is feasible in the topological sense~--~meaning that it could be drawn on the surface, ignoring the underlying graph itself~--~can actually be realised in the embedded graph.
The proof reduces to three core pieces:
(1) The case when $h=g=1,$ that is $\Sigma$ is a disc, requires only that $\lambda(1,1,k) \in \mathbf{O}(k)$.
(2) There exists a theorem by Geelen, Huynh, and Richter \cite{GeelenHR2018Explicit} saying that for every non-separating curve $\gamma$ with both ends on the boundary of some surface, any path of any linkage drawn in the surface can be redrawn to cross $\gamma$ only twice.
(3) One may interpret the surface $\Sigma$ as consisting entirely of the disc $\Delta$ together with a number $s \in \mathbf{O}(g + h)$ many strips attached to $\Delta$.
Using (1) as the base case, the claim now follows by induction on the number $s$ of strips, guaranteed by (3), where (2) acts as the tool to replace a single strip with a bounded number of additional terminals by sacrificing a corresponding amount of the concentric cycles. Combining the steps described above concludes the proof sketch for \cref{th_our_result}; we continue with \cref{thm:intro_folio_irrelevant_informal}.

\paragraph{Combing an annulus.}
We say that a graph $G$ is \emph{partially disc-embedded} if there exists a disc $\Delta$ and a triple $(G_1,G_2,\Gamma)$ where $G = G_1 \cup G_2,$ $G_1$ and $G_2$ intersect in vertices only, and $\Gamma$ is an embedding of $G_1$ in $\Delta$ such that precisely the vertices in $V(G_1) \cap V(G_2)$ are drawn on the boundary of $\Delta$.
Golovach, Stamoulis, and Thilikos \cite{GolovachST2023Combing} proved a useful technical lemma.
Let $G$ be partially disc-embedded as above such that there exists a family $\mathcal{C}$ of $a'(k)$ concentric cycles crossed by $a'(k)$ paths $\mathcal{R}$~--~called the \emph{rails}~--~and let $\circledcirc \subseteq \Delta$ denote the annulus between the innermost and the outermost cycle of $\mathcal{C}$.
Then, if $G$ contains a linkage $\mathcal{L}$ of order at most $k$ such that no terminal of $\mathcal{L}$ lies in $\circledcirc,$ then there also exists a linkage $\mathcal{L}'$ in $G$ of the same pattern as $\mathcal{L}$ such that $\mathcal{L}'$ intersects $\circledcirc$ only in $\mathcal{R}$.
This lemma is known as the \emph{Annulus Combing Lemma}.
The original version of this lemma had no additional assumptions on the terminals of $\mathcal{L}$ and had a quadratic dependency on the Linkage Function.
We provide an independent proof using \zcref{th_our_result} while introducing the additional measure of bidimensionality on the terminals and improving the dependencies to being linear in the Linkage Function.
Notably, our proof of the Annulus Combing Lemma is the only place where we directly use \zcref{th_our_result}.
Once this lemma is proven, all further algorithmic arguments build on it instead of \zcref{th_our_result}.

While this might already be of independent interest, we then proceed by proving a further strengthening of a similar combing lemma from the work of Baste, Sau, and Thilikos \cite{BasteST2023Hitting}.
This lemma allows us, given that the partial disc-embedding is actually an almost embedding and there exists a large wall $W$ drawn in the interior of the disc~--~given that $W$ satisfies a condition known as \textsl{homogeneity}~--~to redraw any rooted minor model of bounded detail in a way such that it avoids the centre of $W$ entirely.
As before, we obtain much better dependencies than the original result in order to prove a version sensitive to the bidimensionality of the roots to adapt the lemma to our setting.
As a side effect, our results in this part show that the highly technical definitions originating from recent controversies about technical details regarding the Flat Wall Theorem (see for example \cite{Arnon2023fixing,SauST2024More}) may be completely avoided and simplified.

\paragraph{Who cares about this vertex?}
At this point, most major pieces towards the proof of \zcref{thm:intro_folio_irrelevant_informal} are in place.
Indeed, we actually prove an algorithmic version of \zcref{thm:intro_folio_irrelevant_informal} in \zcref{sec:irrelevantvertex} and all of our algorithmic main results are derived from it in \zcref{sec:GMAexplanation}.
There are only two remaining pieces to be discussed:
Firstly, we require a way to obtain a clique-minor, a small width tree-decomposition, or a large flat wall $W$ in linear time in order to guarantee the quadratic running time of \zcref{thm:intro_folio_algorithm_informal}.
Moreover, we require this flat wall $W$ to satisfy some additional properties.
One of them is a small connectivity property, and the other is the previously mentioned ``homogeneity''.
To find one of the three required outcomes while also ensuring that this connectivity requirement is met by our wall, we employ a recent result due to Sau, Stamoulis, and Thilikos \cite{SauST2024More}.
Homogeneity is slightly more tricky.
It has been observed in the literature that \cite{SauST2020An,SauST2022kApices2,SauST2023kApices,MorelleSST2024Faster} the original definitions for homogeneity and the original proofs by Robertson and Seymour seem to require an additional exponential blow-up in order to ensure $W$ to be homogeneous.
However, within our proofs we are able to utilise the recent homogeneity framework, due to Gorsky, Seweryn, and Wiederrecht \cite{GorskySW2026Price}, which ensures that the price of homogeneity remains polynomial. 
This last piece of technology allows us to avoid the very last source of possible exponential explosion and, by combining these theorems with our strengthened Annulus Combing Lemma for rooted minors, we are able to show that no rooted minor of small detail in $(G,R)$ cares about the vertices in the very centre of $W$.

\subsection{Organisation}
After the introduction, we start with a more in-depth explanation of the Graph Minor Algorithm in \zcref{sec:GMAexplanation}, which provides more careful definitions and explanations for the variants of the \textsc{Folio}-problem we consider and also illustrates how this allows us to resolve $k$-\textsc{Disjoint Paths} and \textsc{$H$-Minor Checking}.
Once this is taken care of, we provide necessary definitions for the rest of the work (some of which are covered in \zcref{sec:GMAexplanation}) in \zcref{sec:prelims}.

The first proper step into the proofs is taken in \zcref{sec:linkagesinsurfaces}, which starts by introducing tools for a common setting we will encounter throughout our work~--~namely a set of concentric cycles in a disc with a linkage running through it~--~and then moves on to providing a fresh presentation and correction of the crucial core of Mazoit's perspective on vital linkages in graphs embedded on a surface of bounded genus.
This is followed by a pair of sections which shows how to reduce the case in which we are almost embedded on a surface without apices by first showing how to tame vortices in \zcref{sec:linkedness} and then doing the arduous technical core of the proof of \zcref{th_our_result} in \zcref{sec:apexfree} by reducing the apex-free case to graphs embedded on surfaces of bounded genus.
\zcref{sec:noapexnocry} then takes care of getting rid of apices in an almost embedding and building the prerequisites to apply the results of \zcref{sec:linkedness,sec:apexfree}, while \zcref{sec:largeclique} presents a small toolkit to deal with the presence of large clique-minors in our graph.
The proof of our variant of the Vital Linkage Theorem (see \zcref{th_our_result}) is then laid out by combining the results of the previous sections in \zcref{sec:VitalLinkage}.

The second half of our paper begins with \zcref{sec:folio}, in which we provide all the technical tools we need for the \textsc{Folio}-problem that can be derived from \zcref{th_our_result}, including the Annulus Combing Lemma and its variant for rooted minors.
These results are then used to prove the existence of irrelevant vertices for the \textsc{Folio}-problem within graphs of high treewidth in \zcref{sec:irrelevantvertex}, which is capped by showing that this existence-result can be used algorithmically to find such a vertex in linear time.

We close our work by first giving proofs for several lower bounds to our results in \zcref{sec:lowerbounds}, which demonstrate that what we prove is essentially optimal in several ways.
This is then contrasted by some comments and suggestions for future research directions in \zcref{sec:conclusion} to explain that there is indeed still a lot of potential for further developments in this area.

\section{The Graph Minor Algorithm}\label{sec:GMAexplanation}
In this section we provide proofs for our algorithmic main results, namely \zcref{thm:intro_folio_irrelevant_informal}, \zcref{thm:intro_folio_algorithm_informal}, and~--~as a direct consequence~--~\zcref{prop_rs_meta}.
To do this, we provide a small amount of formal definitions and references to later points in the paper for those definitions that are too bulky and technical to be explained at this point.
The proofs of our main theorems are essentially reductions to more powerful and formal variants, namely \zcref{thm:Main1Folio} and \zcref{thm:Main2Folio} whose proof spans mostly \zcref{sec:folio} and \zcref{sec:irrelevantvertex}.

The proof of our most central result, that is the single-exponential bound on the Linkage Function conditioned on the bidimensionality of the annotated vertex set as manifested in the form of \zcref{th_our_result} is confined to the proof of \zcref{thm:VitalLinkage} in \zcref{sec:VitalLinkage} and spans all of \zcref{sec:linkagesinsurfaces,sec:linkedness,sec:apexfree,sec:noapexnocry,sec:largeclique,sec:VitalLinkage}.

We begin by introducing some general notation and some of our core concepts.

\subsection{A small set of preliminary definitions}\label{subsec:smallPrelims}
First, we introduce notation for commonly used concepts.
By $\mathbb{N}$ we denote the set of non-negative integers.
Given any two integers $a,b\in\mathbb{N},$ we write $[a,b]$ for the set $\{z\in\mathbb{N} ~\!\colon\!~ a\leq z\leq b\}.$
Notice that the set $[a,b]$ is empty whenever $a>b.$
For any positive integer $c$ we set $[c]\coloneqq [1,c].$
In order to avoid ambiguity for this notation, we denote closed intervals over the reals by $[x,y]_{\mathbb{R}}$.

In several contexts, we will have a set $\mathcal{G}$ of graphs and want to access the graph their union comprises or want to interact with its edges or vertices.
We use $\bigcup \mathcal{G}$ as for $\bigcup_{G \in \mathcal{G}} G$ and may thus access the vertices and edges as $V(\mathcal{G}) \coloneqq V(\bigcup \mathcal{G})$ and $E(\mathcal{G}) \coloneqq E(\bigcup \mathcal{G})$.

\paragraph{Annotated graphs and minors.}
A \emph{minor-model} $\varphi$ of a graph $H$ in a graph $G$ is a function mapping $V(H)$ to a collection of pairwise disjoint connected subsets of $V(G),$ called the \emph{branch sets}, such that for all $uv\in E(H)$ there is an edge in $G$ between $\varphi(u)$ and $\varphi(v)$.

An \emph{annotated graph} is a pair $(G,R)$ where $G$ is a graph and $R\subseteq V(G)$ is the set of \emph{annotated vertices}, which we often simply refer to as \emph{red vertices}.
If we need easy access to the number of vertices in $R,$ we say that $(G,R)$ is an \emph{$|R|$-annotated graph}.
The \emph{bidimensionality} of an annotated graph $(G,R)$ is the largest integer $k$ such that there exists a minor-model $\varphi$ of the $(k \times k)$-grid where $R \cap \varphi(x) \neq \emptyset$ for all vertices of the grid.
In general, we say that an annotated graph $(H,R_H)$ is a \emph{red minor} of $(H,R_G)$ if there exists a minor-model $\varphi$ of $H$ in $G$ such that $R_G \cap \varphi(v) \neq \emptyset$ for all $v\in R_H$.

We note that, in an annotated graph $(G,R),$ declaring a set $S \subseteq V(G)$ to also be red increases the bidimensionality by no more than $|S|$.

\begin{observation}\label{obs:AddRedVertices}
Let $(G,R)$ be an annotated graph and $S \subseteq V(G),$ then 
\[ \mathsf{bidim}(G,R \cup S) \leq \mathsf{bidim}(G,R) + |S| . \]
\end{observation}

\begin{figure}[ht]
 \centering
 \includegraphics[width=9cm]{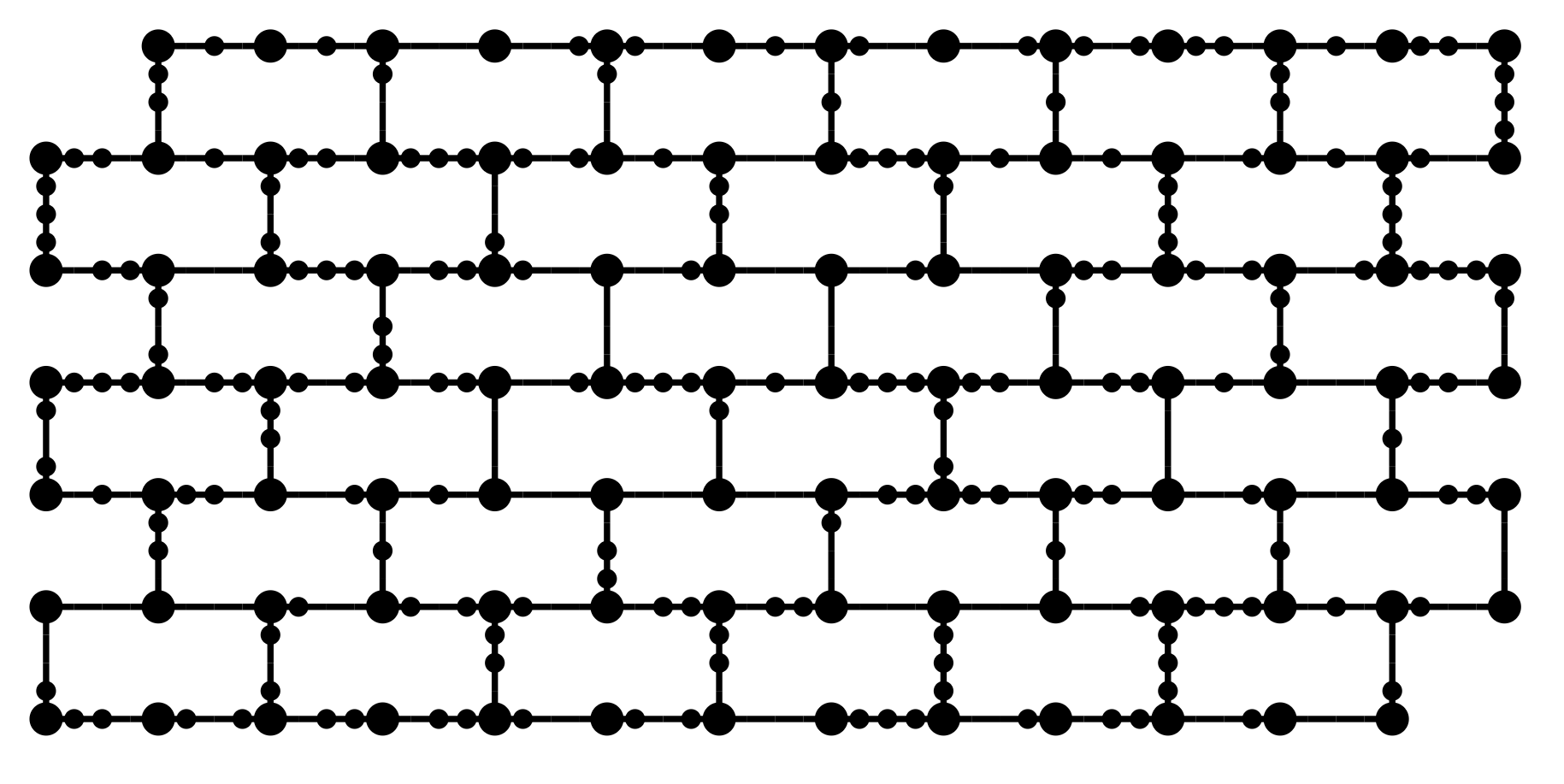}

 \caption{A $7$-wall.}
 \label{fig:wall}
\end{figure}

\paragraph{Grids and Walls.}
Let $n,m\geq 1$ be integers.
The \emph{$(n \times m)$-grid} is the graph with vertex set $[n] \times [m]$ where $\{(i_1,j_1),(i_2,j_2)\}$ is an edge if and only if $|i_1 - i_2| + |j_1 - j_2| =1$.
Given an integer $n,$ the \emph{elementary $n$-wall} is the graph obtained from the $(n \times 2n)$-grid by deleting all edges of the form
\begin{itemize}
 \item $\{(i,j),(i+1,j) \}$ for $i\in[n-1]$ odd and $j\in[2n]$ even, and
 \item $\{(i,j),(i+1,j) \}$ for $i\in[n-1]$ even and $j\in[2n]$ odd.
\end{itemize}
See \zcref{fig:wall} for an illustration.
An \emph{$n$-wall} is any graph that can be obtained from the elementary $n$-wall by subdividing its edges.

One can easily see that any graph that contains the $(n\times 2n)$-grid as a minor must contain an $n$-wall as a subgraph.
Moreover, every $n$-wall contains the $(n\times n)$-grid as a minor.

\paragraph{Treewidth and tree-decompositions.}
A \emph{tree-decomposition} for a graph $G$ is a pair $(T,\beta)$ such that $T$ is a tree, $\beta\colon V(T)\to 2^{V(G)}$ assigns to each node of $T$ a subset of the vertices of $G$ known as a \emph{bag}, $\bigcup_{t\in V(T)}G[\beta(t)] = G,$ and for each $v\in V(G),$ the set $\{ t\in V(T) ~\colon~ v\in\beta(t) \}$ is connected.
The \emph{adhesion} of $(T,\beta)$ is $0$ if $T$ has only one node and $\max_{dt \in E(T)}|\beta(d) \cap \beta(t)|$ otherwise.
The \emph{width} of $(T,\beta)$ is defined as $\max_{t\in V(T)} |\beta(t)|-1$.
The \emph{treewidth} of a graph $G,$ denoted by $\mathsf{tw}(G),$ is the smallest integer $k$ such that $G$ has a tree-decomposition of width at most~$k$.

\paragraph{A rough description of the Flat Wall Theorem.}
The Flat Wall Theorem of Robertson and Seymour \cite{RobertsonS1995Graph} may be seen as a fundamental stepping stone towards the full structure theorem for $H$-minor-free graphs.
Indeed, it is sometimes referred to as the ``weak structure theorem'' as it asserts that every graph of large treewidth must contain a big ``area'' that roughly behaves like a planar graph.
To formally define flat walls requires some amount of machinery that is beyond the scope of this section.
Instead, here we provide some intuition behind the meaning of the Flat Wall Theorem and refer the reader to \zcref{subsec:GraphMinorPrelim} for the formal definitions.

The so-called Two Paths Theorem (see for example \cite{Jung1970Verallgemeinerung,Shiloach1980Polynomial,Seymour1980Disjoint,Thomassen19802Linked}) provides a characterisation for all \textsc{No}-instances of the $2$\textsc{-Disjoint Paths} problem:
Given distinct terminals $s_1,s_2,t_1,t_2$ in a graph $G,$ then there do not exist vertex-disjoint paths $P_1$ and $P_2$ such that $P_i$ connects $s_i$ and $t_i$ in $G$ if and only if $G$ can be ``reduced'' to a planar graph where the four vertices $s_1,s_2,t_1,t_2$ appear on a common face in the order listed.
Here the term ``reduced'' means that one is allowed to hide non-planar parts of the graph behind separators of size at most $3$.
The core observation here is that, in case $s_1,s_2,t_1,t_2$ appear on a common face of some planar graph $H$ in the order listed, then any two paths $P_1$ and $P_2$ as above would be forced to cross each other, thereby contradicting the planarity of $H$.

The Flat Wall Theorem\footnote{Recall that the formal definitions may be found in \zcref{subsec:GraphMinorPrelim} and the variant of the Flat Wall Theorem we are using may be found in \zcref{prop:FindGoodFlatWall}.} now asserts that for any graph $G$ with a large enough wall $W \subseteq G$ one of two things hold:
\begin{enumerate}
 \item either $G$ contains a big clique-minor, or
 \item there exists a small set $A \subseteq V(G)$ and a big subwall $W' \subseteq W - A$ such that, in the subgraph $H$ of $G-A$ consisting of everything that attaches to the ``interior'' of $W',$ there are no two vertex-disjoint paths joining the pairs of opposing corners of $W'$.
\end{enumerate}
In other words, the Two Paths Theorem says that in the second outcome above, the part of $H$ that contains ``most'' of $W'$ can be seen as being planar up to separations of order at most $3$.
We refer to the graph $H$ has the \emph{compass} of the \emph{flat wall} $W'$ in $G-A$.
\medskip

It is an intriguing feature of Robertson and Seymour's Graph Minor Algorithm that most of the effort goes into the proof of its correctness.
That is, in order to show that in any instance $(G,\mathcal{R})$ of large treewidth there exists some vertex $v \in V(G)$ such that $(G,\mathcal{R})$ and $(G - v,\mathcal{R})$ are equivalent instances they employ almost the entire breadth of structure theorems that appeared in their Graph Minors Series.
However, once the correctness is proven, it suffices to find either a clique-minor or a flat wall of large enough order to detect such an \emph{irrelevant} vertex $v$ as above.

\subsection{The folio problem}
\label{subsec:folio}
At the heart of the Graph Minor Algorithm sits the so-called \textsc{Rooted $d$-Folio} problem.
In this problem one is given as input a graph $G$ together with a set of $k$ special vertices $r_1,\dots,r_k$ with labels and an integer $d$.
The question is, what are the minors of $G$ ``rooted'' on the $r_i$ and with ``detail'' at most $d,$ that is, with at most $d$ additional vertices and at most $d$ edges.
This means, instead of finding a single minor of bounded size, or a single pattern of disjoint paths with ends on the roots, the algorithm finds all possible such rooted minors and patterns up to the bound $d + k$ on the number of vertices.

\paragraph{Rooted graphs.}
Let $k\geq 0$ be an integer.
A \emph{$k$-rooted graph} is a pair $(G,\mathcal{R})$ where $\mathcal{R} = \{ r_1,r_2,\dots,r_k \} \subseteq V(G)$ is a multiset of $k$ not necessarily distinct vertices.
We sometimes write $(G,r_1,\dots,r_k)$ instead of $(G,\mathcal{R})$.
The $r_i$ are called the \emph{roots} of $(G,\mathcal{R})$.
A pair $(G,\mathcal{R})$ is a \emph{rooted graph} if it is a $k$-rooted graph for some $k$.

Notice that we make a distinction between a \textsl{rooted} graph $(G,\mathcal{R})$ and an \textsl{annotated} graph $(G,R)$.
One core difference between the two concepts is that for a rooted graph the \textsl{labels} of the roots matter, while for an annotated graph, annotated vertices are interchangeable.
The other key difference is the fact that $\mathcal{R}$ is an \textsl{ordered multiset}.
That is, we allow for vertices to appear more than once in $\mathcal{R}$ and all vertices in $\mathcal{R}$ carry \textsl{distinct} labels.
In fact, as mentioned in the introduction, every rooted graph $(G,\mathcal{R})$ gives rise to a \textsl{uniquely} determined annotated graph $(G,R_{\mathcal{R}})$ where $R_{\mathcal{R}}$ is obtained from $\mathcal{R}$ by forgetting the labels and removing duplicates.
We also extend the notion of bidimensionality of annotated graphs to rooted graphs in this way.
That is, $\mathsf{bidim}(G,\mathcal{R}) \coloneqq \mathsf{bidim}(G,R_{\mathcal{R}})$.
\smallskip

We say that two rooted graphs $(G_1,r^1_1,\dots,r^1_k)$ and $(G_2,r^2_1,\dots,r^2_k)$ are \emph{isomorphic} if there exists an isomorphism $\varphi$ from $G_1$ to $G_2$ such that $\varphi(r^1_i) = r^2_i$ for all $i\in[k]$.
\smallskip

The \emph{detail} of a rooted graph $(G,\mathcal{R})$ is defined as $\mathsf{detail}(G,\mathcal{R}) \coloneqq \max\{ |V(G) \setminus \mathcal{R}|, |E(G)| \}$.
\smallskip

\paragraph{The $d$-folio of a rooted graph.}
A rooted graph $(H,q_1,\dots,q_k)$ is a \emph{rooted minor} of another rooted graph $(G,r_1,\dots,r_k)$ if there exists a minor model $\varphi$ of $H$ in $G$ such that for every $i\in[k]$ it holds that $r_i \in \varphi(q_i)$.

Given an integer $d\geq 0$ and a rooted graph $(G,\mathcal{R}),$ the \emph{$d$-folio} of $G,$ denoted as $d\text{-}\mathsf{folio}(G,\mathcal{R}),$ is the set of all rooted graphs $(H,\mathcal{Q})$ of detail at most $d$ that are rooted minors of $(G,\mathcal{R})$.

\paragraph{The $(k,d)$-folio with respect to a vertex set.}
We now let $(G,R)$ be an annotated graph.
The \emph{$(k,d)$-folio} of $(G,R),$ denoted by $(k,d)\text{-}\mathsf{folio}(G,R),$ is the union of the $d$-folios over all rooted graphs $(G,\mathcal{R})$ where $\mathcal{R} = \{ r_1,\dots,r_k\}$ is an ordered multisubset of $R$.

\paragraph{The $(k,d)$-\textsc{Folio} problem and its rooted sibling.}
The \emph{$(k,d)$-\textsc{Folio} problem} may now be formally defined as the problem that takes as input an annotated graph $(G,R)$ and whose task is the computation of the $(k,d)$-folio of $(G,R)$.

The \emph{\textsc{Rooted $d$-Folio} problem} is closely related to the \textsc{$(k,d)$-Folio} problem.
Here, the input is a rooted graph $(G,\mathcal{R})$ with $|\mathcal{R}| \leq k$ and the task is to determine $d\text{-}\mathsf{folio}(G,\mathcal{R})$.
\medskip

We stress that \textsc{$(k,d)$-Folio} with input $(G,R)$ may be reduced to computing $d\text{-}\mathsf{folio}(G,\mathcal{R})$ for all rooted graphs $(G,\mathcal{R})$ where $\mathcal{R}$ is a set of $k$ roots selected from $R$.

To see that the $k$-\textsc{Disjoint Paths} problem is contained in the \textsc{Rooted} $0$-\textsc{Folio} problem on $2k$ roots, consider the following example.
Let $\mathbf{H} = (H,x_1,x_1,x_2,x_2,\dots,x_k,x_k)$ be the rooted graph with $2k$ roots and vertex set $\{ x_1,\dots,x_k \}$ and where $E(H) = \emptyset$.
Now, if $\mathbf{G} = (G,s_1,t_1,s_2,t_2,\dots,s_k,t_k)$ is a rooted graph with distinct vertices $s_1,\dots,s_k,t_1,\dots,t_k,$ then there exist $k$ pairwise vertex-disjoint paths $P_1,\dots,P_k$ such that $P_i$ has endpoints $s_i$ and $t_i$ if and only if $\mathbf{H}$ is a rooted minor of $\mathbf{G}$.
To see this, notice that the existence of the paths $P_i$ is enough to produce a rooted minor model of $\mathbf{H}$ in $\mathbf{G}$.
For the other direction, if there is a rooted minor model $\varphi$ of $\mathbf{H}$ in $\mathbf{G},$ then for each $i \in[k]$ it must be true that $s_i,t_i \in \varphi(x_i)$ and, since the $\varphi(x_i)$ are connected and pairwise disjoint it follows that each $\varphi(x_i)$ contains one such path $P_i$.

It is easily seen that by definition, $\mathbf{H}$ is a rooted graph with detail $0$ and thus, in order to solve the $k$-\textsc{Disjoint Paths} problem on the instance $\mathbf{G},$ it suffices to solve the \textsc{Rooted} $0$-\textsc{Folio} problem on the rooted graph $(G,\{ s_1,\dots,s_k,t_1,\dots,t_k\})$ and check if $\mathbf{H} \in 0\text{-}\mathsf{Folio}(G,\{ s_1,\dots,s_k,t_1,\dots,t_k\})$.

Hence, in order to prove \zcref{thm:intro_folio_algorithm_informal}, \zcref{prop_rs_meta}, and to obtain a double-exponential \textsf{FPT}-algorithm for $k$-\textsc{Disjoint Paths} it suffices to prove the following theorem.

\begin{theorem}\label{thm:Main1Folio}
There exists an algorithm for \textsc{Rooted $d$-Folio} that takes as input a rooted graph $(G,\mathcal{R})$ with $\mathsf{bidim}(G,\mathcal{R}) \leq b$ and $|\mathcal{R}| \leq k$ and runs in time $2^{2^{\mathbf{poly}(b + d)} \cdot \mathbf{poly}(k)} \cdot |G|^2$.
\end{theorem}

Since the number of possible $k$-tuples with elements from a set $R$ of order $k$ is in $\mathbf{O}(k^k)$ and by our discussion on how \textsc{$(k,d)$-Folio} may be reduced to \textsc{Rooted $d$-Folio}, \zcref{thm:Main1Folio} gives rise to the following corollary. 

\begin{corollary}\label{cor:Main1Folio}
There exists an algorithm for \textsc{$(k,d)$-Folio} that takes as input an annotated graph $(G,R)$ with $\mathsf{bidim}(G,R) \leq b$ and $|R| \leq k$ and runs in time $2^{2^{\mathbf{poly}(b + d)} \cdot \mathbf{poly}(k)} \cdot |G|^2$.
\end{corollary}

When $(G,R)$ is an annotated graph with $|R| \leq k$ then known algorithms are already able to solve \textsc{$(k,d)$-Folio} and its rooted counterpart in single-exponential time depending on the treewidth of $G$.
See for example the algorithm of Adler, Dorn, Fomin, Sau, and Thilikos \cite{AdlerDFST2011Faster}.

\begin{proposition}[Adler, Dorn, Fomin, Sau, and Thilikos \cite{AdlerDFST2011Faster}]\label{prop:FolioTW}
There exists an algorithm for \textsc{Rooted $d$-Folio} that takes as input a rooted graph $(G,\mathcal{R})$ with $\mathsf{tw}(G) \leq w$ and $|\mathcal{R}| \leq k$ and runs in time $2^{\mathbf{poly}(w + k + d)} \cdot |G|$.
\end{proposition}

Hence, the only remaining step towards deriving \zcref{thm:Main1Folio} is an algorithmic version of \zcref{thm:intro_folio_irrelevant_informal}.
In particular, in order to meet the quadratic running time of the algorithm from \zcref{thm:Main1Folio}, we require a way to find a single irrelevant vertex in \textsl{linear time}.
Indeed, the following is the algorithmic core of this paper.

\paragraph{Strongly irrelevant vertices.}
Let $k,d$ be non-negative integers and $(G,R)$ be an annotated graph.
We say that a vertex $v \in V(G) \setminus R$ is \emph{strongly irrelevant} for $(k,d)\text{-}\mathsf{folio}(G,R)$ if for every ordered multiset $\mathcal{R}$ consisting of $k$ roots chosen from among the vertices of $R$ it holds that $d\text{-}\mathsf{folio}(G - v,\mathcal{R}) = d\text{-}\mathsf{folio}(G,\mathcal{R})$.

\begin{theorem}\label{thm:Main2Folio}
There exists a function $f_{\ref{thm:Main2Folio}} \colon \mathbb{N}^3 \to \mathbb{N}$ and an algorithm that takes as input an annotated graph $(G,R)$ of bidimensionality at most $b,$ and either correctly determines that $\mathsf{tw}(G) \leq f_{\ref{thm:Main2Folio}}(k,b,d)$ or finds a vertex $v\in V(G)\setminus R$ such that $v$ is strongly irrelevant for $(k,d)\text{-}\mathsf{folio}(G,R)$ in time $2^{2^{\mathbf{poly}(b + d)}\cdot \mathbf{poly}(k)}\cdot |G|$.
Moreover, $f_{\ref{thm:Main2Folio}}(k,b,d) \in 2^{\mathbf{poly}(b + d)}\cdot \mathbf{poly}(k)$.
\end{theorem}

Here something interesting happens.
Our observation that \textsc{$(k,d)$-Folio} can naturally be reduced to solving \textsc{Rooted $d$-Folio} on all arising rooted graphs also gives rise to a reverse observation when it comes to the irrelevance of some vertex.
If there is a vertex $v \in V(G)$ that is strongly irrelevant for $(k,d)\text{-}\mathsf{folio}(G,R)$ then this means in particular that also $d\text{-}\mathsf{folio}(G,\mathcal{R}) = d\text{-}\mathsf{folio}(G-v,\mathcal{R})$ for all possible ways to select an ordered multiset $\mathcal{R}$ of $k$ roots from $R$.
Hence, in the main body of the paper we will only be concerned with the annotated version, that is with $(k,d)\text{-}\mathsf{folio}(G,R)$.

The linear running time of the algorithm from \zcref{thm:Main2Folio} allows for the following simple algorithm meeting the requirement of \zcref{thm:Main1Folio}.
\medskip

\textbf{Input:} A rooted graph $(G,\mathcal{R})$ with $\mathsf{bidim}(G,R_{\mathcal{R}}) \leq b$ and $|\mathcal{R}| \leq k$.
\begin{description}
 \item[Step 1:] Run the algorithm from \zcref{thm:Main2Folio} on $(G,\mathcal{R})$.
 This yields one of two outcomes:
 \begin{description}
 \item[Case 1.1:] \textsl{The treewidth of $G$ is at most $f_{\ref{thm:Main2Folio}}(k,b,d)$.}
 In this case we \textbf{forward} $(G, \mathcal{R})$ to \textbf{Step 2}.

 \item[Case 1.2:] \textsl{The algorithm returns a vertex $v \in V(G)\setminus R_{\mathcal{R}}$ such that $v$ is strongly irrelevant for $(k,d)\text{-}\mathsf{folio}(G,R_{\mathcal{R}})$.}
 In this case \textbf{forward} $(G-v, \mathcal{R})$ back to the beginning of \textbf{Step 1}.
 \end{description}

 \item[Step 2:] Run the algorithm from \zcref{prop:FolioTW} on $(G,\mathcal{R}),$ then \textbf{return} $d\text{-}\mathsf{folio}(G,\mathcal{R})$ and \textbf{terminate} the algorithm.
\end{description}

Since the instance shrinks by one vertex every time the algorithm enters \textbf{Case 1.2}, the loop in \textbf{Step 1} repeats at most $|G|$ times.
Hence, since the algorithm from \zcref{thm:Main2Folio} runs in linear time, we are sure to enter \textbf{Step 2} after at most $|G|$ iterations and at most $|G|^2$ steps.

\section{Preliminaries}\label{sec:prelims}
There are a variety of established frameworks for the work in the structural theory surrounding the Graph Minors Series and careful definitions are needed throughout.
We opt to mainly work from a set of definitions for annotated graphs that is taken from \cite{GorskyPW2026Quickly} and which are largely based on definitions that evolved from the Graph Minors Series by Robertson and Seymour and were refined in \cite{KawarabayashiTW2018New,KawarabayashiTW2021Quickly,GorskySW2025Polynomial}.

Later on in \zcref{sec:folio}, we will expand our view and consider graphs in which vertices may be marked with many different colours.
As we will need these notions briefly in \zcref{sec:largeclique}, we introduce them here as well, but will be able to avoid working in this more complicated setting for most of the paper.

\subsection{General definitions}
Here, we introduce a series of fundamental definitions required throughout the paper.

\paragraph{Paths, linkages, and separations.}
A \emph{linkage} $\mathcal{L}$ in a graph $G$ is a set of pairwise vertex-disjoint paths.
Let $T \subseteq V(G)$ be the set of endpoints of paths in $\LLL,$ then we call $\LLL$ a \emph{$T$-linkage}.
We define the \emph{pattern of $\LLL$} via $\tau(\LLL) \coloneqq \{\{s,t\} \mid \text{ there is }L\in \LLL \text{ with endpoints } s \text{ and } t\}$. 

We say that a path $P$ in $G$ is \emph{internally disjoint} from a set $X \subseteq V(G)$ if $V(P) \cap X$ does not contain any vertex of $P$ that is not an endpoint.
Given a graph $G$ and two subsets $A, B \subseteq V(G),$ an \emph{$A$-path} in $G$ is a path with both endpoints in $A$ and internally disjoint from $A,$ and an \emph{$A$-$B$-path} is a path with one endpoint in $A,$ the other in $B,$ and internally disjoint from $A \cup B$.
An \emph{$A$-$B$-linkage} in $G$ is a linkage consisting of $A$-$B$ paths.
If $H$ is a subgraph of $G,$ an \emph{$H$-path} is a $V(H)$-path of length at least one with no edge in $E(H)$.

Let $G$ be a graph and $k$ be a positive integer.
A \emph{separation} of $G$ is a tuple $(A,B)$ such that $A,B\subseteq V(G)$ with $A\cup B = V(G)$ and for every edge $uv \in E(G)$ we have $u,v \in A$ or $u,v \in B$.
We call $\Abs{A\cap B}$ its \emph{order}.
The following classic tool will be helpful to us in many ways.

\begin{proposition}[Menger~\cite{Menger1927Zur}]\label{thm:menger}
 Let $k$ be a positive integer, let $G$ be a graph, and let $X,Y \subseteq V(G)$ be two sets with at least $k$ elements each.
 Then there does not exist a separation $(A,B)$ in $G$ of order less than $k$ with $X \subseteq A$ and $Y \subseteq B$ if and only if there exist $k$ pairwise disjoint $X$-$Y$-paths in $G$.
\end{proposition}

\paragraph{Surfaces.}
By a \emph{surface} we mean a compact $2$-dimensional manifold $\Sigma$ with or without boundary.
Let $X \subseteq \Sigma$. The boundary and interior of $X$ will be denoted $\bd(X)$ and $\mathsf{int}(X),$ respectively. The topological closure of $X$ is denoted by $\mathsf{cl}(X)$. If $\Sigma$ is a surface with boundary then we refer to the connected components of $\bd(\Sigma)$ as its \emph{cuffs}.

Given a pair $(\mathsf{h}, \mathsf{c}) \in \mathbb{N} \times \mathbb{N}$ 
we define $\Sigma^{(\mathsf{h}, \mathsf{c})}$ to be the surface without boundary created from the sphere by adding $\mathsf{h}$ handles and $\mathsf{c}$ crosscaps (see \cite{MoharT2001Graphs} for more details). By the surface classification theorem, every compact connected surface without boundary is homeomorphic to $\Sigma^{(\mathsf{h}, \mathsf{c})}$ for some $h, c \geq 0$. 
Note that by Dyck's theorem \cite{Dyck1888Beitraege,FrancisW1999Conways}, two crosscaps are equivalent to a handle in the presence of a third crosscap.

If $\mathsf{c} = 0$ the surface $\Sigma^{(\mathsf{h}, \mathsf{c})}$ is an \emph{orientable} surface, otherwise it is called \emph{non-orientable}.
We let the \emph{genus} of $\Sigma$ be $2\mathsf{h} + \mathsf{c},$ where $\Sigma^{\mathsf{(\mathsf{h},\mathsf{c})}}$ is a surface to which $\Sigma$ is homeomorphic.

Let $\Sigma$ and $\Gamma$ be surfaces. A \emph{homeomorphism} $h \sth \Sigma \rightarrow \Gamma$ is a bijective function from $\Sigma$ to $\Gamma$ which is continuous and has a continuous inverse. A homeomorphism $h \sth \Sigma \rightarrow \Sigma$ on a surface $\Sigma$ with boundary such that $h$ fixes every point on the boundary of $\Sigma$ is called a \emph{$\bd$-homeomorphism}.

Let $\gamma \sth [0,1]_{\mathbb{R}}\to \Sigma$ be continuous and injective on $[0,1[_{\mathbb{R}},$ then we call its image $\gamma([0,1]_{\mathbb{R}})$ a \emph{curve} in $\Sigma$ and $\gamma(0),\gamma(1)$ its endpoints. A curve is \emph{closed} if its endpoints agree. A curve with both endpoints on the boundary of $\Sigma$ is called a \emph{boundary curve}.

Two curves $C, C' \subseteq \Sigma$ are \emph{homeomorphic} if there is a homeomorphism $h \sth \Sigma \rightarrow \Sigma$ which maps $C$ to $C'$. 
$C$ and $C'$ are \emph{homotopic} if there is a continuous function $f \sth \Sigma \times [0,1]_{\mathbb{R}} \rightarrow \Sigma$ such that $f(x, 0) = C(x)$ and $f(x, 1) = C'(x)$ and for every $t\in [0,1]_{\mathbb{R}} $ the function $f_t(x)\coloneqq f(x,t)$ maps the boundary of $\Sigma$ to the boundary of $\Sigma$. We call $f$ a \emph{homotopy}.
Homotopy induces an equivalence relation on the set of curves on a surface $\Sigma$. We refer to the equivalence classes as \emph{homotopy classes} and to the class containing a curve $C \subseteq \Sigma$ its \emph{homotopy type}.

Let $\Sigma$ be a connected surface. A curve $C \subseteq \Sigma$ is \emph{contractible} if one of the components of $\Sigma - C$ is homeomorphic to an open disc. Otherwise it is \emph{non-contractible}. $C$ is \emph{separating} if $\Sigma-C$ is not connected. A non-separating, non-contractible curve is called \emph{essential} or \emph{genus reducing}.

The following is a simple result on the number of homotopy classes of non-contractible curves on a given surface that can be realised simultaneously\footnote{This lemma in different formulation is almost certainly found in standard topology texts. Our use for it is related to \cite{Mazoit2013Single} (see Claim 7 in said paper for comparison.)}. In \cite{JuvanMM1996}, Juvan et al. prove similar results for closed curves disjoint from the boundary. See e.g. Lemma 3.2 in that paper.

\begin{lemma}\label{lem:typecounting}
 Let $\Sigma$ be a connected surface of genus $g$ with $b \geq 1$ boundary
 components. Let $\LLL$ be a set of pairwise disjoint simple curves on
 $\Sigma$ whose endpoints are all on the boundary of $\Sigma$.

 Then the number of homotopy types realised by curves in $\LLL$ is $\leq 1$ in
 case $b=1$ and $g=0,$ it is $\leq 3$ in case $b=2$ and $g=0,$ and otherwise it is bounded by $4b + 6g - 6$.
\end{lemma}
\begin{proof}
 In the special case that $b=1$ and $g=0$ $\Sigma$ is a disc and $\LLL$ is a
 set of simple curves with all endpoints on the boundary. Thus all
 curves in $\LLL$ are null-homotopic and therefore there is only one
 homotopy type that is realised by $\LLL$. 

 The other special case we need to consider is $b=2$ and $g=0$. In this case $\Sigma$ is a cylinder and there are at most three possible homotopy types that can be realised simultaneously: one type for each of the two cuffs realised by curves with both endpoints on that cuff and the third type realised by a curve with both endpoints on different cuffs. Note that there is an unbounded number of pairwise non-homotopic curves with endpoints on different cuffs, which are essentially determined by the number of times they wind around the cylinder. But no two of these can be realised at the same time in $\LLL$.
 
 So we may now assume that $b>1$ or $g>0$.
 Let $L \in \LLL$ be a curve with both endpoints $a, b$ on the same boundary component $C$ of $\Sigma$. We say that $L$ is null-homotopic if $L$ together with one of the two components of $C - \{a,b\}$ bounds a disc on $\Sigma$. Clearly any two null-homotopic curves with endpoints on the same boundary are homotopic. Thus there are at most $b$ null-homotopic curves in $\LLL$ which realise different homotopy types in $\Sigma$.
 
 Let $M \subseteq \LLL$ be a maximal subset of $\LLL$ such that no two curves in $M$ are homotopic and furthermore no curve in $M$ is null-homotopic in the sense above. Let $\Sigma'$ be obtained from $\Sigma$ by closing each cuff $C$ with a disc $\Delta(C)$. For each $C$ let $T(C)$ be the endpoints of curves in $M$ on $C$. For each cuff $C$ choose a point $v = v(C) \in \mathsf{int}\big(\Delta(C)\big)$ and for each $u \in T(C)$ a curve $L_u$ in $\Delta(C)$ between $v$ and $u$ such that for $u \not= u' \in T(C)$ these curves are disjoint except for $v$. We can then extend each curve $L \in M$ with endpoints $a, b$ by the curves $L_a, L_b$ to obtain a new set $M'$ of curves whose endpoints are all contained in $\{ v(C) \sth C$ a cuff of $\Sigma \}$. Without loss of generality~we assume that $M'$ is maximal in the sense that there is no non-null-homotopic curve with endpoints in $\{ v(C) \sth C$ a cuff of $\Sigma \}$ that is internally disjoint from each curve in $M'$ and is not homotopic to a curve in $M'$. 
 
 As $M'$ does not contain any null-homotopic curve, there is a
 triangulation $N$ of $\Sigma'$ such that $M' \subseteq N$. See e.g.~\cite{FarbM2012}. %
 By Euler's formula \cite{Euler1758Elementa}, $|V(N)| - |E(N)| + |F(N)| = 2 - 2g,$ where $V(N)$ denotes the set of vertices of $N,$ $E(N)$ the set of arcs, and $F(N)$ the set of faces of the embedding $N$. By construction, $|V(N)| = b$ and $|E(N)| = |N| \geq |M'|$. Let $f \coloneqq |F(N)|$.
 Thus we obtain $|N| = b + f +2g - 2$.

 On the other hand, as $N$ triangulates $\Sigma',$ we obtain $2 |N| = 3f$. Combining the two equations yields $|N| = 3b + 6g - 6$. 

 Thus, together with the at most $b$ null-homotopic curves we removed above, the maximal number of pairwise non-homotopic curves in $\LLL$ is $\leq 4b + 6g - 6$.
\end{proof}

\paragraph{Embeddings.}
Given a graph $G,$ we say that $G$ has an \emph{embedding} $\psi$ into a surface $\Sigma,$ if $\psi$ is a function with domain $V(G) \cup E(G)$ such that
\begin{enumerate}
 \item $\psi(v)$ is a point of $\Sigma$ for each $v \in V(G),$
 \item $\psi(u) \neq \psi(v)$ for all distinct $u,v \in V(G),$
 \item $\psi(uv)$ is a simple curve in $\Sigma$ with the ends $\psi(u),\psi(v)$ for all $uv \in E(G),$
 \item $\psi(uv) \cap \psi(wx) \subseteq \{ \psi(u), \psi(v) \},$ for all $uv, wx \in E(G),$ and
 \item $\psi(uv) \cap \psi(w) = \emptyset$ for all $uv \in E(G)$ and $w \in V(G) \setminus \{ u,v \}$.
\end{enumerate}
We let $\psi(G)$ be the union of all points and curves in its image and call the connected components of $\Sigma \setminus \psi(G)$ the \emph{faces} of $\psi$.
If every face of $\psi$ is homeomorphic to an open disc, we say that $\psi$ is a \emph{$2$-cell-embedding}.

Given a sequence of cycles $\CCC = \langle C_1,\ldots,C_s\rangle$ in $G$ we call these cycles \emph{concentric (with respect to $\psi$)} if there is a disc $\Delta \subseteq \Sigma$ such that $\psi(\bigcup \CCC) \subseteq \Delta$ and $\psi(C_i)$ bounds a disc $\Delta_i$ such that $\Delta_1 \subsetneq \ldots \subsetneq \Delta_s$.

We will generally work with graphs with a fixed embedding and thus identify the elements of $G$ with the points and curves in $\Sigma$ which $\psi$ maps to each other.
The minimum number of vertices of $G$ that is intersected by any non-contractible curve in $\Sigma$ meeting $G$ only in vertices is called the \emph{representativity} (or \emph{face-width}) of $\psi$.
We consider the representativity of an embedding in the sphere to be infinite.

\subsection{Definitions from Graph Minors}
\label{subsec:GraphMinorPrelim}
We continue by briefly introducing some key concepts for handling the general structure of $H$-minor-free graphs.

\paragraph{Separations and tangles.}
Let $G$ be a graph and $k$ be a positive integer.
Let $t\geq 1$. A sequence of separations $\langle(A_i,B_i)\rangle_{i \in [t]}$ of $G$ is \emph{laminar} if $A_i\subseteq A_j$ and $B_j\subseteq B_i$ for every pair $1\leq i \leq j \leq t$.
We denote by $\mathcal{S}_k(G)$ the collection of all separations $(A,B)$ of order less than $k$ in $G$. 

An \emph{orientation} of $\mathcal{S}_k(G)$ is a set $\mathcal{O}$ such that for all $(A,B)\in\mathcal{S}_k(G)$ exactly one of $(A,B)$ and $(B,A)$ belongs to $\mathcal{O}$. 
A \emph{tangle} of order $k$ in $G$ is an orientation $\mathcal{T}$ of $\mathcal{S}_k(G)$ such that for all $(A_1,B_1),(A_2,B_2),(A_3,B_3)\in\mathcal{T},$ it holds that $G[A_1]\cup G[A_2]\cup G[A_3]\neq G$.
If $\mathcal{T}$ is a tangle and $(A,B)\in\mathcal{T}$ we call $A$ the \emph{small side} and $B$ the \emph{big side} of $(A,B)$.

Let $G$ be a graph and $\mathcal{T}$ and $\mathcal{D}$ be tangles of $G$.
We say that $\mathcal{D}$ is a \emph{truncation} of $\mathcal{T}$ if $\mathcal{D}\subseteq\mathcal{T}$.
\medskip

Let $G$ and $H$ be graphs as well as $\mathcal{T}$ be a tangle in $G$.
We say that a minor-model $\mu$ of $H$ in $G$ is \emph{controlled} by $\mathcal{T}$ if there does not exist a separation $(A,B)\in\mathcal{T}$ of order less than $|V(H)|$ and an $x \in V(H)$ such that $\mu(x)\subseteq A\setminus B$.

Let $\mathcal{T}$ be a tangle of order $k$ in a graph $G$ and let $S\subseteq V(G)$ be a set of size at most $k-2.$
Notice that every separation $(A,B)$ of order at most $k-|S|-1$ in $G-S$ corresponds to a separation $(A\cup S,B\cup S)$ of order at most $k-1$ in $G.$
It follows that $\mathcal{T}_S\coloneqq \{ (A,B) ~\colon~ (A\cup S,B\cup S)\in\mathcal{T} \}$ is a tangle of order $k-|S|$ in $G-S.$
Now suppose that $G-S$ has components $G_1,\dots,G_p.$
Then there exists $i\in[p]$ such that $V(G_i)\subseteq B$ for all separations of order $0$ in $\mathcal{T}_S.$
We call $G_i$ the \emph{$\mathcal{T}$-big component} of $G - S$.

\paragraph{Meshes.}
Let $n,m$ be integers with $n,m\geq 2$.
A \emph{$(n\times m)$-mesh} is a graph $M$ which is the union of paths $M=P_1\cup\cdots\cup P_n\cup Q_1\cup \cdots \cup Q_m$ where
 \begin{itemize}
 \item $P_1,\cdots,P_n$ are pairwise vertex-disjoint, and $Q_1,\cdots,Q_m$ are pairwise vertex-disjoint.
 \item for every $i\in [n]$ and $j\in [m],$ the intersection $P_i\cap Q_j$ induces a path,
 \item each $P_i$ is a $V(Q_1)$-$V(Q_m)$-path intersecting the paths $Q_1,\cdots Q_m$ in the given order, and each $Q_j$ is a $V(P_1)$-$V(P_m)$-path intersecting the paths $P_1,\cdots, P_h$ in the given order. 
 \end{itemize}
We say that the paths $P_1,\cdots,P_n$ are the \emph{horizontal paths}, and the paths $Q_1,\cdots,Q_m$ are the \emph{vertical paths}.
The union $P_{1} \cup P_{n} \cup Q_{1} \cup Q_{m}$ is a cycle called the \emph{perimeter} of $M$.
The unique cycle in the union $P_{i} \cup P_{i+1} \cup Q_{j} \cup Q_{j+1},$ where $i \in [n - 1]$ and $j \in [m - 1],$ is called a \emph{brick} of $M$.
A mesh $M'$ is a \emph{submesh} of a mesh $M$ if every horizontal (vertical) path of $M'$ is a subpath of a horizontal (vertical) path of $M,$ respectively.
We write \emph{$n$-mesh} as a shorthand for an $(n \times n)$-mesh.

Let $r \in \mathbb{N}$ with $r\geq 3,$ let $G$ be a graph, and $M$ be an $r$-mesh in $G$.
Let $\mathcal{T}_M$ be the orientation of $\mathcal{S}_r$ such that for every $(A,B)\in\mathcal{T}_M,$ the set $B\setminus A$ contains the vertex set of both a horizontal and a vertical path of $M$. We call $B$ the \emph{$M$-majority side} of $(A,B)$.
Then $\mathcal{T}_M$ is the tangle \emph{induced} by $M$.
If $\mathcal{T}$ is a tangle in $G,$ we say that $\mathcal{T}$ \emph{controls} the mesh $M$ if $\mathcal{T}_M$ is a truncation of $\mathcal{T}$.

\paragraph{Paintings in surfaces.}
A \emph{painting} in a surface $\Sigma$ is a pair $\Gamma = (U,N),$ where $N \subseteq U \subseteq \Sigma,$ $N$ is finite, $U \setminus N$ has a finite number of arcwise-connected components, called \emph{cells} of $\Gamma,$ such that each cell is homeomorphic to an open disc in $\Sigma,$ $\mathsf{cl}(c)$ is a closed disc in $\Sigma,$ and $\mathsf{cl}(c) \cap \mathsf{cl}(c') \subseteq N$ for two distinct cells. We define $N_\Gamma(c) \coloneqq \mathsf{cl}(c) \cap N$.
If $|N_\Gamma(c)| \geq 4,$ the cell $c$ is called a \emph{vortex}.
We further let $N(\Gamma) \coloneqq N,$ let $U(\Gamma) \coloneqq U,$ and let $C(\Gamma)$ be the set of all cells of $\Gamma$.

Any given painting $\Gamma = (U,N)$ defines a hypergraph with $N$ as its vertices and the set of closures of the cells of $\Gamma$ as its edges.
Accordingly, we call $N$ the \emph{nodes} of $\Gamma$.

\paragraph{$\Sigma$-renditions.}
Let $G$ be a graph and $\Sigma$ be a surface.
A \emph{$\Sigma$-rendition} of $G$ is a triple $\rho = (\Gamma, \sigma, \pi),$ where
\begin{itemize}
 \item $\Gamma$ is a painting in $\Sigma,$
 \item for each cell $c \in C(\Gamma),$ $\sigma(c)$ is a subgraph of $G,$ and
 \item $\pi \colon N(\Gamma) \to V(G)$ is an injection,
\end{itemize}
such that
\begin{description}
 \item[R1] $G = \bigcup_{c \in C(\Gamma)}\sigma(c),$
 \item[R2] for all distinct $c,c' \in C(\Gamma),$ the graphs $\sigma(c)$ and $\sigma(c')$ are edge-disjoint,
 \item[R3] $\pi(N_\Gamma(c)) \subseteq V(\sigma(c))$ for every cell $c \in C(\Gamma),$ and
 \item[R4] for every cell $c \in C(\Gamma),$ we have $V(\sigma(c) \cap \bigcup_{c' \in C(\Gamma) \setminus \{ c \}} (\sigma(c'))) \subseteq \pi(N_\Gamma(c))$.
\end{description}
We write $N(\rho)$ for the set $N(\Gamma),$ let $N_\rho(c) = N_\Gamma(c)$ for all $c \in C(\Gamma),$ and similarly, we lift the set of cells from $C(\Gamma)$ to $C(\rho)$.
If it is clear from the context which $\rho$ is meant, we will sometimes simply write $N(c)$ instead of $N_\rho(c),$ and if the $\Sigma$-rendition $\rho$ for $G$ is understood from the context, we usually identify the sets $\pi(N(\rho))$ and $N(\rho)$ along $\pi$ for ease of notation.

Given a graph $G$ with a $\Sigma$-rendition $\rho$ and a subgraph $G' \subseteq G,$ we define the \emph{rendition $\rho'$ of $G'$ induced by $\rho$} by letting the nodes of $\rho'$ be $N(\rho) \cap V(G'),$ letting the cells of $\rho'$ be the cells $c$ of $\rho$ such that $(\sigma_\rho(c) \cap G') - N(\rho)$ is non-empty, and adjusting $\sigma$ appropriately. 

Note that every embedding $\psi$ of a graph $G$ in a surface $\Sigma$ defines a rendition $\rho_{\psi}$ of $G$ in $\Sigma$ where the graph of every cell corresponds to a unique edge of $G$. The cell itself can be an open disc with the two endpoints of the edge in its boundary.
Thus all definitions for renditions naturally transfer to embeddings.

\paragraph{Blank renditions and renditions of annotated graphs.}
Let $\rho = (\Gamma, \sigma, \pi)$ be a $\Sigma$-rendition of an annotated graph $(G,R)$.
By slightly abusing notation, we call a cell $c$ of $\rho$ a \emph{vortex} if it is either a vortex of $\Gamma$ or $V(\sigma(c)) \cap R \neq \emptyset$.
If $$\pi(N(\rho)) \cup \bigcup \{ V(\sigma(c)) \colon c \in C(\rho) \text{ and } c \text{ is not a vortex}\}$$ is disjoint from $R,$ we call $\rho$ a \emph{blank rendition (of $(G,R)$)}.

\paragraph{Societies.}
Let $\Omega$ be a cyclic ordering of the elements of some set.
Given that $\Omega=\langle v_{1},\ldots,v_{r}\rangle,$ then we denote the set $\{v_{1},\ldots,v_{r}\}$ by $V(\Omega)$.
A \emph{society} is a pair $(G,\Omega),$ where $G$ is a graph and $\Omega$ is a cyclic ordering with $V(\Omega)\subseteq V(G)$.
For a given set $S \subseteq V(\Omega)$ an element
 $s \in S$ is an \emph{endpoint} of $S$ if there exists an element $t \in V(\Omega) \setminus S$ that immediately precedes or succeeds $s$ in $\Omega$.
We call $S$ a \emph{segment} of $\Omega$ if $S$ has two or less endpoints. We call a segment \emph{proper} if $V(S) \neq V(\Omega)$. The \emph{ordering of $S$ induced by $\Omega$} is the restriction of $\Omega$ to $S$.

Let $(G,\Omega)$ be a society and let $\Sigma$ be a surface with one boundary component $B$ homeomorphic to the unit circle.
A \emph{rendition} of $(G,\Omega)$ in $\Sigma$ is a $\Sigma$-rendition $\rho$ of $G$ such that the image under $\pi_{\rho}$ of $N(\rho) \cap B$ is $V(\Omega)$ and $\Omega$ is one of the two cyclic orderings of $V(\Omega)$ defined by the way the points of $\pi_{\rho}(V(\Omega))$ are arranged in the boundary $B$.

\paragraph{Traces of paths and cycles.}
Let $\rho$ be a $\Sigma$-rendition of a graph $G$.
For every cell $c \in C(\rho)$ with $|N_\rho(c)| = 2,$ we select one of the components of $\mathsf{bd}(c) \setminus N(c)$.
This selection will be called a \emph{tie-breaker in $\rho$}, and we assume that every rendition comes equipped with a tie-breaker.

Let $G$ be a graph and $\rho$ be a $\Sigma$-rendition of $G$.
Let $Q$ be a cycle or path in $G$ that uses no edge of $\sigma(c)$ for every vortex $c \in C(\rho)$.
We say that $Q$ is \emph{grounded} if it uses edges of $\sigma(c')$ and $\sigma(c'')$ for two distinct cells $c', c'' \in C(\rho),$ or $Q$ is a path with both endpoints in $N(\rho)$.
If $Q$ is grounded we define the \emph{trace} of $Q$ as follows.
Let $P_1,\dots,P_k$ be distinct maximal subpaths of $Q$ such that $P_i$ is a subgraph of $\sigma(c)$ for some cell $c_i$.
Fix $i \in [k]$.
The maximality of $P_i$ implies that its endpoints are $\pi(n_1)$ and $\pi(n_2)$ for distinct nodes $n_1,n_2 \in N(\rho)$.
If $|N_\rho(c_i)| = 2,$ let $L_i$ be the component of $\mathsf{bd}(c_i)$ selected by the tie-breaker, and if $|N_\rho(c_i)| = 3,$ let $L_i$ be the component of $\mathsf{bd}(c_i) \setminus (N(c_i) \cap V(P_i))$ that is disjoint from $N_\rho(c_i)$.
We define $L_i'$ by pushing $L_i$ slightly so that it is disjoint from all cells in $C(\rho),$ while maintaining that the resulting curves intersect only at a common endpoint.

The \emph{trace} of $Q$ is defined to be $\bigcup_{i\in[k]} L_i'$.
If $Q$ is a cycle, its trace is the homeomorphic image of the unit circle, and otherwise, it is an arc in $\Sigma$ with both endpoints in $N(\rho)$.

\paragraph{Linear decompositions of societies.}
Let $(G,\Omega)$ be a society.
A \emph{linear decomposition} of $(G,\Omega)$ is a labelling $v_1,v_2,\dots,v_n$ of $V(\Omega)$ such that $v_1,v_2,\dots,v_n$ appear in $\Omega$ in the order listed, together with sets ordered as $\langle X_1,X_2,\dots,X_n\rangle$ such that
\begin{enumerate}
 \item $X_i\subseteq V(G)$ and $v_i\in X_i$ for all $i\in[n],$
 \item $\bigcup_{i\in[n]}X_i=V(G)$ and for every $uv\in E(G)$ there exists $i\in[n]$ such that $u,v\in X_i,$ and
 \item for every $x\in V(G)$ the set $\{ i\in[n] ~\!\colon\!~ x\in X_i \}$ forms an interval in $[n]$.
\end{enumerate}
The \emph{adhesion} of a linear decomposition is $\max \{ |X_i\cap X_{i+1}| ~\!\colon\!~ i\in[n-1] \}.$
The \emph{width} of a linear decomposition is $\max \{ |X_i| ~\!\colon\!~ i\in[n] \}.$
We write $\alpha=\left(\langle v_1,\ldots,v_n\rangle,\langle X_1,\ldots,X_n\rangle\right)$ for the respective linear decomposition. Let $1 \leq i \leq j \leq n$ and $S=[i,j]$. Then we call $(\langle v_i,\ldots,v_j\rangle,\langle X_i,\ldots,X_j\rangle)$ the \emph{$S$-restriction of $\alpha$}.

\paragraph{Aligned curves and grounded subgraphs.}
Let $G$ be a graph and let $\rho = (\Gamma, \sigma, \pi)$ be a $\Sigma$-rendition of $G$. 
We say that a $2$-connected subgraph $H$ of $G$ is \emph{grounded (in $\rho$)} if every cycle in $H$ is grounded and no vertex of $H$ is drawn by $\Gamma$ in a vortex of $\rho$.
A curve in $\Sigma$ is called \emph{$\rho$-aligned} if it only intersects $\Gamma$ in nodes.
In turn, a disc in $\Sigma$ is called \emph{$\rho$-aligned} if its boundary is $\rho$-aligned and we call it \emph{vortex-free} if its interior is disjoint from all vortices in $\rho$.
If $H$ is planar, we say that it is \emph{flat in $\rho$} if there exists a vortex-free, $\rho$-aligned disc $\Delta \subseteq \Sigma$ which contains all cells $c \in C(\rho)$ with $E(\sigma(c)) \cap E(H) \neq \emptyset$.

For any $\rho$-aligned disc $\Delta,$ we call the subgraph of $G$ that is the union of $\sigma(c)$ for all cells $c$ of $\rho$ in $\Delta$ the \emph{crop of $G$ by $\Delta$ (in $\rho$)}.
Furthermore, the \emph{restriction $\delta'$ of $\rho$ by $\Delta$} is defined as the $\Delta$-rendition that consists of the restriction of both $\Gamma,$ $\sigma,$ and $\pi$ to $\Delta$.

This allows to define a society associated with $\Delta$ as follows.
Let $V(\Omega_{\Delta})$ be the set of all vertices whose corresponding nodes are drawn in the boundary of $\Delta$ and let $\Omega_{\Delta}$ be the cyclic ordering of $V(\Omega_{\Delta})$ obtained by traversing along the boundary of $\Delta$ in the anticlockwise direction.
Now, let $G_{\Delta}$ be the crop of $G$ by $\Delta$.
We call the society $(G_{\Delta}, \Omega_{\Delta})$ the \emph{$\Delta$-society (in $\rho$)}.
If $\rho$ is clear from the context, we do not mention it.
We also call the restriction of $\rho$ by $\Delta,$ the \emph{restriction of $\rho$ to $(G_{\Delta}, \Omega_{\Delta})$.}

Let $\mathcal{P}$ be an $X$-$Y$-linkage in $G$ such that $X \cap V(G_{\Delta}) \subseteq V(\Omega_{\Delta})$ and $Y \subseteq V(G_{\Delta})$ and assume that each path in $\mathcal{P}$ is grounded in $\rho$.
Then we define the \emph{$\Delta$-truncation (in $\rho$)} of $\mathcal{P}$ to be the $V(\Omega_{\Delta})$-$Y$-linkage in $G_{\Delta}$ which consists of the minimal $V(\Omega_{\Delta})$-$Y$-subpaths of the paths in $\mathcal{P}$.

Let $\rho$ be a rendition of a society $(G, \Omega)$ in a surface $\Sigma$ that is not the sphere.
Given a cycle $C \subseteq G$ that is grounded in $\rho$ we define the \emph{$C$-disc (in $\rho$)} as the unique $\rho$-aligned disc $\Delta \subseteq \Sigma$ bounded by the trace of $C$ in $\rho$.
We also use the terms \emph{$C$-society (in $\rho$)} to denote the $\Delta$-society in $\rho$ and \emph{$C$-truncation (in $\rho$)} to denote the $\Delta$-truncation in $\rho$ of an appropriately defined linkage in $G$.

\paragraph{Transactions and cylindrical renditions.}
Let $(G, \Omega)$ be a society. 
A \emph{transaction} in $(G, \Omega)$ is an $A$-$B$-linkage for disjoint segments $A, B$ of $\Omega$ consisting of $V(\Omega)$-paths.
A rendition of $(G, \Omega)$ in the disc with a unique vortex $c_{0}$ is called a \emph{cylindrical rendition} of $(G, \Omega)$ \emph{around} $c_{0}$.

\paragraph{Nests and radial linkages.}
Let $\rho$ be a rendition of a society $(G, \Omega)$ in a surface $\Sigma$ other than the sphere.
Further, let $\mathcal{C} = \langle C_1,\ldots,C_s\rangle$ be a sequence of pairwise disjoint cycles in $G,$ such that each of them is grounded in $\rho$ and we have $\Delta_1 \subsetneq \ldots \subsetneq \Delta_s \subseteq \Sigma,$ where $\Delta_i$ is the $C_i$-disc for each $i \in [s]$.
We say that the cycles in $\mathcal{C}$ are \emph{concentric (in $\rho$)} and call $C_1$ the \emph{inner} cycle and $C_s$ the \emph{outer cycle} of $\mathcal{C}$ in $\rho$. %

A \emph{nest $\mathcal{C}$ (in $\rho$)} is a set of concentric cycles with the inner cycle $C$ such that every vortex of $\rho$ is found in the $C$-disc.
Moreover, we call a $V(\Omega)$-$V(C)$-linkage $\mathcal{R}$ a \emph{radial linkage (in $\rho$) for $\mathcal{C}$} if all paths in $\mathcal{R}$ are grounded in $\rho$ and internally disjoint from $V(\Omega)$.

If $(G, \Omega)$ is a society with a nest $\mathcal{C}$ in a rendition $\rho$ of $(G, \Omega)$ in a disc, we say that a radial linkage $\mathcal{R}$ for $\mathcal{C}$ is \emph{orthogonal to $\mathcal{C}$} if for all $C \in \mathcal{C}$ and all $R \in \mathcal{R}$ the graph $C \cap R$ is a path.

\paragraph{Inner and outer graphs of a cycle.}
Let $(G, \Omega)$ be a society with a $\Sigma$-rendition $\rho.$
Further, let $C$ be a grounded cycle whose trace bounds a disc $\Delta_C$ and the $\Delta_C$-society $(G', \Omega').$
We call $G' \cup C$ the \emph{inner graph of $C$ (in $\rho$)} and call $G'$ itself the \emph{proper inner graph of $C$ (in $\rho$).}
Let $B = \pi(N(\rho) \cap \mathsf{bd}(\Delta_C)).$
We define the \emph{proper outer graph of $C$ (in $\rho$)} as $G'' \coloneqq G[B \cup (V(G) \setminus V(G'))]$ and call $G'' \cup C$ the \emph{outer graph of $C$ (in $\rho$).} 

\paragraph{Depth of vortices.}
Let $G$ be a graph and $\rho$ be a $\Sigma$-rendition of $G$ with a vortex cell $c_0.$
Notice that $c_0$ defines a society $(\sigma(c_0), \Omega_{c_0}),$ where $V(\Omega_{c_0})$ is the set of vertices of $G$ corresponding to $N_\rho(c_0).$
The ordering $\Omega_{c_{0}}$ is obtained by traversing along the boundary of the closure of $c_0$ in counterclockwise direction.
We call $(\sigma(c_0),\Omega_{c_0})$ as obtained above the \emph{vortex society} of $c_0.$

We define the \emph{depth} of a society $(G, \Omega)$ as the maximum cardinality of a transaction in $(G, \Omega)$.
The \emph{depth} of the vortex $c_{0}$ is thereby defined as the depth of its vortex society.

It easily follows that every society with a linear decomposition of adhesion at most $k$ has depth at most $2k.$
What interests us here is the reverse statement which was shown by Robertson and Seymour in \cite{RobertsonS1990Graphb}.

\begin{proposition}[{\cite{RobertsonS1990Graphb}}]\label{prop:DepthToLinearDec} Let $k$ be a non-negative integer and $(G, \Omega)$ be a society of depth at most $k.$
Then $(G, \Omega)$ has a linear decomposition of adhesion at most $k.$
\end{proposition}

Finally, we define the \emph{width} of a $\Sigma$-rendition of a graph $G$ as the minimum non-negative integer $w$ such that every vortex society of $\rho$ admits a linear decomposition of width at most $w.$
Further, we define the \emph{breadth of $\rho$} as the number of vortex cells of $\rho$ and the \emph{depth of $\rho$} as the maximum depth of its vortex societies.

\paragraph{Flat Meshes.}
Let $n \ge 2$ be an integer.
Let $G$ be a graph, and let $M \subseteq G$ be an $n$-mesh.
We say that $M$ is a \emph{flat mesh} in $G$ if there exists a rendition $\rho$ of $G$ on the sphere with a single vortex $c_0$ such that $M$ is flat in $\rho$ and the trace of the perimeter of $M$ in $\rho$ separates all vertices\footnote{Note that this may mean that non-nodes belonging to the perimeter of the mesh may be drawn in cells outside of the vortex.} in $N(\rho) \cap V(M)$ from $c_0$. 
We say that $\rho$ \emph{witnesses} the flatness of $M.$

\paragraph{Cylindrical meshes.}
In some settings it will be easier to work with a version of meshes that is cylindrical.

Let $m,n$ be positive integers, let $M$ be a graph, and let $C_1, \ldots, C_m$ be cycles and $P_1, \ldots, P_n$ be paths in $M$ such that the following holds for all $i \in [m]$ and $j \in [n]$:
\begin{itemize}
 \item $C_1, \ldots, C_m$ are pairwise vertex-disjoint, $P_1, \ldots, P_n$ are pairwise vertex-disjoint, and $M = C_1 \cup \cdots \cup C_m \cup P_1 \cup \cdots \cup P_n$.

 \item $C_i \cap P_j$ is a path, and if $i \in \{ 1, m \}$ or $j \in \{ 1, n \},$ then $C_i \cap P_j$ has exactly one vertex,

 \item when traversing $C_i$ starting from an endpoint of $P_1\cap C_i,$ then either the paths $P_1, \ldots, P_n$ are encountered in the order listed or the next $P_j$ one encounters is $P_n$ and from here the paths are encountered in the order $P_n,\dots,P_1,$ and

 \item $P_j$ has one end in $C_1$ and the other in $C_m,$ and when traversing $P_j$ starting from its endpoint on $C_1,$ the cycles $C_1, \ldots, C_m$ are encountered in the order listed.
\end{itemize}
If the above conditions hold for $M,$ we call $M$ an \emph{$(n \times m)$-cylindrical mesh}.
The cycles $C_1, \ldots, C_m$ are called the \emph{concentric cycles}, or \emph{cycles}, of $M$ and the paths $P_1, \ldots, P_n$ are called the \emph{radial paths}, or \emph{rails}, of $M$.
We also call $(n \times n)$-cylindrical meshes \emph{$n$-cylindrical meshes}.

Let $n, m \geq 3$ be integers, $G$ be a graph, and $M$ be an $(n \times m)$-cylindrical mesh in $G$.
Similarly to meshes, let $\mathcal{T}_{M}$ be the orientation of $\mathcal{S}_{r},$ where $r = \min \{ n, m \},$ such that for every $(A, B) \in \mathcal{T}_{M},$ the set $B \setminus A$ contains the vertex set of both a concentric cycle and a radial path of $M,$ we call $B$ the \emph{$M$-majority side of $(A, B)$.}
Then $\mathcal{T}_{M}$ is the tangle \emph{induced} by $M$.
If $\mathcal{T}$ is a tangle in $G,$ say that $\mathcal{T}$ controls the mesh $M$ if $\mathcal{T}_{M}$ is a truncation of $\mathcal{T}$.

\subsection{Manipulating vital linkages}
At several points in our proofs we will want to delete vertices, edges, and decompose graphs at separations while preserving the existence of vital linkages.
Most of these operations have been considered already by Robertson and Seymour in \cite{RobertsonS2009Graph} and we will reuse several of their results on these operations.

\begin{proposition}\label{lem:vital_and_separations}
 Let $(G,T)$ be an annotated graph and $\LLL$ a vital $T$-linkage in $G$.
 Let $(A,B)$ be a separation in $G$ and let $S \coloneqq A \cap B$.
 Let $T_B \coloneqq T \cap B$.
 Then $(G[B],T_B \cup S)$ is an annotated graph admitting a vital $(S\cup T_B)$-linkage $\LLL'$ such that each path in $\LLL'$ is a subpath of a path in~$\LLL$.
\end{proposition}

\begin{proposition}\label{lem:vital_and_deletion}
 Let $(G,T)$ be an annotated graph and $\LLL$ a vital $T$-linkage in $G$.
 Let $v \in V(G)\setminus V(T)$ and $G' \coloneqq G - v$.
 Let $u,w \in V(G)$ be the two vertices adjacent to $v$ such that $uv,vw \in E(L)$ for some $L \in \LLL$.
 Let $T' \coloneqq T \cup \{u,w\}$.
 Then $(G',T')$ is an annotated graph admitting a vital $T'$-linkage $\LLL'$ such that each path in $\LLL'$ is a subpath of a path in $\LLL$.

 If $v \in V(G)$ is an isolated vertex, then $v \in T$ and there is some path $L \in \LLL$ consisting solely of~$v$.
 Further, $(G-v,T\setminus\{v\})$ is an annotated graph with a vital linkage $\LLL'$ such that $\LLL' \subseteq \LLL$.
\end{proposition}

Let $G, H$ be graphs with $H \subseteq G$ and let $\mathcal{L}$ be a linkage in $G,$ then we define~$\mathcal{L} \cap H$ to be $\{ P ~\colon~ L \in \mathcal{L} \text{ and } P \text{ is a component of } H \cap L \}$.
The following is an immediate consequence of (2.3) in \cite{RobertsonS2009Graph}\footnote{Note that in Robertson and Seymour's setting separations are defined as tuples of subgraphs that intersect only in vertices.}.

\begin{proposition}\label{prop:vitalitypreservedthroughrestriction}
 Let $G, H$ be graphs with $H \subseteq G$ and let $\mathcal{L}$ be a linkage in $G$.
 If $\mathcal{L}$ is vital in $G,$ then $\mathcal{L} \cap H$ is vital in $H$.
\end{proposition}

A key tool for us will be \cite[Theorem 2.5]{RobertsonS2009Graph}.
This is essential for our later proof of \zcref{thm:no_apices_vital_reduction} analysing the number of ``non-trivial types'' a vital linkage can ``induce'' on the separators given a sequence of laminar separations.
The following is a transcription to our setting.
\begin{proposition}[\!{\cite[Theorem 2.5]{RobertsonS2009Graph}}]\label{thm:Lemma2.5}
 Let $k,d \in \N$ and let $f(k,d) \coloneqq (2k + 1)(k + d + 1)^{2(k+d)} + 1$. Let $\LLL$ be a vital linkage of order $k$ in a graph $G,$ and for $i \in [f(k,d)]$ let $(A_i, B_i)$ be a separation of $G$ of order $d,$ such that
 \begin{enumerate}
 \item $A_i \subseteq A_j$ and $B_j \subseteq B_i$ for $1 \leq i < j\leq n,$ and
 \item for $1\leq i < n,$ there is a linkage $\mathcal{M}_i$ in $G[B_i] \cap G[A_{i+1}]$ of order $d,$ each path of which has one end in $A_i \cap B_i$ and the other in $A_{i+1} \cap B_{i+1}$.
 \end{enumerate} 
 Then there exists $i \in [n]$ such that $\bigcup\LLL \cap G[B_i] \cap G[A_{i+1}] = \bigcup \mathcal{M}_i$.
\end{proposition}

\subsection{Colorful graphs}
To later deal with the folio-problem, we will need definitions for graphs which feature several different colours.
We define these here since we will briefly deal with these graphs and the folio-problem in \zcref{sec:largeclique} as well.
We use the patently incorrect spelling of colour when defining these \textsl{colorful graphs} to stay in line with the name given to them in \cite{ProtopapasTW2025Colorful}.

Let $q$ be a non-negative integer.
A \emph{$q$-colorful graph} is a pair $(G,\chi)$ where $G$ is a graph and $\chi\colon V(G)\rightarrow 2^{[q]}$.
For a vertex $v\in V(G)$ we call $\chi(v)$ the \emph{palette} of $v$.
A $q$-colorful graph $(G,\chi)$ where $\chi(G)\subsetneq[q]$ is called \emph{restricted} and a $q$-colorful graph $(G,\chi)$ with $\chi(G) = \emptyset$ is called \emph{blank}.
A $q$-colorful graph $(G,\chi)$ is said to be \emph{rainbow} if $\chi(v)=[q]$ for all $v\in V(G)$. 
Annotated graphs correspond directly to $1$-colorful graphs.

Let $(G,\chi)$ and $(H,\psi)$ be colorful graphs.
We say that a minor model $\varphi$ of $H$ in $G$ is a \emph{colorful minor model} of $(H,\psi)$ in $(G,\chi)$ if for every $v\in V(H)$ and every $i\in \psi(v)$ there is some vertex $u\in \varphi(v)$ such that $i\in\chi(u)$.
To simplify notation, we may simply say that $(H,\psi)$ is a \emph{colorful minor} of $(G,\chi)$.
Accordingly, in a $1$-colorful graph $G$ that has $H$ as a red minor, $H$ is also a rainbow (colorful) minor.

\subsection{The Local Structure Theorem for annotated graphs}

One of the key results we will use is a variant of the Local Structure Theorem (first proven in \cite{RobertsonS2003Grapha} as (3.1)) recently provided by Gorsky, Protopapas, and Wiederrecht \cite{GorskyPW2026Quickly}.
Their statement of this theorem (see Theorem 5.1 in \cite{GorskyPW2026Quickly}) is quite a bit more technical than what we use here.
We simplify their formulation in three ways.
First, we assume that we are working on an annotated graph with bounded dimensionality, which gets rid of one of the options in their statement.
Secondly, we will not require the complicated structure that their \textsl{layout}-definition provides and thus will simply extract the parts of their structure that are of use to us.
(Accordingly the proposition we state here is quite a bit weaker than what they actually prove.)
Third, we reformulate their first outcome by digging only mildly into their proof (which uses \zcref{thm:colorfulclique} to guarantee this point) to respect a blank $K_t$-minor model instead of the mesh.

\begin{proposition}[Gorsky, Protopapas, and Wiederrecht \cite{GorskyPW2026Quickly}]\label{thm:localstructure}
There exist functions $\mathsf{apex}_{\ref{thm:localstructure}},\mathsf{depth}_{\ref{thm:localstructure}} \colon \mathbb{N}^2 \to \mathbb{N}$ and $\mathsf{mesh}_{\ref{thm:localstructure}} \colon \mathbb{N}^3 \to \mathbb{N}$ such that for all integers $t \geq 1,$ $r \geq 4,$ and $w \geq 3,$ every annotated graph $(G,R)$ with $\mathsf{bidim}(G,T) \leq r,$ and every $\mathsf{mesh}_{\ref{thm:localstructure}}(r, t, w)$-mesh $M \subseteq G$ one of the following holds.
\begin{enumerate}

 \item $(G,R)$ has a blank $K_{\lfloor \nicefrac{5t}{2} \rfloor}$-minor model $\psi$ and a separation $(X,Y)$ of order at most $t-1$ such that the $\mathcal{T}_\psi$-big component of $G - (X \cap Y)$ is blank,
 
 \item $(G, R)$ has a red $K_t$-minor model controlled by $M,$ or
 \item there exists a set $A \subseteq V(G)$ with $|A| \leq \mathsf{apex}_{\ref{thm:localstructure}}(t, r)$ and a surface $\Sigma$ of genus less than $9t^2$ such that $(G - A, R \setminus A)$ has a blank $\Sigma$-rendition $\rho$ with breadth at most $\nicefrac{3}{2}(t-1)(3t-4) + r(r-1)-3,$ depth at most $\mathsf{depth}_{\ref{thm:localstructure}}(t, r),$ and there exists a vortex-free, $\rho$-aligned disc $\Delta \subseteq \Sigma$ such that the restriction of $\rho$ to $\Delta$ contains a flat $w$-wall.
\end{enumerate}
Moreover, it holds that $\mathsf{apex}_{\ref{thm:localstructure}}(t,r), \mathsf{depth}_{\ref{thm:localstructure}}(t,r) \in \mathbf{poly}(t+r)$ and $\mathsf{mesh}_{\ref{thm:localstructure}}(t,r,w) \in \mathbf{poly}(t+r) + \mathbf{O}(t^2w)$.
\end{proposition}

\section{Taming linkages in surfaces}\label{sec:linkagesinsurfaces}

This subsection develops the structural framework needed to control how linkages interact with nested cycles in surface embeddings. The central objective is to enforce a controlled interaction between a family of paths and a system of concentric cycles drawn in a disc, ensuring that paths penetrate the cycles in a well-organised and hierarchical manner. This machinery has been inspired by the result of Frédéric Mazoit in \cite{Mazoit2013Single} 
who proved a single exponential bound 
for the linkage function on surface embeddable graphs. This result is based on an earlier result of Geelen, Huynh, and Richter in \cite{GeelenHR2018Explicit}
on topological linkages. The treatment presented here can be viewed as a streamlined reconstruction of the ideas of 
\cite{Mazoit2013Single}, adapted to our setting and developed with an emphasis on clarity and structural transparency.

\subsection{Preparing the linkage}
\label{sec:linkagepreparation}

This subsection is dedicated to establishing that we can force a set of paths interacting with a set of concentric cycles on a disc to intersect in an orderly manner.
What we present here is a version of a fairly standard approach (see for example \textsl{cheap linkages} in \cite{AdlerKKLST2011Tight} and \textsl{valleys} in \cite{Mazoit2013Single}) adapted to our setting.

\begin{figure}[ht]
 \centering
 \scalebox{0.85}{
 \begin{tikzpicture}[scale=1.25]

 \pgfdeclarelayer{background}
		 \pgfdeclarelayer{foreground}
			
		 \pgfsetlayers{background,main,foreground}

 \begin{pgfonlayer}{background}
 \pgftext{\includegraphics[width=6cm]{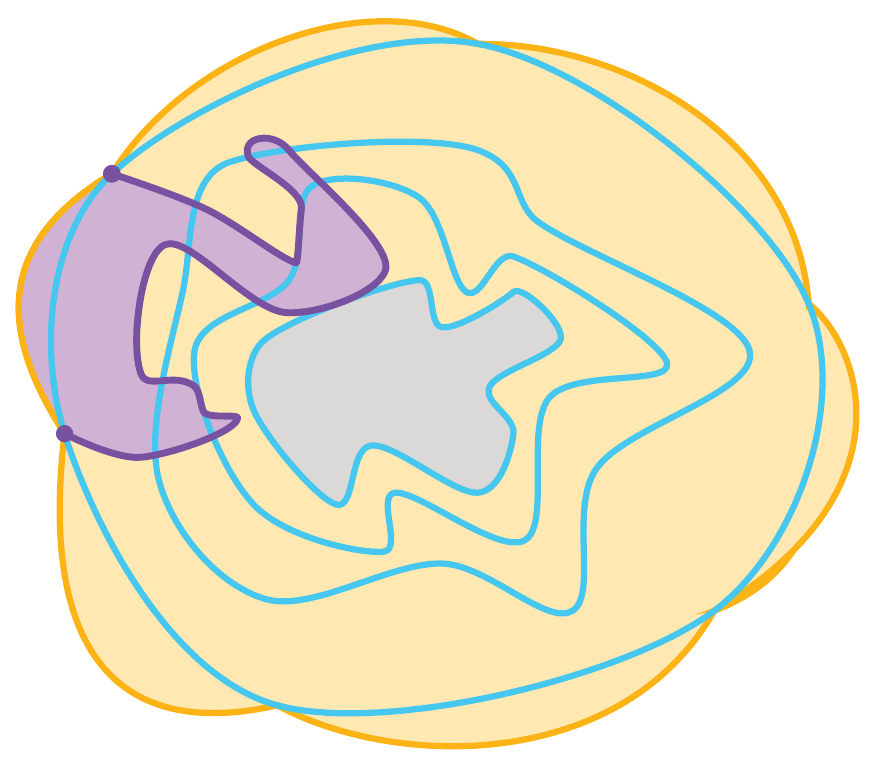}} at (C.center);
 \end{pgfonlayer}{background}
			
 \begin{pgfonlayer}{foreground}

 \node (THEPATH) at (-2.2,-0.2) [draw=none] {$P$};
 \node (DELTA) at (3.05,-0.5) [draw=none] {$\Delta$};
 \node (DISK1) at (-0.2,0) [draw=none] {$\Delta_1$};
 \node (DISKP) at (-2.35,0.5) [draw=none] {$\Delta_P$};
 \node (C1) at (0.55,0.35) [draw=none] {$C_1$};
 \node (C2) at (1.125,0.225) [draw=none] {$C_2$};
 \node (C3) at (1.7,-0.5) [draw=none] {$C_3$};
 \node (C4) at (2.25,0.6) [draw=none] {$C_4$};
 
 \end{pgfonlayer}{foreground}
 
 \begin{pgfonlayer}{foreground}
 \end{pgfonlayer}{foreground}
 
 \end{tikzpicture}
 }
 \caption{A 4-well with a single path drawn in violet. The disc $\Delta_1 \subseteq \Delta$ bounded by $C_1$ forms the ``bottom'' of the well and is illustrated in gray. As drawn, it is possible for the outermost cycle to intersect the boundary of the disc onto which the well is embedded in nodes.}
 \label{fig:wellexample}
\end{figure}
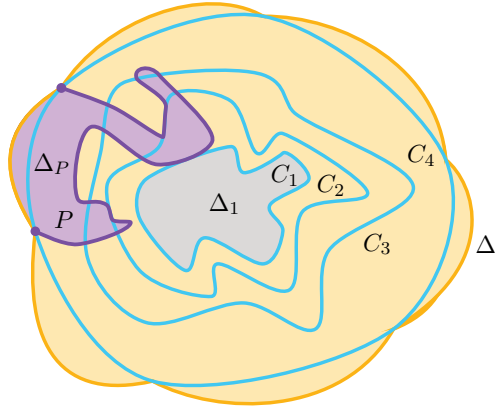

Let $s$ be a positive integer.
An \emph{$s$-well} $\mathcal{W} = (G, \Omega, \rho, \mathcal{C}, \mathcal{P})$ is a society $(G,\Omega)$ with a vortex-free, cylindrical rendition $\rho$ in a disc $\Delta$ such that $G = \bigcup \mathcal{C} \cup \bigcup \mathcal{P}$ for a set of pairwise internally disjoint $V(\Omega)$-paths $\mathcal{P}$ and $\mathcal{C} = \langle C_1, \ldots, C_s\rangle$ is a sequence of concentric cycles in $\Delta$.
Thus for each $P \in \mathcal{P},$ we may let $T$ be the trace of $P$ in $\rho$ and there exists a unique disc $\Delta_P$ in $\Delta - T$ such that the intersection of $\Delta_P$ and the trace of $C_1$ is empty.
We call the closure of $\Delta_P$ the \emph{interior of $P$}.
See \zcref{fig:wellexample} for an example of a 4-well.

We say that an $s$-well $\mathcal{W} = (G, \Omega, \rho, \mathcal{C}, \mathcal{P})$ is \emph{drained} if either $|\mathcal{P}| \leq 1$ or
\begin{itemize}
 \item for each $P \in \mathcal{P}$ we have that for all $i \in [s-1],$ if $P \cap C_i$ is non-empty, then there exists $Q \in \mathcal{P} \setminus \{ P \}$ such that the trace of $Q \cap C_{i+1}$ is non-empty and found in $\Delta_P,$ and
 
 \item for all $i \in [s],$ no $C_i$-path $P'$ in $P$ that is disjoint from $V(\mathcal{C} \setminus \{ C_i \})$ is contained in a cell $c \in C(\rho)$ such that $\sigma(c)$ contains an edge of $C_i$.
\end{itemize}

In drained wells the vortex-free rendition tells us that any path that moves deeply into a drained well implies the existence of many other paths that also go deep.
This is captured in the following observation which we will apply implicitly whenever we work with drained wells.

\begin{observation}\label{obs:fullbucketsindrainedwells}
 Let $\mathcal{W} = (G, \Omega, \rho, \mathcal{C}, \mathcal{P})$ be a drained $s$-well and let $P \in \mathcal{P}$ be a path such that $C_i \cap P$ is non-empty for some $i \in [s]$.
 Then there exist $(s - i) + 1$ distinct paths $P_s, P_{s-1}, \ldots, P_{i+1}, P_i = P$ such that $P_j \cap C_j$ is non-empty for all $j \in [i,s]$ and we have $\Delta_{P_s} \subsetneq \ldots \subsetneq \Delta_{P_{i+1}} $.
\end{observation}

We start by proving that we can drain any given well by rearranging the paths within it without changing the underlying set of cycles that make up the well.

\begin{lemma}\label{lem:drainedwell}
 Let $s$ be a positive integer and let $\mathcal{W} = (G, \Omega, \rho, \mathcal{C}, \mathcal{P})$ be an $s$-well.
 Then there exists an $s$-well $\mathcal{W}' = (G', \Omega, \rho', \mathcal{C}, \mathcal{P}')$ that is drained and such that $G' \subseteq G,$ $\rho'$ is the rendition of $G'$ induced by $\rho,$ $|\mathcal{P}| = |\mathcal{P}'|,$ and for each $P \in \mathcal{P}$ there exists a $P' \in \mathcal{P}'$ with the same endpoints.
\end{lemma}
\begin{proof} 
 Suppose the statement is false and that $\mathcal{W}$ is a minimal counterexample with respect to $|E(\mathcal{C} \cup \mathcal{P})|$.
 Let $\mathcal{C} = \langle C_1, \ldots, C_s\rangle,$ where $\Delta_i$ is the $C_i$-disc in $\rho$ for each $i \in [s],$ and let $P \in \mathcal{P}$ intersect some cycle $C_i$ with $i \in [s-1]$.

 We first establish the second property in the definition of being drained.
 Observe that if there exists a cell $c \in C(\rho)$ such that $\sigma(c)$ contains both an edge of $C_i$ and a $C_i$-path $P' \subseteq P$ with endpoints $u,v,$ then we may simply replace $P'$ in $P$ with the unique $u$-$v$-path in $C_i \cap \sigma(c)$.
 This would contradict the minimality of $\mathcal{W}$.
 
 Towards the first property, suppose that $C_{i+1}$ is not intersected by any path whose trace lies in $\Delta_P$.
 Since $\rho$ is vortex-free, we know that $C_{i+1} \cap P$ is non-empty.
 Then let $a,b$ be the first and last vertex we see of $P$ within $C_{i+1} \cap P$ when traversing $P$ from one endpoint to the other.

 If there exists a grounded $a$-$V(P)$-path $L$ in $C_{i+1} - E(P)$ whose trace is found in $\Delta_P$ (see the left side of \zcref{fig:drainingthewell}), then according to our assumptions on $P$ and $C_{i+1},$ no other path in $\mathcal{P}$ whose trace lies in $\Delta_P$ intersects $L$.
 However, this would allow us to shortcut $P$ along $C_{i+1},$ contradicting the minimality of $\mathcal{W}$. 

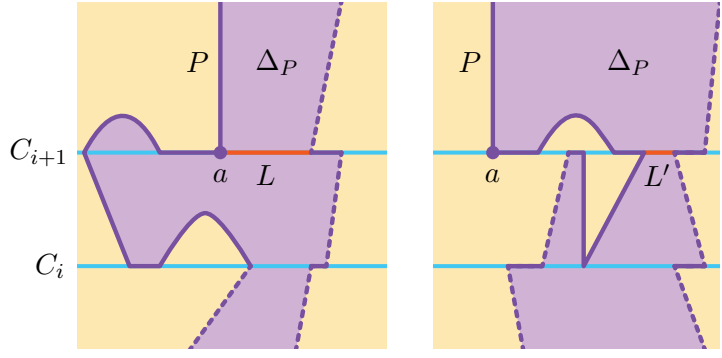
\begin{figure}[ht]
 \centering
 \begin{tikzpicture}[scale=1.5]

 \pgfdeclarelayer{background}
		 \pgfdeclarelayer{foreground}
			
		 \pgfsetlayers{background,main,foreground}

 \begin{pgfonlayer}{background}
 \pgftext{\includegraphics[width=6cm]{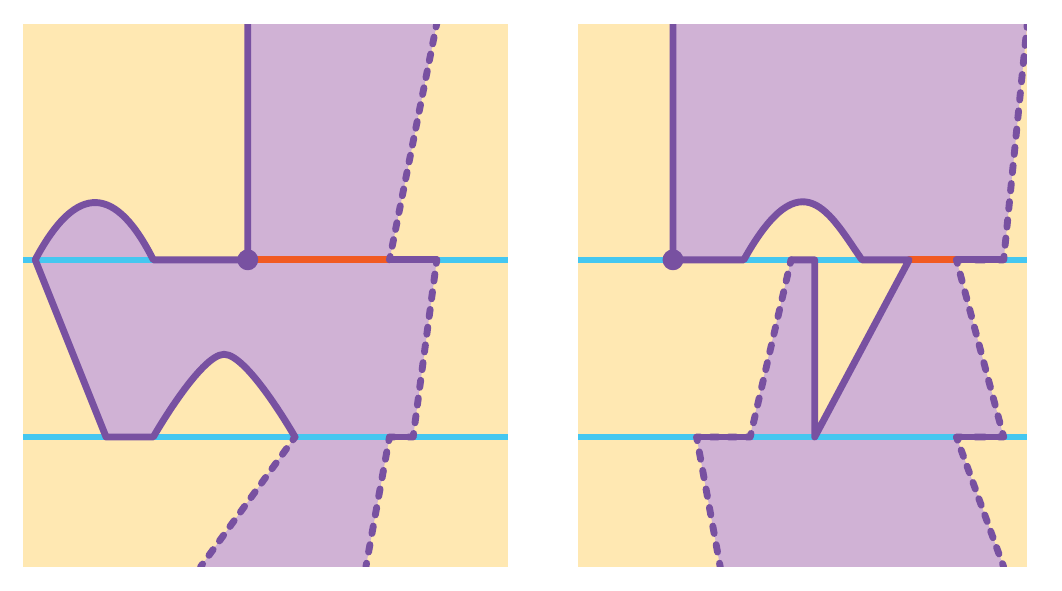}} at (C.center);
 \end{pgfonlayer}{background}
			
 \begin{pgfonlayer}{foreground}

 \node (THEPATH) at (-1.8,1) [draw=none] {$P$};
 \node (THEPATHCOPY) at (0.6,1) [draw=none] {$P$};
 \node (A) at (-1.6,0) [draw=none] {$a$};
 \node (ACOPY) at (0.8,0) [draw=none] {$a$};
 \node (DISKP) at (-1.1,1) [draw=none] {$\Delta_P$};
 \node (DISKPCOPY) at (2,1) [draw=none] {$\Delta_P$};
 \node (CI) at (-3.1,-0.8) [draw=none] {$C_i$};
 \node (CIPLUS) at (-3.2,0.2) [draw=none] {$C_{i+1}$};
 \node (L) at (-1.2,0) [draw=none] {$L$};
 \node (LPRIME) at (2.25,0) [draw=none] {$L'$};
 
 \end{pgfonlayer}{foreground}
 
 \begin{pgfonlayer}{foreground}
 \end{pgfonlayer}{foreground}
 
 \end{tikzpicture}
 \caption{An illustration for \zcref{lem:drainedwell} showing the two kinds of paths $L$ and $L'$ that lead to a contradiction to minimality in the context of the proof of the lemma.}
 \label{fig:drainingthewell}
\end{figure}

 Let $P^* \subseteq P$ be the first $C_{i+1}$-$C_i$-path or $C_{i+1}$-path with a trace in $\Delta_{i+1}$ that we see when traversing $P$ starting in $a$ and heading towards $b$.
 Note that such a path must exist, since $P \cap C_i$ is non-empty and it must be grounded, according to our prior observation.
 Let $a'$ be the endpoints of $P^*$ closest to $a$ on $P$.

 As there does not exist a grounded $a$-$V(P)$-path in $C_{i+1} - E(P)$ whose trace is found in $\Delta_P$ and due to our choice of $P^*$ (see the right side of \zcref{fig:drainingthewell}), there now exists a grounded $a'$-$V(P)$-path $L'$ in $C_{i+1} - E(P)$ whose trace is found in $\Delta_P$ and we arrive at a contradiction exactly as we did earlier with $L$.
\end{proof}

For some applications the above lemma will be sufficient, but we will also need wells that act much more orderly but provide stronger assumptions on the structure of the cycles.

Let $\mathcal{C} = \langle C_1, \ldots, C_s\rangle$ be a sequence of concentric cycles in a $\Delta$-rendition $\rho$ of a graph $G$ in a disc $\Delta$ and let $\Delta_i$ be the $C_i$-disc in $\rho$ for each $i \in [s],$ with $\Delta_0 = \emptyset$.
We say that $\mathcal{C}$ is \emph{tight (in $\rho$)} if for all $i \in [s]$ there does not exist a grounded $C_i$-path whose trace is found in $\Delta_i \setminus \Delta_{i-1}$. 
Given a vortex $c$ with $c \subseteq \Delta_1,$ we say that $\mathcal{C}$ is \emph{tight around $c$ (in $\rho$)} if there does not exist a grounded $C_1$-path whose trace is found in $\Delta_1$.
In general, we call a set of concentric cycles $\mathcal{C}'$ tight in a rendition $\rho'$ (around a vortex) if we can find a disc such that the restriction of $\rho'$ to $\Delta'$ allows us to call $\mathcal{C}'$ tight.
A well $\mathcal{W} = (G', \Omega, \rho', \mathcal{C}', \mathcal{P})$ \emph{tight} if $\mathcal{C}'$ is tight in $\rho'$.

Finding tight concentric cycles is an easy exercise and many articles in the literature provide algorithmic proofs that turn a concentric set of cycles into a tight, concentric set of cycles (see the discussion of tight cycles in \cite{AdlerKKLST2011Tight} and cozy nests in \cite{GorskySW2025Polynomial} for example).
Thus we state the following result without proof.

\begin{proposition}\label{prop:makecyclestight}
 Let $\mathcal{C} = \langle C_1, \ldots, C_s\rangle$ be a set of concentric cycles in a $\Delta$-rendition $\rho$ of a graph $G$ in a disc $\Delta$ with $\Delta_s$ being the $C_s$-disc in $\rho$.
 Then there exists a tight concentric set of cycles $\mathcal{C}' = \langle C_1', \ldots, C_s'\rangle$ such that the trace of each cycle in $\mathcal{C}'$ lies in $\Delta_s$.
\end{proposition}

Under the assumption that we are working with a tight well, we can actually guarantee a much stronger property than simply being drained.
We say that a well $\mathcal{W} = (G, \Omega, \rho, \mathcal{C} = \langle C_1, \ldots, C_s\rangle, \mathcal{P})$ is \emph{dry} if $\mathcal{W}$ is drained and for all $P \in \mathcal{P}$ we have that
\begin{enumerate}
 \item the graph $C_1 \cap P$ is either empty or a single path,

 \item there exists a unique $i \in [s]$ such that $C_i \cap P$ is a single path,

 \item for all $j \in [i-1]$ the graph $C_j \cap P$ is empty, and

 \item for all $j \in [i+1,s]$ the graph $C_j \cap P$ consists of exactly two paths.
\end{enumerate}

It turns out that guaranteeing the existence of a dry well once given a tight well is quite easy, since we already proved \zcref{lem:drainedwell}.

\begin{lemma}\label{lem:drywell}
 Let $s$ be a positive integer and let $\mathcal{W} = (G, \Omega, \rho, \mathcal{C}, \mathcal{P})$ be a tight $s$-well.
 Then there exists a tight, dry $s$-well $\mathcal{W}' = (G', \Omega, \rho', \mathcal{C}, \mathcal{P}')$ such that $G' \subseteq G,$ $\rho'$ is the rendition of $G'$ induced by $\rho,$ $|\mathcal{P}| = |\mathcal{P}'|,$ and for each $P \in \mathcal{P}$ there exists a $P' \in \mathcal{P}'$ with the same endpoints.
\end{lemma}
\begin{proof} 
 Suppose the statement is false and that $\mathcal{W}$ is a minimal counterexample with respect to $|E(\mathcal{C} \cup \mathcal{P})|$.
 Thanks to \zcref{lem:drainedwell}, we may assume that $\mathcal{W}$ is already drained and we let $\mathcal{C} = \langle C_1, \ldots, C_s \rangle,$ where $\Delta_i$ is the $C_i$-disc in $\rho$ for each $i \in [s]$.

 Thanks to the second property of being drained, all paths we discuss from this point forward can be assumed to have a trace.
 Furthermore, we observe that, due to the tightness of $\mathcal{C},$ no path in $\mathcal{P}$ may have a $C_1$-path whose trace intersects the interior of $\Delta_1$.

 Suppose now that there exists a path $P \in \mathcal{P}$ and an $i \in [s]$ such that $P$ contains a $C_i$-path $P'$ with the endpoints $a,b$ whose trace is found outside of $\Delta_i \setminus \mathsf{bd}(\Delta_i)$.
 Then we choose $P$ to maximise $i$ and let $L$ be the $a$-$b$-path in $C_i$ such that the $P' \cup L$-disc is disjoint from $\Delta_i \setminus \mathsf{bd}(\Delta_i)$.
 Due to the maximality of $i$ and the fact that $\rho$ is vortex-free, we conclude that $L$ is internally disjoint from $V(\mathcal{P})$.
 Thus by replacing $P'$ with $L$ in $P$ we produce a contradiction to the minimality of $\mathcal{W}$.

 As a consequence of the observations we have just made, for any path $P \in \mathcal{P}$ such that $C_1 \cap P$ is non-empty, the graph $C_1 \cap P$ must be a path.
 This guarantees the first property of being dry.

 Suppose now that there instead exists a path $P \in \mathcal{P}$ such that for no $i \in [s]$ the graph $C_i \cap P$ consists of a single component.
 We choose $P$ such that the value $i \in [s]$ for which $C_i \cap P$ is non-empty is maximum and for all $Q \in \mathcal{P}$ such that the interior $\Delta_Q$ of $Q$ contains $P$ there does actually exist a $j \in [s]$ such that $C_j \cap Q$ is a path.
 Clearly, the fact that $C_i \cap P$ is comprised of at least two components implies the existence of a $C_i$-path $P'$ in $P$.
 As we have observed earlier, $P'$ must have its trace in $\Delta_i$.
 But since $i$ is maximum and $P$ is chosen to be the ``outermost'' path violating the second property at $i,$ $P'$ must be disjoint from all other cycles in $\mathcal{C}$ and thus the tightness of $\mathcal{C}$ prohibits its existence.
 This ensures that the second property of dryness holds.

 The third property is then easily derived from the fact that $\mathcal{C}$ is concentric and $\rho$ is vortex-free.
 So we may finally suppose that the fourth property is violated and that there exists a $P \in \mathcal{P}$ such that $C_i \cap P$ consists of more than two components for some $i \in [s]$.
 By again choosing $P$ and $i$ to be maximal as in the previous paragraph, we arrive at another contradiction, completing our proof.
\end{proof}

We note at this point that while the ideas behind our considerations in this subsection play an important role in the upcoming discussion of the bounded genus case of the Vital Linkage Theorem, we need somewhat stronger minimality assumptions on the type of ``wells'' that will appear there.
Nonetheless our results in this subsection and the remainder of this section are cut from similar cloth and we will require wells ~--~ both drained and dry ~--~ at several points throughout the paper, as they form a basic type of configuration we encounter.

\subsection{Linkages in graphs on a surface}
\label{sec:mazoitproof}

In 2013 Mazoit proved a single exponential upper bound for the irrelevant vertex theorem for graphs embedded into a fixed surface. More precisely, he proved the following \cite[Theorem 1]{Mazoit2013Single}. Let $G$ be a graph embedded into a surface $\Sigma$ of genus $g \geq 0$. Suppose that $G$ contains a sequence $\langle C_1, \ldots, C_t\rangle$ of pairwise disjoint cycles such that each $C_i$ bounds a disc $\Delta(C_i) \subseteq \Sigma$ with $\Delta(C_i) \subseteq \mathsf{int}\big(\Delta(C_{i+1})\big),$ for all $1 \leq i < t$. 
Let $\LLL \subseteq G$ be a linkage of order $k$ in $G$ such that no path $L \in \LLL$ has an endpoint in $\Delta(C_t)$. In this setting Mazoit proved that if $t \geq f(g,k) \coloneqq \frac{20}{9}\cdot k \cdot c^{3(g-1) + 2k},$ where $c = 3e^{\frac{10}{3e}}$ is a small constant, then there is a linkage $\LLL' \subseteq G$ with $\tau(\LLL) = \tau(\LLL')$ and $\LLL'$ is disjoint from $G \cap \mathsf{int}\big(\Delta(C_1)\big)$.

To prove the theorem, Mazoit first performs some preprocessing steps
to reduce the input instance to a simpler form that he calls a
\emph{reduced instance} which corresponds to what we call a wedding cake below. See \zcref{def:wedding-cake} for details. The technical core of his argument is to
show that in such a reduced instance the vertices deep in the middle
of a separating set of concentric cycles is irrelevant.

Mazoit made his proof available as a preprint on arXiv \cite{Mazoit2013Single}, but to date it has not appeared in a peer-reviewed venue. His proof in \cite{Mazoit2013Single} is at some places very condensed and not always completely polished and therefore occasionally a bit hard to follow. 

Therefore, in this section, we reprove the part of his result that we need in full detail and also adapt it to the notation in the rest of this paper.

Before we state the theorem we prove we first need to introduce some
notation.
\paragraph{Wedding cakes}

We first define the kind of surfaces we work with. 

\begin{definition}
 For a number $g \geq 0,$ a \emph{disc with $g$ strips} is a surface $\Sigma$
 obtained from a closed disc $\Delta$ as follows.
 Let $d^1_1, d_1^2, \dots, d^1_{g},
 d^2_g$ be pairwise disjoint closed and connected subsets of
 $\bd(\Delta)$. Let $S_1, \dots, S_g$ be closed discs and, for all $1 \leq i \leq g,$ let $s_i^1$ and
 $s_i^2$ be disjoint connected closed subsets of $\bd(S_i)$. Then
 $\Sigma$ 
 is obtained from $\Delta$ by gluing, for all $1 \leq i \leq g,$ $s_i^1$ onto $d_i^1$ and $s_i^2$ onto
 $d_i^2$ in any of the two possible orientations.

 We say that $\Sigma$ is \emph{generated} by $(\Delta, S_1, \dots, S_g)$ and
 denote this by $\Sigma \hat{=} (\Delta, S_1, \dots, S_g)$.

 We refer to the discs $S_1, \dots, S_g$ as the \emph{strips} of $\Sigma$
 and to $\Delta$ as its \emph{base}. The curves $s_i^1$ and $s_i^2$ are the \emph{sides} of $S_i$.

 The \emph{segments} of $\Sigma$ are the components of $\bd(\Delta) \cap \bd(\Sigma)$.
 A \emph{pre-segment} of $\Sigma$ is either a segment of $\Sigma$ or a side of
 a strip. 
\end{definition}

In this section we will only work with graphs $G$ embedded into a
surface $\Sigma$ generated by $(\Delta, S_1, \dots, S_g)$. It is therefore convenient for us to assume that $G \subseteq
\Sigma,$ i.\@e.\@~the vertices of $G$ are points on $\Sigma$ and the edges are
simple curves connecting their endpoints. 
Our goal is to show that any vertex in $G$ that is separated from a set of terminals by a sufficient number of concentric cycles will be irrelevant for the existence of a linkage connecting the terminals. In the next definition we formally define what we mean by a set of concentric cycles. 

\begin{definition}
 Let $\CCC \coloneqq \langle C_1, \dots, C_t\rangle$ be a sequence of pairwise disjoint cycles
in $G$. We say that $\CCC$ is a \emph{canonical sequence of concentric cycles}
if each $C_i,$ for $1 \leq i < t,$ bounds a closed disc $\Delta(C_i) \subseteq \Delta$ such that $\Delta(C_i) \subseteq \mathsf{int}\big(\Delta(C_{i+1})\big)$. The number $t$ is the \emph{order} of $\CCC$.
\end{definition}

As mentioned in the introduction of this section, we will only prove the irrelevant vertex theorem for instances of a very special form. In the next definition we introduce notation for linkages that we need to formally define the special form of instances we work with later on.
Recall that the segments of $\Sigma\hat{=}(\delta, S_1, \ldots, S_g)$ are the components of $\bd(\Delta) \cap \bd(\Sigma)$. 

\begin{definition}\label{def:valleys}
Let $\Sigma$ be a disc with $g$ strips generated by $(\Delta, S_1, \dots, S_g)$. 
Let $G \subseteq \Sigma$ be a graph and let $\CCC$ be a canonical sequence of
concentric cycles of order $t$ embedded into $\Delta$ such that $\bd(\Delta) =
C_t$.

Let $\LLL \subseteq G$ be a linkage whose endpoints are all contained in the
 segments of $\Sigma,$ i.\@e.\@\,in $\bd(\Sigma) \cap \bd(\Delta)$.

 A component $L$ of $\LLL \cap \Delta$ is a \emph{valley} if
 there is $d \leq t$ such that the vertices of $L,$ as they appear, are $\langle v_0, \dots, v_{2d}\rangle$
 with
 $v_i, v_{2d-i} \in V(C_{t-i}),$ for all $0 \leq i \leq d$. We call
 $d$ the \emph{depth} of the valley $L$.

 The linkage $\LLL$ is \emph{clean} (w.r.t.~$\Sigma$ and $\CCC$) if each component $L$ of $\LLL
 \cap \Delta$ is a valley and, if $g\geq 1,$ then the two endpoints of $L$ are contained in
 different pre-segments of $\Sigma$. 
\end{definition}

A property of valleys that is implicit in their definition is that
every edge of a valley has its endpoints on different consecutive cycles of $\CCC$.
We will need this property below.

\begin{figure}
 \centering
 \includegraphics{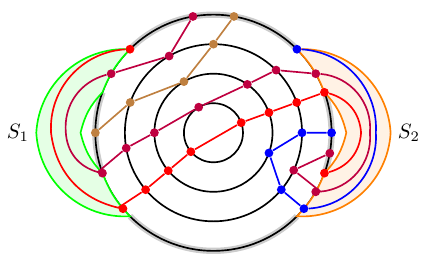}
 \caption{A wedding cake of genus $2,$ order $4,$ and depth $4$. The $4$ segments are marked by the grey areas on the outer cycle. }
 \label{fig:wedding-cake}
\end{figure}

We are now ready to define the special form of instances for the irrelevant
vertex problem we will work with. See \zcref{fig:wedding-cake} for an illustration.

\begin{definition}\label{def:wedding-cake}
 Let $\Sigma$ be a disc with $g$ strips generated by $(\Delta, S_1, \dots,
 S_g)$. Let $G \subseteq \Sigma$ be a graph, $\CCC = \langle C_t, \ldots, C_1\rangle$ be a canonical set of concentric
 cycles of order $t,$ and let $\LLL \subseteq G$ be a linkage of order $k$.
 The tuple $(\Sigma, G, \CCC, \LLL)$ is a \emph{wedding cake} of
 \emph{depth} $t,$ \emph{genus} $g,$ and \emph{order} $k$ if it meets the following
 conditions.
 \begin{enumerate}[(i)]
 \item $\LLL$ is clean (w.r.t.\,$\Sigma$ and $\CCC$), the endpoints of the
 paths in $\LLL$ are contained in the segments of $\Sigma,$ and otherwise $\LLL$ is disjoint from $\bd(\Sigma)$.
 \item $\Delta(C_t) = \Delta$.
 \item $V(G) = V(\LLL) = V(\CCC),$ $E(\CCC) \cap E(\LLL) = \emptyset,$ and $E(G) = E(\CCC) \cup
 E(\LLL)$.
 \end{enumerate}
\end{definition}

Given a wedding cake $(\Sigma, G, \CCC, \LLL),$ we say that a vertex $v \in V(G)$
is \emph{irrelevant} if $G{-}v$ contains a linkage $\LLL'$ with $\tau(\LLL') =
\tau(\LLL)$.

Finally, we fix an orientation of $\bd(\Delta)$ which also orients all
cycles in $\CCC$ in the same way. This allows us to say that a vertex $u$ in a pre-segment $s$ of $\Sigma$ comes \emph{before} or \emph{after} a vertex $v$ on $s$.

Before stating the main theorem of this section we collect a few
useful properties of wedding cakes we use below.

\begin{lemma}\label{lem:pyramid}
 Let $\WWW \coloneqq (\Sigma \hat{=}
 (\Delta, S_1, \dots, S_g), G, \CCC, \LLL)$ be a wedding cake of depth $t,$
 genus $g,$ and order $k$.

 Let $p \subseteq \bd(\Delta)$ be a pre-segment and let $\PPP = \{ P_0, \dots, P_\ell\}$
 be the paths in $\LLL \cap 
 \Delta$ with an endpoint in $p$ ordered so that the endpoints $p_0,
 \dots, p_\ell$ of $P_0, \dots, P_\ell$ appear in this order on $\bd(\Delta)$
 w.r.t.\,to the orientation of $\bd(\Delta)$. 

 Then for each $0 \leq j \leq \lceil \frac{\ell}2 \rceil$ the paths $P_j$ and $P_{\ell-j}$
 both contain a subpath
 linking their endpoint on $p$ and a vertex on $C_{t-j}$.
 \end{lemma}
 \begin{proof}
 Suppose the claim was false and let $j$ be minimal such that $P_j$
 or $P_{\ell-j}$ does not contain a subpath linking their endpoint
 on $p$ with a vertex on $C_{t-j}$.
 W.l.o.g.~we assume that $P_j$ does not contain such a subpath. 

 It is easily seen that $j>0$.
 By our choice of $j,$ the paths $P_{j-1}$ and $P_{\ell + 1 - j}$ both
 contain subpaths $I_{j-1} \subseteq P_{j-1}$ and $I_{\ell + 1 - j} \subseteq P_{\ell + 2
 - j}$ linking their endpoints $p_{j-1}$ and $p_{\ell + 1
 -j}$ on $p,$ resp., to some vertex $q_{j-1}$ and $q_{\ell + 1 - j},$ resp., on
 $C_{t-(j-1)}$. But if $P_j$ is disjoint from $C_{t-j},$ then $p \cup I_{j-1} \cup I_{\ell + 1 - j}$ together
 with the subpath $C' \subseteq C_{t-j}$ between $q_{j-1}$ and $q_{\ell + 1 -
 j}$ bounds a closed disc $D$ that contains $P_j$ with $\bd(D) \cap
 \bd(\Sigma) = p$. But this is a contradiction to the requirement of valleys that
 the other endpoint of $P_j$ is not on $p$.
 \end{proof}

 Given a segment $p$ in a wedding cake we will refer to the set of
 paths $I_1, \dots, I_\ell$ constructed in the previous lemma as a
 \emph{pyramid} over $p$ below. 

We are now ready to state the main technical theorem in Mazoit's
paper. 
\begin{theorem}[see {\cite[Theorem 2]{Mazoit2013Single}}]\label{thm:MazoitEndgame}
 Let $f_{\ref{thm:MazoitEndgame}}(g, k) \coloneqq \frac{20}9k\cdot(3e^{\frac{10}{3e}})^g$ and let $(\Sigma \hat{=}
 (\Delta, S_1, \dots, S_g), G, \CCC, \LLL)$ be a wedding cake of depth $t \geq
 f_{\ref{thm:MazoitEndgame}}(g, k),$ genus $g,$ and order $k$. Then there is a
 linkage $\LLL'$ with $\tau(\LLL') = \tau(\LLL)$ which is disjoint from $\mathsf{int}\big(\Delta(C_1)\big)$.
 That is, every
 vertex of $G$ inside $\mathsf{int}\big(\Delta(C_1)\big)$ is irrelevant for the existence of a
 linkage in $G$ realising $\tau(\LLL)$.
\end{theorem}

In the remainder of this section we prove the previous theorem. So far we have worked with linkages $\LLL$ as sets of paths in $G$. For some of our arguments we will work with \emph{topological linkages} which are simply sets of disjoint curves on $\Sigma$. We introduce the required notation next.

\begin{definition}
 Let $G$ be a graph and let $A \subseteq V(G)$. 
 A \emph{pattern on $A$} is a set $\tau \subseteq \big\{ \{ a,b \} \sth a,b \in A \}$ of pairs of vertices of $A$ (not necessarily disjoint). If $A = V(G)$ we simply say that $\tau$ is a pattern. 
 
 Let $\Sigma$ be a surface. 
 A \emph{topological linkage} $\TTT$ in $\Sigma$ with pattern $\tau$ is a set of
 pairwise 
 disjoint simple curves such that $\TTT$
 contains for each pair $\{a,b\} \in \tau$ a curve with endpoints $a$ and $b$. We
 denote the \emph{pattern} of $\TTT$ by $\tau(\TTT)$.

 A pattern
 $\tau$ is \emph{realisable} (in $G$) if there is a linkage $\LLL$ in
 $G$ with $\tau(\LLL) = \tau$ and it is \emph{topologically feasible}
 (in $\Sigma$) if there is a topological linkage $\TTT$ in $\Sigma$ with
 $\tau(\TTT) = \tau$.
\end{definition}

Obviously if $\tau$ is feasible then it must also be topologically
feasible. The converse is not true in general. We call a
topological linkage $\TTT$ \emph{realisable} in $G$ if there is a
linkage $\LLL \subseteq G$ such that $\tau(\LLL) = \tau(\TTT)$.

Throughout this part we will use $\LLL$ (possibly with subscripts) to denote linkages in $G$ and $\TTT$ (again possibly with subscripts) to denote topological linkages in $\Sigma$.

We first need the following easy lemma which proves a special case to
which we reduce in the sequel.

\begin{lemma}\label{lem:outerplanar-cylinder}
 Let $G$ be a graph embedded on a disc $\Delta$ and let $\CCC \coloneqq \langle C_1,
 \dots, C_t\rangle$ be a canonical sequence of concentric cycles such that $C_t = \bd(\Delta)$. Furthermore, let $\PPP$ be a $V(C_t){-}V(C_1)$-linkage of order $2k$ in $G$ such that $|V(P \cap C_j)| = 1$ for each $P \in \PPP$ and $1 \leq j \leq t$. Let $\tau$ be a pattern of order $k$ on $V(C_t) \cap V(\PPP)$. 
 
 If $t \geq k$ and $\tau$ is topologically feasible then it is realisable in $G$. Furthermore, $\tau$ can be realised in $G \setminus \mathsf{int}(\Delta(C_1))\subseteq \Sigma$.
\end{lemma}
\begin{proof} We prove the lemma by induction on $k$.
 Let $\PPP \coloneqq \{ P_1, \dots, P_{2k} \}$ and, for $1 \leq j \leq 2k,$ let $p_j$ be the endpoint of $P_j$ on $C_t$. As $\tau$ is topologically feasible, there must be a pair $\{a,b\} \in \tau$ such that $C_t$ contains a subpath $C(a,b)$ with ends $a$ and $b$ that does not contain any other vertex occurring in an element of $\tau$. 
 Let $A \subseteq \tau$ be the set of such pairs. Thus, $|A| \geq 1$.

 We construct a pattern $\tau'$ as follows. For $\{a,b\} \in \tau \setminus A$ let $P_a, P_b$ be the paths in $\PPP$ with endpoints $a,b,$ resp., and let $a^\star$ be the vertex in $P_a \cap C_{t-1}$ and let $b^\star$ be the vertex in $P_b \cap C_{t-1}$. Let $P_a^\star$ and $P_b^\star$ be the subpaths of $P_a$ and $P_b$ between $a$ and $a^\star$ and $b$ and $b^\star,$ resp. Then $\tau' \coloneqq \{ \{ a^\star, b^\star\} \sth \{a,b\} \in \tau \setminus A \}$.

 By construction, $\tau'$ is a pattern of order $k'$ for some $k' < k$ and all vertices occurring in $\tau'$ are contained in $V(C_{t-1})$.
 Let $\CCC' \coloneqq \langle C_1, \ldots, C_{t-1}\rangle$ and let $\PPP' \coloneqq \PPP\cap\Delta(C_{t-1})$. Thus, $G' \coloneqq G \cap \Delta(C_{t-1}),$ $\PPP',$ $\CCC',$ and $\tau'$ satisfy the assumptions of the induction hypothesis and hence there is a linkage $\LLL' \subseteq G' \setminus \mathsf{int}(\Delta(C_1))\subseteq \Sigma$ 
 with $\tau(\LLL') = \tau'$. 
 For each $L' \in \LLL'$ with endpoints $a^\star$ and $b^\star$ let $L = P_a^\star \cup L \cup P_b^\star$ be the extension of the path $L' \in \LLL$ to the cycle $C_t$ and let
 $\LLL \coloneqq \{ L \sth L' \in \LLL' \} \cup \{ C(a,b) \sth \{ a,b\}\in A\}$. 
 Then $\LLL \subseteq G\setminus \mathsf{int}(\Delta(C_1)) \subseteq \Sigma$
 and $\tau(\LLL) = \tau$ as required.
\end{proof}

We also need the following lemma in later sections which can be proved in a similar way. 

\begin{lemma}\label{lem:two-sided-cylinder}
 Let $G$ be a graph embedded into a disc $\Delta$ and let $\mathcal{C} \coloneqq \langle C_1,
 \dots, C_t\rangle$ with $t \geq 2k$ be a canonical sequence of concentric cycles such that $C_t = \bd(\Delta)$ and $\mathsf{int}\big(\Delta(C_1)\big) = \emptyset$.
 Let $\mathcal{P}$ be a $V(C_t){-}V(C_1)$-linkage of order $2k$ in $G$ such that $|V(C_i) \cap V(P)| = 1,$ for all $1 \leq i \leq t$ and all $P\in \PPP$. Finally, let $\tau$ be a pattern of order $k$ on $(V(C_t) \cup V(C_1)) \cap V(\mathcal{P})$.

 If $\tau$ is topologically feasible then it is realisable in $G$.
\end{lemma}
\begin{proof}
 A pair $\{a,b\} \in \tau$ is \emph{local} if $\{a,b\} \subseteq V(C_1)$ or $\{a,b\} \subseteq V(C_t)$.
 We prove the lemma by induction on the number of local pairs in $\tau$. 
 
 We first consider the case where $\tau$ does not have a local pair. Fix an orientation of $\Delta$ which implies a consistent orientation of $C_t$ and $C_1$.
 Let $\tau \coloneqq \big\{ \{a_i, b_i \} \sth 1 \leq i \leq k \big\}$ such that all $a_i \in V(C_t)$ and all $b_i \in V(C_1),$ numbered such that $a_1, \ldots, a_k$ occur in this order on $C_t$. As $G$ is planar and $\tau$ is topologically feasible, $b_1, \dots, b_k$ must also appear in this order in $C_1$. 

 Let $P(a_1) \in \PPP$ be the path containing $a_1$ and let $c_1$ be the other endpoint of $P(a_1)$ on $C_1$. Let $S_1$ be the subpath of $C_1$ starting at $c_1$ and ending at $b_1$ in the direction of $C_1$ and let $S_2$ be the other subpath of $C_1$ starting at $b_1$ and ending at $c_1$. Let $P(b_1) \in \PPP$ be the path with endpoint $b_1$.
 
 At least one of $S_1, S_2$ contains $ \geq k_1$ endpoints of paths in $\PPP \setminus \{ P(a_1), P(b_1) \}$. W.l.o.g.~we assume that this is $S_2$. 

 Let $a'$ be the unique vertex in $V(P(a_1)) \cap C_{\lfloor \frac t2 \rfloor}$ and let $b'$ be the unique vertex of $V(P(b_1)) \cap C_{\lfloor \frac t2 \rfloor}$. Let $S$ be the subpath of $C_{\lfloor \frac t2 \rfloor}$ between $a'$ and $b'$ in the same direction as $S_1$. Let $L$ be the union of $S,$ the subpath of $P(a_1)$ between $a_1$ and $a',$ and the subpath of $P(b_1)$ between $b'$ and $b_1$. By construction, $L$ links $a_1$ and $b_1$. 

 Let $\LLL'$ be a linkage of order $k-1$ between $A \coloneqq \{ a_2, \ldots, a_k \}$ and $B \coloneqq \{b_2, \ldots, b_k \}$ in $G - L$. Such a linkage exists, for otherwise, by Menger's theorem, there would be a separation $(A', B')$ of $G - L$ of order $< k-1$ with $A \subseteq A'$ and $B \subseteq B'$. But this is impossible because there are $k-1$ disjoint paths in $\PPP \setminus \{ P(a_1), P(b_1) \}$ with an endpoint in $B$ which intersect each of the cycles $C_j$ with $1 \leq j \leq {\lfloor \frac t2 \rfloor - 1},$ there are $k-1$ paths in $\PPP \setminus \{ P(a_1), P(b_1) \}$ with an endpoint in $A$ which intersect each of the cycles $C_j$ with $t \geq j \geq {\lfloor \frac t2 \rfloor+1}$ and, by our choice of $S_1,$ there are $k-1$ paths in $\PPP \setminus \{ P(a_1), P(b_1) \}$ linking $S_2 \subseteq C_1$ to $C_t$. As $|V(A') \cap V(B')| \leq k-2,$ the union of these paths and cycles contains an $A{-}B$-path disjoint from $V(A')\cap V(B'),$ a contradiction. 

 Thus the $A{-}B$-linkage $\LLL'$ exists. As $G-L$ is a planar graph with $A$ and $B$ contained in disjoint subpaths of the outer cycle, $\LLL'$ must link $a_i$ to $b_i,$ for all $2 \leq i \leq k$. Thus $\LLL \coloneqq \{ L \} \cup \LLL'$ is a linkage realising $\tau$. 
 This concludes the induction base. 

 Now suppose that $\tau$ does contain a local pair $\{a,b\}$. We can now argue exactly as in the proof of \zcref{lem:outerplanar-cylinder}. In particular, as $\tau$ is topologically feasible, $\tau$ must contain a local pair $\{a,b\}$ such that one of the two subpaths of the cycle $C \in \{C_t, C_1\}$ containing $a, b$ between $a$ and $b$ contains no other terminal from $\tau$. Let $P \subseteq \tau$ be the set of all such pairs. Then we can link these pairs directly by these subpaths of the cycle containing them without intersecting any other terminal pair. Let $\tau' \coloneqq \tau \setminus P$. For each terminal $a \in V(C_t)$ occurring in $\tau'$ let $P(a) \in \PPP$ be the path in $\PPP$ containing $a$ and let $a'$ be the unique vertex in $V(P(a) \cap C_{t-1})$. Analogously, for each $b \in V(C_1)$ occurring in $\tau'$ let $P(b) \in \PPP$ be the path containing $b$ and let $b'$ be the unique vertex of $P(b) \cap V(C_2)$. Let $\PPP'$ be the set of subpaths of paths in $\PPP$ between $C_{t-1}$ and $C_2$. Let $G' \coloneqq \PPP' \cup \bigcup \{ C_2, \dots, C_{t-1} \}$. By induction hypothesis, there is a linkage $\LLL'$ in $G'$ realising $\tau'$. This can easily be extended to a linkage $\LLL$ in $G$ realising $\tau$. 
\end{proof}

\subsection{Proof of \zcref{thm:MazoitEndgame}}

We are now ready to start proving \zcref{thm:MazoitEndgame}. The main
idea is as follows. Let $(\Sigma \hat{=} (\Delta, S_1, \dots, S_g), G, \CCC, \LLL)$ be
a wedding cake. We now consider the subpaths of $\LLL$ within a single
strip $S_i,$ i.\@e.\@~the linkage $\LLL \cap S_i$.
If $\LLL$ crosses $S_i$ only a bounded number of times, we can
remove $S_i$ from $\Sigma$ and treat the vertices in $\LLL \cap \bd(S_i)$ as new
terminals. If this can be done for every strip then we arrive at a
planar instance embedded into $\Delta$ with all terminals on $\bd(\Delta)$. This
is the simplest case that we prove in \zcref{lem:outerplanar-routing}.

Otherwise we can take a strip $S_i$ which $\LLL$ crosses
sufficiently often. We then apply the following theorem by Geelen, Huynh, and Richter \cite{GeelenHR2018Explicit}, which is the topological engine of Mazoit's argument.

\begin{theorem}\cite[Theorem 2.6]{GeelenHR2018Explicit}\label{thm:jim}
 Let $\Sigma$ be a surface and let $P$ be a non-separating $\bd(\Sigma)$-path
 in $\Sigma$. For any linkage $\LLL$ in $\Sigma$ whose ends are disjoint from
 $P,$ there is a $\bd$-homeomorphism $\phi \sth \Sigma \rightarrow \Sigma$ such that each
 path of $\phi(\LLL)$ intersects $P$ at most twice. 
\end{theorem}

This will allow us to show that the pattern $\tau(\LLL)$ can be realised by
a topological linkage $\TTT_i$ which crosses $S_i$ only a bounded number of
times. Repeating this for each strip we obtain a topological linkage
$\TTT$ which realises the pattern $\tau(\LLL)$ but crosses each strip only a
bounded number of times. 
The difficult step in Mazoit's proof is to show that $\TTT$ can
be realised by a linkage in $G$ which also crosses each strip only a
bounded number of times. Once this is done we are in the case that we
already discussed above.

We now give the details of this proof strategy.
We start with the base case which essentially settles the case where $\LLL$ crosses
each strip only a bounded number of times.

\begin{lemma}\label{lem:outerplanar-routing}
 Let $(\Sigma \coloneqq \Delta, G, \CCC, \LLL)$ be a wedding cake of depth $t \geq
 2k,$ genus $0,$ and order $k$. Then there is a
 linkage $\LLL'$ with $\tau(\LLL') = \tau(\LLL)$ which is disjoint from $\mathsf{int}\big(\Delta(C_1)\big)$.
 That is, every
 vertex of $G$ in $\mathsf{int}\big(\Delta(C_1)\big)$ is irrelevant for the existence of a
 linkage in $G$ to realise $\tau(\LLL)$.
\end{lemma}
\begin{proof}
 The proof is almost identical to the proof of \zcref{lem:outerplanar-cylinder}. As $\Sigma$ has no strips, the instance is planar with all terminals on the outer cycle $C_t$. Thus, there must be a pair $\{a,b\} \in \tau(\LLL)$ such that $C_t$ contains a subpath $P(a,b)$ with ends $a$ and $b$ which is internally disjoint from $\LLL$ (in fact, by the definition of wedding cakes, $a=b,$ but this does not matter here). As in the proof of \zcref{lem:outerplanar-cylinder} we add $P(a,b)$ to $\LLL'$ for all such pairs and apply induction on the remaining pairs and the sequence $\langle C_1, \ldots, C_{t-1}\rangle$ of cycles.
\end{proof}

We now prove the main lemma that establishes the inductive step. 
For $0 \leq i \leq g$ we define $\Sigma_i \coloneqq \Sigma \setminus \big( (S_1 \cup \dots S_i) \setminus \bd(\Delta)\big)$. That is,
$\Sigma_0 = \Sigma$ and, for $1 \leq i \leq g,$ $\Sigma_i$ is obtained from the disc $\Delta$ by
attaching only the strips $S_{i+1}, \dots, S_g$. For $0 \leq i \leq g,$ a

Given a linkage $\LLL,$ we define the \emph{size} of a strip $S_i,$
denoted by $|S_i|,$ as the number of components in $S_i \cap \LLL$.

Recall that, by definition of a wedding cake, $\LLL$ is clean. 

\begin{lemma}\label{lem:Mazoit-lem3}
 Let $(\Sigma \hat{=} (\Delta, S_1, \dots, S_g), G, \CCC, \LLL)$ be a wedding cake of 
 genus $g \geq 1,$ order $k,$ and depth $t \geq 2k3^g$. If $|S_i| \geq 3k3^i+1,$
 for all $1 \leq i \leq g,$ then every vertex in $\mathsf{int}\big(\Delta(C_1)\big)$ is irrelevant. 
\end{lemma}
\begin{proof}
 Let $1 \leq i \leq g$ and let $s_i \coloneqq |S_i|-1$.
 Let $p$ be one of the two sides of $S_i$
 and let $\JJJ \coloneqq \LLL \cap S_i$. Let $\JJJ = \{ J_0, \dots, J_{s_i}\}$ be ordered in the order in
 which the endpoints $p_0, \dots, p_{s_i}$ of the paths in $\JJJ$ appear on
 $p$.
 By \zcref{lem:pyramid}, the paths in $\Delta \cap \LLL$ with an endpoint in $\{
 p_0, \dots, p_{s_i}\}$ contain subpaths $I_0, \dots, I_{s_i}$ such that $I_j$ links $p_j$ and $C_{t- \min \{ j, {s_i} - j \}}$.

 Now let $p'$ be the other side of $S_i$. For $0 \leq j \leq s_i$ let $p_j'$ be the endpoint of $J_j$ on $p'$. Note that either $p_0', \ldots, p_{s_i}'$ appear in this order on $p'$ or in reverse.
 Again, by
 \zcref{lem:pyramid}, the paths in $\Delta \cap \LLL$ with an endpoint in $p'$
 contain subpaths $I_0', \dots, I_{s_i}'$ such that $I'_j$ links $p'_j$ and
 $C_{t- \min \{ j, {s_i} - j \}}$. Thus, $\SSS_i' \coloneqq \{ J_j \cup I_j \cup I'_j \sth 0 \leq j \leq s_i \}$ is a linkage which contains for each $0 \leq j \leq {s_i}$ a path $P_j$
 containing $p_j$ which has both endpoints on $C_{t- \min \{ j, {s_i} - j \}}$. 
 Let $l_1 \coloneqq k\cdot 3^i$ and $l_2 \coloneqq l_1 + \lfloor \frac12\cdot k \cdot 3^{i-1}\rfloor$. Let $r_2 \coloneqq 2k\cdot 3^i$ and let $r_1 \coloneqq r_2 - \lfloor \frac 12 k \cdot 3^{i-1} \rfloor$. 

 We define $\SSS_i \coloneqq \{ P_j \sth l_1 \leq j \leq r_2 \}$ and
 set 
 \begin{itemize}
 \item
 $\SSS_i^l \coloneqq \{ P_j \sth l_1 \leq j < l_2\},$
 \item $\SSS_i^m \coloneqq \{ P_j \sth l_2 \leq j < r_1 \},$ and
 \item 
 $\SSS_i^r \coloneqq \{ P_j \sth r_1 \leq j < r_2 \}$.
 \end{itemize}

 We also define $\SSS_i^{ll} \coloneqq \{ P_j \sth 0 \leq j < l_1 \}$ and $\SSS_i^{rr} := \{ P_j \sth r_2 < j \leq s_i \}$. See \zcref{fig:pyramids} for an illustration. 
 By construction, $|\SSS_i^l|, |\SSS_i^r| = \frac12 k 3^{i-1}$ and $|\SSS_i^m| = 2k3^{i-1}$.
 Note that $|\SSS_i^{ll}| = |\SSS_i^{rr}| = k\cdot 3^k$ and thus $|\SSS_i^{ll}| = |\SSS_i^l| +|\SSS_i^m| +|\SSS_i^r|$. This will be important below.
 
\begin{figure}
 \centering
 \includegraphics{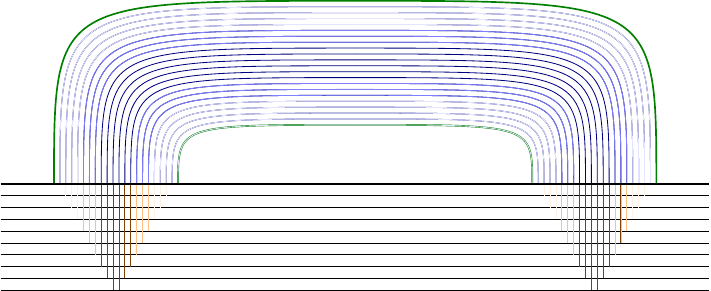}
 \caption{Illustration of the constructions of the sets $\SSS_i^l, \SSS_i^m, \SSS_i^r$. The paths marked in orange are the paths $I_j$ (left) and $I'_j$ (right). The paths in blue are the paths $J_j$. The paths in the middle in the darkest colours are $\SSS_i^m,$ the outermost paths in the very light colours are the paths in $\SSS_i^{ll}$ and $\SSS_i^{rr},$ resp. }
 \label{fig:pyramids}
\end{figure}
 \begin{beautifulclaim}
 There is a topological linkage $\TTT \subseteq \Sigma$ with $\tau(\TTT) = \tau(\LLL)$ such
 that, for all $1 \leq i \leq g,$ $\TTT$ crosses $S_i$ at most $2\cdot k \cdot 3^{i-1}$
 times. Furthermore, $\TTT \cap S_i \subseteq \SSS_i^m,$ for all $1 \leq i \leq g$. 
 \end{beautifulclaim}
 \begin{claimproof}
 We prove by induction on $0 \leq i \leq g$ that there is a topological
 linkage $\TTT_i$ with $\tau(\TTT_i) = \tau(\LLL)$ such that, for all $1 \leq j \leq i,$
 $\TTT_i$ crosses $S_j$ at most $2 \cdot k \cdot 3^{j-1}$ times and $\TTT_j \cap S_j
 \subseteq \SSS^m_j$. Furthermore, $|\TTT_i \cap \Sigma_i| \leq k \cdot 3^i$.

 For $i=0$ we can choose $\TTT_0 = \LLL,$ which obviously satisfies the
 conditions as $\Sigma_0 = \Sigma$.

 Now suppose $\TTT_{i-1}$ has already been defined, for some $1 \leq i \leq 
 g$. By the induction hypothesis, $|\TTT_{i-1} \cap \Sigma_{i-1}| \leq k \cdot 3^{i-1}$.

 Let $R \subseteq S_i \setminus \bd(\Delta)$ be a
 boundary curve in $S_i$ separating its two sides. By 
 \zcref{thm:jim} there is a $\bd$-homeomorphism $\phi \sth \Sigma_{i-1} \rightarrow \Sigma_{i-1}$ that
 maps $\TTT_{i-1}$ to a topological linkage $\TTT_i' = \phi(\TTT_{i-1})$ 
 with the same pattern as $\TTT_{i-1}$ in $\Sigma_{i-1}$ such that each
 path in $\TTT_i'$ crosses $S_i$ at most twice. Thus, in total $\TTT_i'$
 crosses $S_i$ at most $2 \cdot k \cdot 3^{i-1}$ times. As $|\SSS_i^m| = 2 \cdot k
 3^{i-1},$ we may shift $\TTT_i'$ to obtain a linkage $\TTT_i$ with
 $\tau(\TTT_i) = \tau(\TTT'_i)$ such that $\TTT_i \cap S_i = \SSS_i^m$. 

 It follows that $|\TTT_i \cap \Sigma_i| = |\TTT_{i-1} \cap \Sigma_{i-1}| + |\TTT_i \cap S_i| \leq
 k \cdot 3^{i-1} + 2 \cdot k \cdot 3^{i-1} = k \cdot 3^i$ as required.

 We set $\TTT \coloneqq \TTT_g$ which satisfies the requirements of the claim.
 \end{claimproof}

 We show next that the pattern $\tau(\TTT)$ can also be realised by a
 linkage in $G$. For this, we first prove the following claim. 

 Let $s$ be a segment of $\Sigma_i,$ for some $1 \leq i \leq g,$ and let $p_0, \dots, p_r$ be the vertices of $\TTT \cap s$ in the order in which they appear on $s$. Similar to the discussion below \zcref{lem:pyramid} we define a \emph{pyramid} $\PPP(s)$ for $s$ as a set of paths $P_0, \dots, P_r$ such that each vertex of $\TTT \cap s$ is endpoint of exactly one path and for $0 \leq j \leq \lceil \frac r2 \rceil$ the paths $P_j$ and $P_{r-j}$ link a vertex on $C_t$ and a vertex on $C_{t-j}$.
Note that, by construction, if $s$ is a side of $S_i,$ for some $1 \leq i \leq g,$ then $|\TTT \cap s| \leq 2 \cdot 3^{i-1}$. 
 \begin{beautifulclaim}
 Let $p_1, p_2$ be the two segments of $\Sigma_{g-1},$ i.\@e.\@, the two
 components of $\bd(\Delta) - \bd(S_g)$. Then $G \cap \Delta$ contains disjoint
 pyramids $\PPP(p_1)$ and $\PPP(p_2)$ of $p_1$ and $p_2,$ resp., which are disjoint from $\SSS_g$.
 \end{beautifulclaim}
 \begin{claimproof}
 By induction on $0 \leq i < g$ we show that $G$ contains a pyramid $\PPP^i(s)$ for each segment $s$ of $\Sigma_i$. 
 For $i=0$ this follows immediately from \zcref{lem:pyramid} and the fact that $\tau(\LLL) = \tau(\TTT)$. It is easily seen that for different segments $s \not= s'$ these pyramids are disjoint (as $t > k$).

 Now suppose that we have already constructed a pyramid $\PPP^{i-1}(s)$ for each segment of $\Sigma_{i-1}$. Let $s$ be a segment of $\Sigma_i$. 
 We distinguish between the following three cases.

 \noindent\textit{Case 1. } If $s$ is already a segment of $\Sigma_{i-1}$ then we simply set $\PPP^i(s) \coloneqq \PPP^{i-1}(s)$.

 \begin{figure}
 \centering
 \includegraphics[width=0.5\linewidth]{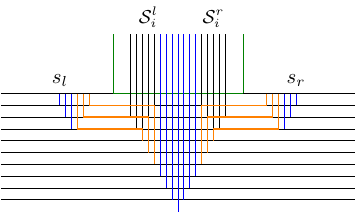}
 \caption{Illustration of the routing in Case 2.}
 \label{fig:mz-lem3-case2}
 \end{figure}
\noindent\textit{Case 2. } The next case we consider is that $s$
 is the union of one of the sides $s_m$ of $S_i$ and the two segments $s_l, s_r$ of $\Sigma_{i-1}$ to the left and to the right of $s_m$. See \zcref{fig:mz-lem3-case2} for an illustration. 
 
 Let $a_0, 
\dots, a_{p-1}, b_1, \dots, b_{q-1},$ and $c_0, \dots, c_{r-1}$ be the vertices of $\TTT \cap s_l,$ $\TTT \cap s_m,$ and $\TTT \cap s_r,$ resp., in the order in which they appear on $\bd(\Delta)$. Recall the definition of the sets $\SSS_i$ as well as $\SSS_i^l, \SSS_i^m,$ and $\SSS_i^r$ from the beginning of this proof. W.l.o.g.~we assume that on the side $s_m$ of $S_i$ the paths in $\SSS^i_l$ appear before the paths in $\SSS^i_m$ and that the paths in $\SSS^i_r$ come last.
By construction of $\TTT,$ $\TTT \cap s_m \subseteq \SSS_i^m \cap s_m$.

 Let $\QQQ_l$ be the components in $\SSS_i^l \cap \Delta$ with an end in $s_m$ and let $\QQQ_m$ and $\QQQ_r$ be defined analogously. Let $\QQQ_l \coloneqq \{ Q^l_0, \dots, Q^l_{|\QQQ_l|-1}\}$ be ordered in the order in which their endpoints appear on $s_m$ and define $\QQQ_r \coloneqq \{ Q^r_0, \dots, Q^r_{|\QQQ_r|-1}\}$ analogously.

 As observed above, $q \leq 2 \cdot k \cdot 3^{i-1} = |\QQQ_m|$ and $|\QQQ_l|, |\QQQ_r| \geq \frac 12\cdot k \cdot 3^{i-1}$. Furthermore, $|\TTT \cap \Sigma_{i-1}| \leq k \cdot 3^{i-1}$ and therefore $p,r 
 \leq k \cdot 3^{i-1}$.
 Let $L_0, \dots, L_{p-1}$ be the paths of $\PPP^{i-1}(s_l)$ such that $L_j$ contains $a_j,$ for $0 \leq j \leq p-1,$ and let $R_0, \dots, R_{r-1}$ be the paths of $\PPP^{i-1}(s_r),$ again numbered such that $R_j$ contains $c_j,$ for $0 \leq j \leq r-1$.
 We now construct a pyramid $\PPP^i(s)$ as follows. 
 \begin{itemize}
 \item For $0 \leq j \leq \lceil \frac {p-1}2 \rceil$ we set $P(a_j) \coloneqq L_j$ and for $0 \leq j \leq \lceil \frac {r-1}2 \rceil$ we set $P(c_{r-j}) \coloneqq R_{r-j}$. 
 \item For $0 \leq j \leq q-1$ we set $P(b_j)$ to be the subpath of the path $P \in \SSS_m$ with endpoint $b_j$ of the right length for it to be a proper part of the pyramid. This is always possible as the paths in $\SSS_m$ link $C_t$ to $C_{t-k\cdot 3^i}$.
 \item For $0 \leq j \leq \lfloor \frac {p-1}2\rfloor$ we construct the path $P(a_{p-j})$ as follows. Let $P'$ be the union of $L_{p-j},$ $Q^l_j,$ and the part of $C_{t-j}$ that connects $L_{p-j}$ and $Q^l_j$ and is disjoint from $\QQQ^i_m$. $P'$ is a path that connects $a_{p-j}$ to $C_{k\cdot 3^i}$ and thus we can take a subpath $P(a_{p-j}) \subseteq P'$ with $a_{p-j}$ as one end and the other end on $C_{t-(p-j)}$.

 The paths $P(c_{r-j}),$ for $0 \leq j \leq \lfloor 
 \frac {r-1}2 \rfloor$ are defined analogously using the paths in $\QQQ^r$ instead.

 The paths constructed in this last step are marked in orange in \zcref{fig:mz-lem3-case2}.
 \end{itemize}
 Now we set $\PPP^i(s) \coloneqq \{ P(u) \sth u \in \{ a_0, \dots, a_{p-1}, b_0, \dots, b_{q-1}, c_0, \dots, c_{r-1} \} \}$. By construction, $\PPP^i(s)$ is a pyramid over $s,$ as required.

 \noindent\textit{Case 3. } The last case we consider is the case where $s$ is the union of both sides $s_1$ and $s_2$ of $S_i$ and the segments $p_1, p_2, p_3,$ where $p_1$ is the segment before $s_1$ on $C_t,$ $p_2$ is the segment between $s_1$ and $s_2,$ and $p_3$ is the segment after $s_2$ on $C_t$. See \zcref{fig:mazoit-lem-3-case-3} for an illustration of this case. Recall the definition of $\SSS_i^l, \SSS_i^m,$ and $\SSS_i^r$ from above. In this case we also need the sets $\SSS_i^{ll}$ and $\SSS_i^{rr}$.

 For $1 \leq j \leq 3$ let $a^j_{0}, \ldots, a^j_{h_j}$ be the vertices of $p_j \cap \TTT$ and for $1 \leq j \leq 2$ let $b^j_{0}, \ldots, b^j_{f_j}$ be the vertices of $s_j \cap V(G)$. By construction, $f_1, f_2 \geq 3k3^i+1$. Furthermore, $h_1 + h_2 + h_3 \leq \TTT \cap \Sigma_{i-1}\leq k\cdot 3^{i-1}$. 

As in the previous case we let $\QQQ^1_l$ be the components of $\SSS_i^l \cap \Delta$ with an end in $s_1$ and we define $\QQQ^1_m, \QQQ^1_r, \QQQ^1_{ll},$ and $\QQQ^1_{rr}$ analogously. In the same way we define $\QQQ^2_l$ as the components of $\SSS_i^l \cap \Delta$ with an end in $s_2$ and we define $\QQQ^2_m, \QQQ^2_r, \QQQ^2_{ll},$ and $\QQQ^2_{rr}$ analogously. W.l.o.g.~we assume that the endpoints of $\QQQ^2_{ll}, \QQQ^2_l, \QQQ^2_m, \QQQ^2_r, \QQQ^2_{rr}$ appear in this order on $s_2$. 

 \begin{figure}
 \centering
 \includegraphics[width=.75\linewidth]{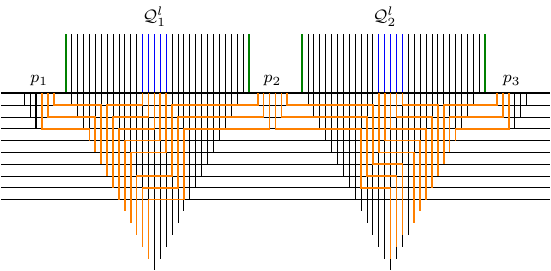}
 \caption{Illustration of Case 3 in \zcref{lem:Mazoit-lem3}. The two orange paths illustrate the rerouting.}
 \label{fig:mazoit-lem-3-case-3}
 \end{figure}

Our goal is to construct a pyramid $\PPP(s)$ as in Case 2 by rerouting some of the paths in $\PPP(p_1)$ along the cycles in $\CCC$ and then through the paths in $\PPP(s_1)$ and likewise for the paths in $\PPP(p_3)$. The main problem here is that we also need to reroute the paths in $\PPP(p_2)$. As these are the innermost paths of $s$ they also need to go deepest into the set $\CCC$ of concentric cycles but as $|\PPP(p_2)| \leq k \cdot 3^{i-1}$ the paths in $\PPP(p_2)$ are not guaranteed to go deep enough themselves. We therefore need to reroute them as well. 
 See \zcref{fig:mazoit-lem-3-case-3} for an illustration. 

The crucial observation is that as $h_1 \leq k \cdot 3^{i-1}$ but $|\QQQ^1_{ll}| = k\cdot 3^i,$ we can easily reroute the paths in $\PPP^{i-1}(p_1)$ with endpoints $a^1_{\lceil\frac 12 h_1\rceil}, \ldots, a^1_{h_1}$ along the cycles $C_{t-h_1}, \ldots, C_{t-\lceil\frac 12 h_1\rceil}$ (in this order) to the paths in $\QQQ^1_{ll}$. More precisely, let $\QQQ^1_{ll} \coloneqq \{ Q^{1,ll}_0, \ldots, Q^{1,ll}_{k3^i-1} \}$ be ordered by their endpoints on $s_1$ and let $\PPP^{i-1}(p_1) \coloneqq \{ L_0, \ldots, L_{|\PPP^{i-1}(p_1)|}\}$. Then for all $0 \leq j \leq \frac12 |\PPP^{i-1}(p_1)|$ we add the path $L_j$ to $\PPP^i(s)$. For $\frac12 |\PPP^{i-1}(p_1)| < j \leq \frac12 |\PPP^{i-1}(p_1)|$ we add the path obtained from the union of $L_{j},$ the subpath of $C = C_{t-|\PPP^{i-1}(p_1)|+j}$ between $L_j$ and $Q^{1,ll}_{j}$ and the subpath of $Q^{1,ll}_j$ between $C$ and $C_{t-j}$. 

The next paths we need to reroute are the paths in $\QQQ^1_m$. Recall that these are precisely $2k3^{i-1}$ paths and we need to route them to the cycles $C_{t-|\PPP^{i-1}(p_1)|+1}$ to $C_{t-|\PPP^{i-1}(p_1)|+2k3^{i-1}}$. 
So far we have only used $\leq k \cdot 3^{i-1}$ paths from $\QQQ^{1,ll}_i$. As $|\QQQ^{1,ll}_i| = k3^i$ there are still $2k3^{i-1}$ paths unused in $\QQQ^{1,ll}_i$. We will use these paths to route the paths of $\QQQ^1_m$ along the cycles $C_{t-k \cdot 3*{i-1}}, \ldots, C_{t-k\cdot 3^i}$ and then along the yet unused paths in $\QQQ^{1,ll}_i$ to the cycles $C_{t-|\PPP^{i-1}(p_1)|+1}$ to $C_{t-|\PPP^{i-1}(p_1)|+2k3^{i-1}}$. 

We route the paths in $\PPP^{i-1}(p_3)$ and $\QQQ^{2,m}_i$ symmetrically. 

The way we rerouted the paths so far ensures that we have not yet used any path in $\QQQ^{1,r}_i \cup \QQQ^{2,l}_i$. Thus we can use these paths to route the paths in $\PPP^{i-1}(p_2)$. This concludes the construction of the pyramid $\PPP^i(s)$ in Case 3. 

This concludes the induction step and thus the proof of the claim.
 \end{claimproof}

 We can now conclude the proof of the lemma. Let $p_1$ and $p_2$ be the two segments of $\Sigma_{g-1}$ and let $s_1$ and $s_2$ be the two sides of $S_g$ such that $p_1, s_1, p_2, s_2$ appear in this order on $\bd(\Delta)$. 
 We proceed as in the proof of the previous claim to find a set of disjoint paths linking $(p_1 \cup p_2) \cap \TTT$ to $C_{t-k\cdot3^g -1}$. By \zcref{lem:outerplanar-cylinder} there is a linkage $\LLL' \subseteq G$ with $\tau(\LLL') = \tau(\TTT)$ which avoids $\mathsf{int}\big(\Delta(C_1)\big)$.
\end{proof}

We are now ready to prove \zcref{thm:MazoitEndgame}.

\begin{proof}[Proof of \zcref{thm:MazoitEndgame}] 
If $g=0$ then we can apply \zcref{lem:outerplanar-routing}. Thus we may assume that $g \geq 1$.

Let $S_1, \dots, S_g$ be the strips of $\Sigma$ ordered such that $|S_i| \leq |S_{i+1}|,$ for all $1 \leq i < g,$ where $|S_i| \coloneqq |\LLL \cap S_i|$. We define $a \coloneqq \frac{10}{3}$.

If $|S_i| \geq k \cdot 3^i+1,$ for all $1 \leq i \leq g,$ then we can apply \zcref{lem:Mazoit-lem3} to obtain the required linkage. 

Otherwise, let $i_1$ be maximal such that $|S_{i_1}| < k \cdot 3^{i_1+1}$. Now suppose $i_j$ has already been defined. If there is some $i > 0$ such that $|S_{i_1 + \dots + i_j + i}| \leq k \cdot a^{j+1} \cdot i_1 \cdot i_2 \cdot \ldots \cdot i_j \cdot 3^{i_1 + \cdots i_{j} + i}$ then we choose $i_{j+1}$ as the maximal such number. Otherwise the construction concludes at this point. 

Let $\ell$ be maximal such that $i_\ell$ exists but $i_{\ell+1}$ does not. As explained above, $\ell \geq 1$. For $1 \leq j \leq\ell$ we define $d_j \coloneqq \sum_{x=1}^{j} i_x$ and $m_j:= \prod_{x=1}^{j} i_x$. Note the disparity between the indices in the upper bound for $|S_{i_1 + \ldots + i_{j}}|$: for $j \leq \ell$ we require $|S_{i_1 + \ldots + i_j}| \leq k \cdot a^j \cdot m_{j-1} \cdot 3^{d_j},$ i.\@e.\@~the sum in the exponent runs from $i_1$ to $i_j$ but the product before the exponential term only runs from $i_1$ to $i_{j-1}$.

We now consider $\tau(\LLL \cap \Sigma),$ which is a pattern of order $k'$ for some $k'$. We first establish a bound on $k'$.
\[
\begin{array}{rl}
 k' = & k\ +\ \big(|S_1| + \ldots + |S_{i_1}|\big)\ +\ \big(|S_{i_1+1}| + \cdots |S_{i_1 + i_2}|\big)\ + \ldots +\ \big(|S_{d_{\ell-1}+1} + \ldots + |S_{d_\ell}|\big) \\[1ex]
 \leq & k + i_1\cdot |S_{i_1}| + \ldots + i_\ell|S_{d_\ell}|\\
 & \textrm{(as $|S_{j'}| \leq |S_{i_j}|$ for all $1 \leq j \leq \ell$ and $j' \leq j$}\\[1ex]
 \leq & k\ +\ i_1 \cdot k \cdot a \cdot 3^{i_1}\ +\ i_2 \cdot k \cdot a \cdot i_1 \cdot 3^{i_1 + i_2}\ + \ldots +\ k \cdot i_\ell \cdot a^\ell \cdot m_{\ell-1} \cdot 3^{d_\ell} \\
 & \textrm{(as $|S_{d_j}| \leq k \cdot a^j \cdot m_{j-1}\cdot 3^{d_j},$ for all $j \leq \ell$)}\\[1ex]
 = & k\ +\ k \cdot a \cdot i_1\cdot 3^{i_1}\ +\ k \cdot a \cdot i_1 \cdot i_2 \cdot 3^{i_1 + i_2}\ +\ \ldots\ +\ k \cdot a^\ell \cdot m_\ell 3^{d_\ell} \\[1ex]
 = & k \cdot a^\ell \cdot m_\ell \cdot 3^{d_\ell} \cdot \Big(
 \frac{1}{a^\ell \cdot m_{\ell}\cdot 3^{d_\ell}} + \frac{1}{a^{\ell-1} \cdot i_2 \cdot \ldots \cdot i_{\ell}\cdot 3^{i_2 + \ldots + i_\ell}} + \ldots + \frac{1}{a \cdot i_\ell \cdot 3^{i_1}} + 1\Big)\\[1ex]
 \leq & k \cdot a^\ell \cdot m_\ell \cdot 3^{d_\ell} \cdot \Big(1 + \frac1{3a} + \frac1{(3a)^2} + \ldots + \frac1{(3a)^\ell}\Big)\\[1ex]
 \leq & k \cdot a^\ell \cdot m_\ell \cdot 3^{d_\ell} \cdot \frac{1}{(3a)^\ell}\cdot \frac{(3a)^{\ell+1} - 1}{3a-1}\\
 \leq & k \cdot \frac1{3^\ell}\frac{(3a)^{\ell+1}-1}{3a-1}\cdot m_\ell \cdot 3^{d_\ell}.
\end{array}
\]

We show next that $|S_{d_\ell+j}| > 3k'3^j,$ for all $j > 0,$ so that we can apply \zcref{lem:Mazoit-lem3}.
The maximality of $\ell$ implies that for all $j > 0$ we have

\begin{align*}
 3 \cdot k' \cdot 3^j & = \frac{3a - 1}{3} \cdot k' \cdot 3^j\qquad \textrm{(as $a = \frac{10}{3}$)}\\
 & < \frac{(3a)^{\ell+1}}{(3a)^{\ell + 1} -1}\cdot\frac{3a-1}{3}\cdot k' \cdot 3^j\\
 & \leq \frac{(3a)^{\ell+1}}{(3a)^{\ell + 1} -1}\cdot\frac{3a-1}{3}\cdot \Big( k \cdot \frac1{3^\ell}\cdot\frac{(3a)^{\ell+1}-1}{3a-1}\cdot m_\ell \cdot 3^{d_\ell}\Big) \cdot 3^j\\
 & = k \cdot a^{\ell+1} \cdot m_\ell \cdot 3^{d_\ell+j} \\
 & \leq 
 |S_{d_\ell + j}|
\end{align*}

We now want to apply \zcref{lem:Mazoit-lem3} to the instance $\big((\Delta, S_{d_\ell+1}, \dots, S_g), G, \CCC, \LLL \setminus (S_{d_\ell+1} \cup \ldots \cup S_g)\big)$. Let $g' \coloneqq g - d_\ell,$ which is the number of strips in the instance. What is left is to calculate the size of $\CCC$ so that the lemma yields an irrelevant vertex. For this we need to show that $t \geq 2k'\cdot 3^{g'}$.

\begin{align*}
 2k'\cdot 3^{g'} & \leq 2\cdot \Big(k\cdot \frac{1}{3^\ell} \cdot \frac{(3a)^{\ell+1}-1}{3a-1} \cdot m_\ell \cdot 3^{d_\ell}\Big)\cdot 3^{g'} \\
 & \leq 2k \cdot a^{\ell+1}\frac{3}{3a-1}\cdot m_\ell \cdot 3^{d_\ell} \cdot 3^{g'}\\
 & = 2k \cdot a \cdot \frac{3}{3a-1}\cdot a^{\ell} \cdot m_\ell \cdot 3^g\\
 & \leq \frac{20}{9} k \cdot \Big(a\cdot \frac{d_\ell}{\ell}\Big)^\ell \cdot 3^g\\
 & \leq \frac{20}{9}k \cdot \Big(\frac{a\cdot g}{\ell}\Big)^\ell 3^g \qquad \textrm{(as $d_\ell \leq g$)}
\end{align*}

The factor $\big(\frac{a\cdot g}{\ell}\big)^\ell$ is maximised for $\ell = \frac{a \cdot g}{e}$. Therefore $2k'\cdot 3^{g'} \leq \frac{20}9k\cdot e^{\frac{10\cdot g}{3e}}\cdot 3^g = f(g, k)$. Thus, \zcref{lem:Mazoit-lem3} implies that all vertices in $\mathsf{int}\big(\Delta(C_1)\big)$ are irrelevant. This concludes our proof. 
\end{proof}

\section{Linking vortices up to some gap}\label{sec:linkedness}

\begin{figure}
 \centering
\includegraphics[width=0.55\linewidth]{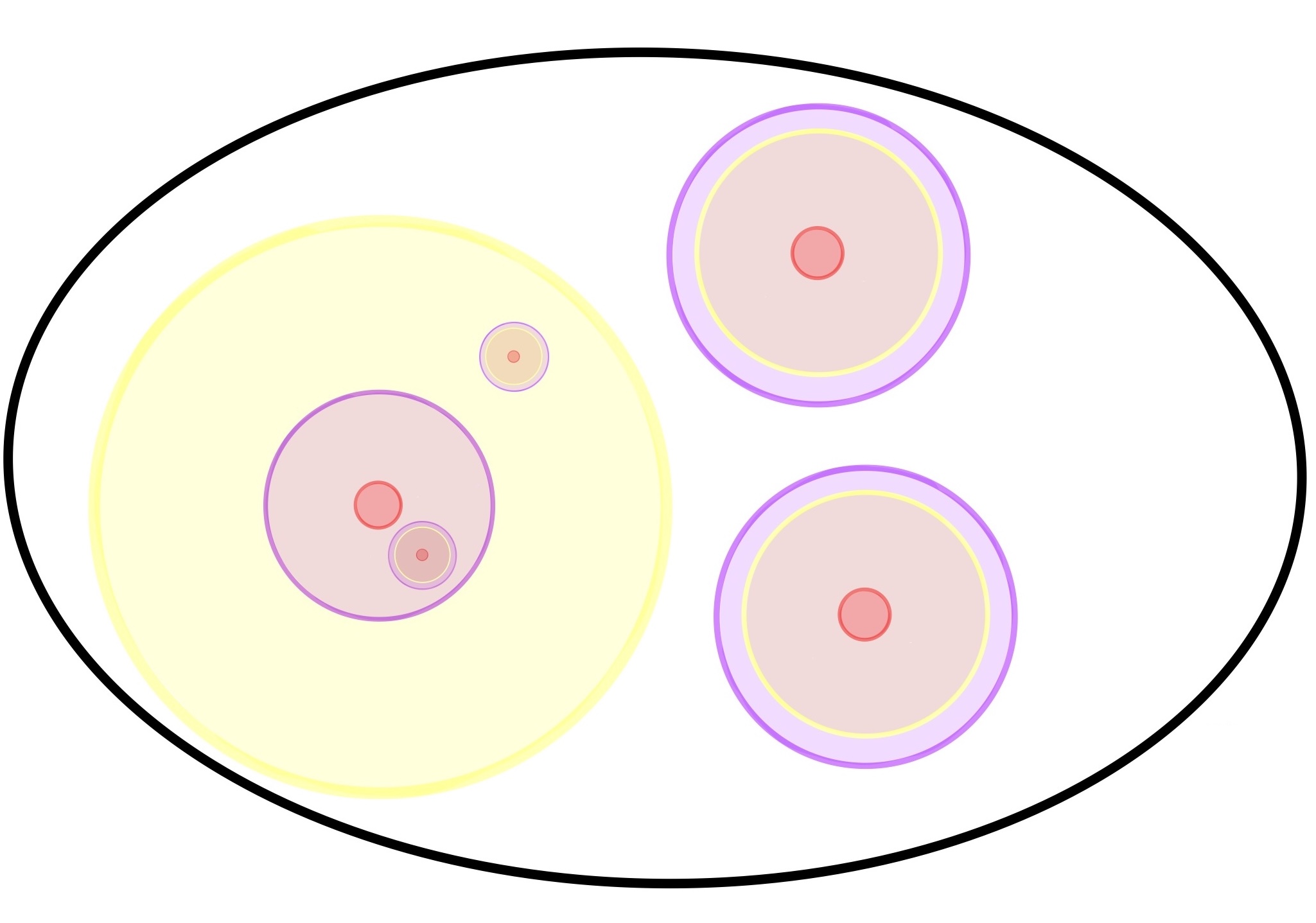}
 \caption{The setting for linking vortices. Red cells denote vortices, each of which come with two layers of concentric cycles: purple areas are schematic representations of $(d+1)$ concentric cycles and yellow areas are schematic representations of the remaining concentric cycles. The cycle whose outer-layer is yellow is the vortex one currently tries to link, where the inner purple layer is used to push cycles deeply into the vortex.}
 \label{fig:vortices_setting}
\end{figure}

This section is dedicated to the proof of a version of \cite[Lemma 15]{DiestelKMW2012Excluded} suitable for our needs: We want to change our rendition slightly in order to guarantee that the graphs of vortices are \emph{$d$-linked up to some bounded gap}. Intuitively speaking, this means that there is a linear decomposition of some segment of the vortex society, that ``covers'' all but a few society vertices, and such that the adhesion of that linear decomposition is $d$ and there is a linkage of order $d$ that ``links'' the first to the last bag of the decomposition. We will need this to apply a variant of \cref{thm:Lemma2.5} in a later section. The main idea to achieve linkedness being that, if we find $d+1$ concentric cycles around a vortex of depth $d,$ we may try to ``push'' the cycles into the vortex as deeply as possible in order to ``calibrate'' regions of the vortex that are of depth lower than $d$ by bloating these up artificially. Unfortunately, we cannot easily guarantee all the assumptions needed for the referenced theorem; more precisely, we may not be able to find ``private'' nests of ``large enough'' size for every vortex. Worse even, it may happen that the nest of one vortex is contained in the disc bounded by a cycle of the nest of another vortex~--~this will be one outcome of what we will call \emph{stacked}~--~without yielding immediate contradictions to any ``distance'' assumptions between vortices; see \cref{fig:vortices_setting} for a schematic representation of the setting. 
It turns out that this is, however, not really a problem when trying to link vortices; it may simply happen that one vortex is ``absorbed'' by another without increasing the depth of the former by too much.

\subsection{Linked linear decompositions}

We will need renditions where vortices are \emph{linked up to some bounded gap}; we transcribe the following definition of \cite{DiestelKMW2012Excluded} to our setting. 

\begin{definition}\label{def:linked_vortex_normal}
 Let $d \in \N$ and $m \geq 1$.
 Let $(H,\Omega)$ be a society.
 Let $\alpha=\left(\langle u_1,\ldots,u_{m}\rangle,\langle U_1,\ldots,U_m\rangle\right)$ be a linear decomposition of $(H,\Omega)$.
 
 Let $1 \leq p \leq q \leq m$ and let $S=[p,q]$. Let $\alpha_S$ be the $S$-restriction of $\alpha$. Let $Z_i \coloneqq (U_i \cap U_{i+1})\setminus V(\Omega_c)$ for $i \in [p,q-1]$. We call $\alpha_S$ \emph{$d$-linked} if 
 \begin{enumerate}
 \item for every $i \in [p,q],$ $u_i \in U_i$ and $U_i \cap V(\Omega) \subseteq \{u_{i-1},u_{i}\},$ where $u_0 \coloneqq u_1$
 \item for every $ i \in [p,q-1],$ $\Abs{Z_i} = d,$ and
 \item there is a $U_{p}$-$U_q$-linkage $\mathcal{M}$ in $G[\bigcup_{i=p}^qU_i]\setminus V(\Omega)$ of order $d$.
 \end{enumerate}
 Note that by definition for every $ i \in [p,q-1],$ it holds $U_i \cap U_{i+1} \subseteq Z_i \cup \{u_i\}$. We say that $\mathcal{M}$ \emph{witnesses that $(H,\Omega)$ is $d$-linked}.

 If additionally the following holds
 \begin{enumerate}[resume]
 \item there is a path $P^*\subseteq G$ disjoint from $\mathcal{M}$ such that $N_{\rho}(P^*) = V(\Omega),$
 \end{enumerate}
then we call $\alpha$ \emph{$d$-stuffed}. We say that $\mathcal{M}\cup\{P^*\}$ \emph{witnesses that $(H,\Omega)$ is $d$-stuffed}. Finally, we say that $\alpha$ is \emph{properly $d$-linked} if $\MMM \cup \{P^*\} \subseteq G[\bigcup_{i=p}^qU_i]$.
\end{definition}
If $(H,\Omega)=(\sigma(c),\Omega)$ is a vortex society for some vortex $c \in C(\rho)$ given some $\Sigma$-rendition of some graph $G$ in a surface $\Sigma,$ then we say that $\mathcal{M}$ \emph{witnesses that $c$ is $d$-linked} and $\mathcal{M}\cup\{P^*\}$ \emph{witnesses that $c$ is $d$-stuffed}, respectively, for simplicity.

While properly $d$-linkedness can in general not be easily guaranteed -- we work with renditions and paths cannot be easily ``tucked'' into vortices -- we will need it for technical reasons in \cref{sec:noapexnocry} where we essentially reduce to an embedded instance. For completeness, we give the definition here. We extend the definition of (properly) linked and stuffed linear decompositions, to linear decompositions with gaps.
\begin{definition}\label{def:linked_vortex}
 Let $d \in \N,$ and $\ell,m \geq 1$.
 Let $(H,\Omega)$ be a society. 
 Let $u_{-(\ell-1)},\ldots,u_0,u_1,\ldots,u_{m}$ be the vertices of $V(\Omega)$ appearing in $\Omega$ in the order listed together with a sequence of sets $\langle U_{-(\ell-1)},\ldots,U_0,U_1,\ldots,U_m\rangle$ such that
 \begin{enumerate}
 \item $U_i\subseteq V(H)$ and $u_i\in U_i$ for all $i\in[-(\ell-1),m],$
 \item $\bigcup_{i\in [-(\ell-1),m]}U_i=V(H)$ and for every $uv\in E(H)$ there exists $i\in[-(\ell-1),m]$ such that~$u,v\in U_i,$
 \item for every $x\in V(H),$ the set $\{ i\in[m] ~\!\colon\!~ x\in X_i \}$ forms an interval in $[m]$ and the set $\{ i\in[-(\ell-1),0])] ~\!\colon\!~ x\in X_i \}$ forms an interval in $[-(\ell-1),0],$
 \item let $X^- \coloneqq \bigcup_{i=-(\ell-1)}^0 U_i $ and $X^+ \coloneqq \bigcup_{i=1}^m U_i,$ then $X^- \cap X^+ \subseteq (U_m \cap U_{-(\ell-1)}) \cup (U_0 \cap U_1)$. 
\end{enumerate}
 
 Let $\alpha=\left(\langle u_{-(\ell-1)},\ldots,u_0,u_1,\ldots,u_{m}\rangle,\langle U_{-(\ell-1)},\ldots,U_0,U_1,\ldots,U_m\rangle\right),$ then we call $\alpha$ a \emph{linear decomposition of $(H,\Omega)$ with gap $\ell$}.
 
 Let $-(\ell-1) \leq p \leq q \leq m$ and let $S=[p,q]$. We call $\alpha_S=\left(\langle u_{p},\ldots,u_q\rangle,\langle U_{p},\ldots,U_q\rangle\right)$ the $S$-restriction of $\alpha$. We call $\alpha_S$ \emph{(properly) $d$-linked} or \emph{$d$-stuffed} if $\alpha_S$ is a (properly) $d$-linked or $d$-stuffed linear decomposition of $(G[\bigcup_{i=p}^q U_i],\Omega'),$ respectively, where $\Omega'$ is the restriction of $\Omega$ to $\{u_i \mid i \in [p,q]\}$. We say that $\alpha_S$ is (properly) $d$-linked, respectively $d$-stuffed, with gap $\ell' \leq \ell$ if $p=-(\ell'-1)$ and $q\geq 1$ and the $[1,q]$-restriction of $\alpha$ is (properly) $d$-linked, respectively $d$-stuffed.
 
 We say that $\alpha$ is \emph{(properly) $d$-linked, respectively $d$-stuffed, with gap $\ell$} if the $[1,m]$-restriction of $\alpha$ is (properly) $d$-linked, respectively $d$-stuffed.
\end{definition}
Note that linear decompositions with gaps are special cases of so called cyclic decompositions; see for example \cite{DiestelJKK2025Canonical}.

Given some vortex $c$ in some $\Sigma$-rendition $\rho,$ we may say that a linear decomposition of $(\sigma(c),\Omega_c)$ is a \emph{$d$-linked linear decomposition with gap $\ell$ of $c$} for simplicity, not stating the exact choice of $\Omega_c$ (there are two up to cyclic rotation).
We say that $(\sigma(c),\Omega_c)$ or $c$ is \emph{(properly) $d$-linked, respectively $d$-stuffed, with gap $\ell$} if it admits a (properly) $d$-linked, respectively $d$-stuffed, linear decomposition with gap $\ell,$ and so on.

\smallskip

Furthermore, we extend the definition of concentric cycles to more general renditions as follows.
\begin{definition}\label{def:concentric_cycles}
 Let $\Sigma$ be a surface, $s \geq 1$ and $\rho$ a $\Sigma$-rendition of some graph $G$. Let $\langle C_1,\ldots,C_s\rangle$ be a sequence of disjoint cycles in $G,$ then we call them \emph{concentric in $\rho$} if there exists a $\rho$-aligned disc $\Delta \subseteq \Sigma$ such that $\langle C_1,\ldots,C_s\rangle$ are concentric in the restriction $\rho^*$ of $\rho$ by $\Delta$. Further, for every $i \in [s]$ we let $\Delta_{C_i}$ denote the respective $C_i$-disc in $\rho^*$ (and hence in $\rho$).
\end{definition} 
Note that whenever $s \geq 2$ the $C_i$-disc is unambiguously defined via the order $\langle C_1,\ldots,C_s\rangle$ even if $\Sigma$ is the sphere.

Finally, we slightly extend the definition of \emph{crop} of a graph by a disc as follows.

\begin{definition}\label{def:extended-crop}
 Let $\Sigma$ be a surface and let $\rho$ be a $\Sigma$-rendition of some graph $G$. Let $C\subseteq G$ be a grounded cycle and let $\Delta_C$ be a $\rho$-aligned disc bounded by the trace of $C$. We define $G[\rho,\Delta_C]$ as \emph{the extended crop of $G$ by $\Delta_C$} obtained from the crop $G_{\Delta_{C}}$ of $G$ by $\Delta_{C},$ by adding to it the cycle $C,$ i.e., $G[\rho,\Delta_C] \coloneqq G_{\Delta_{C}} \cup C$.
 
 If $\Delta_C$ is uniquely defined by $C$ or is clear from context, we write $G[\rho,C]$ for simplicity, and refer to it as the \emph{extended crop of $G$ by $C$}. 
\end{definition}
Note that if $\Sigma$ is not the sphere and the trace of $C$ bound a disc in $\Sigma,$ then the extended crop of $G$ by $C$ is unambiguously defined.

\subsection{Kissing vortices}

As a first lemma will need an adaptation of \cite[Lemma 15]{DiestelKMW2012Excluded}; the main difference being that we do not want to introduce new (major) apices when changing the rendition. Following their proof yields a possible major vertex $u$ say\footnote{Note that all but one of the resulting apices in their proof have all their original neighbours in the graph of the newly constructed vortex, hence they are minor apices.}. This can easily be mitigated by introducing a ``small gap'' where we do not impose as strict restrictions on the behaviour inside the graph of the vortex by simply redrawing $u$ on the resulting vortex boundary and adding a bag to the linear decomposition that takes care of all the formerly introduced apices. While the modification is easy~--~put all the deleted vertices into a single bag~--~we provide a proof for completeness. Note that the main part of the proof simply follows the proof of \cite[Lemma 15]{DiestelKMW2012Excluded}, thus we omit some details by transcribing relevant claims from their setting to ours. Furthermore, note that they work with what they call ``properly attached cells''; this is mitigated in our setting by working with grounded paths and nests.

We define the following, giving a name to a set of paths ``pushed tightly into a vortex''.
\begin{definition}\label{def:kissing}
 Let $\Sigma$ be a surface and $\rho$ a $\Sigma$-rendition of depth $d$ of some graph $G$. Let $c \in C(\rho)$ be a vortex. Let $s \geq 1$ and let $\mathcal{C} = \langle C_1, \ldots, C_{s}\rangle$ be concentric cycles in $\rho$ such that~$c \subseteq \Delta_{C_1}$. Let~$G[\rho,\Delta_{C_s}]$ be the extended crop of $G$ by $\Delta_{C_s}$.

 Let $\BBB=\{B_1,\ldots,B_s\}$ and $\PPP=\{P_1,\ldots,P_s\}$ be linkages in $G[\rho,\Delta_{C_s}]$ such that for distinct $i,j \in [s]$ the paths $P_i, B_j$ are internally disjoint and $\{B_i\cup P_i \mid i\in [s]\}$ is a set of $s$ disjoint cycles where~$B_i \subsetneq C_i$ for every $i \in [s]$. We call $\BBB$ a \emph{base of $\PPP$ in $\CCC$}.
 
 Let $C^* \subseteq G[\rho,\Delta_{C_s}]$ be a grounded cycle in $\rho$ and let $\Delta^*$ be a disc bounded by the trace of $C^*$ containing $c$. We say that $C^*$ \emph{$\rho$-envelopes $\PPP$ at $\BBB$} if $\PPP \cup \BBB$ is contained in the extended crop of $G$ by $\Delta_{C^*}$ and additionally $C^* = P_{s} \cup B_s;$ we call $(P_s,B_s)$ the \emph{perimeter} of $(\PPP,\BBB)$ and $\Delta^*$ the \emph{$C^*$-disc} (in $\rho$ containing $c$). 
 
 We say that $(\PPP,\BBB)$ \emph{$\rho$-kisses $c$ in $\CCC$} if for every linkage $\PPP'$ of order $s$ such that $\BBB$ is a base of $\PPP'$ in $\CCC$ and such that $C'$ $\rho$-envelopes $\PPP$ at $\BBB,$ it holds that $\Delta_{C^*}$ is contained in the respective~$C'$-disc.
\end{definition}

\begin{figure}
 \centering
 \begin{subfigure}[t]{0.5\textwidth}
 \centering
 \includegraphics[width=0.95\linewidth]{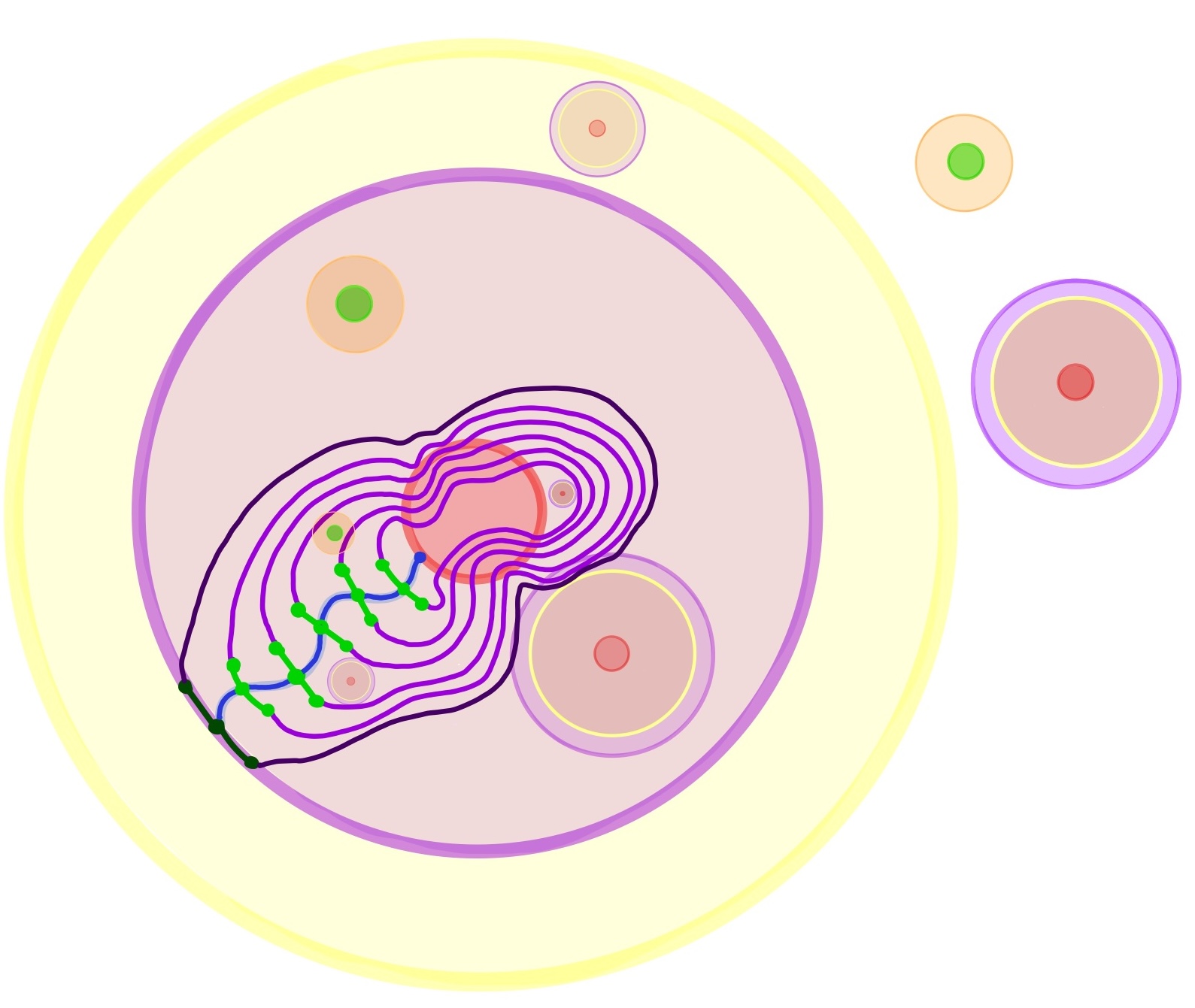}
 \caption{$(\PPP,\BBB)$ kissing a vortex.}
 \end{subfigure}%
 ~ 
 \begin{subfigure}[t]{0.5\textwidth}
 \centering
 \includegraphics[width=0.95\linewidth]{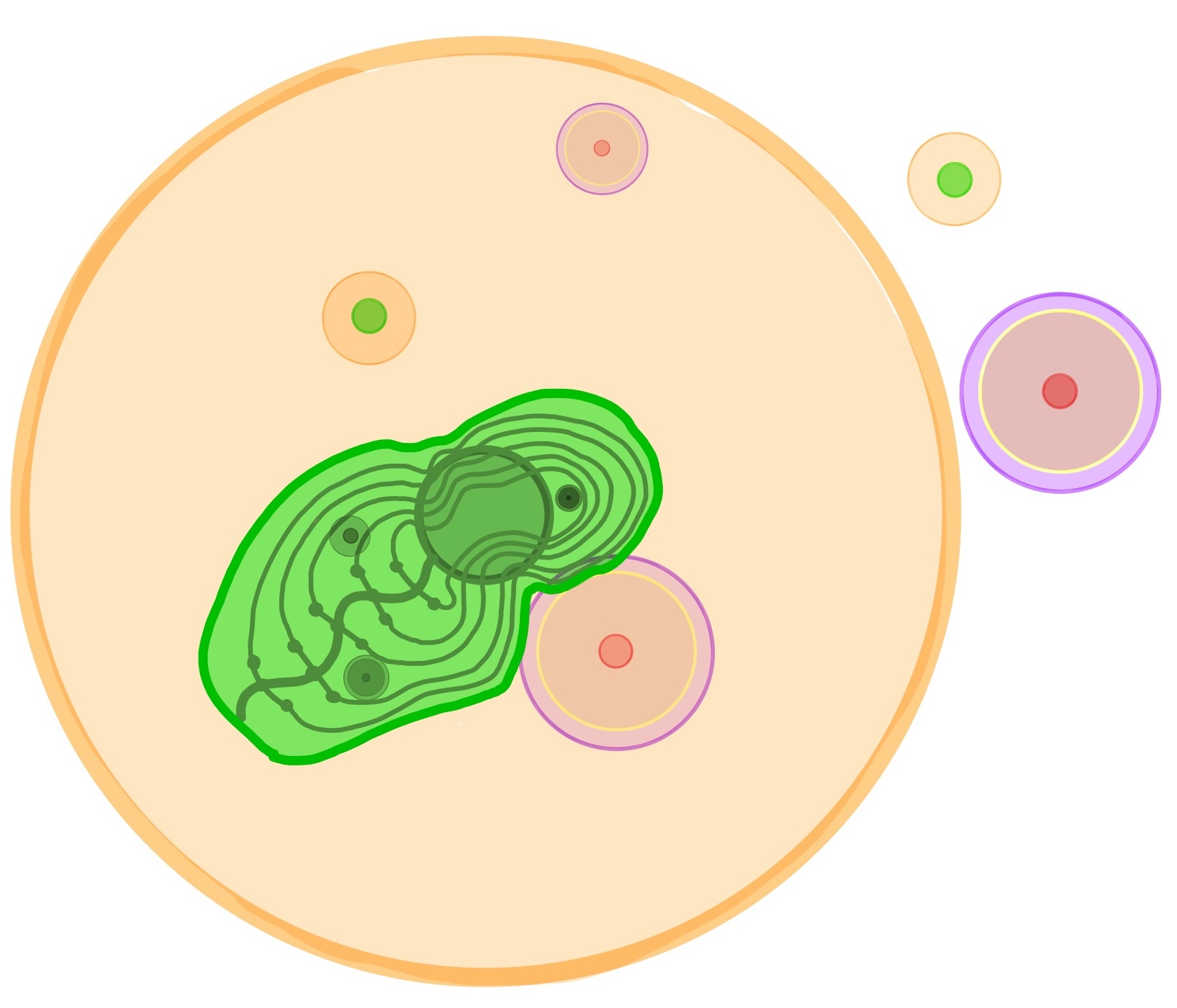}
 \caption{The resulting linked vortex.}
 \end{subfigure}
 \caption{A schematic illustration of the proof of \cref{thm:linking_vortices}. Non-linked vortices are red together with purple regions marking $(d+1)$ concentric cycles and yellow regions marking the remaining concentric cycles. Linked vortices are green together with their remaining set of cycles in orange. All vortices are stacked. On the left hand side a $(d+1)$ linkage $\PPP$ (purple) together with a base $\BBB$ (green) $\rho$-kissing a vortex are depicted. On the right hand side, the resulting linked vortex is depicted, where some vortices have been absorbed.}
 \label{fig:linking_vortices_proof}
\end{figure}
Note that in general, $C$-discs in $\rho$ are only defined if $\Sigma$ is not the sphere: This is because if it is a sphere there are two natural choices. We may however talk about the $C$-disk \emph{containing $X \subseteq \Sigma$} which then unambiguously defines it, i.\@e.\@, if the trace of $C$ bounds a disc $\Delta$ in $\Sigma$ such that $X \subseteq \Delta,$ then we call $\Delta$ the \emph{$C$-disc (in $\rho$) containing $X$}.

We have the following version of \cite[Lemma 15]{DiestelKMW2012Excluded}.
\begin{lemma}\label{lem:linking_vortices}
 Let $d \in \N$. Let $\Sigma$ be a surface and $\rho$ a $\Sigma$-rendition of depth $d$ of some graph $G$.
 Let $c \in C(\rho)$ be a vortex. Let $\Delta$ be a closed $\rho$-aligned disc such that the restriction $\rho_c$ of $\rho$ to $\Delta$ is a cylindrical rendition around $c$ of the $\Delta$-society with a nest $\mathcal{C} = \langle C_1, \ldots, C_{d+1}\rangle$ such that $\CCC$ is tight around $c$.

 Then there exists a $\Sigma$-rendition $\rho'$ of $G$ together with a cell $c^*\in C(\rho')\setminus C(\rho)$ such that $\rho'$ agrees with $\rho$ on $\Sigma \setminus c^*$ and $c \subseteq c^* \subseteq \Delta$ and either $\Abs{N(c^*)}\leq 3$ or the $c^*$-society $(\sigma_{\rho'}(c^*),\Omega_{c^*})$ is $d$-stuffed with gap at most $1$ and the depth of $c^*$ is at most $4d+2$. In the latter case, there exist linkages $\PPP=\{P_1,\ldots,P_{d+1}\}$ and $\BBB=\{B_1,\ldots,B_{d+1}\}$ of order $(d+1)$ in $G$ satisfying the following:
 \begin{enumerate}
 \item $\{P_1,\ldots,P_d\}$ witnesses that $c^*$ is $d$-linked with gap at most $1,$ 
 \item $P_{d+1}$ satisfies $V(P_{d+1}) \cap N(\rho') \subseteq V(\Omega_{c^*}),$ in particular $\PPP$ witnesses that $c^*$ is $d$-stuffed with gap at most $1,$ and
 \item for every $i \in [d+1]$ it holds $B_i \subsetneq C_i$ and $\BBB$ is a base for $\PPP,$
 \end{enumerate} 
 such that $C^* \coloneqq P_{d+1} \cup B_{d+1}$ $\rho$-envelopes $\PPP$ at $\BBB$ (it may not be grounded in $\rho'$) with $V(C^*) \cap N(\rho') = V(\Omega_{c^*})$ and $(\PPP, \BBB)$ $\rho$-kisses $c$ in $\CCC$.
\end{lemma}
\begin{proof}
 Let $c,\CCC,\rho$ and $\Delta$ be as in the lemma. Since the traces of all the cycles of the nest are contained in a common disc $\Delta,$ we may choose a common orientation of their traces, i.\@e.\@, orient all of them counter-clockwise say, and by the choice of nest every cycle $C_i$ comes with a uniquely defined disc $\Delta_i$ containing $c_i$.

 If $T \coloneqq N(\rho) \cap V(C_{d+1}) \leq 3$ simply define $c^* \coloneqq \Delta_{C_{d+1}}\setminus T;$ by construction $c^* \subseteq \Delta$. Let $\rho'$ be obtained from $\rho$ by deleting all cells contained in $c^*$ and adding the cell $c^*$. Finally, let $\rho'$ agree with $\rho$ on all cells in $C(\rho') \cap C(\rho)$ and let $\sigma_{\rho'}(c^*)$ be the graph of the $\Delta_{C_{d+1}}$-society in $\rho$. Then $\rho'$ is as desired; henceforth assume that $N(\rho) \cap V(C_{d+1}) \geq 4$.
 
 Let $(\sigma(c),\Omega_c)$ be a $c$-society fixing a cyclic order $\Omega_c$ visiting the nodes $u_1,\ldots,u_n$ of $c$ in that order. Let $\alpha=(\langle u_1,\ldots,u_n\rangle,\langle U_1,\ldots,U_n\rangle)$ be a linear decomposition of $(\sigma(c),\Omega_c)$ of adhesion at most $d$ which exists by \cref{prop:DepthToLinearDec}. 

 Pick a vertex $v_{d+1} \in V(C_{d+1})$ and let $\gamma$ be a shortest $\rho$-aligned curve~--~where length is taken with respect to the number of nodes it intersects~--~from $v_{d+1}$ to some node in $N(c),$ say $u_i$ for some respective $i \in [n],$ internally disjoint from $c$. Then, using tightness of $\CCC,$ $\gamma$ intersects at most $d+2$ nodes; see the \textcolor{ceruleanblue}{blue} \textsw{squiggly} line on the left hand side of \cref{fig:linking_vortices_proof}. Let $\langle v_{d+1},v_d,\ldots,v_1,u_i\rangle$ be the order of the nodes visited when traversing $\gamma$ starting in $v_{d+1}$. Then $v_j\in V(C_j)$ for every $j \in [d+1]$ where $v_1=u_i$ possibly but $v_j \neq v_\ell$ for distinct $j,\ell \in [d+1]$. Let $F \coloneqq \{v_{d+1},\ldots,v_1,u_i\}$.

 Let $Z_i \coloneqq U_i \cap U_{i-1},$ then $\Abs{Z_i} \leq d$ since the adhesion of $\alpha$ is at most $d$. Let $G'$ be obtained from~$G$ by deleting $F \cup Z_i;$ these are the vertices we will need to ``add back'' later. Since the nest $\CCC$ is grounded in $\rho$ by definition, for $j \in [d+1]$ the vertex $v_j$ has a unique predecessor $x_j \in V(C_j)\cap N(\rho)$ and successor $y_j \in V(C_j)\cap N(\rho)$~--~they are not equal since $\Abs{N(\rho) \cap V(C_j)}\geq 4$~--~with respect to the counter-clockwise orientation of the trace of $C_j,$ such that there are a grounded path $P_j^x \subsetneq C_j$ from $x_j$ to $v_j$ not visiting $y_j$ as an internal vertex, and a grounded path $P_j^y \subsetneq C_j$ from $v_j$ to $y_j$ not containing $x_j$ as an internal vertex, and such that no internal vertex of either path is a node of $\rho$. Fix $X \coloneqq \{x_j \mid j\in [d+1]\}$ and $Y \coloneqq \{y_j \mid j \in [d+1]\}$. Define $\BBB \coloneqq \{P_j^x \circ P_j^y \mid j \in [d+1] \}$ to be the set consisting of the paths obtained by concatenating $P_j^x $ and $ P_j^y$ at $v_j;$ $\BBB$ is an $X$-$Y$-linkage in $G$ by construction where each path has length at least $1$ since $x_j \neq y_j$ for all $j\in [d+1]$. Let $\BBB= \{B_1,\ldots,B_{d+1}\}$ where $B_j$ is the path with endpoints $x_j,y_j$ for $j \in [d+1];$ see the \textcolor{darkspringgreen}{green} paths on the left hand side of \cref{fig:linking_vortices_proof} for a schematic illustration of $\BBB$.
 
 Let $\PPP_0=\{P_1,\ldots,P_{d+1}\}$ be the unique $X$-$Y$-linkage in $G'$~--~and thus in $G$~--~obtained from $\CCC$ by deleting $F$ in the obvious way, in particular $P_j \subsetneq C_j$ and $E(P_j) = E(C_j) \setminus \bigcup_{i=1}^{d+1}E(P_i^x \cup P_i^y)$ for $j \in [d+1]$. Then $\BBB$ is a base for $\PPP_0$ in $G$ by construction and additionally $\{P_j \cup B_j \mid j\in [d+1]\} = \{C_j \mid j \in [d+1]\}$. Thus, $C_{d+1}$ $\rho$-envelopes $\PPP_0$ at $\BBB$ as in \cref{def:kissing}, whence $(P_{d+1},B_{d+1})$ is the perimeter of $(\PPP_0,\BBB)$. 

 We have the following.
 \begin{beautifulclaim}\label{claim:the_paths_are_enveloped}
 Let $\PPP$ be an $X$-$Y$-linkage in $G'$. Then $\BBB$ is a base for $\PPP$ in $G$.
 \end{beautifulclaim}
 \begin{claimproof}
 This is immediate from the fact that $\BBB$ is internally disjoint from $V(G')$ and every path in $\BBB$ shares its endpoints with a unique path in $\PPP$.
 \end{claimproof}

 The following summarises Equation $(4)$ and its implications as given in the proof of \cite[Lemma 15]{Diestel2010Graph} transcribed to our setting; again the paths being grounded in $\rho$ take over the role of properly attached cells in their setting.
 \begin{beautifulclaim}\label{stolen_uno}
 For every set $\PPP$ of $d + 1$ disjoint $X$–$Y$-paths in $G'$ whose traces are contained in $\Delta,$ the path $P$ starting in $y_{d+1}$ is grounded in $\rho$ and ends in $x_{d+1}$. 
 \end{beautifulclaim}

 Combining \cref{stolen_uno} and \cref{claim:the_paths_are_enveloped} we derive the following.

 \begin{itemize}
 \item[$(\star)$] Let $\PPP$ be an $X$-$Y$-linkage in $G'$ such that the trace of every path in $\PPP$ is contained in $\Delta$. Let $P \in \PPP$ be the unique path starting in $y_{d+1}$ and ending in $x_{d+1}$. Then $C(P) \coloneqq P\cup B_{d+1}$ is a cycle that $\rho$-envelopes $\PPP$ at $\BBB,$ in particular $(P,B_{d+1})$ is the perimeter of $(\PPP,\BBB)$. 
 \end{itemize}
 
 Note that $\PPP_0$ and $P_{d+1}$ satisfy \cref{stolen_uno}, and given the respective cycle $C(P_{d+1})$ by \cref{stolen_uno}, the $C(P_{d+1})$-disc $\Delta(P_{d+1})$ in $\rho$ containing $c$ satisfies $c \subseteq \Delta(P_{d+1})$ and the trace of every path in $\PPP_0$ is contained in $\Delta(P_{d+1})$. Let $\PPP, P$ and $C(P)$ be as in $(\star)$ and let $\Delta(P)$ be the $\rho$-aligned $C(P)$-disc. We define $G'[\rho,C(P)]$ as the extended crop of $G'$ by $C(P)$. As in the proof of \cite[Lemma 15]{DiestelKMW2012Excluded} let $\PPP$ and $P^* \in \PPP$ be chosen to satisfy \cref{stolen_uno} -- hence~$(\star)$ -- minimizing $G'[\rho,C(P^*)];$ see the left hand side of \cref{fig:linking_vortices_proof} for a schematic illustration of such a linkage $\PPP,$ forgetting about the remaining vortices in the figure.

 By our choice of $\PPP$ and $P^*,$ together with $(\star)$ and the \cref{def:kissing} of $\rho$-kissing, we derive the following.

 \begin{itemize}
 \item[$(\star\star)$] $(\PPP,\BBB)$ $\rho$-kisses $c$ at $\CCC$ and $(P^*,B_{d+1})$ is the perimeter of $(\PPP,\BBB)$. 
 \end{itemize}

 Let $\Omega'=\langle y_{d+1},w_1,\ldots,w_k,x_{d+1}\rangle$ be a cyclic ordering of $V(P^*)\cap N(\rho) = \{w_1,\ldots,w_k,x_{d+1},y_{d+1}\}$ respecting the orientation of the trace of $P^*$ when starting at $y_{d+1}$ and ending at $x_{d+1}$ for some respective $k \in \N;$ note that $(G'[\rho,C(P^*)],\Omega')$ defines a society. Then Diestel, Kawarabayashi, Müller, and Wollan. \cite[Lemma 15]{DiestelKMW2012Excluded} prove that one can slightly shrink $\Delta(P^*)$ to obtain a new disc $\Delta' \subseteq \Delta(P^*) \subseteq \Delta$ satisfying the following (transcribed to our setting).

 \begin{beautifulclaim}\label{stolen_two}
 There is a rendition $\rho'$ of $G'=G-(F\cup Z)$ of depth at most $2d+2,$ that agrees with~$\rho$ outside of $\Delta',$ and such that $\Delta'\setminus (V(P^*)\cap N(\rho))$ is a vortex $c'$ of $\rho'$ with $\sigma_{\rho'}(\Delta') \subseteq G[\rho,C(P^*)]$ and $P^* \cap \sigma_{\rho'}(\Delta') = V(P^*) \cap N(\rho')$. Further $(\sigma_{\rho'}(\Delta'),\Omega')$ is a $c'$-society that admits a $d$-stuffed linear decomposition $\alpha'=\left(\langle y_{d+1},w_1,\ldots,w_k,x_{d+1}\rangle,\langle Y_{d+1},W_1,\ldots,W_k,X_{d+1}\rangle\right),$ witnessed by $\PPP,$ where the $d$-linkedness is witnessed by $\PPP\setminus\{P^*\}$.
 \end{beautifulclaim}
 Note that the depth is at most $2d+2$ since the adhesion of the linear decomposition is $\leq d+1$. See the right hand side of \cref{fig:linking_vortices_proof} for a schematic illustration of the resulting vortex.

We continue with the outcome of \cref{stolen_two}. Note that when constructing $\Delta'$ in their proof, they ``push'' certain cells of $\rho$ out of the new vortex that a priori intersect $\Delta(P^*)$. Then, the reason why the path $P^*$ visiting the society vertices exists in their setting is guaranteed by the fact that said cells are properly attached, whereas in our setting this follows from the fact that the path was grounded in $\rho$ to begin with (this is simply a matter of choice of definitions for the setting). Note that after this construction, the path $P^*$ is in general not grounded in $\rho'$ and may be partially contained in the graph of the new vortex. 

Up to here nothing~--~besides transcribing to our setting~--~inherently changed when comparing to \cite[Lemma 15]{DiestelKMW2012Excluded}; we are left to ``add back $F \cup Z$'' to $G'$ by allowing a gap.

 We change $\rho'$ slightly to a rendition of $G,$ adding back $F \cup Z$ as follows. Note that by \cref{stolen_two} the nodes of $x_{d+1}$ and $y_{d+1}$ in $\rho'$ are the same points in $\Sigma$ as the respective nodes in $\rho$. By construction of $\Delta'$ the curve $\gamma\subsetneq \Sigma$ with nodes $F$ intersects the boundary of $\Delta' \subseteq \Delta(P^*) \subseteq \Delta(P_{d+1})\subseteq \Delta$ in $\Sigma$. Let $\xi$ be the intersection point of $\gamma$ and $\Delta'$ closest to the original node of $\pi_{\rho}^{-1}(v_{d+1})$ in $\Sigma,$ possibly equal to it. Fix $\pi_{\rho'}(\xi) = v_{d+1}$ and denote the resulting vortex via $c^*\coloneqq c' \setminus \{\xi\}$. Fix $$\sigma_{\rho^*}(c^*)\coloneqq G[V(\sigma_{\rho'}(c)) \cup (Z \cup F)],$$ noting that, by construction of $\Delta',$ all cells of $\rho$ that had a vertex of $F \setminus \{v_{d+1}\}$ as a node are completely contained in $\Delta',$ and $Z \subseteq V(\sigma_{\rho}(c))$ by definition. Finally, since all the cells of $\rho$ that had $v_{d+1}$ as a node are either completely contained in $\Delta'$ or agree with the cells of $\rho',$ this yields a valid $\Sigma$-rendition $\rho^*$ of $G$ with the desired properties as in the lemma. 
 
 Note that the depth of $\rho^*$ may only be marginally higher than that of $\rho'$: Since the depth of $\rho'$ is at most $2d+2$ and $\Abs{Z \cup F}\leq 2d,$ the depth of $c^*$ is at most $4d+2$ and so is the depth of $\rho^*$. 

 Let $\Omega^* \coloneqq \langle v_{d+1},y_{d+1},w_1,\ldots,w_k,x_{d+1}\rangle,$ then $\Omega^*$ is a cyclic order of the nodes of $c^*$ by the above construction, in particular $(\sigma_{\rho^*}(c^*),\Omega^*)$ is a $c^*$-society and $P^*$ is a path satisfying $V(P^*) \cap N(\rho^*) \subseteq V(\Omega_c^*)$ concluding $ii)$ of the lemma.
 
 Define $\alpha^* = \left(\langle v_{d+1},y_{d+1},w_1,\ldots,w_k,x_{d+1}\rangle,\langle v_{d+1},y_{d+1},w_1,\ldots,w_k,x_{d+1}\rangle\right)$ by setting $V_{d+1} \coloneqq F \cup Z \cup Y_{d+1} \cup X_{d+1}$. Then $\alpha^*$ is a linear decomposition of $(\sigma_{\rho^*}(c^*),\Omega^*)$ with gap $1$ by construction, and since $\left(\langle y_{d+1},w_1,\ldots,w_k,x_{d+1}\rangle,\langle Y_{d+1},W_1,\ldots,W_k,X_{d+1}\rangle\right)$ is a $d$-stuffed linear decomposition of $(G[Y_{d+1} \cup X_{d+1} \cup \bigcup_{i=1}^k W_i],\Omega')$ where $\Omega'$ is the restriction of $\Omega$ to $\{y_{d+1},w_1,\ldots,w_k,x_{d+1}\},$ $\alpha^*$ is $d$-stuffed with gap $1$ by \cref{def:linked_vortex} and this is still witnessed by $\PPP$ concluding $i)$ of the lemma.

 By \cref{claim:the_paths_are_enveloped} we derive $iii)$ of the lemma (this only depends on $\rho$ not $\rho^*$), and by $(\star)$ $C^* = P^* \cup B_{d+1}$ is a cycle grounded in $\rho$ with $V(C^*) \cap N(\rho^*) = V(\Omega_{c^*})$ $\rho$-enveloping $\PPP$ at $\BBB$ and $(\PPP,\BBB)$ $\rho$-kisses $c$ at $\CCC$ by $(\star\star)$.

 This concludes the proof.
\end{proof}

\subsection{Linking stacked vortices}

In light of the discussion above, we make precise what we mean by ``stacked vortices''; see \cref{fig:vortices_setting} for a schematic illustration.

\begin{definition}\label{def:stacked}
 Let $s_1,s_2 \geq 1$. Let $\Sigma$ be a surface, and let $\rho$ be a $\Sigma$-rendition of a graph $G$. Let $\Delta_0^1,\Delta_0^2$ be two distinct vortices in $\rho$.
 Let $\mathcal{C}_1 = \langle C_1^1, \ldots, C_{s_1}^1 \rangle$ and $\mathcal{C}_2 = \langle C_1^2, \ldots, C_{s_2}^2 \rangle$ be two sequences of concentric cycles in $\rho$ around $\Delta_0^1$ and $\Delta_0^2,$ respectively. Let $\Delta_j^i$ for $i \in [2]$ and $j \in [s_i]$ be the $C_j^i$-disc in $\rho$ containing $\Delta_0^i$ such that $\Delta_0^i \subsetneq \Delta_1^i \subsetneq\ldots \subsetneq \Delta_{s_i}^i,$ and $\CCC_i$ is tight around $\Delta_0^i$ for both $i \in [2]$. Then we say that $(\mathcal{C}_1, \Delta_0^1)$ and $(\mathcal{C}_2, \Delta_0^2)$ are \emph{stacked} with respect to each other if we either have that $\Delta_{s_1}^1$ and $\Delta_{s_2}^2$ are disjoint, or there exists an $i \in [2]$ and a $j \in [s_i]$ such that $\Delta_{s_{3-i}}^{3-i}$ is found in $\Delta_j^i - (\mathsf{bd}(\Delta_j^i) \cup \Delta_{j-1}^i),$ i.\@e.\@, in the interior of the annulus bounded by traces of $C_j^i$ and $C_{j-1}^i$.
\end{definition}

The following is the main theorem of this section.
The main idea is that, given two vortices $(\mathcal{C}_1, \Delta_0^1)$ and $(\CCC_2, \Delta_0^2)$ as in \cref{def:stacked} with enough concentric cycles~--~meaning $b(d+1)$ many where $b$ is the number of vortices~--~and after contracting the graph of the first vortex together with the $(b-1)(d+1)$ inner cycles of $\CCC_1$ into a single vertex $x$ and subsequently applying \cref{lem:linking_vortices} resulting in a new linked vortex $\tilde \Delta_0^2,$ we can guarantee that $x$ will not be part of the adhesion of the ``$d$-linked part'' of $\tilde \Delta_0^2$. This is possible since \cref{lem:linking_vortices} results in a ``tight'' set $\PPP$ of $d$ paths witnessing linkedness~--~this is made precise by the notion of $\rho$-kissing~--~which, due to the presence of $d+1$ concentric cycles around the contracted vertex $x$ are ``deflected'' from $x$ by the choice of ``tightness'' for $\PPP$. See the left hand side of \cref{fig:linking_vortices_proof} for an illustration of this, where $c_2$ may be thought of as the focused vortex with the purple region close around it~--~these are the $(d+1)$ cycles needed to link it~--~and $c_1$ is a vortex inside the purple region of $c_2,$ where the yellow region is to be thought of as the $(b-1)(d-1)$ cycles that are ``protected'' by the purple region of $(d+1)$-layers needed to deflect the paths. After linking $c_2,$ there are more than $(d+1)$ cycles around it left~--~see the right hand side \cref{fig:linking_vortices_proof} where these cycles are marked by the orange region~--~that can be used to deflect paths in a next iteration step.

More precisely, the whole vortex $\Delta_0^1$ will either be contained in one of the bags of the resulting linear decomposition of $\tilde {\Delta}_0^2,$ i.\@e.\@, it was \emph{absorbed} by the second vortex, or $\Delta_0^1$ together with the inner $(b-1)(d+1)$ cycles of $\mathcal{C}_1$ will lie outside of $\tilde \Delta_0^2$.
Note that once a vortex is absorbed, it will no longer be a cell, and thus no longer be a vortex in the resulting rendition; in particular no vortex is absorbed more than once.
If the first vortex was not absorbed, we continue linking it in a second iteration.
Note that this does not result in an exponential blow-up, as the number of cycles needed to link a vortex is only $(d+1)$ and similarly only $(d+1)$ concentric cycles are needed to guarantee that a vortex can be safely absorbed.
Therefore, we have enough cycles left to continue linking $\Delta_0^1$ possibly absorbing other vortices, and after having linked $\tilde \Delta_0^2,$ it has enough concentric cycles left around it so that it can be safely absorbed in a later step. For technical reasons we will work with $(b+1)(d+1)$ many concentric cycles, in order to easily witness that the resulting rendition remains \emph{stretched}, i.\@e.\@, vortex boundaries are disjoint, and since we will need ``disjoint'' $d$-stuffed linear decompositions, this easily guarantees that the resulting linkages witnessing this are too. 

We continue with a formal formulation and proof of the above.

\begin{theorem}\label{thm:linking_vortices}
 Let $d,b \in \N$. Let $\Sigma$ be a surface and $\rho$ a $\Sigma$-rendition of depth $d$ of some graph $G$.
 Let $\{c_1,\ldots,c_b\} \subseteq C(\rho)$ be the (possibly empty) set of all vortex cells.
 For $i\in [b]$ let $\mathcal{C}_i = \langle C_1^i,\ldots,C_{(b+1)\cdot(d+1)}^i\rangle$ be disjoint sequences of $(b+1)\cdot(d + 1)$ concentric cycles around $c_i,$ with~$\Delta_i$ being the closed $C_1^i$-disc in $\rho$ containing $c_i,$ such that $(\mathcal{C}_i,c_i)$ and $(\mathcal{C}_j,c_j)$ are stacked for all distinct $i,j \in [b]$.
 
 Then there exists a $\Sigma$-rendition $\rho'$ of $G$ together with a set of $b'$ cells $\{c_{i_1}^*,\ldots,c_{i_{b'}}^*\}$ for $1 \leq i_1 < \ldots i_{b'} \leq b$ such that for every $ j \in [b'],$ $c_{i_j}^*$ has either $\leq 3$ nodes or it is $d$-stuffed with gap at most~$1$ together with linkages $\PPP_{i_j},\BBB_{i_j}$ of order $d+1,$ such that 

 \begin{enumerate}
 \item for all $\ell \in [b]$ there exists $j \in [b']$ such that $c_\ell \subseteq c_{i_j}^*\subseteq \Delta_{i_j},$ 
 \item for every $j \in [b'],$ $\PPP_{i_j}$ witnesses that $c_{i_j}^*$ is $d$-stuffed with gap $1,$ 
 \item for all distinct $j,k \in [b'],$ $(\PPP_{i_j},\BBB_{i_j})$ $\rho$-kisses $c_{i_j}$ in $\langle C_1^{i_j},\ldots,C_{d+1}^{i_j}\rangle,$ and for all $X_j \in \{\PPP_{i_j},\BBB_{i_j}\}$ and $X_k\in \{\PPP_{i_k},\BBB_{i_k}\},$ $X_j$ and $X_k$ are disjoint,
 \end{enumerate}
 and the restrictions of $\rho$ and $\rho'$ to $\Sigma\setminus\bigcup_{j=1}^{b'} c_{i_j}^*$ agree. In particular, every vortex of $C(\rho')$ is part of $\{c_{i_1}^*,\ldots,c_{i_b}^*\}$.
\end{theorem}
\begin{proof}
If there are no vortices in $C(\rho)$ then there is nothing to show. Let $\{c_1,\ldots,c_b\} \subseteq C(\rho)$ be the set of vortices for some $b \geq 1$. Note that then every vortex comes with at least $2$ concentric cycles around it, and hence the respective $C_j^i$-discs are unambiguously defined for $j \in[(b+1)\cdot(d+1)]$ and $i \in [b]$.
Let $\mathcal{V}_0 \subseteq\{c_1,\ldots,c_b\}$ be the set of those vortices that are \emph{not} $d$-stuffed with gap at most $1$ and let $\mathcal{W}_0\subseteq \{c_1,\ldots,c_b\}$ be the set of those vortices that are $d$-stuffed with gap $1;$ clearly $\mathcal{W}_0,\mathcal{V}_0$ partition $\{c_1,\ldots,c_b\}$. The proof is by induction on $\Abs{\mathcal{V}_0}$ and without loss of generality let it be $b;$ no vortex is $d$-stuffed with gap $1$. Let~$\rho_0 \coloneqq \rho$

For $0 \leq i \leq t\leq b$ for some respective $t,$ we will inductively construct $\rho_i,\mathcal{V}_i,\mathcal{W}_i$ maintaining the following induction invariants:
\begin{description}
 \item[$\rho_i$:] We have that $\rho_i$ is a $\Sigma$-rendition of $G$ such that $\rho_i$ agrees with $\rho$ on $\Sigma \setminus (\bigcup\mathcal{V}_i\cup\bigcup \mathcal{W}_i),$ and the set of vortices of $\rho_i$ is contained in $\mathcal{V}_i \cup \mathcal{W}_i$. 
 \item[$\mathcal{V}_i$:] It holds $\mathcal{V}_i \subsetneq \mathcal{V}_{i-1}$ is a strict subset, in particular $\Abs{\mathcal{V}_i} < \Abs{\mathcal{V}_{i-1}},$ where $\mathcal{V}_i \subseteq C(\rho_i)$ is the set of those vortices that are \emph{not} $d$-stuffed with gap at most $1$ in $\rho_i$. Further for every vortex $c\in \mathcal{V}_i$ there is a set of $(b+1-i)(d+1)$ tight concentric cycles $\CCC_{\VVV}(c)$ in $\rho_i$ around $c$. See the left hand \cref{fig:linking_vortices_proof} for an illustration where $\VVV_i$ is the set of vortices highlighted in \textcolor{amaranth}{red}, and the set $\CCC_{\VVV}(c)$ is highlighted in \textcolor{Mustard}{yellow}.
 \item[$\mathcal{W}_i$:] The set $\mathcal{W}_i \subseteq C(\rho_i)$ is a set of cells of $\rho_i$ that either have at most $3$ nodes or they are $d$-stuffed vortices $c_j^* \in \WWW_i$ with gap at most $1$ witnessed by some linkage $\PPP_j$ of order $(d+1)$ together with a linkage $\BBB_j$ of order $(d+1)$ such that $(\PPP_j,\BBB_j)$ $\rho_i$-kisses $c_j^*$ at $\langle C_1^j,\ldots,C_{d+1}^j\rangle$ for some respective $j \in [b]$. Further for every $c^* \in \mathcal{W}_i$ there is some vortex $c \in C(\rho)$ such that $c \subseteq c^*$ and there is a set of $(b-i+1)(d+1)$ (not necessarily tight) concentric cycles $\CCC_{\WWW}(c)$ in $\rho_i$ around $c^*$. See the left hand \cref{fig:linking_vortices_proof} for an illustration where $\WWW_i$ is the set of vortices highlighted in \textcolor{darkspringgreen}{green}, and the set $\CCC_{\WWW}(c)$ is highlighted in \textcolor{orangepeel}{orange}.
 \item[$(\star)$:] For every vortex $c \in C(\rho)$ it holds that either $c \in \mathcal{V}_i$ or $c \subseteq c^*$ for some $c^* \in \mathcal{W}_i$ and $\mathcal{W}_i, \mathcal{V}_i$ partition the set of vortices of $\rho_i$. Furthermore, for any pair of two distinct vortices $c,c' \in C(\rho_i)$ it holds that $(\CCC(c),c)$ and $(\CCC(c'),c')$ are stacked, where $\CCC(c) \in \{\CCC_{\VVV}(c),\CCC_\WWW(c)\}$ and $\CCC(c') \in \{\CCC_{\VVV}(c'),\CCC_\WWW(c')\},$ respectively; note that for $i=b$ there are still $(d+1)$ concentric cycles left, hence the notion of stacked is well-defined.
\end{description}

Let $\CCC_1,\ldots,\CCC_b$ be as in the theorem, then $\VVV_0,\WWW_0,\rho_0$ satisfy the induction invariant; thus assume that we constructed $\VVV_{i-1},\WWW_{i-1},\rho_{i-1}$ as in the theorem for some $1 \leq i\leq b$~--~note that for $i=b$ we are done as then $\Abs{\VVV_b} = 0$~--~and let $1 \leq \Abs{\VVV_{i-1}} ,$ otherwise we are done again. Note that~$\Abs{\VVV_{i-1}} \leq b-i+1$.

Let $b^* \coloneqq \Abs{\VVV_{i-1}} + \Abs{\WWW_{i-1}},$ then $b^* \leq b$ and for simplicity let $\{c_1,\ldots,c_{b^*}\} = \VVV_{i-1} \cup \WWW_{i-1}$ with $c_1 \in \VVV_{i-1};$ we will ``stuff'' $c_1$ in order to construct $\VVV_i,\WWW_i,\rho_i$.

By assumption $(\star),$ there exist $(\CCC_1,c_1), \ldots, (\CCC_{b^*},c_{b^*})$ that are pairwise stacked, and~--~after possibly renaming the cycles for vortices in $\WWW_{i-1},$ as they may not be the respective most inner set of the tight concentric cycles following the induction invariant~--~let $\CCC_{j} = \langle C_1^j,\ldots,C_{(b-i+2)(d+1)}^{j}\rangle)$. Let $\Delta_j^i$ be the $C_j^i$-disc in $\rho$ containing $c_i$ for $i \in [b^*]$ and $j\in[(b-i+2)(d+1)]$.
 By \cref{def:stacked} of stacked, we derive that for every $j \in [b^*]$ either $\Delta_{(b-i+2)(d+1)}^j \subseteq \Delta_{(d+1)}^1$ or they are disjoint; recall that~$b-i+2 \geq 1$.

 Let $I \subseteq [b^*]\setminus\{1\}$ be maximal with $\Delta_{(b-i+2)(d+1)}^i \subseteq \Delta_{d+1}^1;$ note that for all $j \in [b^*]\setminus (I\cup \{1\})$ the following construction does not interfere with the vortices or the respective sets of concentric cycles around them guaranteed by the fact that they are stacked using the induction invariant. For every $j \in I$ let $G[\rho,C_{(b-i+1)(d+1)}^j]$ denote the extended crop of $G$ by $C_{(b-i+1)(d+1)}^j$. 
 Furthermore, let $G^*$ be obtained from $G$ by contracting for every $j \in I$ the graph $G[\rho,C_{(b-i+1)(d+1)}^j]$ to a single vertex $x^j$. Let $\rho_i^*$ be the resulting rendition obtained from $\rho_{i-1}$ by deleting all the cells contained in $\bigcup_{j \in I}\Delta_{(b-i+1)(d+1)}^j$ and identifying all the nodes of $N(\rho) \cap C_{(b-i+1)(d+1)}^j$ to $x^j$ for all $j\in I$. Then $\rho_i^*$ is a rendition of $G^*$ and since we did not create any new vortices by the above deletion and subsequent contraction, the restriction $\rho^*$ of $\rho_i^*$ to $\Delta_{d+1}^1$ is a cylindrical rendition around $c_1$ of the $\Delta_{d+1}^1$-society in $\rho_i^*$. Finally, since for all $j \in I,$ $(\CCC_j,c_j)$ and $(\CCC_1,c_1)$ were stacked we derive that $\CCC_1'\coloneqq \langle C_1^1,\ldots,C_{d+1}^1\rangle$ is still a set of tight concentric cycles around $c_1$ in $\rho_i^*$ by the induction invariant on $\VVV_{i-1},$ in particular it is a nest in $\rho^*$. And further 
 \begin{itemize}
 \item[$(\star\star)$] for every $j \in I$ the sequence $\langle C_{(b-i+1)(d+1)+1}^j,\ldots,C_{(b-i+2)(d+1)}^j\rangle$ consists of $(d+1)$ concentric cycles around $x^j$ in $\rho_i^*$\footnote{Concentric cycles around a vortex are defined as concentric cycles around the cell consisting solely of $x$.}.
 \end{itemize}
To illustrate $(\star\star)$: The \textcolor{mediumslateblue}{purple} outer regions of vortices depicted on the right hand side of \cref{fig:linking_vortices_proof} are schematic illustrations of the respective sequences $\langle C_{(b-i+1)(d+1)+1}^j,\ldots,C_{(b-i+2)(d+1)}^j\rangle,$ where the \textcolor{Mustard}{yellow} regions have been contracted to vertices $x^i$ in the above process.

 Thus, we may apply \cref{lem:linking_vortices} to $(\CCC_1',c_1),$ $\Delta_{d+1}^1$ and $\rho_i^*$ deriving the following: There exists a $\Sigma$-rendition $\rho_i'$ of $G^*$ together with a cell $c_1^* \in C(\rho_i')\setminus C(\rho_i^*)$ such that $\rho_i'$ agrees with $\rho_i^*$ on $\Sigma \setminus c_1^*$ and $c_1 \subseteq c_1^* \subseteq \Delta_{d+1}^1$ and one of two cases occur, i.\@e.\@, $\Abs{N(c^*)} \leq 3$ or not. We start with the former; let $I^*\subseteq I$ be maximal such that for $j \in I^*$ it holds $x^j \in \sigma_{\rho_i'}(c_1^*)$ (intuitively this represents the set of vortices that have been \emph{absorbed} by $c_1^*$ in the process).
 \begin{description}
 \item[$\Abs{N(c_1^*)} \leq 3$:] In this case define $\VVV_i \coloneqq \VVV_{i-1} \setminus \{ c_1, c_j \mid j \in I^*\}$ and $\WWW_i \coloneqq (\WWW_{i-1}\setminus \{c_j\mid j\in I^*\})\cup\{c_1^*\}$. Finally let $\rho_i$ be obtained from $\rho_i'$ by ``uncontracting'' every $x^j$ for $j\in I\setminus I^*$ in the obvious way such that $\rho_i$ agrees with $\rho_{i-1}$ on $\Sigma \setminus c_1^*$. This reintroduces the vortices $c_j$ together with their ``nests'' $\langle C_1^j,\ldots,C_{(b-i+1)(d+1)}^j\rangle $ of possibly tight~--~depending on whether they are in $\VVV_{i-1}$ or $\WWW_{i-1}$~--~concentric cycles around $c_j$ for every $j\in I\setminus I^*$. 

 By construction $\VVV_i,\WWW_i$ and $\rho_i$ satisfy the inductive invariants concluding the proof for this case.
 \end{description}
 
 Thus we continue with the case that $\Abs{N(c_1^*)}> 3,$ in particular $c_1^*$ is a vortex and by \cref{lem:linking_vortices} there is a $d$-stuffed linear decomposition $\alpha$ with gap at most $1$ of the $c_1^*$-society $(\sigma_{\rho_i'}(c_1^*),\Omega_{c_1^*})$. Further there exist two linkages $\PPP=\{P_1,\ldots,P_{d+1}\}$ and $\BBB=\{B_1,\ldots,B_{d+1}\}$ of order $(d+1)$ in $G^*$ satisfying the following:
 \begin{enumerate}
 
 \item $\{P_1,\ldots,P_d\}$ witnesses that $c_1^*$ is $d$-linked with gap at most $1,$ and
 \item $P_{d+1}$ satisfies $V(P_{d+1}) \cap N(\rho') \subseteq V(\Omega_{c^*}),$ in particular $\PPP$ witnesses that $c^*$ is $d$-stuffed with gap at most $1,$ and
 \item for every $i \in [d+1]$ it holds $B_i \subsetneq C_i^1$ and $\BBB$ is a base for $\PPP$
 \end{enumerate} 
 such that the cycle $C_1^* \coloneqq P_{d+1} \cup B_{d+1}$ $\rho_i^*$-envelopes $\PPP$ at $\BBB$ with $V(C_1^*) \cap N(\rho_i') = V(\Omega_{c_1^*})$ and $(\PPP, \BBB)$ $\rho_i^*$-kisses $c_1$ in $\CCC_1'$. Let $X,Y \subseteq V(G^*)$ be disjoint such that $\PPP$ and $\BBB$ are $X$-$Y$-linkages; these exist since $\BBB$ is a base for $\PPP$ in $\CCC$. 

 Let $\alpha=\left(\langle u_0,u_1,\ldots,u_m\rangle,\langle U_0,U_1,\ldots,U_m\rangle\right)$ and for $p \in [m-1]$ let $Z_p \coloneqq (U_p \cap U_{p+1})\setminus V(\Omega_{c_1^*}),$ then $\Abs{Z_p} = d$.

 \begin{beautifulclaim}\label{claim:absorb_vortices}
 It holds $(V(\Omega_{c_1^*}) \cup \bigcup_{p=1}^{d}Z_p) \cap \{x^j\mid j \in I\} = \emptyset$.
 \end{beautifulclaim}
 \begin{claimproof}
 Note that $\{x^j \mid j\in I\} \cap V(\bigcup \BBB) = \emptyset$ since $\{x^j \mid j\in I\} \cap V(\bigcup \CCC_1') = \emptyset$ and $\BBB$ is a base for $\PPP$ at $\CCC_1'$. Thus if $V(\Omega_{c_1^*}) \cap \{x^j \mid j \in I\} \neq \emptyset$ implies $V(P_{d+1}) \cap \{x^j \mid j \in I\} \neq \emptyset$.

 Since $(\PPP,\BBB)$ $\rho_i^*$-kisses $c_1$ in $\CCC_1',$ we derive by \cref{def:kissing} that for every linkage $\PPP'$ of order $(d+1)$ such that $\BBB$ is a base of $\PPP'$ in $\CCC_1'$ and such that $C'$ $\rho_i^*$-envelopes $\PPP$ at $\BBB,$ it holds that $\Delta_{C_1^*}$ is contained in the $C'$-disc. 

 Let $j \in I$ and assume towards a contradiction that the claim is wrong. We have two cases to consider.
\begin{description}
 \item[Case $x^j \in V(\Omega^*)$:] As discussed above this implies $x^j \in V(P_{d+1})$. Let $H\coloneqq G^*[\rho_i^*,C_1^*]$ be the extended crop of $G^*$ by $C_1^*$ in $\rho_i^*,$ and let $H' = H - x^j$. Note that if there is an $X$-$Y$-linkage $\PPP'=\{P_1',\ldots,P'_{d+1}\}$ in $H'$ after deleting the internal vertices of $\BBB,$ then $(\PPP',\BBB)$ $\rho_i^*$-kisses $c_1$ in $\CCC_1'$ with perimeter $(P'_{d+1},B_{d+1});$ this follows analogously to $(\star\star)_{\ref{lem:linking_vortices}}$ in the proof of \cref{lem:linking_vortices} using \cref{stolen_uno}.
 By definition of $\rho_i^*$-kissing this is impossible since for the cycle $C' = P'_{d+1}\cup B_{d+1} \neq C_1^*$ the $\rho_i^*$-aligned $C'$-disc $\Delta_{C'}$ is contained in the $C_1^*$-disc as a contradiction to the fact that $(\PPP,\BBB)$ $\rho_i^*$-kisses $c_1^*$ in $\CCC_1'$.
 
 However, Menger's Theorem \cref{thm:menger} guarantees the existence of such a linkage $\PPP'$ as a contradiction. To see this, assume the contrary. Then there exists a separation $(A,B)$ in $H'$ of order at most $d$ such that $X \subseteq A$ and $Y \subseteq B;$ let $S\coloneqq V(A)\cap V(B)$. Since~$P_{d+1}$ is a grounded path in $\rho_i^*$ and since by $(\star\star)$ we derive that $\bigcup_{\ell=(b-i+1)(d+1)+1}^{(b-i+2)(d+1)}C_\ell^j$ is grounded in $\rho_i^*$ satisfying $ x^j\in \Delta_{(b-i+1)(d+1)+1}^j \subsetneq \ldots \subsetneq \Delta_{(b-i+2)(d+1)}^j,$ we derive that $\Abs{V(P_{d+1}) \cap V(C_{\ell}^j)} \geq 2$ for all $\ell \in [(b-i+1)(d+1)+1,(b-i+2)(d+1)]$.
 
 Now $\PPP\setminus \{P_{d+1}\}$ is still a linkage of size $d$ in $H,$ in particular $S \subseteq \bigcup_{i=1}^d V(P_i)$. By the pigeonhole principle, there is a $t \in [(b-i+1)(d+1)+1,(b-i+2)(d+1)]$ such that $V(C_t^j) \cap S = \emptyset$. Let $P_{d+1}^{\mathsf{in}},P_{d+1}^{\mathsf{out}} \subsetneq P_{d+1}$ be disjoint paths of shortest possible length such that $P_{d+1}^\mathsf{in}$ has one endpoint in $Y$ and the other in $V(C_{t}^j),$ and $P_{d+1}^\mathsf{out}$ has one endpoint in $V(C_{t}^j)$ and the other in $X$. Clearly $x^j \notin V(P_{d+1}^\mathsf{in}) \cup V(P_{d+1}^\mathsf{out})$. But then $P_{d+1}^\mathsf{in} \cup P_{d+1}^\mathsf{out}\cup C_t^j$ is a connected subgraph in $H'-S$ with one vertex in $X$ and one in $Y$ as a contradiction to the assumption that $(A,B)$ was a separation of order $d$ separating $X$ from $Y$.

 \item[Case $x^j\in Z_p$:] Let $p \in [m]$. By $i)$ every path in $\{P_1,\ldots,P_d\}$ contains exactly one vertex of $Z_p$ for every $p \in [m]$. If $x^j \in Z_p,$ then $x^j 
 \in V(P_j)$ for some $j \in [d]$. Again by Mengers Theorem \cref{thm:menger}, analogously to the case $x^j\in V(\Omega_{c_1^*}),$ we derive the existence of an $X$-$Y$-linkage $\PPP'$ of order $(d+1)$ in $H' \coloneqq H-x^j$ . Then this is a contradiction to the fact that $(Z_p \cup \{u_p\})\setminus\{x^j\}$ witnesses a separation $(A,B)$ of order $d$ in $H'$. 
\end{description}
This concludes the proof of the claim.
 \end{claimproof}

 By \cref{claim:absorb_vortices} we derive that for every $j \in I$ if $x^j \in U_p$ for $p \in [m]$ then for all $p' \in [m] \setminus \{p\}$ we have $x^j \notin U_{p'}$. Note that a priori it may happen that $x^j \in U_0 \cap U_1$ or $x^j \in U_0 \cap U_m,$ however, this is of no importance for the statement the theorem, as we are working with linkedness up to exactly that gap. Recall the definition of $I^*$. Then we construct $\VVV_i,\WWW_i$ and $\rho_i$ equivalently to the case $\Abs{N_{\rho_i'}(c_1^*)} \leq 3;$ we repeat it for ease of readability.

\begin{description}
 \item[$\Abs{N_{\rho_i'}(c_1^*)} \geq 4$:]D $\VVV_i \coloneqq \VVV_{i-1} \setminus \{ c_1, c_j \mid j \in I^*\}$ and $\WWW_i \coloneqq (\WWW_{i-1}\setminus \{c_j\mid j\in I^*\})\cup\{c_1^*\}$. Finally let $\rho_i$ be obtained from $\rho_i'$ by ``uncontracting'' every $x^j$ for $j\in I\setminus I^*$ in the obvious way such that $\rho_i$ agrees with $\rho_{i-1}$ on $\Sigma \setminus c_1^*$. This reintroduces the vortices $c_j$ together with a sequence $\langle C_1^j,\ldots,C_{(b-i+1)(d+1)}^j\rangle$ of $(b-i+1)(d+1)$ tight concentric cycles around $c_j$ for every $j\in I\setminus I^*$. 
 \end{description}

 We have the following.
 \begin{beautifulclaim}
 $\VVV_i,\WWW_i$ and $\rho_i$ satisfy the induction invariants.
 \end{beautifulclaim}
 \begin{claimproof}
 By construction $\rho_i$ is a $\Sigma$-rendition of $G$ such that $\rho_i$ agrees with $\rho_{i-1}$ on $\Sigma \setminus c_1^*,$ and hence since $\rho_{i-1}$ satisfies the induction invariant by assumption, we derive that $\rho_i$ agrees with $\rho$ on $\Sigma \setminus (\bigcup\mathcal{V}_i\cup\bigcup \mathcal{W}_i)$. Similarly by construction the set of vortices of $\rho_i$ is contained in~$\mathcal{V}_i \cup \mathcal{W}_i$: To see this note that for all the cells $c_j$ for $j \in I^*$ the graph $\sigma_{\rho_{i-1}}(c_j)$ is contained in $\sigma_{\rho_i}(c_1^*),$ in particular $c_j \subseteq c_1^*\subseteq \Sigma$ whence it is not a vortex in $\rho_i$. And further $\VVV_i \cup \WWW_i$ partition $(\VVV_{i-1}\cup \WWW_{i-1})\setminus\{c_j \mid j \in I^*\},$ hence this follows by assumption that $\VVV_{i-1},\WWW_{i-1},\rho_{i-1}$ satisfy the induction invariants. This proves the invariant on $\rho_i$.

 By definition $\VVV_{i}$ is a strict subset of $\VVV_{i-1}$ and $\VVV_i$ consists exactly of the vortices in $\rho_i$ that are not yet $d$-stuffed with gap $1,$ for those vortices are part of $\WWW_i$ and both sets are disjoint and partition the vortices of $\rho_i$ as derived above. But every cell $c \in \WWW_i$ is either in $\WWW_{i-1}$ or is equal to $c_1^*;$ each of them are $d$-stuffed with gap at most $1$. For $c_1^*$ this follows by construction, and for all others by the fact that $\WWW_{i-1}$ satisfies the induction invariant together with the fact that~$\rho_i$ agrees with $\rho_{i-1}$ on $\Sigma \setminus c_1^*$. By the same arguments we derive that $\mathcal{W}_i \subseteq C(\rho_i)$ is a set of cells of $\rho_i$ that either have at most $3$ nodes or they are $d$-stuffed vortices with gap at most $1$. By construction of $c_1^*,$ $(\PPP,\BBB)$ $\rho_i^*$-kisses $c_1$ in $\CCC_1'$ and in particular it $\rho_i$-kisses $c_1$ in $\CCC_1',$ and for the other vortices in $\WWW_i$ this remains true by the induction invariant and the construction of $\rho_i$ using \cref{claim:absorb_vortices} (the linkages $\PPP,\BBB$ do not intersect the linkages respective witnessing that the remaining vortices in $\WWW_i$ are $d$-stuffed with gap $1$).

 For every $c \in \VVV_i$ we derive that there is a set of $(b-i+1)(d+1)$ tight concentric cycles $\CCC_\VVV(c)$ in $\rho_i$ around $c$: This follows from the fact that $\VVV_i = \VVV_{i-1} \setminus \{c_1,c_j \mid j \in I^*\}$ and the fact that for every $j \in I \setminus I^*$ the sequence $\langle C_1^j,\ldots,C_{(b-i+1)(d+1)}^j\rangle$ is still a sequence of tight concentric cycles around $c_j$ by the previous construction; see the construction in the case $\Abs{N_{\rho_i'}(c_1^*)} \geq 4$. Further, for all $c \in \VVV_{i-1}\setminus I$ this follows from the fact that $\rho_i$ agrees with $\rho_{i-1}$ on $\Sigma^*\setminus c_1^*$ and our initial choice of $I$ using that the original families of vortices were stacked. This proves the induction invariant for $\VVV_i$.

 Using analogous arguments, for every $c^* \in \mathcal{W}_i$ there is some vortex $c \in C(\rho)$ such that $c \subseteq c^*$ and there is a set of $(b-i+1)(d+1)$ (not necessarily tight) concentric cycles $\CCC_{\WWW}(c)$ in $\rho_i$ around $c,$ disjoint from the respective linkages witnessing that it is $d$-stuffed with gap at most $1,$ that is disjoint from the respective linkages of $(\PPP_c,\BBB_c)$ that $\rho_i$-kisses $c$. Note that this is true for $c_1^*$ by the fact that $c_1^* \subseteq \Delta_{(d+1)}^1\subsetneq \Delta_{(d+2)}^1 \subsetneq \ldots \subsetneq \Delta_{(b-i+2)(d+1)}^1,$ hence $\langle C_{d+2}^1,\ldots,C_{(b-i+2)(d+1)}^1\rangle$ is a sequence of $(b-i+1)(d+1)$ concentric cycles around $c_1^*$ satisfying the above. This proves the induction invariant for $\WWW_i$.

 Finally $(\star)$ holds true by the above and the fact that the sequences of concentric cycles are subsequences of $\CCC_1,\ldots,\CCC_b$ which were stacked in the first place.
 \end{claimproof}

Since $(\VVV_i,\WWW_i,\rho_i)$ maintain the inductive invariants, let $t \leq b $ be minimal such that $\Abs{\VVV_t}= 0;$ then $\rho_t,\WWW_t$ satisfy the theorem where the vortices in $\WWW_t$ are $d$-stuffed with gap $1$ and satisfy $i),ii),iii)$ of the theorem by construction, concluding the proof.
\end{proof}
Note that a careful analysis of the above proof shows that depth of the resulting rendition $\rho'$ is at most $4d+2$~--~using \cref{lem:linking_vortices}, but as we do not need explicit bounds on the depth~--~we only need $d$-linkedness up to some bounded gap~--~we omitted a thorough proof.

In light of \cref{thm:linking_vortices} we define the following.

\begin{definition}\label{def:properly_stuffed_and_linked_rendition}
 Let $\Sigma$ be a surface and let $\rho$ be a $\Sigma$-rendition of some graph $G$. Let $d,\ell \in \N$. We call $\rho$ \emph{(properly) $d$-linked, respectively $d$-stuffed with gap $\ell$}, if every vortex of $\rho$ is (properly) $d$-linked, respectively $d$-stuffed with gap $\ell,$ and for every pair of distinct vortices $c_1,c_2 \in C(\rho),$ the paths $\PPP_1,\PPP_2$ witnessing that they are (properly) $d$-linked, respectively $d$-stuffed with gap $\ell,$ are disjoint.
\end{definition}

We derive the following from \cref{thm:linking_vortices}.

\begin{corollary}\label{cor:linked_vortices}
 Let $d,b \in \N$. Let $\Sigma$ be a surface and $\rho$ a $\Sigma$-rendition of depth $d$ of some graph $G$.
 Let $\{c_1,\ldots,c_b\} \subseteq C(\rho)$ be the (possibly empty) set of all vortex cells.
 For $i\in [b]$ let $\mathcal{C}_i = \{C_1^i,\ldots,C_{(b+1)\cdot(d+1)}^i\}$ be disjoint sets of $(b+1)\cdot(d + 1)$ concentric cycles around $c_i,$ with $\Delta_i$ being the closed $C_1^i$-disc in $\rho$ containing $c_i$ for each $i \in [b],$ and $(\mathcal{C}_i,c_i)$ and $(\mathcal{C}_j,c_j)$ are stacked for all distinct $i,j \in [b]$.
 
 Then there exists a $\Sigma$-rendition $\rho'$ of $G$ that is $d$-stuffed with gap at most $1$ such that for every vortex $c \in C(\rho)$ there is a vortex $c^* \in C(\rho')$ with $c \subseteq c^*$.
\end{corollary}

\section{The apex-free case}\label{sec:apexfree}

Throughout this section we work with renditions in surfaces, in particular we are \emph{apex-free}. The main goal of this section is to prove that, given an instance $(G,T)$ with a vital $T$-linkage~$\LLL$ and a blank $\Sigma$-rendition of $(G,T),$ then there exists another instance $(G',T')$ with a vital $T'$-linkage $\LLL'$ such that $\Abs{T'} \in \mathbf{poly}(\Abs{T}) \cdot 2^{\mathbf{poly}(b)}$ where $b = \mathsf{bidim}(G,T)$ and such that $(G,T)$ can be embedded in a surface $\Sigma'$ with all terminals drawn on boundary components of $\Sigma',$ where the surfaces $\hat\Sigma'$ and $\hat\Sigma$ obtained by ``capping their holes with discs'' are homeomorphic. We will formulate and prove the relevant results in a more general setting as it is more natural.

We need specific renditions throughout this section; we call a rendition $\rho$ \emph{stretched} if for any pair of distinct vortices $c,c' \in C(\rho)$ it holds $N(c) \cap N(c') = \emptyset$. Recall the \cref{def:concentric_cycles} of concentric cycles in $\rho,$ and that, given a sequence $\langle C_1,\ldots,C_s\rangle$ of concentric cycles in $\rho,$ we let $\Delta_{C_i}$ denote the respective $C_i$-disc in $\rho$ satisfying $\Delta_{C_1}\subsetneq \ldots \subsetneq \Delta_{C_s}$.

\begin{definition}\label{def:locus}
 Let $\Sigma$ be a surface. Let $m\geq 1$. Let $(G,T)$ be an annotated graph. Let $\rho$ be a blank and stretched $\Sigma$-rendition of $(G,T)$ and let $C^* \subseteq C(\rho)$ be its set of vortices. Let $M=\langle C_1,\ldots,C_m\rangle$ be a sequence of concentric cycles in $\rho$ such that $\Delta_{C_m}$ is vortex-free\footnote{If $m=1$ and there are two choices for $\Delta_{c_m}$ -- in particular $C^* = \emptyset$ -- fix one of them.}. Then we call $(\rho,M)$ an \emph{$m$-locus (of $(G,T)$ in $\Sigma$)}. We refer to~$C_i$ as a \emph{cycle of $M$} for every $i \in [m]$.

We call $(G,T),\LLL,(\rho,M)$ \emph{exhausted} if~$$G = \bigcup \LLL \cup \bigcup_{i=1}^m C_i \cup \bigcup_{c \in C^*}\sigma(c)\mbox{~and~}{E(\bigcup \LLL) \cap E(\bigcup_{i=1}^m C_i) = \emptyset}.$$
\end{definition}
Note that due to the definition of $m$-locus, no vertex of $T$ is contained in the interior of $\Delta_{C_m}$. Furthermore, recall that by definition of concentric cycles, each cycle of $M$ is grounded in $\rho,$ in particular it is edge-disjoint from the graph of every vortex.

\begin{definition}
 Let $\Sigma$ be a surface, $G$ a graph and $\rho$ a $\Sigma$-rendition of $G$. Let $c\in C(\rho)$ be a cell with $N(c)= \{a_c,b_c\}$ such that $\sigma(c)$ consists of the single edge $a_cb_c$. Then we call $c$ \emph{dead}\footnote{Such cells are sometimes referred to as ``trivial'' in the literature.}.
 
 We call $\rho$ a \emph{$\Sigma$-skeleton of $G$} if every non-vortex cell is dead.
\end{definition}

Given a $\Sigma$-rendition, we will classify different homotopy types of paths whose traces have both endpoints in vortex-boundaries: These will be loops, if both endpoints lie on the society of a common vortex, and \emph{links} otherwise. 

\begin{definition}\label{def:links}
 Let $\Sigma$ be a surface and $\rho$ a $\Sigma$-rendition of some graph $G$. Let $c,c' \in C(\rho)$ be vortices. Let $P \subseteq G$ be grounded such that both its endpoints are in $N(c)\cup N(c')$. If $c \neq c'$ we call $P$ a \emph{$c$-$c'$ link} or simply \emph{link}. If $c = c',$ we call $P$ a \emph{loop (at $c$)}.
 
 Given a loop $P$ at $c,$ let $\gamma\subsetneq \Sigma$ be the trace of $P$. If there exists an arc-wise connected component $\tau \subsetneq \mathsf{cl}(c)\setminus \{u,v\}$ such that $\gamma\cup \tau$ bounds a disc $\Delta_P(c) \subseteq \Sigma$ that contains no vortex of $\rho$ and is disjoint from the closure of all vortices distinct from $c,$ we call $P$ \emph{simple}.

 Let $C^* \subseteq C(\rho)$ be the set of vortices. We call two links (or loops) \emph{$\rho$-homotopic} if their traces are homotopic in $\Sigma \setminus \bigcup_{c \in C^*}\mathsf{int}(c),$ respectively.
\end{definition}
Note that given a simple loop $P$ the disc $\Delta_P(c)$ is unique.

The following is the main theorem of this section; note that if a rendition is stretched then no path can be simultaneously a link and a loop.
\begin{theorem}\label{thm:no_apices_vital_reduction}
 Let $t,\xi,g, k \in \N$ and $d\geq 1$. Let $\Sigma$ be a surface of genus $g$. Then there exist functions $\mathsf{terminals}\colon \N^4\to \N,$ $\mathsf{cycles}\colon \N^3\to \N$ $\mathsf{fat}\colon\N^2\to \N$ satisfying the following.
 
 Let $(G,T)$ be a $k$-annotated graph such that 
 \begin{itemize}
 \item there exists a $\mathsf{cycles}(k,d,t)$-locus $(\rho,M)$ of $(G,T)$ in $\Sigma$ such that $\rho$ is $d$-stuffed with gap at most $d$ admitting at most $\xi$ vortices, and
 \item there exists a vital $T$-linkage in $G$.
 \end{itemize}
 Then there exists a $\mathsf{terminals}(k,\xi,d,g)$-annotated graph $(G',T'),$ and a $t$-locus $(\rho',M')$ of $(G',T')$ in $\Sigma$ together with a vital $T'$-linkage $\LLL'$ in $G'$ such that 
 \begin{enumerate}
 \item $\rho'$ is a $\Sigma$-skeleton with at most $\xi$ vortices, and $(G',T'),\LLL',(\rho',M')$ is exhausted,
 \item all vortices $c \in C(\rho')$ satisfy $\sigma(c) = (N(c),\emptyset),$
 \item let $L\in \LLL',$ then $L$ is either a non-simple loop or a link; in particular it is internally disjoint from $N_{\rho'}(c)$ for every vortex cell $c \in C(\rho'),$
 \item let $P \in \LLL'$ be a non-simple loop or a link in $G',$ and let $\LLL'(P) \subseteq \LLL'$ be the maximum size linkage of paths $\rho$-homotopic to $P$. Then $\Abs{\LLL'(P)}\leq \mathsf{fat}(k,d)$.
 \end{enumerate}
 Moreover, $\mathsf{terminals}(k,\xi,d,g)\in {\mathbf{poly}(k+\xi+g)}\cdot 2^{\mathbf{poly}(d)},$ $\mathsf{cycles}(k,d,t) \in {\mathbf{poly}(k+t)}\cdot2^{\mathbf{poly}(d)},$ and $\mathsf{fat}(k,d)\in {\mathbf{poly}(k)}\cdot 2^{\mathbf{poly}(d)},$ 
\end{theorem}

\subsection{Exhausting the instance}
In order to prove \zcref{thm:no_apices_vital_reduction} we need several technical steps; the first one~--~marking the main result of this subsection~--~allows us to reduce an instance $(G,T),\LLL,(\rho,M)$ where $(\rho,M)$ is some $(m+1)$-locus of $(G,T)$ in a surface $\Sigma,$ to that of an exhausted instance $(G',T),\LLL',(\rho',M')$ on the same set of annotations where $(\rho',M')$ is an $m$-locus of $(G',T')$ in $\Sigma;$ recall \cref{def:locus}.
 In particular, the instance we are left with comes with a $\Sigma$-skeleton, and thus it is ``quasi-embedded'', meaning that it is essentially embedded up-to the vortex cells. The main techniques employed to prove this result expand on existing techniques \cite{RobertsonS2009Graph,CavallaroKK2024EdgeDisjoint,AdlerKKLST2017Irrelevant}. Unfortunately, we cannot easily transcribe the proposed results from the literature as we need to take care of further technicalities, and thus we give a full proof dedicated to our setting and needs for completion. The main idea of the proof being that no two paths in $\LLL$ can use non-node vertices of the graph of a common non-vortex cell, for the graphs of non-vortex cells can be separated from the rest by separations of order~$\leq 3$. This can be exploited in order to ``throw away'' cells with~$\leq 3$ nodes via simply replacing them with single edges respecting the linkage $\LLL$. Continuing this argumentation inductively by a careful analysis of the graph's structure inside the cells with respect to the cells' societies, we can reduce to $\Sigma$-skeletons. In a next instance, one then ``throws away'' vertices and contracts edges until left with an exhausted instance maintaining a large locus; recall the definitions from \zcref{sec:linkagepreparation}.

 There is one more technicality that we need to add to the above sketch, which we briefly discussed in \cref{sec:linkedness}: When given a $d$-stuffed vortex with gap $1,$ one of the $(d+1)$ paths witnessing that the vortex is $d$-stuffed, may partially ``live'' outside of the vortex, due to the fact that we work with renditions. We need to carefully ``tuck'' that path into the graph of the vortex when exhausting the instance, as we want the adhesion of our linear decomposition to match the size of the linkage inside the graph of the vortex ``linking'' the separators; compare with \cref{thm:Lemma2.5}. To this extent, recall \cref{def:linked_vortex,def:properly_stuffed_and_linked_rendition} of properly $d$-linked renditions with gap at most $\ell$.

\begin{figure}
 \centering
 \begin{subfigure}[t]{0.32\textwidth}
 \centering
 \includegraphics[width=0.75\linewidth]{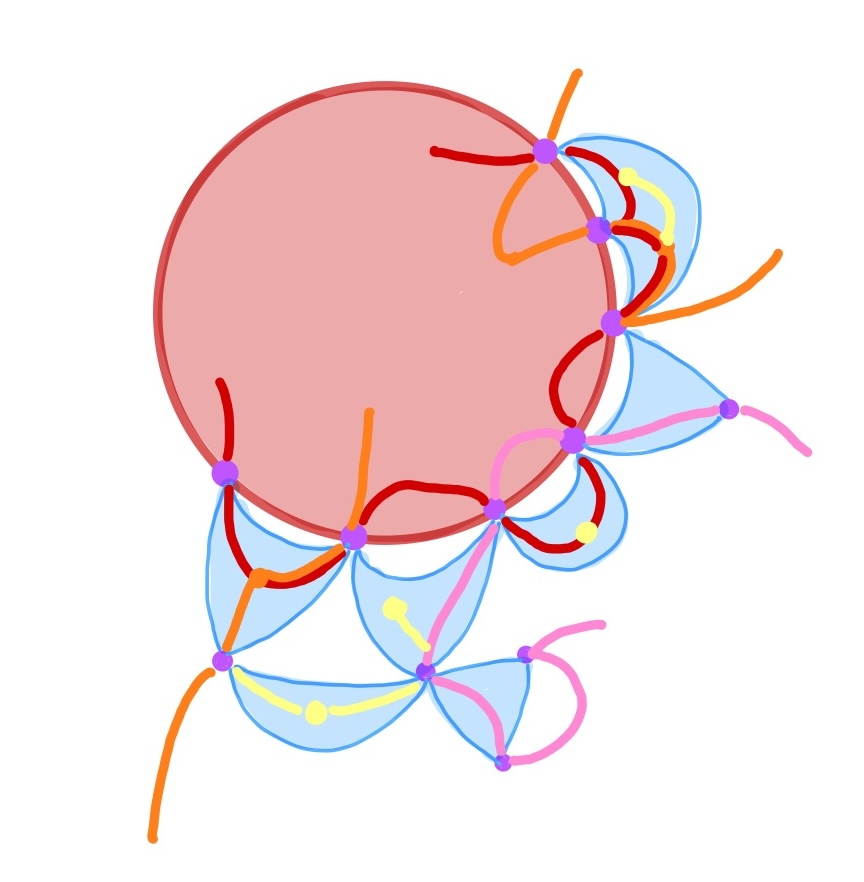}
 \caption{$P_i^*$ and forbidden configurations.}
 \end{subfigure}%
 ~ 
 \begin{subfigure}[t]{0.32\textwidth}
 \centering
 \includegraphics[width=0.7\linewidth]{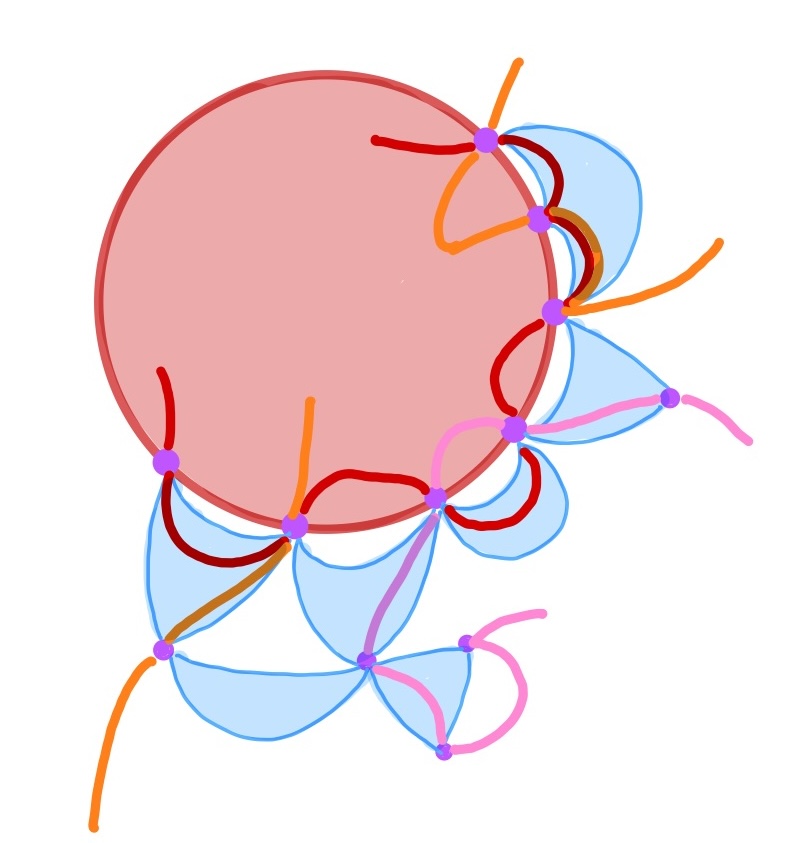}
 \caption{$P_i^*$ after contraction.}
 \end{subfigure}
 \begin{subfigure}[t]{0.32\textwidth}
 \centering
 \includegraphics[width=0.7\linewidth]{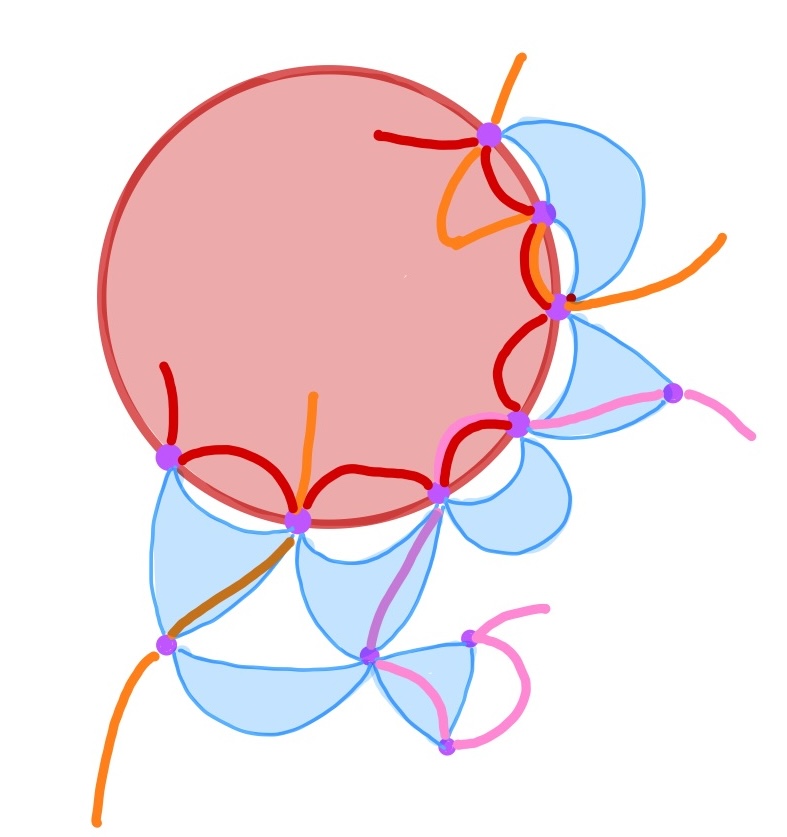}
 \caption{Tucking in $P_i^*$.}
 \end{subfigure}%
 \caption{A schematic illustration of the proof of \cref{thm:rendition_to_skeleton}. The red regions depict a vortex $c_i$. The figures depict a path $P_i^*$ in red together with paths of a vital linkage in orange and pink in some rendition, where blue regions denote cells. Yellow elements in cells denote ``forbidden'' configurations violating vitality of $\LLL$. Figure $(a)$ illustartes the initial setting. Figure $(b)$ illustrates the setting after contracting the subpaths of $P_i^*$ in cells to edges. Figure $(c)$ illustrates how to ``tuck'' $P_i^*$ into the vortex.}
 \label{fig:lemma_skeleton}
\end{figure}
\begin{theorem}\label{thm:rendition_to_skeleton}
 Let $\Sigma$ be a surface. Let $m,d,\ell \geq 1$ and let $(G,T)$ be an annotated graph admitting a vital $T$-linkage $\LLL$. Let $(\rho,M)$ be an $(m+1)$-locus of $(G,T)$ such that $\rho$ is $d$-stuffed with gap at most $\ell$. Let $C^*$ be the set of vortices of $\rho$.

 Then there exists a graph $G'$ that is a minor of $G,$ such that $(G',T)$ is an annotated graph admitting a vital $T$-linkage $\LLL'$ with $\tau(\LLL') =\tau(\LLL)$ and an $m$-locus $(\rho',M')$ where $M'$ is tight, and $\rho'$ is a $\Sigma$-skeleton of $G'$. Furthermore, $(G',T),\LLL',(\rho',M')$ is exhausted and $\rho'$ is properly $d$-linked with gap at most $\ell$ such that $C^*$ is the set of vortices of $\rho'$ and for every $c \in C^*$ it holds $N_{\rho'}(c) = N_{\rho}(c)$ and $\bigcup_{c \in C^*}\sigma_{\rho'}(c) \cap \bigcup M' = \emptyset$.
\end{theorem}
\begin{proof}
 Throughout this proof, whenever we say rendition we mean a $\Sigma$-rendition, as $\Sigma$ is clear from context. Let $(G,T),$ $\LLL,$ $(\rho,M)$ be as in the theorem. Let $M=\langle C_1,\ldots,C_{m+1}\rangle$ and let $C^*=\{c_1,\ldots,c_b\}$ be the set of vortices of $\rho$ (possibly empty) for some $b \in \N$ with respective vortex societies $(\sigma_{\rho}(c_i),\Omega_{c_i})$ for $i \in [b]$. 

 For every $i\in [b],$ let $\PPP_i$ be a linkage of order $(d+1)$ witnessing that $c_i$ is $d$-stuffed with gap at most $\ell$ where $\PPP_i,\PPP_j$ are disjoint for distinct $i,j \in [b]$. Let $P_i^* \in \PPP_i$ be the unique path with its nodes in $V(\Omega_{c_i})$ and $\bigcup(\PPP_i\setminus \{P_i^*\}) \subseteq \sigma_{\rho}(c_i)$ by \cref{def:linked_vortex}; call such a path $P_i^*$ a \emph{$c_i$-hook}. Since $M$ is an $(m+1)$-locus, for every $i \in [m]$ it holds $V(C_i) \cap \bigcup_{j=1}^b\sigma_{\rho}(c_j) = \emptyset;$ let $M' \coloneqq \langle C_1,\ldots,C_m\rangle$. Furthermore for every $i \in [b]$ and $c_i$-hook $P_i^*$ let $\PPP(P_i^*)$ be the set of maximal subpaths of $P_i^*$ with both ends on $V(\Omega_{c_i})$ such that they are completely contained in the graph of some single non-vortex cell. By construction, for every $i \in [b],$ $\PPP(P_i^*)$ is a set of internally disjoint paths. Note that if $\PPP_i^* = \emptyset,$ then $P_i^* \subseteq \sigma(c_i),$ and hence $c_i$ is properly $d$-linked with gap $\ell$ by \cref{def:linked_vortex}. 
 
 As a first step, let $G'$ be obtained from $G$ by deleting all edges $e \in E(G)$ that are not part of $E(\bigcup M'),$ not part of $E(\bigcup \LLL),$ not part of $E(\bigcup_{i=1}^bP_i^*)$ and are not in $\bigcup_{i=1}^b\sigma_{\rho}(c_i)$. Then $G' \subseteq G,$ $V(G')=V(G),$ $\LLL$ is still a vital $T$-linkage in $G',$ and $(\rho,M')$ is an $m$-locus of $(G',T)$ in $\Sigma$. For notational simplicity, let $G = G'$ and $M = M';$ where the outercycle of $M$ is by the above construction disjoint from $\bigcup_{i=1}^b\sigma_{\rho}(c_i)$. We have the following.
 \begin{itemize}
 \item[$(\star)$] $G = \bigcup \LLL \cup \bigcup M \cup \bigcup_{i=1}^b P_i^* \cup \bigcup_{i=1}^b\sigma(c_i),$
 \end{itemize}
and $\rho$ is a rendition of $G,$ $d$-stuffed with gap at most $\ell$ witnessed by the same linkages $\PPP_1,\ldots,\PPP_b$.
 
We call a cell $c\in C(\rho)$ \emph{dying} if it is not a vortex, it is not dead and either $E(\sigma(c)) = \emptyset$ or $V(\sigma(c)) =N_\rho(c)$. Note that technically $E(\sigma(c)) = \emptyset$ is allowed for vortices, hence the additional assumption.
 
\begin{beautifulclaim}\label{thm:rendition_to_skeleton_claim_vital_cell_interaction}
 Let $c \in C(\sigma)$ be a non-vortex cell.
 Then, either
 \begin{enumerate}
 \item $c$ is dying, or
 \item there is a unique path $L(c) \in \LLL$ such that $E(L(c)) \cap E(\sigma(c)) \neq \emptyset$ and~$\Abs{V(L(c)) \cap N_\rho(c)} \geq 2$. 
 \end{enumerate}
 \end{beautifulclaim}
 \begin{claimproof}
 By assumption $c$ is not a vortex cell. If $c$ is dead $ii)$ is true and if it is dying $i)$ is true. Henceforth, assume that $E(\sigma(c)) \neq \emptyset$ as well as $V(\sigma(c)) \neq N_\rho(c),$ whence there exist a vertex $u \in V(\sigma(c))\setminus N_\rho(c)$ and an edge $e \in E(\sigma(c))$ incident to $u$. Note that such an edge must exist since $\LLL$ is vital implying that $u \in V(L)$ for some $L \in \LLL$ and since $\rho$ is a blank rendition and $c$ is not a vortex cell, $u$ must be adjacent to some edge. Further, since $u$ is not a node of the cell and $\rho$ is a rendition, $E(\sigma(c))$ must contain at least one edge adjacent to $u$. Since $\rho$ is blank and $\Abs{N(c)} \leq 3$~--~for it is not a vortex~--~there is at most one path $L \in \LLL$ such that $E(L) \cap E(\sigma(c)) \neq \emptyset$. To see this note that both terminals of $L$ lie in $V(G) \setminus V(\sigma(c))$ since $\rho$ is blank, in particular $N(c)$ separates $T$ from $\sigma(c)$ with $E(\sigma(c)) \cap E(\sigma(c')) = \emptyset$ for all cells $c' \neq c$. Thus, since $\Abs{N(c)} \leq 3$ and paths in $\LLL$ are pairwise vertex-disjoint, at most one path in $\LLL$ may use an edge of $\sigma(c)$. Similarly, since $N_\rho(c)$ separates $u$ from $T,$ at least one path of $\LLL$ satisfies $E(L) \cap E(\sigma(c)) \neq \emptyset$. Combining both observations, such a path exists and is unique, and since both endpoints of the path are in $T,$ it must have at least two common vertices with $N_{\rho}(c)$. This concludes the proof; see \cref{fig:lemma_skeleton} for a schematic illustration of the general setting.
 \end{claimproof}

Throughout the proof we will write $L(c)$ for the unique path satisfying \zcref{thm:rendition_to_skeleton_claim_vital_cell_interaction} whenever applicable. 

 \begin{beautifulclaim}\label{claim:all_P^*_are_edges}
 There is a minor $G' \subseteq G$ such that $(G',T)$ is an annotated graph admitting a vital $T$-linkage $\LLL'$ such that $\tau(\LLL') = \tau(\LLL)$. Further there is an $m$-locus $(\rho',M)$ of $(G,T)$ such that $\rho'$ is $d$-stuffed with gap at most $\ell$ and vortices $\{c_1,\ldots,c_b\}$ (but it may not be equal to $\rho$ on these). Let $\PPP_1',\ldots,\PPP_b'$ be respective linkages witnessing that $\rho'$ is $d$-stuffed with gap $\ell$ and with $c_i$-hooks $P_i'$ respectively. Then every path in $\PPP(P_i')$ consists of a single edge.
 \end{beautifulclaim}
 \begin{claimproof}
 The proof is by induction on $\{P_1^*,\ldots,P_b^*\}$. Note that $P_1^*,\ldots,P_b^*$ are disjoint, and further no two paths in $\bigcup_{i=1}^b \PPP(P_i^*)$ are contained in the graph of the same cell by definition of $\PPP(P_i^*)$. Thus, we may construct a minor $G'$ of $G$ inductively by contracting the paths in $\bigcup_{i=1}^b \PPP(P_i^*)$ as follows. Let $P \in \PPP(P_i^*)$ for some $i \in [b]$ and let $c \in C(\rho)$ be the respective cell with $P \subseteq \sigma(c)$. If $N(c) \subseteq V(P),$ then, since $V(P) \cap N(\rho) \subseteq V(\Omega_{c_i}),$ we can add the whole graph $\sigma_{\rho}(c)$ to $\sigma_{\rho}(c_i),$ and the resulting rendition is still $d$-stuffed with gap at most $\ell$ as witnessed by the same linkages. Henceforth, assume that $\Abs{N(P) \cap N(c)} = 2$ and $\Abs{N(c)}=3;$ let $N(c) = \{u_1,u_2,u_3\}$ and let the endpoints of $P$ be $u_1,u_2,$ say. If $c$ is dying, then $P$ is an edge, for $\Abs{N(P) \cap N(c)} = 2$ and we are done; hence assume that $c$ is not dying. 
 
 By \cref{thm:rendition_to_skeleton_claim_vital_cell_interaction} we derive that either $E(P) \cap E(L) = \emptyset$ for every $L \in \LLL$ or $E(P) \cap E(L(c)) \neq \emptyset$ and it is disjoint from all other paths in $\LLL$. If the former holds true $P$ must again be an edge, for otherwise there is $v \in V(P) \setminus N(c)$ which is not part of any $L \in \LLL$ contradicting vitality of $\LLL$. We derive that $P \cap L(c) \neq \emptyset;$ note that their intersection must be a subpath of $L(c),$ i.\@e.\@, a single component, for else $\LLL$ was not vital as we could reroute $L(c)$ using $P$. Contract every edge in $e \in E(P) \cap E(L(c))$. Then this results in a graph $G'$ for which the resulting linkage $\LLL'$ has the same pattern as $\LLL$ since $T \subseteq V(G')$ as the graphs of vortices have not been altered. We claim that $\LLL'$ is still a vital linkage: Otherwise, let $\LLL^*$ be a linkage in $G'$ with the same pattern as $\LLL'$ but distinct from it. Uncontracting $P$ now yields a linkage $\tilde \LLL^*$ in $G$ with the same pattern as $\LLL$ but distinct from it as a contradiction.
 
 Finally, we claim that after the above contraction $P$ consists of a single edge. If not, then there is $u \in V(P) \setminus \{u_1,u_2\}$ with $u \in \sigma(c)$. Note that $u \neq u_3$ for $\Abs{V(P) \cap N(c)} = 2$. By \cref{thm:rendition_to_skeleton_claim_vital_cell_interaction} there is a single path of $\LLL$ that has an edge in $\sigma(c)$. Thus $u \in V(L(c))$ still. Let $u_i,u_j \in \{u_1,u_2,u_3\}$ be the two endpoints of $L(c)$. Then $\{u_i,u_j\} \cap \{u_1,u_2\} \neq \emptyset,$ say $u_i = u_1$ for simplicity. Hence, there is a $\{u\}$-$\{u_1\}$-path $L' \subseteq L(c)$ as well as a $\{u\}$-$\{u_1\}$-path $P' \subseteq P$. Since $E(P) \cap E(L(c)) = \emptyset$ we can reroute $L(c)$ along $P'$ as a contradiction to vitality. See \cref{fig:lemma_skeleton} $(a)$ and $(b)$ for a schematic illustration depicting most cases of the steps discussed above.

 Repeating the above for all $i \in [b]$ and all paths in $\PPP(P_i^*)$ we conclude the proof, noting that we never alter the graphs of cells containing edges of $\bigcup M$.
 \end{claimproof}

For simplicity, assume that $(G,T)$ $\LLL,$ $(\rho,M)$ and $\PPP_1,\ldots,\PPP_b$ as well as $P_1^*,\ldots,P_b^*$ are such that they satisfy \cref{claim:all_P^*_are_edges}. We are left to ``tuck'' the respective edges as given by $\PPP(P_i^*)$ ``into the vortices'' to make the rendition properly $d$-linked with gap at most $\ell$. To do this, note that for all $i \in [b]$ and all $P \in \PPP(P_i^*),$ we have $P = uv$ for consecutive $u,v \in V(\Omega_{c_i})$. 
Define $\rho'$ from $\rho$ by adding $P$ to the graph $\sigma_{\rho}(c_i)$ and removing it from the respective graph of the cell of $\rho$ it was part of. This results in an $m$-locus $(\rho',M)$ of $(G,T)$ such that $\rho'$ is properly $d$-linked with gap at most $d$. For notational simplicity, let $\rho = \rho'$. We summarise the above as follows:

 \begin{itemize}
 \item[$(\star\star)$] $G = \bigcup \LLL \cup \bigcup M \cup \bigcup_{i=1}^b\sigma(c_i)$ and $\rho$ is properly $d$-linked with gap at most $\ell$.
 \end{itemize}
 
 With this in mind, we continue to reduce $\rho$ to a $\Sigma$-skeleton. 
 
 The following is an easy, but slightly technical claim; it essentially states that whenever we can refine $\rho$ by splitting one cell into two cells without creating new vortices, then we do so. Note that this does not affect $G$ or $\LLL$ in any way but is a pure refining of the rendition at hand.

\begin{beautifulclaim}\label{thm:rendition_to_skeleton_claim_vital_3cells_are_2connected}
 There is an $m$-locus $(\rho',M)$ of $(G,T)$ satisfying $(\star\star)$ such that for every non-vortex cell $c \in C(\rho)$ with $N_\rho(c) =\{u_1,u_2,u_3\}$ for distinct vertices and every separation $(A,B)$ in $\sigma(c)$ with $A \cap \{u_1,u_2,u_3\} \neq \emptyset$ and $B \cap \{u_1,u_2,u_3\} \neq \emptyset,$ it holds $\Abs{A \cap B} \geq 2$.
\end{beautifulclaim}
\begin{claimproof}
 The proof is by inductively refining $\rho$ as follows. Let $c \in C(\sigma)$ be a cell such that $N_\rho(c) =\{u_1,u_2,u_3\}$ for distinct vertices and let $(A,B)$ be an order $\leq 1$ separation in $\sigma(c)$ with $A \cap \{u_1,u_2,u_3\} \neq \emptyset$ and $B \cap \{u_1,u_2,u_3\} \neq \emptyset$: call such a cell \emph{bad}. We prove how to refine the rendition by reducing the number of bad cells inductively.
 
 First note that if the separation is of order $0,$ then $A, B$ partition $\{u_1,u_2,u_3\}$. Thus, we may define $\rho'$ from $\rho$ by removing the cell $c$ and adding two new disjoint cells $c_A,c_B \subsetneq c \subsetneq \Sigma$ by letting $\sigma_{\rho'}(c_A) = \sigma_{\rho}(c)[A]$ and $\sigma_{\rho'}(c_B) = \sigma_{\rho}(c)[B]$ where $c_A$ has as nodes $A \cap \{u_1,u_2,u_3\}$ and $c_B$ has as nodes $B\cap \{u_1,u_2,u_3\};$ both have $\leq 2$ nodes. Then $(\rho',M)$ is an $m$-locus of $G$ in $\Sigma$ and $\PPP_1,\ldots,\PPP_b$ still witness that it is $d$-stuffed with gap at most $\ell;$ the details are easy and left to the reader. In particular the number of bad cells in $\rho'$ is strictly less than the ones in $\rho$ by construction.

 Thus let the separation be of order $1$ and let $x \in A\cap B$ be the unique vertex (possibly $x \in \{u_1,u_2,u_3\}$). Then construct $\rho'$ from $\rho$ by deleting $c$ and introducing new disjoint cells $c_A,c_B \subsetneq c \subsetneq \Sigma$ with nodes $A\cap\{u_1,u_2,u_3\} \cup \{x\}$ and $B\cap\{u_1,u_2,u_3\} \cup \{x\}$ respectively, by letting $\sigma_{\rho'}(c_A) = \sigma_\rho(c)[A]$ and $\sigma_{\rho'}(c_B) = \sigma_{\rho}(c)[B],$ and drawing $x$ as a node inside $c$ away from its boundary. One easily verifies that $(\rho',M)$ is an $m$-locus of $(G,T)$ in $\Sigma$ (after possibly fixing a tiebreaker if $\Abs{N_{\rho'}(c_A)} = 2$ or $\Abs{N_{\rho'}(c_B)} = 2$). To see this note that if some cycle $C$ of $M$ has an edge in $E(\sigma(c)),$ then the trace of $C$ in $\rho'$ is still equal to its old trace except for its subsets of the boundary of $c_A$ and $c_B$. Since $c_A,c_B \subseteq c$ by construction, one can easily choose a tiebreaker (if necessary) such that the trace of $C$ in $\rho'$ is as desired.
 
 Note however, that both cells $c_A$ and $c_B$ may have exactly three nodes each, and the nodes may still be separable by separators of order $1$. In particular the number of bad cells may have increased. But one may now inductively continue with the cells $c_A$ and $c_B$ by refining the rendition $\rho'$ further at possible $1$ separations. Since in each step we further decompose the graph $\sigma(c),$ this procedure eventually stops, and we effectively reduced the number of bad cells.

 By repeating the above for every bad cell of $\rho$ we conclude the proof (note that $(\star\star)$ is still satisfied by construction).
\end{claimproof}

Henceforth, we assume that $(G,T)$ and $(\rho,M)$ are as given by \zcref{thm:rendition_to_skeleton_claim_vital_3cells_are_2connected}; $\LLL$ is still a vital linkage in $G$. 

\begin{beautifulclaim}\label{thm:rendition_to_skeleton_claim_vital_3cell_interaction}
 Let $c \in C(\rho)$ be a non-vortex cell with $N_\rho(c) =\{u_1,u_2,u_3\}$ for distinct vertices. If $c$ is not dying, then $\Abs{V(L(c)) \cap N_\rho(c)} = 2$.
\end{beautifulclaim}
\begin{claimproof}
 Assume that $\Abs{V(L(c)) \cap N_\rho(c))} = 3,$ i.\@e.\@, the path visits all three vertices. Fix an orientation of the path via a start and endpoint of the path and let $\langle v_0,\ldots,v_m\rangle$ be a linear order on its vertices with respect to that orientation. Assume that $\langle u_1,u_2,u_3\rangle$ are visited in that order by $L(c)$ and let $L'(c)\subseteq L(c)$ be the respective subpath starting in $u_1$ and ending in $u_3$. By \zcref{thm:rendition_to_skeleton_claim_vital_3cells_are_2connected} we derive that there is a $\{u_1\}$-$\{u_3\}$-path $P$ in $\sigma(c) - u_2$. In particular $P \neq L(c)'$. But since no path of $\LLL^* \coloneqq \LLL \setminus \{L(c)\}$ uses an edge or vertex in $V(\sigma(c))$~--~for $L(c)$ visits all three vertices $\{u_1,u_2,u_3\}$ by assumption~--~the path $P$ is vertex-disjoint from every path in $\LLL^*$. Let $L$ be obtained from $L(c)$ by replacing $L'(c)$ with $P,$ then $L$ is a path and it has the same endpoints as $L(c)$. By construction then $\LLL^* \cup \{L\}$ is a linkage in $G$ with pattern $\tau(\LLL)$ contradicting the fact that $\LLL$ was vital.
\end{claimproof}

Combining the above claims, one easily kills dying cells as follows.
\begin{beautifulclaim}\label{thm:rendition_to_skeleton_claim_killing_dying_cells}
 There exists a graph $G' \subseteq G,$ such that $(G',T)$ is an annotated graph admitting a vital $T$-linkage $\LLL'$ with $\tau(\LLL') =\tau(\LLL)$ and an $m$-locus $(\rho',M)$ where $\rho'$ has no dying cells, and it agrees with $\rho$ on all cells of $\rho$ that are not dying, in particular the instance satisfies $(\star\star)$.
\end{beautifulclaim}
\begin{claimproof}
 Let $c \in C(\rho)$ be a dying cell. If $E(\sigma(c)) = \emptyset$
 then simply delete $c$ from the rendition; this does not have any effect on $M,$ i.\@e.\@, for every cycle $C$ of $M,$ $C\subseteq G'$ remains grounded in the resulting rendition, and no cell containing an edge of $E(C)$ is altered. Henceforth assume that $V(\sigma(c)) =N_\rho(c)$ and $E(\sigma(c)) \neq \emptyset$. If $\Abs{N_\rho(c)} = 2,$ then $c$ is dead as a contradiction to it being dying; thus let $N_\rho(c) = \{u_1,u_2,u_3\}$ for distinct vertices. By \zcref{thm:rendition_to_skeleton_claim_vital_3cells_are_2connected} we derive that $\sigma(c)$ is isomorphic to $K_3$. 
 We define $\rho'$ from $\rho$ by deleting $c$ and adding three new cells $c_1,c_2,c_3 \notin C(\sigma)$ and defining $\sigma_{\rho'}(c_j)$ to consist of the edge $\{u_i,u_k\}$ for $\{i,j,k\}=\{1,2,3\}$. One easily verifies that $\rho'$ is a blank rendition of $G$ still satisfying the assumptions of the claim, in particular, every cycle of $M$ remains grounded in the new rendition. Note that this may affect at most one cycle of $M,$ say $C,$ since the cycles are pairwise disjoint. Furthermore, the cycles remain concentric, since the trace of $C$ in the new rendition remains effectively equal when choosing the new cells carefully and setting the tie-breaker for the trace of the new cells accordingly. Finally, $c_1,c_2$ and $c_3$ are dead by definition.

 Repeating this argument for every dying cell concludes the proof.
\end{claimproof}

Next, we kill cells away from cycles in $M$.
\begin{beautifulclaim}\label{thm:rendition_to_skeleton_claim_outside_wall_dead}
 There is an $m$-locus $(\rho',M)$ of $(G,T)$ such that $\rho'$ agrees with $\rho$ on all cells that contain an edge of a cycle in $M$. Further, all other non-vortex cells of $\rho'$ are dead, and the instance satisfies $(\star\star)$.
 \end{beautifulclaim}
 \begin{claimproof}
 We proceed by induction, showing how to ``kill'' cells away from $M$. Let $c \in C(\rho)$ be a non-vortex cell that is not dead and such that it does not contain an edge of a cycle in $M$. By \zcref{thm:rendition_to_skeleton_claim_killing_dying_cells} we may assume that it is not dying. 
 
 Thus, by $(\star\star)$ some path of $\LLL$ uses an edge of $\sigma(c)$. Since $c$ is not dying, \zcref{thm:rendition_to_skeleton_claim_vital_cell_interaction} implies that there is a unique path $L(c) \in \LLL$ using an edge in $\sigma(c)$. Combining \zcref{thm:rendition_to_skeleton_claim_vital_3cells_are_2connected} and \zcref{thm:rendition_to_skeleton_claim_vital_3cell_interaction} we derive that $\Abs{V(L(c)) \cap N_\rho(c)} = 2$. Let $L'(c)$ be the maximum length subpath of $L(c)$ contained in $\sigma(c),$ and let $u$ and $v$ be its endpoints. Then $u,v \in N_\rho(c)$. Since no cycle of $ M$ admits an edge in $\sigma(c)$ by assumption, we derive by $(\star\star)$ that $\sigma(c) = L(c)$. Thus, replace the cell $c$ by introducing dead cells, one for each edge of $L(c)$.

 By induction we can proceed as above until all non-vortex cells not containing an edge of a cycle of $M$ are dead, never producing new non-dead cells, i.\@e.\@, the induction stops after $r \in \N$ steps where $r$ is the number of non-dead cells in $C(\rho)$. Clearly the instance still satisfies $(\star\star)$ for we did not alter the rendition on vortices. This concludes the proof of the claim. 
 \end{claimproof}
 Combining \zcref{thm:rendition_to_skeleton_claim_killing_dying_cells} and \zcref{thm:rendition_to_skeleton_claim_outside_wall_dead} we may assume that every non-vortex cell of $\rho$ whose graph has no edge in common with a cycle in $M$ is dead~--~by switching to the respective subgraph and rendition~--~and there are no dying cells left.
 Let $M=\langle C_1,\ldots,C_m\rangle$ and fix $\CCC=\{C_1,\ldots,C_m\}$.
 \begin{beautifulclaim}\label{thm:rendition_to_skeleton_claim_inside_wall_dead}
 Let $c \in C(\rho)$ be a non-dead cell such that $E(\sigma(c)) \cap E(\bigcup \CCC) \neq \emptyset$. There exists a graph $G'$ that is a minor of $G,$ such that $(G',T)$ is a $k$-annotated graph having a vital $T$-linkage $\LLL'$ with $\tau(\LLL') =\tau(\LLL)$ and an $m$-locus $(\rho',M)$ where $\rho'$ agrees with $\rho$ on all cells $c' \neq c,$ and all other cells of $\rho'$ are dead; in particular the instance satisfies $(\star\star)$.
 \end{beautifulclaim}
 \begin{claimproof}
 Let $c$ be as in the claim and let $e \in E(\sigma) \cap E(\bigcup \CCC)$ and let $C \in \CCC$ be the respective cycle with $e \in E(C)$. There are two cases to consider: $\Abs{N_\rho(c)} \in \{2,3\}$.
 \begin{description}
 \item[Case 1; $|N_{\rho}(c)| = 2$:] Let $N_\rho(c) = \{u_1,u_2\}$ for distinct vertices. 
 Every edge of $E(\bigcup \CCC)$ is part of a grounded cycle. Thus, since $\CCC$ is a set of disjoint cycles $\sigma(c)$ contains edges of exactly one cycle in $\CCC,$ hence $C,$ and in particular it must contain a subpath of length at least $2$ of $C$ since $e \in E(C) \cap E(\sigma(c));$ let $P_C \subseteq C$ be maximal such that $E(P_c) \subseteq E(\sigma(c))$. Since $C$ is a grounded cycle we derive that $P_C$ must be a path with endpoints $u_1,u_2$.
 
 Further, since $c$ is not dying, \zcref{thm:rendition_to_skeleton_claim_vital_cell_interaction} implies that $L(c)$ exists. Let $L'(c)$ be the maximum subpath of $L(c)$ contained in $\sigma(c)$ with endpoints $u_1,u_2,$ then $L'(c)$ by vitality of $\LLL$ we derive that $L'(c) = P_C,$ else we could reroute $L(c)$ by replacing $L'(c)$ with $P_C$ contradicting that $\LLL$ was the unique linkage of its pattern. Furthermore, since $\LLL$ is vital we derive that $V(\sigma(c)) \subseteq V(L'(c))$. But then $E(\sigma(c)) = E(L'(c)),$ for any other edge could be used to reroute $L'(c)$ inside $\sigma(c)$ again refuting vitality of $\LLL$. Now we can refine the rendition $\rho$ by deleting the cell $c$ and introducing for each edge $e=xy \in E(\sigma(c)$ a new cell $c_e \subseteq c$ such that $\sigma(c_e)$ consists exactly of the edge $e$ with nodes $x,y$ in the obvious way: clearly $C$ is grounded in $\rho'$. Since each cell $c_e$ is contained in $c,$ no matter the tiebreaker for each cell, $M$ remains a sequence of concentric cycles in $\rho'$.

 \item[Case 2; $|N_{\rho}(c)| = 3$:] Let $N_\rho(c) = \{u_1,u_2,u_3\}$ for distinct vertices. Combining \zcref{thm:rendition_to_skeleton_claim_vital_3cells_are_2connected} and \zcref{thm:rendition_to_skeleton_claim_vital_3cell_interaction} we derive that $L(c)$ is the unique path with an edge in $E(\sigma(c))$ and $\Abs{V(L(c))\cap N_\rho(c)} = 2;$ without loss of generality assume that $u_1,u_2 \in V(L(c))$ and let $L'(c)$ be the maximum subpath of $L(c)$ contained in $\sigma(c)$ with endpoints $u_1,u_2$. 
 
 Since $M$ is a sequence of concentric cycles in $\rho,$ at most one cycle of $M,$ say $C,$ has an edge in $E(\sigma(c))$ and $C$ intersects $\sigma(c)$ in a supath $P_C \subsetneq C$. By $(\star)$ we derive that $\sigma(c) = L'(c) \cup P_C$. Further every vertex $u \in V(\sigma(c))\setminus\{u_3\}$ satisfies $u \in L'(c)$ since $\LLL$ is vital. By \zcref{thm:rendition_to_skeleton_claim_vital_3cells_are_2connected} we derive the existence of three (not necessarily disjoint) paths $P_1,P_2,P_3\subseteq \sigma(c)$ of length at least $2$ such that $P_i$ has $u_j,u_k$ as endpoints for $\{i,j,k\} = \{1,2,3\}$. Since $u_3 \notin V(L'(c))$ we may set $P_3 = L'(c)$ and we know that $P_j \not\subseteq L'(c)$ for $j=1,2$ since both contain $u_3$. Let $P_1',P_2'\subseteq P_C$ be maximal such that $E(P_i') \cap E(L'(c)) = \emptyset$ for $i=1,2;$ these paths exist since $\sigma(c) = L'(c) \cup P_C$. Now both $P_1'$ and $P_2'$ have length at least $2,$ whence they must intersect $L'(c)$. 

 If there exists an edge $e \in E(\sigma(c)) \setminus (E(L'(c))\cup E(P_1')\cup E(P_2')),$ then $e \in E(P_C)$ and $e$ has both its endpoints on $V(L(c)),$ since $P_c$ has two incidences with $u_3,$ and the respective edges are in $E(P_1')$ and $E(P_2')$ respectively. But then we can reroute $L'(c)$ along $e,$ and hence $L(c),$ keeping the resulting path edge-disjoint from the remaining paths in $\LLL$. This is a contradiction to the unique pattern of $\LLL$. Henceforth, we may assume that no such edge $e$ exists, whence $\sigma(c) = P_1' \cup P_2' \cup L'(c)$. Analogously we derive that no edge of $P_1'$ or $P_2'$ has both its incidences with vertices of $L'(c)$ for else we could reroute.
 
 Thus, again applying \zcref{thm:rendition_to_skeleton_claim_vital_3cells_are_2connected} we derive that $P_1'$ must end in $u_1,$ say, and $P_2'$ must end in $u_2,$ say (or vice-versa but then rename the paths accordingly) and $\{P_1',P_2',L'(c)\}$ is an internally disjoint linkage. In particular it is a subdivision of $K_3$ and admits a planar embedding. We may thus refine $\rho$ accordingly by deleting $c$ and introducing new cells, one for each edge in $E(P_1')\cup E(P_2') \cup E(L'(c))$ keeping $M$ a set of concentric circles in the resulting rendition $\rho'$ using analogous argumentation as in the previous proofs.
 \end{description} 
 As we did not alter any graphs of vortices, we derive that $(\star\star)$ still holds. This concludes the proof of the claim.
 \end{claimproof}
 
 Recall that according to our previous assumptions, all non-dead cells $c \in C(\rho)$ satisfy $E(M)\cap E(\sigma(c)) \neq \emptyset$. Thus, repeatedly applying \zcref{thm:rendition_to_skeleton_claim_inside_wall_dead} for every non-dead cell that is left, and using the fact that being a minor is a transitive relation we derive the existence of a $k$-annotated graph $(G',T)$ admitting a vital $T$-linkage $\LLL'$ with $\tau(\LLL') = \tau(\LLL)$ and an $m$-locus $(\rho',M)$ of $(G',T)$ in $\Sigma$ such that $\rho'$ is a $\Sigma$-skeleton and $(G',T),\LLL',(\rho',M)$ satisfies $(\star\star)$. 
 
 If $M=\langle C_1,\ldots,C_m\rangle$ is not tight, then \zcref{prop:makecyclestight} implies the existence of a set of $m$-concentric cycles $M^*=\langle C_1^*,\ldots,C_m^*\rangle$ such that $M^*$ is still a set of concentric cycles in $\rho$ where the trace of each $C_i^*$ lies in $\Delta_{C_m},$ in particular the cycles are grounded, their discs are vortex-free and $E(\bigcup_{i=1}^m C_i^*) \cap \bigcup_{c \in C^*} \sigma(c) = \emptyset$. In particular $(\rho',M^*)$ is an $m$-locus of $(G',T)$ in $\Sigma$ by \zcref{def:locus}. Now $(G',T),\LLL',(\rho',M^*)$ satisfies $(\star\star)$ and $M^*$ is tight.
 
 \begin{beautifulclaim}\label{thm:rendition_to_skeleton_claim_last}
 Let $(G,T),\LLL,(\rho,M)$ satisfy $(\star\star)$ and such that $M$ is tight. Then there exists a minor $G'$ of $G$ and an $m$-locus $(\rho',M')$ of $(G',T)$ in $\Sigma$ together with a vital $T$-linkage $\LLL'$ with $\tau(\LLL') = \tau(\LLL),$ such that the set $C^*$ of vortex cells of $\rho$ is the set of vortex cells of $\rho',$ and $\rho'$ agrees with $\rho$ on $C^*$. Further $\rho'$ is a $\Sigma$-skeleton, properly $d$-linked rendition with gap at most $\ell,$ and $(G',T),\LLL',(\rho',M')$ is exhausted.
 \end{beautifulclaim}
 \begin{claimproof}
 The proof is by iteratively contracting edges $e \in E(\bigcup \LLL) \cap E(\bigcup_{i=1}^mC_i);$ let $e$ be such an edge and $e \in E(C_i)$ for some $i \in [m],$ contract the edge $e$ to a single vertex and adapt the rendition accordingly in the obvious way since each edge corresponds to a single cell for $\rho$ is a $\Sigma$-skeleton. Let $C_i'$ be the respective cycle obtained from $C_i$ by contracting $e$. Then $\rho'$ is still a $\Sigma$-skeleton, $M'=\langle C_1,\ldots,C_{i-1},C_i',C_{i+1},\ldots,C_m\rangle$ is still a set of concentric cycles in $\rho',$ the resulting graph is by construction a minor of $G$ and $\LLL'$ is still a vital linkage with the same pattern. Finally $M'$ is tight, for otherwise $M$ was not---simply reroute the cycles of $M$ accordingly by ``decontracting'' the edge if necessary---as a contradiction to the above assumption.
 
 Since the resulting instance satisfies $(\star\star),$ by repeating the above iteratively reducing the amount of edges in $ E(\bigcup \LLL) \cap E(\bigcup_{i=1}^mC_i)$ with each iteration concludes the proof of the claim.
 \end{claimproof}
 Applying \zcref{thm:rendition_to_skeleton_claim_last} to $(G',T),\LLL',(\rho',M^*)$ we conclude the proof. Note that $C^*$ is still the set of vortices of $\rho'$ and $N_{\rho'}(c) = N_{\rho}(c)$ for all $c \in C^*$ by construction.
\end{proof}

\subsection{Taming the instance}

\begin{figure}[ht]
 \centering
 \begin{tikzpicture}

 \pgfdeclarelayer{background}
		\pgfdeclarelayer{foreground}
			
		\pgfsetlayers{background,main,foreground}

 \begin{pgfonlayer}{background}
 \pgftext{\includegraphics[width=11cm]{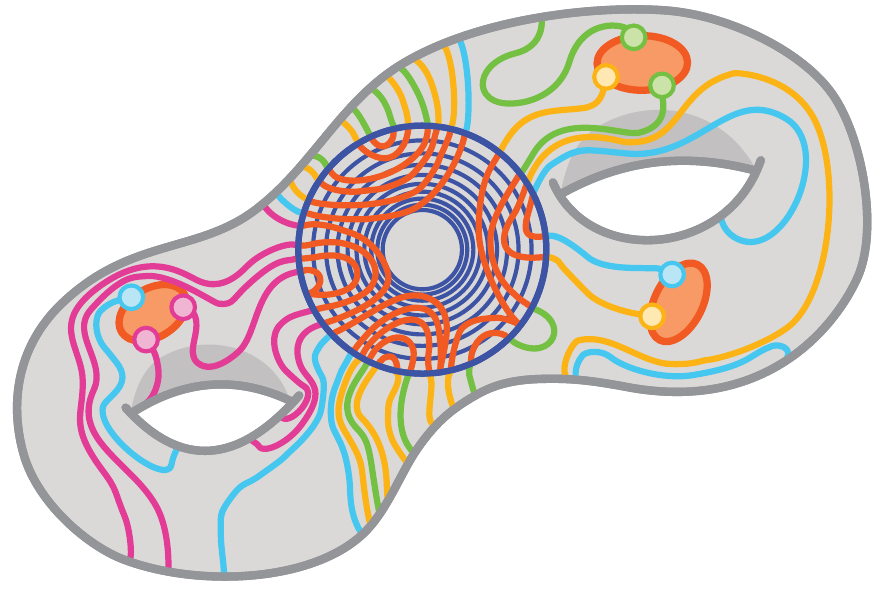}} at (C.center);
 \end{pgfonlayer}{background}
			
 \begin{pgfonlayer}{main}
 \node (C) [v:ghost] {};
 
 \end{pgfonlayer}{main}
 
 \begin{pgfonlayer}{foreground}
 \end{pgfonlayer}{foreground}

 \end{tikzpicture}
 \caption{A schematic representation of the interaction of segments of a vital linkage with a large set of concentric cycles in a disc. The four segments are highlighted in different colours, and the vortices are marked by red regions. The restriction of the paths to the cycles is highlighted by red paths exhibiting a tamed behaviour. The four segments however exhibit ``untamed'' behaviour, where a single segment produces many ``stacked'' paths in the disc.}
 \label{fig:in_vs_out_bucket}
\end{figure}

In this subsection, we gather various technical lemmas that leverage \zcref{lem:drywell} on dry wells, in order to tame the behaviour of the vital linkage $\LLL$ inside the $C_m$-disc $\Delta_{C_m}$ of the outer cycle $C_m$ of $M,$ given that $(G,T),\LLL,(\rho,M)$ is exhausted. Let the \emph{segments} of $\LLL$ be the maximal subpaths with ends in vortex boundaries that are otherwise disjoint from the graphs in vortices. We will adapt several ideas of \cite{CavallaroKK2024EdgeDisjoint} to our setting, in order to lift the derived structure of the paths ``inside'' $\Delta_{C_m}$ to the segments of $\LLL,$ which requires a non-trivial amount of work: Note that while the paths inside of $\Delta_{C_m}$ may exhibit tamed behaviour, this is a priori not guaranteed for the segments of $\LLL;$ see \cref{fig:in_vs_out_bucket}.

 To this extent, we define the following, which is essentially the restriction of the linkage $\LLL$ to a vortex-free $\rho$-aligned disc.

\begin{definition}\label{def:linkage_restriction_to_a_disc}
 Let $\Sigma$ be a surface, and let $\rho$ be a $\Sigma$-rendition of some annotated graph $G$. Let $\Delta\subseteq \Sigma$ be a vortex-free $\rho$-aligned disc. Let $\LLL$ be a linkage in $G$ such that every endpoint of a path in $\LLL$ is contained in the graph of some cell $c$ not contained in $\Delta$. Let $(G_\Delta,\Omega_\Delta)$ be the $\Delta$-society in $\rho$. We define $\LLL(\Delta)$ as the set of paths obtained from $\LLL$ as follows. For every $L\in \LLL$ let $P\subseteq L$ be a maximal subpath such that $P$ has its endpoints in $V(\Omega_\Delta),$ the trace $\gamma_P$ of $P$ is contained in $\Delta,$ and no interior point of $P$ is part of $V(\Omega_\Delta),$ then $P \in \LLL(\Delta)$. We call $\LLL(\Delta)$ the \emph{restriction of $\LLL$ to $\Delta$}.
\end{definition}
Note that the above is \emph{not} the same as $\LLL \cap G_\Delta$ where $G_\Delta$ denotes the crop of $G$ by $\Delta,$ since $\LLL \cap G_\Delta$ may contain paths with internal vertices on the boundary of $\Delta$.

One easily verifies that $\LLL(\Delta)$ is a set of pairwise internally disjoint paths; see the \textcolor{amaranth}{red} paths in \cref{fig:in_vs_out_bucket} for an illustration. We extend the definition of patterns of linkages to sets of pairwise internally disjoint paths $\PPP$ in the obvious way, i.\@e.\@, $\tau(\PPP)\coloneqq\{\{s,t\}\mid P\in \PPP \text{ has endpoints } s \text{ and } t\}$.
Note that by definition, every path $P\in \LLL(\Delta)$ has all its vertices in graphs of cells of the rendition $\rho'$ obtained by restricting $\rho$ to $\Delta$. The following is an easy observation following from the definition of $\LLL(\Delta)$ (using the notation of \zcref{def:linkage_restriction_to_a_disc}).

\begin{observation}\label{obs:sharing_ends_in_restricted_linkage}
 Let $P_1,P_2\in \LLL(\Delta)$ be distinct such that $V(P_1)\cap V(P_2)\neq \emptyset,$ then they intersect in a single vertex $x$ with $x \in V(\Omega_\Delta)$ and $P_1,P_2 \subsetneq L$ for some common $L \in \LLL$. Further, $x$ is an endpoint of both, say that $P_1$ ends in $x$ and $P_2$ starts in $x,$ and the concatenation $P_1\circ P_2$ is a subpath of $P \in \LLL$.
\end{observation}
\begin{proof}
 This is clear by construction: There exist $L_1, L_2 \in \LLL$ (possibly equal) such that $P_1 \subseteq L_1$ and $P_2 \subseteq L_2$. If $L_1 = L_2$ we derive that $\Abs{V(P_1) \cap V(P_2)} \leq 1$ and if $V(P_1) \cap V(P_2) \neq \emptyset$ then their concatenation is clearly a subpath of $L_1$. By definition of $\LLL(\Delta),$ if they intersect they must do so in a vertex of $V(\Omega_\Delta)$. If $L_1 \neq L_2$ then $V(P_1) \cap V(P_2) = \emptyset$ since $\LLL$ is a linkage. This concludes the proof.
\end{proof}

By definition of wells and exhausted instances, one can easily extract the former from the latter as follows. Given an exhausted instance $(G,T),\LLL,(\rho,M)$ with outer cycle $C_m$ of $M$ and let $\Delta_{C_m}$ be the respective $C_m$-disc. 
Recall \cref{def:extended-crop} and define $G(M,\LLL) \coloneqq G[\rho,C_m]$ as the extended crop of $G$ by $C_m$ and let $\Omega_{C_m}$ be a cyclic orientation of $V(C_m)$ respecting the cycle, i.e., two consecutive vertices with respect to $\Omega_{c_m}$ share an edge in $C_m$. Then $(G(M,\LLL),\Omega_{C_m})$ is a society and it admits a $\Delta_{C_m}$-rendition $\rho_{\Delta_m}$. By construction $\bigcup M \subseteq G(M,\LLL)$ and $\bigcup \LLL(\Delta_{C_m}) \subseteq G(M,\LLL)$ where $\LLL(\Delta_{C_m})$ is a set of internally disjoint paths with its endpoints in $V(\Omega_{C_m});$ the resulting graph and $\LLL(\Delta_{C_m})$ are essentially given as in \cref{fig:in_vs_out_bucket}, where the \textbf{thick} cycle highlights the outer cycle $C_m$ and the \textcolor{amaranth}{red} paths illustrate $\LLL(\Delta_{C_m})$.

\begin{observation}\label{obs:well_from_locus}
 Let $\Sigma$ be a surface. Let $m \geq 1$. Let $(G,T),\LLL,(\rho,M)$ be exhausted where $(G,T)$ is an annotated graph, $\LLL$ is some vital $T$-linkage in $G$ and $(\rho,M)$ is an $m$-locus of $(G,T)$ in $\Sigma$. Then $\WWW \coloneqq (G(M,\LLL),\Omega_{\Delta_{C_m}},\rho_{\Delta_m},M,\LLL(\Delta_{C_m}))$ is an $m$-well. And if $M$ is tight, then $\WWW$ is tight.
\end{observation}

Note that the reason we add back $C_m$ in the construction of $G(\LLL,M)$ is purely for aesthetic reasons: Otherwise, we could only guarantee an $(m-1)$-well obtained from $M$ by additionally deleting the outer cycle, for it may not be part of the initial $\Delta_{C_m}$-society. 

Thus, given an annotated graph $(G,T)$ with $m$-locus $(\rho,M)$ and a vital $T$-linkage $\LLL,$ we write $\left(G(M,\LLL),\Omega_{\Delta_{C_m}},\rho_{\Delta_{C_m}},M,\LLL(\Delta_{C_m})\right)$ to mean the $m$-well obtained by \zcref{obs:well_from_locus}. Note that it is uniquely defined up to reversing $\Omega_{\Delta_{C_m}}$. It is not difficult to see that if $(G,T),\LLL,(\rho,M)$ is exhausted and $M$ is tight, then $(G(M,\LLL),\Omega_{\Delta_{C_m}},\rho_{\Delta_{C_m}},M,\LLL(\Delta_{C_m}))$ is dry.

\begin{lemma}\label{lem:no_apex_vital_from_exhausted_to_dry}
 Let $\Sigma$ be a surface. Let $m \geq 1$ and let $(G,T)$ be an annotated graph admitting a vital $T$-linkage $\LLL$. Let $(\rho,M)$ be an $m$-locus of $(G,T)$. If $(G,T),\LLL,(\rho,M)$ is exhausted and $M$ is tight, then $\left(G(M,\LLL),\Omega_{\Delta_{C_m}},\rho_{\Delta_{C_m}},M,\LLL({\Delta_{C_m}})\right)$ is dry.
\end{lemma}
\begin{proof}
 Towards a contradiction assume the contrary. By assumption of the theorem, in particular $M$ is tight, \zcref{lem:drywell} implies the existence of an $m$-well $(G',\Omega_{\Delta_m},\rho',\CCC,\LLL')$ where $G'\subseteq G(M,\LLL)$ with $G'=\bigcup\CCC \cup \bigcup\LLL'$ and $\rho'$ is the rendition of $G'$ induced by $\rho$. Further $\tau(\LLL') = \tau(\LLL(\Delta_{C_m}))$ but $\LLL'\neq \LLL(\Delta_{C_m})$ as we assumed that the well was not dry. Since their patterns agree, there is a bijection $\alpha: \LLL(\Delta_{C_m})\to \LLL'$ such that for $P \in \LLL(\Delta_{C_m}),$ the path $P$ and $\alpha(P)$ have agreeing endpoints, each of which are nodes of $V(\Omega_{\Delta_m}),$ and their traces are otherwise disjoint from~$\mathsf{cl}(\Sigma\setminus \Delta_{C_m})$.

 Let $\LLL=\{L_1,\ldots,L_t\},$ say, and for every $i \in [t]$ let $\{P_{L_i}^1,\ldots,P_{L_i}^{k_i}\} \subseteq \LLL(\Delta_{C_m})$ for some $k \geq 0$ be maximal such that $P_{L_i}^p\subseteq L$ for every $p \in [k_i]$ and every $i \in [t]$. 
 By definition of $\LLL(\Delta_{C_m})$ and \zcref{obs:sharing_ends_in_restricted_linkage} we derive the following.
 \begin{beautifulclaim}\label{lem:no_apex_vital_from_exhausted_to_dry_claim}
 For every $i\in[t],$ $\{P_{L_i}^1,\ldots,P_{L_i}^k\}$ is a set of internally disjoint paths where each path has its endpoints in $V(\Omega_{\Delta_{C_m}})$. Further, for every $j\in [t]\setminus\{i\}$ and any $(p,q)\in [k_i]\times[k_j],$ the paths $P_{L_i}^p$ and $P_{L_j}^q$ are disjoint.
 \end{beautifulclaim}
 For every $i \in [t]$ define the following.
 Let ${L_{i,0}}^* \coloneqq L_i$. For $p \in [k_i],$ let $L^*_{i,p}$ be obtained from $L^*_{i,p-1},$ by replacing $P_{L_i}^p$ by $\alpha(P_{L_i}^p)$. Since $P_{L_i}^p$ and $\alpha(P_{L_i}^p)$ have the same endpoints, we derive that $L^*_{i,p}$ is a path in $G$ with the same endpoints as $L_i$ (assuming by induction that $L^*_{i,p-1}$ satisfied the same). Similarly, using \zcref{lem:no_apex_vital_from_exhausted_to_dry_claim} one maintains that $L_{i,p}^*$ and $L_{j,q}^*$ are vertex disjoint for every iteration where $q\in[ k_j]$ for some $j \in [t]$. Finally, let $\LLL^*\coloneqq \{L^*_{i,k_i}\mid 1\leq i \leq t\}$. Then $\LLL^*$ is a linkage in $G$ with $\tau(\LLL^*)=\tau(\LLL)$ but $\LLL^*\neq \LLL$ since $\LLL(\Delta_{C_m})\neq \LLL'$ by assumption. This is a contradiction to $\LLL$ being vital, concluding the proof. 
\end{proof}
\zcref{lem:no_apex_vital_from_exhausted_to_dry} brings to light a clear pattern for how a vital linkage interacts with a sequence of concentric cycles $M$ inside the disc bounded by its outer cycle. It turns out that we can leverage the above, to tame the behaviour of $\LLL$ ``outside'' of the cycles of $M$ as we discuss next. To this extent, we introduce \emph{segments} of a linkage given a blank rendition, which we define as follows; see \cref{fig:in_vs_out_bucket} for a schematic illustration.

\begin{definition}\label{def:segments}
 Let $\Sigma$ be a surface and $\rho$ a blank $\Sigma$-rendition of an annotated graph $(G,T)$. Let $P$ be a path in $G$. A \emph{segment of $P$ (in $\rho$)} is a grounded subpath $P' \subseteq P$ with both endpoints in vortex-boundaries such that no internal vertex of $P'$ is part of a vortex-boundary.
\end{definition}
Note that by definition, every segment is a loop or a link; recall \zcref{def:links}. The following is an easy observation.

\begin{observation}\label{obs:segments}
 Let $\Sigma$ be a surface and $\rho$ a blank $\Sigma$-rendition of some annotated graph~$(G,T)$. Let $\LLL$ be a $T$-linkage and $\SSS$ the set consisting of all segments of paths in $\LLL$. Then $\SSS$ is a set of internally disjoint paths and every endpoint of a path in $\SSS$ is part of at most two paths in~$\SSS$.
\end{observation}

Recall that, given a disc $\Delta,$ a sequence of concentric cycles $\langle C_1,\ldots,C_m\rangle$ in some $\Delta$-rendition $\rho$ and a grounded path $P$ whose trace $\gamma_P$ is contained in $\Delta$ with its endpoints on the boundary of $\Delta$ but internally disjoint from the boundary of $\Delta,$ we denote by $\Delta_P$ the unique disc of $\Delta - \gamma_P,$ after taking the topological closure of the component, that does not contain the trace of $C_1$.

We refine on this definition as follows.
\begin{definition}\label{def:disc_of_restricted_linkage}
 Let $\Sigma$ be a surface. Let $m \geq 1$. Let $(\rho,M$) be an $m$-locus of some annotated graph $(G,T)$ such that $\rho$ is a $\Sigma$-skeleton and $M=\langle C_1,\ldots,C_m\rangle$. Let $\LLL$ be a vital $T$-linkage. Let $P \in \LLL(\Delta_m)$ with trace $\gamma_P$. Then we define $\Delta^M_P$ to be the unique disc of $\Delta_{C_m} - \gamma_P,$ after taking the topological closure of the component, that does not contain the trace of $C_1$. We write $\Delta_P$ for simplicity if $M$ is clear from context.
\end{definition}
Note that $\Delta_P$ is $\rho$-aligned by construction. The following is an easy observation.

\begin{observation}\label{obs:discs_of_paths_in_locus_are_wellbehaved}
 Let $\Sigma$ be a surface. Let $m \geq 1$. Let $(\rho,M$) be an $m$-locus of some annotated graph $(G,T)$ such that $\rho$ is a $\Sigma$-skeleton and $M=\langle C_1,\ldots,C_m\rangle$. Let $\LLL$ be a vital $T$-linkage. Let $P,P' \in \LLL(\Delta_{C_m})$. Then either $\Delta_P$ and $\Delta_{P'}$ are internally disjoint or one is contained in the other.
\end{observation}
\begin{proof}
 By \zcref{obs:sharing_ends_in_restricted_linkage}, $P$ and $P'$ are internally disjoint. Since $\rho$ is a $\Sigma$-rendition, if the trace $\gamma_P$ of $P$ is not internally disjoint from $\Delta_{P'}$ then $\gamma_P \subsetneq \Delta_{P'}$. To see this note that $\gamma_P$ is internally disjoint from the trace of $C_m$ as well as from the trace of $P'$. One easily verifies that $\Delta_{P} \subseteq \Delta_{P'}$. Switching the roles of $P$ and $P'$ concludes the proof. 
\end{proof}

The following lemma proves that, given an exhausted instance $(G,T),\LLL,(\rho,M)$ such that $M$ is tight, no subpath of $L\in \LLL$ whose trace is contained in $\Delta_{C_m}$ may intersect $C_m$ in more than two components, i.\@e.\@, the traces of the paths do not ``bounce off'' the outer cycle while staying inside the disc. In particular, together with \zcref{obs:sharing_ends_in_restricted_linkage} this implies that $\LLL(\Delta_{C_m})$ is a linkage (as opposed to internally disjoint).
\begin{lemma}\label{lem:no_apex_dry_no_bounce_ins}
 Let $\Sigma$ be a surface. Let $m \geq 2$. Let $(\rho,M$) be an $m$-locus of some annotated graph $(G,T)$ such that $M= \langle C_1,\ldots,C_m\rangle$ is tight and $\rho$ is a $\Sigma$-skeleton. Let $\LLL$ be a vital $T$-linkage. Let $L \in \LLL$ and let $S\subseteq L$ be a segment. Let $P\subseteq S$ be a subpath such that its trace is contained in $\Delta_{C_m}$ and its endpoints lie on the boundary of $\Delta_{C_m}$. Then $P\cap C_m$ consists of at most $2$ components, each of which is a single vertex of $C_m$.
 
 In particular, there is a unique path $P_{\Delta_m} \in \LLL(\Delta_m)$ that is a subpath of $P$ and $\LLL(\Delta_m)$ is a linkage. 

\end{lemma}
\begin{proof}
 Note that the first part of the lemma implies the second, i.\@e.\@, if $P\cap C_m$ consists of at most two components then, since $P$ was arbitrary, we derive that $\LLL(\Delta_m)$ is a linkage, since, a priori, if two paths of $\LLL(\Delta_m)$ intersect, then they belong to a common segment and intersect at a common endpoint by \cref{obs:sharing_ends_in_restricted_linkage}. Since $P \cap C_m$ consists of at most two components this is impossible; we are left to prove the former.

 \begin{figure}[ht]
 \centering
 \begin{tikzpicture}
 \pgfdeclarelayer{background}
		 \pgfdeclarelayer{foreground}
			
		 \pgfsetlayers{background,main,foreground}

 \begin{pgfonlayer}{background}
 \pgftext{\includegraphics[width=12cm]{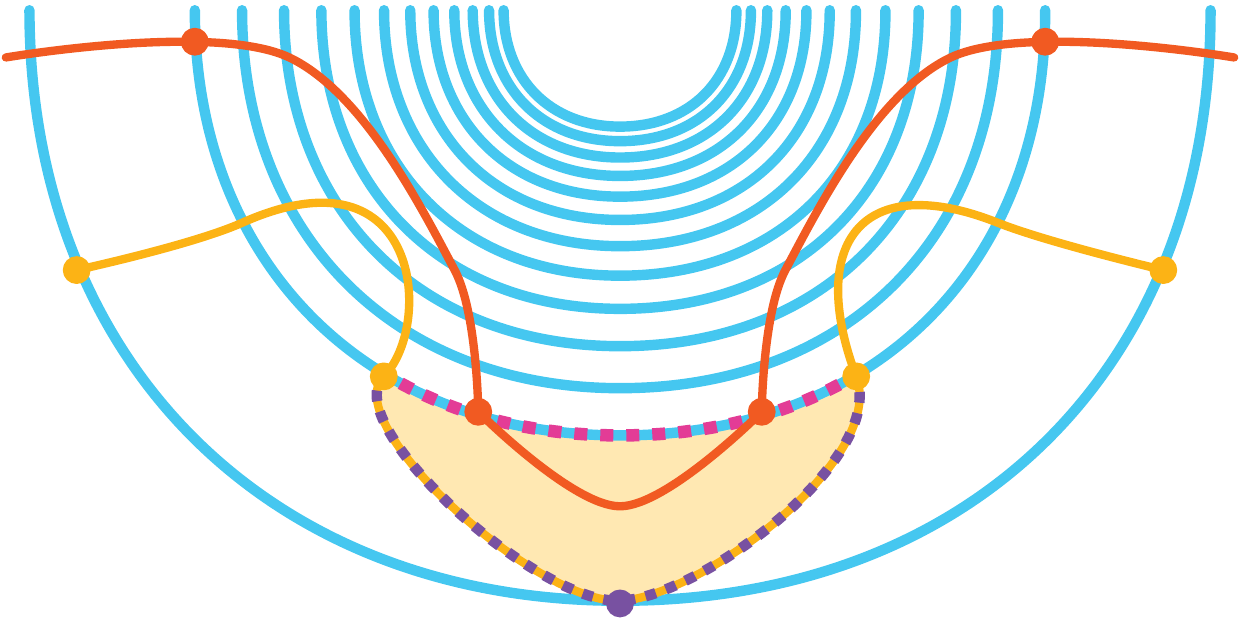}};
 \end{pgfonlayer}{background}
			
 \begin{pgfonlayer}{foreground}

 \node (X) at (0,-2.5) [draw=none] {$x$};
 \node (Y1) at (-1.6,-1.2) [draw=none] {$y_1$};
 \node (Y2) at (1.65,-1.2) [draw=none] {$y_2$};
 \node (X1) at (-5.6,0.4) [draw=none] {$x_1$};
 \node (X2) at (5.65,0.4) [draw=none] {$x_2$};
 \node (W1) at (-4.4,2.3) [draw=none] {$w_1$};
 \node (W2) at (4.4,2.3) [draw=none] {$w_2$};
 \node (US) at (-2.7,-0.7) [draw=none] {$u_s$};
 \node (UT) at (2.7,-0.7) [draw=none] {$u_t$};
 \node (PPRIME) at (0.1,-1.55) [draw=none] {$P'$};
 \node (P) at (-6.2,2.55) [draw=none] {$P$};
 \node (CM) at (3.2,-2.3) [draw=none] {$C_m$};
 \node (CMMINUS1) at (3.6,-0.1) [draw=none] {$C_{m-1}$};
 
 \end{pgfonlayer}{foreground}
 \end{tikzpicture}
 \caption{A schematic illustration of the proof of \cref{lem:no_apex_dry_no_bounce_ins}. The path $P$ with intersections $x_1,u_s,x,u_t,x_2$ is drawn in {yellow}, each of the intersection points are highlighted. The path $P_s^t$ is marked by a purple dashed line, whereas $C_{m-1}^{s,t}$ is highlighted with a {pink} dashed line, with the disc $\Delta^*$ marked in {yellow}. Further, the path $P'$ is highlighted in dark orange.}
 \label{fig:no_bounce_inside}
 \end{figure}
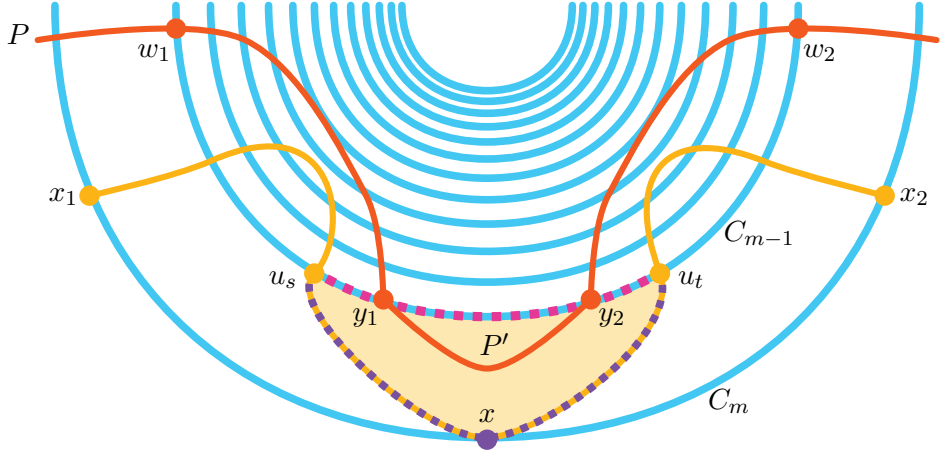
 
 Towards a contradiction assume the lemma is wrong and let $P\subseteq S$ be subpath as in the theorem such that its intersection with $C_m$ results in more than two components in particular $\Abs{V(P)}\geq 3;$ note that since the instance is exhausted, $P$ shares no edge with $C_m$ whence $P\cap C_m$ consists of isolated vertices. Further assume that $P$ is as short as possible satisfying the above. Since $\rho$ is a $\Sigma$-skeleton, every edge that is not part of the graph of a vortex is grounded. In particular, $P$ is grounded and we derive that $P$ starts in a vertex of $V(C_m),$ say $x_1,$ and ends in a vertex of $V(C_m),$ say $x_2$. Let $p \geq 2$ be the length of $P$~--~it cannot be $1$ by the previous assumption~--~and let $\langle x_1,u_1,\ldots,u_{p-1},x_2\rangle$ be the vertices of $P$ as they appear in $P$. Since we chose $P$ as short as possible, we derive that $x_1u_1, u_{p-1}x_2\notin E(C_m)$ for otherwise, $\langle u_1,\ldots,x_2\rangle$ or $\langle x_1,\ldots,u_{p-1}\rangle$ would be the vertices of subpaths of $P$ that are shorter than $P,$ contradicting the lemma and hence our assumption.
 
 By \zcref{def:linkage_restriction_to_a_disc} of $\LLL(\Delta_{C_m})$ and the assumptions on $P,$ we derive that $P$ is the concatenation of paths $P_1,\ldots,P_\ell \in \LLL(\Delta_{C_m})$ for some $\ell\geq 1,$ i.\@e.\@, $P=P_1\circ \ldots \circ P_\ell$. Note that for all $i \in [\ell],$ $P_i$ has its endpoints on $V(C_m),$ since $\rho$ is a $\Sigma$-skeleton and the trace of $P_i$ has its endpoints on the trace of $C_m$.

 \begin{beautifulclaim}\label{lem:no_apex_dry_no_bounce_claim_P1_P2}
 It holds $\ell=2$ and $\Abs{V(P) \cap V(C_m)} = 3$. Further $\Abs{V(P_i)}\geq 2$ and $E(P_i)\cap E(C_m) = \emptyset$ for $i=1,2$.
 \end{beautifulclaim}
 \begin{claimproof}
 Clearly $\ell > 1,$ since otherwise $P$ would be internally disjoint from $V(C_m)$ and thus $P\cap C_m = \{x_1,x_2\}$. Next note that for $i\in [\ell-1],$ no $P_i$ consists of a single vertex by \zcref{def:linkage_restriction_to_a_disc} of $\LLL(\Delta_m),$ since otherwise $P_i\circ P_{i+1}$ would be a longer path with both endpoints on the boundary of $\Delta_{C_{m-1}}$ that is otherwise disjoint from it refuting maximality of $P_i$. Similarly, $P_\ell$ is not a single vertex for else $P_{\ell-1}$ would be longer again contradicting its existence.

 This implies that $\ell\leq 2$ since otherwise the concatenation $P'=P_1 \circ P_2 \subseteq P$ would be a shorter path satisfying the assumptions of the lemma and $P' \cap C_m$ consists of at least $3$ components by \zcref{obs:sharing_ends_in_restricted_linkage}.

 To see that $\Abs{V(P_m) \cap V(C_m)} = 3,$ note that we already established that $P_1,P_2$ have their endpoints on $V(C_m)$. Since neither of the two paths consist of a single vertex, and neither has an edge in common with $E(C_m),$ we derive that $\Abs{V(P_1) \cap V(C_m)}= 2$ as well as $\Abs{V(P_2) \cap V(C_m)}= 2$. Finally, since they share exactly one common endpoint we derive that $\Abs{V(P) \cap V(C_m)} = 3$.
 \end{claimproof}

 Thus let $x \in V(P) \cap V(C_m)\setminus\{x_1,x_2\}$. Let $\Omega_{\Delta_{C_m}}=\langle v_1,\ldots,v_m\rangle$ be a linear order on $V(C_m)$ such that $v_1=x_1$ and such that $v_iv_{i+1}\in E(C_m)$ for every $i\in [m-1]$ and $v_mv_1\in E(C_m)$. Without loss of generality assume that $x_1,x,x_2$ appear in that order in $\Omega_{\Delta_{C_m}},$ else reverse the order, i.\@e.\@, choose $\Omega_{\Delta_{C_m}}=\langle v_1,v_m,\ldots,v_2\rangle$ which is possible since $C_m$ is a cycle.

 By \zcref{lem:no_apex_vital_from_exhausted_to_dry} and the assumption of the theorem, we derive that $\WWW\coloneqq (G(M,\LLL),\Omega_{\Delta_{C_m}},\rho_{\Delta_{C_m}},M,\LLL(\Delta_{C_m}))$ is dry. 

 \begin{beautifulclaim}\label{lem:no_apex_dry_no_bounce_claim_touch_m-1}
 Let $i=1,2,$ then $1\leq \Abs{V(P_i)\cap V(C_{m-1})} \leq 2$.
 \end{beautifulclaim}
 \begin{claimproof}
 Let $i=1,2$ be arbitrary. Since $\WWW$ is dry and since $P_i\cap C_m$ consists of two components, none of which is a single vertex or subpath of $C_m$ as given by \zcref{lem:no_apex_dry_no_bounce_claim_P1_P2}, we derive by $ii)$ of the definition of dryness that there exists $j \in [m-1]$ such that $P_i \cap C_j$ is a single path, in particular a single vertex since our instance is exhausted. By $iv)$ of the definition of dryness we derive the claim, again noting that since the instance is exhausted each intersection with $C_{m-1}$ is a single vertex.
 \end{claimproof}

 By \zcref{lem:no_apex_dry_no_bounce_claim_touch_m-1} we may choose $u_{s} \in V(P_1)\cap V(C_{m-1})$ for $ s \in [p-1]$ maximal. Similarly, we may choose $u_t \in V(P_2)\cap V(C_{m-1})$ for $t \in [p-1]$ minimal. Since $P=P_1\circ P_2$ and since $V(P_1)\cap V(P_2) = \{x\}$ with $x \notin V(C_{m-1})$ we derive that $s<t$ and in particular $u_s,x,u_t$ appear in that order when traversing $P$ from $x_1$ to $x_2$. Let $C_{m-1}^{s,t} \subsetneq C_{m-1}$ be the unique subpath with endpoints $u_s,u_t$ and let $P_s^t\subseteq P$ be the unique subpath whose vertices are ordered via $\langle u_s,\ldots,u_t\rangle$. Let $C^* = P_{s}^t\cup C_{m-1}^{s,t};$ see \cref{fig:no_bounce_inside} for a schematic illustration of the above defined objects.

 \begin{beautifulclaim}\label{lem:no_apex_dry_no_bounce_claim_C*_cycle}
 $C^*= P_{s}^t\cup C_{m-1}^{s,t}$ is a cycle satisfying $V(C^*) \cap V(C_m) = \{x\}$. Further, the trace of $C^*$ bounds a vortex-free disc $\Delta^*$ satisfying $\Delta^* \subsetneq \Delta_{C_m}$. 
 \end{beautifulclaim}
 \begin{claimproof}
 Clearly $C^*$ is a cycle since $P_{s}^t$ and $C_{m-1}^{s,t}$ are internally disjoint by construction. Further, by construction of $C^*$ we have $V(C^*) \cap V(C_m) = \{x\}$. Since both the traces of $P_{s}^t$ and $C_{m-1}^{s,t}$ are contained in $\Delta_{C_m}$ we derive that the trace of $C^*$ is and hence the claim.
 \end{claimproof}
 Let $\Delta^*$ be the ``$C^*$-disc'' as given by \zcref{lem:no_apex_dry_no_bounce_claim_C*_cycle}; see the \textcolor{Mustard}{yellow} region in \cref{fig:no_bounce_inside}.

 \begin{beautifulclaim}
 There exists $P' \in \LLL(\Delta_m)$ with $V(P')\cap V(C_{m-1}^{s,t}) \neq \emptyset$.
 \end{beautifulclaim}
 \begin{claimproof}
 Since $\LLL$ is vital the above is true unless $C_{m-1}^{s,t}=(u_s,u_t)$ is a single edge $u_su_t \in E(C_{m-1})$. In that case there is a path $P^*$ whose vertices are ordered 
 via $\langle x_1,u_1,\ldots,u_s,u_t,\ldots,u_{p-1},x_2\rangle$ with the same endpoints as $P$ but $P^*\neq P$ for $x \notin V(P^*)$. Finally replacing $P$ with $P^*$ in $L$ results in a new linkage $\LLL^*$ in $G$ with $\tau(\LLL^*)=\tau(\LLL)$ but $\LLL^*\neq \LLL;$ this contradicts the vitality of~$\LLL$.
 \end{claimproof}

 Let $P'$ be as in the claim. We claim that $P'$ violates $iv)$ in the definition of dryness of $\WWW;$ see the \textcolor{orangepeel}{orange} path in \cref{fig:no_bounce_inside}.

 \begin{beautifulclaim}\label{lem:no_apex_dry_no_bounce_claim_end}
 $P'\cap C_{m-1}$ consists of at least $3$ components.
 \end{beautifulclaim}
 \begin{claimproof}
 To see this, note that since $\rho$ is a $\Sigma$-skeleton, we derive that $P'$ has both its endpoints on $V(C_m)$ since its trace does and every vertex of $P'$ is a node of the rendition; call them $w_1$ and $w_2$ respectively. Let $y_1,y_2 \in V(P')\cap V(C_{m-1}^{s,t})$~--~possibly equal~--~be such that $P_1^*,P_2^* \subseteq P'$ start in $y_1,y_2,$ respectively, and end in $w_1,w_2$ respectively, and such that both are internally disjoint from each other as well as internally disjoint from $C_{m-1}^{s,t}$. Clearly, these paths exist by choosing $y_1,y_2$ as the closest vertices to $u_s$ and $u_t$ along $C_{m-1}^{s,t}$ in $V(P') \cap V(C_{m-1}),$ respectively. 

 Let $i=1,2$ be arbitrary. We claim that $P_i^*$ and $C_{m-1}$ are not internally disjoint. To see this, assume the contrary towards a contradiction and recall that $\Delta^*$ is a disc. Let $\gamma_i$ be the trace of~$P_i^*$. We are to show that $\gamma_i$ has an interior point on the boundary of $\Delta_{C_{m-1}}$ which is essentially the trace of $C_{m-1}$. If not, this implies that $\gamma_i$ is internally disjoint from $\Delta_{C_{m-1}},$ since it has one endpoint in $V(C_{m})$. This in turn means that the trace of $\gamma_i$ is contained in $\mathsf{cl}(\Delta_{C_m}\setminus \Delta_{C_{m-1}})$. Since one end of $\gamma_i$ lies on the trace of $C_{m-1}^{s,t}$ we derive that $\gamma_i$ has one end on the boundary of $\Delta^*$ and since further $\Delta^* \subseteq \mathsf{cl}(\Delta_{C_m}\setminus \Delta_{C_{m-1}})$ but also $\gamma_i\subseteq \mathsf{cl}(\Delta_{C_m}\setminus \Delta_{C_{m-1}})$ we derive that $\gamma_i \subseteq \Delta^*$. To see this, note that we assumed that $P_i^*$ is internally disjoint from $C_{m-1}$ and hence its trace is internally disjoint from the trace of $C_{m-1}$. Since further $\LLL(\Delta_m)$ is a set of internally disjoint paths by \zcref{obs:sharing_ends_in_restricted_linkage} we have that $P_i^*$ is internally disjoint from $P_s^t$ and combining both with \zcref{lem:no_apex_dry_no_bounce_claim_C*_cycle} we derive that $\gamma_i$ is internally disjoint from the boundary of $\Delta^*$ but also from $\Delta_{C_{m-1}}$. Since $y_i$ is an interior point of the trace of $C_{m-1}^{s,t},$ we derive that $\gamma_i$ has an interior point in $\mathsf{cl}(\Delta_{C_m}\setminus \Delta_{C_{m-1}})$ and thus we derive that it is contained in $ \mathsf{cl}(\Delta_{C_m}\setminus \Delta_{C_{m-1}})$ as claimed.
 
 Finally, again applying \zcref{lem:no_apex_dry_no_bounce_claim_C*_cycle} we derive that $V(C^*)\cap V(C_m) = \{x\}$. But $\gamma_i \subseteq \Delta^*$ together with the fact that $P_i^*$ ends in $w_i \in V(C_m)$ imply that $w_i = x;$ this is a contradiction to \zcref{obs:sharing_ends_in_restricted_linkage} since then three distinct paths of $\LLL(\Delta_m),$ namely $P_i^*,P_1,$ and $P_2$ share a common end. Thus, our initial assumption was wrong whence $P_i^*$ and $C_{m-1}$ are not internally disjoint.

 Since $i=1,2$ was arbitrary, $P'\cap C_{m-1}$ consists of at least three components, one of which contains $y_1$ (possibly $= y_2$).
 \end{claimproof}

 \zcref{lem:no_apex_dry_no_bounce_claim_end} is a contradiction to the fact that $\WWW$ is dry. In particular, our assumption that $P$ intersects $C_m$ in more than two components was wrong, concluding the proof. 
\end{proof}

While \zcref{lem:no_apex_dry_no_bounce_ins} tames the behaviour of segments inside $\Delta_{C_m},$ the next lemma tames their behaviour outside of $\Delta_{C_m}$. Intuitively speaking, the following lemma proves that given an exhausted instance $(G,T),\LLL,(\rho,M)$ such that $M$ is tight, the trace of a subpath of $L\in \LLL$ that lies outside of $\Delta_{C_m}$ with both endpoints on its boundary has to be an essential curve. That is, paths do not ``jump'' on the outer cycle $C_m$ without really ``using'' the surface $\Sigma$.

\begin{figure}[ht]
 \centering
 \begin{tikzpicture}
 \pgfdeclarelayer{background}
		 \pgfdeclarelayer{foreground}
			
		 \pgfsetlayers{background,main,foreground}

 \begin{pgfonlayer}{background}
 \pgftext{\includegraphics[width=10cm]{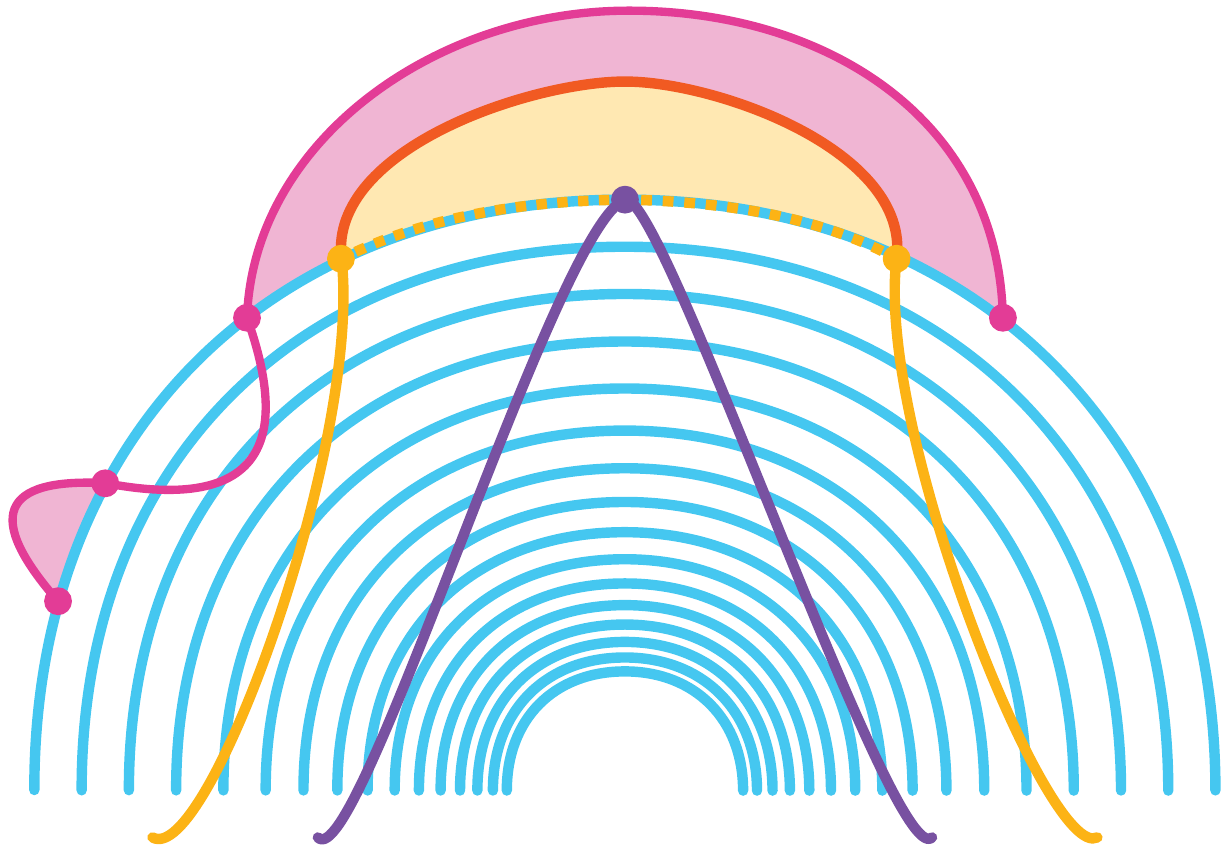}};
 \end{pgfonlayer}{background}
			
 \begin{pgfonlayer}{foreground}

 \node (X) at (0.1,2.1) [draw=none] {$x$};
 \node (U) at (-2.45,1.45) [draw=none] {$u$};
 \node (V) at (2.55,1.45) [draw=none] {$v$};
 \node (S) at (-3.95,-3.2) [draw=none] {$S$};
 \node (SPRIME) at (-2.65,-3.2) [draw=none] {$S'$};
 \node (PPRIME) at (-3.1,1.9) [draw=none] {$P'$};
 \node (P) at (-0.25,2.5) [draw=none] {$P$};
 
 \end{pgfonlayer}{foreground}
 \end{tikzpicture}
 \caption{A schematic illustration of the proof of \cref{lem:no_apex_dry_no_bounce_out}. The $\Delta_{C_m}$-bumps are highlighted by pink and orange regions. The path $P$ with vertices $u,v$ is highlighted in red with corresponding disc $\Delta^P$ highlighted in {yellow}. Furthermore, the vertex $x$ of the segment $S'$ on $C_P \cap C_m$ is highlighted in purple.}
 \label{fig:outer-bounce}
\end{figure}
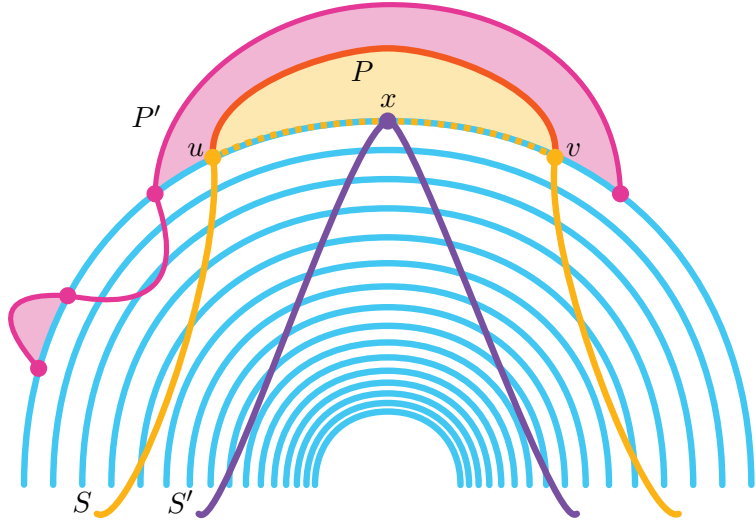

\begin{definition}\label{def:boundary_jump}
 Let $\Sigma$ be a surface and $m \geq 1$. Let $(\rho,M$) be an $m$-locus of some annotated graph $(G,T)$ with $M = \langle C_1,\ldots,C_m\rangle$ such that $\rho$ is a $\Sigma$-skeleton. Let $P \subseteq G$ be a path with at least one edge and both endpoints $u,v$ in $V(C_m)$ and such that its trace is contained in $\mathsf{cl}(\Sigma \setminus \Delta_{C_m})$. If there exists a $\rho$-aligned vortex-free disc $\Delta \subseteq \Sigma$ such that $\Delta_{C_m} \subsetneq \Delta$ and such that the trace of $P$ is contained in $\Delta,$ we call $P$ \emph{a $C_m$-boundary jump}. Let $C_m^1,C_m^2 \subsetneq C_m$ be internally disjoint subpaths such that they intersect exactly in $\{u,v\}$ and satisfy $C_m^1 \cup C_m^2 = C_m$. Then for $i=1,2,$ $C_i' \coloneqq C_m^i \cup P$ is a cycle whose trace is contained in $\Delta$ and for exactly one $i \in \{1,2\}$ it holds that the trace of $C_i'$ bounds a vortex-free disc $\Delta' \subseteq \Delta$ that is internally disjoint from $\Delta_{C_m}$. We define $\Delta^P$ to be $\Delta',$ if it exists, and we call $\Delta^P$ \emph{the $\Delta_{C_m}$-bump by $P$}. 
\end{definition}
Note that $\Delta^P$ is $\rho$-aligned by construction.
See the \textcolor{lightsalmonpink}{pink} and \textcolor{orangepeel}{orange} highlighted regions in \cref{fig:outer-bounce} for an illustration of $\Delta_{C_m}$-bumps.

\begin{lemma}\label{lem:no_apex_dry_no_bounce_out}
 Let $\Sigma$ be a surface. Let $m \geq 2$. Let $(\rho,M$) be an $m$-locus of some annotated graph $(G,T)$ such that $M= \langle C_1,\ldots,C_m\rangle$ is tight and $\rho$ is a $\Sigma$-skeleton. Let $\LLL$ be a vital $T$-linkage and let $(G,T),\LLL,(\rho,M)$ be exhausted. Let $P\subseteq G$ be a path on at least one edge with both its endpoints on $V(C_m)$. Then $P$ is not a $C_m$-boundary jump.
\end{lemma}
\begin{proof}
 Towards a contradiction assume that there exist $C_m$-boundary jumps. Let $\PPP$ be the set of all $C_m$-boundary jumps and let $\{\Delta^P \mid P \in \PPP\}$ be the set of $\Delta_{c_m}$-bumps by $P$ for each $P \in \PPP$. 

 \begin{beautifulclaim}\label{lem:no_apex_dry_no_bounce_out_claim_discs_of_jumps_disjoint}
 Let $P,P' \in \PPP$ be distinct, then either $\Delta^{P}$ and $\Delta^{P'} $ are internally disjoint or one is contained in the other.
 \end{beautifulclaim}
 \begin{claimproof}
 Since $\rho$ is a $\Sigma$-skeleton, and in particular a $\Sigma$-rendition, the traces of $P$ and $P'$ are internally disjoint whence $P$ and $P'$ are internally disjoint. Let $C_P,C_{P'} \subsetneq C_m$ be the unique subpaths of $C_m$ such that the trace of $C_P \cup P$ bounds $\Delta^P$ and the trace of $C_{P'}\cup P'$ bounds $\Delta^{P'}$. Then either $C_P$ and $C_P'$ are internally disjoint or one is a subpath of the other, say $C_P \subsetneq C_{P'}$ in that case. Finally by definition of $C_m$-boundary jump, the traces of $P$ and $P'$ are both internally disjoint from the traces of $C_{P}$ and $C_{P'}$. Thus, combining all of the above, if $\Delta^{P}$ and $\Delta^{P'} $ are \emph{not} internally disjoint, then this means that two paths of $\{P,P',C_P,C_{P'}\}$ are \emph{not} internally disjoint. Since this is only possible for $C_{P}$ and $C_{P'}$ we derive that $C_{P} \subsetneq C_{P'}$. Finally, we derive that $\Delta^P \subsetneq \Delta^{P'},$ concluding the proof; see \cref{fig:outer-bounce} for a schematic illustration of possible arrangements of $\Delta_{C_m}$-bumps.
 \end{claimproof}
 Applying \zcref{lem:no_apex_dry_no_bounce_out_claim_discs_of_jumps_disjoint} we derive the existence of a $C_m$-boundary jump $P \in \PPP$ such that for all $P' \in \PPP$ we derive that either $\Delta^P$ and $\Delta^{P'}$ are internally disjoint, or $\Delta^P \subseteq \Delta^{P'}$.

 Let $u,v$ be the distinct endpoints of $P$ in $V(C_m);$ they exist since $P$ has at least one edge. Let $P_{C_m}$ be the subpath of $C_m$ with ends $u,v$ such that the trace of the cycle $C^*=P\cup P_{C_m}$ bounds~$\Delta^P$. Note that $P_{C_m}$ contains at least one edge. By our choice of $P,$ and the fact that $\rho$ is a $\Sigma$-skeleton, $\Delta^P$ only intersects or contains cells that correspond to edges of $E(C^*)$. Let $\SSS(\LLL)$ be the segments of $\LLL$.

 \begin{beautifulclaim}
 There is a segment $S \in \SSS(\LLL)$ such that $P \subseteq S$. In addition, there is a segment $S' \in \SSS(\LLL)$ such that $V(S') \cap (V(P_{C_m})\setminus \{u,v\}) \neq \emptyset$.
 \end{beautifulclaim}
 \begin{claimproof}
 Since $P$ is a $C_m$-boundary jump we derive that the trace of $P$ is internally disjoint from all vortices of $\rho$ and from $\Delta_{C_m}$. In particular $E(P) \cap E(\bigcup_{i=1}^m C_i) = \emptyset$. Since $(G,T),\LLL,(\rho,M)$ is exhausted we derive that $E(P) \subseteq E(\bigcup \LLL)$ and since its trace is disjoint from the vortices of $\rho$ there is $S \in \SSS(\LLL)$ with $P \subseteq S$ as desired. Furthermore, since $\LLL$ is vital, we derive that $P_{C_m}$ contains at least three vertices. To see this, note that otherwise $P_{C_m}$ consists solely of the edge $uv \in E(C_m),$ but then we could reroute $S$ by replacing $P$ with $uv$. This results in a new $T$-linkage $\LLL^*$ with the same pattern as $\LLL$ but $\LLL \neq \LLL^*$ as a contradiction to the vitality of $\LLL$. 
 Let $x \in V(P_{C_m}) \setminus \{u,v\}$. Again, since $\LLL$ is vital, there exists $L \in \LLL$ such that $x \in V(L)$. Thus there exists $S' \in \SSS(\LLL)$ -- possibly $S'=S$ -- with $x \in V(S')$ whence $V(S') \cap V(P_{C_m}) \neq \emptyset$ concluding the proof.
 \end{claimproof}

 Let $S' \in \SSS(\LLL)$ with $V(S') \cap (V(P_{C_m})\setminus \{u,v\})\neq \emptyset$ as in the previous claim, and let $x \in V(S') \cap (V(P_{C_m})\setminus \{u,v\});$ see the \textcolor{mediumslateblue}{purple} path in \cref{fig:outer-bounce} for an illustration.
 
 We claim that $S'$ witnesses a contradiction to \zcref{lem:no_apex_dry_no_bounce_ins}. To see this, let $P^* \subseteq S'$ be a maximal subpath such that $x \in V(P^*),$ the endpoints $u',v'$ (possibly equal) of $P^*$ are part of $V(C_m)$ and such that the trace of $P^*$ is contained in $\Delta_{C_m}$.

 \begin{beautifulclaim}\label{lem:no_apex_dry_no_bounce_out_claim_triplets}
 $\Abs{\{u',x,v'\}} \geq 3,$ in particular $P^* \cap C_m$ consists of at least $3$ components.
 \end{beautifulclaim}
 \begin{claimproof}
 If $u' = v'$ then $u' = x$ and $v'=x$ and $P^*$ consists of a single vertex. Thus it suffices to prove that $u' \neq x$ and $v' \neq x$. Towards a contradiction assume that $u' = x$: the case $x=v'$ is analogous. Let $S_1',S_2' \subseteq S'$ be the two internally disjoint subpaths agreeing in $x$ such that $S_1' \cup S_2' = S'$. By our choice of $P^*$~--~the maximality assumptions~--~this implies that the traces $\gamma_{S_1'},\gamma_{S_2'}$ of $S_1',S_2'$ must be internally disjoint from $\Delta_{C_m}$. In particular we derive that the trace $\gamma_{S'}$ of $S'$ is contained in $\mathsf{cl}(\Sigma\setminus \Delta_{C_m})$ with both endpoints on vortex boundaries. Recall that for $C^* = C_{P_m} \cup P$ we established that $\Delta^P$ only contains cells that correspond to edges of $E(C^*)$. Thus $\gamma_{S_1'},\gamma_{S_2'}$ can intersect $\Delta^{P}$ only in the trace of $C^*$. Since $x \in V(S_1') \cap V(S_2')$ and since $V(P_{C_m}) \geq 3$ we derive that there is at least one edge $e^* \in E(C_m)$ that is adjacent to $x;$ one easily verifies that there are actually two such edges using \zcref{lem:no_apex_dry_no_bounce_ins} but this is irrelevant. Thus $e^* \in E(S_1') \cup E(S_2')$. This is a contradiction to the fact that $E(\bigcup \LLL) \cap E(C_m) = \emptyset;$ thus our assumption was wrong and $u' \neq x$ concluding the proof.
 \end{claimproof}
 
 Finally, $P^*$ satisfies the assumptions of \zcref{lem:no_apex_dry_no_bounce_ins} whence we derive that $P^* \cap C_m$ consists of at most $2$ components each of which is a single point on $C_m$. In particular $\Abs{V(P^*) \cap V(C_m)} \leq 2$. This contradicts \zcref{lem:no_apex_dry_no_bounce_out_claim_triplets}, concluding the proof.
\end{proof}

With \zcref{lem:no_apex_dry_no_bounce_ins,lem:no_apex_dry_no_bounce_out} at hand, we finally derive that, given a $\rho$-aligned disc $\Delta$ containing some segment $S$ of a path $L\in \LLL,$ the segment does not ``interact'' with the concentric cycles in a ``nested'' way, essentially proving a version of ``dryness'' for segments; compare to the original setting in \cref{fig:in_vs_out_bucket}. This marks the main result of this section.

\begin{figure}[ht]
 \centering
 \begin{tikzpicture}
 \pgfdeclarelayer{background}
		 \pgfdeclarelayer{foreground}
			
		 \pgfsetlayers{background,main,foreground}

 \begin{pgfonlayer}{background}
 \pgftext{\includegraphics[width=10cm]{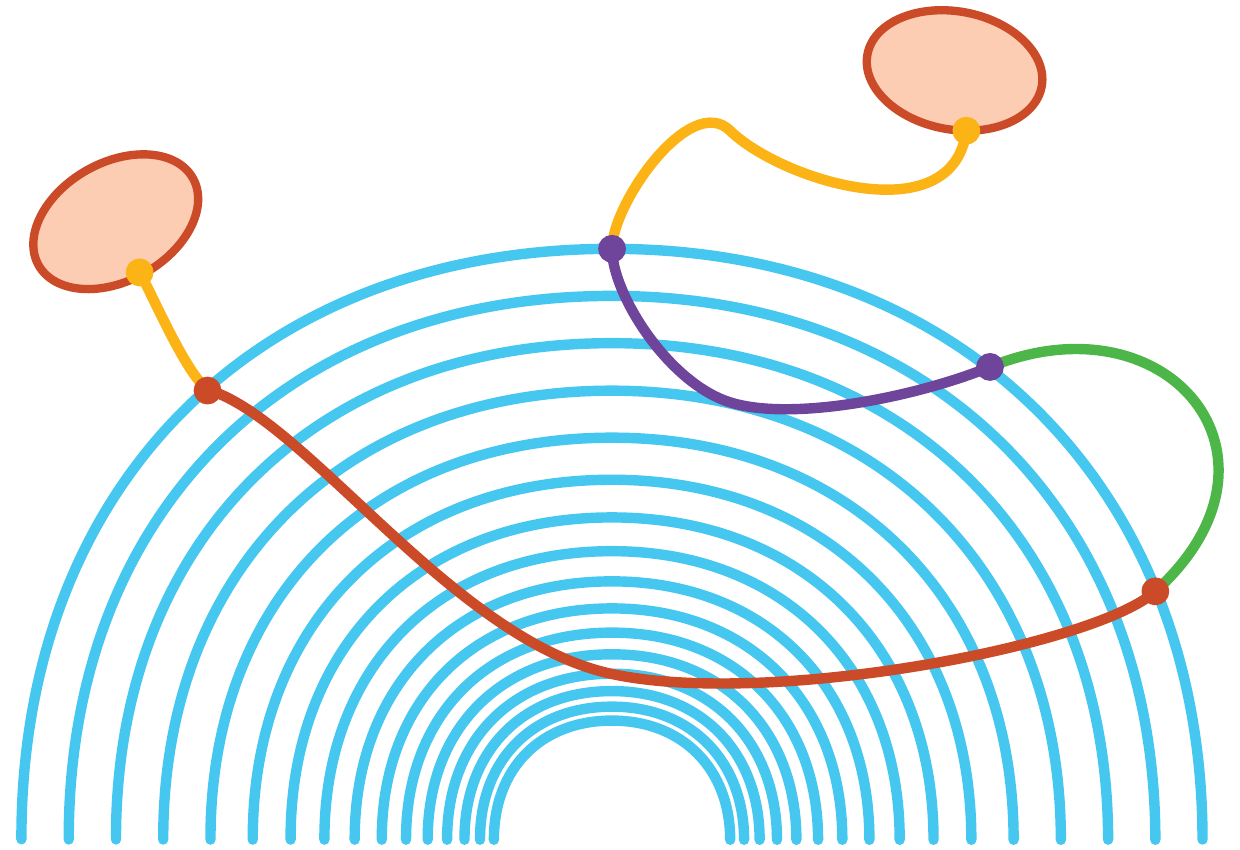}};
 \end{pgfonlayer}{background}
			
 \begin{pgfonlayer}{foreground}

 \node (U1) at (-3.7,0.3) [draw=none] {$u_1$};
 \node (U2) at (4.7,-1.35) [draw=none] {$u_2$};
 \node (V1) at (3,0.75) [draw=none] {$v_1$};
 \node (V2) at (-0.35,1.65) [draw=none] {$v_2$};
 \node (P1) at (0.05,-2.35) [draw=none] {$P_1$};
 \node (P2) at (1.05,0.4) [draw=none] {$P_2$};
 \node (P) at (5.6,0) [draw=none] {$P^* = P$};
 \node (S) at (-3.35,0.9) [draw=none] {$S$};
 
 \end{pgfonlayer}{foreground}
 \end{tikzpicture}
 \caption{A schematic illustration of the proof of \cref{thm:no_apex_dry_then_clean_segments} in the case that $P^* = P$ and $\gamma_P \subseteq \mathsf{cl}(\Sigma \setminus \Delta_{C_m})$. The segment $S$ is drawn in orange with the path $P_1$ drawn in dark orange and $P_2$ drawn in purple. The path $P^*=P$ is highlighted in green.}
 \label{fig:no_jumps_final}
\end{figure}
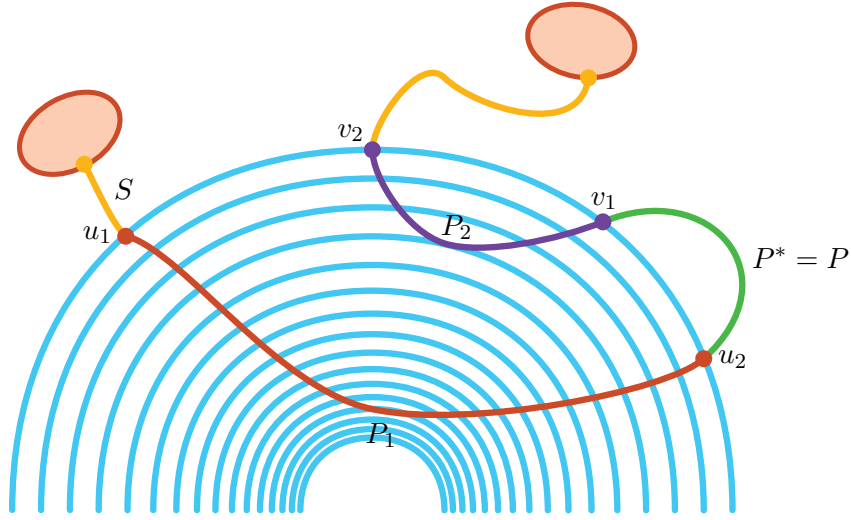

\begin{theorem}\label{thm:no_apex_dry_then_clean_segments}
 Let $\Sigma$ be a surface and $m \geq 2$. Let $(\rho,M)$ be an $m$-locus of some annotated graph $(G,T)$ in $\Sigma$ such that $M= \langle C_1,\ldots,C_m\rangle$ is tight and $\rho$ is a $\Sigma$-skeleton. Let $\LLL$ be a vital $T$-linkage in $G$ and let $(G,T),\LLL,(\rho,M)$ be exhausted. Let $\SSS(\LLL)$ be the set of segments of $\LLL$. Let $\Delta\subseteq \Sigma$ be a vortex-free $\rho$-aligned disc such that for each segment $S\in \SSS(\LLL)$ its trace is either contained in $\Delta$ or is internally disjoint from $\Delta$.
 Let $S \in \SSS(\LLL)$ be such that its trace is contained in $\Delta$. Let $\LLL(\Delta_{C_{m}})$ be the restriction of $\LLL$ to $\Delta_{C_{m}}$. Let $P_1,P_2\in \LLL(\Delta_m)$ be distinct such that $P_1,P_2 \subsetneq S$ and both contain at least one edge each.
 
 Then $\Delta_{P_1} \cap \Delta_{P_2} = \emptyset$.
\end{theorem}
\begin{proof}
 Recall that since $\rho$ is a $\Sigma$-skeleton, all paths in $\LLL(\Delta_m)$ are grounded, and every cell represents a single edge of the graph. Since $(G,T),\LLL,(\rho,M)$ is exhausted and $M$ is dry, \zcref{lem:no_apex_vital_from_exhausted_to_dry} implies that $$\WWW \coloneqq \left(G(M,\LLL),\Omega_{\Delta_{C_m}},\rho_{\Delta_{C_m}},M,\LLL(\Delta_{C_m})\right)$$ is dry.

 Towards a contradiction, assume the lemma to be wrong. Henceforth, let $P_1,P_2 \in \LLL(\Delta_{m})$ be distinct with $P_1,P_2\subsetneq S$ as in the lemma and assume that $\Delta_{P_1}\cap \Delta_{P_2} \neq \emptyset$. Further, let $\Delta$ be as in the theorem, in particular it is vortex-free. Note that since the trace of $S$ is contained in $\Delta$ we derive that the traces of $P_1$ and $P_2$ are contained in $\Delta$. By \zcref{obs:sharing_ends_in_restricted_linkage} we derive that $P_1,P_2$ are internally disjoint, and combined with \zcref{lem:no_apex_dry_no_bounce_ins} we derive that they are in fact completely disjoint. In particular, $\Delta_P$ and $\Delta_{P'}$ are not internally disjoint.
 
 \begin{beautifulclaim}\label{thm:no_apex_dry_then_clean_segments_stack}
 It holds $\Delta_{P_1}\subsetneq \Delta_{P_2}$ or $\Delta_{P_2}\subsetneq \Delta_{P_1}$.
 \end{beautifulclaim}
 \begin{claimproof}
 This follows immediately from \zcref{obs:discs_of_paths_in_locus_are_wellbehaved} together with the fact that they are not internally disjoint.
 \end{claimproof}
 
 By \zcref{thm:no_apex_dry_then_clean_segments_stack} we may assume without loss of generality that $\Delta_{P_1}\subsetneq \Delta_{P_2};$ otherwise rename $P_1$ and $P_2$ accordingly.
 Denote the ends of $P_i$ by $u_i,v_i$ for $i=1,2;$ by assumption of the lemma $P_i$ contains at least one edge thus its endpoints are distinct. Since $P_1,P_2$ are disjoint paths, we derive that $u_1,u_2,v_1,v_2$ are pairwise distinct. Fix an orientation of $S$ such that $P_1$ is visited prior to $P_2$. Then, there exists a subpath $P^*$ with endpoints $x_1,x_2$ where $x_i \in \{u_i,v_i\}$ such that the concatenation $P_1\circ P^* \circ P_2$ is a subpath of $S$. In particular, the trace of $P^*$ is contained in $\Delta;$ see \cref{fig:no_jumps_final} for a schematic illustration.

 We claim that $P^*$ witnesses the existence of a path violating either \zcref{lem:no_apex_dry_no_bounce_ins} or \zcref{lem:no_apex_dry_no_bounce_out}, in either case contradicting the assumption that the theorem does not hold. To this extent, let $P \subseteq P^*$ with one endpoint $x_1,$ and the other endpoint $x \in V(C_m)\setminus \{x_1\}$ such that $P$ is otherwise internally disjoint from $V(C_m)$. To see that this exists, note that $E(P^*) \cap E(C_m) = \emptyset$ by our assumption that $(G,T),\LLL,(\rho,M)$ is exhausted and further $\Abs{V(P^*) \cap V(C_m)} \geq 2$. By construction, the trace $\gamma_P$ of $P$ satisfies the following:
 \begin{itemize}
 \item[$(\star)$] either $\gamma_P \subseteq \Delta_{C_m}$ or $\gamma_P \subseteq \mathsf{cl}(\Sigma\setminus \Delta_{C_m})$.
 \end{itemize} 

 \begin{beautifulclaim}
 $\gamma_P$ is not contained in $\Delta_{C_m}$.
 \end{beautifulclaim}
 \begin{claimproof}
 Assume the contrary. Let $Q=P_1 \cup P \subseteq S,$ then $Q$ is a path and its trace is contained in $\Delta_{C_m}$ with both its endpoints on the boundary of $\Delta_{C_m}$. But $u_1,u_2,x$ are pairwise distinct, and $\{u_1,u_2,x\} \subseteq V(Q) \cap V(C_m)$. In particular $\Abs{V(Q) \cap V(C_m)} \geq 3$ as a contradiction to \zcref{lem:no_apex_dry_no_bounce_ins}.
 \end{claimproof}

\begin{beautifulclaim}
 $\gamma_P$ is not contained in $\mathsf{cl}(\Sigma\setminus \Delta_{C_m})$.
 \end{beautifulclaim}
 \begin{claimproof}
 Assume the contrary. Then by \zcref{def:boundary_jump} together with the fact that $\gamma_P \subseteq \Delta$ with $\Delta$ vortex-free, we derive that $P$ is a $C_m$-boundary jump; note that it admits at least one edge by construction. This is a contradiction to \zcref{lem:no_apex_dry_no_bounce_out}.
 \end{claimproof}
 Combining both claims we derive a contradiction to $(\star),$ concluding the proof.
\end{proof}

\subsection{Controlling vortices}
This subsection is devoted to the proof of technical lemmas needed to show that, given an exhausted instance $(G,T),\LLL,(\rho,M)$ in some surface $\Sigma$ with a ``large enough'' locus with respect to the depth of the vortices, the number of non-simple loops and links in $\SSS(\LLL)$ is bounded in $d$ and the genus of the surface $\Sigma,$ whereas simple loops are guaranteed to not penetrate ``deep'' into $M$. The main result consists of a rather technical analysis on the interaction of $\LLL$ with vortices, combining \cref{thm:Lemma2.5,thm:linking_vortices,thm:no_apex_dry_then_clean_segments}; recall the relevant definitions. 

We will need \emph{$T$-normal linear decompositions} which we define as follows.
\begin{definition}
 Let $(G,T)$ be an annotated graph and let $\alpha = \left(\langle u_1,\ldots,u_n\rangle,\langle U_1,\ldots,U_n\rangle\right)$ be a linear decomposition (possibly with gap) of some society $(H,\Omega)$ with $H \subseteq G$. Let $1 \leq i < j \leq n$ and let $\alpha'$ be the $[i,j]$-restriction of $\alpha$. We call $\alpha'$ \emph{$T$-normal} if for every $t \in T$ either $t \in \bigcap_{p=i}^j U_p$ or $t \notin \bigcup_{p=i}^j U_i$. 
\end{definition}

\begin{figure}[ht]
 \centering
 \begin{tikzpicture}
 \pgfdeclarelayer{background}
		 \pgfdeclarelayer{foreground}
			
		 \pgfsetlayers{background,main,foreground}

 \begin{pgfonlayer}{background}
 \pgftext{\includegraphics[width=10cm]{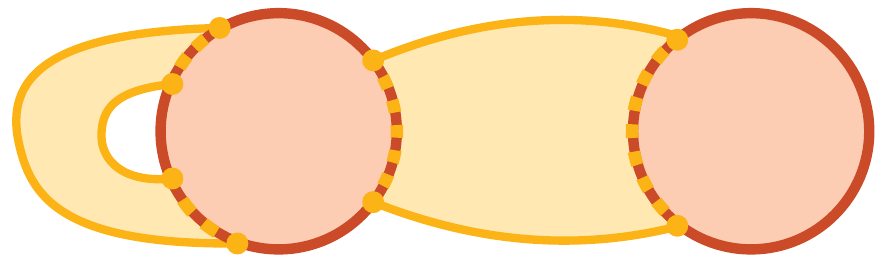}};
 \end{pgfonlayer}{background}
			
 \begin{pgfonlayer}{foreground}

 \end{pgfonlayer}{foreground}
 \end{tikzpicture}
 \caption{Two examples of $(F_1,F_2)$-strips, highlighted in orange, one between two vortices, drawn in red, and one on a single vortex.}
 \label{fig:strips}
\end{figure}

For ease of notation we define what we mean by ``strips'' bounded by segments of vortices; see \cref{fig:strips} for an illustration.

\begin{definition}
 Let $\Sigma$ be a surface and let $\rho$ be a $\Sigma$-rendition of some graph $G$. Let $c_1,c_2 \in C(\rho)$ be vortices (possibly equal) with respective vortex societies $(\sigma(c_1),\Omega_{c_1})$ and $(\sigma(c_2),\Omega_{c_2})$. For $i=1,2,$ let $F_i$ be a segment of $\Omega_{c_i}$ with endpoints $s_i,t_i$ such that $F_1 \cap F_2 = \emptyset$. Let $\gamma_1,\gamma_2$ be $\rho$-aligned curves such that for $i=1,2$ the curve $\gamma_i$ has endpoints $s_i,t_i$ and such that $\gamma_1 \cup \gamma_2$ together with the boundary segments $\eta_1,\eta_w$ of $c_1,c_2$ that contain exactly the nodes of $F_1,F_2,$ respectively, bound a vortex-free disc $\Delta$. Then we call $\Delta$ an \emph{$(F_1,F_2)$-strip (guarded by $\gamma_1,\gamma_2)$}. 
\end{definition}
Note that the $(F_1,F_2)$-strip is $\rho$-aligned as the boundary of every cell is $\rho$-aligned by definition. 

The following is a folklore topological fact transcribed to our setting.
\begin{observation}\label{obs:strip_from_homotopic_paths}
 Let $\Sigma$ be a surface and let $\rho$ be a $\Sigma$-rendition of some graph $G$. Let $c_1,c_2 \in C(\rho)$ be vortices (possibly equal) with respective vortex societies $(\sigma(c_1),\Omega_{c_1})$ and $(\sigma(c_2),\Omega_{c_2})$. Let $\SSS(c_1,c_2)$ be a set of at least two $\rho$-homotopic paths, each with one endpoint in $N(c_1)$ and the other in $N(c_2)$. For $i=1,2$ let $F_i$ be segments of $\Omega(c_i)$ such that $\SSS(c_1,c_2)$ is an $F_1$-$F_2$-linkage. Then there exist $S,S' \in \SSS(c_1,c_2)$ with traces $\gamma,\gamma'$ and disjoint segments $F_i^* \subseteq F_i$ such that $\gamma,\gamma'$ guard an $(F_1^*,F_2^*)$-strip $\Delta,$ and such that $\SSS(c_1,c_2)$ is an $F_1^*$-$F_2^*$-linkage, and the trace of each path in $\SSS(c_1,c_2)$ is contained in $\Delta.$
\end{observation}

The following is the main technical lemma of this section. It is a powerful tool that, given a vital linkage $\LLL,$ leverages \zcref{thm:Lemma2.5} in order to ``control'' the number of homotopic non-simple loops and homotopic links of $\LLL,$ as well as to quantify ``how much'' simple loops can interact with the concentric cycles of some large locus given that the vortices of the instance are properly $d$-linked; recall \zcref{def:segments,def:linked_vortex_normal} and the fact that if a linear decomposition is properly $d$-linked, then the adhesion is $ (d+1)$ and it admits a linkage of order $(d+1)$ between the bags. This ``off-by-one'' definition is the reason why the following lemma is stated as is.

\begin{figure}[ht]
 \centering
 \begin{tikzpicture}

 \pgfdeclarelayer{background}
		\pgfdeclarelayer{foreground}
			
		\pgfsetlayers{background,main,foreground}
			
 \begin{pgfonlayer}{main}
 \node (C) [v:ghost] {};

 \node(L) [v:ghost] at (-3.5,0) {
 \begin{tikzpicture}

 \pgfdeclarelayer{background}
		 \pgfdeclarelayer{foreground}
			
		 \pgfsetlayers{background,main,foreground}

 \begin{pgfonlayer}{background}
 \pgftext{\includegraphics[width=4.5cm]{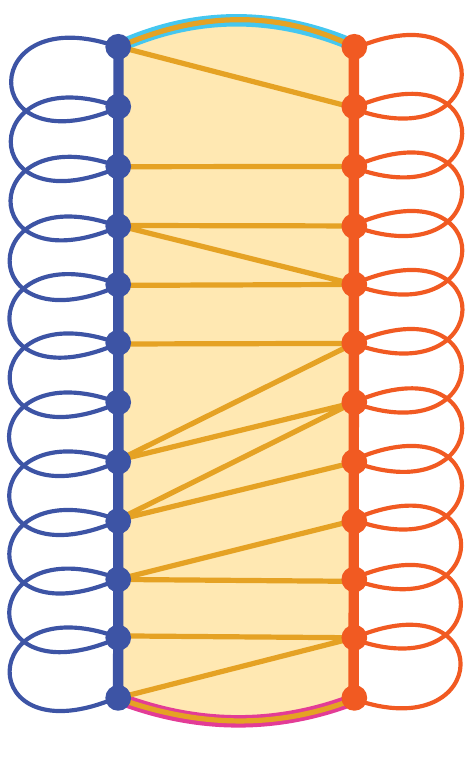}} at (C.center);
 \end{pgfonlayer}{background}
			
 \begin{pgfonlayer}{main}
 \node (C) [v:ghost] {};
 
 \end{pgfonlayer}{main}
 
 \begin{pgfonlayer}{foreground}
 \end{pgfonlayer}{foreground}

 \end{tikzpicture}
 };

 \node(M) [v:ghost] at (3.5,0) {
 \begin{tikzpicture}

 \pgfdeclarelayer{background}
		 \pgfdeclarelayer{foreground}
			
		 \pgfsetlayers{background,main,foreground}

 \begin{pgfonlayer}{background}
 \pgftext{\includegraphics[width=4.5cm]{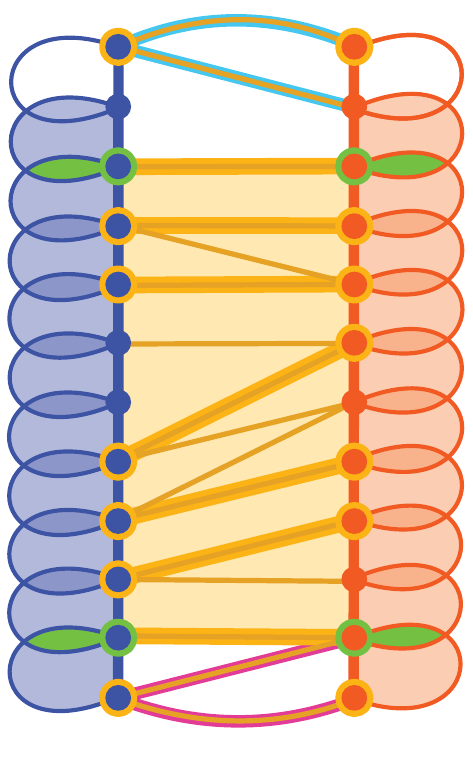}} at (C.center);
 \end{pgfonlayer}{background}
			
 \begin{pgfonlayer}{main}
 \node (C) [v:ghost] {};
 
 \end{pgfonlayer}{main}
 
 \begin{pgfonlayer}{foreground}
 \end{pgfonlayer}{foreground}

 \end{tikzpicture}
 };

 \node (i) [v:ghost] at (-3.5,-4.2) {(a) The segments and $\Delta$.};
 \node (ii) [v:ghost] at (3.5,-4.2) {(b) The graph $G^*$ and $T^*$.};

 \node (PuD) [v:ghost] at (-3.5,3.9) {$P^{\Delta}$};
 \node (PlD) [v:ghost] at (-3.5,-3.6) {$P_{\Delta}$};

 \node (UL) [v:ghost] at (-4.7,3.6) {$U_{\ell}$};

 \node (iL) [v:ghost] at (-4.25,3.6) {\small $i_{\ell}$};
 \node (iL-1) [v:ghost] at (-4.15,2.87) {\small $i_{\ell-1}$};
 \node (iL-2) [v:ghost] at (-4.15,2.25) {\small $i_{\ell-2}$};

 \node (i5) [v:ghost] at (-4.25,-1.55) {\small $i_{5}$};
 \node (i4) [v:ghost] at (-4.25,-2.02) {\small $i_{4}$};
 \node (i3) [v:ghost] at (-4.25,-2.35) {\small $i_{3}$};
 \node (i2) [v:ghost] at (-4.25,-2.82) {\small $i_{2}$};
 \node (i1) [v:ghost] at (-4.25,-3.4) {\small $i_{1}$};

 \node (WL) [v:ghost] at (-2.4,3.6) {$W_{\ell}$};

 \node (jL) [v:ghost] at (-2.8,3.2) {\small $j_{\ell}$};
 \node (jL-1) [v:ghost] at (-2.8,2.81) {\small $j_{\ell-1}$};
 \node (jL-2) [v:ghost] at (-2.8,2.25) {\small $j_{\ell-2}$};

 \node (j5) [v:ghost] at (-2.7,-1.3) {\small $j_{5}$};
 \node (j4) [v:ghost] at (-2.7,-1.9) {\small $j_{4}$};
 \node (j3) [v:ghost] at (-2.7,-2.25) {\small $j_{3}$};
 \node (j2) [v:ghost] at (-2.7,-2.7) {\small $j_{2}$};
 \node (j1) [v:ghost] at (-2.7,-3.4) {\small $j_{1}$};

 \node (xPuD) [v:ghost] at (3.5,3.9) {$P^{\Delta}$};
 \node (xPlD) [v:ghost] at (3.5,-3.6) {$P_{\Delta}$};

 \node (xWL) [v:ghost] at (4.7,3.6) {$W_{\ell}$};

 \node (xjw) [v:ghost] at (4.25,2.25) {\small $j_{w}$};

 \node (xj4) [v:ghost] at (4.35,0.25) {\small $j_{4}$};
 \node (xj3) [v:ghost] at (4.35,-0.75) {\small $j_{3}$};
 \node (xj2) [v:ghost] at (4.35,-1.55) {\small $j_{2}$};
 \node (xj1) [v:ghost] at (4.35,-2.5) {\small $j_{1}$};
 \node (xj0) [v:ghost] at (4.35,-3.4) {\small $j_{0}$};

 \node (xYL) [v:ghost] at (2.4,3.6) {$W_{\ell}$};

 \node (xiw) [v:ghost] at (2.8,2.25) {\small $i_{w}$};

 \node (xi4) [v:ghost] at (2.8,-0.5) {\small $i_{4}$};
 \node (xi3) [v:ghost] at (2.8,-1.1) {\small $i_{3}$};
 \node (xi2) [v:ghost] at (2.8,-1.75) {\small $i_{2}$};
 \node (xi1) [v:ghost] at (2.8,-2.5) {\small $i_{1}$};
 \node (xi0) [v:ghost] at (2.7,-3.3) {\small $i_{0}$};

 \end{pgfonlayer}{main}
 
 \begin{pgfonlayer}{foreground}
 \end{pgfonlayer}{foreground}

 \end{tikzpicture}
 \caption{A schematic illustration of the proof of \cref{lem:bounded_jumps_between_vortices}. The left hand figure highlights the setting of the proof, where the segments $\SSS$ are drawn in orange, $F_1$ is drawn in blue and $F_2$ is drawn in red with respective schematic representations of linear decompositions. The right hand figure depicts the remaining set $\SSS^*$ after sorting the segments, where $G^*$ is given by combining the orange, blue, and red highlighted regions, whereas $T^*$ is schematically highlighted in green. The set of paths $\tilde\PPP$ is highlighted in light blue and light purple. }
 \label{fig:lemma_vortexjumps}
\end{figure}
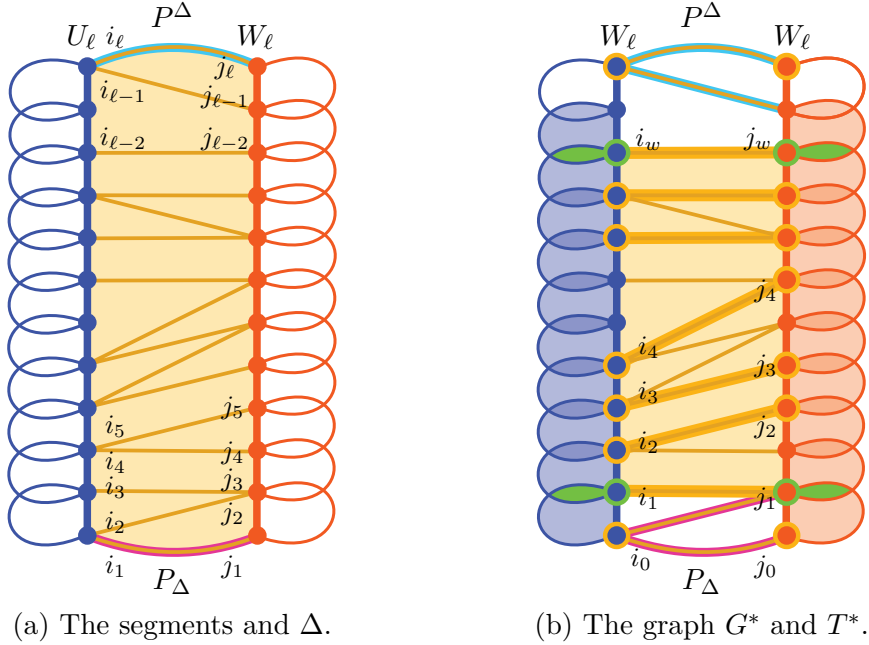

\begin{lemma}\label{lem:bounded_jumps_between_vortices}
 Let $\ell',k',d_1,d_2\geq 1$. Let $\Sigma$ be a surface. There exists a function $g_{\ref{lem:bounded_jumps_between_vortices}}:\N \to \N$ such that the following holds. Let $\rho$ be a blank and stretched $\Sigma$-rendition of some annotated graph~$(G,T)$. Let $c_1,c_2 \in C(\rho)$ be vortices (possibly equal). Let $\alpha_1=\left(\langle u_1,\ldots,u_{\ell'}\rangle,\langle U_1,\ldots,U_{\ell'}\rangle\right)$ be a linear decomposition of the vortex society $(\sigma(c_1),\Omega_{c_1})$ of adhesion at most $d_1,$ and let $\alpha_2=(\langle w_1,\ldots,w_{k'}\rangle,\langle W_1,\ldots,W_{k'}\rangle)$ be a linear decomposition of the vortex society $(\sigma(c_2),\Omega_{c_2})$ of adhesion at most $d_2$. Let $1 \leq \ell_1<\ell_2 \leq \ell'$ and $1 \leq k_1 < k_2 \leq k'$ be such that the $[\ell_1,\ell_2]$-restriction of $\alpha_1$ is $T$-normal and properly $(d_1-1)$-linked and the $[k_1,k_2]$-restriction of $\alpha_2$ is $T$-normal and properly $(d_2-1)$-linked. Further assume that for every $i \in [k_1,k_2]$ and every $j \in [\ell_1,\ell_2]$ it holds $U_i \cap W_j = \emptyset$. 
 
 Let $\LLL$ be a vital $T$-linkage such that $G=\bigcup \LLL \cup \sigma(c_1) \cup \sigma(c_2)$. Let $\SSS$ be a set of $\rho$-homotopic segments of $\LLL$ such that every $P \in \SSS$ has one endpoint in $\{u_{\ell_1},\ldots,u_{\ell_2}\}$ and the other in $\{w_{k_1},\ldots,w_{k_2}\}$. Then $\Abs{\SSS}< g_{\ref{lem:bounded_jumps_between_vortices}}(d_1+d_2)$.

 Furthermore, $g_{\ref{lem:bounded_jumps_between_vortices}}(d) \in 2^{\mathbf{poly}(d)}$.
\end{lemma}
\begin{proof}
Let $f \coloneqq f_{\ref{thm:Lemma2.5}}$ and let $d \coloneqq d_1 + d_2 $. Finally fix $g(d) \coloneqq 12f(d,d)+14;$ we claim that $g$ satisfies the theorem; $g(d) \in 2^{\mathbf{poly}(d)}$ by definition. Note that if $c_1=c_2$ we may still have $\Omega_{c_1} \neq \Omega_{c_2};$ then the second is the ``reverse'' order of the first.

Towards a contradiction, assume that for our choice of $g$ the theorem is wrong, and let $\SSS$ be as in the theorem with $\Abs{\SSS}\geq 12f(d,d)+14$. Let $\SSS(\LLL)$ be the set of all segments of $\LLL;$ in particular~$\SSS \subseteq \SSS(\LLL)$. 

Without loss of generality, by renaming the vertices accordingly with respect to the cyclic order~$\Omega_{c_i}$ for $i=1,2,$ respectively, we may assume that $\ell_1=1$ and let $\ell \coloneqq \ell_2$ as well as $k_1 = 1$ and let $k \coloneqq k_2$. 
Note that if $c_1=c_2$ then $d_1=d_2,$ but it may still be true that $\alpha_1 \neq \alpha_2,$ in particular the linear decompositions are allowed to have differing bags. Furthermore, note that since $\rho$ is stretched, the additional assumption $\{u_i \mid i \in [\ell_1,\ell_2]\} \cap \{w_i \mid i \in [k_1,k_2]\} = \emptyset$ is only relevant if $c_1=c_2$. It turns out that the additional assumption makes the proof in the single vortex case analogous to the one with two vortices. This is intuitively easy to see, since there is no tangible way to differentiate between both: If the graph $\sigma(c)$ were to consist of two disjoint components one could easily split the single vortex into two vortices. 

In addition, without loss of generality, we may assume that $u_1,w_1,u_\ell,w_k$ are endpoints of some paths in $\SSS$ since otherwise we may restrict the segments of $\Omega_{c_1}$ and $\Omega_{c_2}$ accordingly. Let $U \coloneqq \{u_1,\ldots,u_\ell\}$ and $W \coloneqq \{w_1,\ldots,w_k\}$.

Since $\LLL$ is a linkage, $\SSS$ is a set of pairwise internally disjoint paths. And since the paths in $\SSS$ are grounded and pairwise $\rho$-homotopic, there exists a $\rho$-aligned vortex-free disc $\Delta\subseteq \Sigma$ that intersects $N(\rho)$ exactly in all of $U \cup W$ such that for every $P \in \SSS$ its trace is contained in $\Delta$. This follows from \cref{obs:strip_from_homotopic_paths}, but we provide an exemplary construction for completion. To see this note that the traces of any two $\rho$-homotopic paths $P,P' \in \SSS$ bound a disc $\Delta_{P,P'}$ together with the respective segments of vortex boundaries between them by definition of homotopy, call these segments $F_1,F_2$ for $c_1,c_2$ respectively. That is, given the traces $\gamma_{P},\gamma_{P'}$ of $P,P',$ there exists an $(F_1,F_2)$-strip $\Delta_{P,P'}$ guarded by $\gamma_P,\gamma_{P'};$ in particular that disc is $\rho$-aligned by construction. Call two paths $P,P' \in \PPP$ \emph{consecutive} if the trace of no other path $P \in \SSS$ is contained in $\Delta_{P,P'}$. Note here that there may still be $P^* \in \SSS(\LLL)\setminus \SSS$ such that the trace of $P$ is contained in $\Delta_{P,P'}$ but then $P^*$ needs to be a simple loop with both ends either in $U$ or both in $W;$ if $c_1=c_2$ then $P$ and $P'$ may also be (simple) loops, but they still have endpoints in the different segments. Finally, one can construct $\Delta$ by gluing the discs $\Delta_{P,P'}$ of pairwise consecutive paths carefully~--~possibly omitting one such disc if the resulting surface would be no disc~--~until the resulting disc contains the traces of all paths in $\SSS$.

Note that since $\Delta$ is a disc and all the paths in $\SSS$ have their traces in $\Delta$ with both endpoints on its boundary, we derive that~--~after possibly reversing the order of $\langle w_1,\ldots,w_k \rangle$ and renaming~--~for any two paths $P_1,P_2 \in \SSS$ with ends $u_{i_1},w_{j_1}$ and $u_{i_2},w_{j_2}$ we have 
\begin{itemize}
 \item[$(\star)$] if $1 \leq i_1 \leq i_2 \leq \ell$ then $1 \leq j_1 \leq j_2 \leq k$ and vice versa.
\end{itemize}
Let $P_\Delta, P^\Delta$ be the paths in $\SSS$ with endpoints $u_1,w_1$ and $u_\ell,w_k,$ respectively, and let $\tilde \PPP\subsetneq \SSS$ be the at most six paths with endpoints in $u_1,w_1,u_\ell,w_k$ including $P^\Delta,P_\Delta$\footnote{Due to topological restrictions these are actually at most $4$ paths.}. Let $\SSS' \coloneqq \SSS\setminus \tilde \PPP,$ in particular $\Abs{\SSS'} \geq 12 f(d,d)+8,$ then for every path in $\SSS'$ its trace is contained in the interior of $\Delta,$ in particular $\bigcup \SSS'$ is contained in the crop of $G$ by $\Delta;$ see the left hand figure in \cref{fig:lemma_vortexjumps}, where the paths $\tilde{\PPP}$ are highlighted by \textcolor{mediumslateblue}{purple} and \textcolor{ceruleanblue}{blue} tubes in the right hand figure. 

Since $\LLL$ is a linkage and every path in $\SSS$ is a segment of some path in $\LLL,$ an endpoint of any segment in $\SSS$ is part of at most $2$ paths in $\SSS$. By the pigeonhole principle, there is a set $\SSS'' \subseteq \SSS'$ of size $\geq 6f(d,d) + 4$ such that no two paths in $\SSS''$ end in the same vertex of $U$. Applying the pigeonhole principle again, we derive that there exists a set $\SSS^*\subseteq \SSS''$ of size $\geq 3f(d,d)+2$ such that $\SSS^*$ is a $U$-$W$-linkage; see the paths highlighted by an \textcolor{orangepeel}{orange} tube in the right hand \cref{fig:lemma_vortexjumps}. 

Without loss of generality, assume that $\ell \leq k$: the other case is analogous. Let $1\leq i_0< \ldots < i_{\Abs{\SSS^*}-1}\leq \ell$ be such that $U^*=\{u_{i_j}\mid 0\leq j \leq \Abs{\SSS^*}-1\}\subseteq U$ is the set of vertices that are endpoints of paths in $\SSS^*$. Let $1\leq j_0<\ldots <j_{\Abs{\SSS^*}-1}\leq k$ be such that $W^*=\{w_{j_i}\mid 0\leq i \leq \Abs{\SSS^*}-1\}\subseteq W$ is the set of vertices that are endpoints of paths in $\SSS^*$. Since $\SSS^*$ is a $U$-$W$-linkage, $(\star)$ implies that for every $P \in \SSS$ there is $p \in [0,\Abs{\SSS^*}-1]$ such that $P$ has $u_{i_p},w_{i_p}$ as its endpoints. 
Fix $\omega = 3f(d,d)$ and let $\{P_0,P_1\ldots,P_\omega,P_{\omega+1}\} \subseteq \SSS^*$ be such that $P_p$ has $u_{i_p},w_{i_p}$ as its endpoints which exists since $\Abs{\SSS^*} \geq \omega+2;$ let $\PPP \coloneqq \{P_1,\ldots,P_\omega\}$. For every $i \in [\omega]$ let $\gamma_i$ be the trace of~$P_i$. 
Let $\Delta^*$ be the closure of the component of $\Delta$ after deleting the trace of $P_0$ and $P_{\omega+1}$ that contains the traces of $P_1,\ldots,P_\omega$. By construction $\Delta^*$ is $\rho$-aligned; we derive the following.
\begin{beautifulclaim}\label{lem:bounded_jumps_between_vortices_claim_Delta_i}
 Let $i \in [\omega]$ and let $\gamma_i$ be the trace of~$P_i$. Then $\gamma_i \subseteq \Delta^*,$ it has both its endpoints on the boundary of $\Delta^*$ and is otherwise disjoint from the boundary of $\Delta^*$. The traces of $P_0$ and~$P_{\omega+1}$ are both contained in the boundary of $\Delta^*$ and the crop $G_{\Delta^*}$ of $G$ by $\Delta^*$ contains all paths~$P_1,\ldots,P_{\omega}$.
\end{beautifulclaim}

Next let $G' \coloneqq G_{\Delta^*} \cup \bigcup_{j=0}^{i_{\omega+1}}U_j \cup \bigcup_{i=0}^{j_{\omega+1}}W_i$ and finally let $G^*$ be obtained from $G'$ by deleting all the internal vertices of $P_0$ and $P_{\omega+1};$ see the right hand \cref{fig:lemma_vortexjumps}. 

Note that by our choice of $\SSS^* \subseteq \SSS'$ we have $\tilde \PPP \cap \SSS^* = \emptyset,$ in particular $i_0-1 \geq 1$ and $i_{\omega+1}+1 \leq k$ as well as $j_0-1 \geq 1$ and $j_{\omega+1}+1 \leq \ell$. Therefore, the following is well-defined
$$ T^* \coloneqq (U_{i_0}\cap U_{i_0-1}) \cup (U_{i_{\omega+1}}\cap U_{i_{\omega+1}+1}) \cup (W_{j_0}\cap W_{j_0-1}) \cup (W_{j_{\omega+1}}\cap W_{j_{\omega+1}+1})$$
and satisfies $\Abs{T^*} \leq 2d_1 + 2d_2 = 2d;$ see the \textcolor{myGreen}{green} highlighted regions in the right hand \cref{fig:lemma_vortexjumps}. We have the following.

\begin{beautifulclaim}\label{lem:bounded_jumps_between_vortices_claim_separation_G*}
 Let $\bar V \coloneqq (V(G)\setminus V(G^*)) \cup T^*$. Then $(V(G^*),\bar V)$ is a separation of $G$ of order $2d$ satisfying $V(G^*)\cap \bar V = T^*$.
\end{beautifulclaim}
\begin{claimproof}
 By definition $V(G^*) \cup \bar V = V(G)$ and $V(G^*)\cap \bar V = T^*$. As derived above $\Abs{T^*} =2d$. We are left to prove that it is indeed a separation of $G^*$. To this extent, let $uv \in E(G)$ be an edge: We need to prove that either $\{u,v\} \subseteq V(G^*) \setminus \bar V$ or $\{u,v\} \subseteq \bar V \setminus V(G^*)$ or $\{u,v\} \cap T^* \neq \emptyset$. 
 
 Since $G = \bigcup\LLL \cup \sigma(c_1) \cup \sigma(c_2)$ and, since either $V(\sigma(c_1)) \cap V(\sigma(c_2)) = \emptyset,$ for $\rho$ is stretched, or $V(\sigma(c_1)) = V(\sigma(c_2))$ if $c_1=c_2,$ we derive that either $\{u,v\} \subseteq V(\sigma(c_1))$ or $\{u,v\} \subseteq V(\sigma(c_2))$ or $e \in E(L)$ for some $L \in \LLL$. Let $i \in [i_0,i_{\omega+1}]$ and $i' \in [k'] \setminus [i_0,i_{\omega+1}]$ as well as $j \in [j_0,j_{\omega+1}]$ and $j' \in [\ell'] \setminus [j_0,j_{\omega+1}]$. Since $\alpha_1$ and $\alpha_2$ are linear decompositions, we derive that if $u \in U_i$ and $v \in U_{i'}$ or $u \in W_j$ and $v \in W_{j'},$ then $\{u,v\} \cap T^* \neq \emptyset$ which we wanted to prove. In particular we may assume that $\{u,v\} \not \subseteq V(\sigma(c_i))$ for $i=1,2$.

 Thus, let $e \in E(L),$ hence $e\in E(P)$ for some $P \in \SSS(\LLL),$ or $\{u,v\} \subseteq V(\sigma(c_i))$ for some $i=1,2;$ the latter we already excluded, therefore, let $e \in E(P)$. Let $C_{\Delta^*} \subseteq C(\rho)$ denote all the cells that are contained in $\Delta$ and let $C_{\bar\Delta^*}\subseteq C(\rho)$ be all the cells that are contained in $\mathsf{cl}(\Sigma \setminus \Delta^*);$ by definition of $\Delta^*,$ the sets $C_{\Delta^*},C_{\bar\Delta^*}$ form a partition of $C(\rho)$. 
 Let now $c_e \in C(\rho)$ be such that $e \in E(\sigma(c_e)),$ in particular $u,v \in V(\sigma(c_e))$. 
 By the above analysis we know that $c_e \notin \{c_1,c_2\}$. If none of $u,v$ are nodes of $c_e,$ then, since $C_{\Delta^*},C_{\bar\Delta^*}$ partitions $C(\rho)$ and since $c_e \notin \{c_1,c_2\},$ we either derive $\{u,v\} \subseteq V(G^*)\setminus \bar V$ or $\{u,v\} \subseteq \bar V \setminus V(G^*);$ in both cases we are done.
 Thus at least one of $u,v$ is a node of $c_e$. Furthermore, note that $c_e$ must have at least one node on the boundary of $\Delta^*$ since, otherwise, either $\{u,v\} \subseteq V(G^*)\setminus \bar V$ or $\{u,v\} \subseteq \bar V \setminus V(G^*)$ by the definition of crop by $\Delta^*;$ in particular, we are done.
 
 We are left to prove that if $e=uv\in E(P)$ where $u$ or $v$ lie on the boundary of $\Delta^*$ with $u \in V(G^*),$ say, and $w \in \bar V$ say, then $\{u,w\} \cap T^* \neq \emptyset$. To see this, note that by the construction of $\Delta^*$ and the fact that $\SSS$ is a $U$-$W$-linkage we derive that either $P\in \{P_0,P_{\omega+1}\},$ or its trace $\gamma_P$ is internally disjoint from the traces of $P_0,P_{\omega+1}$ and satisfies $\gamma_P \subseteq\Delta^*$ or $\gamma_P \subseteq\mathsf{cl}(\Sigma\setminus \Delta^*)$. 
 
 If $P = P_0,$ then by construction the only vertices of $V(G^*)$ in common with $P_0$ are $u_{i_0},w_{j_0} \in T^*;$ an analogous argument works for $P= P_{\omega +1}$. Hence we may assume $P \notin \{P_0,P_{\omega+1}\}$ which implies that either $\gamma_u \subseteq\Delta^*$ or $\gamma_P \subseteq\Sigma\setminus \Delta^*$. 
 
 Assume for now that $\gamma_P \subseteq\Delta^*$. Then again, since $\SSS(\LLL)$ is a set of internally disjoint paths by \zcref{obs:segments} and since it has its ends in $U \cup W,$ we derive that $\gamma_P$ is disjoint from the traces of all others paths in $\SSS$ except for possibly their endpoints, i.\@e.\@, except for possibly nodes in $$\{u_{i_0},u_{i_0+1},\ldots,u_{i_{\omega+1}-1},u_{i_\omega+1}\} \cup \{w_{j_0},w_{j_0+1},\ldots,w_{j_{\omega+1}-1},w_{j_\omega+1}\}.$$ Thus $\gamma_P$ is disjoint from $\Sigma \setminus \Delta^*$ except for possibly nodes in $u_{i_0},u_{i_{\omega+1}},w_{j_0},w_{j_{\omega+1}}$ since the remaining nodes are not part of $\bar V^*$. Thus, we derive that $\{u,v\} \cap T^* \neq \emptyset$ as desired. An analogous argument works if $\gamma_P \subseteq\Sigma \setminus \Delta^*$ since its trace is disjoint from the traces of $P_0$ and $P_{\omega+1}$ we derive that it must be disjoint from all of the interior of $\Delta^*$ and hence from $\Delta^*$~--~its boundary consists of segments of the vortex cell and the trace of $P_0$ and $P_{\omega+1}$~--~except for possibly nodes in $$\{u_{i_0},u_{i_0+1},\ldots,u_{i_{\omega+1}-1},u_{i_\omega+1}\} \cup \{w_{j_0},w_{j_0+1},\ldots,w_{j_{\omega+1}-1},w_{j_\omega+1}\}$$ on its boundary. Then again we derive that $\gamma_P$ is disjoint from $\Delta^*$ except for possibly nodes in $u_{i_0},u_{i_{\omega+1}},w_{j_0},w_{j_{\omega+1}} \subseteq T^*$ as desired, concluding the proof.
\end{claimproof}

We derive the following. 

\begin{beautifulclaim}\label{lem:bounded_jumps_between_vortices_claim_G^*_annotated}
 $(G^*,T^*)$ is an annotated graph admitting a vital $T^*$-linkage $\LLL^*$ of order $\leq d$ such that every path in $\LLL^*$ is a subpath of a path in $\LLL$.
\end{beautifulclaim}
\begin{claimproof}
 Note that by the assumption of the theorem that $\alpha_1$ and $\alpha_2$ are $T$-normal, we derive that every $v \in T \cap (\bigcup_{p=1}^{k'}U_p \cup \bigcup_{q=1}^{\ell'}W_q)$ satisfies $v \in T^*$. Thus the claim follows at once by combining \zcref{lem:bounded_jumps_between_vortices_claim_separation_G*} with \zcref{lem:vital_and_separations}, noting that $\Abs{T^*} \leq 2d$ is easily seen to imply that $\LLL^*$ is of order at most $d$: Note that $L \in \LLL^*$ may consist of a single terminal, but then $\Abs{T^*} < 2d$. 
\end{claimproof}

Next we define a sequence of laminar separations of equal order in $G^*;$ recall that $\omega = 3f(d,d)$. Towards this, note that by \zcref{lem:bounded_jumps_between_vortices_claim_Delta_i} for every $i \in [\omega]$ and trace $\gamma_i$ of $P_i \in \PPP$ we derive that $\Delta^* \setminus \gamma_i$ consists of exactly two non-empty components, one containing the trace of $P_0$ and the other that of $P_{\omega+1}$. With this in mind, we define the following; let $p \in [\omega]$.
\begin{itemize}
 \item[$\Delta_p$:] Let $\Delta_p \subseteq\Delta^*$ be the closure of the component of $\Delta^*\setminus \gamma_p$ that contains the traces of $P_0,P_1,\ldots,P_p$ and let $\bar\Delta_p$ be the closure of the other component. If $p=\omega+1$ fix~$\Delta_{\omega+1}=\Delta^*$. By construction $\Delta_p$ is $\rho$-aligned.
 \item[$G_p$:] Let $G_{\Delta_p}$ be the crop of $G$ by $\Delta_p$ in $\rho$ and let $G_p' \coloneqq G_{\Delta_p} \cup \bigcup_{j=0}^{i_{p}}U_j \cup \bigcup_{i=0}^{j_{p}}W_i$ and finally let $G_p$ be obtained from $G_p'$ by adding all the edges and vertices of $P_p$ if they were not already part of the crop. 
\end{itemize}

For every $1\leq 3t \leq \omega$ let 
 \begin{align*}
 A_t&\coloneqq V(G_{i_{3t}})\\
 B_t&\coloneqq (V(G^*) \setminus V(G_{i_{3t}})) \cup U_{i_{3t}+1} \cup W_{j_{3t}+1}.
 \end{align*} 
 We have the following.
\begin{beautifulclaim}\label{lem:bounded_jumps_between_vortices_claim_laminar_family}
 $\langle(A_t,B_t)\rangle_{1 \leq t \leq f(d,d)}$ is a sequence of laminar separations of $G^*,$ each of order $d$ and for each $1< t \leq f(d,d)$ there is a linkage $M_t$ in $G^*[B_{t-1}]\cap G^*[A_{t}]$ of order $d$ that is edge-disjoint from $\PPP$ and such that each path in $M_t$ has one end in $A_{t-1}\cap B_{t-1}$ and the other in $A_{t}\cap B_{t}$.
 \end{beautifulclaim}
 \begin{claimproof}
 Notice that for every $1\leq t \leq f(d,d),$ $V(A_t)\cup V(B_t) = V(G^*)$. By construction we derive that
 \begin{equation}\label{lem:bounded_jumps_between_vortices_eq}
 A_t\cap B_t = (U_{i_{3t}}\cap U_{i_{3t+1}}) \cup (W_{j_{3t}}\cap W_{j_{3t+1}}).
 \end{equation}

 Using similar arguments as in \zcref{lem:bounded_jumps_between_vortices_claim_separation_G*} one easily verifies that $(A_t,B_t)$ is indeed a separation of $G^*$: note again that since $G = \LLL \cup \sigma(c_1) \cup \sigma(c_2),$ every edge adjacent to a vertex in some path in $\SSS(\LLL)$ must be contained in said path. Thus, any edge with one end in $A_t$ and the other in $B_t$ that is not contained in one of the bags of the linear decompositions, must be contained in a cell with at least one node on the boundary of $\Delta_i$. In particular, it must be part of a path in $\SSS(\LLL)$ whose trace intersects the boundary of $\Delta_i$. Therefore it is either contained in $E(P_i)$ or in $E(P_0)$ or it is part of some path that intersects them at their endpoints~--~$\SSS(\LLL)$ is a set of internally disjoint paths~--~which are in $A_t \cap B_t$ by construction.

 Recall that by assumption of the theorem we have that for every $i \in [k]$ and $j \in [\ell]$ it holds that $U_i \cap W_j= \emptyset$. Further, by assumption of the theorem the, respective $[k]$- and $[\ell]$-restrictions of the linear decompositions are properly $(d_1-1)$- and $(d_2-1)$-linked. Combining both we derive that $\Abs{A_t\cap B_t} = \Abs{U_{i_{3t}}\cap U_{i_{3t+1}} \cup W_{j_{3t}}\cap W_{j_{3t+1}}} = \Abs{U_{i_{3t}}\cap U_{i_{3t+1}}} + \Abs{W_{j_{3t}}\cap W_{j_{3t+1}}} = d_1+d_2 =d$.
 
 Since $(\langle u_1,\ldots,u_\ell \rangle,\langle U_1,\ldots,U_\ell \rangle)$ is properly $(d_1-1)$-linked, the \cref{def:linked_vortex_normal} of properly linked segments of linear decompositions implies that for every $1< p <p' < \ell$ there is a $d_1$-linkage -- recall the ``off-by-one'' error -- $\LLL^U_{p,p'} \subseteq G[U_p\cup \ldots \cup U_{p'}]$ such that each path in $\LLL^U_{p,p'}$ has exactly one end in $U_{p-1}\cap U_p$ and the other in $U_{p'} \cap U_{p'+1}$. Analogously, for the properly $(d_2-1)$-linked linear decomposition $(\langle w_1,\ldots,w_k\rangle,\langle W_1,\ldots,W_k\rangle)$ and every $1 < q < q'< k$ there is a $d_2$-linkage $\LLL^W_{q,q'}\subseteq G[W_q\cup \ldots \cup W_{q'}]$ such that each path in $\LLL^W_{q,q'}$ has exactly one end in $W_{q-1}\cap W_q$ and the other in $W_{q'} \cap W_{q'+1}$. Again since for every $i \in [k]$ and $j \in [\ell]$ it holds that $U_i \cap W_j= \emptyset,$ both linkages are disjoint whence $\LLL^{q,q'}_{p,p'}\coloneqq \LLL^U_{p,p'} \cup \LLL^W_{q,q'}$ is a $d$-linkage in $G^*$.
 
 Finally, let $1< t \leq f(d,d)$. By definition of $(A_{t-1},B_{t-1}),(A_{t},B_{t})$ and the above analysis, we derive that $U_{i_{3(t-1)+1}},U_{i_{3(t-1)+2}},U_{i_{3t}} \subseteq B_{t-1}\cap A_{t}$ as well as $W_{j_{3(t-1)+1}},W_{j_{3(t-1)+2}},W_{j_{3t}} \subseteq B_{t-1}\cap A_{t}$. Thus we derive that $M_t\coloneqq \LLL^{j_{3(t-1)+1},j_{3t}}_{i_{3(t-1)+1},i_{3t}} \subseteq G[B_{t-1}] \cap G[A_{t}]$. Now $M_t$ is by the above observation a $d$-linkage such that each path has exactly one end in $(U_{i_{3(t-1)+1}}\cap U_{i_{3(t-1)+2}}) \cup (W_{j_{3(t-1)+1}}\cap W_{j_{3(t-1)+2}})$ and the other in $(U_{i_{3(t-1)+2}}\cap U_{i_{3t}}) \cup (W_{j_{3(t-1)+2}}\cap W_{j_{3t}})$. Thus, using \zcref{lem:bounded_jumps_between_vortices_eq}, we derive that $M_t$ has one end in $A_{t-1}\cap B_{t-1}$ and the other in $A_{t+1}\cap B_{t+1}$ as desired. Note that it has no edge in $E(\bigcup \SSS(\LLL))\cap E(G^*)$ and in particular no edge in $E(\bigcup \PPP)$ by construction.
 \end{claimproof}

 Let $\langle (A_t,B_t)\rangle_{1 \leq t \leq f(d,d)}$ be as in \zcref{lem:bounded_jumps_between_vortices_claim_laminar_family}. We have the following.

 \begin{beautifulclaim}\label{lem:bounded_jumps_between_vortices_claim2}
 Let $1\leq t < f(d,d)$. Then $V(P_{i_{3t}-1})$ is contained in $B_{t-1}\cap A_t$.
 \end{beautifulclaim}
 \begin{claimproof}
 This is immediate by construction of the sequence $\langle(A_t,B_t)\rangle_{1 \leq t \leq f(d,d)}$. In particular, since $G[A_t]$ contains the crop of $G^*$ by $\Delta_{i_{3t}},$ it contains $P_{i_{3t}-1}$ and since $\Delta_{i_{3(t-1)}}$ is disjoint from the trace of $P_{i_{3t}-1}$ we derive that $V(P_{i_{3t}-1}) \subseteq B_t$.
 \end{claimproof}

 By \zcref{lem:bounded_jumps_between_vortices_claim_G^*_annotated} $\LLL^*$ is a vital $T^*$-linkage of order $\leq d$ in $G^*$. By our choice of $f(d,d) = f_{\ref{thm:Lemma2.5}}(d,d),$ and \zcref{lem:bounded_jumps_between_vortices_claim_laminar_family}, \zcref{thm:Lemma2.5} implies that there is $1 \leq i<\omega$ such that $\bigcup\LLL^* \cap G^*[B_i] \cap G^*[A_{i+1}] = \bigcup M_i$. But $\LLL^*$ satisfies $P_{i_{3t}-1} \subseteq L$ for some $L \in \LLL^*,$ whereas $M_i$ is edge-disjoint from $\PPP;$ this is a contradiction concluding the proof.
\end{proof}

It is a small but crucial observation that the function $f_{\ref{lem:bounded_jumps_between_vortices}}$ does \emph{not} depend on the size $\Abs{T}$ but solely on the depth of the vortices.

\smallskip
Our main goal is to generalise \cref{lem:bounded_jumps_between_vortices} to a more general setting. The first step in that direction is an easy combinatorial result that guarantees the existence of $T$-normal segments in linear decompositions of sufficient length.

\begin{lemma}\label{lem:T-normal_segments_exist}
 Let $\ell^*,\ell,k\geq 1$ and $t \geq 2$. Let $\Sigma$ be a surface. Let $\rho$ be a $\Sigma$-rendition of some $k$-annotated graph $(G,T)$. Let $c_1 \in C(\rho)$ be a vortex. Let $\alpha=(\langle u_1,\ldots,u_{\ell^*}\rangle,\langle U_1,\ldots,U_{\ell^*}\rangle)$ be a linear decomposition of $c$. Let $J \subseteq[\ell^*]$ be a set of indices such that $\Abs{J} =\ell$. If $\ell > (2k+1)(t+1),$ then there exist $ j_1,j_2 \in J$ with $\Abs{[j_1,j_2]\cap J}\geq t,$ and such that the $[{j_1},{j_2}]$-segment of $\alpha$ is $T$-normal.
\end{lemma}
\begin{proof}
 Let $T= \{t_1,\ldots,t_k\}$ be non-empty for some $k \geq 1$. For every $t \in [k]$ let $1\leq s_t \leq f_t \leq \ell^*$ be such that $t \in \bigcap_{p=s_t}^{f_t} U_p$ and $t\notin U_q$ for all $q \notin [s_t,f_t];$ this exists since $\alpha$ is a linear decomposition. Let $X = \{s_t,f_t \mid t \in [k]\}$. Then $\Abs{X} \leq 2k$. Let $X=\{x_1,\ldots,x_m\}$ for some $ m \in [2k]$ such that $x_1 \leq \ldots,\leq x_m$. We derive the existence of $(m+1)$ intervals $I_0,\ldots,I_{m}$ of $[\ell^*]$ satisfying the following
 \begin{enumerate}
 \item $\Abs{I_0 \cap X} = 1$ as well as $\Abs{I_m \cap X}=1$ and $I_0 = [0,x_1]$ (where $x_1=0$ is possible) as well as $I_m = [x_m,\ell]$ (where $x_m=\ell$ is possible)
 \item for every $i \in [1,m-1]$ it holds $\Abs{I_i \cap X} = 2$ and $I_i=[x_i,x_{i+1}]$.
 \end{enumerate}
 In particular, note that for $0 \leq i < j \leq m $ and every $x \in I_i$ and $y \in I_j$ it holds $x \leq y$ and if $I_i \cap I_j \neq \emptyset,$ then $j=i+1$ and $I_i \cap I_{i+1} = \{x_{i+1}\}$.
 
 Note that $\ell > (2k+1)(t+1)$ implies $\ell>(m+1)(t+1)$. By the pigeonhole principle, there is $p \in [0,m]$ such that $I_p$ is an interval of length at least $(t+2)$ containing at least $(t+2)$ elements of $J$. In particular since $I_p=[x_p,x_{p+1}]$ we derive the $x_{p+1}\geq x_p+t+2$. Hence $I \coloneqq [x_{p}+1,x_{p+1}-1]$ is a non-empty interval containing at least $t$ elements of $J$. Thus, there exist $x_p+1 \leq j_1<j_2\leq x_{p+1}$ with and $j_1,j_2 \in J$ and $\Abs{[j_1,j_2]\cap J}\geq t$.
 \begin{beautifulclaim}
 Let $i \in I$ and $t \in T \cap U_i,$ then $t \in \bigcap_{j \in I}U_j$.
 \end{beautifulclaim}
 \begin{claimproof}
 First note that by $i)$ and $ii),$ we derive $I \cap X = \emptyset$. In particular $s_t,f_t \notin I$. By construction of $s_t,f_t$ we derive that $t \in U_j$ if and only if $j \in [s_t,f_t]$. Since $t \in U_i$ we derive that $s_t \leq i \leq f_t$ and in particular since $I\cap X = \emptyset$ we derive $s_t \leq x_{p}+1\leq x_{p+1}-1 \leq f_t$ and finally $t \in U_j$ for all $j \in I$. This concludes the proof.
 \end{claimproof}

 By the above claim, the $[j_1,j_2]$-segment of $\alpha$ is $T$-normal, concluding the proof.
\end{proof}

As a second step, in order to generalise \cref{lem:bounded_jumps_between_vortices}, we want to guarantee that the respective $[k_1,k_2]$-segment and $[\ell_1,\ell_2]$-segment, as given by the lemma, are indeed disjoint. This is a given if $\rho$ is stretched and $c_1 \neq c_2,$ but it is not if $c_1 = c_2$. As briefly discussed in the proof of \cref{lem:bounded_jumps_between_vortices}, this is the only relevant distinction between these two cases: Note that if $c_1 = c_2,$ then there may be a single vertex that is part of all the bags of its linear decomposition.

We define the following.
\begin{definition}
 Let $\ell \geq 1$. Let $\Sigma$ be a surface. Let $\rho$ be a $\Sigma$-rendition of some graph $G$. Let $c \in C(\rho)$ be a vortex. Let $\alpha=(\langle u_1,\ldots,u_{\ell}\rangle,\langle U_1,\ldots,U_{\ell}\rangle )$ be a linear decomposition of the vortex society $(\sigma(c),\Omega_{c})$ (possibly with some gap). Let $F = [i,j]$ for some $1 \leq i \leq j$. We define~$G_\alpha[F] \coloneqq \bigcup_{u_i \in F}U_i,$ omitting the subscript $\alpha$ if clear from context.
\end{definition}

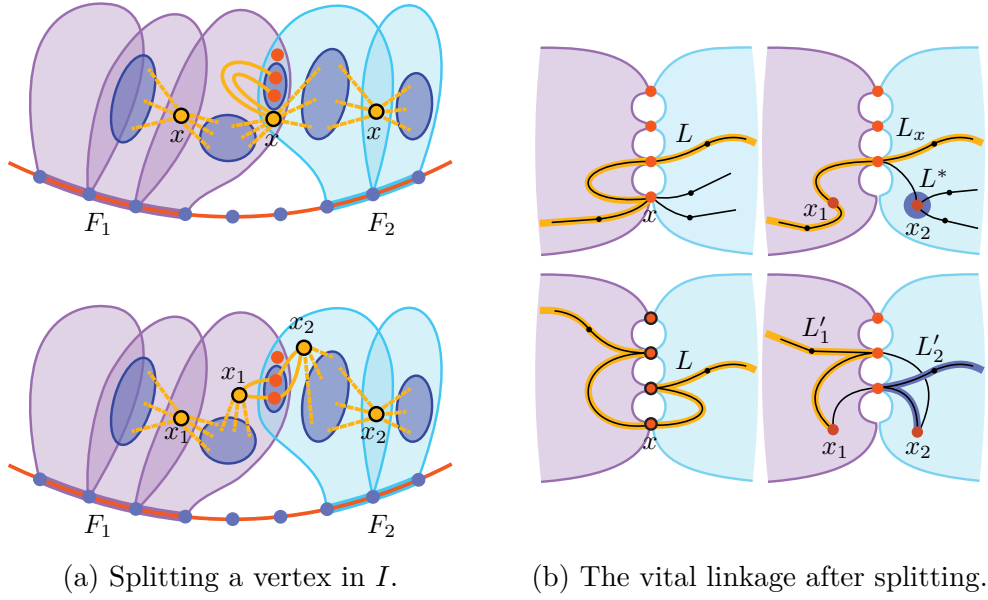
\begin{figure}[ht]
 \centering
 \begin{tikzpicture}

 \pgfdeclarelayer{background}
		\pgfdeclarelayer{foreground}
			
		\pgfsetlayers{background,main,foreground}
			
 \begin{pgfonlayer}{main}
 \node (C) [v:ghost] {};

 \node(L) [v:ghost] at (-3.5,2) {
 \begin{tikzpicture}

 \pgfdeclarelayer{background}
		 \pgfdeclarelayer{foreground}
			
		 \pgfsetlayers{background,main,foreground}

 \begin{pgfonlayer}{background}
 \pgftext{\includegraphics[width=6cm]{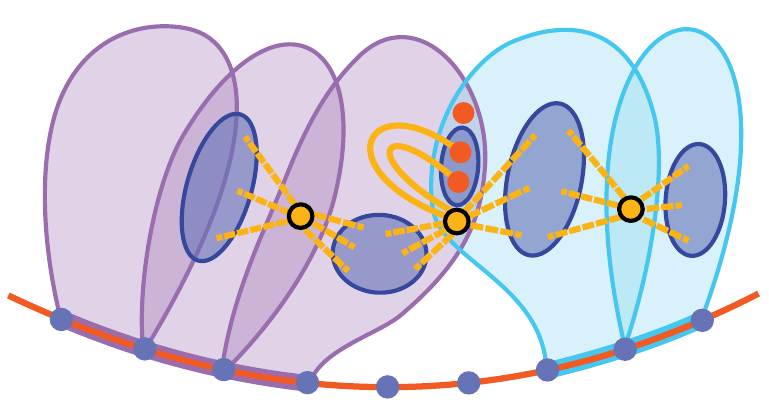}} at (C.center);
 \end{pgfonlayer}{background}
			
 \begin{pgfonlayer}{main}
 \node (C) [v:ghost] {};
 
 \end{pgfonlayer}{main}
 
 \begin{pgfonlayer}{foreground}
 \end{pgfonlayer}{foreground}

 \end{tikzpicture}
 };

 \node(M) [v:ghost] at (-3.5,-2) {
 \begin{tikzpicture}

 \pgfdeclarelayer{background}
		 \pgfdeclarelayer{foreground}
			
		 \pgfsetlayers{background,main,foreground}

 \begin{pgfonlayer}{background}
 \pgftext{\includegraphics[width=6cm]{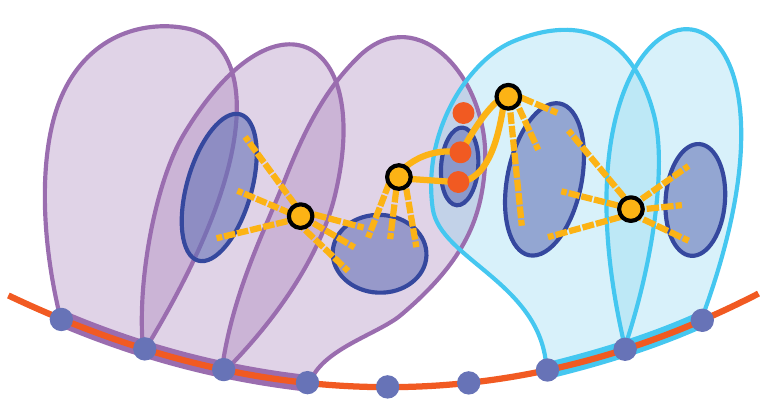}} at (C.center);
 \end{pgfonlayer}{background}
			
 \begin{pgfonlayer}{main}
 \node (C) [v:ghost] {};
 
 \end{pgfonlayer}{main}
 
 \begin{pgfonlayer}{foreground}
 \end{pgfonlayer}{foreground}

 \end{tikzpicture}
 };

 \node(M) [v:ghost] at (2,1.5) {
 \begin{tikzpicture}

 \pgfdeclarelayer{background}
		 \pgfdeclarelayer{foreground}
			
		 \pgfsetlayers{background,main,foreground}

 \begin{pgfonlayer}{background}
 \pgftext{\includegraphics[width=3cm]{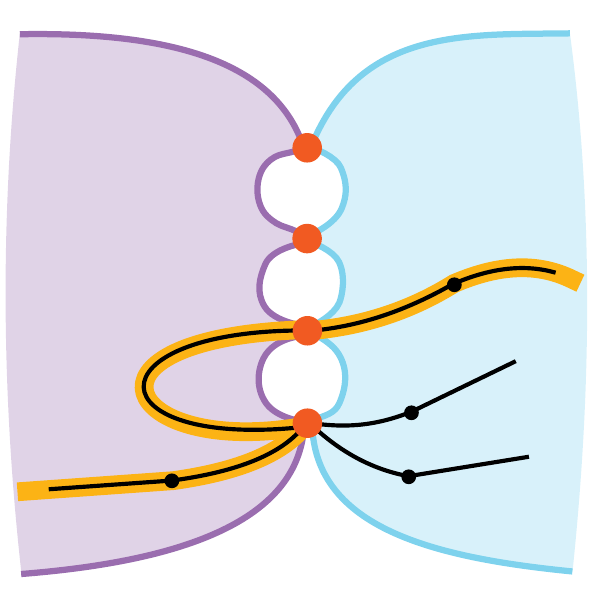}} at (C.center);
 \end{pgfonlayer}{background}
			
 \begin{pgfonlayer}{main}
 \node (C) [v:ghost] {};
 
 \end{pgfonlayer}{main}
 
 \begin{pgfonlayer}{foreground}
 \end{pgfonlayer}{foreground}

 \end{tikzpicture}
 };

 \node(M) [v:ghost] at (5,1.5) {
 \begin{tikzpicture}

 \pgfdeclarelayer{background}
		 \pgfdeclarelayer{foreground}
			
		 \pgfsetlayers{background,main,foreground}

 \begin{pgfonlayer}{background}
 \pgftext{\includegraphics[width=3cm]{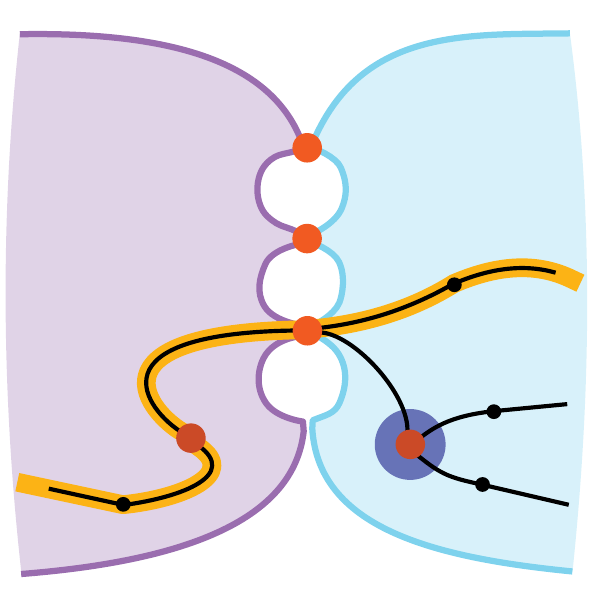}} at (C.center);
 \end{pgfonlayer}{background}
			
 \begin{pgfonlayer}{main}
 \node (C) [v:ghost] {};
 
 \end{pgfonlayer}{main}
 
 \begin{pgfonlayer}{foreground}
 \end{pgfonlayer}{foreground}

 \end{tikzpicture}
 };

 \node(M) [v:ghost] at (2,-1.5) {
 \begin{tikzpicture}

 \pgfdeclarelayer{background}
		 \pgfdeclarelayer{foreground}
			
		 \pgfsetlayers{background,main,foreground}

 \begin{pgfonlayer}{background}
 \pgftext{\includegraphics[width=3cm]{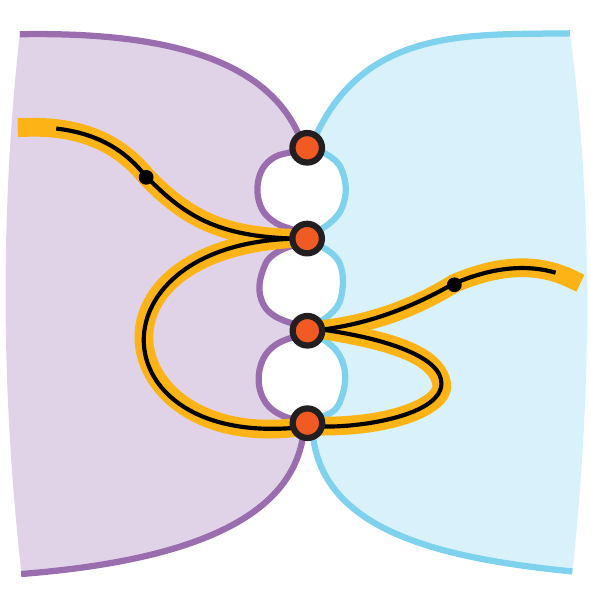}} at (C.center);
 \end{pgfonlayer}{background}
			
 \begin{pgfonlayer}{main}
 \node (C) [v:ghost] {};
 
 \end{pgfonlayer}{main}
 
 \begin{pgfonlayer}{foreground}
 \end{pgfonlayer}{foreground}

 \end{tikzpicture}
 };

 \node(M) [v:ghost] at (5,-1.5) {
 \begin{tikzpicture}

 \pgfdeclarelayer{background}
		 \pgfdeclarelayer{foreground}
			
		 \pgfsetlayers{background,main,foreground}

 \begin{pgfonlayer}{background}
 \pgftext{\includegraphics[width=3cm]{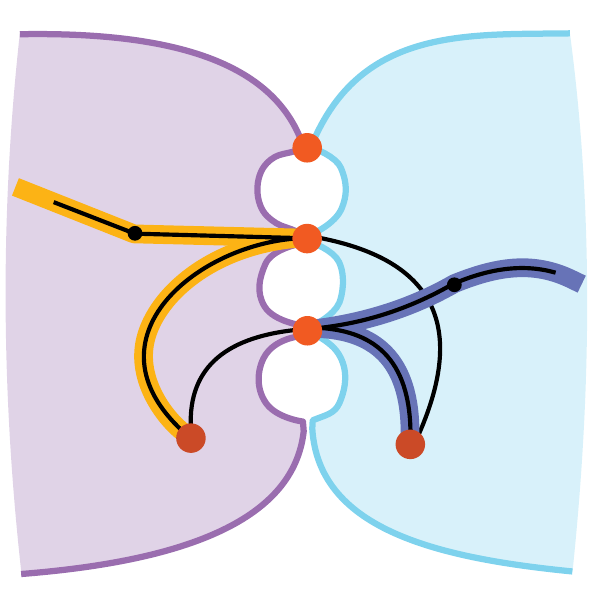}} at (C.center);
 \end{pgfonlayer}{background}
			
 \begin{pgfonlayer}{main}
 \node (C) [v:ghost] {};
 
 \end{pgfonlayer}{main}
 
 \begin{pgfonlayer}{foreground}
 \end{pgfonlayer}{foreground}

 \end{tikzpicture}
 };

 \node (i) [v:ghost] at (-3.5,-4.2) {(a) Splitting a vertex in $I$.};
 \node (ii) [v:ghost] at (3.5,-4.2) {(b) The vital linkage after splitting.};

 \node (x1) [v:ghost] at (-4.20,1.65) {\small $x$};
 \node (x2) [v:ghost] at (-2.95,1.6) {\small $x$};
 \node (x3) [v:ghost] at (-1.6,1.7) {\small $x$};

 \node (F11) [v:ghost] at (-5.25,0.5) {\small $F_1$};
 \node (F21) [v:ghost] at (-1.5,0.5) {\small $F_2$};

 \node (x11) [v:ghost] at (-4.2,-2.3) {\small $x_1$};
 \node (x12) [v:ghost] at (-3.45,-1.5) {\small $x_1$};
 \node (x21) [v:ghost] at (-2.55,-0.85) {\small $x_2$};
 \node (x22) [v:ghost] at (-1.6,-2.3) {\small $x_2$};

 \node (F12) [v:ghost] at (-5.25,-3.5) {\small $F_1$};
 \node (F22) [v:ghost] at (-1.5,-3.5) {\small $F_2$};

 \node (Ax) [v:ghost] at (2.05,0.6) {\small $x$};
 \node (AL) [v:ghost] at (2.5,1.7) {\small $L$};

 \node (Bx1) [v:ghost] at (4.25,0.65) {\small $x_1$};
 \node (BLx) [v:ghost] at (5.5,1.75) {\small $L_x$};
 \node (Bx2) [v:ghost] at (5.6,0.4) {\small $x_2$};
 \node (BLstar) [v:ghost] at (5.8,1.1) {\small $L^*$};

 \node (Cx) [v:ghost] at (2.05,-2.4) {\small $x$};
 \node (CL) [v:ghost] at (2.5,-1.3) {\small $L$};

 \node (Dx1) [v:ghost] at (4.5,-2.5) {\small $x_1$};
 \node (DL1prime) [v:ghost] at (4.25,-0.9) {\small $L_1'$};
 \node (Dx2) [v:ghost] at (5.6,-2.5) {\small $x_2$};
 \node (DL2prime) [v:ghost] at (5.75,-1.15) {\small $L_2'$};

 \end{pgfonlayer}{main}
 
 \begin{pgfonlayer}{foreground}
 \end{pgfonlayer}{foreground}

 \end{tikzpicture}
 \caption{A schematic illustration of the proof of \cref{lem:cut_open_vortices}. Figure $(a)$ highlights the operation of splitting $x,$ and the resulting $\alpha_x$. Light purple regions are the bags indexed by $F_1$ and light blue regions are the bags indexed by $F_2$. Dark blue regions highlight vertices in bags that are neighbors of $x$. Figure $(b)$ depicts two exemplary constructions of the vital linkage after splitting $x$ (from left to right), where black lines are edges in $G$.}
 \label{fig:cut_vortices_proof}
\end{figure}

\begin{lemma}\label{lem:cut_open_vortices}
 Let $d,m \geq 1$. Let $\Sigma$ be a surface. Let $\rho$ be a blank and stretched $\Sigma$-rendition of some $k$-annotated graph $(G,T)$. Let $\LLL$ be a vital $T$-linkage. Let $c \in C(\rho)$ be a vortex. Let $\alpha$ be a linear decomposition with gap $d'+1 \leq d$ of the respective vortex society $(\sigma(c),\Omega_{c})$ where the vertices $u_{-d'},\ldots,u_m$ of $V(\Omega_c)$ appear in this order in $\Omega_c$. Let $F_1,F_2 \subsetneq [m]$ be proper disjoint intervals such that the $F_i$-restriction of $\alpha$ is properly $(d-1)$-linked and such that $G_\alpha[F_1] \cap G_\alpha[F_2]$ is disjoint from $V(\Omega_c)$. Let $\SSS(\LLL)$ be the set of segments of $\LLL$ and let $\SSS(F_1,F_2)\subseteq \SSS(\LLL)$ be a set of $\rho$-homotopic paths with one endpoint in $\{u_i \mid i \in F_1\}$ and the other in $\{u_j \mid j \in F_2\}$ and let~$\omega \coloneqq \Abs{\SSS(F_1,F_2)}$.

 Then there exists a $(k+2d)$-annotated graph $(G',T')$ and a blank and stretched rendition $\rho'$ of $G'$ together with a vital $T'$-linkage $\LLL'$ such that 
 \begin{enumerate}
 \item $C(\rho') = C(\rho)$ and further $\rho'$ agrees with $\rho$ away from $c$. In particular $G'- V(\sigma_{\rho'}(c)) = G- V(\sigma_{\rho}(c)),$ and
 \item Let $\SSS(\LLL')$ denote the segments of $\LLL'$. There exists a vortex society $(\sigma_{\rho'}(c), \Omega'_{c})$ of $c$ such that $\Omega_c' = \Omega_c,$ together with a linear decomposition $\alpha'$ with gap $d'+1\leq d$ such that the vertices $u_{-d'},\ldots,u_m$ appear in this order in $\Omega_c'$. Further, the $F_i$-restriction of $\alpha'$ is properly~$(d-1)$-linked and $G_{\alpha'}[F_1] \cap G_{\alpha'}[F_2] = \emptyset$ and $\SSS(F_1,F_2)\subseteq \SSS(\LLL')$ is a set of $\rho'$-homotopic paths with one endpoint in $\{v_i \mid i \in F_1\}$ and the other in $\{u_j \mid j \in F_2\}$.
 \end{enumerate}
\end{lemma}
\begin{proof}
 Let $\mathbf{G} =\left((G,T),\rho,\LLL,\alpha,F_1,F_2,\SSS(\LLL)\right),$ each entry of the tuple of $\mathbf{G}$ being as in the theorem. Let $(\sigma(c),\Omega_c)$ be a $c$-society. Let $\alpha=\left(\langle u_{-d'},\ldots,u_m\rangle,\langle U_{-d'},\ldots,U_m\rangle \right)$ as in the theorem: note that $m\geq 3$ since the segments are disjoint. Since $F_1 \cap F_2 = \emptyset$ and the $F_i$-restriction of $\alpha$ is $d$-linear for $i=1,2$ we derive that $I \coloneqq G_{\alpha}[F_1] \cap G_{\alpha}[F_2] \subseteq V(G)$. Further, by assumption of the lemma we have 
 \begin{itemize}
 \item[$(\star)$] $I \cap V(\Omega_c) = \emptyset$.
 \end{itemize}
 
 Assume without loss of generality that $F_1=[\ell_1,\ell_2]$ and $F_2 = [k_1,k_2]$ for some $1 \leq \ell_1 < \ell_2 < k_1 < k_2 \leq m$ by possibly reversing the order of $\Omega_c$. For every $i \in [\ell_1,\ell_2-1] \cup [k_1,k_2-1]$ let $Z_i \coloneqq U_i \cap U_{i+1}$. By \cref{def:linked_vortex_normal} of properly $(d-1)$-linkedness we have the following:
 \begin{itemize}
 \item for every $i \in [\ell_1,\ell_2] \cup [k_1,k_2],$ $u_i \in U_i$ and $U_i \cap V(\Omega_c) \subseteq \{u_{i-1},u_{i}\},$ if $\ell_1=1$ let $u_{i-1} = u_{\ell_1}$
 \item for every $i \in [\ell_1,\ell_2-1] \cup [k_1,k_2-1]$ it holds $\Abs{Z_i} = d,$ and
 \item there is a $U_{\ell_1}$-$U_{\ell_2}$-linkage $M_1$ in $G[F_1]\setminus V(\Omega_c)$ of order $d$ and a $U_{k_1}$-$U_{k_2}$-linkage $M_2$ in $G[F_2]\setminus V(\Omega_c)$ of order $d$.
 \end{itemize}
 
 Let $x \in I$. We define the \emph{split of $\mathbf{G}$ at $x$} to be the tuple $\mathbf{G}_x =\left((G_x,T_x),\rho_x,\LLL_x,\alpha_x,F_1,F_2,\SSS(\LLL_x)\right)$ obtained as follows.
 
 \begin{description}
 \item[$G_x$:] Define $G_x$ to be the graph obtained from $G$ by adding two fresh vertices $x_1,x_2 \notin V(G)$. Then for $i=1,2$ and every edge $xv \in G_{\alpha}[F_i]$ add the edge $x_iv$ to $G$. Finally delete $x$ from $G$. Let $V_c \coloneqq V(\sigma_\rho(c)) \setminus \{x\},$ then $V_c \subseteq V(G_x);$ see \cref{fig:cut_vortices_proof} $(a)$ for a schematic illustration of the operation where $I$ is highlighted in \textcolor{amaranth}{red}.
 
 \item[$\rho_x$:] Let $C(\rho_x) = C(\rho)$ and let $\rho_x$ agree with $\rho$ on all cells away from $c$. Note that by $(\star)$ we derive $x \notin N(c)$. Finally fix $\sigma_{\rho_x}(c) \coloneqq G_x[V_c \cup \{x_1,x_2\}]$. Then $\rho_x$ is a $\Sigma$-rendition of $G_x$ as we did not alter the nodes on the boundary of $c$.

 \item[$\LLL_x$ and $T_x$:] Since $\LLL$ is vital in $G,$ there is a unique path $L \in \LLL$ with $x \in V(L)$. Let $\LLL^- = \LLL \setminus \{L\}$. Let $ux, xw \in E(L)$ be the edges adjacent to $x;$ possibly $x$ is a terminal of $L,$ in that case there is only one such edge and the arguments are analogous. There are four cases to consider.
 
 If $ux_1, x_1w \in E(G_x)$ then let $L_x$ be obtained from $L$ by replacing $ux$ and $xw$ with $ux_1$ and $x_2w;$ clearly $L_x$ is a path in $G_x$ with the same endpoints as $L$. In that case add a path $L^*$ consisting solely of the vertex $x_2$. By construction $\LLL_x \coloneqq \LLL^- \cup \{L_x, L^*\}$ is a vital linkage in $G_x$ and $T_x \coloneqq \{v \mid v \text{ is an endpoint of } L \in \LLL_x\}$ with $T_x = T \cup \{x_2\}$ whence $\Abs{T_x} \leq \Abs{T} + 1$. See the top two illustrations in \cref{fig:cut_vortices_proof} $(b),$ giving a schematic illustration of the construction of $L_x$ and $L^*$ in which the path $L$ uses an edge between vertices of $I$ highlighting the construction in said ``corner case''.

 If the first case is not true but $ux_2, x_2w \in E(G_x)$ then construct $\LLL_x$ as above by symmetry.

 Thus we are left with $ux_1 \in E(G_x)$ and $x_2w \in E(G_x)$ or $ux_2 \in E(G_x)$ and $x_1w \in E(G_x);$ assume the former since the cases are symmetrical. Note that in that case $x$ is not a terminal, hence at least to edges $ux,xw \in E(G)$ exist. With this in mind let $L_1,L_2 \subsetneq L$ be subpaths on at least one edge such that $L_1 \cup L_2 = L$ and $V(L_1 \cap L_2) = \{x\}$. Without loss of generality let $ux \in E(L_1)$ and $xw \in E(L_2)$. Finally let $L_1'\subsetneq G_x$ be obtained from $L_1$ by replacing $ux$ with $ux_1$ and let $L_2'\subsetneq G_x$ be obtained from $L_2$ by replacing $xw$ with $x_2w$. Then $\LLL_x \coloneqq \LLL^- \cup \{L_1',L_2'\}$ is a linkage in $G_x$ by construction. Further it is vital in $G_x$ for else $\LLL$ was not vital in $G;$ note that all the newly introduced edges cannot be used by paths different from $L_1',L_2'$ since they are all incident with $\{x_1,x_2\}$. Finally let $T_x = \{v \mid v \text{ is an endpoint of } L \in \LLL_x\},$ then $T_x = T\cup \{x_1,x_2\}$ whence $\Abs{T_x} \leq \Abs{T} +2$. See the bottom two illustrations in \cref{fig:cut_vortices_proof} $(b),$ giving a schematic illustration of the construction of $L_1'$ and $L_2',$ in which the path $L$ uses two edges between vertices of $I,$ highlighting the construction in said ``corner case''.
 \item[$\alpha_x:$] We define $\alpha_x=\left(\langle u_{-d'},\ldots,u_m \rangle,\langle U_{-d'}',\ldots,U_m'\rangle \right)$ as follows. For $i\in[\ell_1,\ell_2]$ with $x \in U_i$ define $U_i' \coloneqq (U_i\setminus\{x\}) \cup \{x_1\}$ and if $x \notin U_i$ define $U_i' \coloneqq U_i$. For $i\in[k_1,k_2]$ with $x \in U_i$ let $U_i' \coloneqq U_i\setminus\{x\} \cup \{x_2\}$ and if $x \notin U_i$ define $U_i' \coloneqq U_i$. For all $i \notin [\ell_1,\ell_2] \cup [k_1,k_2]$ let $U_i' \coloneqq U_i \setminus\{x\}$. By construction $\alpha_x$ is a linear decomposition with gap $d'+1 \leq d$. To see this note that $G_{x}[F_i] \cong G_\alpha[F_i]$ for $i=1,2$ by construction and every edge adjacent to $x_i$ is by construction part of $G_x[F_i];$ see \cref{fig:cut_vortices_proof} for a schematic representation, where the neighbors of $x$ and $x_1,x_2$ are highlighted by dark blue regions.
 \end{description}
Let $\SSS(\LLL_x)$ be the segments of $\LLL_x$ in $\rho$. By construction of $\rho_x$ and the fact that every terminal of $T_x \setminus T$ is contained in $\sigma(c) \setminus V(\Omega_c)$ by construction, we derive that $\SSS(\LLL_x) = \SSS(\LLL);$ recall \cref{def:segments}. In particular $\SSS(F_1,F_2) \subseteq \SSS(\LLL_x)$ and the paths in $\SSS(F_1,F_2)$ are $\rho_x$-homotopic, since $\rho_x$ agrees with $\rho$ on $\mathsf{cl}(\Sigma\setminus c)$.

Finally, for $i=1,2$ the $F_i$-restriction of $\alpha_x$ is still properly $(d-1)$-linked. To see this note again that by construction $G_{x}[F_i] \cong G_\alpha[F_i]$ and for every $j \in F_i$ the bags $U_j'$ and $U_j$ agree up to renaming $x$ to $x_i$. In particular, we have
 \begin{itemize}
 \item for every $i \in [\ell_1,\ell_2] \cup [k_1,k_2],$ $u_i \in U_i'$ and $U_i' \cap V(\Omega_c) \subseteq \{u_{i-1},u_{i}\};$ if $\ell_1=1$ let $u_{i-1} = u_{\ell_1},$
 \item for every $i \in [\ell_1,\ell_2-1] \cup [k_1,k_2-1]$ it holds $U_i' \cap U_{i+1}'$ and $U_i\cap U_{i+1}$ agree up to renaming $x$ to $x_1$ if $i \in [\ell_1,\ell_2-1]$ and $x_2$ if $i \in [k_1,k_2],$ thus $\Abs{U_i' \cap U_{i+1}'} = \Abs{U_i \cap U_{i+1}}=d,$
 \item there is a $U_{\ell_1}'$-$U_{\ell_2}'$-linkage $M_1'$ in $G_x[F_1]\setminus V(\Omega_c)$ of order $d$ and a $U_{k_1}'$-$U_{k_2}'$-linkage $M_2'$ in $G_x[F_2]\setminus V(\Omega_c)$ of order $d$ using that $G_x[F_i] \cong G[F_i]$ and the construction of $\alpha_x$.
 \end{itemize}

Repeating the above construction inductively splitting $\mathbf{G}$ at all $x \in I,$ we conclude the proof.
\end{proof}

Finally, we combine \cref{lem:bounded_jumps_between_vortices,lem:T-normal_segments_exist,lem:cut_open_vortices}, proving the main technical theorem of this subsection.

\begin{theorem}\label{thm:the_vortex_thing_refined}
 Let $\Sigma$ be a surface and $k^*,p,\ell,d_1,d_2\geq 1$. There exists a function $f_{\ref{thm:the_vortex_thing_refined}}:\N^2\to \N$ such that the following holds. Let $\rho$ be a blank and stretched $\Sigma$-rendition of some $k^*$-annotated graph $(G,T)$. Let $\LLL$ be a vital $T$-linkage in $G$.
 
 Let $c_1,c_2 \in C(\rho)$ be vortices and $\alpha_1=\left(\langle u_{-d'},\ldots,u_{\ell}\rangle,\langle U_{-d'},\ldots,U_\ell\rangle \right)$ a linear decomposition of the vortex society $(\sigma(c_1),\Omega_{c_1})$ with gap at most $d'+1 \leq d_1,$ and $\alpha_2=\left( \langle w_{-d''},\ldots,w_{p}\rangle,\langle W_{-d''},\ldots,W_{p}\rangle\right)$ with gap at most $d''+1 \leq d_2$ (if $c_1=c_2$ let $\alpha_1=\alpha_2$).
 Further, let $F_1 \subseteq [\ell]$ and $F_2\subseteq [p]$ be intervals such that $\{u_i \mid i \in F_1\} \cap \{u_j \mid j \in F_2\} = \emptyset$ and such that the $F_i$-restriction of $\alpha_i$ is properly $(d_i-1)$-linked. Let $\SSS(\LLL)$ be the segments of $\LLL$ and let $\SSS(F_1,F_2)=\{S_1,\ldots,S_{\omega^*}\} \subseteq \SSS(\LLL)$ be a linkage of $\rho$-homotopic segments with one endpoint in $\{u_i \mid i \in F_1\}$ and one in $\{u_j \mid j \in F_2\}$ for some $\omega^*\geq 1$. Then $\Abs{\SSS(F_1,F_2)} \leq f_{\ref{thm:the_vortex_thing_refined}}(k^*,\max(d_1,d_2))$. 
 
 Moreover $f_{\ref{thm:the_vortex_thing_refined}}(k^*,d) \in \mathbf{poly}(k^*)\cdot 2^{\mathbf{poly}(d)}$.
\end{theorem}
\begin{proof}
 Fix $d \coloneqq \max(d_1,d_2)$. Let $f(k^*,d) = (2(k^*+2d)+1)^2(g_{\ref{lem:bounded_jumps_between_vortices}}(2d)+1)+8;$ we claim that $f$ satisfies the theorem. Clearly $f(k^*,d) \in \mathbf{poly}(k^*)\cdot 2^{\mathbf{poly}(d)}$ by definition. Towards a contradiction, assume that $\omega^* \geq f(k^*,d)$. By assumption on $\SSS(F_1,F_2)$ and the fact that $d \geq 1$ we derive that $\Abs{F_1},\Abs{F_2} \geq 4$.\footnote{This is a technical detail, note that the segments can indeed be assumed to be huge for otherwise there is nothing to show.}
 
 Let $\alpha_1,\alpha_2$ as well as $F_1,F_2$ and $\SSS(F_1,F_2)$ be as in the theorem. Note that if $c_1 \neq c_2,$ since $\rho$ is stretched, we derive that $G_{\alpha_1}[F_1] \cap G_{\alpha_2}[F_2] = \emptyset$. However, if $c_1 = c_2,$ this may not be true a priori; we deal with that first. In that case, by assumption of the theorem, $\alpha_1=\alpha_2$ whence $d_1=d_2;$ let $F_1=[\ell_1,\ell_2]$ and $F_2=[k_1,k_2]$ be the respective segments where $\ell_i,k_i \geq 1$ for $i=1,2$. Since $\{u_i \mid i\in F_1\} \cap \{u_j \mid j \in F_2\} = \emptyset$ we may assume without loss of generality that $\ell_1< \ell_2 < k_1 < k_2$. Let $\tilde F_1= [\ell_1+1,\ell_2-1]$ be obtained from $F_1$ by discarding its endpoints and similarly $\tilde F_2 = [k_1+1,k_2-1]$. Since the $F_i$-restriction of $\alpha_i$ is properly $(d_i-1)$-linked, we derive by \cref{def:linked_vortex_normal} that $(U_{\ell_2-1} \cap U_{k_1+1}) \cap V(\Omega_{c_1}) = \emptyset$ as well as $(U_{\ell_1+1} \cap U_{k_2-1}) \cap V(\Omega_{c_1}) = \emptyset$. In particular $G_{\alpha_1}[\tilde F_1] \cap G_{\alpha_1}[\tilde F_2]$ is disjoint from $V(\Omega_{c_1}) = V(\Omega_{c_2})$ by construction. Note that $\SSS(\tilde F_1, \tilde F_2) \subseteq \SSS( F_1, F_2),$ defined in the obvious way, satisfies 
 $$\Abs{\SSS(\tilde F_1, \tilde F_2)} \geq \Abs{\SSS( F_1, F_2)} - 8 \geq (2(k^*+2d)+1)^2\cdot(g_{\ref{lem:bounded_jumps_between_vortices}}(2d)+1),$$
 since $\SSS( F_1, F_2)$ is a set of internally disjoint paths by \cref{obs:segments}, and thus at most $2$ paths in $\SSS( F_1, F_2)$ may use $u_i$ for every $i\in \{\ell_1,\ell_2,k_1,k_2\}$.
 
 \cref{lem:cut_open_vortices} implies the existence of a $(k^*+2d)$-annotated graph $(G',T')$ with a vital linkage $\LLL'$ together with a $\Sigma$-rendition $\rho'$ and a linear decomposition $\alpha'$ with gap at most $d' \leq d_1$ of the vortex $c \in C(\rho') =C(\rho)$ and disjoint proper segments $F_1',F_2',$ such that the $F_i'$-restriction to $\alpha'$ is properly $(d_1-1)$-linked for $i=1,2$. Further $G_{\alpha'}[F_1'] \cap G_{\alpha'}[F_2'] = \emptyset$ and there is a set of $\rho'$-homotopic segments of $\SSS(F_1',F_2')$ satisfying $\Abs{\SSS(F_1',F_2')} \geq (2(k^*+2d)+1)^2\cdot(g_{\ref{lem:bounded_jumps_between_vortices}}(2d)+1)$. In particular it suffices to prove that $\Abs{\SSS(F_1',F_2')} \geq (2(k^*+2d)+1)^2\cdot(g_{\ref{lem:bounded_jumps_between_vortices}}(2d)+1)$ is impossible, which in turn contradicts the assumption that $\omega^* \geq f(k^*,d)$ as desired. Since $(G',T'),\LLL',\rho',\alpha',F_1',F_2',\SSS(F_1',F_2')$ still satisfy the assumptions of the lemma, we may set $$k \coloneqq k^*+2d,~\omega \coloneqq (2(k^*+2d)+1)^2(g_{\ref{lem:bounded_jumps_between_vortices}}(2d)+1),$$ and, for notational simplicity, continue with the instance $(G,T),\rho,\alpha_1,\alpha_2,c_1,c_2,F_1,F_2,\SSS(F_1,F_2)$ by renaming, with the additional assumption that
 \begin{itemize}
 \item[$(\cap)$] $G_{\alpha_1}[F_1] \cap G_{\alpha_2}[F_2] = \emptyset$.
 \end{itemize}
 Due to the identification $k=k^* + 2d,$ and the choice $\omega = (2(k^*+2d)+1)^2(g_{\ref{lem:bounded_jumps_between_vortices}}(2d)+1)$ we are left to refute the assumption that $\omega \geq f(k,d)$.

 From here on, the proofs for $c_1=c_2$ and $c_1 \neq c_2$ are analogous. Let $i \in \{1,2\}$.

 Since $\SSS(F_1,F_2)$ is a set of $\rho$-homotopic paths, \cref{obs:strip_from_homotopic_paths} implies that there exist $S_1,S_{\omega} \in \SSS(F_1,F_2)$ with traces $\gamma_1,\gamma_{\omega}$ such that $\gamma_1,\gamma_{\omega}$ guard an $(F_1,F_2)$-strip $\Delta(F_1,F_2)$ that is vortex-free, and for every $S \in \SSS(F_1,F_2)$ its trace is contained in $\Delta(F_1,F_2)$. Further, by \cref{def:segments} of segments, the paths in $\SSS(F_1,F_2)$ are disjoint from vortices other than $c_1,c_2,$ hence the boundary of $\Delta(F_1,F_2)$ is disjoint from the boundary of all vortices that are not $c_1,c_2$. Note that given the disc $\Delta(F_1,F_2),$ the four endpoints of $S_1$ and $S_2$ agree with $\{u_{\ell_1},u_{\ell_2},w_{k_1},w_{k_2}\}$.

 Let $\SSS^*=\{S_1,\ldots,S_{\omega}\}$ be an enumeration of the paths in $\SSS(F_1,F_2)$ such that the following holds: Let the endpoints of $S_i$ be $\{u_{i_u},w_{j_w}\}$ for respective $i_u \in [\ell_1,\ell_2]$ and $j_w \in [k_1,k_2],$ then~$1_u \leq \ldots \leq \omega_u$.

 Since the paths in $\SSS^*$ are $\rho$-homotopic and each trace of a path in $\SSS^*$ is contained in the disc $\Delta(F_1,F_2),$ we derive that either $1_w \leq \ldots \leq \omega_w$ appear in this order or in reverse order, i.\@e.\@, $\omega_w \leq \ldots \leq 1_w;$ by reversing the cyclic order of $\Omega_{c_2}$ we may assume without loss of generality that $1_w \leq \ldots \leq \omega_w$ appear in this order. Note here that if $c_1 = c_2$ this results in possibly fixing $\Omega_{c_1} \neq \Omega_{c_2},$ for they may have reversed orders. In particular now $\alpha_1 \neq \alpha_2$ is possible (but the bags still agree up to bijection).

 For simplicity, let us assume that $\alpha_1=\left(\langle u_{-d'},\ldots,u_{\ell}\rangle,\langle U_{-d'},\ldots,U_{\ell}\rangle \right)$ and $\alpha_2=\left(\langle w_{-d''},\ldots,w_{p}\rangle,\langle W_{-d'},\ldots,W_{p}\rangle\right)$ be the resulting decompositions after the previous possible reversing of the orders. By the above analysis we derive the following:
 \begin{itemize}
 \item[$(\star)$] Let $S_j \in \SSS^*,$ then $S_j$ has endpoints $u_{j_u},w_{j_w}$ for respective $j_u \in [\ell],$ $j_w \in [p]$ and $1_u \leq \ldots \leq {\omega}_u$ as well as ${1_w} \leq \ldots \leq \omega_w$.
 \end{itemize}

Furthermore, 
\begin{itemize}
 \item[$\boxtimes$] if $c_1 = c_2,$ then $d' = d''$ and $k = p$. Further, given $u \in N(c_1)$ and $i \in [\ell+d']$ and $j \in [p + d'']$ with $u=u_i$ and $w=w_j$ we derive that $U_i = W_j$.
\end{itemize}

Fix $d \coloneqq \max{(d_1,d_2)}$. Let $I_u\coloneqq \{1_u,\ldots,\omega_u\}$ and $I_w\coloneqq \{1_w,\ldots,\omega_w\}$.

Let $\alpha_u^*$ be the $[1_u,\omega_u]$-restriction of $\alpha_1$ and let $\alpha_w^*$ be the $[1_w,\omega_w]$-restriction of $\alpha_2$ respectively. Recall that $\omega \geq (2k+1)^2(g_{\ref{lem:bounded_jumps_between_vortices}}(2d)+1)$. Since $1_u \leq \ldots \leq \omega_u$ by construction, applying \zcref{lem:T-normal_segments_exist} to $\alpha_u^*$ and $I_u,$ we derive the existence of $i,j \in [k]$ such that $\Abs{[i_u,j_u]\cap I_u} \geq (2k+1)g_{\ref{lem:bounded_jumps_between_vortices}}(2d)$ and such that the $[i_u,j_u]$-segment of $\alpha_u^*$ is $T$-normal. Let $\SSS_u \coloneqq \{S_j \in \SSS^* \mid j\in [i_u,j_u]\},$ then $\Abs{\SSS_u} \geq (2k+1)(g_{\ref{lem:bounded_jumps_between_vortices}}(2d)+1)$. 
 
Since $S_m$ has endpoints $\{u_{m_u},w_{m_w}\}$ and since for $S_m,S_{m'}$ it holds $m_u\leq m'_u \iff m_w \leq m'_w$ we derive that $\Abs{I_w \cap [i_w,j_w]}\geq (2k+1)(g_{\ref{lem:bounded_jumps_between_vortices}}(2d)+1)$. Again, applying \zcref{lem:T-normal_segments_exist} to $J_w \coloneqq I_w \cap [i_w,j_w]$ we derive the existence of $s,t \in [i,j]$ such that $\Abs{[s_w,t_w]\cap J_w} \geq g_{\ref{lem:bounded_jumps_between_vortices}}(2d)+1$ and such that the $[s_w,t_w]$-segment of $\alpha_w^*$ is $T$-normal. Finally fix $J_u \coloneqq [s_u,t_u]$.
 
Summarizing we derive that for $\SSS \coloneqq \{S_j \in \SSS^* \mid j\in[s,t] \},$ it holds \begin{itemize}
 \item[$(\star\star)$] $\Abs{\SSS} \geq g_{\ref{lem:bounded_jumps_between_vortices}}(2d)+1,$
 \end{itemize}
 and the $[s_u,t_u]$-restriction $\alpha_u$ of $\alpha_1$ as well as the $[s_w,t_w]$-restriction $\alpha_w$ of $\alpha_2$ are by construction $T$-normal as well as properly $(d_1-1)$-linked and properly $(d_2-1)$-linked, respectively. 
 
$$ \text{Let $U \coloneqq \bigcup_{j=s_u}^{t_u}U_j$ and let $W \coloneqq \bigcup_{j=s_w}^{t_w}W_j$}.$$ 
By construction $U \cup W \subseteq V(\sigma(c))$ and $\SSS$ is a $U$-$W$-linkage.
 
 Let $\LLL^* \subseteq \LLL$ be maximal such that for every path $L \in \LLL^*$ it holds $V(L) \cap (U \cup W) \neq \emptyset$. 
 \begin{beautifulclaim}\label{thm:the_vortex_thing_refined_stopit}
 Every path in $\SSS$ is a subpath of some path in $\LLL^*$. Further every path in $\LLL^*$ has at least one vertex in two of $\{U_{s_u},U_{t_u},W_{s_w},W_{t_w}\}$.
 \end{beautifulclaim}
 \begin{claimproof}
 Since $\rho$ is a blank rendition, all terminals of paths in $\LLL^*$ are part of the graph of some vortex. If an endpoint of $L \in \LLL^*$ is contained in $U$ then it is contained in $\bigcap_{j=s_u}^{t_u}U_j $ since $\alpha_u$ is $T$-normal, and similarly if an endpoint of $L$ is contained in $W$ then it is contained in $\bigcap_{j=s_w}^{t_w}W_j $ since $\alpha_w$ is $T$-normal. If a path $L \in \LLL^*$ has no endpoint in $U$ but intersects $U$ then this implies that $V(L) \cap (U_{s_u} \cup U_{t_u}) \neq \emptyset$ for $\alpha_u$ is a linear decomposition, and an analogous result holds if it has no endpoint in $W$ but intersects $W$. Summarising the above, every path $L \in \LLL^*$ admits at least one vertex in two of $\{U_{s_u},U_{t_u},W_{s_w},W_{t_w}\};$ note that the path may have the same vertex in $U_{s_u} \cap U_{t_u}$ but no vertex in $W_{s_w} \cup W_{t_w}$ if it is disjoint from $\SSS$ for example.
 
 Finally, let $P \in \SSS$. $P$ has one endpoint in $U$ and the other in $W$. Since $\SSS \subseteq \SSS(\LLL)$ by construction we derive that $P \subseteq L$ for some $L \in \LLL$. In particular $L \in \LLL^*$ by definition of the latter. This concludes the proof.
 \end{claimproof}

 Recall that $\SSS=\{S_s,\ldots,S_t\}$ with $t-s \geq g_{\ref{lem:bounded_jumps_between_vortices}}(2d)+1 $ such that $S_j$ has endpoints $u_{j_u},w_{j_w}$ and $s_u \leq \ldots \leq {t_u}$ as well as $s_w \leq \ldots \leq {t_w}$ by $(\star)$ for every $j \in [s,t]$. Furthermore, recall that the paths in $\SSS$ are $\rho$-homotopic. 
 
 Let $\gamma_s,\gamma_{t}$ be the traces of $S_s$ and $S_t$. Let $\Delta \coloneqq \Delta([s_u,t_u],[s_w,t_w])$ be the $([s_u,t_u],[s_w,t_w])$-strip guarded by $\gamma_1,\gamma_t$ as given by \cref{obs:strip_from_homotopic_paths}. Then by construction $\Delta$ contains the traces of all paths in $\SSS$.
\begin{description}
 \item[$G_1$:] Let $G_1\subsetneq G$ be obtained from $G$ by deleting all edges $E(\bigcup_{i=1}^m C_i) \setminus E(\bigcup \LLL),$ then $(G_1,T)$ is still an annotated graph and $\LLL$ is a vital linkage in $G_1$. Let $\rho_1$ be the restriction of $\rho$ to $G_1$. Since we did not alter the graphs of vortices, every vortex cell of $\rho$ is one of $\rho_1$ and their societies agree~--~up to a choice of ordering~--~in particular $\alpha_u$ and $\alpha_w$ are still $T$-normal and properly $(d_1-1)$-linked and properly $(d_2-1)$-linked respectively. 
 
 \item[$G_2$:] Let $G_2\subseteq G_1$ be obtained from $G_1$ by deleting all the paths in $\LLL^- \coloneqq \LLL \setminus \LLL^*$. Let $T_2$ be the set of endpoints of paths in $\LLL^*$ then $(G_2, T_2)$ is still an annotated graph. One easily verifies that $\LLL_2 \coloneqq \LLL^*$ is a vital $T_2$-linkage in $G_2$ by construction; otherwise $\LLL$ was not vital for $(G_1,T)$. Let $\rho_2$ be the restriction of $\rho_1$ to $G_2$. Then $c_1,c_2$ are still vortices in $\rho_2$. Note that $V(\bigcup \LLL^-) \cap (U \cup W) = \emptyset$. In particular, $\alpha_u$ and $\alpha_w$ remain both $T_2$-normal as well as properly $(d_1-1)$-linked and properly $(d_2-1)$-linked given the $c_1$-society and $c_2$-society in $\rho_2,$ respectively. Note that every vortex of $\rho_2$ is a vortex of $\rho_1,$ but it may have fewer vortices. 

 \item[$G_3:$] Let $C^*\subseteq C(\rho_2)\setminus\{c_1,c_2\}$ be the set of vortex cells distinct from $c_1,c_2$. Let $G_3\subseteq G_2$ be obtained from $G_2$ by deleting all vertices in $\bigcup_{c \in C^*}V(\sigma(c))$ and let $\LLL_3 = \LLL_2 \cap G_3$ be the resulting linkage. Then $\LLL_3$ is vital by construction and \cref{lem:vital_and_deletion}. Note that the set of endpoints $T_3$ of paths in $\LLL_3$ are part of cells that are outside of $\Delta$ since the boundary of $\Delta$ is disjoint from the boundary of vortex cells in $C^*$. Furthermore, by construction it holds that $G_3 = V(\sigma(c_1))\cup V(\sigma(c_2))\cup \LLL^*$. Let $\rho_3$ be the restriction of $\rho_2$ to $G_3$. Since we did not alter $V(\sigma_{\rho_2}(c_1)) \cup V(\sigma_{\rho_2}(c_2)),$ and since the rendition $\rho$ was stretched~--~in particular $\rho_2$ is~--~we derive that $\alpha_u$ and $\alpha_w$ are still $T_3$-normal as well as properly $(d_1-1)$-linked and properly $(d_2-1)$-linked, respectively. 

 \item[$G^*$:] Let $G^*$ be obtained from $G_3$ by deleting every vertex that is part of a non-vortex cell contained in $\mathsf{cl}(\Sigma \setminus \Delta)$. Let $\PPP^* = \LLL_3 \cap G^*,$ then again $\PPP^*$ is vital by \cref{lem:vital_and_deletion}. Let $\rho^*$ be the restriction of $\rho_3$ to $G^*$. Let $T^*$ be the endpoints of paths in $\PPP^*,$ then $T^* \subseteq V(\sigma_{\rho^*}(c_1))\cup V(\sigma_{\rho^*}(c_2))$. Further by construction $T^* \cap (U \cup W) \subseteq T$. In particular $\alpha_u,\alpha_w$ are $T^*$-normal and since we did not delete any vertices that are part of $V(\sigma_{\rho_3}(c_1))\cup V(\sigma_{\rho_3}(c_2)),$ they are both still properly $(d_1-1)$-linked and properly $(d_2-1)$-linked respectively.
 Further $G^*= \sigma_{\rho^*}(c_1)\cup \sigma_{\rho^*}(c_2) \cup \bigcup \PPP^* $ by construction. 
\end{description}
 Finally let $K_u = [s_u,t_u]$ and $K_w = [s_w,t_w]$. Then these are still segments of the $c_1$- and $c_2$-societies in $\rho^*$ (choosing the respective cyclic orientation) and by the above the $K_u$- and $K_w$-restrictions of the respective linear decompositions are still properly $(d_1-1)$-linked and properly $(d_2-1)$-linked as well as $T^*$-normal, respectively. Note that we did not alter $G_{\alpha_x}[K_x]$ throughout the above construction for $x \in \{u,w\}$. In particular, by $(\cap)$ we derive that $G^*_{\alpha_u}[K_u] \cap G^*_{\alpha_2}[K_w] = \emptyset$. Let now $\alpha_1^*,\alpha_2^*$ be linear decompositions of the vortex societies $(\sigma_{\rho^*}(c_1),\Omega_{c_1}^*)$ and $(\sigma_{\rho^*}(c_2),\Omega_{c_2}^*)$ respectively, such that $\alpha_u$ is the $K_u$-restriction of $\alpha_1^*$ and $\alpha_2$ is the $K_w$-restriction of $\alpha_2^*;$ these clearly exist. Note that we may have $\alpha_1^* \neq \alpha_2^*$ in the case that $c_1=c_2,$ simply because they need to respect the inherent order given by $\alpha_u,\alpha_w$\footnote{Although not needed, one may guarantee that they satisfy $(\boxtimes),$ i.\@e.\@, their bags agree, but $\alpha_1^*$ and $\alpha_2^*$ may be defined with respect to different orderings.}.
 
 Finally, $(G^*,T^*),\rho^*,\PPP^*,\SSS,c_1,c_2,d_1,d_2,\alpha_1^*,\alpha_2^*,\alpha_u,\alpha_w$ satisfy the assumptions of \zcref{lem:bounded_jumps_between_vortices}; note here that the size of $T^*$ is irrelevant and may be very large regarding the above construction. By \zcref{lem:bounded_jumps_between_vortices} we conclude that $\Abs{\SSS^*} \leq g_{\ref{lem:bounded_jumps_between_vortices}}(2d)$ as a contradiction to $(\star \star),$ concluding the proof.
\end{proof}

In a final iteration, we combine \zcref{lem:no_apex_vital_from_exhausted_to_dry,lem:T-normal_segments_exist} and \zcref{thm:no_apex_dry_then_clean_segments,lem:bounded_jumps_between_vortices} to derive $iii)$ and $iv)$ of \zcref{thm:no_apices_vital_reduction}.

\subsection{Taming simple loops}

In a first instance, we prove the main ingredient to prove $iii)$ of \zcref{thm:no_apices_vital_reduction}, i.\@e.\@, given an exhausted tuple $(G,T), \LLL, (\rho,M)$ in some surface $\Sigma,$ such that the locus is ``large enough'', simple loops cannot intersect cycles ``deeply nested'' in $M,$ whence we will be able to remove them completely from the instance.

We outsource the following lemma, as we will use the same technique to tame non-simple loops and links. The lemma is an easy combinatorial filtering result.
\begin{lemma}\label{lem:scheissortiereneh}
 Let $(G,\Omega)$ be a society and $d,\ell \geq 1$. Let $S_1,S_2 \subsetneq V(\Omega)$ be disjoint non-empty segments. Let $\mathsf{gap}_1,\mathsf{gap}_2 \subseteq V(\Omega)$ be segments (not necessarily disjoint) with $\Abs{\mathsf{gap}_i} \leq \ell$ for $i=1,2$ and such that either $\mathsf{gap}_1=\mathsf{gap}_2$ or $S_i \cap \mathsf{gap}_j =\emptyset$ for $\{i,j\} = \{1,2\}$. Let $u_1,\ldots,u_m$ be a labelling of $V(\Omega)$ such that they appear in that order in $\Omega$ for $m \coloneqq \Abs{V(\Omega)}$. Let $\SSS$ be a set of $16d + 4\ell$ internally disjoint paths in $G$ such that every path in $\SSS$ has one endpoint in $S_1$ and the other in $S_2,$ and no vertex in $S_1 \cup S_2$ is part of more than two paths in $\SSS$. 

 Then there exist disjoint non-empty intervals $F_1,F_2 \subsetneq [m]$ such that 
 \begin{enumerate}
 \item for $i=1,2$ it holds $\{u_j \mid j\in F_i\} \subseteq S_i$ and $\Abs{F_i} \geq d,$
 \item $\{u_i \mid i \in F_1\cup F_2\} \cap (\mathsf{gap}_1 \cup \mathsf{gap}_2)= \emptyset,$ and
 \item there is $\SSS' \subseteq \SSS$ of size $d$ such that every path in $\SSS'$ has one endpoint in $\{u_i \mid i \in F_1\}$ and the other in $\{u_i \mid i \in F_2\}$.
 \end{enumerate}
\end{lemma}
\begin{proof}
We start by defining linear orders on $S_1,S_2$ ``induced'' by $\Omega$.
\begin{description}
 \item[$\preceq_i$:] Note first that $\Abs{\SSS} \geq 16d+4\ell$ implies that $\Abs{S_i} \geq 2,$ and further since the segments are proper (for they are disjoint), there exist two distinct endpoints $s_i,t_i \in S_i$ for $i=1,2$. Let $i=1,2,$ then since $S_i$ is a proper segment we derive the existence of $x \in V(\Omega)\setminus V(S_i)$ such that either $s_i, x, t_i$ or $t_i,x,s_i$ appear in the order listed in $\Omega$. Without loss of generality assume the latter and call $s_i$ the \emph{start} of $S_i$ and $t_i$ the \emph{end} of $S_i$. Then, again since $S_i$ is a segment, we may define a linear order $ x\preceq_i y$ on $S_i$ by saying $x \preceq_i y \iff s_i,x,y,t_i$ appear in the order listed in $\Omega$. In particular $s_i \preceq_i t_i$.
\end{description}

Next we filter $\SSS$ by removing all the paths that have an endpoint in one of $\mathsf{gap}_1,\mathsf{gap}_2$
\begin{description}
 \item[$\SSS^*$:] By assumption of the lemma $\SSS$ is a set of internally disjoint paths, such that at most $2$ paths of $\SSS$ may use the same vertex $x \in \mathsf{gap}_1\cup\mathsf{gap}_2 ,$ i.\@e.\@, a total of at most $4\ell$ paths use a vertex in $\mathsf{gap}_1\cup\mathsf{gap}_2$. Let $\SSS^* \subseteq \SSS$ be the set of those paths that do not end in a vertex of $\mathsf{gap}_1\cup\mathsf{gap}_2,$ then $\Abs{\SSS^*} \geq \Abs{\SSS} - 4\ell = 16d$. 
\end{description}

Since $S_i$ and $\mathsf{gap}_i$ are both segments of $\Omega,$ we derive that $S_i \setminus \mathsf{gap}_i$ is the union of at most two segments $S_i^1 \cup S_i^2$ for $i=1,2$. 

By construction and assumption of the lemma we further derive that $S_i^p \cap \mathsf{gap}_j = \emptyset$ for all $i,j,p \in \{1,2\},$ since either $F_i \cap \mathsf{gap}_j = \emptyset$ by assumption or $S_i^p \cap \mathsf{gap}_j = \emptyset$ by construction, for $\mathsf{gap}_1 = \mathsf{gap}_2$. Let $\SSS^*_1,\SSS^*_2 \subseteq \SSS^*$ form a partition of $\SSS^*$ such that every path in $\SSS^*_i$ has an endpoint in $S_1^i$ for $i=1,2$. By the pigeonhole principle there is $i \in \{1,2\}$ such that $\Abs{\SSS^*_i} \geq 8d;$ without loss of generality let it be $i=1$ for the other case is analogous.

Similarly there is a partition $\SSS'_1,\SSS'_2 \subseteq \SSS^*_1$ such that every path in $\SSS'_i$ has an endpoint in $S_2^i$ for $i=1,2$. By the pigeonhole principle there is $i \in \{1,2\}$ such that $\Abs{\SSS'_i} \geq 4d;$ without loss of generality let it be $i=1$ for the other case is analogous.

Finally since $S_i^1$ is a segment of $\Omega,$ we derive that there are at most two disjoint intervals $F_i^1,F_i^2 \subseteq [m]$ -- one may be empty -- such that $S_i^1 = \{u_j \mid j \in F_i^1 \cup F_i^2\}$ for $i=1,2$. By two more repeated applications of the pigeonhole principle, there exist $i,j \in \{1,2\}$ such that $\PPP_i^j \coloneqq \{S \in \SSS'_1 \mid S \text{ has one endpoint in } F_1^i \text{ and the other in } F_2^j\}$ satisfies $\Abs{\PPP_i^j} \geq d$.

Then $\SSS' \coloneqq \PPP_i^j$ and $F_1 \coloneqq F_1^i$ as well as $F_2 = F_2^i$ are as desired concluding the proof.
\end{proof}

The following is the main result of this section.

\begin{theorem}\label{thm:simple_loops_not_deep}
 Let $\Sigma$ be a surface and $k,d \geq 1$. There exists a function $f_{\ref{thm:simple_loops_not_deep}}:\N^2\to \N$ such that the following holds. Let $m> f_{\ref{thm:simple_loops_not_deep}}(k,d)$ and $(G,T)$ be a $k$-annotated graph with an $m$-locus $(\rho,M)$ of $(G,T)$ in $\Sigma$ such that $\rho$ is a $\Sigma$-skeleton and
 \begin{itemize}
 \item $M$ is tight,
 \item every vortex of $C(\rho)$ is properly $(d-1)$-linked with gap at most $d,$
 \item there exists a vital $T$-linkage $\LLL$ in $G,$ and
 \item $(G,T),\LLL,(\rho,M)$ is exhausted.
 \end{itemize}
 Let $M =\langle C_1,\ldots,C_m\rangle$ and let $\SSS(\LLL)$ be the segments of $\LLL$. Then for every simple loop $S \in \SSS(\LLL)$ and every $ i \in[m-f_{\ref{thm:simple_loops_not_deep}}(k,d)]$ it holds $V(S)\cap V(C_i) =\emptyset$. 
 
 Moreover $f_{\ref{thm:simple_loops_not_deep}}(k,d)\in \mathbf{poly}(k) \cdot 2^{\mathbf{poly}(d)}$.
\end{theorem}
\begin{proof}
 We claim that $f(k,d) =16(f_{\ref{thm:the_vortex_thing_refined}}(k,d)+1)+4d$ satisfies the theorem. Let $m > f(k,d)$. For notational simplicity, let $\omega^* \coloneqq f(k,d)$. Let $(G,T),\LLL,(\rho,M)$ be as in the theorem, in particular it is exhausted and $M=\langle C_1,\ldots,C_m\rangle$ is tight. Recall \zcref{def:linkage_restriction_to_a_disc} as well as \zcref{obs:well_from_locus} and the respective definition of the respective $m$-well $$\WWW\coloneqq \left(G(M,\LLL),\Omega_{\Delta_{C_m}},\rho_{\Delta_{C_m}},M,\LLL({\Delta_{\CCC_m}})\right).$$ 
 Since by assumption of the theorem $(G,T),\LLL,(\rho,M)$ is exhausted and $M$ is tight, \zcref{lem:no_apex_vital_from_exhausted_to_dry} implies that $\WWW$ is dry.
 
 Let $S \in \SSS(\LLL)$ be a simple loop. Let $c \in C(\rho)$ be the unique vortex cell such that $S$ is a simple loop at $c$. Let $\Delta(S)$ denote the unique disc bounded by $S$ and a segment of $c$ such that $\Delta(S)$ is vortex-free and disjoint from the closure of each of vortex distinct from $c;$ this exists by \zcref{def:links} and it is $\rho$-aligned by construction. Towards a contradiction, assume that $V(S) \cap V(C_i) \neq \emptyset$ for some $1 \leq i \leq m -\omega^*;$ note that $1<m-\omega^*$ by assumption. Since $V(S) \cap V(C_i)$ is not empty, there is a path $P \in \LLL(\Delta_{C_m})$ such that $P$ is a subpath of $S$ and $V(P) \cap V(C_i) \neq \emptyset$. Recall \zcref{def:disc_of_restricted_linkage}. Then, since $\WWW$ is dry~--~in particular it is drained~--~\zcref{obs:fullbucketsindrainedwells} implies that there exist $m-i+1$ distinct paths $P_{m},\ldots,P_i \in \LLL(\Delta_{C_{m}})$ such that 
 \begin{enumerate}
 \item $V(P_j)\cap V(C_j) \neq \emptyset$ for all $j \in [i,m],$ and
 \item $\Delta_{P_m} \subsetneq \ldots \subsetneq \Delta_{P_i}$.
 \end{enumerate}
 By \zcref{obs:segments}, $\SSS(\LLL)$ is a set of internally disjoint paths. In particular for every $S' \in \SSS(\LLL)$ its trace is either contained in $\Delta(S)$ or disjoint from it. Hence, since $\Delta(S)$ is vortex-free and $\rho$-aligned by the above construction, and further $\Delta_{P_t} \subseteq \Delta_{P_j}\subseteq \Delta(S)$ for every $j \geq t \geq i ,$ \zcref{thm:no_apex_dry_then_clean_segments} implies that there is an injection $\beta:\{P_i,\ldots,P_{m}\} \to \SSS$ such that for every $j \in [i,m]$ it holds $P_j \subseteq \beta(P_j)$. Let $S_j \coloneqq \beta(P_j)$ for every $j \in [i,m]$. 

 \begin{beautifulclaim}\label{thm:simple_loops_not_deep_claim_niceAndcozy}
 For every $j \in [i,m]$ the segment $S_j$ is a simple loop at $c$ with corresponding vortex-free disc $\Delta(S_j)$ and $\Delta(S_{m}) \subsetneq \ldots \subsetneq \Delta(S_{i}) \subseteq \Delta(S)$.
 \end{beautifulclaim}
 \begin{claimproof}
 Let $j \in [i,m]$ be arbitrary, then since $S_j \in \SSS(\LLL)$ it is either a link or a loop. Since $\Delta(S)$ is a disc, and $S_j$ and $S$ are either internally disjoint or equal by \zcref{obs:segments}~--~they may only be equal for $j=m$~--~we derive that the trace of $S_j$ is either contained in the interior of $\Delta(S)$ or it is disjoint from it. Since $P_j \subseteq S_j$ and $\Delta_{P_j} \subseteq \Delta(S)$ we derive that the trace of $S_j$ is either equal to that of $S$ if $j=i$ and $S_i=S,$ or contained in the interior of $\Delta(S)$. Since $S_j$ is a segment of $\LLL$ its trace has both its endpoints on vortex boundaries. Since for every vortex $c' \neq c,$ $\Delta(S)$ is disjoint from the boundary of $c',$ we derive that the trace of $S_j$ has both its endpoints in $c$. Finally $S_j$ is a simple loop at $c$ by \zcref{def:links}.

 Similarly to \zcref{obs:discs_of_paths_in_locus_are_wellbehaved} one easily verifies that $\Delta(S_{m}) \subsetneq \ldots \subsetneq \Delta(S_{i}) \subseteq \Delta(S)$ leveraging the fact that $\Delta_{P_m} \subsetneq \ldots \subsetneq \Delta_{P_{i}}\subseteq \Delta(S)$ where $\Delta_{P_j} \subseteq \Delta(S_j)$ for all $j \in [i,m-1]$.
 \end{claimproof}

 Let $\SSS \coloneqq \{S_i,\ldots,S_{m}\}$ be as in \zcref{thm:simple_loops_not_deep_claim_niceAndcozy}. By our previous assumption $i<m -\omega^*$ whence
 \begin{itemize}
 \item[$(\star)$] $\Abs{\SSS} = m-i+1 \geq \omega^*+1 > 16(f_{\ref{thm:the_vortex_thing_refined}}(k,d)+1)+4d$.
 \end{itemize} 

 Since $c$ is properly $(d-1)$-linked with gap at most $d,$ \zcref{def:linked_vortex} implies the existence of a linear decomposition $\alpha=\left(\langle u_{-d'},\ldots,u_0,u_1,\ldots,u_{\ell}\rangle,\langle U_{-d'},\ldots,U_0,U_1,\ldots,U_\ell\rangle \right)$ of $(\sigma(c),\Omega_c),$ with gap $d'$ for some $d' < d,$ such that the $[\ell]$-restriction of $\alpha$ is properly $(d-1)$-linked. Note that by \cref{obs:segments} $\SSS$ is a set of internally disjoint paths with endpoints in $V(\Omega_c)$ and every vertex of $V(\Omega_c)$ is part of at most two paths of $\SSS$. Let $\mathsf{gap}\coloneqq\{u_i \mid i \in [-d',0]\}$. By \cref{thm:simple_loops_not_deep_claim_niceAndcozy} we have $\Delta(S_m) \subsetneq \ldots \subsetneq \Delta(S_i),$ hence by ``planarity'' -- recall that $\rho$ is a $\Sigma$-skeleton -- there exist disjoint segments $T_1,T_2 \subsetneq V(\Omega_c)$ of $\Omega_c$ such that every path in $\SSS$ has one endpoint in $T_1$ and the other in $T_2,$ and the endpoints of the segments are part of such paths. Let $\Delta(T_1,T_2)$ be a $(T_1,T_2)$-strip guarded by the traces of $S_i$ and $S_m$ which contains all the traces of $\SSS;$ see also \cref{obs:strip_from_homotopic_paths}. Let $\Omega^*$ be the cyclic order of $\Omega_c$ restricted to $T_1 \cup T_2$. Let $G^*$ be the crop of $G$ by $\Delta(T_1,T_2)$.

 By applying \cref{lem:scheissortiereneh} to $(G^*,\Omega^*),$ $\SSS,$ $\mathsf{gap},$ $S_1,$ $S_2$ together with $(\star)$ we derive the existence of disjoint intervals $F_1,F_2 \subsetneq [\ell]$ and a subset $\SSS^* \subseteq \SSS$ with one endpoint in $\{u_j \mid j \in F_1\}$ and the other endpoint in $\{u_j \mid j \in F_2\}$ such that 
 \begin{itemize}
 \item[$(\star\star)$] $\Abs{\SSS^*} \geq f_{\ref{thm:the_vortex_thing_refined}}(k,d)+1$. 
 \end{itemize}
 Without loss of generality assume that for every $i \in F_1$ and $j \in F_2$ it holds $i<j;$ otherwise rename the intervals accordingly.

Let $S_1,\ldots,S_p$ be an enumeration of $\SSS^*$ such that $\Delta(S_1) \subsetneq \ldots \subsetneq \Delta(S_p)$. Further, for every $j \in [p]$ denote the endpoints of $S_j$ by $u_{i_u},w_{i_w}$ where $i_u \in F_1$ and $i_w \in F_2$. 
 
Again, since $\Delta(S_1) \subsetneq \ldots \subsetneq \Delta(S_{p}),$ each of the paths being a simple loop, we derive that $1_u \leq \ldots \leq p_u < p_w \leq \ldots \leq 1_w$.

 By construction the $F_i$-restriction of $\alpha$ is properly $(d-1)$-linked. Let $c_1,c_2 \coloneqq c$. Then $G,\rho,c_1,c_2,\SSS^*,\alpha,F_1,F_2$ satisfy the assumptions of \zcref{thm:the_vortex_thing_refined}. Thus, said lemma implies that $\Abs{\SSS}\leq f_{\ref{thm:the_vortex_thing_refined}}(k,d)$ as a contradiction to $(\star\star),$ concluding the proof.
 
\end{proof}

Note that while \cref{thm:simple_loops_not_deep} bounds how ``deep'' simple loops may interact with a locus, it does not bound the overall number of homotopic simple loops, simply because $\rho$-homotopic simple loops at some common vortex $c$ may bound disjoint discs. There is no way to ``efficiently'' bound the number of such simple loops, but fortunately we will not have too.

\subsection{Controlling homotopic links and non-simple loops}

In a second instance, we prove $iv)$ of \zcref{thm:no_apices_vital_reduction}, i.\@e.\@, given an exhausted tuple $(G,T), \LLL, (\rho,M)$ in some surface $\Sigma,$ such that the locus is ``large enough'', there are not ``too many'' links or non-simple loops of the same homotopy type. The proof is easier than that of \cref{thm:simple_loops_not_deep} as it does not rely and is independent of the size of the locus.

\begin{theorem}\label{thm:links_not_many}
 Let $\Sigma$ be a surface and $k,d,m \geq 1$. There exists a function $f_{\ref{thm:links_not_many}}:\N^2\to \N$ such that the following holds. Let $(G,T)$ be a $k$-annotated graph with a blank and stretched~$\Sigma$-skeleton $\rho$. Let $C^*\subseteq C(\rho)$ be the set of vortices such that
 \begin{itemize}
 \item every vortex of $C(\rho)$ is properly $(d-1)$-linked with gap at most $d,$
 \item there exists a vital $T$-linkage $\LLL$ in $G$.
 \end{itemize}
 Let $S$ be a segment of $\LLL$ that is not a simple loop. Let $\SSS(S)$ be the set of all segments of $\LLL$ that are $\rho$-homotopic to $S$. Then $\Abs{\SSS(s)} < f_{\ref{thm:links_not_many}}(k,d)$. 
 
 Moreover, $f_{\ref{thm:links_not_many}}(k,d)\in \mathbf{poly}(k) \cdot 2^{\mathbf{poly}(d)}$.
\end{theorem}
\begin{proof}
 Choose $f(k,d) =16(f_{\ref{thm:the_vortex_thing_refined}}(k,d)+1)+4d;$ we claim that $f$ satisfies the theorem.
 
 Let $S$ be as in the theorem, i.\@e.\@, for some $L \in \LLL$ the path $S \subseteq L$ is a segment of $L$ that is not a simple loop. Let $c_1,c_2 \in C(\rho)$ be vortices (possibly equal if $S$ is a non-simple loop) such that $S$ has its endpoints on the nodes of these vortices. Let $\SSS(\LLL)$ be the set of segments of $\LLL$ and let $\SSS(S)\subseteq \SSS(\LLL)$ be as in the theorem. Towards a contradiction, assume that $\Abs{\SSS(S)} \geq f(k,d)$.

 Now let $$\alpha_1\coloneqq \left(\langle u_{-d'},\ldots,u_0,u_1,\ldots,u_{\ell}\rangle,\langle U_{-d'},\ldots,U_0,U_1,\ldots,U_{\ell}\rangle\right)$$ be a properly $(d-1)$-linked linear decomposition of the $c_1$-society $(\sigma(c_1),\Omega_{c_1})$ with gap $d' \leq d$ and let $$\alpha_2 \coloneqq \left(\langle w_{-d''},\ldots,w_0,w_1,\ldots,w_{p}\rangle,\langle W_{-d''},\ldots,W_0,W_1,\ldots,W_{p}\rangle\right)$$ be a properly $(d-1)$-linked linear decomposition of the $c_2$-society $(\sigma(c_1),\Omega_{c_1})$ with gap $d'' \leq d;$ if $c_1=c_2$ let $\alpha_1 = \alpha_2$.

 Since $\SSS(s)$ is a set of $\rho$-homotopic paths, \cref{obs:strip_from_homotopic_paths} implies the existence of disjoint non-empty segments $T_i \subsetneq V(\Omega_{c_i})$ and two paths $S_1,S_{\omega}$ with respective traces $\gamma_1,\gamma_\omega$ such that~$\SSS(s)$ is a $T_1$-$T_2$-linkage and the trace of each path of $\SSS(s)$ is contained in the $(T_1,T_2)$-strip $\Delta(T_1,T_2);$ let $\Abs{T_1}=s$ and $\Abs{T_2}=t$

 Let $\mathsf{gap}_1^*\coloneqq\{u_i \mid i \in [-d',0]\}$ and let $\mathsf{gap}_2^* \coloneqq \{w_i \mid i \in [-d'',0]\},$ and if $c_1 = c_2$ then $\mathsf{gap}_1=\mathsf{gap_2}$. Choose an orientation of the boundary of the disc $\Delta(T_1,T_2)$ and let $\Omega^*$ be a cyclic order of $T_1 \cup T_2$~--~they are part of the boundary of $\Delta(T_1,T_2)$~--~such that it agrees with said orientation. By construction there exist $1_u,\ldots,s_u \in [-d',\ell]$ and $1_w,\ldots,t_w \in [-d'',p]$ such that $u_{1_u},\ldots,u_{t_u},w_{1_w},\ldots,w_{s_w}$ appear in the order listed in $\Omega^*$.
 
 Let $G^*$ be the crop of $G$ by $\Delta(T_1,T_2)$ and let $\mathsf{gap}_i \coloneqq \mathsf{gap}_i^* \cap (T_1 \cup T_2)$ for $i =1,2$. Then $(G^*,\Omega^*),$ $T_1,$ $T_2,$ $\mathsf{gap}_1,$ $\mathsf{gap}_2,$ and $\SSS(s)$ satisfy \cref{lem:scheissortiereneh}. Thus we derive the existence of disjoint intervals $F_1 \subseteq [1_u,s_u]$ and $F_2 \subseteq [1_w,t_w]$~--~note that $T_1 = \{u_{1_u},\ldots,u_{t_u}\}$ and $T_2 = \{w_{1_w},\ldots,w_{s_w}\}$~--~of size at least $f_{\ref{thm:the_vortex_thing_refined}}(k,d)+1$ such that 
 \begin{enumerate}
 \item $\{u_i \mid i \in F_1 \cup F_2\} \cap (\mathsf{gap}_1 \cup \mathsf{gap}_2) = \emptyset,$ and
 \item there is $\SSS' \subseteq \SSS(S)$ of size $f_{\ref{thm:the_vortex_thing_refined}}(k,d)+1$ such that every path in $\SSS'$ has one endpoint in~$\{u_i \mid i \in F_1 \}$ and the other in $\{u_i \mid i \in F_2\}$.
 \end{enumerate}

Finally, $i)$ implies that $F_1 \subseteq [\ell]$ and $ii)$ implies that $F_2 \subseteq [p]$
 
 Thus $(G,T),\LLL,c_1,c_2,\SSS',\alpha_1,\alpha_2,F_1,F_2$ satisfy the assumptions of \cref{thm:the_vortex_thing_refined}, and hence~$\Abs{\SSS'}\leq f_{\ref{thm:the_vortex_thing_refined}}(k,d)$ as a contradiction to $ii)$. This concludes the proof.
\end{proof}

\subsection{Embedding the instance}
\label{sec:mygifttoMazoit}

We are ready to prove \cref{thm:no_apices_vital_reduction} by combining \cref{thm:rendition_to_skeleton,thm:simple_loops_not_deep,thm:links_not_many}.

For simplicity we extend the definition of vitality of linkages to sets of internally disjoint paths as follows. Let $G$ be a graph and $\SSS$ a set of internally disjoint paths in $G$. Recall that $\tau(\SSS) \coloneqq \{ (u,v) \mid \text{a path in $\SSS$ had endpoints $u$ and $v$} \}$ is a multiset. Then we call $\SSS$ \emph{vital} in $G$ if for every set of internally disjoint paths $\SSS^*$ with $\tau(\SSS^*) = \tau(\SSS)$ it holds $\SSS^* = \SSS$. 

\begin{proof}[Proof of \cref{thm:no_apices_vital_reduction}]
 Let $t,\xi,g,k\in \N$ and $d\geq 1$. Let $\Sigma$ be a surface of genus $g$. Let $\mathsf{cycles}(k,d,t) \coloneqq f_{\ref{thm:simple_loops_not_deep}}(k,d+1) +t + 2,$ then $\mathsf{cycles}(k,d,t) \in \mathbf{poly}(k+t)\cdot 2^{\mathbf{poly}(d)}$ by definition\footnote{In fact it is linear in $t$.}.
 
 Let $(G,T)$ be a $k$-annotated graph such that
 \begin{enumerate}[label=(\roman*)]
 \item there exists a $\mathsf{cycles(k,d,t)}$-locus $(\rho,M)$ of $(G,T)$ in $\Sigma$ such that $\rho$ has at most $\xi$ vortices, and is $d$-stuffed with gap at most $d,$ and
 \item there exists a vital $T$-linkage $\LLL$ in $G$.
 \end{enumerate}
In a first step we ``exhaust the instance'' using \cref{thm:rendition_to_skeleton}, in particular we derive the existence of a $k$-annotated graph $(G',T)$ such that $G'$ is a minor of $G$ together with 
\begin{enumerate}[label=(\arabic*)]
 \item a vital $T$-linkage $\LLL'$ with $\tau(\LLL') = \tau(\LLL),$
 \item a $(\mathsf{cycles(k,d,t)}-1)$-locus $(\rho',M')$ of $(G',T)$ in $\Sigma$ where $M'$ is tight and $\rho'$ is a $\Sigma$-skeleton of $G',$
 \item Let $C^*\subseteq C(\rho)$ be the set of vortices of $\rho,$ then $C^*$ is also the set vortices of $\rho'$ satisfying $N_{\rho'}(c) = N_{\rho}(c)$ for all $c \in C^*,$ $\rho'$ is properly $d$-linked with gap at most $d,$ $\bigcup_{c\in C^*} \sigma_{\rho'}(c) \cap \bigcup M' = \emptyset,$ and
 \item $(G',T),\LLL',(\rho',M')$ is exhausted.
\end{enumerate}
Let $\omega \coloneqq \mathsf{cycles(k,d,t)}-1$ and let $M' \coloneqq \langle C_1,\ldots,C_\omega\rangle$ be the respective concentric cycles. Let~$\SSS(\LLL')$ denote the segments of $\LLL'$.
Combining $(1)-(4)$ and $\omega >f_{\ref{thm:simple_loops_not_deep}}(k,d+1)+t,$ \cref{thm:simple_loops_not_deep} implies that 
\begin{itemize}
 \item[$(\star)$] for every simple loop $S$ of $\SSS(\LLL')$ and for every $i \in [t]$ it holds $V(S) \cap V(C_i) = \emptyset$.
\end{itemize}

Let $\mathsf{fat}(k,d) \coloneqq f_{\ref{thm:links_not_many}}(k,d+1),$ then $\mathsf{fat}(k,d) \in \mathbf{poly}(k)\cdot 2^{\mathbf{poly}(d)} $. Let $S$ be a segment of $\LLL'$ that is not a simple loop, then we define $\SSS(S)\subseteq \SSS(\LLL')$ to be the set of all segments of $\LLL'$ that are $\rho'$-homotopic to $S$. Combining $(1)-(4)$~--~the size of the locus is not needed~--~and the definition of $\mathsf{fat}(k,d)$ we derive that $\Abs{\SSS(S)}< \mathsf{fat}(k,d) $ by \cref{thm:links_not_many}. Summarizing the above we derive
\begin{itemize}
 \item[$(\star\star)$] Let $S \in \SSS(\LLL')$ be a segment of $\LLL'$ that is not a simple loop, then $\Abs{\SSS(S)}< \mathsf{fat}(k,d)$.
\end{itemize}

By \cref{lem:typecounting} we derive that there are at most $\delta(\xi,g) \coloneqq \max(4\xi + 6g - 6,3)$ different $\rho'$-homotopy types for segments of $\LLL';$ cap all boundary components of $\Sigma$ (if they exist) then delete the vortices to produce the respective boundary components, then this follows by definition of $\rho'$-homotopy.

In particular together with $(\star\star)$ we conclude the following: Let $\SSS^* \subseteq \SSS(\LLL')$ be maximal such that every segment in $\SSS^*$ is not a simple loop.

\begin{itemize}
 \item[$(\star\star\star)$] There are at most $\delta(\xi,g) \cdot \mathsf{fat}(k,d)$ segments of $\SSS(\LLL')$ that are not simple loops. In particular $\Abs{\SSS^*} \leq \delta(\xi,g) \cdot \mathsf{fat}(k,d)$.
\end{itemize}

Note that every path in $\SSS^*$ is a link or a non-simple loop. Let $E' \coloneqq E(\bigcup\SSS^*) \cup E(\bigcup_{i=1}^t C_i)$. Let~$G^*$ be obtained from $G$ by deleting all edges in $E(G) \setminus E'$ and subsequently deleting all isolated vertices. We derive the following; recall that $C^*$ are the vortices of $\rho'$.
\begin{beautifulclaim}\label{claim:vortices_are_Edgeless}
 Let $v \in V(G^*) \cap V(\bigcup_{c \in C^*} \sigma_{\rho'}(c)),$ Then $v \in V(\bigcup\SSS^*) \cap \bigcup_{c \in C^*}N(c)$. Furthermore, $E(\bigcup_{c \in C^*} \sigma_{\rho'}(c)) = \emptyset$. 
\end{beautifulclaim}
\begin{claimproof}
 To see this note that $\bigcup_{i=1}^t C_i$ is disjoint from $V(\bigcup_{c \in C^*} \sigma_{\rho}(c))$ by (3), we derive that $v$ must be adjacent to some edge in $ E(\bigcup\SSS^*) $. Since $\SSS^* \subseteq \SSS(\LLL')$ are segments, \cref{def:segments} implies that only their endpoints may be part of $\bigcup_{c \in C^*} \sigma_{\rho}(c)$ and if so they are nodes of the vortex cells. This concludes the proof of the first part; the second part follows immediately by construction and the fact that no edge of a graph of a vortex cell is part of $E'$.
\end{claimproof}

Let $T^*$ be the set of endpoints of paths in $\SSS^*,$ then $T^* \subseteq \bigcup_{c \in C^*}N(c)$ by the previous claim. Since $\rho'$ is a $\Sigma$-skeleton the graph of every non-vortex cell of $\rho'$ is isomorphic to $K_2$. Let $\rho^*$ be the rendition of $G^*$ obtained by deleting all the cells corresponding to edges of $E(G)$ that were deleted in the construction of $G^*$. Let $M^*\coloneqq\langle C_1,\ldots,C_t\rangle$. We derive the following by construction.
\begin{beautifulclaim}\label{claim:das_Ende_naht_neuer_locus}
 $(\rho^*,M^*)$ is a $t$-locus of $G^*$ in $\Sigma$ and $\rho^*$ is a $\Sigma$-skeleton.
\end{beautifulclaim}

By \cref{obs:segments} $\SSS^*$ is a set of internally disjoint paths and every endpoint of $\SSS^*$ is part of at most two paths in $\SSS^*$. We have the following.
\begin{beautifulclaim}\label{claim:das_Ende_naht_S*_vital}
 $\SSS^*$ is vital in $G^*$ and no path of $\SSS^*$ is a simple loop.
\end{beautifulclaim}
\begin{claimproof}
 If not, then there exists a set $\PPP$ of internally disjoint paths such that $\tau(\SSS^*) = \tau(\PPP)$ but $\PPP \neq \SSS^*,$ in particular we have $\Abs{\SSS^*} = \Abs{\PPP}$ since their ``patterns'' are multisets. Let $\beta: \SSS^* \to \PPP$ be a bijection such that $S \in \SSS^*$ and $\beta(S)$ have the same endpoints.

 Let $\tilde \LLL'$ be obtained from $\LLL'$ by replacing for each $L \in \LLL'$ and $S \in \SSS^*$ with $S \subseteq L,$ the subpath $S$ of $L$ by $\beta(S)$. By construction this is a linkage and it verifies $\tau(\tilde \LLL') = \tau(\LLL')$ but $\tilde \LLL' \neq \LLL'$ refuting the vitality of $\LLL'$.

 Finally no path of $\SSS^*$ is a simple loop by construction.
\end{claimproof}

Let $\mathsf{bad}(\SSS^*) \subseteq T^*$ be the maximum set consisting of all those vertices in $T^*$ that are adjacent to two paths in $\SSS^*$. Let $x \in \mathsf{bad}(\SSS^*)$. By construction of $G^*,$ $\deg(x)=2;$ let $ux,xv \in E(G^*)$ be the two edges. We construct \emph{the split of $G^*,$ $\SSS^*$ and $\rho^*$ at $x$} denoted by $G^*_x,$ $\SSS^*_x,$ $\rho^*_x$ as follows.
\begin{description}
 \item[$G^*_x$:] Let $G^*_x$ be obtained from $G^*$ by introducing two fresh vertices $x_1,x_2\notin V(G^*),$ deleting $x$ from $G^*$ and adding the edges $ux_1$ and $x_2v$. 
 \item[$\SSS^*_x$] Let $S_1,S_2 \in \SSS^*$ be the two paths that have $x$ as an endpoint; without loss of generality let $ux \in E(S_1)$ and $xv \in E(S_2)$. Let $\tilde S_1$ be obtained from $S_1$ by replacing the edge $ux$ with $ux_1$ and let $\tilde S_2$ be obtained from $S_2$ by replacing the edge $xv$ with $vx_2$. Clearly $\tilde S_1$ and $\tilde S_2$ are paths in $G^*_x$. Moreover they are disjoint, for otherwise $S_1\cup S_2$ is a cycle in $G^*$ and hence in $G',$ as opposed to $\LLL'$ being a linkage. In particular $\SSS^*_x$ is an internally disjoint linkage. 
 \item[$\rho^*_x$:] Since $x \in T^*,$ it is a node of some vortex cell $c \in C^*$. Let $c_u,c_v \in C(\rho^*)$ be the cells whose graphs contain $ux$ and $xv$ respectively; they are distinct since $\rho^*$ is a $\Sigma$-skeleton, and they are not vortices by \cref{claim:vortices_are_Edgeless}. Define $\rho^*_x$ to agree with $\rho^*$ on cells $C(\rho^*)\setminus\{c_u,c_v,c\},$ and subsequently drawing $x_1,x_2$ on the boundary of $c$~--~thus changing the vortex cell $c$ to a new cell $c_x$ with one more node~--~and slightly shift the cells $c_u,c_v$ so that they are disjoint, the first containing $x_1$ and the second $x_2$ as a node. This is clearly feasible: One way to see this is by taking a small disc $\Delta_\epsilon \subsetneq \Sigma$ centred at the node $x,$ such that the disc contains no other node of $\rho^*$. Let $\zeta_1,\zeta_2$ be the two intersection points of the boundary of $\Delta_\epsilon$ and the boundary of $c$. Finally define $\pi_{\rho^*_x}(\zeta_i) = x_i$~--~possibly renaming $\zeta_1,\zeta_2$~--~such that the curve $\gamma_{ux}$ between $u$ and $x_1$ on $\Sigma$ and the curve $\gamma_{xv}$ between $v$ and $x_2$ do not intersect. Then let $\tilde c_u',\tilde c_v'$ be disjoint discs, internally disjoint from all cells $C(\rho^*) \setminus \{c_u,c_v\},$ the first containing $\gamma_{ux}$ and intersecting it exactly at its endpoints $u,\zeta_1,$ the second containing $\gamma_{xv}$ and intersecting it exactly at its endpoints $v,\zeta_2$. Finally delete $c_u,c_v$ and let $c_u'\coloneqq \tilde c_u' \setminus \{u,\zeta_1\}$ and $c_v'\coloneqq \tilde c_v' \setminus \{v,\zeta_2\}$. Define $\sigma_{\rho_x}(c_u')$ as the graph consisting solely of the edge $ux_1$ and $\sigma_{\rho_x}(c_v')$ as the graph consisting solely of the edge $x_2v$. This concludes the construction.
\end{description}

By construction of $G^*_x,\SSS^*_x,$ and $\rho^*_x$ we derive the following.
 
\begin{beautifulclaim}\label{claim:das_Ende_naht_S*_split_vital}
 Let $T^*_x$ be the set of endpoints of $\SSS^*_x$. Then $(\rho^*_x,M^*)$ is a $t$-locus of $(G^*_x,T^*_x)$ in $\Sigma$ such that $\rho^*_x$ is a $\Sigma$-skeleton. Furthermore $\Abs{T^*_x} = \Abs{T^*} + 1,$ $\Abs{\mathsf{bad}(\SSS^*_x)} < \mathsf{bad}(\SSS^*)$ and $\SSS^*_x$ is vital in $G^*_x$ admitting no simple loop. 

 Let $P_1,P_2 \in \SSS^*$ then there exist $\tilde P_1,\tilde P_2 \in \SSS^*_x$ where $\tilde P_i = P_i$ if $P_i\notin \{S_1,S_2\}$ and $\tilde P_i = \tilde S_j$ if $P_i=S_j$ for $i,j \in \{1,2\}$ (and vice versa). Furthermore, $P_1,P_2$ are $\rho^*$-homotopic if and only if $\tilde P_1,\tilde P_2$ are $\rho^*_x$-homotopic.
\end{beautifulclaim}
\begin{claimproof}
 By construction $\rho^*_x$ is a blank and stretched $\Sigma$-skeleton of $(G^*_x,T^*_x)$. Clearly $\bigcup_{i=1}^t C_i \subseteq G^*_x$ since $T^* \cap V(\bigcup_{i=1}^t C_i) = \emptyset$. Together with \cref{claim:das_Ende_naht_neuer_locus} we derive that $(\rho^*_x,M^*)$ is a $t$-locus of $(G^*,T^*)$ in $\Sigma$.
 
 The claims $\Abs{T^*_x} = \Abs{T^*} + 1$ and $\Abs{\mathsf{bad}(\SSS^*_x)} < \mathsf{bad}(\SSS^*)$ follow by construction of $\SSS^*_x,$ since $T^*_x = (T^* \cup \{x_1,x_2\}) \setminus\{x\}$. Further $x_1,x_2 \notin \mathsf{bad}(\SSS^*_x)$ whence $\mathsf{bad}(\SSS^*_x) \subsetneq \mathsf{bad}(\SSS^*)$ deriving said claim.

 To see that $\SSS^*_x$ is vital in $G^*_x,$ note that $\SSS^*$ is vital in $G^*$ by \cref{claim:das_Ende_naht_S*_vital}. Using analogous arguments as usual we derive that if $\SSS^*_x$ is not vital, then $\SSS^*$ was not via ``gluing'' the unique paths of $\SSS^*_x$ containing $x_1,x_2$ by identifying both with $x$ in $G^*$.
 
 Furthermore note that if $\SSS^*_x$ contains a simple loop $S,$ then $\beta^{-1}(S)$ must be a simple loop. And finally note that the trace of $\tilde S_i$ is homotopic to the trace of $S_i$ by construction, deriving the last part of the claim (the first half of which follows by construction). This concludes the proof.
\end{claimproof}

Define $\mathsf{terminals}(k,\xi,d,g) \coloneqq 4 \delta(\xi,g) \cdot \mathsf{fat}(k,d),$ then by definition we have $\mathsf{terminals}(k,\xi,d,g) \in \mathbf{poly}(k+\xi+g) \cdot 2^{\mathbf{poly}(d)}$\footnote{The dependency on $\xi$ and $g$ is actually linear.}.

Iteratively splitting $G^*,\SSS^*$ and $\rho^*$ for every vertex in $\mathsf{bad}(\SSS^*)$ and applying \cref{claim:das_Ende_naht_S*_split_vital} we derive the following.

\begin{beautifulclaim}\label{claim:das_Ende_naht_FIN}
 There exists an annotated graph $(\tilde G^*,\tilde T^*)$ and $t$-locus $(\tilde \rho^*,M^*)$ of $(\tilde G^*,\tilde T^*)$ in $\Sigma$ with a vital $\tilde T^*$-linkage $ \tilde \SSS^*$ such that 
 \begin{enumerate}
 \item[$(a)$] $\tilde \rho^*$ is a $\Sigma$-skeleton, and $(G^*,\tilde T^*),\tilde \SSS^*,(\tilde \rho^*,M^*)$ is exhausted,
 \item[$(b)$] all vortices $c \in C(\tilde \rho^*)$ satisfy $\sigma(C) = (N(c),\emptyset),$
 \item[$(c)$] let $L\in \tilde \SSS^*,$ then $L$ is either a non-simple loop or a link,
 \item[$(d)$] let $P \in \tilde\SSS^*$ be a non-simple loop or a link in $G^*,$ and let $\tilde\SSS^*(P) \subseteq \tilde\SSS^*$ be the maximum size linkage of paths $\rho$-homotopic to $P$. Then $\Abs{\tilde\SSS^*(P)}\leq \mathsf{fat}(k,d)$.
 \item[$(e)$] $\Abs{\tilde T^*} \leq \mathsf{terminals}(k,\xi,d,g)$
 \end{enumerate}
\end{beautifulclaim}
\begin{claimproof}
 Let $\tilde G^*,\tilde\SSS^*,\tilde \rho^*$ be obtained by iteratively splitting $G^*,\SSS^*,\rho^*$ at the vertices $\mathsf{bad}(\SSS^*)$. Then $(b)$ is clear by construction, since there are no edges left in the graphs of vortices, and we deleted all isolated vertices when constructing $G^*$.
 
 By \cref{claim:das_Ende_naht_S*_split_vital} we derive that $\tilde\SSS^*$ is a vital linkage in $\tilde G^*$ (note that there are no ``bad'' vertices left, i.\@e.\@, $\mathsf{bad}(\tilde\SSS^*) = \emptyset$). Since there were no simple loops in $\SSS^*,$ and splitting never produced new simple loops by \cref{claim:das_Ende_naht_S*_split_vital}, we derive that $\tilde\SSS^*$ admits no simple loop, in particular we derive $(c)$. 
 
 Let $\tilde T^*$ be the set of endpoints of paths in $\tilde\SSS^*,$ again by repeated application of \cref{claim:das_Ende_naht_S*_split_vital} we derive that $(\tilde \rho^*,M^*)$ is a $t$-locus of $(\tilde G^*,\tilde T^*)$ in $\Sigma$ and $\tilde \rho^*$ is a $\Sigma$-skeleton. Altogether we derive $(a),$ where exhaustedness follows by construction. 

 To see $(d)$ note that \cref{claim:das_Ende_naht_S*_split_vital} implies that if two paths of $\SSS^*$ are $\rho^*$-homotopic, then the unique resulting paths following the splitting construction remain $\tilde \rho^*$ homotopic.

 Finally $\Abs{ \tilde T^*} \leq 2 \Abs{T^*}$ and $\Abs{T^*} \leq 2 \Abs{\SSS^*}$ which by $(\star\star\star)$ implies $$\Abs{\tilde T^*}\leq 4\delta(\xi,g) \cdot \mathsf{fat}(k,d).$$ In particular $\Abs{\tilde T^*} \leq \mathsf{terminals}(k,\xi,d,g)$ concluding $(e)$.

 Altogether this concludes the proof of the claim.
\end{claimproof}

Finally the instance resulting from \cref{claim:das_Ende_naht_FIN} is as desired, concluding the proof.
\end{proof}

We have the following corollary; given a surface $\Sigma$ possibly with boundary we define $\hat \Sigma$ to be the surface without boundary obtained from $\Sigma$ by capping all the holes of $\Sigma$ with discs. 
\begin{corollary}\label{cor:noapexreducetomazoit}
 There exist two functions $\mathsf{terminals}_{\ref{cor:noapexreducetomazoit}}\colon \N^4\to \N$ and $\mathsf{cycles}_{\ref{cor:noapexreducetomazoit}}\colon \N^3\to \N$ such that the following holds.

 Let $t,\xi,g,k\in \N$ and $d\geq 1$.
 Let $\Sigma$ be a surface of genus $g$ and let $(G,T)$ be a $k$-annotated graph such that 
 \begin{itemize}
 \item there exists a $\mathsf{cycles}_{\ref{cor:noapexreducetomazoit}}(k,d,t)$-locus $(\rho,M)$ of $(G,T)$ in $\Sigma$ such that $\rho$ is $d$-stuffed with gap at most $d$ and has at most $\xi$ vortices, and
 
 \item there exists a vital $T$-linkage in $G$.
 \end{itemize}
 Then there exists a surface $\Sigma^*,$ a $\mathsf{terminals}_{\ref{cor:noapexreducetomazoit}}(k,\xi,d,g)$-annotated graph $(G',T'),$ a vital $T'$-linkage $\LLL',$ a disc $\Delta \subseteq \Sigma^*$ disjoint from $\mathsf{bd}(\Sigma^*),$ a sequence of cycles $\CCC = \langle C_1, \ldots, C_t\rangle$ and an embedding $\psi : V(G)\cup E(G) \to \Sigma^*$ with $\psi(C_i) \subseteq \Delta$ for every $i \in [t]$ such that the cycles $\CCC$ are concentric with respect to $\psi,$ satisfying the following:
 \begin{enumerate}
 \item $\hat \Sigma^* \cong \hat \Sigma$ and $\Sigma^*$ has at most $\xi$ boundary components,
 \item $G' = \bigcup \mathcal{C} \cup \mathcal{L}',$
 \item $\psi(T')$ is contained in $\mathsf{bd}(\Sigma^*),$ and
 \item let $L\in \LLL',$ then $L$ is internally disjoint from $\mathsf{bd}(\Sigma^*)$.
 \end{enumerate}

 Moreover, $\mathsf{terminals}_{\ref{cor:noapexreducetomazoit}}(k,\xi,d,g)\in {\mathbf{poly}(k+\xi+g)}\cdot 2^{\mathbf{poly}(d)}$ and $\mathsf{cycles}_{\ref{cor:noapexreducetomazoit}}(k,d,t) \in {\mathbf{poly}(k+t)}\cdot2^{\mathbf{poly}(d)}$. 
\end{corollary}

\section{Eliminating apices and building nests}\label{sec:noapexnocry}
While \zcref{thm:localstructure} already guarantees both a large wall grounded in the rendition together with private nests for each vortex, those nests are not large enough for us to apply \zcref{thm:linking_vortices}.
This is fundamentally a consequence of the way \zcref{thm:localstructure}~--~and more generally the main results in \cite{GorskySW2025Polynomial}~--~is proven and how the vortex nests are created.
The goal of this section is to show that we are able to surround the vortices with disjoint nests of essentially arbitrary size, with the depth of their vortices being our goal size.

The main tools in this section are inspired by results due to Diestel, Kawarabayashi, M{\"u}ller, and Wollan \cite{DiestelKMW2012Excluded}.
Since their results produce several annoying problems for us, which we will explain in detail in a moment, we have to do quite a bit of work to fix them up for our purposes.

\paragraph{Representativity and distance in a rendition.}
We start by generalising the notion of representativity to renditions in a natural way and introduce a measure for distance between nodes of a rendition.

Let $\Sigma$ be a surface and let $\rho$ be a $\Sigma$-rendition of a graph $G$.
Then we let $G_\rho$ be the graph on the vertex set $N(\rho)$ and with the following edges
\begin{itemize}
 \item for each non-vortex cell $c \in C(\rho),$ add all edges between distinct vertices in $N(c),$ and

 \item for each vortex, with $(H,\Omega)$ being its society and with $\Omega = \langle v_1, \ldots, v_k \rangle,$ we introduce the edges $v_iv_{i+1}$ for all $i \in [k-1]$ and the edge $v_kv_1$.\footnote{We delete all parallel edges after this process to turn $G_\rho$ into a simple graph.}
\end{itemize}
We then define an embedding $\psi_\rho$ of $G_\rho$ on $\Sigma$ by drawing $N(\rho)$ onto $\Sigma$ exactly as $\rho$ does and for each edge $uv \in E(G_\rho)$ choosing a cell $c \in C(\rho)$ such that $c \in N(c),$ letting $T$ be the closure of a component of $\mathsf{bd}(c) - \{ u, v \}$ for which $T$ is disjoint from $N(c) \setminus \{ u, v \}$ but contains both $u$ and $v,$ and setting $\psi(uv) = T$.
The \emph{representativity of $\rho$} is then simply the representativity of $\psi_\rho$.

Furthermore, given two distinct vertices $u,v \in N(\rho),$ the \emph{distance between $u$ and $v$ in $\rho$} is the minimum number of nodes of $\rho$ that a $\rho$-aligned curve ending in $u$ and $v$ contains.
Given two vortices $c_1,c_2$ with vortex societies $(H_1,\Omega_1),(H_2,\Omega_2),$ the \emph{distance between $c_1$ and $c_2$ in $\rho$} is the minimum distance between a vertex in $V(\Omega_1)$ and a vertex in $V(\Omega_2)$ in $\rho$.

Given an embedding $\psi_\rho$ derived from a rendition $\rho$ of a graph $G$ in a surface $\Sigma,$ we want to translate a curve $T$ in $\Sigma$ that meets $\psi_\rho$ only in the vertices of $G$ into what should be a $\rho$-aligned curve.
It is easy to see that $T$ can be shifted in $\Sigma$ such that it meets none of the cells of $\rho$.
However, if $T$ passes through the interior of a vortex $c$ of $\rho,$ there may be no way to simply shift $T$ to avoid $c$.
Nonetheless, we want to use such a curve $T$ to cut apart our surface in this section.

We call a curve in $\Sigma$ that meets $\rho$ only in nodes and the interior of vortices of $\rho$ a \emph{$\rho$-conformal curve}.
The following is immediate from our discussion.

\begin{observation}\label{obs:alginedvsconformal}
 Given an embedding $\psi_\rho$ of $G$ with representativity at most $r$ derived from a rendition $\rho$ of $G$ in a surface $\Sigma,$ then for any curve $T$ that meets $\psi_\rho$ only in $V(G)$ there exists a $\rho$-conformal curve $T'$ such that $T \cap V(G) = T' \cap V(G),$ the curve $T'$ intersects none of the non-vortex cells of $\rho,$ and $T'$ intersects a vortex of $\rho$ if and only if $T$ intersects said vortex.
\end{observation}

We will need a small technical result on curves witnessing low representativity from \cite{DiestelKMW2012Excluded}.
For a graph $G$ with an embedding $\psi$ in a surface $\Sigma,$ a face $f$ of $\psi,$ and a closed curve $C,$ let $\mathsf{comp}(C,f)$ denote the number of components of $C \cap f$.

\begin{proposition}[Diestel et al.\ \cite{DiestelKMW2012Excluded}]\label{prop:shortgenusreducingcurves}
 Let $G$ be a graph with an embedding $\psi$ in a surface $\Sigma$ of positive genus and let $F$ be the set of faces of $\psi$.
 For a positive integer $r,$ consider all genus-reducing curves $C$ that meet $\psi$ only in $\psi(V(G))$ and such that $|\psi(V(G)) \cap C| < r$.
 Let $C$ be chosen such that $\sum_{f \in F} \mathsf{comp}(C,f)$ is minimal.

 Then $\mathsf{comp}(C,f) \leq 1$ for all $f \in F$.
\end{proposition}

What we would love to do at this point is to simply import the Lemmas 14 and 16 from \cite{DiestelKMW2012Excluded} which first fuse together vortices that are close together and then increase the representativity of the surface our rendition is working in.
However, both of their lemmas introduce vertices that must be deleted~--~commonly called \emph{apices}~--~which have neighbours outside of the vortices of the established rendition.
Such apices are called \emph{major apices}.\footnote{Notably we are being a bit vague here by what we mean by ``neighbours outside of vortices''. Our requirements are actually slightly more lax than what is usually required. In the broader context of Graph Minor Structure Theory an apex is considered \textsl{major} if it sees any node of the rendition or has a neighbour in a non-vortex cell of said rendition. In our case, we will actually allow apices that have neighbours on the boundary of vortices, which are nodes of the rendition.}
This causes immense problems for us, because we desperately need all newly introduced apices to only have neighbours in vortices.

As such, we will have to go through a much more considerable effort than Diestel et al.\ to establish our versions of their lemmas.
In particular, we make the following adjustments:
We will only work on graphs which already host a linkage whose vitality we are trying to preserve.
In the course of both of our modifications, we will actually change the graph, the rendition, and the linkage we are working on to produce a new instance with more convenient properties for us.
The vortices we create throughout will be more complicated than what Diestel et al.\ find, but in exchange we will have more control over the apices that get created.

\subsection{Joining close vortices}
Our first goal will be to join together vortices in a rendition which are not very far apart.
This is dealt with in Lemma 14 of \cite{DiestelKMW2012Excluded}, where the nodes on the curve that shows that the two vortices are not far from each other are simply taken to be apices.
For us this is not an option and we will have to make some more involved modifications to our graph.

\begin{lemma}\label{lem:joinvortices}
 Let $b,d,r$ be positive integers.
 Let $\Sigma$ be a surface, let $\rho$ be a $\Sigma$-rendition with depth $d$ and breadth $b$ of a graph $G,$ let $\mathcal{L}$ be a linkage in $G,$ and let $c_1,c_2$ be two distinct vortices in $\rho$ such that there exists a $\rho$-aligned $\mathsf{bd}(c_1)$-$\mathsf{bd}(c_2)$-curve $T$ that is otherwise disjoint from all vortices of $\rho$ such that $|T \cap N(\rho)| \leq r$.

 Then there exists a set $A \subseteq V(G)$ with $|A| \leq 2d,$ a set $B$ with $|B| = 2|T \cap N(\rho)|,$ a graph $G',$ a linkage $\mathcal{L}'$ in $G',$ a $\Sigma$-rendition $\rho'$ of $G'$ with breadth $b-1,$ and a $\rho'$-aligned disc $\Delta$ such that
 \begin{enumerate}
 \item $V(G') = (V(G) \cup B) \setminus (A \cup (T \cap N(\rho))),$

 \item $c_1 \cup c_2 \subseteq \Delta,$

 \item $\rho$ and $\rho'$ agree on $G - (A \cup (T \cap N(\rho)) \cup \{ \sigma_\rho(c) - N_\rho(c) ~\colon~ c \in C(\rho) \text{ and } (c \setminus \bd(c)) \cap \Delta \neq \emptyset \}),$
 
 \item $\Delta$ is a vortex of $\rho'$ with the vortex society $(H,\Omega)$ of depth at most 2 plus the sum of the depths of the vortex societies of $c_1$ and $c_2,$

 \item $|\mathcal{L}'| \leq |\mathcal{L}| + r,$
 
 \item if $S$ are the endpoints of $\mathcal{L},$ then the endpoints of $\mathcal{L}'$ lie in $S \cup V(H),$ and

 \item if $\mathcal{L}$ is vital in $G,$ then $\mathcal{L}'$ is vital in $G'$.
 \end{enumerate}
\end{lemma}
\begin{proof}
 We let $(G_1, \Omega_1)$ be the $c_1$-society and $(G_2, \Omega_2)$ be the $c_2$-society both with depth at most $d,$ where $\Omega_1 = \langle v_1,\ldots,v_n\rangle$ and $\Omega_2 = \langle w_1, \ldots, w_m\rangle$.
 By slightly adjusting $T$ we may assume that the endpoints of $T$ are vertices $v_h$ and $w_\ell,$ so that $T \cap c_1 = \{v_h\}$ and $T \cap c_2 = \{w_{\ell}\}$ for some indices $k \in [n]$ and $\ell \in [m]$. 
 Let $T' = N(\rho) \cap T$ and let $\psi = \langle v_h = x_1, x_2, \ldots, x_{r'} = w_\ell\rangle$ be the linear ordering of the vertices in $T'$ obtained by traversing $T$ starting from $v_h$ towards $w_\ell$. 

 We now fatten $T$ slightly to a disc $D$ such that we maintain that $c \not\subseteq D$ for all $c \in C(\rho)$ and $\mathsf{bd}(D) \setminus (c_1 \cup c_2)$ has exactly two components $Y_1,Y_2$.
 Then we obtain a closed disc $\Delta \coloneqq \mathsf{cl}( D \cup c_1 \cup c_2 )$ such that $\Delta \cap N(\rho) \setminus T' = ( V(\Omega_1) \cup V(\Omega_2) ) \setminus \{ v_k, w_\ell \}$.
 By reindexing if necessary, we may assume that the orientations of $\mathsf{bd}(\Delta)$ induced by $V(\Omega_1) \setminus \{ v_k \}$ and $V(\Omega_2) \setminus \{ w_\ell \}$ agree.

 At this point we must deviate from the proof of Lemma 14 in \cite{DiestelKMW2012Excluded} as they now use $\Delta$ as the new vortex and $T'$ as part of the apices.
 This however will mean that $T'$ has many neighbours outside of vortices, namely all neighbouring nodes.
 To avoid this, we will now copy $T'$ twice and place the copies onto $Y_1$ and $Y_2,$ adjusting $\rho$ mildly to fit this new graph.

 Let $T_1,T_2$ be two disjoint copies of $T',$ let $\pi_i \colon T_i \rightarrow T'$ be a bijection identifying these copies for each $i \in [2],$ and let $B = T_1 \cup T_2$ for later reference.
 Further, let $\Gamma$ be the painting in $\rho$.
 We adjust $\Gamma$ and $G$ as follows for both $i \in [2]$.
 The nodes corresponding to $T'$ are removed from $\Gamma$ and replaced by $T_1 \cup T_2$.
 Similarly, to define $G^*$ we first remove $T'$ and then add $T_1 \cup T_2,$ with additional edges to be determined.
 We then place $T_i$ onto $Y_i$ according to $\pi_i$ in the order prescribed by $\psi$.
 For each $c \in C(\rho) \setminus \{ c_1, c_2 \}$ and $i \in [2]$ such that $c \cap T_i \neq \emptyset$ we note that $c \cap T_{3-i} = \emptyset$ and accordingly, we may deform $c$ slightly such that its closure intersects $T_i$ only in $\pi_i^{-1}(N_\rho(c))$.
 Furthermore, for each $v \in \sigma_\rho(c)$ such that there exists a $u \in N_\rho(c) \cap T',$ we introduce the edge $v\pi_i^{-1}(u)$ to $G^*$.
 This finalises our definition of $G^*$ and we can now also add $\Delta$ as a new vortex to $\Gamma$.
 The resulting painting lets us derive a rendition $\rho^*$ of $G^*$ in $\Sigma$ for which $V(\sigma_{\rho^*}(\Delta)) = V(\sigma_\rho(c_1) \cup \sigma_\rho(c_2)) \cup T_1 \cup T_2$.
 In particular, we make our choices such that $v_{h-1}, v_h, \pi_1^{-1}(x_2)$ (where we let $v_{h-1}$ be $v_h$ if $h=1$) are seen in this order on $\mathsf{bd}(\Delta)$ without interruptions by other nodes in $N(\rho^*)$.

 We then modify the linkage $\mathcal{L}$ in $G$ into a linkage $\mathcal{L}^*$ in $G^*$ as follows (see \zcref{fig:pathsplitting} for an illustration):
 If $V(P) \cap T' = \emptyset,$ we add $P$ to $\mathcal{L}^*$.
 Suppose instead that $P$ contains a vertex $v \in T'$.
 Then let $P_v$ be the component of $P[T']$ that contains $v$.
 \begin{itemize}
 \item If $P_v = P$ we remove $P$ from $\mathcal{L}$ and add the paths $P_1,P_2$ derived by replacing each $x \in V(P)$ with $\pi_i^{-1}(x)$ for each $i \in [2]$.

 \item Should we have $P_v \neq P,$ but $P_v$ still contains an endpoint of $P,$ let $u$ be the other endpoint of $P_v$ with $uu' \in E(P) \setminus E(P_v)$.
 Then there exists a unique $i \in [2]$ such that $\pi_i^{-1}(u)u' \in E(G^*)$ and given this $i,$ we replace all vertices in $P_v$ with their copies from $T_i$ and repeat the process we are outlining here on the new path until it is disjoint from $T'$.
 Furthermore, we add the path derived from $P_v$ by replacing all vertices in $P_v$ with their copies from $T_{3-i}$ to $\mathcal{L}^*$.

 \item Finally, let $u$ and $w$ be the endpoints of $P_v$ with $uu', ww' \in E(P) \setminus E(P_v)$.
 If there exists an $i \in [2]$ such that $\pi_i^{-1}(u)u', \pi_i^{-1}(w)w' \in E(G^*),$ replace all vertices in $P_v$ with their copies from $T_i$ and repeat this process on the resulting path, and add the path derived from $P_v$ by replacing all of its vertices with their copies from $T_{3-i}$ to $\mathcal{L}^*$.
 
 Otherwise, let $i,j \in [2]$ be distinct such that $\pi_i^{-1}(u)u', \pi_j^{-1}(w)w' \in E(G^*)$.
 We then split $P$ into two paths $P_u,P_w$ with $P_u \cap P_w = P_v$ such that $u$ is an endpoint of $P_u$ and $w$ is an endpoint of $P_w$.
 Afterwards, we replace all vertices of $P_v$ in $P_u$ by their copies from $T_i$ and proceed analogously from $P_w$.
 The process is then repeated on $P_u$ and $P_w$.
 \end{itemize}
 Though the definition is somewhat involved, it is not hard to see that if $\mathcal{L}$ was vital in $G,$ then $\mathcal{L}^*$ is vital in $G^*$.
 In particular, this increases the order of $\mathcal{L}^*$ by at most $r$ compared to $\mathcal{L}$.

 \begin{figure}[ht]
 \centering
 \scalebox{0.8}{
 \begin{tikzpicture}[scale=1.25]

 \pgfdeclarelayer{background}
		 \pgfdeclarelayer{foreground}
			
		 \pgfsetlayers{background,main,foreground}

 \begin{pgfonlayer}{background}
 \pgftext{\includegraphics[width=16cm]{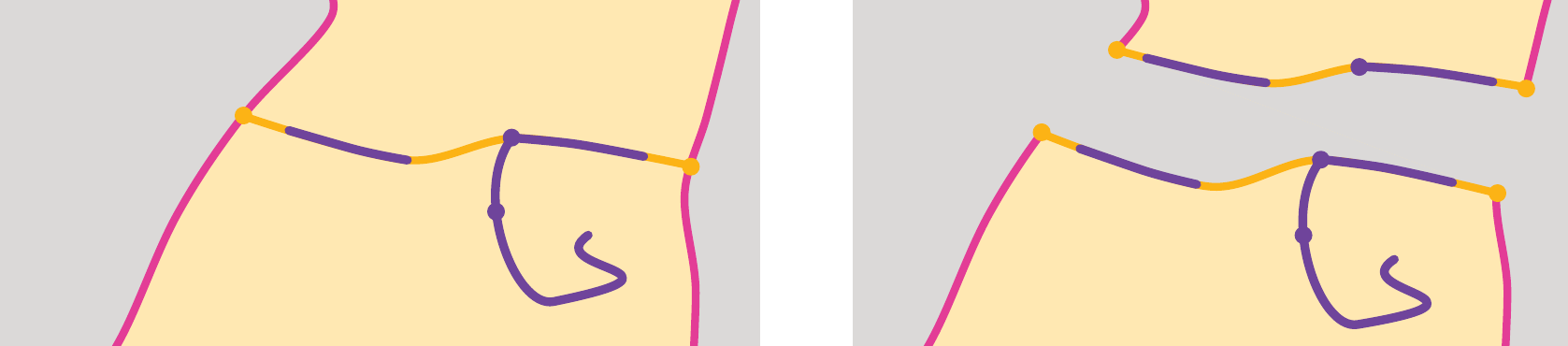}};
 \end{pgfonlayer}{background}
			
 \begin{pgfonlayer}{foreground}

 \node (T) at (-5.35,0.3) [draw=none] {$T$};
 \node (VH) at (-5.75,0.7) [draw=none] {$v_h$};
 \node (P) at (-4.5,0) [draw=none] {$P$};
 \node (UB) at (-2.8,0.575) [draw=none] {$u$};
 \node (UPRIMEB) at (-3.125,-0.3) [draw=none] {$u'$};
 \node (WL) at (-0.65,0) [draw=none] {$w_\ell$};
 \node (A1) at (-4.4,0.55) [draw=none] {a)};
 \node (B1) at (-2,0.525) [draw=none] {b)};

 \node (T1) at (4.8,-0.325) [draw=none] {$T_1$};
 \node (T2) at (5.35,1.2) [draw=none] {$T_2$};
 \node (P1) at (3.45,-0.15) [draw=none] {$P_1$};
 \node (P2) at (4.3,1.25) [draw=none] {$P_2$};
 \node (PIVH1) at (2,0.5) [draw=none] {$\pi^{-1}_1(v_h)$};
 \node (PIVH2) at (2.75,1.3) [draw=none] {$\pi^{-1}_2(v_h)$};
 \node (UPRIMEBAGAIN) at (5.1,-0.55) [draw=none] {$u'$};
 \node (PIUB1) at (5.475,0.375) [draw=none] {$\pi^{-1}_1(u)$};
 \node (PIWL1) at (7.5,0.05) [draw=none] {$\pi^{-1}_1(w_\ell)$};
 \node (PIWL2) at (7.475,0.6) [draw=none] {$\pi^{-1}_2(w_\ell)$};
 \node (A2) at (3.85,0.5) [draw=none] {a)};
 \node (B2) at (6.25,0.6) [draw=none] {b)};

 \end{pgfonlayer}{foreground}
 \end{tikzpicture}
 }\vspace{0.25cm}
 \scalebox{0.8}{
 \begin{tikzpicture}[scale=1.25]

 \pgfdeclarelayer{background}
		 \pgfdeclarelayer{foreground}
			
		 \pgfsetlayers{background,main,foreground}

 \begin{pgfonlayer}{background}
 \pgftext{\includegraphics[width=16cm]{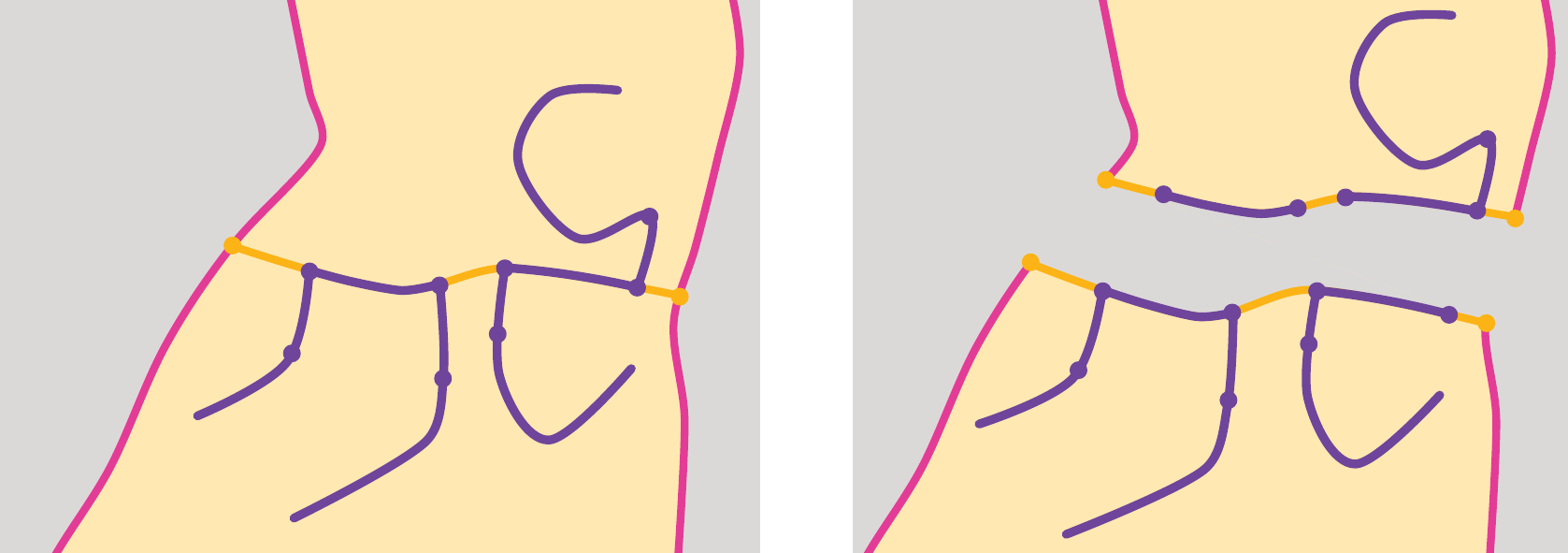}};
 \end{pgfonlayer}{background}
			
 \begin{pgfonlayer}{foreground}

 \node (UC) at (-4.8,0.25) [draw=none] {$u$};
 \node (UPRIMEC) at (-5.2,-0.7) [draw=none] {$u'$};
 \node (WC) at (-3.5,0.125) [draw=none] {$w$};
 \node (WPRIMEC) at (-3.7,-1) [draw=none] {$w'$};
 \node (UD) at (-2.85,0.3) [draw=none] {$u$};
 \node (UPRIMED) at (-2.675,-0.55) [draw=none] {$u'$};
 \node (WD) at (-1.5,-0.325) [draw=none] {$w$};
 \node (WPRIMED) at (-1.325,0.85) [draw=none] {$w'$};
 \node (C1) at (-4.15,0.15) [draw=none] {c)};
 \node (D1) at (-1.65,0.2) [draw=none] {d)};

 \node (PIUC1) at (3.35,0.15) [draw=none] {$\pi^{-1}_1(u)$};
 \node (UCPRIMEAGAIN) at (3.1,-1.1) [draw=none] {$u'$};
 \node (PIWC1) at (4.375,-0.075) [draw=none] {$\pi^{-1}_1(w)$};
 \node (WCPRIMEAGAIN) at (4.3,-1.25) [draw=none] {$w'$};
 \node (PIUD1) at (5.5,0.1) [draw=none] {$\pi^{-1}_1(u)$};
 \node (UDPRIMEAGAIN) at (5.6,-0.6) [draw=none] {$u'$};
 \node (PIWD2) at (7.15,0.375) [draw=none] {$\pi^{-1}_2(w)$};
 \node (WDPRIMEAGAIN) at (7.25,1.65) [draw=none] {$w'$};
 \node (C2) at (4.1,0.475) [draw=none] {c)};
 \node (D2) at (6.225,0.4) [draw=none] {d)};
 
 \end{pgfonlayer}{foreground}
 \end{tikzpicture}
 }
 \caption{An illustration of how to split paths of a vital linkage at a curve. Following our arguments in the proof of \zcref{lem:joinvortices}, the part of the figure labelled a) shows how to proceed if $P_v = P$. The construction shown with the label b) corresponds to the case in which exactly one endpoint of $P$ is found in $P_v$. The remaining two cases on the bottom show how to proceed if neither endpoint of $P$ is found in $P_v$. Note that $P$ may intersect $T$ in several disjoint segments, which is not depicted here and is the cause for us having to repeat the splitting process.}
 \label{fig:pathsplitting}
\end{figure}
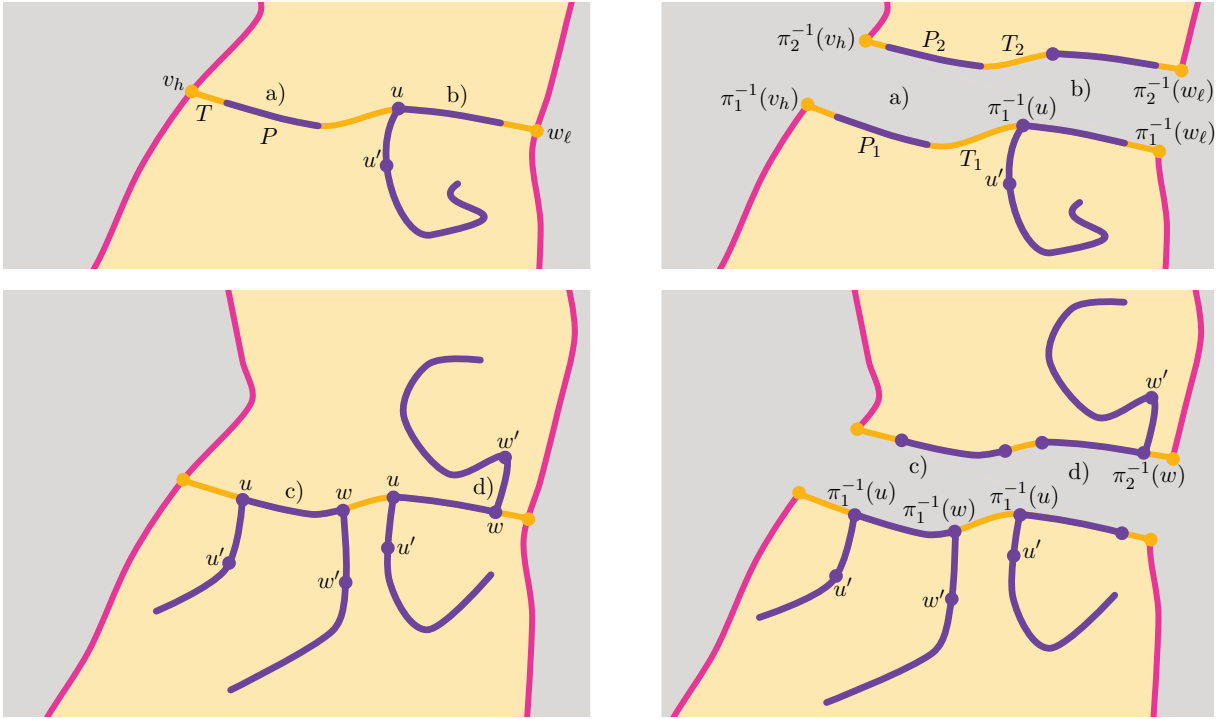

 To check the depth of $\rho^*$ we first note that for all vortices that $\rho$ and $\rho^*$ share, which is every vortex of $\rho^*$ except for $\Delta,$ the depth cannot have changed since we did not modify that part of the graph.
 Thus we now must describe how to derive an appropriate society for $\Delta$ from the information we have and explain with it has bounded depth.

 The graph $\sigma_{\rho^*}(\Delta)$ is derived by taking $\sigma_\rho(c_1) \cup \sigma_\rho(c_2),$ replacing $v_h,v_\ell$ with its two copies from $T_1 \cup T_2$ whilst retaining the original neighbourhood of $v_h$ and $v_\ell,$ and then adding $(T_1 \cup T_2) \setminus \{ \pi_1^{-1}(v_h), \pi_2^{-1}(v_h), \pi_1^{-1}(v_\ell), \pi_2^{-1}(v_\ell) \}$ as isolated vertices. 
 The cyclic ordering defining the corresponding society is then simply
 \begin{align*}
 \Omega' &\coloneqq \langle v_{h+1}, \ldots, v_n,v_1, \ldots, v_{h-1},\pi_1^{-1}(v_h),x_2^1, \ldots, x_{r'-1}^1,\pi_1^{-1}(w_\ell), \\
 &\qquad w_{\ell+1}, \ldots, w_m,w_1, \ldots, w_{\ell-1}, \pi_2^{-1}(w_\ell),x_{r'-1}^2, \ldots, x_2^2,\pi_2^{-1}(v_h)\rangle.
 \end{align*}
 Since we effectively only added a single vertex to both $\sigma_\rho(c_1)$ and $\sigma_\rho(c_2)$ it is easy to see that the depth of $(\sigma_\rho(\Delta), \Omega')$ is at most the sum of the depths of $c_1$ and $c_2$ plus 2, as the largest transaction in $(\sigma_\rho(\Delta), \Omega')$ at worst consists of a large transaction in $(G_1,\Omega_1),$ another large transaction in $(G_2,\Omega_2),$ and two extra paths that may appear due to us copying $v_h$ and $v_\ell$.

 The first four items of our lemma are now easy to confirm according to our construction with $G^*$ being the final graph we return together with the rendition $\rho^*$.
 Note that our construction also immediately shows that the breadth of $\rho^*$ is one less than the breadth of $\rho$.
 The fact that all new endpoints that $\mathcal{L}^*$ added lie in $V(\sigma_{\rho^*}(\Delta))$ is also an easy consequence of our construction.
 Finally, the last point was already explicitly guaranteed earlier during the construction of $\mathcal{L}^*$.
\end{proof}

\subsection{Pumping up the representativity}
Our second goal of this section is to pump up the representativity of the rendition we are dealing with.
This process is taken care of by Lemma 16 in \cite{DiestelKMW2012Excluded}.
As noted earlier, we will also have to modify this lemma and its proof quite a bit to make it work for our purposes.
In the statement we discuss the existence of two surfaces independent of what kind of genus-reducing curve we find, but we will only need more than one surface if the genus-reducing curve is separating.
For the other cases we will let the pairs of objects be equal.

\begin{lemma}\label{lem:pumprepresentativity}
 Let $k,t,r,r',b,d,w$ with $r' \geq r$ and $w \geq 2r$ be positive integers.
 Let $\Sigma$ be a surface of positive genus, let $G$ be a graph with a $\Sigma$-rendition $\rho$ with depth $d,$ breadth $b,$ and the property that any two vortices have distance at least $r',$ and let $\mathcal{L}$ be a linkage in $G$.
 Furthermore, there exists a closed, $\rho$-aligned disc $\Delta_M \subseteq \Sigma$ and a flat $w$-mesh $M$ in $G$ such that the trace of $M$ is found in $\Delta_M$ and $\Delta_M$ is disjoint from all vortices in $\rho$.

 If the representativity of $\rho$ is less than $r,$ then there exist graphs $G_1,G_2,$ surfaces $\Sigma_1, \Sigma_2,$ linkages $\mathcal{L}_1,\mathcal{L}_2,$ a $\Sigma_1$-rendition $\rho_1$ of $G_1,$ and a $\Sigma_2$-rendition $\rho_2$ of $G_2,$ such that for both $i \in [2]$:
 \begin{enumerate}
 \item the genus of $\Sigma_i$ is less than the genus of $\Sigma,$

 \item the breadth of $\rho_i$ is at most $b+2,$

 \item the depth of $\rho_i$ is at most $d+2,$

 \item $|\mathcal{L}_i| \leq |\mathcal{L}| + r + 6d,$

 \item if $\mathcal{L}$ is vital, then $\mathcal{L}_i$ is vital, and

 \item there exists a $j \in [2],$ a $(w - 2r)$-mesh $M' \subseteq M,$ and a closed, $\rho_j$-aligned disc $\Delta' \subseteq \Sigma_j,$ such that $M'$ is a flat mesh in $\rho_j,$ the trace of $M'$ in $\rho_j$ is contained in $\Delta',$ and all vortices of $\rho_j$ are disjoint from $\Delta'$.
 \end{enumerate}
\end{lemma}
\begin{proof}
 As $\rho$ has representativity less than $r,$ we may let $T$ be a closed, $\rho$-conformal, genus-reducing curve such that $|T \cap N(\rho)|$ is minimal.
 Note that thanks to \zcref{prop:shortgenusreducingcurves}, we know that $T$ intersects each vortex of $\rho$ in at most one component.
 Furthermore, since $r' \geq r,$ we know that $T$ intersects at most one vortex $c$ of $\rho$ and if $T$ intersects $N_\rho(c),$ then either $|N_\rho(c) \cap T| = 1$ and $T$ does not intersect $c$ itself~--~it thus just ``grazes'' $c$~--~or $|N_\rho(c) \cap T| = 2$ and $T \cap c$ consists of a single component.
 We note that for the sake of this proof, we can ignore the first option.

 We now have a few cases to distinguish:
 First of all, we are interested in whether $T$ is a separating curve or a non-separating curve.
 If $T$ is non-separating, we distinguish further between the cases in which $\Sigma - T$ has two open holes and the case in which it has only one open hole.

 \paragraph{$T$ is separating.}
 We first suppose that there exists a vortex $c$ in $\rho$ such that $T \cap c$ is non-empty and thus consists of a single component.
 The other case in which $T$ intersects no vortices essentially just uses a subset of the arguments we will now present and thus we will not discuss it explicitly.
 Let $\langle X_1, \ldots, X_k\rangle$ be a linear decomposition of the vortex society $(H,\Psi)$ of $c$ with at most $d,$ with $\Psi=\langle v_1, \ldots, v_k\rangle$.
 Furthermore, let $N_\rho(c) \cap T = \{ v_\ell, v_h \},$ chosen such that $\ell < h$.

 Let $Y^\ell = X_\ell \cap X_{\ell+1}$ and, if $h \neq k$ we let $Y^h = Y_h \cap Y_{h+1}$ and otherwise $Y^h = \{ v_h \}$.
 Then there exists a separation $(A_1, A_2)$ in $G$ such that
 \begin{itemize}
 \item $A_1 \cap A_2 = (T \cap N(\rho)) \cup Y^\ell \cup Y^h,$
 \item for all cells $c' \in C(\rho) \setminus \{ c \},$ we have $\sigma_\rho(c') \subseteq G[A_1]$ or $\sigma_\rho(c') \subseteq G[A_2],$
 \item for $i \in [\ell+1, h],$ we have $X_i \subseteq A_i,$ and
 \item for $i \in [\ell] \cup [h+1,k],$ we have $X_i \subseteq A_2$.
 \end{itemize}
 
 Let $\Sigma_1,\Sigma_2$ be the result of capping the boundaries of the two components of $\Sigma - T$ with closed discs $D_1,D_2,$ chosen such that the restriction $\rho_i'$ of $\rho$ to $G[A_i] - ( V(H) \setminus N_\rho(c) )$ is a $\Sigma_i$-rendition for both $i \in [2]$.
 We then let $\Delta_i = D_i \cup (c \cap \Sigma_i)$ be a new vortex we introduce to $\rho_i'$ to create $\rho_i$ for each $i \in [2]$.
 Further, let $\langle u_1 = v_h, \ldots, u_l = v_\ell\rangle$ be an ordering of the nodes on $T$ by traversing $T$ from $v_h$ to $v_\ell$. 
 
 The vortex society of $\Delta_1$ is then simply $(H_1 \coloneqq H[\bigcup_{i = \ell+1}^h X_i] \cup (T \cap N(\rho), \emptyset), \Omega_1),$ where
 \[ \Omega_1 \coloneqq \langle v_\ell, v_{\ell+1}, \ldots, v_h, u_2, \ldots, u_{l-1}\rangle . \]
 Similarly, the vortex society of $\Delta_2$ is $(H_2 \coloneqq H[\bigcup_{i \in [\ell] \cup [h+1,k]} X_i] \cup (T \cap N(\rho), \emptyset), \Omega_2),$ where 
 \[ \Omega_2 \coloneqq \langle v_\ell, v_{\ell - 1}, \ldots, v_1, v_k, v_{k-1}, \ldots, v_h, u_2, \ldots, u_{l-1} \rangle. \]
 For both of our new vortices, we observe that their depth is at most the depth of $(H,\Psi)$.

 We have therefore already guaranteed the first, second, and third points of our statement.
 The last point is then easily verified by noting that $T$ cannot see more than the outermost $r$ vertical and horizontal paths of the mesh $M$ and thus we can choose a $(w-2r)$-submesh $M' \subseteq M$ that has the appropriate properties and will necessarily be found either in $\rho_1$ or $\rho_2$.

 Finally, we let $\mathcal{L}_1 = \mathcal{L} \cap G[A_1]$ and let $\mathcal{L}_2 = \mathcal{L} \cap G[A_2]$.
 Accordingly, due to $(A_1,A_2)$ being a separation in $G,$ $|\mathcal{L}_1| \leq |\mathcal{L}| + r + 2d$ and $|\mathcal{L}_2| \leq |\mathcal{L}| + r + 2d,$ since $|S| \leq 2d$.
 By \zcref{prop:vitalitypreservedthroughrestriction}, we can therefore also ensure the fourth and fifth item of our statement.
 
 \paragraph{$T$ is non-separating and creates two open holes.}
 We let $c,$ $\langle X_1, \ldots, X_k\rangle,$ $(H,\Psi),$ $\Psi = \langle v_1, \ldots, v_k\rangle,$ and $N_\rho(c) \cap T = \{ v_\ell, v_h \},$ with $\ell < h,$ be defined as above.

 Choose some $\epsilon$ such that the $\epsilon$-region $r_{\epsilon(p),\Sigma} \subseteq \Sigma$ around each point $p \in T$ is a disc.
 Let $B_1,B_2$ be the two boundaries of $\Sigma - T$.
 We create two disjoint copies $T_1,T_2$ of $T$ and glue $T_i$ onto $B_i$ for both $i \in [2]$ to define the new surface $\Sigma'$ such that for both $i \in [2]$ there exists a bijection $\pi_i \colon T \rightarrow T_i$ with $r_{\epsilon,\Sigma'}(p) \setminus T_i \subseteq r_{\epsilon,\Sigma}(\pi^{-1}_i(p))$ for all $p \in T_i$.
 (See \zcref{fig:cuttube} for an illustration of this procedure.)
 Note that the genus of $\Sigma'$ is less than the genus of $\Sigma$.
 Since $T$ is $\rho$-aligned and the cells of $\rho$ explicitly do not contain the nodes of $\rho,$ we know that $T$ does not touch the cells of $\rho$ and thus all cells of the painting $\Gamma$ in $\rho$ are found in $\Sigma'$.
 Let $N = T \cap N(\rho)$ be the nodes of $\rho$ found in $T$ and let $N_1,N_2$ be two disjoint copies of $N$ such that $N_i \subseteq T_i$ and $v = \pi_i^{-1}(v)$ for all $v \in N$ and $i \in [2]$.
 We then extend the drawing $\Gamma'$ by replacing $N$ with $N_1 \cup N_2$ in the natural way according to $\pi_1$ and $\pi_2$.
 This also gives a natural definition of a graph $G'$ associated with this modification.

 \begin{figure}[ht]
 \centering
 \begin{tikzpicture}
 \pgfdeclarelayer{background}
		 \pgfdeclarelayer{foreground}
			
		 \pgfsetlayers{background,main,foreground}

 \begin{pgfonlayer}{background}
 \pgftext{\includegraphics[width=14cm]{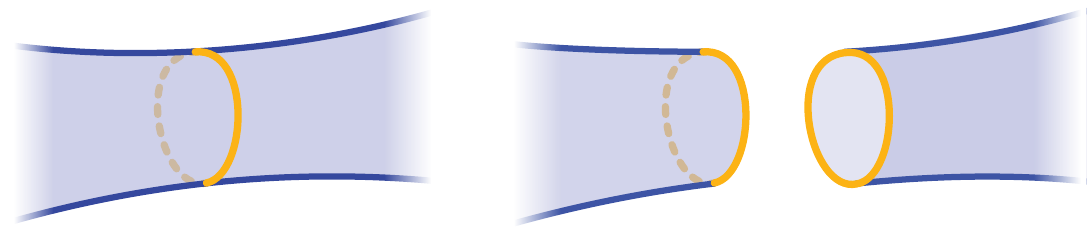}};
 \end{pgfonlayer}{background}
			
 \begin{pgfonlayer}{foreground}

 \node (T) at (-4.5,1.15) [draw=none] {$T$};
 \node (T1) at (2.05,1.15) [draw=none] {$T_1$};
 \node (T2) at (3.9,1.15) [draw=none] {$T_2$};
 
 \end{pgfonlayer}{foreground}
 \end{tikzpicture}
 \caption{An illustration of the modification of $\Sigma$ made to build two new boundaries when $T$ is a non-separating curve that creates two new boundaries upon being removed.}
 \label{fig:cuttube}
\end{figure}

 We now need to discuss how to modify $\rho$ into a rendition $\rho'$ of $G'$ on the surface $\Sigma''$ we get from $\Sigma'$ by capping the two holes bounded by $T_1$ and $T_2$ in $\Sigma'$ with disjoint open discs $\Delta_1',\Delta_2'$.
 For both $i \in [2],$ let the closure $\Delta_i$ of $\Delta_i' \cup c$ be the two new vortices which we wish to introduce.
 For all cells in $\rho$ other than $c,$ we can retain their placement and the definition of what part of the graph is assigned to them (up to the modifications $\pi_1$ and $\pi_2$ prescribe) to define $\rho'$.
 Note in particular that such cells are by definition disjoint from $N$ and thus are still found within $\Sigma'$.
 This only leaves us to define what should be placed into $\Delta_1$ and $\Delta_2$.

 We can now proceed quite similarly to the construction in our previous case and in the proof of \zcref{lem:joinvortices}.
 There exists a separation $(A_1,A_2)$ in $H$ such that $A_1 = \bigcup_{i=\ell+1}^h X_i$ and $A_2 = \bigcup_{i=1}^\ell X_i \cup \bigcup_{i=h+1}^k X_i,$ with the set $A_1 \cap A_2$ consisting of the sets $Y^\ell = X_\ell \cap X_{\ell+1}$ and if $h \neq k,$ $Y^h = X_h \cap X_{h+1}$ and $Y^h = \{ v_h \}$ otherwise.
 Further, we let $$S^1 \coloneqq (Y^\ell \cup Y^h) \setminus \{ v_i ~\colon~ i \in [\ell+,h-1] \},$$ let $$S^2 \coloneqq (Y^\ell \cup Y^h) \setminus \{ v_i ~\colon~ i \in [\ell-1] \cup [h+1,k] \},$$ and let $$S \coloneqq (Y^\ell \cup Y^h) \setminus N(\rho).$$
 We then let $F$ be the set of edges in $H$ with at least one edge incident to a vertex in $Y^\ell \cup Y^h,$ together with their counterparts translated through $\pi_1$ and $\pi_2$ if $v_h$ and $v_\ell$ are found in $Y^\ell \cup Y^h$.
 
 This allows us to define the graph $H_1$ as the graph derived from $(H[A_1 \setminus S^1] - F) \cup (\pi_1^{-1}(N) \setminus \{ \pi_1^{-1}(v_\ell), \pi_1^{-1}(v_h) \}, \emptyset)$ by replacing $v_h$ by $\pi_1^{-1}(v_h)$ and replacing $v_\ell$ by $\pi_1^{-1}(v_\ell)$.
 Analogously, we define $H_2$ by taking $$(H[A_2 \setminus S^2] - F) \cup (\pi_2^{-1}(N) \setminus \{ \pi_2^{-1}(v_\ell), \pi_2^{-1}(v_h) \}, \emptyset)$$ and replacing $v_h$ by $\pi_2^{-1}(v_h)$ and replacing $v_\ell$ by $\pi_2^{-1}(v_\ell)$.

 Let $\langle u_1 = v_h, u_2, \ldots, u_l = v_\ell\rangle$ be an ordering of the nodes on $T$ by traversing $T$ from $v_h$ to $v_\ell$ whilst avoiding the interior of $c$.
 Then we thus have
 \begin{align*}
 \Omega_1 &= \langle \pi_1^{-1}(v_\ell), v_{\ell+1}, \ldots, v_{h-1}, \pi_1^{-1}(v_h), \pi_1^{-1}(u_2), \ldots, \pi_1^{-1}(u_{l-1})\rangle \text{ and} \\
 \Omega_2 &= \langle \pi_2^{-1}(v_\ell), v_{\ell-1}, \ldots, v_1, v_k, v_{k-1}, \ldots, \pi_2^{-1}(v_h), \pi_2^{-1}(u_2), \ldots, \pi_2^{-1}(u_{l_1})\rangle .
 \end{align*}
 Thus the appropriate vortex societies are $(H_1, \Omega_1)$ and $(H_2, \Omega_2)$.
 This leads to us ultimately defining a $\Sigma''$-rendition $\rho'$ of the graph $(G' - S) - F$.
 As before the depth of these new societies is at most the depth of $(H, \Psi)$ plus 1 due to the copied vertex.
 The first three points of our statement and the last one hold by construction for similar reasons as in our last case, where we set the pair of objects asked for in the statement to simply be two copies of what we defined here.

 We now want to split up $\mathcal{L}$ along $N$ to derive a new linkage appropriate for $G'$.
 For this sake we perform the exact same modifications presented in the proof of \zcref{lem:joinvortices}, which we therefore do not repeat here, to derive a linkage $\mathcal{L}^*$ from $\mathcal{L}$ that lives in $G',$ is vital in $G'$ if $\mathcal{L}$ was vital in $G$ and has at most $r$ more paths than $\mathcal{L}$.
 This confirms the remaining two points of our statement.
 
 Furthermore, taking $\mathcal{L}' \coloneqq \mathcal{L}^* \cap ((G' - S) - F),$ we quickly note that $|\mathcal{L}'| \leq |\mathcal{L}| + r + 6d$.
 Thus, using \zcref{prop:vitalitypreservedthroughrestriction}, we can easily confirm the remaining two points of our statement.

 \paragraph{$T$ is non-separating and creates a single open hole.}
 The last remaining case is largely analogous to the previous one once we have appropriately split apart our surface.
 Once again assuming that $T$ intersects some vortex $c$ in a non-trivial component, we let $(H,\Psi),$ $\Psi = \langle v_1, \ldots, v_k \rangle,$ and $N_\rho(c) \cap T = \{ v_\ell, v_h \},$ with $\ell < h,$ be defined as above.

 Let $B$ be the unique boundary of $\Sigma - T$.
 We again choose $\epsilon$ as above and now define $\Sigma'$ by first splitting $T$ at $v_\ell$ to create a curve $L$ with endpoints $s,t,$ with $s = v_\ell$ such that the immediate successor of $s$ in $L$ lies in $c,$ creating two copies $L_1,L_2$ of $L,$ gluing these together at their endpoints to form the closed curve $T'',$ and attaching $T''$ onto $B$ to form a new surface $\Sigma''$ such that for both $i \in [2]$ there exist a bijection $\pi_i \colon L \rightarrow L_i$ with $r_{\epsilon,\Sigma'}(p) \setminus (L_1 \cup L_2) \subseteq r_{\epsilon,\Sigma}(\pi^{-1}_i(p))$ for all $p \in L_i$.
 Once again the genus of $\Sigma'$ is less than the genus of $\Sigma$.
 We then define the new graph $G'$ analogously to what we discussed above.

 Following this, we cap $\Sigma'$ with a disc $\Delta'$ bounded by $T'',$ which yields the surface $\Sigma'',$ and let $\Delta = \Delta' \cup c$ be the new vortex we want to introduce.
 Again for all cells in $\rho$ other than $c,$ we can retain their placement and what $\sigma$ assigns to them (up to the modifications $\pi_1$ and $\pi_2$ prescribe) to define $\rho',$ leaving only the contents of $\Delta$ up in the air.

 We can now follow the construction of the vortex society from our previous case quite closely.
 We define $H^*$ by taking $H,$ replacing $v_\ell$ with its copies $\pi_1^{-1}(v_\ell),\pi_2^{-1}(v_\ell)$ and $v_h$ with its copies $\pi_1^{-1}(v_h),\pi_2^{-1}(v_h)$.

 Now the appropriate vortex society is $(H^* \cup (\pi_1^{-1}(N) \cup \pi_2^{-1}(N), \emptyset), \Psi'),$ where $\Psi'$ is derived from $\Psi$ and $L_1$ and $L_2,$ as follows:
 If the order of $N$ according to $L$ is $\langle u_1 = v_h, u_2, \ldots, u_l = v_\ell\rangle,$ then we have
 \begin{align*}
 \Psi' &= \langle \pi_1^{-1}(v_\ell), v_{\ell+1}, v_{\ell+2}, \ldots, v_{h-1}, \pi_1^{-1}(v_h) = \pi_1^{-1}(w_2), \pi_1^{-1}(w_3), \ldots, \pi_1^{-1}(w_l), \\
 &\qquad \pi_2^{-1}(v_\ell), v_{\ell-1}, v_{\ell-2}, \ldots, v_1, v_k, v_{k-1}, \ldots, v_{h+1}, \pi_2^{-1}(v_h) = \pi_2^{-1}(w_2), \ldots, \pi_2^{-1}(w_l) \rangle .
 \end{align*}
 Note that the depth of the vortex society of $\Delta$ is at most the depth of the vortex society of $c$ plus 2.
 This ultimately defines $\rho'$ as a $\Sigma'$-rendition of $G'$ and the first four points as well as the last point of our statement holds for reasons analogous to our previous arguments.
 To capture the fourth and fifth points of the statement, we then modify $\mathcal{L}$ as in the previous cases and thus arrive at the conclusion of our proof.
\end{proof}

\subsection{Building nests in renditions with large distances and high representativity}
Now that we can ensure that the representativity of our renditions is high and that any two vortices in them are at large distance with respect to each other, we need some tools to build nests around each of these vortices that will end up being stacked with respect to each other.
We start with what can be considered the base case in which we have a rendition in the sphere.

\begin{lemma}\label{lem:sphererenditionnestbuilder}
 Let $d,k,s$ be positive integers with $s$ being even.
 Let $G$ be a graph with a sphere-rendition $\rho$ such that any two vortices of $\rho$ lie at distance at least $sd+3$ and let $c_1, \ldots, c_k$ be the vortices of $\rho$.
 Further, let there exists a disc $\Delta_0$ such that there exists a flat $(2sd)$-wall $W$ in $\rho$ whose trace is contained in $\Delta_0$ and $\Delta_0$ is disjoint from all vortices of $\rho$.
 
 Then there exist $k$ sets of $\nicefrac{sd}{2}$ tight, concentric cycles $\mathcal{C}_1, \ldots, \mathcal{C}_k,$ with $\Delta_i$ being the closed disc bounded by the trace of the inner cycle of $\mathcal{C}_i$ in $\rho$ for each $i \in [k],$ such that $c_i \subseteq \Delta_i$ and $(\mathcal{C}_i,c_i)$ and $(\mathcal{C}_j,c_j)$ are stacked for all distinct $i,j \in [k]$.
\end{lemma}
\begin{proof}
 Due to the presence of $W,$ we can definitely find at least $k$ sets of $s' \coloneqq \nicefrac{sd}{2}$ concentric cycles $\mathcal{C}_1, \ldots, \mathcal{C}_k,$ with $\Delta_i$ defined as in our statement, in $\rho$ such that $c_i \subseteq \Delta_i$.
 In particular, we may choose all of these sets such that they are tight in $\rho$.
 For each $i \in [k]$ let $\mathcal{C}_i = \langle C_1^i, \ldots, C_{s'}^i\rangle,$ which $C_1^i$ being the innermost cycle of $\mathcal{C}_i$.

 Suppose that there exist distinct $i,j \in [k]$ and $h_i,h_j \in [s']$ such that $V(C_{h_i}^i) \cap V(C_{h_j}^j) \neq \emptyset$.
 Then let $v \in V(C_{h_i}^i) \cap V(C_{h_j}^j)$ be a node of $\rho$ and note that the tightness of both $\mathcal{C}_i$ and $\mathcal{C}_j$ implies that there exist $\rho$-aligned curves $\gamma_i,\gamma_j$ such that for both $x \in \{ i,j \}$ the curve $\gamma_i$ ends in $v$ and $N_\rho(c_x)$ and intersects at most $h_x+1$ nodes of $\rho$.
 Thus $\gamma_i \cup \gamma_j$ contain a $\rho$-aligned curve intersecting at most $h_i+h_j+2 \leq 2\nicefrac{sd}{2}+2$ nodes of $\rho$ and thus contradict the fact $c_i$ and $c_j$ lie at distance at least $sd + 3$ in $\rho$.

 Therefore, for any two distinct $i,j \in [k]$ we have $V(\mathcal{C}_i) \cap V(\mathcal{C}_j) = \emptyset$.
 The fact that all our sets are now pairwise stacked then follows directly from standard considerations for the topology of the sphere.
\end{proof}

Now that we have proven this base case, we note that the following result due to Mohar and Thomassen lifts the above lemma to graphs with a rendition in a surface of positive genus and high representativity.

\begin{proposition}[Mohar and Thomassen \cite{MoharT2001Graphs}]
 Let $d$ be a positive integer.
 Let $G$ be a graph with an embedding $\psi$ in a surface of positive genus and representativity at least $2d+2$.
 Then for every face $F$ of $\psi$ there exists a set of $d$ cycles $C_1, \ldots, C_d$ in $G$ such that the $C_d$-disc $\Delta$ contains $f$ and $C_1, \ldots, C_{d-1}$ are concentric in $\Delta$.
\end{proposition}
 
\begin{corollary}\label{lem:generalrenditionnestbuilder}
 Let $d,k,s$ be positive integers with $s$ being even.
 Let $G$ be a graph with a $\Sigma$-rendition $\rho$ such that $\Sigma$ has positive genus, any two vortices of $\rho$ lie at distance at least $sd+3,$ and $\rho$ has representativity at least $sd + 2,$ and let $c_1, \ldots, c_k$ be the vortices of $\rho$.

 Then there exist $k$ sets of $\nicefrac{sd}{2}$ tight, concentric cycles $\mathcal{C}_1, \ldots, \mathcal{C}_k,$ with $\Delta_i$ being the closed disc bounded by the trace of the inner cycle of $\mathcal{C}_i$ in $\rho$ for each $i \in [k],$ such that $c_i \subseteq \Delta_i$ and $(\mathcal{C}_i,c_i)$ and $(\mathcal{C}_j,c_j)$ are stacked for all distinct $i,j \in [k]$.
\end{corollary}

\subsection{Reduction to graphs without apices}
Before we can prove the main result of this section, we need a way to produce a large mesh if our graph has high treewidth.
The standard reference for this is the grid theorem \cite{RobertsonS86GraphminorsV,ChekuriC2016Polynomial,ChuzhoyT2021Tighter}.
Since we also exclude a clique as a minor, we can actually use slightly more optimised bounds from \cite{GorskySTW2026Catching}.

\begin{proposition}[Gorsky, Stamoulis, Thilikos, and Wiederrecht \cite{GorskySTW2026Catching}]\label{prop:ratcatcher}
 There exists a function $\mathsf{grid}_{\ref{prop:ratcatcher}} \colon \mathbb{N}^2 \rightarrow \mathbb{N}$ with $\mathsf{grid}(t,k) \in \mathbf{O}(t^2k + t^{2304})$ such that the following holds.
 
 Let $k,t$ be positive integers.
 Any graph excluding both a $(k \times k)$-grid and $K_t$ as a minor has treewidth at most $\mathsf{grid}_{\ref{prop:ratcatcher}}(t,k)$. 
\end{proposition}

Having now gathered all of our tools, we are finally ready to gather them up to prove the main theorem of this section which reduces an instance of a graph with a vital linkage that is already excluding a clique-minor of a certain size to an instance which has a rendition with some convenient properties.
Most importantly, the rendition is for the entire instance and does not require the deletion of apices.

\begin{theorem}\label{thm:reducetonoapices}
 There exist functions $\mathsf{clique}_{\ref{thm:reducetonoapices}}, \mathsf{breadth}_{\ref{thm:reducetonoapices}}, \mathsf{depth}_{\ref{thm:reducetonoapices}} \colon \mathbb{N} \rightarrow \mathbb{N},$ $\mathsf{paths}_{\ref{thm:reducetonoapices}} \colon \mathbb{N}^2 \rightarrow \mathbb{N}$ and $\mathsf{wall}_{\ref{thm:reducetonoapices}} \colon \mathbb{N}^3 \rightarrow \mathbb{N},$ with $\mathsf{breadth}_{\ref{thm:reducetonoapices}}(t), \mathsf{clique}_{\ref{thm:reducetonoapices}}(t), \mathsf{depth}_{\ref{thm:reducetonoapices}}(t), \mathsf{paths}_{\ref{thm:reducetonoapices}}(t,s), \mathsf{wall}_{\ref{thm:reducetonoapices}}(t,s,r) \in \mathbf{poly}(t + s + r)$.

 Let $r,t,s$ be positive integers.
 Let $(G,R)$ be an annotated graph with $\mathsf{bidim}(G,R) \leq t$ that does not contain a $K_{\mathsf{clique}_{\ref{thm:reducetonoapices}}(t)}$-minor and let $\mathcal{L}$ be an $R$-linkage in $G$ with a pattern partitioning $R$.

 Then either $G$ has treewidth in $\mathsf{wall}_{\ref{thm:reducetonoapices}}(t,s,r),$ or there exists a surface $\Sigma$ with a blank $\Sigma$-rendition $\rho$ with breadth $b \leq \mathsf{breadth}_{\ref{thm:reducetonoapices}}(t)$ and depth at most $\mathsf{depth}_{\ref{thm:reducetonoapices}}(t)$ of an annotated graph $(H,R')$ containing a linkage $\mathcal{L}',$ and there exist $b$ vortices $c_1, \ldots, c_b \in C(\rho)$ such that
 \begin{enumerate}
 \item the genus of $\Sigma$ is in $9t^4,$ and $|\mathcal{L}'| - |\mathcal{L}| \leq \mathsf{paths}_{\ref{thm:reducetonoapices}}(t,s),$

 \item there exist $b$ sets of $s(\mathsf{depth}_{\ref{thm:reducetonoapices}}(t) + 1)$ tight, concentric cycles $\mathcal{C}_1, \ldots, \mathcal{C}_b,$ with $\Delta_i$ being the closed disc bounded by the trace of the inner cycle of $\mathcal{C}_i$ in $\rho$ for each $i \in [b],$ such that $c_i \subseteq \Delta_i,$ and $(\mathcal{C}_i,c_i)$ and $(\mathcal{C}_j,c_j)$ are stacked for all distinct $i,j \in [b],$

 \item the vortices in $\rho$ pairwise lie at distance at least $2s(\mathsf{depth}_{\ref{thm:reducetonoapices}}(t) + 1) + 3,$

 \item there exists a disc $\Delta_0$ such that, if we let $\Delta_i'$ be the trace of the outer cycle of $\mathcal{C}_i$ for each $i \in [b_1],$ we have $\Delta_0 \cap \Delta_i' = \emptyset$ for all $i \in [b_1]$ and there exists a flat $r$-wall in $\rho$ whose trace is contained in $\Delta_0,$ and
 
 \item $R'$ contains exactly the endpoints of $\mathcal{L}'$ and if $\mathcal{L}$ is vital in $G$ then $\mathcal{L}'$ is vital in $H$.
 \end{enumerate}
\end{theorem}
\begin{proof}
 Our functions are chosen as follows:
 \begin{align*}
 \mathsf{clique}_{\ref{thm:reducetonoapices}}(t) &\coloneqq \lfloor \nicefrac{5}{2} \rfloor t^2 + 1, \\
 \mathsf{breadth}_{\ref{thm:reducetonoapices}}(t) &\coloneqq \nicefrac{3}{2}(t^2-1)(3t-4) + t(t-1)-3 + 2\mathsf{apex}_{\ref{thm:localstructure}}(t^2, t) + 9t^4, \\
 \mathsf{depth}_{\ref{thm:reducetonoapices}}(t) &\coloneqq (\mathsf{depth}_{\ref{thm:localstructure}}(t^2, t)+2)(\mathsf{breadth}_{\ref{thm:reducetonoapices}}(t) - 1) + 18t^4, \\
 \mathsf{paths}_{\ref{thm:reducetonoapices}}(t,s) &\coloneqq 2\mathsf{apex}_{\ref{thm:localstructure}}(t^2, t) + 18st^4(\mathsf{depth}_{\ref{thm:reducetonoapices}}(t) + 1) \\
 &\qquad \qquad + \ \mathsf{breadth}_{\ref{thm:reducetonoapices}}(t) \cdot (2s(\mathsf{depth}_{\ref{thm:reducetonoapices}}(t) + 1) + 8\mathsf{depth}_{\ref{thm:reducetonoapices}}(t)), \text{ and}\\
 \mathsf{wall}_{\ref{thm:reducetonoapices}}(t,s,r) &\coloneqq \mathsf{grid}_{\ref{prop:ratcatcher}}(\mathsf{clique}_{\ref{thm:localstructure}}(t), \mathsf{mesh}_{\ref{thm:localstructure}}(t,\mathsf{clique}_{\ref{thm:reducetonoapices}}(t), ( r \\
 &\qquad \qquad + \ (9t^4(4s(\mathsf{depth}_{\ref{thm:reducetonoapices}}(t) + 1)+4)+ \mathsf{breadth}_{\ref{thm:reducetonoapices}}(t)(4s(\mathsf{depth}_{\ref{thm:reducetonoapices}}(t) + 1) + 6) \\
 &\qquad \qquad + \ 2s(\mathsf{depth}_{\ref{thm:reducetonoapices}}(t) + 1) ))\sqrt{\mathsf{apex}_{\ref{thm:localstructure}}(t^2,t) + 1} ) .
 \end{align*}
 
 For the rest of our proof we may assume that $\mathcal{L}$ is vital.
 This makes a few of our constructions easier to analyse.
 We may assume that $G$ has treewidth at least $\mathsf{wall}_{\ref{thm:reducetonoapices}}(t,s,r),$ as we are otherwise done.
 Thus, \zcref{prop:ratcatcher} tells us that there exists an $r'$-mesh in $G,$ where
 \begin{align*}
 r' &\geq \mathsf{mesh}_{\ref{thm:localstructure}}(t,\mathsf{clique}_{\ref{thm:reducetonoapices}}(t), ( r + (9t^4(4s(\mathsf{depth}_{\ref{thm:reducetonoapices}}(t) + 1)+4) \\
 &\qquad \qquad + \ \mathsf{breadth}_{\ref{thm:reducetonoapices}}(t)(4s(\mathsf{depth}_{\ref{thm:reducetonoapices}}(t) + 1) + 6) + 2s(\mathsf{depth}_{\ref{thm:reducetonoapices}}(t) + 1) ))\sqrt{\mathsf{apex}_{\ref{thm:localstructure}}(t^2,t) + 1} ) .
 \end{align*}
 This allows us to apply \zcref{thm:localstructure}.
 We cannot reach the first item of \zcref{thm:localstructure}, since $G$ already excludes a $K_{\mathsf{clique}_{\ref{thm:reducetonoapices}}(t)}$-minor, and we cannot reach the second item, since finding a red $K_{t^2+1}$-minor contradicts $\mathsf{bidim}(G,R) \leq t$.
 Thus \zcref{thm:localstructure} returns a set $A \subseteq V(G)$ with $|A| \leq \mathsf{apex}_{\ref{thm:localstructure}}(t^2, t)$ and a surface $\Sigma$ of genus less than $9t^4$ such that $(G - A, R \setminus A)$ has a blank $\Sigma$-rendition $\rho$ with breadth at most $\nicefrac{3}{2}(t^2-1)(3t-4) + t(t-1)-3,$ depth at most $\mathsf{depth}_{\ref{thm:localstructure}}(t^2, t),$ and there exists a vortex-free, $\rho$-aligned disc $\Delta \subseteq \Sigma$ such that the restriction of $\rho$ to $\Delta$ contains a flat $r_0$-mesh $W_0,$ where
 \begin{align*}
 r_0 &= (r+9t^4(4s(\mathsf{depth}_{\ref{thm:reducetonoapices}}(t) + 1)+4) + \mathsf{breadth}_{\ref{thm:reducetonoapices}}(t)(4s(\mathsf{depth}_{\ref{thm:reducetonoapices}}(t) + 1) + 6) \\
 &\qquad \qquad + \ 2s(\mathsf{depth}_{\ref{thm:reducetonoapices}}(t) + 1) )\sqrt{\mathsf{apex}_{\ref{thm:localstructure}}(t^2,t) + 1} . 
 \end{align*}

 We first modify the graph to get rid of the apices in $A$ at the cost of introducing new vortices.
 For each $a \in A$ let $u_1^a,u_2^a$ be its up to two neighbours in $\bigcup \mathcal{L}$.
 (We will ignore the vertices in $A$ which constitute a trivial path in $\mathcal{L}$ and simply delete them from $\mathcal{L}$ and $G$ for simplicity.
 Clearly this does not affect the vitality of $\mathcal{L}$ in $G$.)
 Let $L_a \in \mathcal{L}$ be the path with $a \in A$.
 First, we delete $a$ from $G$ and instead introduce $a_1,a_2$ (with $a_1=a_2$ if $u_1^a = u_2^a$), making them adjacent to $u_1^a,u_2^a$ respectively.
 Then we let $L_a^i$ be the component of $L_a - a$ that contains $u_a^i$ for each $i \in [2],$ remove $L_a$ from $\mathcal{L},$ and add $L_a^i \cup ( \{ a_i,u_i^a \}, \{ a_iu_i^a \})$ for each $i \in [2]$.
 These modifications are performed iteratively.
 The identity of $u_1^a,$ $u_2^a,$ and $L_a$ may change several times over the course of this process, but this is not dangerous as long as we make these modifications iteratively.
 
 Let $G'$ and $\mathcal{L}_1$ be the resulting graph and linkage.
 We then modify $\rho$ into a $\Sigma$-rendition of $G'$ as follows.
 For each $a \in A$ and $i \in [2],$ if $u_i^a \not\in N(\rho)$ there exists a unique $c \in C(\rho)$ such that $u_i^a \in V(\sigma_\rho(c)),$ then we simply add $a_i$ and its incident edge to $\sigma_\rho(c)$.
 Note that if $c$ is a vortex, this does not increase its depth, as we may simply add $a_i$ to one of the bag of its linear decomposition which also includes $u_i^a$.
 Thus we may instead assume that $u_i^a \in N(\rho)$.
 In this case we introduce a new cell $c_a \subseteq \Sigma$ that intersects no other cell of $\rho$ (or any of the cells we have added thus far) and intersects $N(\rho)$ exactly in $u_i^a$.
 We then update $\rho$ to let $\sigma_\rho(c_a) = ( \{ a_i,u_i^a \}, \{ a_iu_i^a \})$ and note that $c_a$ is a vortex despite having very short length, since it contains a terminal of our linkage and thus one of the coloured vertices of our annotated graph.

 We let $\rho_1$ be the resulting rendition, noting that this is indeed a $\Sigma$-rendition of $G'$ whose depth is still at most $d$ and the fact that $\mathcal{L}_1$ is vital in $G'$ is implied by the vitality of $\mathcal{L}$ in $G$.
 Furthermore, $\rho_1$ is blank and we have added at most $2|A|$ ``vortices'' to $\rho_1$ in comparison to $\rho,$ which means that the breadth of $\rho_1$ is safely within $\mathbf{poly}(t + b')$ and similarly, if we let $R' = (R \cap V(G')) \cup \{ a_1,a_2 ~\colon~ a \in A \},$ we can easily note that $\mathsf{bidim}(G',R') \in \mathbf{poly}(t + b')$.
 Unfortunately, some of these new vortices may have landed in our $r_0$-wall $W_0$.
 So to sort these out, we use the fact that $r_0 = r_1\sqrt{\mathsf{apex}_{\ref{thm:localstructure}}(t^2,t) + 1}$ to find an $r_1$-wall $W_1 \subseteq W_0$ contained in a $\rho'$-aligned disc that is entirely disjoint from all of our vortices, by dividing our wall up into $\mathsf{apex}_{\ref{thm:localstructure}}(t^2,t) + 1$ pairwise disjoint $r_1$-walls, where
 \[ r_1 = r + (9t^4(4s(\mathsf{depth}_{\ref{thm:reducetonoapices}}(t) + 1)+4) + \mathsf{breadth}_{\ref{thm:reducetonoapices}}(t)(4s(\mathsf{depth}_{\ref{thm:reducetonoapices}}(t) + 1) + 6) + 2s(\mathsf{depth}_{\ref{thm:reducetonoapices}}(t) + 1) ) . \]

 \paragraph{Finding a rendition with high distance and representativity.}
 Now our goal is to get to an instance of a graph with a vital linkage embedded on some surface of high representativity (or on the sphere) and with high distance between any two vortices.
 We will accomplish this as follows.
 First, we increase the distance between any two vortices \zcref{lem:joinvortices}.
 Then, we check if the representativity of our rendition is high enough and if it is not, we apply \zcref{lem:pumprepresentativity}.
 However, whilst this reduces the genus of the instance we are interested in, this does introduce up to two more vortices that may now actually be close to some of the existing vortices.
 Thus, whenever we are forced to apply \zcref{lem:pumprepresentativity}, we return to the first step and iterate.
 Since \zcref{lem:pumprepresentativity} decreases the genus of the surface we work on, this will terminate.

 In total, \zcref{lem:pumprepresentativity} will produce a number of vortices proportional to the genus of $\Sigma$ and \zcref{lem:joinvortices} decreases the number of vortices by one with each application whilst only adding at most 2 to the depth of the rendition.
 Thus, if $b_1$ is the breadth of $\rho_1$ and $g$ is the genus of $\Sigma,$ we will at worst generate $b_2 \coloneqq b_1 + 2g$ vortices in this process.
 In particular, since each application of \zcref{lem:joinvortices} fuses two of these vortices whose depth simply gets added together with 2, the final depth of our vortices is at most $d_2 \coloneqq (d+2)(b_2-1) + 2g,$ since \zcref{lem:joinvortices} can only be applied at most $b_2 - 1$ times and \zcref{lem:pumprepresentativity} can only applied $g$ times.

 Let us now be more precise about the numbers with which we apply both of these lemmas.
 First, consider \zcref{lem:joinvortices}, for which we try to ask for distance $2s(d_2 + 1) + 3$ between the vortices of the current rendition.
 Any modification that results from this lemma does not change the surface, results in a vital linkage with at most $4s(d_2 + 1) + 6$ new endpoints (and thus the annotated graph created by taking these endpoints as the red vertices still has bounded bidimensionality), all of these endpoints are found in vortices, and outside of the addition of the new fused vortex the rendition is not changed.
 Thus, if the old rendition contained a flat wall of order $w',$ our new rendition preserves the existence of a $(w' - (4s(d_2+1) + 6))$-wall.

 On the other hand, we apply \zcref{lem:pumprepresentativity} asking for representativity $2s(d_2 + 1) + 2$.
 If this cuts apart our surface into two components, we choose the $w'$-wall in our original instance to orient ourselves, since the new instance still contains a flat $(w' - (4s(d_2 + 1) + 4))$-mesh which is draw on either of the two new surfaces.
 This allows us to choose the surface containing the remains of our flat wall.
 We then note that this also creates a vital linkage with at most $2s(d_2 + 1) + 2 + 6d_2$ endpoints (again ensuring that the bidimensionality of the corresponding annotated graph is still bounded), all of which are found in vortices.
 Furthermore, both the breadth and the depth of our rendition increases at most by 2.

 Thus, after this process we will arrive at a graph $G_2$ containing a vital linkage $\mathcal{L}_2$ with a $\Sigma_2$-rendition $\rho_2$ such that
 \begin{enumerate}
 \item the genus of $\Sigma_2$ is at most the genus of $\Sigma,$

 \item $|\mathcal{L}_2| \leq |\mathcal{L}_1| + g(2s(d_2 + 1)+2) + b_2(2s(d_2 + 1) + 2 + 6d_2),$

 \item all endpoints of $\mathcal{L}_2$ are found in vortices of $\rho_2,$

 \item the depth of $\rho_2$ is at most $d + d_2,$

 \item the breadth of $\rho_2$ is at most $b_2,$

 \item there exists a flat $(r_1 - (g(4s(d_2 + 1)+4) + b_2(4s(d_2 + 1) + 6)))$-wall $W_2$ in $\rho_2,$

 \item the vortices of $\rho_2$ pairwise lie at distance at least $2s(d_2 + 1) + 3,$ and

 \item if $\Sigma_2$ has positive genus, the representativity of $\rho_2$ is at least $2s(d_2 + 1) + 2$.
 \end{enumerate}
 From this we can easily derive an annotated graph $(G_2,R_2)$ by letting $R_2$ be the set of endpoints of the paths in $\mathcal{L}_2$.

 \paragraph{Building nests and making the rendition blank.}
 By simply applying \zcref{lem:sphererenditionnestbuilder} or \zcref{lem:generalrenditionnestbuilder}, we can now find the desired pairwise disjoint nests of order $s(d_2 + 1)$ around the vortices of $\rho_2$.
 These nests may intersect $W_2,$ but this at worst causes us to lose the outermost $2s(d_2+1)$ horizontal and vertical paths of $W_2$ and thus we can find an $r_2$-wall $W_3 \subseteq W_2,$ with $r_2 = (r_1 - (g(4s(d_2 + 1)+4) + b_2(4s(d_2 + 1) + 6))) - 2s(d_2 + 1),$ whose perimeter is disjoint from the trace of all cycles of all nests we have just constructed.
 We note that 
 \[ r_1 - (9t^4(4s(\mathsf{depth}_{\ref{thm:reducetonoapices}}(t) + 1)+4) + \mathsf{breadth}_{\ref{thm:reducetonoapices}}(t)(4s(\mathsf{depth}_{\ref{thm:reducetonoapices}}(t) + 1) + 6) + 2s(\mathsf{depth}_{\ref{thm:reducetonoapices}}(t) + 1) ) \geq r .\]

 Finally, we still have to turn $\rho_2$ into a blank rendition of $(G_2,R_2)$.
 Note that the only problem here is that some vertices of $R_2$ may be placed on the boundary of one of our vortices.
 Thus for each $v \in R_2$ with $N_{\rho_2}(c)$ for some vortex $c$ of $\rho_2,$ we create a copy $v' \not\in V(G_2),$ make $v$ adjacent to $v'$ and add $v'$ to the vortex $c$.
 Note that this does not increase the depth of $c$ and thus does not increase the depth of $\rho_2$.
 Furthermore, we then elongate the path in $\mathcal{L}_2$ that uses $v$ by adding $v'$ and the edge $vv'$ to it and replacing $v$ in $R_2$ with $v'$.
 This clearly preserves the vitality of the linkage.
 Repeating this for every vertex in $R_2$ yields a blank rendition of a graph with a vital linkage having the desired properties.
 This completes our proof.
\end{proof}

\section{What to do with a clique?}\label{sec:largeclique}
One of the most fundamental parts of our proof and the Graph Minor Algorithm we will present later is the problem of dealing with the existence of large clique-minors in our graph, since we would really like them to be absent to apply the tools from Graph Minor Structure Theory.
Luckily there is a well-established set of tools which allow us to see that these cliques are actually quite helpful for us.
This gives us permission to focus on graphs that exclude a clique-minor.

The following is (6.1) from \cite{RobertsonS1995Graph}, where we note that the original statement confusingly discusses \textsl{digraphs}.
However the proof by Robertson and Seymour in particular shows the following statement for undirected graphs.
To get the algorithmic part of the statement derives from algorithm (6.2) in \cite{RobertsonS1995Graph}.

\begin{proposition}[Robertson and Seymour \cite{RobertsonS1995Graph}]\label{thm:cliqueirrelevantvertex}
 Let $G$ be a graph, let $\ell, d \in \mathbb{N}$ be non-negative integers, and let $Z \subseteq V(G)$ with $|Z| \leq \ell$.
 Further let $k = \lfloor \nicefrac{5}{2} \cdot \ell \rfloor + 3d^2 + 1,$ let $\varphi$ be a minor model of $K_k$ in $G,$ and let $(A,B)$ be a separation in $G$ such that
 \begin{enumerate}
 \item \label{item:cliqueright} $A \cap \varphi(v) = \emptyset$ for some $v \in V(K_k),$
 \item \label{item:terminalleft} $Z \subseteq A,$
 \item \label{item:smallsep} subject to \zcref{item:cliqueright} and \zcref{item:terminalleft}, $(A,B)$ is of minimum order, and
 \item subject to \zcref{item:cliqueright}, \zcref{item:terminalleft}, and \zcref{item:smallsep}, $A$ is maximal.
 \end{enumerate}
 Then for any $v \in B \setminus A,$ $v$ is irrelevant to the $d$-folio of $G$ relative to $Z$.

 In particular, given $G,$ $Z,$ and $\varphi$ as above, such a vertex $v$ can be found in $\mathbf{poly}(\ell + d)||H||$-time.
\end{proposition}

We will also require the following result from \cite{ProtopapasTW2025Colorful} which refines some of the tools from \cite{RobertsonS1995Graph}.

\begin{proposition}[Protopapas, Thilikos, and Wiederrecht \cite{ProtopapasTW2025Colorful}]\label{thm:colorfulclique}
 Let $q,t,k \in \mathbb{N}$ be positive integers with $k \geq \lfloor \nicefrac{3}{2} \cdot q t \rfloor + t$.
 Let $(G,\chi)$ be a $q$-colorful graph such that $G$ contains a minor model $\varphi$ of $K_k$.
 Then one of the following is true:
 \begin{enumerate}
 \item There exists a colorful minor model $\psi$ of a rainbow $K_t$ in $(G,\chi)$ such that $\mathcal{T}_\psi$ is a truncation of $\mathcal{T}_\varphi,$ or

 \item There exists a set $S \subseteq V(G)$ of size at most $qt - 1$ such that the $\mathcal{T}_\varphi$-big component of $G - S$ is restricted.
 \end{enumerate}
 Furthermore, there exists an algorithm that takes as input $t,$ $(G,\chi),$ and $\varphi$ as above and finds one of the two outcomes in time $\mathbf{poly}(q+t)|\!|G|\!|$.
\end{proposition}

We combine these two results into a simple lemma that will later on allow us to efficiently find irrelevant vertices for the folio-problem in \zcref{sec:folio}.

\begin{lemma}\label{lem:bigcliqueirrelevantvertexforfolio}
 Let $k,b,d \in \mathbb{N}$ be positive integers, let $(G,Z)$ be an annotated graph with bidimensionality less than $b$ and $|Z| \leq k,$ and let $\varphi$ be a $K_t$-minor-model in $G$ with $t \geq \lfloor \nicefrac{5}{2} \rfloor b^2 + 3d^2 + 1$.
 
 Then $G$ contains an irrelevant vertex for the $(k,d)$-folio of $(G,Z),$ which can be identified in $\mathbf{poly}(k+b+d)|G|$-time.
\end{lemma}
\begin{proof}
 We view $(G,Z)$ as a colorful graph $(G,\chi)$ with one colour which is assigned by $\chi$ exactly to the vertices in $Z$.
 If there existed a colorful minor model $\psi$ of a rainbow $K_{b^2}$ in $(G,\chi),$ then this would directly contradict $\mathsf{bidim}(G,Z) < b$.
 Thus, as $t \geq \lfloor \nicefrac{5}{2} \rfloor b^2 + 1 \geq \lfloor \nicefrac{3}{2} \cdot b^2 \rfloor + b^2,$ \zcref{thm:cliqueirrelevantvertex} tells us that there exists a set $S \subseteq V(G)$ of order $b^2 - 1$ such that the $\mathcal{T}_\varphi$-big component of $G - S$ is restricted and thus contains no terminals from $Z$.

 In particular, this means that there exists a separation $(A,B)$ in $G$ with $A \cap B = S$ such that $A \cap \varphi(v) = \emptyset$ for some $V(K_t)$ and we have $Z \subseteq A$.
 This allows us to apply \zcref{thm:cliqueirrelevantvertex}, which confirms that any vertex in $B$ is irrelevant, completing our proof.
\end{proof}

Of course to do something with a clique-minor, we first need to find one.
To do this efficiently, we take a tool from Reed and Wood \cite{ReedW2009Lineartime}, who state that their algorithm actually takes a running time dependent on both the number of vertices and the number of edges of $G$.
However as Kawarabayashi, Kobayashi, and Reed \cite{KawarabayashiKR2012Disjoint} point out, by only considering $2^{t-3}|G|$ edges of the graph, this algorithm ends up finding the clique-minor in a running time only depending on $|G|$ and the size of the clique we search for.

\begin{proposition}[Reed and Wood \cite{ReedW2009Lineartime}]\label{prop:quicklyfindingaclique}
 Let $t$ be a positive integer and let $G$ be a graph with $|\!|G|\!| \geq 2^{t-3}|G|$.
 Then one can find a $K_t$-minor model in $G$ in $\mathbf{O}(t|G|)$-time. 
\end{proposition}

\section{The Vital Linkage Theorem}\label{sec:VitalLinkage}
In the previous sections we have shown that
\begin{itemize}
 \item graphs with large cliques contain irrelevant vertices for vital linkages (see \zcref{sec:largeclique}),
 
 \item graphs with high treewidth and excluding a clique-minor yields a graph with a well-structured rendition without apices (see \zcref{sec:noapexnocry}), 

 \item graphs hosting a vital linkage that have a well-structured rendition without apices can be reduce down until they no longer have vortices without introducing many terminals (see \zcref{sec:apexfree}), and

 \item graphs containing a vital linkage embedded on a surface with bounded genus and a bounded number of holes hosting terminals of the linkage cannot have large treewidth (see \zcref{sec:mazoitproof}).
\end{itemize}
We will now apply these insights in the given order to prove our first main theorem, which gives the first explicit bound for the Linkage Function, showing that it is polynomial in the number of terminals and singly exponential only in the bidimensionality of these terminals.

\begin{theorem}\label{thm:VitalLinkage}
There exists a function $\beta_{\ref{thm:VitalLinkage}} \colon \mathbb{N}^2 \to \mathbb{N}$ such that for every annotated graph $(G,T)$ with $\mathsf{bidim}(G,T) \leq b,$ if there exists a vital $T$-linkage of order at most $k$ in $G,$ then $\mathsf{tw}(G) \leq \beta_{\ref{thm:VitalLinkage}}(k,b)$.
Moreover, $\beta_{\ref{thm:VitalLinkage}}(k,b) \in 2^{\mathbf{poly}(b)} \cdot \mathbf{poly}(k)$.
\end{theorem}
\begin{proof}
 For future use, we let
 \begin{align*}
 r &\coloneqq 2\mathsf{cycles}_{\ref{cor:noapexreducetomazoit}}(b + 2\mathsf{paths}_{\ref{thm:reducetonoapices}}(b,\mathsf{breadth}_{\ref{thm:reducetonoapices}}(b)), \mathsf{depth}_{\ref{thm:reducetonoapices}}(b), \\
 &\qquad f_{\ref{thm:MazoitEndgame}}(\mathsf{breadth}_{\ref{thm:reducetonoapices}}(b) + 9b^4, \mathsf{terminals}_{\ref{cor:noapexreducetomazoit}}(k \\
 &\qquad + \ 2\mathsf{paths}_{\ref{thm:reducetonoapices}}(b,\mathsf{breadth}_{\ref{thm:reducetonoapices}}(b)), \mathsf{breadth}_{\ref{thm:reducetonoapices}}(b), \mathsf{depth}_{\ref{thm:reducetonoapices}}(b), 9b^4))) .
 \end{align*}
 We will prove this theorem by setting $\beta_{\ref{thm:VitalLinkage}}(k,b) \coloneqq \mathsf{wall}_{\ref{thm:reducetonoapices}}(b,\mathsf{breadth}_{\ref{thm:reducetonoapices}}(b),r)$.
 Let $(G_0,T_0)$ be a counterexample to our statement, meaning that there exists a vital $T_0$-linkage of order at most $k$ in $G_0,$ we have $\mathsf{bidim}(G_0,T_0) \leq b,$ and $\mathsf{tw}(G_0) > \beta_{\ref{thm:VitalLinkage}}(k,b)$.

 We now first apply \zcref{lem:bigcliqueirrelevantvertexforfolio} to exclude the presence of a $K_t$-minor in $G_0$ with $t = \lfloor \nicefrac{5}{2} \rfloor b^2 + 1,$ whilst noting that the linkage problem is simply the folio problem with detail 0.
 This follows, since the presence of an irrelevant vertex for the linkage would imply that $G_0$ is in fact not a minimal counterexample.
 Together with the large treewidth of $G_0,$ this allows us to now apply \zcref{thm:reducetonoapices}, which yields an annotated graph $(G_1,T_1),$ a vital $T_1$-linkage $\mathcal{L}_1$ in $G_1,$ and a blank, stretched $\Sigma_1$-rendition $\rho_1$ with breadth $b_1 \leq \mathsf{breadth}_{\ref{thm:reducetonoapices}}(b)$ and depth at most $\mathsf{depth}_{\ref{thm:reducetonoapices}}(b)$ on a surface $\Sigma_1$ of genus at most $9b^4,$ such that
 \begin{enumerate}
 \item $\mathsf{bidim}(G_1,T_1) \leq b + 2\mathsf{paths}_{\ref{thm:reducetonoapices}}(b,\mathsf{breadth}_{\ref{thm:reducetonoapices}}(b))$ and $|\mathcal{L}_1| \leq k + \mathsf{paths}_{\ref{thm:reducetonoapices}}(b,\mathsf{breadth}_{\ref{thm:reducetonoapices}}(b)),$

 \item there exist $b_1$ sets of $\mathsf{breadth}_{\ref{thm:reducetonoapices}}(b)(\mathsf{depth}_{\ref{thm:reducetonoapices}}(b) + 1)$ concentric cycles $\mathcal{C}_1, \ldots, \mathcal{C}_{b_1},$ with $\Delta_i$ being the closed disc bounded by the trace of the inner cycle of $\mathcal{C}_i$ in $\rho_1$ for each $i \in [b_1],$ such that $c_i \subseteq \Delta_i$ and $(\mathcal{C}_i,c_i)$ and $(\mathcal{C}_j,c_j)$ are stacked for all distinct $i,j \in [b_1],$ 
 
 \item the vortices in $\rho_1$ pairwise lie at distance at least $4(\mathsf{depth}_{\ref{thm:reducetonoapices}}(t) + 1) + 3,$ and

 \item there exists a disc $\Delta_0$ such that, if we let $\Delta_i'$ be the trace of the outer cycle of $\mathcal{C}_i$ for each $i \in [b_1],$ we have $\Delta_0 \cap \Delta_i' = \emptyset$ for all $i \in [b_1]$ and there exists a flat $r$-mesh $M_1$ in $\rho_1$ whose trace is contained in $\Delta_0$.
 \end{enumerate}
 Thanks to the second property, we may now apply \zcref{thm:linking_vortices} to our rendition $\rho_1$ to find a blank, stretched $\Sigma_1$-rendition $\rho_1'$ of $(G_1,T_1)$ with breadth at most $b_1$ such that property i) above is preserved, and each vortex of $\rho_1'$ either has a boundary consisting of at most 3 nodes, or is $\mathsf{depth}_{\ref{thm:reducetonoapices}}(b)$-stuffed with gap at most 1.
 Furthermore, $M_1$ continues to be flat in $\rho_1$ and $\Delta_0$ is disjoint from all vortices of $\rho_1'$.

 These constructions allow us to observe that by iteratively stripping the perimeter from $M_1,$ we find a set of $\nicefrac{r}{2}$ concentric cycles $\mathcal{C}_M$ forming an $\nicefrac{r}{2}$-locus.
 We now let 
 \[ s \coloneqq f_{\ref{thm:MazoitEndgame}}(b_1 + 9b^4, \mathsf{terminals}_{\ref{cor:noapexreducetomazoit}}(k + 2\mathsf{paths}_{\ref{thm:reducetonoapices}}(b,\mathsf{breadth}_{\ref{thm:reducetonoapices}}(b)), b_1, \mathsf{depth}_{\ref{thm:reducetonoapices}}(b), 9b^4)) . \]
 Note that we have
 \[ \nicefrac{r}{2} \geq \mathsf{cycles}_{\ref{cor:noapexreducetomazoit}}(b + 2\mathsf{paths}_{\ref{thm:reducetonoapices}}(b,\mathsf{breadth}_{\ref{thm:reducetonoapices}}(b)), \mathsf{depth}_{\ref{thm:reducetonoapices}}(b), s) . \]
 We may therefore apply \zcref{cor:noapexreducetomazoit} to find an annotated graph $(G_2,T_2)$ with
 \[ |T_2| \leq \mathsf{terminals}_{\ref{cor:noapexreducetomazoit}}(k + 2\mathsf{paths}_{\ref{thm:reducetonoapices}}(b,\mathsf{breadth}_{\ref{thm:reducetonoapices}}(b)), b_1, \mathsf{depth}_{\ref{thm:reducetonoapices}}(b), 9b^4) \in \mathbf{poly}(k+b) \cdot 2^{\mathbf{poly}(b)}, \]
 a vital $T_2$-linkage $\mathcal{L}_2,$ an embedding $\psi_2$ of $G_2$ in a surface $\Sigma_2$ with at most $b_1$ boundary components, and containing a set of concentric cycles $\mathcal{C} = \langle C_1, \ldots, C_s\rangle$.
 Furthermore, we have $\hat{\Sigma}_2 \cong \hat{\Sigma}_1,$ there exists a disc $\Delta_2 \subseteq \Sigma_2$ such that $\psi(C_i) \subseteq \Delta_2$ for all $i \in [s],$ $G_2 = \bigcup \mathcal{C} \cup \bigcup \mathcal{L}_2,$ $\psi(T') \subseteq \mathsf{bd}(\Sigma_2),$ and for all $L \in \mathcal{L}_2$ the curve $\psi(L)$ is internally disjoint from $\mathsf{bd}(\Sigma_2)$.

 Our goal is to now adjust $\psi_2$ and $\Sigma_2$ to fit the requirements of \zcref{thm:MazoitEndgame}.
 To accomplish this, we simply cut out an open disc in $\Sigma_2$ that does not intersect any point used by $\psi_2$.
 This allows us to homeomorphically deform both $\Sigma_2$ and $\psi_2$ into a surface that can be described as the union of a closed disc $\Delta_3$ with $\mathsf{bd}(\Delta_3) = \psi(C_s)$ and a set of at most $g_3 \coloneqq b_1 + 9b^4$ strips $S_1, \ldots, S_{g_3}$.
 Now we are ready to apply \zcref{thm:MazoitEndgame} to $(\Sigma_2, G_2, \mathcal{C}, \mathcal{L}_2)$ and since $s \geq f_{\ref{thm:MazoitEndgame}}(g_3,|T_2|),$ this leads to a contradiction to the vitality of our linkage.
 We note that $g_3 \in \mathbf{poly}(b)$.
 Thus
 \begin{align*}
 f_{\ref{thm:MazoitEndgame}}(g_3,|T_2|) \in \mathbf{poly}(k+b) \cdot 2^{\mathbf{poly}(b)}& \text{ and} \\
 \mathsf{cycles}_{\ref{cor:noapexreducetomazoit}}(b + 2\mathsf{paths}_{\ref{thm:reducetonoapices}}(b,\mathsf{breadth}_{\ref{thm:reducetonoapices}}(b)), \mathsf{depth}_{\ref{thm:reducetonoapices}}(b), s) \in \mathbf{poly}(k+b) \cdot 2^{\mathbf{poly}(b)} &,
 \end{align*}
 which confirms the bound claimed in our statement as $b \leq k$.
\end{proof}

\section{Tackling the Folio Problem}\label{sec:folio}
The previous sections were all concerned almost entirely with proving \zcref{thm:VitalLinkage}.
However, establishing good bounds for the Linkage Function is merely the first step towards reigning in the galactic parametric dependencies in the running time of the Graph Minor Algorithm.
Indeed, in order to prove \zcref{thm:Main1Folio} we need to prove \zcref{thm:Main2Folio} which is the statement that guarantees the existence of an \textsl{irrelevant vertex} in all instances of large enough treewidth.
In this section and the following one we take on the task to deal with the algorithmic consequences of \zcref{thm:VitalLinkage} for the $(k,d)\text{-}\textsc{Folio}$ problem.

The main goal of this section is to show that within a large ``homogeneous'' flat wall $W$ avoiding the set $R$ of annotated vertices, it is always possible to redraw a minor model invading the compass of $W$ in a way such that it avoids the centre of the wall.

\subsection{Quickly combing an annulus}
\label{sec:combing}
A key result in the algorithm theory of graph minors is the so-called ``Annulus Combing Lemma'' of Golovach, Stamoulis, and Thilikos \cite{GolovachST2023Combing}.
This result is a slight strengthening of similar lemmas originally proven by Robertson and Seymour \cite{RobertsonS2012Graph} and essentially says the following:

Given a $\Sigma$-rendition $\rho$ of a graph together with some disc $\Delta \subseteq \Sigma$ such that the restriction of $\rho$ to $\Delta$ is a vortex-free rendition containing a large grounded cylindrical wall $W$ together with a pattern of order $k,$ say $\Pi,$ in $G$ and a $\Pi$-linkage $\mathcal{L},$ then there exists a still big cylindrical subwall $W'$ of $W,$ another $\Pi$-linkage $\mathcal{L}'$ and a subset $\mathcal{R}$ of radial paths in $W'$ such that the intersection of $\mathcal{L}'$ with the compass of $W'$ is a subset of the paths in $\mathcal{R}$.
Here, the \emph{compass} of a cylindrical wall $W$ grounded in $\rho$ is the union of all graphs $\sigma(c)$ over the cells $c$ contained in the annulus defined by the traces of the outer- and inner-most cycle of $W$.
Similarly, the \emph{compass} of a wall $W^*$ grounded in $\rho$ is the union of all graphs $\sigma(c)$ over the set of all cells $c$ contained in the disc defined by the trace of the perimeter of $W^*$.
Notice that, whenever we have a graph $G$ with an embedding $\Gamma$ in some surface $\Sigma,$ one may obtain a trivial $\Sigma$-rendition of $G$ from $\Gamma$ by defining, for each edge $xy \in E(G)$ a cell $c_{xy}$ obtained by choosing an $\varepsilon$-disc $\Delta_{xy}$ around the drawing of $xy$ with $\varepsilon > 0$ such that $x,y \in \mathsf{bd}(\Delta_{xy}),$ $\Delta_{xy}$ intersects $\Gamma$ only in $x,$ $y,$ and $xy,$ and the curve $xy \setminus \{ x,y\}$ is contained in the interior of $\Delta_{xy}$.

The notions of ``large'' and ``big'' in the statement above depend on $k$ and the Linkage Function from \zcref{thm:VitalLinkage}.
In fact, we prove a slightly more general statement for which we require some further definitions.

\begin{figure}[ht]
 \centering
 \begin{tikzpicture}

 \pgfdeclarelayer{background}
		\pgfdeclarelayer{foreground}
			
		\pgfsetlayers{background,main,foreground}

 \begin{pgfonlayer}{background}
 \pgftext{\includegraphics[width=7.5cm]{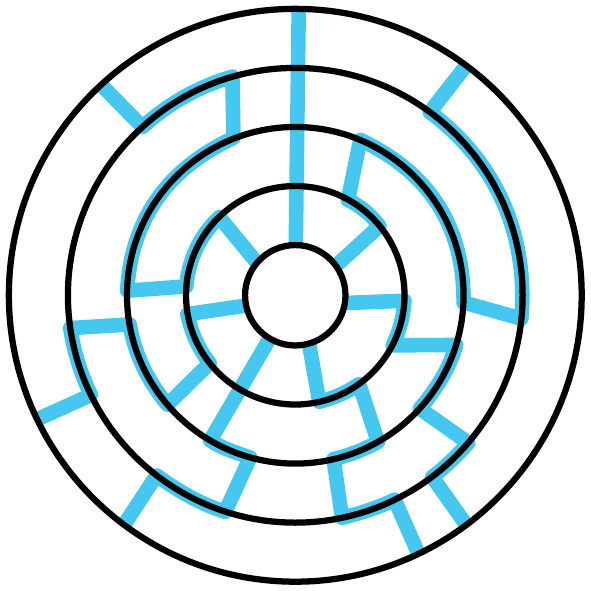}} at (C.center);
 \end{pgfonlayer}{background}
			
 \begin{pgfonlayer}{main}
 \node (C) [v:ghost] {};
 
 \end{pgfonlayer}{main}
 
 \begin{pgfonlayer}{foreground}
 \end{pgfonlayer}{foreground}

 \end{tikzpicture}
 \caption{A railed annulus with $5$ circles and $7$ rails.}
 \label{fig:RailedAnnulus}
\end{figure}

\paragraph{A partially embedded railed annulus.}
Let $w,r \geq 0$ be two integers.
A \emph{$(w,r)$-railed annulus} is a graph $A_{w,r}$ such that there exist pairwise vertex-disjoint cycles $C_1,\dots, C_w,$ called the \emph{circles}, and pairwise vertex-disjoint paths $R_1,\dots R_r,$ called the \emph{rails}, such that
\begin{enumerate}
 \item $A_{w,r} = \bigcup_{i\in [w]} C_i \cup \bigcup_{i \in [r]} R_i$ is a planar graph,
 \item $R_i$ has one endpoint on $C_1,$ the other on $C_w$ for all $i\in[r],$ and is otherwise disjoint from $C_1\cup C_w,$
 \item for all $i\in[r],$ $j\in[w],$ $R_i \cap C_j$ is a path, and
 \item for each $i\in[r],$ $R_i$ visits the cycles $C_1, \dots,C_w$ in the order listed when traversing along $R_i$ starting from its endpoint on $C_1$.
\end{enumerate}
We call $C_1$ and $C_w$ the \emph{boundary circles} of $A_{w,r}$.
We say that a graph $A$ is a \emph{railed annulus} if there exist $r$ and $w$ such that $A$ is a $(w,r)$-railed annulus.\
See \zcref{fig:RailedAnnulus} for an example.

\begin{figure}[ht]
 \centering
 \begin{tikzpicture}

 \pgfdeclarelayer{background}
		\pgfdeclarelayer{foreground}
			
		\pgfsetlayers{background,main,foreground}
			
 \begin{pgfonlayer}{main}
 \node (C) [v:ghost] {};

 \node(L) [v:ghost,position=180:7cm from C] {
 \begin{tikzpicture}

 \pgfdeclarelayer{background}
		 \pgfdeclarelayer{foreground}
			
		 \pgfsetlayers{background,main,foreground}

 \begin{pgfonlayer}{background}
 \pgftext{\includegraphics[width=5cm]{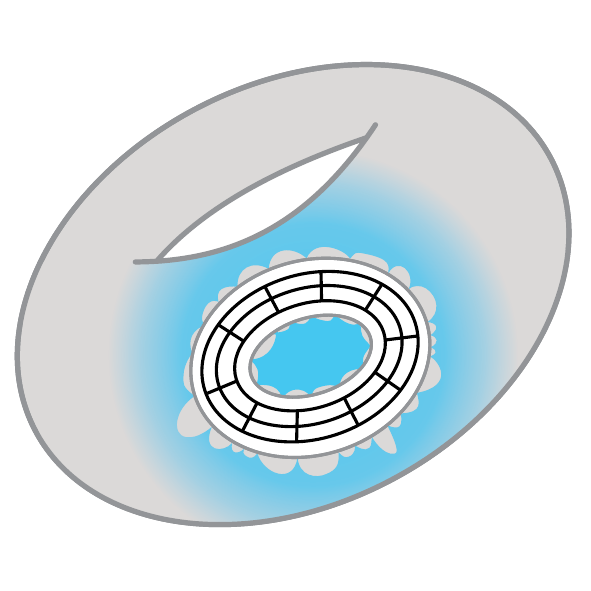}} at (C.center);
 \end{pgfonlayer}{background}
			
 \begin{pgfonlayer}{main}
 \node (C) [v:ghost] {};
 
 \end{pgfonlayer}{main}
 
 \begin{pgfonlayer}{foreground}
 \end{pgfonlayer}{foreground}

 \end{tikzpicture}
 };

 \node(M) [v:ghost,position=0:0cm from C] {
 \begin{tikzpicture}

 \pgfdeclarelayer{background}
		 \pgfdeclarelayer{foreground}
			
		 \pgfsetlayers{background,main,foreground}

 \begin{pgfonlayer}{background}
 \pgftext{\includegraphics[width=5cm]{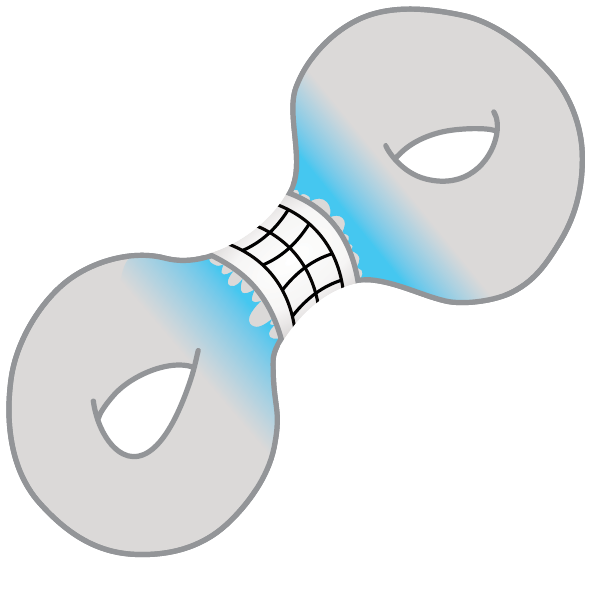}} at (C.center);
 \end{pgfonlayer}{background}
			
 \begin{pgfonlayer}{main}
 \node (C) [v:ghost] {};
 
 \end{pgfonlayer}{main}
 
 \begin{pgfonlayer}{foreground}
 \end{pgfonlayer}{foreground}

 \end{tikzpicture}
 };

 \end{pgfonlayer}{main}
 
 \begin{pgfonlayer}{foreground}
 \end{pgfonlayer}{foreground}

 \end{tikzpicture}
 \caption{Two separating annuli in different surfaces. The separating annulus on the left bounds a disc in its interior, while the separating annulus on the right is genus-reducing for its surface.
 In this paper we only care about the situation on the left, however, \zcref{thm:CombedAnus} holds also in the situation on the right.}
 \label{fig:SeparatingAnnuli}
\end{figure}
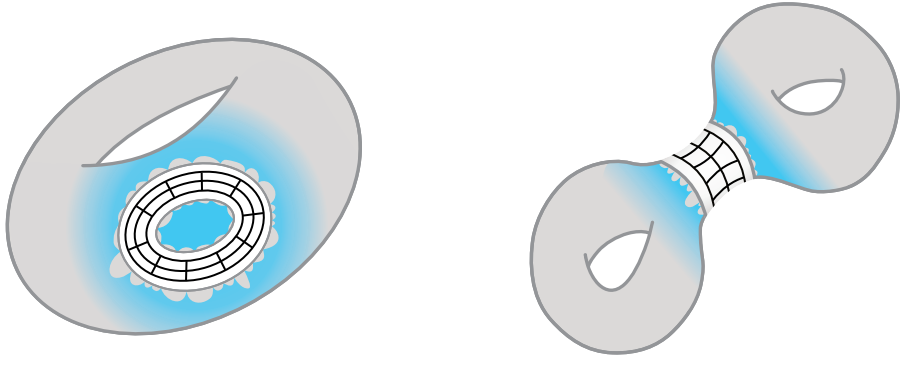

Now let $\Sigma$ be a surface and $G$ be a graph with a $\Sigma$-rendition $\rho$.
We say that a railed annulus $A \subseteq G$ is \emph{separating} in $\rho$ if
\begin{enumerate}
 \item $A$ is grounded in $\rho,$
 \item if $A$ has at least two circles and at least two rails, then there exists an annulus $\circledcirc \subseteq \Sigma$ such that the traces of the boundary circles of $A$ coincide with the two boundary components of $\circledcirc,$ and the trace of every rail of $A$ is contained in $\circledcirc,$
 \item $\circledcirc$ is disjoint from the closures of the vortices of $\rho,$ and
 \item for every circle $C$ of $A,$ the trace of $C$ is a separating curve in $\Sigma$.
\end{enumerate}
See \zcref{fig:SeparatingAnnuli} for some examples.
We call the annulus $\circledcirc$ above the \emph{domain} of $A$ in $\rho$.
The \emph{crop of $G$ by $\circledcirc$} is the union of all graphs $\sigma(c)$ whose cell $c$ of $\rho$ is contained in $\circledcirc$.
Finally, if $(G,T)$ is an annotated graph with a $\Sigma$-rendition $\rho$ and a separating railed annulus $A$ with domain $\circledcirc,$ we say that $A$ is \emph{blank} if the crop of $G$ by $\circledcirc$ does not contain any vertices from $T$.

Let $p\geq 0$ and $q\geq 2$ be integers.
Let $A = A_{2p+q,r}$ be a railed annulus.
The \emph{$p$-trimming} of $A$ is the railed annulus $A_{q,r}$ obtained from $A$ by removing the first and last $p$ circles and restricting the rails to their minimal $C_{p+1}$-$C_{p+q}$-subpaths.

\paragraph{Combed linkages.}
Let $\Sigma$ be a surface, $G$ be a graph with a $\Sigma$-rendition $\rho,$ and let $A \subseteq G$ be a separating railed annulus in $G$ with domain $\circledcirc$ in $\rho$.
Let $H$ be the crop of $G$ by $\circledcirc$.
A linkage $\mathcal{L}$ in $G$ is \emph{combed in $A$} if the intersection of $\bigcup\mathcal{L}$ and $H$ is a subgraph of the union of the rails of $A$.

The main theorem of this subsection now reads as follows.

\begin{theorem}\label{thm:CombedAnus}
There exist functions $f_{\ref{thm:CombedAnus}},h_{\ref{thm:CombedAnus}},g_{\ref{thm:CombedAnus}}\colon \mathbb{N}^2 \to \mathbb{N}$ such that for every annotated graph $(G,T)$ with $\mathsf{bidim}(G,T) \leq b,$ every surface $\Sigma,$ and every $\Sigma$-rendition $\rho$ of $G,$ the following holds:

If $A = A_{a,r}$ is a blank and separating railed annulus with $a \geq f_{\ref{thm:CombedAnus}}(k,b)$ and $r \geq g_{\ref{thm:CombedAnus}}(k,b)$ in $G$ and $\rho$ and $\mathcal{L}$ is a $T$-linkage of order at most $k,$ then there exists a $T$-linkage $\mathcal{L}^{\star}$ in $G$ such that
\begin{enumerate}
 \item $\tau(\mathcal{L}) = \tau(\mathcal{L}^{\star}),$
 \item $\mathcal{L}^{\star}$ is combed in the $h_{\ref{thm:CombedAnus}}(k,b)$-trimming of $A,$ and
 \item if $\circledcirc$ is the domain of $A$ and $G'$ is the crop of $G$ by $\circledcirc,$ then $\bigcup\mathcal{L}^{\star} - G' \subseteq \bigcup\mathcal{L} - G'$.
\end{enumerate}

Moreover, $f_{\ref{thm:CombedAnus}}(k,b), h_{\ref{thm:CombedAnus}}(k,b), g_{\ref{thm:CombedAnus}}(k,b) \in 2^{\mathbf{poly}(b)} \cdot \mathbf{poly}(k)$.
\end{theorem}

We will need the following classic theorem for lower bounds on the treewidth of graphs (see for example \cite{SeymourT1993GraphSearching,RobertsonS2012Graph}).

Let $G$ be a graph.
A \emph{bramble} in $G$ is a collection of connected subgraphs $\mathcal{B}$ of $G$ such that for all $B_1,B_2 \in \mathcal{B}$ either $V(B_1) \cap V(B_2) \neq \emptyset$ or there exists an edge with one end in $B_1$ and the other in $B_2$.
A \emph{hitting set} for a bramble $\mathcal{B}$ is a set $S \subseteq V(G)$ such that $S \cap V(B) \neq \emptyset$ for all $B \in \mathcal{B}$.
The \emph{order} of a bramble $\mathcal{B}$ is the minimum size of a hitting set for $\mathcal{B}$.

\begin{proposition}[Reed \cite{Reed1997Tree}]\label{prop:brambles}
Let $k \geq 0$ be an integer and $G$ be a graph.
If $G$ has a bramble of order $k+1,$ then $\mathsf{tw}(G) \geq k$.
\end{proposition}

\begin{proof}[Proof of \zcref{thm:CombedAnus}.]
Let us start by defining our functions, to this end let $\beta$ be the Linkage Function from \zcref{thm:VitalLinkage}.
We set
\begin{align*}
 f_{\ref{thm:CombedAnus}}(k,b) & \coloneqq 16\beta(k,b) + 16\\
 h_{\ref{thm:CombedAnus}}(k,b) & \coloneqq 8\beta(k,b) + 8\text{, and}\\
 g_{\ref{thm:CombedAnus}}(k,b) & \coloneqq 6\beta(k,b) + 6.
\end{align*}
The claimed bounds now follow directly from \zcref{thm:VitalLinkage}.
\smallskip

Let $\mathcal{C} = \langle C_1,\dots, C_a\rangle$ be the circles of $A$ and let $G''' \coloneqq \bigcup\mathcal{L} \cup \bigcup\mathcal{C}$.
In the following we assume that $\mathcal{L}$ is chosen to minimise $|E(\mathcal{L}) \setminus E(\mathcal{C})|$ among all linkages $\mathcal{Q}$ in $G$ with $\tau(\mathcal{Q}) = \tau(\mathcal{L})$.
This is a fair assumption since the assertion allows us to make changes to the original linkage.
\smallskip

Let $\circledcirc$ be the domain of $A$ in $\Sigma$.
Let $\mathcal{S}$ be the set of all $(V(C_1) \cup V(C_a))$-subpaths of the paths in $\mathcal{L}$ which are contained in the crop of $G$ by $\circledcirc$.
A path in $\mathcal{S}$ is called a \emph{river} if it is a $V(C_1)$-$V(C_a)$-path, a \emph{mountain} if it is a $V(C_1)$-path, and a \emph{valley} if it is a $V(C_a)$-path.
See \zcref{fig:MountainsAndValleys} for an illustration.

Let now $G''$ be obtained from $G$ by iteratively picking a vertex $v$ of degree exactly $2$ and contracting one of its two incident edges.
In a second step, let now $G'$ be obtained from $G''$ by contracting every remaining edge which belongs to some circle of $A$ \textsl{and} some path in $\mathcal{L}$.
Notice that $G'$ still contains a sequence $\mathcal{C}' =\langle C'_1,\dots, C'_a \rangle$ of vertex-disjoint cycles such that $C_i'$ is obtained from $C_i$ for each $i\in[a]$.
Moreover, if $x$ is an endpoint of some path in $\mathcal{L},$ then $x$ is not a vertex of any circle of $A$ since $A$ is blank and thus, $x$ has degree $1$ in $G''$.
It follows that $G'$ contains a linkage $\mathcal{L}'$ such that $\tau(\mathcal{L}') = \tau(\mathcal{L})$.
This is true in particular since the construction step to go from $G''$ to $G'$ only contracted edges that already belonged to $\mathcal{L}$ and thus ensured that the paths in $\mathcal{L}$ remain paths.

Finally, notice that, since $G'$ is a minor of $G,$ we may inherit a $\Sigma$-rendition $\rho'$ of $G'$ from the restriction of $\rho$ to $G'''$.

\begin{figure}[ht]
 \centering
 \begin{tikzpicture}

 \pgfdeclarelayer{background}
		\pgfdeclarelayer{foreground}
			
		\pgfsetlayers{background,main,foreground}

 \begin{pgfonlayer}{background}
 \pgftext{\includegraphics[width=10cm]{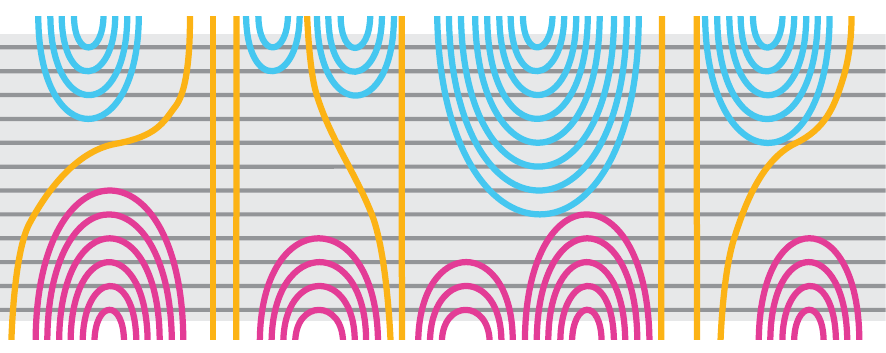}} at (C.center);
 \end{pgfonlayer}{background}
			
 \begin{pgfonlayer}{main}
 \node (C) [v:ghost] {};
 
 \end{pgfonlayer}{main}
 
 \begin{pgfonlayer}{foreground}
 \end{pgfonlayer}{foreground}

 \end{tikzpicture}
 \caption{An illustration of a collection of concentric cycles together with segments of a linkage. In \textcolor{HotMagenta}{magenta} we depict the mountains while the valleys are drawn in \textcolor{CornflowerBlue}{blue} and the rivers are drawn in \textcolor{ChromeYellow}{yellow}.}
 \label{fig:MountainsAndValleys}
\end{figure}

\begin{beautifulclaim}\label{claim:VitalAnus}
The linkage $\mathcal{L}'$ is vital in $G'$.
\end{beautifulclaim}

\begin{claimproof}
Suppose there is a vertex $v\in V(G')$ not used by $\mathcal{L}'$.
Then there exists $i\in[a]$ such that $v \in V(C_i')\setminus V(\mathcal{L})$.
But then, $v$ must have degree $2$ in $G'$ which is ruled out by the construction of $G'$.
Hence, $V(\mathcal{L}') = V(G')$.

Indeed, we may make a stronger observation:
Let $v\in V(C'_i)$ for some $i\in[\ell],$ then $v$ has degree $4$ in $G'$ since $v$ has degree at least $3,$ may be incident with at most four edges of $G',$ and two of those edges are required to belong to $\mathcal{L}'$.
Since $\mathcal{L}'$ cannot share an edge with $C_i'$ as otherwise we would have contracted this edge, it follows that $v$ cannot be of degree $3$.

Now suppose that $\mathcal{L}'$ is not unique in $G'$.
Then there exists a linkage $\mathcal{Q}'$ in $G'$ with the same pattern such that $\mathcal{L}' \neq \mathcal{Q}'$.
Suppose further that there exists an edge $e$ in $E(\mathcal{L}') \setminus E(\mathcal{Q}')$ such that $e \notin E(\mathcal{Q}')$.
Then, since $G'$ is a minor of $G$ and all endpoints of the paths in $\mathcal{L}$ still have degree $1$ in $G',$ there would be a linkage $\mathcal{Q}$ in $G$ of the same pattern as $\mathcal{L}$ such that $E(\mathcal{Q}) \subseteq (E(\mathcal{L}) \cup E(\mathcal{C})) \setminus \{ e \}$.
Hence, $|E(\mathcal{Q}) \setminus E(\mathcal{C})| < |E(\mathcal{L}) \setminus E(\mathcal{C})|$ which is a contradiction to our choice of $\mathcal{L}$.
Therefore, $\mathcal{L}'$ must be unique in $G'$ which finally implies that $\mathcal{L}'$ is a vital linkage.
\end{claimproof}

Let $T' \coloneqq V(G') \cap T$.

\begin{beautifulclaim}\label{claim:BidimensionalitySurvives}
We have that $\mathsf{bidim}(G',T') \leq b$.
\end{beautifulclaim}

\begin{claimproof}
The claim follows directly from the fact that bidimensionality is monotone under taking red minors and $(G',T')$ is indeed a red minor of $(G,T)$.
\end{claimproof}

We say that two rivers $J_1$ and $J_2$ are \emph{consecutive}, if deleting their traces from $\circledcirc$ leaves one component $\zeta',$ whose closure is homeomorphic to a disc, that does not intersect the trace of any other river in $\mathcal{S}$.
Note that we may define what it means for two rails of $A'$ to be \emph{consecutive} in the same way.
We call the closure $\zeta = \mathsf{cl}(\zeta')$ of such a disc a \emph{slice} of $\circledcirc$.
Notice that the rivers in $\mathcal{S}$ decompose $\circledcirc$ into $|\mathcal{S}|$ slices $\zeta_1,\dots,\zeta_{\ell}$ where $\ell$ is the number of rivers in $\mathcal{S}$.
Moreover, for every mountain or valley $P$ in $\mathcal{S}$ there exists a unique $i\in[\ell]$ such that the trace of $P$ is contained in $\zeta_i$.
For each $i\in [\ell]$ let us denote by $\mathcal{M}_i$ the set of all mountains whose trace is contained in $\zeta_i$ and by $\mathcal{V}_i$ the set of all valleys whose trace is contained in $\zeta_i$.

In what follows, we wish to argue that rivers and valleys cannot enter deeply into the circles.
For this, the goal is to apply~--~among others~--~\zcref{lem:drainedwell} which requires us to confine the curves of mountains and valleys into wells.
Since the arguments for dealing with valleys is analogous to those for mountains, we only discuss the mountain case below.
\smallskip

\textbf{Draining a well in a mountain.}
Fix any $i\in[\ell]$.
Let $G_i$ be the graph $\bigcup \mathcal{M}_i \cup \mathcal{C}'$.
Then $G_i$ is a planar graph with a plane embedding $\Gamma_i$ inherited from $\rho$ through \zcref{thm:rendition_to_skeleton} and we may declare $\Omega_i$ to be the cyclic ordering of the vertices of $C'_1$ obtained by traversing along $C'_1$ in clockwise order.
Then $(G_i,\Omega_i,\rho,\mathcal{C}',\mathcal{M}_i)$ is a well and by \zcref{lem:drainedwell} we may assume this well to be drained.
We call $(G_i,\Omega_i,\rho,\mathcal{C}',\mathcal{M}_i)$ the \emph{well associated with $\mathcal{M}_i$}.
This is, because the new collection of paths returned by \zcref{lem:drainedwell} must still have all of its traces contained in $\zeta_i$ and is therefore guaranteed to stay disjoint from the rivers of $\mathcal{S}$.

In particular, any mountain in $\mathcal{M}_i$ encloses a disc with the boundary of $\zeta_i$ which cannot contain the trace of any valley from $\mathcal{V}_i$.
This observation allows us to apply the same treatment to the valleys and assume that their respective wells are also drained.
\smallskip

\textbf{Levelling mountains and filling valleys.}
With the assumption that each of the wells associated with $\mathcal{M}_i$ or $\mathcal{V}_i$ is drained we are now ready to prove our first major claim.

\begin{beautifulclaim}\label{claim:FlatMountains}
Let $M$ be a mountain and $V$ be a valley in $\mathcal{S}$.
Then $M$ is disjoint from $C'_{2\beta(k,b) + 3}$ and $V$ is disjoint from $C'_{a - 2\beta(k,b) - 2}$.
\end{beautifulclaim}

\begin{claimproof}
As before, the arguments for mountains and valleys are analogous so it suffices to only discuss the mountains.

Suppose there exist $i\in[\ell]$ and a mountain $M \in \mathcal{M}_i$ which intersects $C'_{2\beta(k,b)+3}$.
Then, since the well $(G_i,\Omega_i,\rho,\mathcal{C}',\mathcal{M}_i)$ is drained, we know that there exists a sequence $\langle M_1,\dots,M_{2\beta(k,b)+3} \rangle$ of pairwise vertex-disjoint mountains in $\mathcal{M}_i$ such that for each $j\in [2\beta(k,b)+3],$ $M_j$ intersects $C'_j$.
Notice that, whenever a mountain $M_j$ intersects the circle $C'_j,$ then $M_j$ also intersects $C'_{j'}$ for all $j' \leq j$.

Now let us define, for each $j \in [2\beta(k,b) + 3]$ the graph $B_j$ to be $M_j \cup C_j'$ and let $\mathcal{B} \coloneqq \{ B_j ~\colon~ j \in [2\beta(k,b) + 3] \}$.
See \zcref{fig:BrambleMountain} for an illustration of this construction.
Notice that every member of $\mathcal{B}$ is connected, and any two members must intersect by the remark above.
Moreover, it follows from the fact that the $M_j$ are pairwise vertex-disjoint and the $C_j'$ are also pairwise vertex disjoint, that no vertex of $G'$ may be contained in more than two distinct members of $\mathcal{B}$.
Hence, $\mathcal{B}$ is a bramble of order at least $\beta(k,b) + 2$.
Also recall that by \zcref{claim:BidimensionalitySurvives} we know that $\mathsf{bidim}(G',T') \leq b$.
Therefore, \zcref{prop:brambles} says that $\mathsf{tw}(G') \geq \beta(k,b) + 1$.
However, by \zcref{claim:VitalAnus} we have that $\mathcal{L}'$ is a vital $T$-linkage in $G'$ and by \zcref{thm:VitalLinkage} this means that $\mathsf{tw}(G') \leq \beta(k,b)$ which is absurd.
\end{claimproof}

\begin{figure}[ht]
 \centering
 \begin{tikzpicture}

 \pgfdeclarelayer{background}
		\pgfdeclarelayer{foreground}
			
		\pgfsetlayers{background,main,foreground}

 \begin{pgfonlayer}{background}
 \pgftext{\includegraphics[width=10cm]{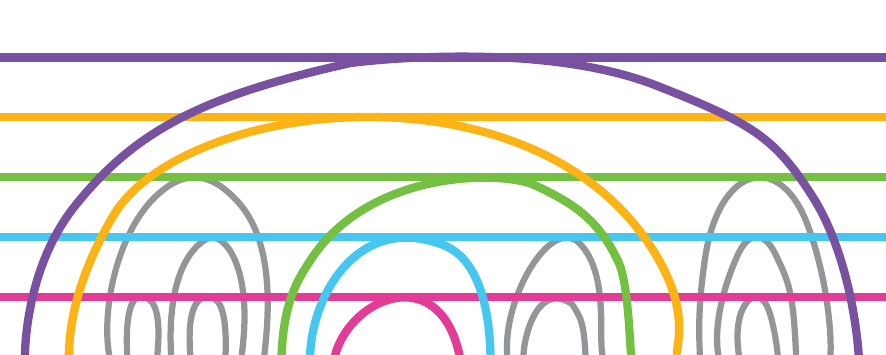}} at (C.center);
 \end{pgfonlayer}{background}
			
 \begin{pgfonlayer}{main}
 \node (C) [v:ghost] {};
 
 \end{pgfonlayer}{main}
 
 \begin{pgfonlayer}{foreground}
 \end{pgfonlayer}{foreground}

 \end{tikzpicture}
 \caption{An illustration of the construction of the bramble in the proof of \zcref{claim:FlatMountains}. Each horizontal line corresponds to one of the circles of the annulus and the pairs consisting of a single circle and one mountain that form the bramble elements are colour coded. The mountains in \textcolor{DarkGray}{gray} illustrate how other parts of the linkage might interact with the mountains we select for our bramble.}
 \label{fig:BrambleMountain}
\end{figure}

\textbf{Bounding the rivers.}
Now we know that neither mountains, nor valleys are able to penetrate deep into our annulus.
The only members of $\mathcal{S}$ we need to be concerned with now are the rivers.
These must necessarily intersect all cycles in $\mathcal{C}',$ however, with an argument similar to the one above, we can find that their number is bounded.
Let $\mathcal{W} \subseteq \mathcal{S}$ be the set of all rivers in $\mathcal{S}$.

\begin{beautifulclaim}\label{claim:FewRivers}
We have that $|\mathcal{W}| \leq 2\beta(k,b) + 2$.
\end{beautifulclaim}

\begin{claimproof}
Suppose $|\mathcal{W}| \geq 2\beta(k,b) + 3$.
As before, we define a bramble of large order to reach a contradiction~--~via the Linkage Function~--~to the fact that $\mathcal{L}'$ is a vital linkage in $G'$.

To this end let $\mathcal{W} = \{ W_1,\dots,W_{w}\}$ with $w \geq 2\beta(k,b) + 3$ and let $w' \coloneqq 2\beta(k,b) + 3$.
Now, for each $i \in [w'],$ let $Q_i \coloneqq W_i \cup C_i'$ and let $\mathcal{Q} \coloneqq \{ Q_i ~\colon~ i\in[w'] \}$.
Then, since the $W_i$ are pairwise vertex disjoint and so are the $C_i',$ but every $W_i$ intersects every $C_j',$ it follows that $\mathcal{Q}$ is a bramble in $G'$ such that no vertex of $G'$ may belong to more than two members of $\mathcal{Q}$.
Hence, the order of $\mathcal{Q}$ is at least $\beta(k,b) + 2$.
As before, \zcref{prop:brambles} now implies that $\mathsf{tw}(G') \geq \beta(k,b) + 1$.
On the other hand, \zcref{claim:VitalAnus} says that $\mathcal{L}'$ is a vital $T$-linkage in $G'$ and thus $\mathsf{tw}(G') \leq \beta(k,b)$ by \zcref{thm:VitalLinkage} which is impossible.
\end{claimproof}

\textbf{Straightening rivers.}
Finally, we reach the point where the actual combing can take place.
Notice that $G'$ still contains an annulus $A''$ with $q$ circles and $r$ rails which is obtained from $A$ through the contractions that transformed $G'''$ into $G'$.

Let $A'$ be the $(2\beta(k,b) + 2)$-trimming of $A'',$ let $\circledcirc'$ be the domain of $A',$ and let $\mathcal{R}$ be the set of rails of $A'$.
Notice that the set of circles of $A'$ is the set $$\mathcal{C}'' \coloneqq \langle C_{2\beta(k,b) + 3}',\dots, C_{a - 2\beta(k,b) - 2}' \rangle.$$
Let us set $a' \coloneqq a - 4\beta(k,b) - 4$ and adjust indices to set $\mathcal{C}'' = \langle C_1'',\dots, C_{a'}'' \rangle$ such that the order of the indices is preserved.

We require two additional railed annuli as follows.
Let $A_1'$ be the railed annulus obtained from $A'$ by first restricting its circles to the cycles $C_1'',\dots,C_{6\beta(k,b) + 6}'',$ and then taking, as rails, from each $R \in \mathcal{R}$ its minimal $C_1''$-$C_{6\beta(k,b) + 6}''$-subpath.
Let $\circledcirc_1'$ be the domain of $A_1'$.
Similarly, let $A_2'$ be the railed annulus obtained from $A'$ by taking as circles the cycles $C_{a' - 6\beta(k,b) - 5}'',\dots, C_{a'}''$ and by taking as rails from each $R \in \mathcal{R}$ its minimal $C_{a' - 6\beta(k,b) - 5}''$-$C_{a'}''$-subpath.
Let $\circledcirc_2'$ be the domain of $A_2'$.

Finally, let $\mathcal{S}'$ be the collection of all $(V(C_1'') \cup V(C_{a'}''))$-subpaths of the paths in $\mathcal{L}'$ which are contained in the crop of $G'$ by $\circledcirc'$.
As before, we may partition $\mathcal{S}'$ into mountains, valleys, and rivers.
However, by our choice of $A'$ and \zcref{claim:FlatMountains} it follows that every member of $\mathcal{S}'$ is a river.
Moreover, $|\mathcal{S}'| \leq 2\beta(k,b) + 2$ by \zcref{claim:FewRivers}.

Let now for each $i\in[2]$ be $\mathcal{S}'_i$ be the restriction of the paths in $\mathcal{S}'$ to the crop of $G'$ by $\circledcirc'_i$~--~let us denote the crop of $G'$ by $\circledcirc'$ by $G_i'$.
Moreover, let $X_1$ be the set of endpoints of the paths in $\mathcal{S}'_1$ on $C_1''$ and let $X_2$ be the set of endpoints of the paths in $\mathcal{S}'_2$ on $C_{a'}''$.
We say that a set of $\mathcal{D}$ rails is \emph{consecutive} if $\mathcal{D} = \{ D_1,\dots,D_{|\mathcal{D}|}\}$ such that $D_i$ and $D_{i+1}$ are consecutive for all $i\in [|\mathcal{D}|-1]$.
Let $Y_1$ be the set of endpoints of $2\beta(k,b) + 2$ consecutive rails of $A_1'$ on $C_{6\beta(k,b) + 6}''$.
Finally, let $Y_2$ be the set of endpoints of the rails of $A_2'$ on $C_{a' - 6\beta(k,b) - 5}''$.

\begin{beautifulclaim}\label{claim:RailTheAnnulus1}
There exists a linkage $\mathcal{F}_1$ of size $|\mathcal{S}'| = |X_1|$ between $X_1$ and $Y_1$ in $G_1'$.
\end{beautifulclaim}

\begin{claimproof}
Suppose there does not exist a $X_1$-$Y_1$-linkage of order $X_1$ in $G_1'$.
Then, by Menger's theorem \zcref{thm:menger}, there exists a set $Z_1 \subseteq V(G_1')$ of size at most $|X_1| - 1$ such that in $G_1' - Z_1$ there is no $X_1$-$Y_1$-path.

Since $|X_1| \leq 2\beta(k,b) + 2$ as a consequence of \zcref{claim:FewRivers} and $|Y_1| = |\mathcal{R}| = r \geq 2\beta(k,b) + 2$ there must exist a rail $R$ of $A_1'$ which completely avoids $Z_1$.
Moreover, there must also exist a river $S \in \mathcal{S}_1'$ which is disjoint from $Z_1$ and finally, there must exist a circle $C$ of $A_1'$ which is disjoint from $Z_1$.
Now, however, $R \cup S \cup C$ is a connected graph which contains a vertex of $X_1$~--~an endpoint of $S$~--~and a vertex of $Y_1$~--~an endpoint of $R$.
This is a contradiction to the choice of $Z_1$ and thus the claim must be true.
\end{claimproof}

Now let $Y_2'$ be the subset of $Y_2$ consisting of precisely those vertices in $Y_2$ that belong to members of $\mathcal{R}'$ that contain an endpoint of some path of $\mathcal{F}_1$.
Then $|Y_2'| = |\mathcal{F}_1| = |X_2|$.
Moreover, let $G_2''$ be the subgraph of $G_2'$ consisting only of the circles of $A_2',$ the rails of $A_2'$ with vertices in $Y_2'',$ and the paths in $\mathcal{S}_2'$.

\begin{beautifulclaim}\label{claim:RailTheAnnulus2}
There exists a linkage $\mathcal{F}_2$ of size $|\mathcal{S}'| = |X_2|$ between $X_2$ and $Y_2'$ in $G_1''$.
Moreover, for every path $S \in \mathcal{S}',$ if $x \in Y_1$ is the endpoint of the path in $\mathcal{F}_1$ containing an endpoint of $S$ and $y \in Y_2'$ is the endpoint of the path in $\mathcal{F}_2$ containing the other endpoint of $S,$ then there exists a rail of $A'$ which contains both $x$ and $y$.
\end{beautifulclaim}

\begin{claimproof}
The proof of this claim follows along two steps.
First we route the paths in $\mathcal{S}'$ starting from $X_2$ onto the rails containing the end of $\mathcal{F}_1$ in $Y_1,$ and then we use the remaining infrastructure in $A_2'$ in order to guarantee that the vertices of $X_2$ are connected to the rails in the way required by the claim.

Let now $B_1$ be the railed annulus obtained from $A_2'$ by taking the cycles $C_{a'-2\beta(k,b)-1}'',\dots,C_{a'}''$ together with the restriction of the rails of $A_2'$ to these cycles.
Let $\circledcirc''_1$ be the domain of $B_1$ and let $H_1$ be the crop of $G'$ by $\circledcirc''_1$.

Moreover, let $B_2$ be the railed annulus obtained from $A_2'$ by first deleting $B_1$ and then iteratively removing vertices of degree $1$.
Let $\circledcirc''_2$ be the domain of $B_2$ and let $H_2$ be the crop of $G'$ by $\circledcirc''_2$.
Notice that $B_2$ has $4\beta(k,b) + 4$ cycles.

We denote by $\mathcal{R}^{\star}$ the set of rails of $A'$ containing vertices from $Y_2'$.
Recall that these rails are precisely those which contain the endpoints of the paths from $\mathcal{F}_1$ in $Y_1$.
Now let $Z_1$ be the set of all ends of those rails of $B_1$ on $C_{a'-2\beta(k,b)-1}''$ which are subpaths of the paths in $\mathcal{R}^{\star}$.
Similarly, let $Z_2$ be the set of ends endpoints of those rails of $B_2$ on the cycle $C_{a'-2\beta(k,b)-2}''$ which are subpaths of the paths in $\mathcal{R}^{\star}$.

Finally, let $\mathcal{F}_2''$ be the collection of all $Z_1$-$Z_2$-subpaths of the paths in $\mathcal{R}^{\star}$.
\medskip

\textbf{Re-routing in $B_1$.}
We now proceed with the first step.
That is, we use arguments analogous to those from the proof of \zcref{claim:RailTheAnnulus1} to find an $X_2$-$Z_1$-linkage $\mathcal{F}_2'$ of size $|X_2|$ in $H_1$.
As before, the existence of this linkage is guaranteed by Menger's theorem \zcref{thm:menger} and the fact that $|X_2| = |Z_1| \leq 2\beta(k,b) + 2$.
\smallskip

Notice that, it is possible to combine the linkages $\mathcal{F}_2'$ and $\mathcal{F}_2''$ along the rails in $\mathcal{R}^{\star}$ to create an $X_2$-$Z_2$-linkage of size $|X_2|$.
All that is left is now to ensure that the second part of our claim is satisfied.
That is, we need to route the vertices in $Z_2$ to the correct vertices in $Y_2'$ as specified by $\mathcal{F}_1$.
See \zcref{fig:SplitAnnulus} for an illustration of the current situation.
\medskip

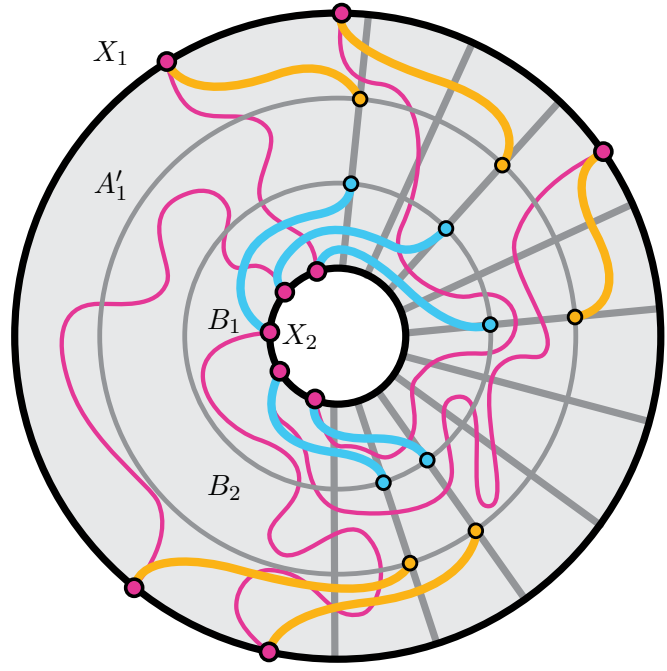
\begin{figure}[ht]
 \centering
 \begin{tikzpicture}

 \pgfdeclarelayer{background}
		\pgfdeclarelayer{foreground}
			
		\pgfsetlayers{background,main,foreground}

 \begin{pgfonlayer}{background}
 \pgftext{\includegraphics[width=9cm]{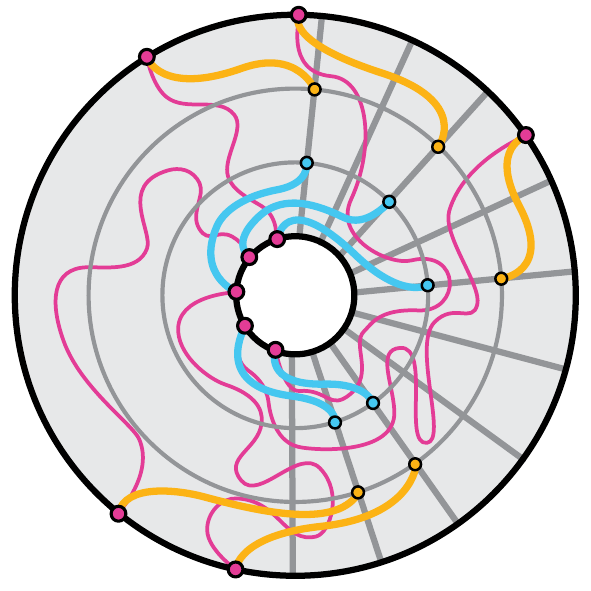}} at (C.center);
 \end{pgfonlayer}{background}
			
 \begin{pgfonlayer}{main}
 \node (C) [v:ghost] {};

 \node (A1) [v:ghost] at (-3,2) {$A_1'$};
 \node (B2) [v:ghost] at (-1.5,-2) {$B_2$};
 \node (B1) [v:ghost] at (-1.5,0.2) {$B_1$};

 \node (X1) [v:ghost] at (-3,3.75) {$X_1$};
 \node (X2) [v:ghost] at (-0.5,0) {$X_2$};
 
 \end{pgfonlayer}{main}
 
 \begin{pgfonlayer}{foreground}
 \end{pgfonlayer}{foreground}

 \end{tikzpicture}
 \caption{An illustration of the situation in the proof of \zcref{thm:CombedAnus} after rerouting the paths in $A_1'$ and $B_1$.
 The figure depicts the situation in the proof of \zcref{claim:RailTheAnnulus2} right before we begin to fix the connections in $B_2$.
 The \textcolor{HotMagenta}{magenta} paths illustrate the original rivers, while the \textcolor{ChromeYellow}{yellow} paths depict the linkage $\mathcal{F}_1$ and the \textcolor{CornflowerBlue}{blue} paths depict the linkage $\mathcal{F}_2'$.}
 \label{fig:SplitAnnulus}
\end{figure}

\textbf{Fixing the connections in $B_2$.}
Recall that $Z_2$ denotes the set of endpoints of the rails of $B_2$ on the cycle $C_{a'-2\beta(k,b)-2}''$ which are subpaths of the paths in $\mathcal{R}^{\star}$.
Moreover, $Y_2'$ denotes the set of endpoints of the rails of $B_2$ on the cycle $C_{a'-6\beta(k,b)-5}''$.
Similarly, $X_2'$ denotes the set of the endpoints of the paths in $\mathcal{F}_2'$ on the cycle $C_{a'}''$ and $Z_1$ denotes the set of endpoints of the paths in $\mathcal{F}_2'$ on the cycle $C_{a'-2\beta(k,b)-1}''$ while $Y_1'$ denotes the endpoints of the paths in $\mathcal{F}_1$ on the cycle $C_{2\beta(k,b) + 2}''$.

Let us now set $r^{\star} \coloneqq |\mathcal{R}^{\star}|$ and $\mathcal{R}^{\star} = \{ R^{\star}_1,\dots,R^{\star}_{r^{\star}} \}$ such that, when restricting the rails of $A'$ to $\mathcal{R}^{\star},$ the rails $R^{\star}_i$ and $R^{\star}_{i+1}$ are consecutive.
We derive from this ordering indices for the elements of $Y_1',$ $Y_2',$ $Z_1,$ and $Z_2$ as follows.
We number $Y_1' = \{ y^1_1,\dots,y^1_{r^{\star}}\}$ and $Y_2' = \{ y^2_1,\dots,y^2_{r^{\star}} \}$ such that $y^j_i$ lies on $R^{\star}_i$ for both $j\in[2]$ and all $i\in[r^{\star}]$.
We also number the vertices in $X_1 = \{ x^1_1,\dots,x^1_{r^{\star}}\}$ such that $x^1_i$ is the endpoint of the path in $\mathcal{F}_1$ that contains $y^1_i$ for each $i\in[r^{\star}]$.
Then, we number the vertices $X_2 = \{x^2_1,\dots,x^2_{r^{\star}}\}$ such that $x^1_i$ and $x^2_i$ belong to the same path, denoted by $S_i,$ of $\mathcal{S}'$ for each $i\in [r^{\star}]$.
Finally, we number $Z_1 = \{ z^1_1,\dots, z^1_{r^{\star}}\}$ and $Z_2 = \{ z^2_1,\dots, z^2_{r^{\star}}\}$ such that $z^1_i$ and $z^2_i$ belong to the path $R^{\star}_i$.

Let us now define the bijection $\pi \colon [r^{\star}] \to [r^{\star}]$ such that for each $i\in[r^{\star}],$ $z^2_i$ lies on $R^{\star}_{\pi(i)}$.

It is now entirely possible that $i \neq \pi(i),$ however, in order to complete our proof we need to find a $Z_2$-$Y_2'$-linkage in $B_2$ which connects $z^2_{\pi(i)}$ to $y^2_i$ for each $i\in[r^{\star}]$.
However, we claim that within $\circledcirc_2'',$ this linkage is topologically feasible.
To see this notice that we know we may use the paths from $\mathcal{F}_1$ and $\mathcal{F}_2'\cup\mathcal{F}_2''$ to redraw the curves obtained from the traces of the paths in $\mathcal{S}'$ to pass through the points of $Y_2'$ and $Z_2$.
It follows now that for each of those curves, the remaining curve between the two points in $Y_2'$ and $Z_2$ may now be confined to be contained entirely in $\circledcirc_2''$ and this is possible while keeping the collection of these curves being a topological linkage.

With the pattern of order $r^{\star}$ defined as $\big\{ \{ y^2_i,z^2_{\pi(i)} \} ~\colon~ i\in[r^{\star}] \big\}$ being topologically feasible in $\circledcirc,$ and the fact that $B_2$ is a railed annulus with at least $2r^{\star}$ circles and at least $2r^{\star}$ rails, it follows from \zcref{lem:two-sided-cylinder} that it is realisable in $B_2,$ leading to a linkage $\mathcal{F}_2'''$.
Now, $\mathcal{F}_2 \coloneqq \mathcal{F}_2' \cup \mathcal{F}_2'' \cup \mathcal{F}_2'''$ is our desired linkage.
\end{claimproof}

Let $Y_1'$ denote the subset of $Y_1$ consisting only of endpoints of paths from $\mathcal{F}_1$.
To complete the proof of the theorem, we now observe that $\mathcal{R}^{\star}$ contains a $Y_1'$-$Y_2'$-linkage $\mathcal{F}_3$ of size $|Y_1'| = |Y_2'| = |\mathcal{R}^{\star}|$ that is internally vertex-disjoint from both $\mathcal{F}_1$ and $\mathcal{F}_2$.
By combining $\mathcal{F}_1,$ $\mathcal{F}_2,$ and $\mathcal{F}_3$ into a single linkage $\mathcal{S}^{\star},$ we have no found a set of rivers that is combed in the $(6\beta(k,b) + 6)$-trimming $A^*$ of $A'$.
This means that, by replacing $\mathcal{S}'$ with $\mathcal{S}^{\star}$ we may obtain a linkage $\mathcal{L}^*$ which is combed in $A^*$ and which satisfies $\tau(\mathcal{L}^*) = \tau(\mathcal{L}')$.
Since $G'$ is a minor of $G$ and $A^*$ is a minor of the $(8\beta(k,b) + 8)$-trimming $A^{\star}$ of $A,$ it follows that $\mathcal{L}^*$ corresponds to a linkage $\mathcal{L}^{\star}$ in $G$ which is combed on $A^{\star}$ and has the same pattern as $\mathcal{L}$.
With this, the proof is complete.
\end{proof}

\subsection{Redrawing topological minors models in partially disc-embedded graphs}
\label{sec:redraw}
With \zcref{thm:CombedAnus} we now have a way to redraw a linkage onto the rails of some subannulus.
Next we need to show that this can be applied in a broader setting as follows:
Suppose some rooted graph $(H,\mathcal{J})$ is a rooted minor of our rooted graph $(G,\mathcal{R})$ and part of $G$ is embedded in a disc $\Delta$ containing, in particular, a large flat subwall of $G$.
Inside of this subwall we may find a big railed annulus $A$ separating an interior area of $\Delta$ from everything of $G$ not drawn into $\Delta$.
From here, the plan is to proceed in two steps.
First, redraw the minor model of $(H,\mathcal{R})$ such that it intersects the domain of $A$ only in few rails and nowhere else.
Second, use $A$ to crop the interior part of the minor model and then redraw it within $A$ itself, thereby avoiding the centre of the wall.

\subsection{Combing a topological minor}
\label{subsec:TopoCombing}
Our first target is a generalisation of \zcref{thm:CombedAnus} to work also with fixed topological minors and not just with vertex-disjoint paths.
As we are going to work with topological minors throughout most of this section, we now spend some time defining the necessary notions.

\paragraph{Topological minors.}
We say that an annotated graph $(J,T)$ is a \emph{tm-pair} if all vertices in $V(J)\setminus T$ have degree precisely $2$ in $J$ and no component of $J-T$ is a path.
Given a graph $G$ and a vertex $v\in V(G)$ of degree $2,$ the operation of \emph{dissolving} $v$ is defined as deleting $v$ and introducing an edge between the two neighbours of $v$ if they were not already adjacent.
We denote by $\mathsf{Diss}(J,T)$ the graph obtained from the tm-pair $(J,T)$ by dissolving all vertices of $V(J\setminus T$.
A \emph{tm-pair} of a graph $G$ is a tm-pair $(J,T)$ such that $J \subseteq G$.
We now say that a graph $H$ is a \emph{topological minor} of a graph $G$ if there exists a tm-pair $(J,T)$ in $G$ such that $\mathsf{Diss}(J,T)$ is isomorphic to $H$.
In this case we say that $(J,T)$ \emph{models} $H$ and $J$ is a \emph{subdivision} of $H$.
The vertices of $T$ are called the \emph{branch vertices} and the paths in $J$ with both endpoints in $T$ but internally vertex-disjoint from $T$ are called the \emph{subdivision paths} of $(J,T),$ finally, the vertices in $V(J)\setminus T$ are called the \emph{subdivision vertices}.
\smallskip

For our railed annuli, we already had a glimpse of the idea what a ``partial'' embedding might look like:
In the case of annuli, we only really needed a sphere-rendition of our graph $G$ with two vortices such that the railed annulus is grounded and the trace of each of its circles separates the two vortices in the sphere.
Now, we go one step further and focus only on graphs that are ``partially'' embedded in a disc.

\paragraph{Partially disc-embedded graphs.}
Let $\Delta$ be a closed disc.
We say that a graph $G$ is \emph{$\Delta$-embedded} (via $\Gamma$) if $\Gamma$ is an embedding of $G$ in $\Delta$ such that the intersection of $\mathsf{bd}(\Delta)$ and $\Gamma$ is a subset of $V(G)$.
A \emph{partial $\Delta$-embedding} of a graph $G$ is a triple $(G_1,G_2,\Gamma)$ where $G=G_1\cup G_2,$ $G_1$ and $G_2$ intersect in vertices only, and $G_1$ is $\Delta$-embedded via $\Gamma$.
We say that $G_1$ is the \emph{$\Delta$-embedded part} of $G$ and we usually identify $G_1$ and $\Gamma$.
When we say that a graph $G$ is \emph{partially $\Delta$-embedded}, we implicitly assume that $G$ comes with a partial $\Delta$-embedding.
A graph $G$ is said to be \emph{partially disc embedded} if there exists a disc $\Delta$ such that $G$ is \emph{partially $\Delta$-embedded}.

In a more general setting, we say that a sphere-rendition $\rho$ of a graph $G$ is a \emph{imperfect rendition in a disc} if $\rho$ has a unique vortex $c_0$.
In such a situation, we refer to the closure of the set obtained from the sphere by deleting $c_0$ as the \emph{disc of $\rho$}.
We sometimes say that $\rho$ is an \emph{imperfect rendition in a disc $\Delta$} to say that $\rho$ is an imperfect rendition in a disc with disc $\Delta$.
Notice that every partial $\Delta$-embedding of some graph $G$ gives rise to an imperfect rendition of $G$ in a disc.

\paragraph{A railed annulus in a disc.}
Let $G$ be a graph with an imperfect rendition in a disc, $\Delta$ be the disc of $\rho,$ and $A \subseteq G$ be a separating railed annulus in $G$ with domain $\circledcirc \subseteq \Delta$.
In such a situation we assume that the circles of $A$ are numbered as $C_1,\dots,C_s$ such that for all $i < j \in [r],$ the $C_i$-disc $\Delta_i$ is contained in the $C_j$-disc $\Delta_j$.

Now let $b\geq 0$ be an integer and let $(G,R)$ be an annotated graph and $\rho$ be an imperfect rendition in a disc $\Delta$ for $G,$ a set $X\subseteq V(G)$ with $R\subseteq X$ is \emph{split in $\Delta$} by a railed annulus $A$ with domain $\circledcirc,$ if
\begin{enumerate}
 \item $A$ is blank,
 \item there exists a set $D \subseteq R$ such that, if $G'$ is the crop of $G$ by $\Delta,$ then $D = G' \cap R,$ and
 \item if $\Delta' \subseteq \Delta$ is the unique component of $\Delta \setminus \circledcirc$ which is homeomorphic to a disc, then $D$ is a subset of the vertices of the crop of $G$ by $\Delta'$.
\end{enumerate}
We say that the set $D$ is the set of \emph{trapped terminals}.

\paragraph{Controlling topological minors.}
Let $(J,T)$ be a tm-pair in an annotated graph $(G,R)$.
The \emph{detail} of $(J,T),$ denoted by $\mathsf{detail}_R(J,T),$ is the integer $d \coloneqq |T \setminus R|$.
Let $\rho$ be an imperfect rendition in a disc $\Delta$ for $G,$ then $(J,T)$ is \emph{carefully split in $\Delta$} by a railed annulus $A$ if $(G,R\cup T)$ is split in $\Delta$ where the set of trapped terminals is a subset of $T \setminus R$.

\begin{figure}[ht]
 \centering
 \begin{tikzpicture}

 \pgfdeclarelayer{background}
		\pgfdeclarelayer{foreground}
			
		\pgfsetlayers{background,main,foreground}

 \begin{pgfonlayer}{background}
 \pgftext{\includegraphics[width=8cm]{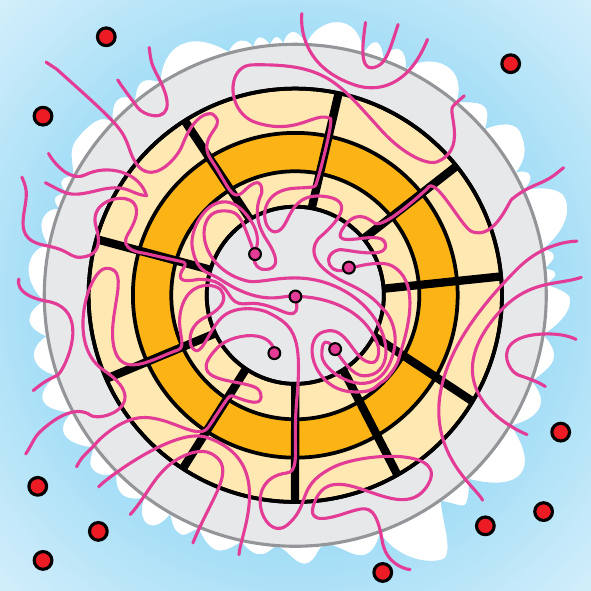}} at (C.center);
 \end{pgfonlayer}{background}
			
 \begin{pgfonlayer}{main}
 \node (C) [v:ghost] {};
 
 \end{pgfonlayer}{main}
 
 \begin{pgfonlayer}{foreground}
 \end{pgfonlayer}{foreground}

 \end{tikzpicture}
 \caption{A sketch of a tm-pair $(J,T)$ with terminals in \textcolor{BostonUniversityRed}{red} -- those are the ones belonging to $R$ -- and \textcolor{HotMagenta}{magenta} with a partial rendition in a disc $\Delta$ depicted in \textcolor{DarkGray}{gray}.
 The \textcolor{HotMagenta}{magenta} curves illustrate the paths of $J$ between the vertices in $T$.
 The \textcolor{ChromeYellow}{yellow} area is the domain of a railed annulus $A$ which carefully splits $(T,J)$.
 The darker shade of \textcolor{ChromeYellow}{yellow} marks the domain of a subannulus $A'$ of $A$ that confines $(T,J)$.}
 \label{fig:CarefulSplit}
\end{figure}

Let now $(J,T)$ be a tm-pair in $(G,R),$ $\rho$ be an imperfect rendition in a disc $\Delta$ for $G,$ and $A$ be a railed annulus in $G$.
We say that $A$ \emph{confines} $(J,T)$ if
\begin{enumerate}
 \item $(J,T)$ is carefully split by $A,$ and
 \item if $\circledcirc$ is the domain of $A$ in $\Delta$ and $G'$ is the crop of $G$ by $\circledcirc,$ then $J \cap G'$ is a collection of subpaths of the rails of $A$.
\end{enumerate}
See \zcref{fig:CarefulSplit} for an illustration of these definitions.
\medskip

What follows is an analogue of the so-called ``Model Combing Lemma'' by Golovach, Stamoulis, and Thilikos \cite{GolovachSG2023Hitting}.
This may be seen as a strengthening of \zcref{thm:CombedAnus} where we replace the linkage with a topological minor model.

\begin{lemma}\label{lemma:CombingAMinor}
There exist functions $f_{\ref{lemma:CombingAMinor}},g_{\ref{lemma:CombingAMinor}},h_{\ref{lemma:CombingAMinor}} \mathbb{N}^{3} \to \mathbb{N}$ such that for all integers $k,b,d \in \mathbb{N},$ graphs $H,$ and annotated graphs $(G,R),$ if
\begin{itemize}
 \item $\mathsf{bidim}(G,R) \leq b,$
 \item $(J,T)$ is a tm-pair in $G$ modelling $H,$
 \item $|T \cap R| \leq k,$
 \item $\mathsf{detail}(J,T) \leq d,$
 \item $G$ has a partial $\Delta$-embedding $(G_1,G_2,\Gamma)$ where $\Delta$ is a closed disc,
 \item $V(G_1) \cap R = \emptyset,$ and
 \item $A = A_{s,r}$ is a railed annulus with $s \geq f_{\ref{lemma:CombingAMinor}}(k,b,d)$ and $r \geq g_{\ref{lemma:CombingAMinor}}(k,b,d)$ such that $A$ carefully splits $(J,T)$.
\end{itemize}
Let $A'$ be the $h_{\ref{lemma:CombingAMinor}}(k,b,d)$-trimming of $A$.
Then there exists a tm-pair $(J',T')$ modelling $H$ in $G$ such that,
\begin{enumerate}
 \item $T' = T,$
 \item $A'$ confines $(J',T'),$ and
 \item if $\circledcirc$ is the domain of $A$ and $G'$ is the crop of $G$ by $\circledcirc,$ then $J - G' \subseteq J' - G'$.
\end{enumerate}
Moreover, $f_{\ref{lemma:CombingAMinor}}(k,b,d),g_{\ref{lemma:CombingAMinor}}(k,b,d),h_{\ref{lemma:CombingAMinor}}(k,b,d) \in 2^{\mathbf{poly(b + d)}} \cdot \mathbf{poly}(k)$.
\end{lemma}

Before we dive into the proof, we need a small auxiliary lemma on the bidimensionality of vertices in partially disc-embedded graphs.

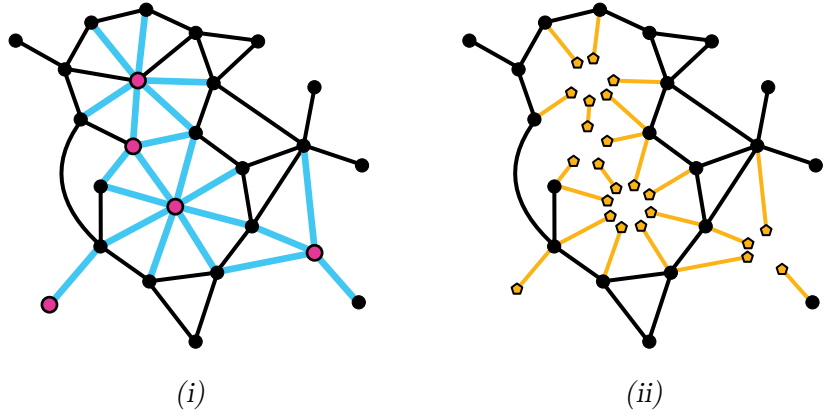
\begin{figure}[ht]
 \centering
 \begin{tikzpicture}

 \pgfdeclarelayer{background}
		\pgfdeclarelayer{foreground}
			
		\pgfsetlayers{background,main,foreground}
			
 \begin{pgfonlayer}{main}
 \node (C) [v:ghost] {};

 \node(L) [v:ghost] at (-3,0) {
 \begin{tikzpicture}

 \pgfdeclarelayer{background}
		 \pgfdeclarelayer{foreground}
			
		 \pgfsetlayers{background,main,foreground}

 \begin{pgfonlayer}{background}
 \pgftext{\includegraphics[width=5cm]{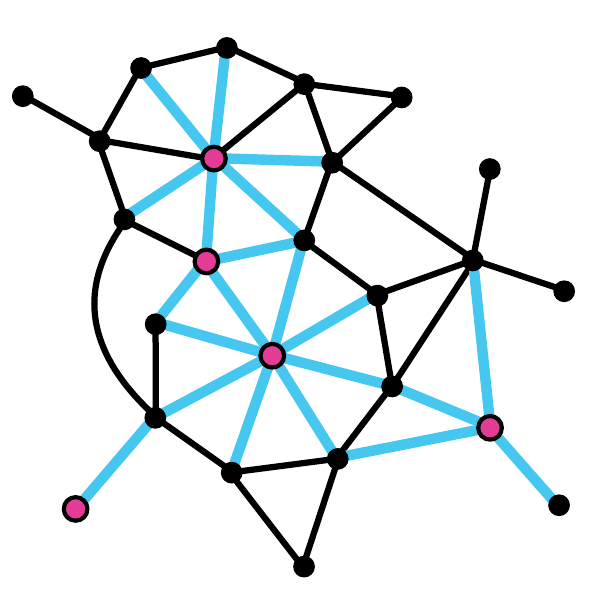}} at (C.center);
 \end{pgfonlayer}{background}
			
 \begin{pgfonlayer}{main}
 \node (C) [v:ghost] {};
 
 \end{pgfonlayer}{main}
 
 \begin{pgfonlayer}{foreground}
 \end{pgfonlayer}{foreground}

 \end{tikzpicture}
 };

 \node(M) [v:ghost] at (3,0) {
 \begin{tikzpicture}

 \pgfdeclarelayer{background}
		 \pgfdeclarelayer{foreground}
			
		 \pgfsetlayers{background,main,foreground}

 \begin{pgfonlayer}{background}
 \pgftext{\includegraphics[width=5cm]{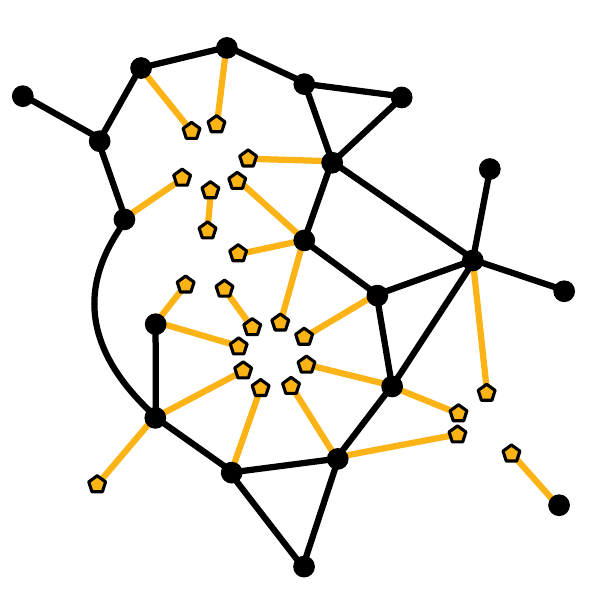}} at (C.center);
 \end{pgfonlayer}{background}
			
 \begin{pgfonlayer}{main}
 \node (C) [v:ghost] {};
 
 \end{pgfonlayer}{main}
 
 \begin{pgfonlayer}{foreground}
 \end{pgfonlayer}{foreground}

 \end{tikzpicture}
 };

 \node (i) [v:ghost] at (-3,-3) {\textsl{(i)}};
 \node (ii) [v:ghost] at (3,-3) {\textsl{(ii)}};

 \end{pgfonlayer}{main}
 
 \begin{pgfonlayer}{foreground}
 \end{pgfonlayer}{foreground}

 \end{tikzpicture}
 \caption{(i) A graph $G$ with an edge set $F$ -- depicted in \textcolor{CornflowerBlue}{blue} -- and vertex set $S$ -- depicted in \textcolor{HotMagenta}{magenta}. (ii) The $F$-splitting of $G$ at $S$ where the residue of $S$ is depicted in \textcolor{ChromeYellow}{yellow}.}
 \label{fig:ResidualVertices}
\end{figure}

Let $G$ be a graph $v\in V(G)$ be a vertex and $F \subseteq E(G)$ be a set of edges incident with $v$.
The \emph{$F$-splitting} of $v$ in $G$ is the graph $G_v$ obtained from $G$ by first subdividing each edge $e \in F$ with a new vertex $x_e$ and then deleting $v$.
Similarly, if $S \subseteq V(G)$ is a set of vertices and $F \subseteq E(G)$ is a set of edges such that each edge in $F$ is incident with some vertex of $S,$ the \emph{$F$-splitting} of $S$ in $G$ is the graph obtained from $G$ by iteratively splitting the vertices in $S$ at their respective edges in $F$.
In case some edge $e$ in $F$ is incident with two vertices of $S,$ this edge will receive two subdivision vertices.
That is, after splitting the first endpoint $u$ of $e = uv$ in the process above and introducing the subdivision vertex $x_e,$ the edge $x_ev$ must be added to $F,$ replacing $e$.
We call the set $X$ of all newly introduced vertices $x_e$ for the $F$-splitting of $S$ in $G$ the \emph{residue of $S$} and we call the elements of $X$ the \emph{residual vertices}.
See \zcref{fig:ResidualVertices} for an illustration.

\begin{lemma}\label{lemma:IncreaseBidimdiscEmbedded}
Let $\Delta$ be a closed disc, $(G,R)$ be an annotated graph and $(G_1,G_2,\Gamma)$ be a partial $\Delta$-embedding of $G$.
Let $ S \subseteq V(G_1)\setminus R$ be a set of vertices not drawn on $\mathsf{bd}(\Delta)$ and let $F\subseteq E(G)_1$ be a set of edges such that each edge in $F$ is incident with some vertex of $S$.
Finally, let $G'$ be the $F$-splitting of $S$ in $G$ and let $X$ be the residue of $S$ in $G'$.
Then, $\mathsf{bidim}(G',R \cup X) \leq \mathsf{bidim}(G, R \cup X) + 2 \leq \mathsf{bidim}(G,R) + |X| +2$.
\end{lemma}

\begin{proof}
Let $s \in S$ be any vertex and let $x_1,\dots,x_p$ be the residual vertices introduced when we split at $s$ in the creating of $G'$.
Notice that we inherit a natural partial $\Delta$-embedding $(G_1',G_2',\Gamma')$ from $(G_1,G_2,\Gamma)$ where $G_2 = G_2'$ and $G_1'$ is the $F$-splitting of $S$ in $G_1$.
Then all of $x_1,\dots,x_p$ sit on a common face in $\Delta$.
It follows that, if $(H,Y)$ is a red-minor of $(G',R\cup X)$ where $H$ is planar and $\varphi$ is a red minor-model of $(H,Y)$ in $G',$ if $P = \{ v \in V(H) ~\colon~ \text{there exists } i\in[p] \text{ s.\@t.\@ } x_i \in \varphi(v) \},$ and $\Gamma_H$ is the drawing of $H$ obtained from $G',$ then there exists a face of $\Gamma_H$ which contains all of $P$.
Hence, in any red grid-minor $J$ of $(G',R \cup X)$ there must exist a facial cycle of $J$ containing all involved vertices among the $x_i$.
Thus, the number of $4$-cycles in any red grid-minor of $(G',R \cup X)$ may exceed the maximum number of $4$-cycles in a red grid-minor of $(G,R)$ by at most $|S|$.
Notice that any grid has at most one facial cycle of length more than $4$~--~its ``outer cycle''~--~and deleting this outer cycle reduces the order of the grid by $2$.
So after removing the outer cycle of any red grid-minor in $(G',R \cup X),$ the total number of $4$-cycles remaining cannot exceed $\mathsf{bidim}(G,R)^2$ by more than $S$.
Thus, $\mathsf{bidim}(G', R \cup X) \leq \mathsf{bidim}(G,R) + |S| + 2$ as desired.
\end{proof}

\begin{proof}[Proof of \zcref{lemma:CombingAMinor}]
We first fix our functions as follows:
\begin{align*}
 f_{\ref{lemma:CombingAMinor}}(k,b,d) & \coloneqq f_{\ref{thm:CombedAnus}}(k + (k+d)^2,b + d +2),\\
 g_{\ref{lemma:CombingAMinor}}(k,b,d) & \coloneqq g_{\ref{thm:CombedAnus}}(k + (k+d)^2,b + d +2)\text{, and}\\
 h_{\ref{lemma:CombingAMinor}}(k,b,d) & \coloneqq h_{\ref{thm:CombedAnus}}(k + (k+d)^2,b + d +2).\\
\end{align*}
With \zcref{lemma:IncreaseBidimdiscEmbedded} we are now able to simply reduce the claim to an application of \zcref{thm:CombedAnus}.
For this recall that $(J,T)$ is a tm-pair in $G$ modelling a graph $H,$ we have that $\mathsf{detail}(J,T) \leq d,$ $(G_1,G_2,\Gamma)$ is a partial $\Delta$-embedding of $G,$ and $A$ is a railed annulus in $G$ with $s \geq f_{\ref{lemma:CombingAMinor}}(k,b,d)$ circles and $r \geq g_{\ref{lemma:CombingAMinor}}(k,b,d)$ rails.

Let now $D$ be the set of trapped terminals for $(G, R \cup T),$ $A,$ and $(G_1,G_2,\Gamma)$.
Then $D \subseteq T \setminus R$ by definition.
Moreover, since the vertices in $D$ are trapped (by $A$) it follows that no vertex of $D$ is drawn on $\mathsf{bd}(\Delta)$.
Let now $F\subseteq E(G) \cap J$ be the set of all edges of $J$ that are incident with some vertex in $D$.
It follows that $|F| \leq (k + d)^2$.

Now consider the $F$-splitting $G'$ of $D$ in $G$ and let $X$ be the set of residual vertices in $G'$.
Notice that $(J,T)$ naturally defines a pattern of order $|E(H)|$ in $G$ by assigning to each subdivision path in $J$ its two endpoints as terminals.
Similarly, in $G',$ we may define a pattern of order $|E(H)|$ on the set $X \cup (T\setminus D)$ which is realised by the linkage $\mathcal{L}$ consisting of the paths in $G'$ obtained from subdivision paths of $(J,T)$ by replacing, in each subdivision path $P,$ the endpoints of $P$ that belong to $D$ with the respective residual vertices.
Finally, we may observe that $\mathsf{bidim}(G', R \cup (T \setminus D) \cup X ) \leq \mathsf{bidim}( G, R \cup (T \setminus D) ) + |D| + 2$ by \zcref{lemma:IncreaseBidimdiscEmbedded}.
Moreover, by \zcref{obs:AddRedVertices} we have that $$\mathsf{bidim}( G, R \cup (T \setminus D) ) + |D| + 2 \leq \mathsf{bidim}( G,R ) + |T\setminus D| + |D| + 2 \leq b + d +2.$$

Now, we apply \zcref{thm:CombedAnus} to $(G', R \cup (T\setminus D) \cup X)$ and $\mathcal{L}$ and $A$.
This yields a linkage $\mathcal{L}'$ in $G'$ with $\tau(\mathcal{L}) = \tau(\mathcal{L}')$ such that $\mathcal{L}'$ is combed in the $h(k,b,d)$-trimming $A'$ of $A$.
Moreover, if $\circledcirc$ is the domain of $A,$ and $G''$ is the crop of $G'$ by $\circledcirc,$ then $\bigcup \mathcal{L}' - G'' \subseteq \bigcup \mathcal{L} - G''$.
Notice that, if we let $\mathcal{L}''$ be the collection of paths obtained from the paths in $\mathcal{L}'$ by replacing each vertex of $X$ with the corresponding vertex of $D,$ then $\mathcal{L}''$ is a collection of internally vertex-disjoint paths in $G$ with all endpoints in $T$.
Indeed, it follows from the construction of $G',$ that $(\bigcup \mathcal{L}'',T)$ is indeed a tm-pair in $G$ modelling $H$ and our proof is complete.
\end{proof}

Notice that, if one was not interested in the bounds on the functions involved or in the explicit dependency on the bidimensionality, \zcref{lemma:CombingAMinor} is precisely the same as the Model Combing Lemma of Golovach, Stamoulis, and Thilikos (Theorem 6 in \cite{GolovachSG2023Hitting}).
Indeed, when considering the bounds abstractly, the main difference is that our proof depends only linearly on the Linkage Function, while the original proof had a quadratic dependency.
This difference is only visible in the proof though as this subtle difference in the degrees is hidden in the $\mathbf{poly}$-notation.

Our main application of \zcref{lemma:CombingAMinor} is to argue that another redrawing-type theorem from the literature~--~this time due to Baste, Sau, and Thilikos \cite{BasteST2023Hitting}~--~is applicable in our setting while maintaining the close-to-optimal bounds to the Linkage Function guaranteed by \zcref{thm:VitalLinkage}.
Before we continue, we need to introduce a hand full of additional definitions in order to state the result of Baste et al.

\paragraph{Contracting respectfully.}
Let $G$ and $H$ be graphs and $Q \subseteq V(H)$ be a set of vertices.
We say that a map $\phi \colon V(H) \to 2^{V(G)}$ is a \emph{$Q$-respecting contraction-mapping of $H$ to $G$} if
\begin{itemize}
 \item $\bigcup_{x \in V(H)} \phi(x) = V(G),$
 \item for all $x \neq y \in V(H)$ it holds that $\phi(x) \cap \phi(y) = \emptyset,$
 \item for each $x \in V(H),$ $\phi(x)$ is connected,
 \item for each $xy \in E(H),$ $\phi(x) \cup \phi(y)$ is connected, and
 \item for each $x \in Q$ it holds that $|\phi(x)| = 1$.
\end{itemize}
 One may interpret the idea behind $Q$-respectful contraction-mappings as an interpolation between topological minors and minors.
 For some special vertices~--~namely those in $Q$~--~we do not allow any form of further contraction, but for those outside $Q$ normal minor rules apply.

\paragraph{Dispersing vertices in a wall.}
Let $G$ be a planar graph with an embedding $\Gamma$ in a disc $\Delta$.
The \emph{outer face} of $\Gamma$ is the unique face that contains $\mathsf{bd}(\Delta)$.
Given two vertices $x$ and $y$ of $G,$ we define their \emph{face-distance} in $\Gamma$ to be the smallest integer $d$ such that there exists a curve $\gamma$ between $x$ and $y$ in $\Delta$ such that $\gamma$ is internally disjoint from the outer face of $\Gamma,$ no internal point of $\gamma$ is a vertex of $\Gamma,$ and $\gamma$ intersects at most $d$ faces of $\Gamma$.

We say that a graph $G$ is \emph{plane} if $G$ comes identified with an embedding $\Gamma$ of $G$ in a disc $\Delta$.
In case $G$ is a plane graph, the \emph{face distance} between two of its vertices is the face distance between those vertices in $\Gamma$.
Given two vertex sets $X,Y \subseteq V(G),$ the \emph{face distance} between $X$ and $Y$ is the minimum integer $d$ such that there are $x\in X$ and $y\in Y$ with face distance at most $d$.
Let $H_1, H_2$ be two subgraphs of $G$.
The \emph{face distance} of $H_1$ and $H_2$ is the face distance of $V(H_1)$ and $V(H_2)$.

Given an integer $d \in \mathbb{N}$ and a vertex $x \in V(G),$ we denote by $\mathsf{F}_G^d(x)$ the set of all vertices of $G$ that are within face distance at most $d$ of $x$.
\smallskip

Let $G$ be a plane graph, $s \in \mathbb{N},$ and a tm-pair $(J,T)$ of $G,$ we say that $(J,T)$ is \emph{safely $s$-dispersed} in $G$ if
\begin{enumerate}
 \item any two distinct vertices $t,t \in T$ are within face distance at least $2s + 1$ in $G,$ and
 \item for every $t \in T$ with $\mathsf{deg}_J(t) = d,$ the graph $J[\mathsf{F}^{s}_G(t) \cap V(J)]$ consists of $d$ paths with $t$ as a unique common endpoint.
\end{enumerate}

\paragraph{Layers and central subwalls.}
Let $r \geq 3$ be an integer and $W$ be an $r$-wall.
Notice that $W$ may be understood as an $(r \times r)$-mesh with a perimeter.
We define the layers of $W$ recursively as follows:
The first \emph{layer} of $W$ is its perimeter.
For $i \in [2,\lfloor \nicefrac{r}{2} \rfloor],$ the \emph{$i$th layer} of $W$ is the $(i-1)$th layer of the $(r-1)$-wall obtained from $W$ by removing its perimeter and iteratively removing vertices of degree $1$.

Given an integer $q \in [3,r]$ such that $r-q$ is even, the \emph{$q$-central subwall} of $W$ is the $q$-wall obtained from $W$ by removing the first $\nicefrac{r-q}{2}$ layers of $W$ and then iteratively deleting vertices of degree at most $1$.

Next, we need a lemma that allows us to redraw planar parts of a tm-pair onto an annulus while also guaranteeing big distances between certain branch vertices and maintaining a fixed interface.

\begin{proposition}[Baste, Sau, and Thilikos \cite{BasteST2023Hitting}]\label{prop:RedrawOntoWall}
There exists a function $f_{\ref{prop:RedrawOntoWall}} \colon \mathbb{N}^3 \to \mathbb{N}$ such that the following holds:

Let $c,r,r',\ell \in \mathbb{N},$ $r' \leq r,$ let $H$ be a $\Delta$-embedded graph with $\ell + r'$ vertices and let $Z = \{ z_1,\dots,z_{r'}\} \subseteq V(H)$ such that
\begin{itemize}
 \item the vertices of $H$ have degree at most $3,$
 \item $Z$ is an independent set of $H,$
 \item all vertices of $Z$ have degree $1$ in $H,$
 \item $\mathsf{bd}(\Delta) \cap H = Z,$ and
 \item $\langle z_1,\dots,z_{r'}\rangle$ is the cyclic ordering of the vertices of $Z$ as they appear in the boundary of $\Delta$.
\end{itemize}
Let also $G$ be a $\Delta'$-embedded graph, $A$ be a railed annulus with $x$ circles $C_1,\dots,C_{x}$ and $y$ rails $R_1,\dots,R_y$ such that $x, y \geq f_{\ref{prop:RedrawOntoWall}}(c,r,\ell) \geq r$.
Let further $w_i$ denote the endpoint of $R_i$ on $C_1$ for each $i\in[r]$ and let $I \coloneqq \{ i_1,\dots,i_{r'}\} \subseteq [r]$.
Then $A$ contains a tm-pair $(J,T)$ such that
\begin{enumerate}
 \item $\mathsf{Diss}(J,T) = H,$
 \item there exists a bijection $\sigma \colon V(H) \to T$ inducing an isomorphism between $\mathsf{Diss}(J,T)$ and $H$ such that for each $j \in [r']$ it holds that $\sigma(z_j) = w_{i_j},$
 \item the tm-pair $(J,T)$ is safely $c$-dispersed in $A,$ and
 \item none of the vertices of $T \setminus \{ w_{i_1},\dots, w_{i_{r'}}\}$ is within face distance less than $c$ from some vertex in $V(C_1) \cup V(C_r)$.
\end{enumerate}
Moreover, it holds that $f_{\ref{prop:RedrawOntoWall}}(c,r,\ell) \in \mathcal{O}(cr (\ell + r))$.
\end{proposition}

\begin{figure}[ht]
 \centering
 \begin{tikzpicture}

 \pgfdeclarelayer{background}
		\pgfdeclarelayer{foreground}
			
		\pgfsetlayers{background,main,foreground}
			
 \begin{pgfonlayer}{main}
 \node (C) [v:ghost] {};

 \node(L) [v:ghost] at (-3.5,0) {
 \begin{tikzpicture}

 \pgfdeclarelayer{background}
		 \pgfdeclarelayer{foreground}
			
		 \pgfsetlayers{background,main,foreground}

 \begin{pgfonlayer}{background}
 \pgftext{\includegraphics[width=6cm]{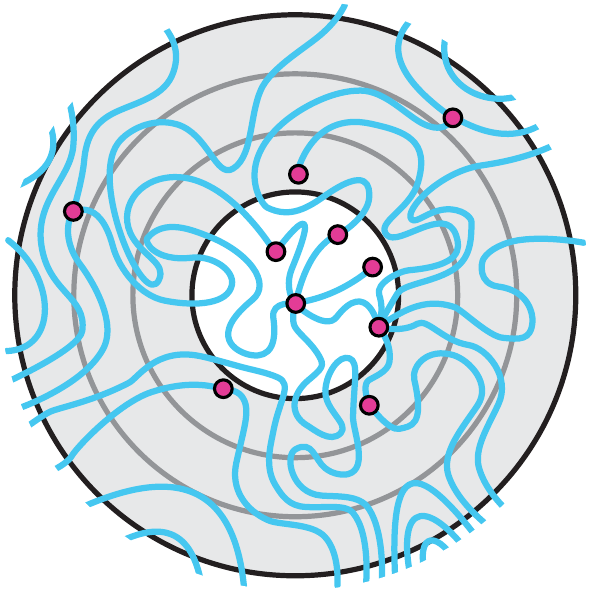}} at (C.center);
 \end{pgfonlayer}{background}
			
 \begin{pgfonlayer}{main}
 \node (C) [v:ghost] {};
 
 \end{pgfonlayer}{main}
 
 \begin{pgfonlayer}{foreground}
 \end{pgfonlayer}{foreground}

 \end{tikzpicture}
 };

 \node(M) [v:ghost] at (3.5,0) {
 \begin{tikzpicture}

 \pgfdeclarelayer{background}
		 \pgfdeclarelayer{foreground}
			
		 \pgfsetlayers{background,main,foreground}

 \begin{pgfonlayer}{background}
 \pgftext{\includegraphics[width=6cm]{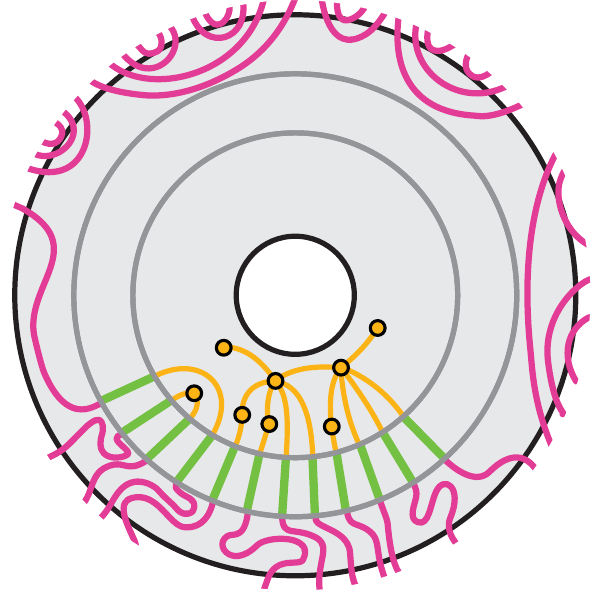}} at (C.center);
 \end{pgfonlayer}{background}
			
 \begin{pgfonlayer}{main}
 \node (C) [v:ghost] {};
 
 \end{pgfonlayer}{main}
 
 \begin{pgfonlayer}{foreground}
 \end{pgfonlayer}{foreground}

 \end{tikzpicture}
 };

 \node (i) [v:ghost] at (-3.5,-3.5) {\textsl{(i)}};
 \node (ii) [v:ghost] at (3.5,-3.5) {\textsl{(ii)}};

 \end{pgfonlayer}{main}
 
 \begin{pgfonlayer}{foreground}
 \end{pgfonlayer}{foreground}

 \end{tikzpicture}
 \caption{(i) The input for \zcref{thm:RedrawAndDisperseTM}: A tm-pair partially drawn in the disk-embedded part of a partially disk-embedded graph together with a railed annulus separating the branch vertices -- depicted in \textcolor{HotMagenta}{magenta}. (ii) The outcome of \zcref{thm:RedrawAndDisperseTM} applied to the input depicted in (i): the tm-pair now coincides with the rails of the central part of the annulus -- the \textcolor{AppleGreen}{green} paths -- while the inner part of the tm-pair -- depicted in \textcolor{ChromeYellow}{yellow} has been fully redrawn into the inner part of the annulus and the outer part -- indicated in \textcolor{HotMagenta}{magenta} -- has been tamed to avoid the central and the inner layer of our annulus.}
 \label{fig:RedrawModel}
\end{figure}
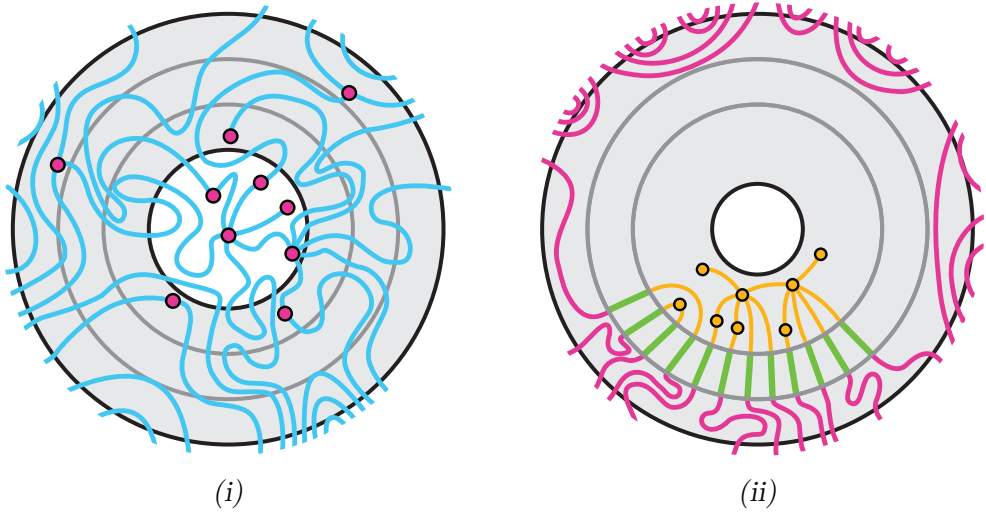

We are finally ready to state our variant of the tm-pair redrawing lemma of Baste, Sau, and Thilikos \cite{BasteST2023Hitting} needed for confining and dispersing partially disc embedded tm-pairs inside a railed annulus.
See \zcref{fig:RedrawModel} for an illustration.

\begin{theorem}\label{thm:RedrawAndDisperseTM}
There exist three functions $f_{\ref{thm:RedrawAndDisperseTM}} \colon \mathbb{N}^5 \to \mathbb{N},$ $g_{\ref{thm:RedrawAndDisperseTM}} \colon \mathbb{N}^3 \to \mathbb{N},$ and $h_{\ref{thm:RedrawAndDisperseTM}} \colon \mathbb{N}^4 \to \mathbb{N}$ such that the following holds:

Let $s,b,d,k,z \in \mathbb{N}$ with $z \geq 3$ being odd.
Let $\Delta$ be a disc and $(G,R)$ be an annotated graph satisfying $\mathsf{bidim}(G,R) \leq b$ with a partial $\Delta$-embedding $(G_1,G_2,\Gamma)$ such that
\begin{itemize}
 \item $V(G_1) \cap R = \emptyset,$
 \item $G_1$ contains an $f_{\ref{thm:RedrawAndDisperseTM}}(s,z,b,d,k)$-wall $W$ whose perimeter separates its interior from $\mathsf{bd}(\Delta)$ in $\Delta,$
 \item $\mathcal{C} = \{ C_1,\dots,C_{g_{\ref{thm:RedrawAndDisperseTM}}(d,b,k) + h_{\ref{thm:RedrawAndDisperseTM}}(s,d,k,b)} \}$ are the first $g_{\ref{thm:RedrawAndDisperseTM}}(d,b,k) + h_{\ref{thm:RedrawAndDisperseTM}}(s,d,k,b)$ layers of $W$ indexed in the natural order, and
 \item $\Delta_1,\dots,\Delta_{g_{\ref{thm:RedrawAndDisperseTM}}(d,b,k) + h_{\ref{thm:RedrawAndDisperseTM}}(s,d,k,b)} \subseteq \Delta$ denote the open discs defined by these layers.
\end{itemize}

If $(J,T)$ is a tm-pair of $G$ with $\mathsf{detail}(J,T) \leq d$ and $|R \cap T| \leq k$ where $Q \subseteq T$ is a set containing vertices of degree at most $3$ in $J,$ then there exists a tm-pair $(J',T')$ of $G$ and an integer $q \in [g_{\ref{thm:RedrawAndDisperseTM}}(d,b,k)]$ such that
\begin{enumerate}
 \item $J' - (V(G_1) \cap \Delta_q)$ is a subgraph of $J - (V(G_1) \cap \Delta_q),$
 \item if $\circledcirc \subseteq \Delta$ denotes the annulus with boundary components $C_q$ and $C_{q + h_{\ref{thm:RedrawAndDisperseTM}}(s,d,k,b)-1}$ and $G'$ denotes the crop of $G_1$ by $\circledcirc,$ then $V(G') \cap (T \cup T') = \emptyset,$
 \item $(J' \cap (G_1 \cap \mathsf{cl}( \Delta_{q + h_{\ref{thm:RedrawAndDisperseTM}}(s,d,k,b)} )), T' \cap (V(G_1) \cap \mathsf{cl}( \Delta_{q + h_{\ref{thm:RedrawAndDisperseTM}}(s,d,k,b)} )))$ is a tm-pair of $W$ that is safely $s$-dispersed in $W$ and none of the vertices of $T' \cap (V(G_1) \cap \mathsf{cl}( \Delta_{q + h_{\ref{thm:RedrawAndDisperseTM}}(s,d,k,b)} ))$ is within face distance less than $s$ in $W$ from some vertex of $C_{q + h_{\ref{thm:RedrawAndDisperseTM}}(s,d,k,b)} \cup C_{q + g_{\ref{thm:RedrawAndDisperseTM}}(d,b,k) + h_{\ref{thm:RedrawAndDisperseTM}}(s,d,k,b)},$
 \item $J' \cap (G_1 \cap \Delta_{g_{\ref{thm:RedrawAndDisperseTM}}(d,b,k)} + h_{\ref{thm:RedrawAndDisperseTM}}(s,d,k,b)) = \emptyset,$
 \item the compass of the $z$-central subwall of $W$ is contained in $\Delta_{g_{\ref{thm:RedrawAndDisperseTM}}(d,b,k) + h_{\ref{thm:RedrawAndDisperseTM}}(s,d,k,b)},$ and
 \item there is a $Q$-respecting contraction mapping of $\mathsf{Diss}(J,T)$ to $\mathsf{Diss}(J',T')$.
\end{enumerate}
Moreover, it holds that $f_{\ref{thm:RedrawAndDisperseTM}}(s,z,d,b,k) \in 2^{\mathbf{poly}(b + d)} \cdot \mathbf{poly}(s + z + k)$ and $g_{\ref{thm:RedrawAndDisperseTM}}(d,b,k) \in 2^{\mathbf{poly}(b + d)} \cdot \mathbf{poly}(k)$ and,$h_{\ref{thm:RedrawAndDisperseTM}}(s,d,k,b) \in 2^{\mathbf{poly}(b + d)} \cdot \mathbf{poly}(s + k)$.
\end{theorem}

\begin{proof}
We open by defining our functions.
\begin{align*}
 f_{\ref{thm:RedrawAndDisperseTM}}(s,z,b,d,k) \coloneqq ~ & 2((d+2) \cdot f_{\ref{lemma:CombingAMinor}}(k,b,d) + f_{\ref{prop:RedrawOntoWall}}(s,k,d) + 2g_{\ref{lemma:CombingAMinor}}(k,b,d))\\
 & + (3s + 1) \cdot (g_{\ref{lemma:CombingAMinor}}(k,b,d)) + f_{\ref{prop:RedrawOntoWall}}(s,f_{\ref{lemma:CombingAMinor}}(k,b,d),d+2d^2) + z + 1\\
 g_{\ref{thm:RedrawAndDisperseTM}}(d,b,k) \coloneqq ~ & (d+2) \cdot f_{\ref{lemma:CombingAMinor}}(k,b,d) + 2g_{\ref{lemma:CombingAMinor}}(k,b,d)\\
 h_{\ref{thm:RedrawAndDisperseTM}}(s,d,k,b) \coloneqq ~ & f_{\ref{prop:RedrawOntoWall}}(s,f_{\ref{lemma:CombingAMinor}}(k,b,d),d+2d^2)
\end{align*}
With this, the claimed bounds on the functions follow directly from \zcref{lemma:CombingAMinor} and \zcref{prop:RedrawOntoWall}.

Let us now define $C_1,\dots, C_{(d+2) \cdot f_{\ref{lemma:CombingAMinor}}(k,b,d) + h_{\ref{thm:RedrawAndDisperseTM}}(s,d,k)+ 2g_{\ref{lemma:CombingAMinor}}(k,b,d)}$ the first $(d+2) \cdot f_{\ref{lemma:CombingAMinor}}(k,b,d) + h_{\ref{thm:RedrawAndDisperseTM}}(s,d,k,b)+ 2g_{\ref{lemma:CombingAMinor}}(k,b,d)$ layers of $W$.
for each $i\in[d+1]$ define $\circledcirc_i \subseteq \Delta$ to be the annulus with boundary components $C_{(i-1) \cdot (f_{\ref{lemma:CombingAMinor}}(k,b,d)) + 1}$ and $C_{i \cdot f_{\ref{lemma:CombingAMinor}}(k,b,d)}$ and let $H_i$ be the crop of $G$ by $\circledcirc_i$.
We say that $C_{(i-1)(\cdot f_{\ref{lemma:CombingAMinor}}(k,b,d)) + 1}$ is the \emph{outer cycle} of $\circledcirc_i$ and $C_{i \cdot \cdot f_{\ref{lemma:CombingAMinor}}(k,b,d)}$ is the \emph{inner cycle} of $\circledcirc_i$.

Notice that there exist $(3s + 1) \cdot (g_{\ref{lemma:CombingAMinor}}(k,b,d)) + f_{\ref{prop:RedrawOntoWall}}(s,f_{\ref{lemma:CombingAMinor}}(k,b,d),d+2d^2)$ vertical paths of $R_1'' \dots R_{(3s + 1) \cdot (g_{\ref{lemma:CombingAMinor}}(k,b,d)) + f_{\ref{prop:RedrawOntoWall}}(s,f_{\ref{lemma:CombingAMinor}}(k,b,d),d+2d^2)}''$ of $W$ such that for each $i \in [(3s + 1) \cdot (g_{\ref{lemma:CombingAMinor}}(k,b,d)) + f_{\ref{prop:RedrawOntoWall}}(s,f_{\ref{lemma:CombingAMinor}}(k,b,d),d+2d^2)]$ and $j \in [(d+2) \cdot f_{\ref{lemma:CombingAMinor}}(k,b,d) + h_{\ref{thm:RedrawAndDisperseTM}}(s,d,k,b)+ 2g_{\ref{lemma:CombingAMinor}}(k,b,d)],$ $R''_i \cap C_{j}$ has precisely two components and each such component is a subpath of both $R''_{i}$ and $C_{j}$.
From each $R''_i$ select $R_i'$ to be a minimal $V(C_1)$-$V(C_{(d+2) \cdot f_{\ref{lemma:CombingAMinor}}(k,b,d) + h_{\ref{thm:RedrawAndDisperseTM}}(s,d,k,b)}+ 2g_{\ref{lemma:CombingAMinor}}(k,b,d))$-subpath of $R_i''$ and assume, without loss of generality, that the paths $$R_1',\dots,R_{(3s + 1) \cdot (g_{\ref{lemma:CombingAMinor}}(k,b,d)) + f_{\ref{prop:RedrawOntoWall}}(s,f_{\ref{lemma:CombingAMinor}}(k,b,d),d+2d^2)}'$$ are numbered according to the clockwise order of the appearance of their endpoints in $C_1$.
Then, for every $j \in [g_{\ref{lemma:CombingAMinor}}(k,b,d)]$ let
\begin{align*}
 a_j \coloneqq 3js + 1,
\end{align*}
Now, for each $j \in [g_{\ref{lemma:CombingAMinor}}(k,b,d)]$ let
\begin{align*}
 R_j \coloneqq R'_{a_j}
\end{align*}
and notice that, in $W,$ any two of the $R_j$'s are at face distance at least $3s$ from each other.
Moreover, $\bigcup_{i \in[(d+2) \cdot f_{\ref{lemma:CombingAMinor}}(k,b,d) + h_{\ref{thm:RedrawAndDisperseTM}}(s,d,k,b)]}C_i \cup \bigcup_{j \in [g_{\ref{lemma:CombingAMinor}}(k,b,d)]}$ forms a railed annulus $A$ with $g_{\ref{lemma:CombingAMinor}}(k,b,d)$ rails and $(d+2) \cdot f_{\ref{lemma:CombingAMinor}}(k,b,d) + h_{\ref{thm:RedrawAndDisperseTM}}(s,d,k,b)$ circles.

For each $i \in[d+1]$ let $A_i$ denote the railed annulus with $f_{\ref{lemma:CombingAMinor}}(k,b,d)$ circles and $g_{\ref{lemma:CombingAMinor}}(k,b,d)$ rails obtained by cropping $A$ by $\circledcirc_i$.

\paragraph{Combing the model.}
Since $\mathsf{detail}(J,T) \leq d$ and $V(G_1) \cap R = \emptyset,$ there must exist $i\in[d+1]$ such that $V(H_i) \cap T = \emptyset$.
Let $\Delta_i$ denote the disc in $\Delta$ bounded by the outer cycle of $\circledcirc_i$.
Moreover, let $G_1'$ be the crop of $G_1$ by $\Delta_i$ and let $G_2'$ be the graph obtained from $G$ by deleting all vertices and edges drawn on the interior of $\Delta_i$.
Finally, let $\Gamma'$ be the restriction of $\Gamma$ to $\Delta_i$.
Then $(G_1',G_2',\Gamma')$ is a $\Delta_i$-embedded graph.
Notice now that
\begin{itemize}
 \item $\mathsf{bidim}(G,R) \leq b,$
 \item $(J,T)$ is a tm-pair in $G$ modelling $\mathsf{Diss}(J,T),$
 \item $|T \cap R| \leq k,$
 \item $\mathsf{detail}(J,T) \leq d,$
 \item $G$ has a partial $\Delta_i$-embedding $(G_1',G_2',\Gamma)$ where $\Delta'$ is a closed disc,
 \item $V(G_1') \cap R = \emptyset,$ and
 \item $A_i$ is a railed annulus with $s \geq f_{\ref{lemma:CombingAMinor}}(k,b,d)$ circles and $r \geq g_{\ref{lemma:CombingAMinor}}(k,b,d)$ rails such that $A_i$ carefully splits $(J,T)$.
\end{itemize}
Hence, we may now apply \zcref{lemma:CombingAMinor} to $(J,T)$ and $A_i$ in $G$.
As a result we get that, if $A_i'$ denotes the $h_{\ref{lemma:CombingAMinor}}(k,b,d)$-trimming of $A$ and $\circledcirc_i'$ denotes the domain of $A_i',$ then there exists a tm-pair $(\widetilde{J},\widetilde{T})$ modelling $\mathsf{Diss}(J,T)$ in $G$ such that,
\begin{enumerate}
 \item $\widetilde{T} = T,$
 \item $A_i'$ confines $(\widetilde{J},\widetilde{T}),$ and
 \item $J - H_i' \subseteq \widetilde{J} - H_i'$ where $H_i'$ denotes the crop of $G_1$ by $\circledcirc_i'$.
\end{enumerate}
For each $j \in [g_{\ref{lemma:CombingAMinor}}(k,b,d)]$ let $\widehat{R}_j$ denote the subpath of $R_j$ in $A_i'$.
We say that $\widehat{R}_j$ is \emph{used} if $\widehat{R}_j \subseteq J'$ and we denote by $I \subseteq [g_{\ref{lemma:CombingAMinor}}(k,b,d)]$ the set of all $j \in [g_{\ref{lemma:CombingAMinor}}(k,b,d)]$ such that $\widehat{R}_j$ is used.
For every $j \in I$ we denote by $u_j$ the endpoint of $R_j$ on $C_{(d+2) \cdot f_{\ref{lemma:CombingAMinor}}(k,b,d)}$.

\paragraph{Redrawing the model.}
Now let $\widehat{T}' \subseteq \widetilde{T}$ be the set of all terminals of $(\widetilde{J},\widetilde{T})$ trapped by $A_i'$.
Furthermore, let $\widehat{J}$ be the crop of $\widetilde{J}$ by the disc $\widehat{\Delta}$ whose boundary coincides with the inner cycle of $\circledcirc_i'$~--~that is the unique boundary component of $\circledcirc'$ which is separated from $\mathsf{bd}(\Delta)$ by the other boundary component of $\circledcirc'_i$.
Let $\pi\colon [|I|] \to I$ be a bijection such that $x < y \in I$ implies that $\pi^{-1}(x) < \pi^{-1}(y)$.
For each $j \in [|I|]$ let $z_j$ be the end of $\widehat{R}_{\pi(j)}$ on the inner cycle of $\circledcirc'_i$ and let $Z \coloneqq \{ z_j ~\colon~ j\in[|I|]\}$.
Finally set $\widehat{T} \coloneqq \widehat{T}' \cup Z$.
Now notice that $(\widehat{J},\widehat{T})$ is a tm-pair in $G$.
Indeed, if we let $G_1''$ be the crop of $G_1$ by $\widehat{\Delta},$ then $G_1''$ is $\widehat{\Delta}$-embedded and $(\widehat{T},\widehat{J})$ is a tm-pair in $G_1''$.
\smallskip

Let us briefly pause and observe that $\Delta(\widehat{J})$ might be larger than $3$.
Indeed, currently, we have that $\widehat{J}$ itself is a $\widehat{\Delta}$-embedded graph intersecting\footnote{Notice that this claim is not entirely true, but we may assume to slightly shift all other vertices and edges on the boundary away from it as they play a less crucial role.} $\mathsf{bd}(\widehat{\Delta})$ precisely in $Z$.
In the next step we create a new tm-pair $(J^{\star},T^{\star})$ from $(\widehat{J},\widehat{T})$ by iteratively ``splitting'' vertices while maintaining 
\begin{enumerate}
 \item that $(J^{\star},T^{\star})$ is $\widehat{\Delta}$-embedded,
 
 \item that $J^{\star}$ intersects $\mathsf{bd}(\widehat{\Delta})$ precisely in $Z,$
 
 \item that $|T^{\star} \setminus Z| \leq |\widehat{T}| + 2|E(\mathsf{Diss}(\widehat{J},\widehat{T}))|,$ and

 \item $\mathsf{Diss}(\widehat{J},\widehat{T})$ is a minor of $J^{\star}$.
\end{enumerate}
Before we start the process, notice that each vertex of $Z$ has degree $1$ in $\widehat{J}$ by construction.
Now, set $\widehat{J}_{0} \coloneqq \widehat{J}$ and $\widehat{T}_0 \coloneqq \widehat{T}$.
As long as there exists $v \in \widehat{J}_i$ with $\mathsf{deg}_{\widehat{J}_i}(v) \geq 4$ we partition $N_{\widehat{J}_i}$ into sets $X$ and $Y$ such that $|X| = 2$ and $|Y| = \mathsf{deg}_{\widehat{J}_i}(v) - 2$.
Then we delete $v$ and instead introduce two vertices $v_1$ and $v_2$ together with edges such that $v_1$ is adjacent exactly to the vertices in $X$ and $v_2,$ and $v_2$ is adjacent exactly to the vertices in $Y$ and $v_1$.
We then denote the resulting graph by $\widehat{J}_{i+1}$.
By our choice of $v$ it must be true that $v \in \widehat{T}_i$ and now we replace $v \in \widehat{T}_i$ by $v_1$ and $v_2$ to create the set $\widehat{T}_{i+1}$.
Notice that in each step the number of vertices increases by $1,$ the number of edges increases by $1,$ and the degree of some vertex decreases by $1$.
Hence, after a bounded number of iterations, we reach a point where $J^{\star} = \widehat{J}_i$ has maximum degree $3$ and $|T^{\star} \setminus Z| \leq |\widehat{T}| + 2|E(\mathsf{Diss}(\widehat{J},\widehat{T}))|$ where $T^{\star} \coloneqq \widehat{T}_i$ as desired.

Notice also that the set $Q \subseteq T$ was a set of vertices of degree at most $3$ in $J$ and so we still have that $Q \cap \widehat{T}' \subseteq T^{\star}$ as $J - H_i' \subseteq \widetilde{J} - H_i'$.
\smallskip

Now require two additional railed annuli.
One in order to apply \zcref{prop:RedrawOntoWall} and the other in order to reconnect the redrawn model into the used rails of $A_i'$.
Recall that $C_{i \cdot f_{\ref{lemma:CombingAMinor}}(k,b,d)}$ is the inner cycle of $\circledcirc'$.

Let $B_1$ denote the railed annulus with circles $C_{i \cdot f_{\ref{lemma:CombingAMinor}}(k,b,d) + 1},\dots,C_{i \cdot f_{\ref{lemma:CombingAMinor}}(k,b,d) + 2g_{\ref{lemma:CombingAMinor}}(k,b,d))}$ and rails $P_1,\dots,P_{(3s + 1) \cdot (g_{\ref{lemma:CombingAMinor}}(k,b,d)) + f_{\ref{prop:RedrawOntoWall}}(s,f_{\ref{lemma:CombingAMinor}}(k,b,d),d+2d^2)}$ obtained by restricting the $R_j'$ to their minimal $C_{i \cdot f_{\ref{lemma:CombingAMinor}}(k,b,d) + 1}$-$C_{2g_{\ref{lemma:CombingAMinor}}(k,b,d))}$-subpaths.
Let $\circledcirc^1$ denote the domain of $B^1$.

Similarly, let $B_2$ denote the railed annulus with circles 
\begin{align*}
C_{i \cdot f_{\ref{lemma:CombingAMinor}}(k,b,d) + 2g_{\ref{lemma:CombingAMinor}}(k,b,d)) + 1},\dots,C_{i \cdot f_{\ref{lemma:CombingAMinor}}(k,b,d) + 2g_{\ref{lemma:CombingAMinor}}(k,b,d)) + f_{\ref{prop:RedrawOntoWall}}(s,f_{\ref{lemma:CombingAMinor}}(k,b,d),d+2d^2)}
\end{align*}
and rails $Q_1,\dots,Q_{(3s + 1) \cdot (g_{\ref{lemma:CombingAMinor}}(k,b,d)) + f_{\ref{prop:RedrawOntoWall}}(s,f_{\ref{lemma:CombingAMinor}}(k,b,d),d+2d^2)}$ obtained by restricting the $R_j'$ to their minimal $C_{ \cdot f_{\ref{lemma:CombingAMinor}}(k,b,d) + 2g_{\ref{lemma:CombingAMinor}}(k,b,d)) + 1}$-$C_{i \cdot f_{\ref{lemma:CombingAMinor}}(k,b,d) + 2g_{\ref{lemma:CombingAMinor}}(k,b,d)) + f_{\ref{prop:RedrawOntoWall}}(s,f_{\ref{lemma:CombingAMinor}}(k,b,d),d+2d^2)}$-subpaths.

We have now gained the right to apply \zcref{prop:RedrawOntoWall} to $J^{\star}$ and the annulus $B_2$.

Let $w_j$ denote the endpoint of $Q_j$ on $C_{i \cdot f_{\ref{lemma:CombingAMinor}}(k,b,d) + 2g_{\ref{lemma:CombingAMinor}}(k,b,d)) + 1}$ for each $j\in[i \cdot f_{\ref{lemma:CombingAMinor}}(k,b,d) + 2g_{\ref{lemma:CombingAMinor}}(k,b,d)) + f_{\ref{prop:RedrawOntoWall}}(s,f_{\ref{lemma:CombingAMinor}}(k,b,d),d+2d^2)]$ and recall our set $I$ of indices.
We now gain from \zcref{prop:RedrawOntoWall} that $B_2$ contains a tm-pair $(J^{\dagger},T^{\dagger})$ such that
\begin{enumerate}
 \item $\mathsf{Diss}(J^{\dagger},T^{\dagger}) = J^{\star},$
 \item there exists a bijection $\sigma \colon V(J^{\star}) \to T^{\dagger}$ inducing an isomorphism between $\mathsf{Diss}(J^{\dagger},T^{\dagger})$ and $J^{\star}$ such that for each $j \in [|I|]$ it holds that $\sigma(z_{\pi^{-1}(j)}) = w_{j},$
 \item the tm-pair $(J^{\dagger},T^{\dagger})$ is safely $s$-dispersed in $B_2,$ and
 \item none of the vertices of $T^{\dagger} \setminus \{ w_i ~\colon~ j \in I \}$ is within face distance less than $s$ from some vertex in $V(C_{i \cdot f_{\ref{lemma:CombingAMinor}}(k,b,d) + 2g_{\ref{lemma:CombingAMinor}}(k,b,d)) + 1}) \cup V(C_{i \cdot f_{\ref{lemma:CombingAMinor}}(k,b,d) + 2g_{\ref{lemma:CombingAMinor}}(k,b,d)) + f_{\ref{prop:RedrawOntoWall}}(s,f_{\ref{lemma:CombingAMinor}}(k,b,d),d+2d^2)})$.
\end{enumerate}
See (ii) in \zcref{fig:RedrawModel}: Here, the \textcolor{ChromeYellow}{yellow} part illustrates the tm-pair $(J^{\dagger},T^{\dagger})$ now rooted on the rails of the annulus.

\paragraph{Reconnecting the model.}
The final step is to simply call upon \zcref{lem:outerplanar-cylinder} to connect the vertices $z_{\pi^{-1}(j)}$ to the vertices $w_j$ within $B_1$.
As this linkage is clearly topologically feasible within the domain of $B_1,$ \zcref{lem:outerplanar-cylinder} yields the existence of the desired linkage $\mathcal{L}$.

Let $(J^{\leftmoon},T^{\leftmoon})$ be obtained by adding the paths in $\mathcal{L}$ to $J^{\dagger}$ and replacing the vertices $w_j$ in $T^{\dagger}$ with the vertices $z_{\pi^{-1}(j)}$ for all $j \in I$.

Let $G_2''$ be obtained from $G$ by deleting all vertices and edges of $G_1''$ drawn in the interior of $\widehat{\Delta}$.
Now let $J^{\rightmoon}$ be the intersection of $\widetilde{J}$ with $G_2''$ and let $T^{\rightmoon}$ be the intersection of $T$ with the vertices of $G_2''$.
Finally, let $J^{\fullmoon} \coloneqq J^{\leftmoon} \cup J^{\rightmoon}$ and $T^{\fullmoon} \coloneqq T^{\leftmoon} \cup T^{\rightmoon}$.
It is easy to see that $(J^{\fullmoon},T^{\fullmoon})$ is now a tm-pair of $G$ such that no vertex or edge of $J^{\fullmoon}$ is drawn in or intersects the interior of the disc in $\Delta$ bounded by $C_{i \cdot f_{\ref{lemma:CombingAMinor}}(k,b,d) + 2g_{\ref{lemma:CombingAMinor}}(k,b,d)) + f_{\ref{prop:RedrawOntoWall}}(s,f_{\ref{lemma:CombingAMinor}}(k,b,d),d+2d^2)}$
Moreover, $\mathsf{Diss}(J)$ is a minor of $J^{\fullmoon}$ and, in particular, there exists a $Q$-respecting contracting mapping of $\mathsf{Diss}(J,Y)$ to $\mathsf{Diss}(J^{\fullmoon},T^{\fullmoon})$ by construction.
From here it is straightforward to verify that $(J^{\fullmoon},T^{\fullmoon})$ is indeed the tm-pair required by our assertion.
\end{proof}

\section{An irrelevant vertex for folio}\label{sec:irrelevantvertex}
In this section we finally prove the core statement of the Graph Minor Algorithm, that is \zcref{thm:Main2Folio}.
Towards this goal we need a number of additional definition and intermediate results.
In broad strokes, what is left to do at this point is to first, explaining homogeneous flat walls in the abstract and then describing the set of colours which defines the homogeneity.
Second, we explain how \zcref{thm:RedrawAndDisperseTM} may be used to show that any minor model invading deep into a flat homogeneous wall may be redrawn in a way that fully avoids its centre.
Third, we require a way to efficiently find a flat wall whose compass has bounded treewidth and which allows us to efficiently compute the colours required in the first and second step.
Combining the outcomes of all three steps with the tools from \zcref{sec:largeclique} for handling clique-minors then finally results in a proof of \zcref{thm:Main2Folio}.

\subsection{Homogeneous walls}
\label{subsec:HomoMuralis}
In order to describe how minor models may invade the compass of a flat wall, we need to take into account the way such minor models may interact with the graphs contained in the cells of the rendition witnessing the flatness of our wall.
Moreover, since most flat walls come under the condition that a small apex set has already been removed, we also need to take into account how the graphs inside the cells can be enhanced through the additional boundary represented by those apex vertices.
In particular, to really say that the centre of some flat wall is irrelevant, we must ensure that these interactions with minor models are essentially interchangeable everywhere in the wall.

Please note that below we talk about meshes, but since every wall is a mesh and every $2r$-mesh contains an $r$-wall, the two definitions are interchangeable up to a factor of $2$.
\smallskip

Let $H$ be a subgraph of a graph $G$.
An \emph{$H$-bridge} in $G$ is a connected subgraph $J$ of $G$ such that $E(J) \cap E(H) = \emptyset$ and either $E(J)$ consists of a unique edge with both ends in $H,$ or 
$J$ is constructed from a component $C$ of $G - V(H)$ and the non-empty set of edges $F \subseteq E(G)$ with one end in $V(C)$ and the other in $V(H),$ by taking the union of $C,$ the endpoints of the edges in $F,$ and $F$ itself.
Notice that the $H$-bridges induce a partition of $E(G)\setminus E(H)$.
The vertices in $V(J) \cap V(H)$ are called the \emph{attachments} of $B$ and the set $V(J)\setminus V(H)$ is called the \emph{interior} of $J$.

\paragraph{Homogeneous meshes.}
Let $r\geq 3$ be an integer and $M$ be an $r$-mesh in $G$ and let $C_P$ be the perimeter of $M$.
Let $C\subseteq M$ be any cycle of $M$.
The \emph{compass} of $C$ is the union of $C$ together with all $C$-bridges in $G$ that are entirely contained in the compass of $M$.
We denote the compass of $C$ with respect to $\rho$ by $\mathsf{compass}_M(C)$ and we drop the $\rho$ in the subscript if it is understood from the context.
The \emph{interior} of a cycle $C\subseteq M,$ denoted by $\mathsf{int}_M(C)$ is the graph $\mathsf{compass}_{\rho}(C)-C$.
We also drop the $M$ in the subscript if the wall is understood from the context.
\smallskip

Recall the definition of $q$-colorful graphs.
We allow for $q$ to be $0$.
In this case, for every $0$-colorful graph $(G,\chi)$ it holds that $\chi(v) = \emptyset$ for all $v\in V(G)$.
Moreover, in such a case we do not distinguish between the graph $G$ and the $0$-colorful graph $(G,\chi)$.

Given an arbitrary $q$-colorful graph $(G,\chi),$ a set $X\subseteq V(G)$ and a subgraph $H\subseteq G,$ we define
\begin{align*}
 \chi(X) \coloneq \bigcup_{v\in X} \chi(v),
\end{align*}
as well as $\chi(H) \coloneqq \chi (V(H))$.
Note that this should not be confused with the chromatic number of $H,$ which is of no concern to us.
Finally, we write $(H,\chi)$ for the \emph{$q$-colorful subgraph} of $G,$ where we implicitly restrict $\chi$ to the vertex set of $H$.
\medskip

Let $q$ be a non-negative integer and $(G,\chi)$ be a $q$-colorful graph.
Let $r\geq 3$ be an integer and $M$ be an $r$-mesh in $G$ such that $M$ is flat in $G$ witnessed by a sphere-rendition $\rho$.
We say that $M$ is \emph{homogeneous} in $(G,\chi)$ if for every cycle $C\subseteq M$ it holds that
\begin{align*}
 \chi( \mathsf{int} (C) ) = \chi( \mathsf{compass} (M) ).
\end{align*}

We stress that the original notions of homogeneous flat walls were introduced by Robertson and Seymour \cite{RobertsonS1995Graph,RobertsonS2012Graph} however, as pointed out by Sau, Stamoulis, and Thilikos, the price of homogeneity \cite{SauST2020An,SauST2022kApices2,SauST2023kApices} that the exponential price for obtaining a homogeneous wall is a major obstacle along the path towards an efficient graph minor theory.
It was shown by Gorsky, Seweryn, and Wiederrecht \cite{GorskySW2026Price} that in $q$-colorful graphs it is possible to find a large homogeneous flat wall within another huge flat wall by paying a price polynomial in $q$ and the order of the target wall.

\begin{proposition}[Gorsky, Seweryn, and Wiederrecht \cite{GorskySW2026Price}]\label{prop:HomoMuralis}
There exists a function $f_{\ref{prop:HomoMuralis}}\colon \mathbb{N}^2\to\mathbb{N}$ with $f(q,k) \in \mathcal{O}(q^4 \cdot k^6)$ such that for all non-negative integers $q$ and $k$ and every $q$-colorful graph $G$ with a flat $f(q,k)$-wall $W_0$ there exists a flat $k$-wall $W_1\subseteq W_0$ such that the tangle of $W_1$ is a truncation of the tangle of $W_0$ and $W_1$ is homogeneous.

Moreover, there exists an algorithm that takes as input a colorful graph $(G,\chi)$ and a flat wall $W_0$ and computes the flat wall $W_1$ as above in time $\mathbf{poly}(q+k)|\!|G|\!|$.
\end{proposition}

The main problem now is to define the colours in a way that allows us to model our desired properties.

\paragraph{Augmenting cells.}
Let $a,r \in \mathbb{N}$ with $r \geq 3$.
Let $G$ be a graph, $A \subseteq V(G)$ be a set of at most $a$ vertices, and let $W$ be an $r$-wall in $G$ that is flat in $G-A$ witnessed by the sphere rendition $\rho$.
We denote by $C(W)$ the collection of all cells contained in the disc defined by the trace of the perimeter of $W$ which contains all non-perimeter vertices of $W$ which are of degree $3$ in $W$.

For each cell $c \in C(W)$ we consider a labelling $\lambda_c \colon N(c) \to [3]$ such that the set of labels assigned by $\lambda_c$ to $N(c)$ is one of $[1],$ $[2],$ or $[3]$.
In addition, we consider a bijection $\alpha_A \colon A \to [a]$.

The labellings $\lambda_c$ and $\alpha_A$ will be used to fix labels for the boundaries of the annotated graphs $(G[V(\sigma(c)) \cup A],N(c) \cup A)$ to translate them into rooted graphs.
For this to fully serve out purposes, we require some additional definitions.

\begin{figure}[ht]
 \centering
 \begin{tikzpicture}

 \pgfdeclarelayer{background}
		\pgfdeclarelayer{foreground}
			
		\pgfsetlayers{background,main,foreground}
			
 \begin{pgfonlayer}{main}
 \node (C) [v:ghost] {};

 \node(L) [v:ghost,position=180:5cm from C] {
 \begin{tikzpicture}

 \pgfdeclarelayer{background}
		 \pgfdeclarelayer{foreground}
			
		 \pgfsetlayers{background,main,foreground}

 \begin{pgfonlayer}{background}
 \pgftext{\includegraphics[width=3.5cm]{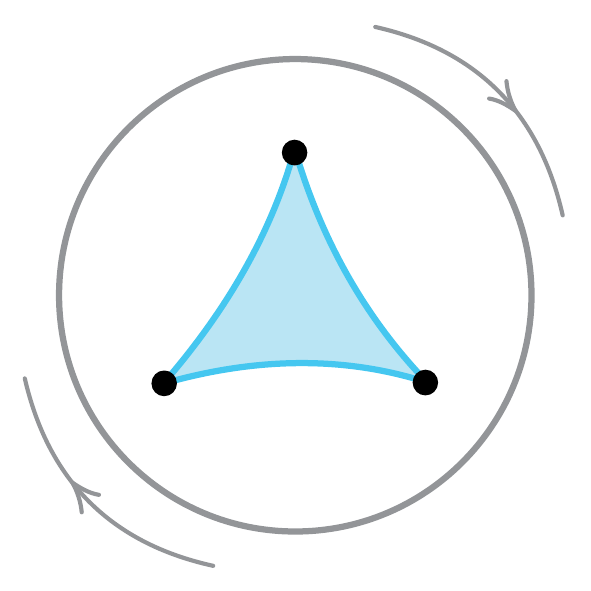}} at (C.center);
 \end{pgfonlayer}{background}
			
 \begin{pgfonlayer}{main}
 \node (C) [v:ghost] {};
 
 \end{pgfonlayer}{main}
 
 \begin{pgfonlayer}{foreground}
 \end{pgfonlayer}{foreground}

 \end{tikzpicture}
 };

 \node(M) [v:ghost,position=0:0cm from C] {
 \begin{tikzpicture}

 \pgfdeclarelayer{background}
		 \pgfdeclarelayer{foreground}
			
		 \pgfsetlayers{background,main,foreground}

 \begin{pgfonlayer}{background}
 \pgftext{\includegraphics[width=3.5cm]{Figures_optlink_arXiv/SpinningTriangle.pdf}} at (C.center);
 \end{pgfonlayer}{background}
			
 \begin{pgfonlayer}{main}
 \node (C) [v:ghost] {};
 
 \end{pgfonlayer}{main}
 
 \begin{pgfonlayer}{foreground}
 \end{pgfonlayer}{foreground}

 \end{tikzpicture}
 };

 \node (Llabel) [v:ghost,position=270:2.3cm from L] {(i)};

 \node (Mlabel) [v:ghost,position=270:2.3cm from M] {(ii)};

 \node (1) [v:ghost,position=90:23mm from Llabel] {};

 \node (L1) [v:ghost,position=90:11mm from 1] {$x_1$};
 \node (L2) [v:ghost,position=200:11mm from 1] {$x_2$};
 \node (L3) [v:ghost,position=340:11mm from 1] {$x_3$};

 \node (2) [v:ghost,position=90:23mm from Mlabel] {};

 \node (R1) [v:ghost,position=90:11mm from 2] {$x_1$};
 \node (R2) [v:ghost,position=200:11mm from 2] {$x_3$};
 \node (R3) [v:ghost,position=340:11mm from 2] {$x_2$};
 
 \end{pgfonlayer}{main}
 
 \begin{pgfonlayer}{foreground}
 \end{pgfonlayer}{foreground}

 \end{tikzpicture}
 \caption{Two labellings of the nodes of some $3$-cell $c$ that cannot be translated into each other by rotation..}
 \label{fig:SpinningTriangle}
\end{figure}

Let $c \in C(W),$ we define the ordering $\Omega_c \coloneqq \langle x_1,\dots,x_{\ell}\rangle$ with $\ell \leq 3$ of the vertices of $N(c)$ such that the following holds:
\begin{enumerate}
 \item $\langle x_1,\dots,x_q\rangle$ is a counter clockwise cyclic ordering of the vertices of $N(c)$ as they appear on the boundary of $c$ in $\rho$ and
 \item for each $i\in[q]$ we have that $\lambda_c(x_i) = i$.
\end{enumerate}
Notice that the cyclic ordering above is significant only in the case where $|N(c)| = 3$ in the sense that $\langle x_1,x_2,x_3\rangle$ remains invariant under shifting, i.\@e.\@, $\langle x_1,x_2,x_\rangle)$ is the same as $\langle x_2,x_3,x_1\rangle,$ but not under inversion, i.\@e.\@, $\langle x_1,x_2,x_3\rangle$ is not the same as $\langle x_3,x_2,x_1\rangle$.
See \zcref{fig:SpinningTriangle} for an illustration.
The second condition, on the other hand, is necessary for us to align the topological information of $\rho$ with the abstract information of the labelling $\lambda_c$ which will give rise to our rooted graphs.

For each cell $c \in C(W)$ we now fix $\pi_c \colon A \cup N(c) \to [a + |N(c)|]$ such that $\pi_c(x) = \alpha_A(x)$ for all $x \in A$ and $\pi_c(x) - a = \lambda_c(x)$ for all $x \in N(c)$.
Finally, this allows us to define the \emph{augmented cell} $c$ as the rooted graph
\begin{align*}
 \mathbf{R}_c \coloneqq (G[V(\sigma(c) \cup A)],r_1,\dots,r_{a + |N(c)|})
\end{align*}
where $r_i = \pi_c^{-1}(i)$ for all $i \in [ a + |N(c)| ]$.
We denote the graph $G[V(\sigma(c) \cup A)]$ by $R_c$.

\paragraph{Palettes of augmented cells.}
With each $c \in C(W)$ we now associate a set of graphs.

Let $d \geq 0$ be an integer.
We denote by $d\text{-}\mathsf{palette}(c)$ the set $d\text{-}\mathsf{folio}(\mathbf{R}_c)$ to which we refer as the \emph{$d$-palette} of $c$.

Notice that $| d\text{-}\mathsf{palette}(c) | \in 2^{\mathbf{O}((a + d)^2)}$.
In the following, we will treat the set of all rooted graphs of detail at most $d$ with at most $a + 3$\, roots, i.\@e.\@, the largest possible choice for a set $d\text{-}\mathsf{palette}(c),$ denoted by $\mathsf{Folio}(d,a+3),$ as the set of colours for our colorful graph.
Hence, in the following we will have that $q \in 2^{(a + d)^2}$ for our colorful graph.

\paragraph{A $d$-folio homogeneous wall.}
Let $a,r \in \mathbb{N}$ with $r \geq 3$.
Let $G$ be a graph $A \subseteq V(G)$ a set of at most $a$ vertices and let $W$ be an $r$-wall in $G$ that is flat in $G-A$ witnessed by the sphere rendition $\rho$.
We define the \emph{$(d,A,W,\rho)$-folio augmentation} of $G,$ denoted by $(G^{(d,A,W,\rho)},\phi),$ as the colorful graph obtained from $G$ by introducing, for each $c \in C(W),$ a vertex $x_c$ adjacent to exactly the vertices of $N(c),$ setting $\phi(v) \coloneqq \emptyset$ for all $v\in V(G)$ and $\phi(c) \coloneqq d\text{-}\mathsf{palette}(c)$ for every $c \in C(W)$.
Notice that $\rho$ may be extended to a sphere rendition $\rho'$ of $G^{(d,A,W,\rho)} - A$ such that $\rho$ is the restriction of $\rho'$ to $G - A \subseteq G^{(d,A,W,\rho)} - A$ simply by setting $\sigma_{\rho'} \coloneqq \sigma(c) + x_c$ for all $c\in C(W)$ and $\sigma_{\rho'}(c) \coloneqq \sigma(c)$ for all $c \in C(\rho) \setminus C(W)$.
We call $\rho'$ the \emph{$\rho$-augmentation} for $G^{(d,A,W,\rho)}$.

Notice that $W$ is still flat in $G^{(d,A,W,\rho)} - A$ as witnessed by $\rho'$.

We say that $W$ is \emph{$(d,A)$-folio homogeneous} in $G$ if $W$ is homogeneous in $(G^{(d,A,W,\rho)} - A,\phi)$.

\paragraph{A flat wall of small treewidth.}
Let $G$ be a graph and $W \subseteq G$ be a flat wall in $G$ witnessed by the sphere rendition $\rho$.
The \emph{treewidth of $W$ in $G$} is the smallest integer $w$ such that the treewidth of the compass of $W$ in $\rho$ is at most $w$.

\paragraph{Suspended walls.}
Let $r\geq 4$.
Given a flat wall $W$ in a graph $G$ witnessed by a sphere rendition $\rho,$ we say that a cell $c\in C(W)$ is an \emph{inner cell} of $W$ if the second and third layer of $W$ separate $N(c)$ from the perimeter of $W$ in $G$.
\smallskip

Let $r \geq 4$ be an integer.
Let $G$ be a graph and $W \subseteq G$ be an $r$-wall in $G$.
We say that $W$ is a \emph{suspended $r$-wall}~--~or simply a \emph{suspended wall}~--~in $G$ if there exists a sphere-rendition $\rho$ of $G-A$ such that $W$ is flat in $\rho$ and for every cell $c\in C(W)$ the following tow properties hold:
\begin{enumerate}
 \item for every non-empty set $S \subseteq N(c)$ there exists a connected subgraph of $\sigma(c)$ containing all vertices from $S$ and avoiding all vertices from $N(c)\setminus S$ and
 \item there exist $|N(c)|$ vertex-disjoint paths from $N(c)$ to the perimeter of $W$.
\end{enumerate}
We say that $\rho$ \emph{witnesses} that $W$ is suspended in $G$.
\medskip

We are now ready to state the main result of this subsection.

\begin{theorem}\label{thm:FolioHomoWall}
There exists a function $f_{\ref{thm:FolioHomoWall}} \colon \mathbb{N}^{3} \to \mathbb{N}$ such that for all integers $r\geq 3$ and $d \geq 0$ the following holds.
Let $G$ be a graph, $A \subseteq V(G)$ with $a \coloneqq |A|,$ and $W_0$ be a flat $f_{\ref{thm:FolioHomoWall}}(a,d,r)$-wall of in $G-A$ with a sphere rendition $\rho$ witnessing the flatness of $W$ in $G-A$.
Then there exists an $r$-wall $W_1 \subseteq W_0$ such that
\begin{enumerate}
 \item $W_1$ is flat in $G-A$ witnessed by $\rho,$
 \item the tangle of $W_1$ in $G-A$ is a truncation of the tangle of $W_0$ in $G-A,$ and
 \item $W_1$ is $(d,A)$-folio homogeneous in $G$.
\end{enumerate}
In particular, if $W_0$ is suspended in $G-A$ witnessed by $\rho,$ then so is $W_1$.
Moreover, $f_{\ref{thm:FolioHomoWall}}(a,d,r) \in 2^{\mathbf{poly}(a + d)} \cdot \mathbf{poly}(r)$ and there exists an algorithm that takes as input $G,$ $A,$ $W_0,$ and $\rho$ as above and computes $W_1$ in time $2^{\mathbf{poly}(w + a + d)} \cdot \mathbf{poly}(r) \cdot (|G| + |\!|G|\!|)$ where $w$ denotes the treewidth of $W_0$ in $G$.
\end{theorem}

\begin{proof}
Let us set $q = q(d,a) \coloneqq |\mathsf{Folio}(d,a+3)|$.
We then set $f_{\ref{thm:FolioHomoWall}}(a,d,r) \coloneqq f_{\ref{prop:HomoMuralis}}(q,r) \in 2^{\mathbf{poly( d + a )}} \cdot \mathbf{poly}(r)$.

Notice that $|C(W_0)| \in \mathbf{O}(|G|)$.
Moreover, $\sum_{c \in C(W_0)} |\sigma(c)| \in \mathbf{O}(|G|)$.
It follows, that we may use the algorithm from \zcref{prop:FolioTW} to compute the $(d,A,W,\rho)$-folio augmentation $(G^{(d,A,W,\rho)},\phi)$ of $G$ in time $2^{\mathbf{poly}(w + a + d)} \cdot |G|$.
We then apply the algorithm from \zcref{prop:HomoMuralis} to the wall $W_0$ in $(G^{(d,A,W,\rho)} - A,\phi)$ to obtain the desired $r$-wall $W_1$ in $(G^{(d,A,W,\rho)} - A,\phi)$.
It is now straightforward to verify that $W_1$ is indeed the desired $(d,A)$-folio homogeneous wall in $G$.
\end{proof}

\subsection{Avoiding the centre of a homogeneous wall}
\label{subsec:CleanCenter}
We have now almost reached the core proof if this section.
Our goal is to show that given a small apex set $A$ in a graph $G$ together with a large enough suspended wall $W$ in $G-A$ that is $(d,A)$-folio homogeneous in $G,$ any tm-pair with small intrusion to the interior of $W$ can be replaced with an equivalent tm-pair avoiding the compass of the centre of $W$.

Before we can dive into this part of the proof, we need a way to bridge the setting of \zcref{thm:RedrawAndDisperseTM} and the setting of flat walls.
This is because a flat wall comes together with a rendition, while \zcref{thm:RedrawAndDisperseTM} requires a partially disc-embedded graph.
To deal with this we adopt the solution of Baste, Sau, and Thilikos \cite{BasteST2023Hitting} -- see \zcref{fig:RedrawModel} for an illustration how tm-pairs may be redrawn to avoid the central part of a homogeneous wall.

\paragraph{Levellings.}
Let $r \geq 4$ be an integer, $G$ be a graph, and $W \subseteq G$ be a flat $r$-wall in $G$ as witnessed by the sphere rendition $\rho$.
Let further $\Delta$ denote the disc bounded by the perimeter of $W$ and $W'$ be the $(r-1)$-central subwall of $W$ such that $W'$ is grounded in the restriction $\rho'$ of $\rho$ to $\Delta$.
Let further $G' \subseteq G$ be the crop of $G$ by $\Delta$.

\begin{figure}[ht]
 \centering
 \begin{tikzpicture}

 \pgfdeclarelayer{background}
		\pgfdeclarelayer{foreground}
			
		\pgfsetlayers{background,main,foreground}
			
 \begin{pgfonlayer}{main}
 \node (C) [v:ghost] {};

 \node(L) [v:ghost] at (-3.5,0) {
 \begin{tikzpicture}

 \pgfdeclarelayer{background}
		 \pgfdeclarelayer{foreground}
			
		 \pgfsetlayers{background,main,foreground}

 \begin{pgfonlayer}{background}
 \pgftext{\includegraphics[width=6cm]{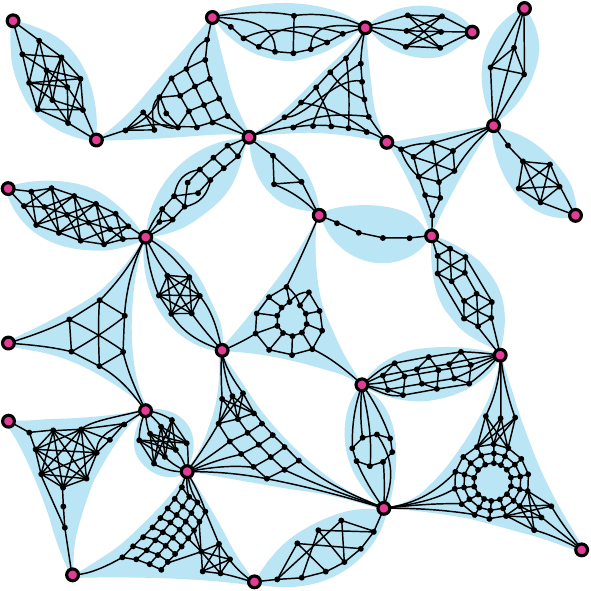}} at (C.center);
 \end{pgfonlayer}{background}
			
 \begin{pgfonlayer}{main}
 \node (C) [v:ghost] {};
 
 \end{pgfonlayer}{main}
 
 \begin{pgfonlayer}{foreground}
 \end{pgfonlayer}{foreground}

 \end{tikzpicture}
 };

 \node(M) [v:ghost] at (3.5,0) {
 \begin{tikzpicture}

 \pgfdeclarelayer{background}
		 \pgfdeclarelayer{foreground}
			
		 \pgfsetlayers{background,main,foreground}

 \begin{pgfonlayer}{background}
 \pgftext{\includegraphics[width=6cm]{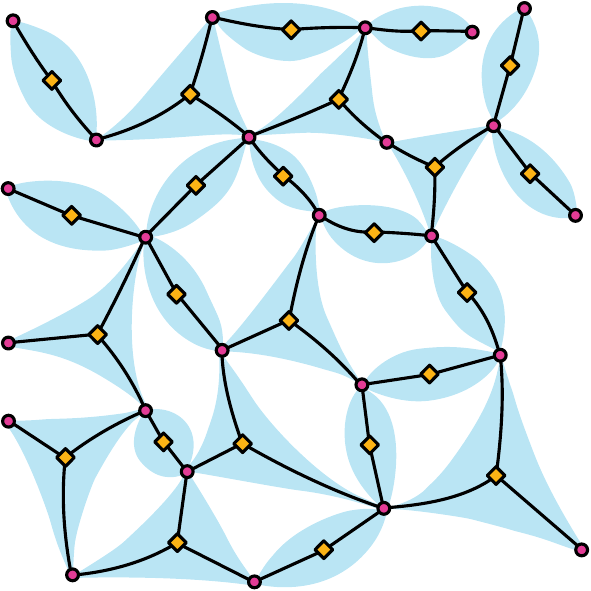}} at (C.center);
 \end{pgfonlayer}{background}
			
 \begin{pgfonlayer}{main}
 \node (C) [v:ghost] {};
 
 \end{pgfonlayer}{main}
 
 \begin{pgfonlayer}{foreground}
 \end{pgfonlayer}{foreground}

 \end{tikzpicture}
 };

 \node (i) [v:ghost] at (-3.5,-3.5) {\textsl{(i)}};
 \node (ii) [v:ghost] at (3.5,-3.5) {\textsl{(ii)}};

 \end{pgfonlayer}{main}
 
 \begin{pgfonlayer}{foreground}
 \end{pgfonlayer}{foreground}

 \end{tikzpicture}
 \caption{A snapshot of a $\Sigma$-rendition $\rho$ of some graph where (i) depicts some ground vertices -- in \textcolor{HotMagenta}{magenta} -- together with some cells $c$ -- the \textcolor{CornflowerBlue}{blue} discs -- and the graphs $\sigma(c)$ associated with those cells.
 In (ii), the result of the levelling of $\rho$ in the same snapshot is depicted. Now each graph $\sigma(c)$ has been replaced with a single cell-vertex -- depicted in \textcolor{ChromeYellow}{yellow}.}
 \label{fig:Levelling}
\end{figure}
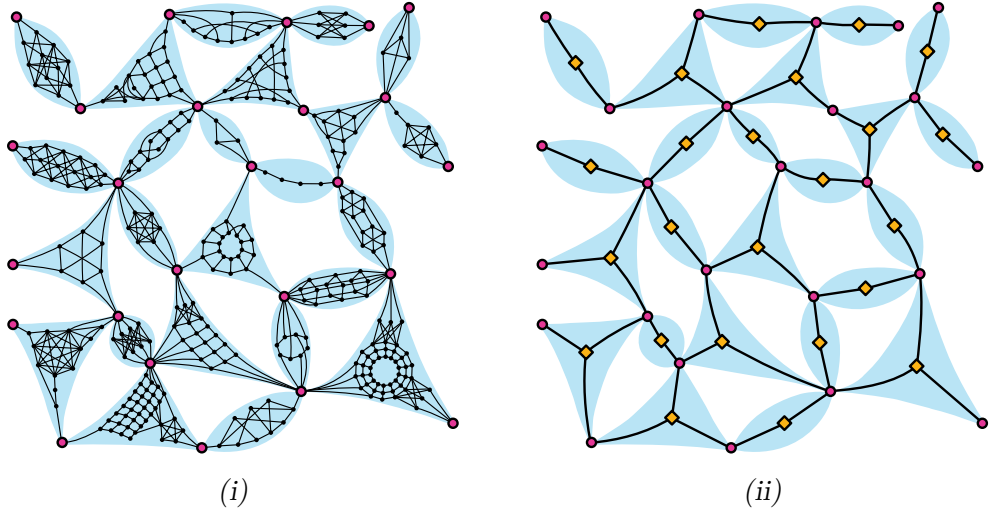

We define the \emph{$(W,\rho)$-levelling} of $G,$ denoted by $\mathsf{Level}(G,W,\rho)$ as the bipartite graph 
\begin{align*}
(N(G') \cup \{ v_c ~ \colon~ c \in C(\rho') \}, \{ xv_c ~\colon~ c \in C(\rho') \text{ and } x \in N(c) \}).
\end{align*}
See \zcref{fig:Levelling} for an illustration.
That is, the vertex set of $\mathsf{Level}(G,W,\rho)$ consists of the nodes contained in $\Delta,$ referred to as the \emph{node-vertices}, on one side, one fresh vertex $v_c,$ called a \emph{cell-vertex}, for each cell contained in $\Delta$ on the other side, and edges $xv_c$ between nodes and cell-vertices whenever $x$ is a node on the boundary of $c$.
Notice also that $\mathsf{Level}(G,W,\rho)$ is the incidence graph of the hypergraph with vertex set $N(\rho')$ and with hyperedges $\{ N(c) ~\colon~ c\in C(\rho') \}$.
Moreover, it is easy to see that $G'$ has an embedding $\Gamma'$ in $\Delta$ such that $\mathsf{\Delta}$ intersects $\Gamma'$ only in the nodes of $W$ that belong to the perimeter of $W$.

The \emph{partial $(W,\rho)$-levelling} of $G,$ denoted by $\mathsf{PLevel}(G,W,\rho),$ is the partially disc-embedded graph with vertex set $( V(G) \setminus V(G') ) \cup V(\mathsf{Level}(G,W,\rho))$ obtained as the union of $G'$ and the graph $G''$ which is obtained from $G$ by deleting all edges of $G',$ all vertices of $G' \setminus \mathsf{perimeter}(W),$ where $\mathsf{perimeter}(W)$ denotes the perimeter of $W$.
It follows that $(G',G'',\Gamma')$ is a partial $\Delta$-embedding of $\mathsf{PLevel}(G,W,\rho),$ we refer to $(G',G'',\Gamma')$ as the \emph{natural $\Delta$-embedding of $\mathsf{PLevel}(G,W,\rho)$}.
The \emph{node-vertices} of $\mathsf{PLevel}(G,W,\rho)$ are the node-vertices of $G',$ the \emph{cell-vertices} of $\mathsf{PLevel}(G,W,\rho)$ are the cell-vertices of $G',$ and the \emph{outer-vertices} of $\mathsf{PLevel}(G,W,\rho)$ are the vertices of $G'' - G'$.

A key observation is that the cell-vertices of $\mathsf{PLevel}(G,W,\rho)$ have maximum degree $3$.

\begin{observation}\label{obs:CellDegree}
Let $r \geq 4$ be an integer, $G$ be a graph, and $W \subseteq G$ be a flat $r$-wall in $G$ as witnessed by the sphere rendition $\rho$.
Let further $\Delta$ be the disc bounded by the perimeter of $W$ that contains all nodes of $W$ and $G^{\star} \coloneqq \mathsf{PLevel}(G,W,\rho)$ be the partial $(W,\rho)$-levelling of $G$.
Then, for every cell-vertex $v_c$ of $G^{\star}$ it holds that $\mathsf{deg}_{G^{\star}}(v) \leq 3$. 
\end{observation}

\paragraph{Levelling a tm-pair.}
Let $r \geq 4$ be an integer, $G$ be a graph, and $W \subseteq G$ be a flat $r$-wall in $G$ as witnessed by the sphere rendition $\rho$.
Let further $\Delta$ be the disc bounded by the perimeter of $W$ that contains all nodes of $W$ and $G^{\star} \coloneqq \mathsf{PLevel}(G,W,\rho)$ be the partial $(W,\rho)$-levelling of $G$.

Now let $\mathcal{J} = (J,T)$ be a tm-pair in $G$.
We denote by $\mathsf{Level}(G,W,\rho,\mathcal{J})$ the \emph{$(W,\rho)$-levelling of $\mathcal{J}$} defined as follows:
Let $(G_1,G_2,\Gamma)$ be the natural $\Delta$-embedding of $G^{\star}$ and let $J_1 \coloneqq J \cap G_1$ as well as $J_2'$ be the intersection of $J$ with the crop of $G$ by $\Delta$.
The vertex set of $J_2'$ consists of two types of vertices, those in $N(\rho)$ and those that belong to $\sigma(c) - N(c)$ for some cell $c \subseteq \Delta$.
Similar to the $(W,\rho)$-levelling of $G,$ we may now define $J_2$ to be the bipartite graph with vertex set 
\begin{align*}
(N(\rho) \cap V(J_2')) \cup \{ v_c ~\colon~ c \in C(\rho) \text{, } c \subseteq \Delta \text{, and } V(J_2') \cap V(\sigma(c) - N(c)) \neq \emptyset\}
\end{align*}
and edge set
\begin{align*}
\{ xv_c ~\colon~ x \in N(\rho) \cap V(J_2') \text{, } x \in N(c) \text{, } c \subseteq \Delta \text{, and } V(J_2') \cap V(\sigma(c) - N(c)) \neq \emptyset\}.
\end{align*}
We set $\widehat{J} \coloneqq J_1 \cup J_2$.
Moreover, let 
\begin{align*}
\widehat{T} \coloneqq (T \cap V(J_1)) \cup (T \cap N(\rho)) \cup \{ v_c \in V(J_2) \}.
\end{align*}
That is, we remove all vertices of $T$ which are drawn inside of cells from $\Delta$ and replace them with the corresponding cell-vertices of $J_2$.
The resulting pair $(\widehat{J}, \widehat{T})$ is now a tm-pair of $G^{\star} \coloneqq \mathsf{PLevel}(G,W,\rho)$ and we define $\mathsf{Level}(G,W,\rho,\mathcal{J}) \coloneqq (\widehat{J}, \widehat{T})$.

\begin{figure}[ht]
 \centering
 \begin{tikzpicture}

 \pgfdeclarelayer{background}
		\pgfdeclarelayer{foreground}
			
		\pgfsetlayers{background,main,foreground}

 \begin{pgfonlayer}{background}
 \pgftext{\includegraphics[width=9cm]{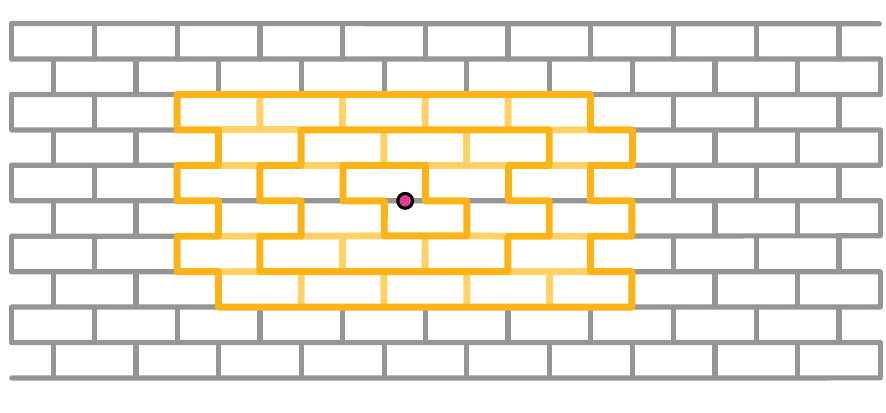}} at (C.center);
 \end{pgfonlayer}{background}
			
 \begin{pgfonlayer}{main}
 \node (C) [v:ghost] {};

 \node(v) [v:ghost] at (-0.55,0.25) {$x$};
 
 \end{pgfonlayer}{main}
 
 \begin{pgfonlayer}{foreground}
 \end{pgfonlayer}{foreground}

 \end{tikzpicture}
 \caption{The $3$-central subwall around the vertex $x$. The three layers of the subwall are depicted in a slightly darker shade of \textcolor{ChromeYellow}{yellow} than the rest of the subwall.}
 \label{fig:CentredWalls}
\end{figure}

Let $W$ be a wall, $s \in \mathbb{N}$ be an integer, and $x \in V(W)$ be a vertex at face-distance at least $s+1$ from the perimeter of $W$ in $W$.
The \emph{$s$-central subwall around $x$} of $W$ is the smallest subwall $W'$ of $W$ with $s$ layers such that the innermost layer of $W'$ separates $s$ from the perimeter of $W$.
See \zcref{fig:CentredWalls} for an illustration.

\begin{theorem}\label{thm:irrelevantVertex}
There exists a function $f_{\ref{thm:irrelevantVertex}} \colon \mathbb{N}^4 \to \mathbb{N}$ such that for all integers $d,k,b,r \in \mathbb{N},$ if $(G,R)$ is an annotated graph of $\mathsf{bidim}(G,R) \leq b$ with $A \subseteq V(G)$ containing an $f_{\ref{thm:irrelevantVertex}}(b,d,k,r)$-wall such that $W$ is suspended and $(4d,A)$-folio homogeneous in $G - A$ as witnessed by the sphere rendition $\rho$ of $G - A$ such that the compass of $W$ does not contain a vertex of $R,$ then the following holds:

If $\Delta_r$ is the disc bounded by the perimeter of the $r$-central subwall $W'$ of $W$ which contains all nodes of $W'$ in $\rho$ and $$X \coloneqq \bigcup_{\substack{c \in C(\rho) \\ c \subseteq \Delta_r}} V(\sigma(c)),$$ then $v$ is strongly irrelevant for $(k,d)\text{-}\mathsf{folio}(G,R)$.

Moreover, $f_{\ref{thm:irrelevantVertex}}(b,d,k,r) \in 2^{\mathbf{poly}(b + d)} \cdot \mathbf{poly}(k + r)$.
\end{theorem}

Before we proceed with the proof of \zcref{thm:irrelevantVertex} we require some technical lemmas on locally redrawing paths inside of suspended walls.
First we need to ensure that any three paths entering a large enough wall and passing through it can always be connected into specific branch vertices of the wall while maintaining their endpoints in the perimeter and avoiding the central part of the wall.

\begin{lemma}\label{lemma:TamePathsComingIntoWall}
Let $s \geq 14$ and $t \geq 13$ be integers, $\Delta$ be a closed disc, and $W$ be a $\Delta$-embedded $s$-wall such that the perimeter of $W$ separates $\mathsf{bd}(\Delta)$ from all non-perimeter vertices of $W$.
Let now $q\in[3]$ be an integer and $v \in V(W)$ be a vertex contained in the $(t - 4)$-central subwall $W'$ of $W$.
If there exist internally vertex disjoint paths $P_1,\dots,P_q$ in $W$ such that
\begin{itemize}
 \item $v$ is an endpoint of $P_i$ for all $i\in [q],$
 \item $P_i$ has its other endpoint $x_i$ on the perimeter of $W$ for all $i \in [q],$ and
 \item the $P_i$ are indexed such that their edges incident with $v$ appear in $\Delta$ in counter-clockwise order.
\end{itemize}
Then there exist pairwise vertex disjoint paths $Q_1,\dots,Q_q$ in $W$ such that
\begin{enumerate}
 \item there are vertices $y_i,$ $i \in [q],$ on the perimeter of $W'$ appearing in counter-clockwise order, each with a neighbour in the interior of the disc $\Delta'$ bounded by the perimeter of $W',$
 \item $Q_i$ has endpoints $x_i$ and $y_i$ for each $i \in [q],$ and
 \item $Q_i$ is internally disjoint from $\Delta'$ for each $i \in [q],$
\end{enumerate}
\end{lemma}

\begin{proof}
This is a simple application of \zcref{thm:menger}.
Let $C_1,\dots,C_4$ be the first four layers of $W',$ let $Y$ denote the set of all vertices of $C_4$ which have degree $3$ in $W$ and a neighbour in the interior of the disc bounded by $C_4$.
Finally, let $\mathcal{F}$ be the collection of all minimal $V(C_1)$-$V(C_4)$-paths in $W$.
Notice that $|\mathcal{F}| \geq 4$ since $t - 8 \geq 5$.

Now let $G'$ be the graph consisting of the union of the $P_i,$ $i\in[q],$ together with the union of the $C_i,$ $i\in[4],$ and the union of all paths in $\mathcal{F}$.
Let $\mathcal{Q}$ be an $X \coloneqq \{ x_i ~\colon~ i\in[q] \}$-$Y$-linkage in $G'$ of maximum order.
We claim that $|\mathcal{Q}| = q$ which would imply the assertion.
So suppose that $|\mathcal{Q}| < q$.
Then \zcref{thm:menger} gives a separator $S$ with $|S| < q \leq 3$ between $X$ and $Y$ in $G'$.
However, by our choices we now know that $S$ must miss one of the $P_i,$ one of the cycles $C_j$ and at least $3$ of the paths from $\mathcal{F},$ say $F_1,\dots,F_3$.
Since $|S| \leq 2$ we know that on $C_4$ there must be some vertex $y$ of $Y$ and $h \in[3]$ such that there exists a subpath of $C_4$ between the endpoint of $F_h$ on $C_4$ and $y$ in $G - S$.
As this is enough to certify that $G' - S$ still contains an $X$-$Y$-path we have reached a contradiction.
\end{proof}

Next, we wish to show that any cell deep in a suspended wall can be linked to a set of prescribed branch vertices on the perimeter of the same wall.
Together with \zcref{lemma:TamePathsComingIntoWall}, this acts as the foundation allowing us to replace cell-vertices with actual cells when translating between tm-pairs in the levelling of some suspended wall and tm-pairs in the actual graph.
\smallskip

Let $r$ and $s$ be integers with $2r \geq s + 3$ and let $W$ be an $r$-wall.
We say a brick $B$ of $W$ is \emph{$s$-central} if is a brick of the $s$-central wall of $W$.

\begin{lemma}\label{lemma:LinkingSuspendedCell}
Let $r$ and $s$ be integers with $2r \geq s + 25$.
Let $G$ be a graph with an $r$-wall $W$ suspended in $G$ as witnessed by the sphere rendition $\rho$.

If $c \in C(\rho)$ is a cell contained in the compass of an $s$-central brick of $W$ together with the labelling $\lambda_c$ of $N(c)$ and $x_1,\dots,x_{|N(c)|}$ are distinct vertices of degree $3$ on the perimeter of $W$ numbered in counter-clockwise order according to their appearance on the perimeter of $W,$ then there exist $|N(c)|$ pairwise vertex-disjoint paths $P_1,\dots,P_{|N(c)|}$ such that for each $i \in [|N(c)|],$ $P_i$ has endpoints $x_i$ and $\lambda_c(i)$.
\end{lemma}

\begin{proof}
Let $A \subseteq W$ be the annulus defined on the first $8$ layers of $W$ together with all minimal subpaths of $W$ with one end on the perimeter of $W$ and the other on the $8$th layer of $W$.
Let further $W'$ denote the $7$-central subwall of $W$.

Then $A$ and $W'$ coincide precisely in the perimeter of $W'$.

\begin{beautifulclaim}\label{claim:linkCellInside}
There exists a linkage $\mathcal{L}_1$ of order $|N(c)|$ between $N(c)$ and the set $Y$ of vertices on the perimeter of $W'$ that have neighbours in $W$ which do not belong to $W'$.
Moreover, $\mathcal{L}_1$ is entirely contained in the compass of the $7$th layer of $W$.
\end{beautifulclaim}

\begin{claimproof}
We have assumed that $W$ is suspended in $G$ as witnessed by $\rho$.
Hence, there exist $|N(c)|$ pairwise vertex disjoint paths between $N(c)$ and the perimeter of $W$.
It follows that there also exists a linkage $\mathcal{P}$ of order $|N(c)|$ between $N(c)$ and the perimeter of $W'$ which is entirely contained in the compass of the $7$th layer of $W$.

Let $C$ now denote the $4$th layer of $W'$.
Then there exists a linkage $\mathcal{Q}$ of order $5$ between $Y$ and $C$ in $W'$ which is internally disjoint from the compass of $C$.
Now a similar argument as the one used in the proof of \zcref{lemma:TamePathsComingIntoWall} together with \zcref{thm:menger} yields the existence of the desired linkage $\mathcal{L}_1$.
\end{claimproof}

Now consider the layers $2$ to $7$ of $W$ together with all minimal subpaths of $W$ between these two layers.
This results in an annulus $A'$ with $6$ circles and at least $6$ rails.
Recall that the vertices $x_1,\dots,x_{|N(c)|}$ are numbered in counter-clockwise order along the perimeter of $W$ and observe that the endpoints of the paths in $\mathcal{L}_1$ on the perimeter of $W'$ are forced to follow the labelling $\lambda_c$ which is also based on the counter-clockwise orientation of $\Delta$.
Let now $\circledcirc$ be the domain of $A$ and $\circledcirc' \subseteq \circledcirc$ be the domain of $A'$.
We may now find a drawing of the graph $G''$ that is the union of $A' \cup \bigcup\mathcal{L}_1$ together with three pairwise vertex-disjoint subpaths of $W$ between the vertices $x_i$ and the $2$nd layer of $W$ in $\circledcirc$ such that $A'$ is drawn in $\circledcirc',$ the $x_i$ lie on one boundary component of $\circledcirc$ while each path in $\mathcal{L}_1$ has one endpoint on the other boundary component of $\circledcirc$.
It follows from our discussion above and this construction that the topological linkage with pattern $$\{ \{ x_i,\lambda_c^{-1}(i) \} ~\colon~ i\in[|N(c)|] \}$$ is topologically feasible in $\circledcirc$.
Hence, an application of \zcref{lem:two-sided-cylinder} yields that this linkage is also realisable in $G''$.
Let us denote this linkage by $\mathcal{L}_2$.
We now obtain the desired paths $P_i$ by combining the linkages $\mathcal{L}_1$ and $\mathcal{L}_2$.
\end{proof}

We immediately obtain the following corollary.

\begin{corollary}\label{cor:LinkingSuspendedCell}
Let $r$ and $s$ be integers with $2r \geq s + 25$.
Let $G$ be a graph with an $r$-wall $W$ suspended in $G$ as witnessed by the sphere rendition $\rho$ and let $G^{\star}$ be the $(W,\rho)$-levelling of $G$ where $W^{\star}$ denotes the $(W,\rho)$-levelling of $W$ in $G$.

If $c \in C(\rho)$ is a cell contained in the compass of an $s$-central brick of $W$ and $x_1,\dots,x_{|N(c)|}$ are distinct vertices of degree $3$ on the perimeter of $W^{\star}$ numbered in counter-clockwise order according to their appearance on the perimeter of $W^{\star},$ then there exist $|N(c)|$ pairwise internally vertex-disjoint paths $P_1,\dots,P_{|N(c)|}$ such that for each $i \in [|N(c)|],$ $P_i$ has endpoints $x_i$ and $v_c$.
\end{corollary}

\begin{figure}[ht]
 \centering
 \begin{tikzpicture}

 \pgfdeclarelayer{background}
		\pgfdeclarelayer{foreground}
			
		\pgfsetlayers{background,main,foreground}

 \begin{pgfonlayer}{background}
 \pgftext{\includegraphics[width=9cm]{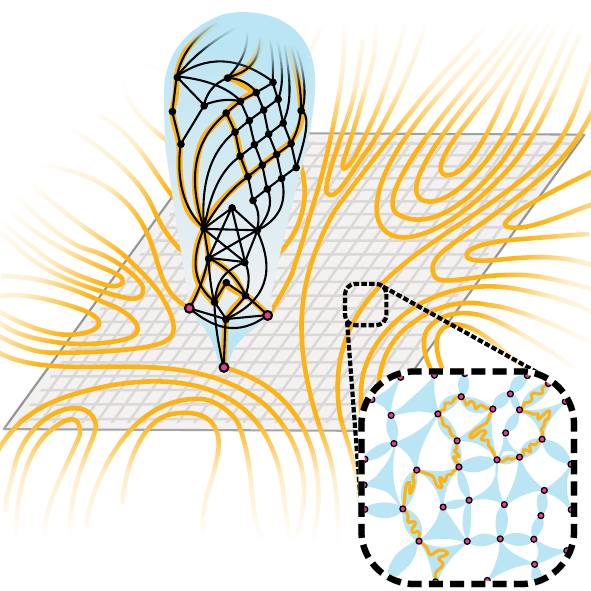}} at (C.center);
 \end{pgfonlayer}{background}
			
 \begin{pgfonlayer}{main}
 \node (C) [v:ghost] {};
 
 \end{pgfonlayer}{main}
 
 \begin{pgfonlayer}{foreground}
 \end{pgfonlayer}{foreground}

 \end{tikzpicture}
 \caption{An illustration of how a tm-pair $(J,T)$ may interact with the compass of a flat wall.
 The \textcolor{ChromeYellow}{yellow} curves indicate the paths of $J$ between the vertices of $T$.
 For a single cell $c,$ the figure illustrates the possible behaviour of $(J,T)$ inside $\sigma(c)$ and for one of the paths $P$ of $J,$ its behaviour inside the cells of the rendition is indicated in the bottom right of the figure.}
 \label{fig:PathsCrossingRendition}
\end{figure}

Consult \zcref{fig:PathsCrossingRendition} to gain some further intuition on how tm-pairs may interact with flat walls in graphs.
Our main goal is now to combine most tools of this section and the previous one into a single argument that allows us to redraw any possible minor model into an ``equivalent'' one that voids the centre of a given flat and homogeneous wall.

\begin{proof}[Proof of \zcref{thm:irrelevantVertex}.]
We begin by setting the stage through defining the function $f_{\ref{thm:irrelevantVertex}}$.
\begin{align*}
 f_{\ref{thm:irrelevantVertex}}(b,d,k,r) \coloneqq f_{\ref{thm:RedrawAndDisperseTM}}(67,r,b,4d,k) + (b+2)^2 + 2,
\end{align*}
then the claimed bounds follow directly from \zcref{thm:RedrawAndDisperseTM}.
\smallskip

Let us select an ordered multiset $\mathcal{R}$ of $k$ roots from $R$ and let $(H,\mathcal{R}_H) \in (k,d)\text{-}\mathsf{folio}(G,\mathcal{R})$ be any graph.
Recall that $\mathsf{detail}(H,\mathcal{R}_H) = \max \{ |V(H) \setminus \mathcal{R}_H|, |E(H)| \}$.
It follows that, since $(H,\mathcal{R}_H)$ is a rooted minor of detail at most $d$ in $(G,R)$ where $\mathcal{R}_H$ is a labelled multisubset of $R$ with $|\mathcal{R}_H| \leq k,$ there exists a tm-pair $(J,T)$ in $(G,R)$ such that $|T \setminus R| \leq 4d$ and there exists a labelled multiset $\mathcal{T}$ of $T\cap R$ with $|\mathcal{T}| \leq k$ such that $(H,\mathcal{R}_H)$ is a rooted minor of $(J,\mathcal{T})$.

Recall that $\rho$ is a sphere rendition of $G - A$ witnessing that $W$ is suspended and flat in $G - A$.
Let $\Delta$ denote the disc bounded by the trace of the perimeter of $W$ which contains all nodes of $W$.
Let further $(G_1,G_2,\Gamma)$ be the natural partial $\Delta$-embedding of $G^{\star} \coloneqq \mathsf{PLevel}(G,W,\rho)$ and let $(\widetilde{J},\widetilde{T}) = \mathsf{Level}(G,W,\rho,(J,T))$ be the $(W,\rho)$-levelling of $(J,T)$.

Let $Q$ denote the set of all cells $c \in C(\rho)$ such that $v_c \in V(\widetilde{J})$.
For each $c \in C$ let $\mathbf{R_c} = (H_c,\mathcal{U}_c)$ be the augmented cell $c$ and recall that we associated a labelling $\lambda(c)$ of $N(c)$ with $c$ which is reflected in the labelling $\mathcal{U}_c$.
Next, consider the graph 
\begin{align*}
J_c \coloneqq (H_c \cap J) \cup \{ u \in (N(c) \cup A)\setminus V(J) ~\colon~ \text{there is } u'\in \mathcal{U}_c \text{ with } u = u'\}
\end{align*}
and denote by $\mathbf{J}_c$ the rooted graph $(J,\mathcal{U}_c)$.
Notice that, since the detail of $(J,T)$ is at most $4d$ we have that $\mathbf{J}_c \in 4d\text{-}\mathsf{folio}(\mathbf{R}_c)$.
Finally, observe that every vertex of $Q$ has degree at most $3$ in $\widetilde{J}$.

\paragraph{Redrawing the levelling.}
We now apply \zcref{thm:RedrawAndDisperseTM} to the tm-pair $(\widetilde{J},\widetilde{T})$ in $G^{\star}$ and the wall $W$ aiming for a safely $67$-dispersed redrawing of $(\widetilde{J},\widetilde{T})$ avoiding the compass of the $r$-central subwall of $W$.

Now note the following:
\begin{itemize}
 \item by our assumption that the compass of $W$ does not contain a vertex from $R$ we have that $R \cap V(G_1) = \emptyset,$
 \item let $W'$ be the subwall of $W$ obtained by removing the first layer of $W$ and let $W^{\star} \subseteq G^{\star}$ be the $(W,\rho)$-levelling\footnote{We omit denoting the branch vertices of $W'$ here as those are precisely all of its vertices that are of degree $3$ in $W$.} of $W',$ then $G_1$ contains an $f_{\ref{thm:irrelevantVertex}}(b,d,k,r) \coloneqq f_{\ref{thm:RedrawAndDisperseTM}}(67,r,b,4d,k)$-wall $W'$ whose perimeter separates its interior from $\mathsf{bd}(\Delta),$
 \item we may denote by $C^{\star}_1,\dots,C^{\star}_{g_{\ref{thm:RedrawAndDisperseTM}}(4d,b,k) + h_{\ref{thm:RedrawAndDisperseTM}}(67,4d,k,b)}$ the first $g_{\ref{thm:RedrawAndDisperseTM}}(4d,b,k) + h_{\ref{thm:RedrawAndDisperseTM}}(67,4d,k,b)$ layers of $W^{\star}$ indexed in the natural order, and
 \item we denote by $\Delta^{\star}_1,\dots,\Delta^{\star}_{g_{\ref{thm:RedrawAndDisperseTM}}(4d,b,k) + h_{\ref{thm:RedrawAndDisperseTM}}(67,4d,k,b)} \subseteq \Delta$ the open discs defined by these layers.
\end{itemize}
Then \zcref{thm:RedrawAndDisperseTM} provides us with a tm-pair $(J^{\star},T^{\star})$ of $G^{\star}$ and an integer $q \in [g_{\ref{thm:RedrawAndDisperseTM}}(4d,b,k)]$ such that
\begin{enumerate}
 \item $J^{\star} - (V(G_1) \cap \Delta^{\star}_q)$ is a subgraph of $\widetilde{J} - (V(G_1) \cap \Delta^{\star}_q),$
 \item if $\circledcirc \subseteq \Delta$ denotes the annulus with boundary components $C^{\star}_q$ and $C^{\star}_{q + h_{\ref{thm:RedrawAndDisperseTM}}(67,4d,k,b)-1}$ and $G'$ denotes the crop of $G_1$ by $\circledcirc,$ then $V(G') \cap (\widetilde{T} \cup T^{\star}) = \emptyset,$
 \item $(J^{\star} \cap (G_1 \cap \mathsf{cl}( \Delta^{\star}_{q + h_{\ref{thm:RedrawAndDisperseTM}}(67,4d,k,b)} )), T^{\star} \cap (V(G_1) \cap \mathsf{cl}( \Delta^{\star}_{q + h_{\ref{thm:RedrawAndDisperseTM}}(67,4d,k,b)} )))$ is a tm-pair of $W^{\star}$ that is safely $67$-dispersed in $W^{\star}$ and none of the vertices of $T^{\star} \cap (V(G_1) \cap \mathsf{cl}( \Delta^{\star}_{q + h_{\ref{thm:RedrawAndDisperseTM}}(67,4d,k,b)} ))$ is within face distance less than $67$ in $W^{\star}$ from some vertex of $C^{\star}_{q + h_{\ref{thm:RedrawAndDisperseTM}}(67,4d,k,b)} \cup C^{\star}_{q + g_{\ref{thm:RedrawAndDisperseTM}}(4d,b,k) + h_{\ref{thm:RedrawAndDisperseTM}}(67,4d,k,b)},$
 \item $J^{\star} \cap (G_1 \cap \Delta^{\star}_{g_{\ref{thm:RedrawAndDisperseTM}}(4d,b,k) + h_{\ref{thm:RedrawAndDisperseTM}}(67,4d,k,b)}) = \emptyset,$
 \item the compass of the $r$-central subwall of $W^{\star}$ is contained in $\Delta_{g_{\ref{thm:RedrawAndDisperseTM}}(4d,b,k) + h_{\ref{thm:RedrawAndDisperseTM}}(67,4d,k,b)},$ and
 \item there is a $Q$-respecting contraction mapping $\phi$ of $\mathsf{Diss}(\widetilde{J},\widetilde{T})$ to $\mathsf{Diss}(J^{\star},T^{\star})$.
\end{enumerate}

\paragraph{Reconstructing $(H,\mathcal{R}_H)$.}
Our goal is to now use the homogeneity of $W$ and the tm-pair $(J^{\star},T^{\star})$ in order to reconstruct a rooted minor model of $(H,\mathcal{R}_H)$ which still avoids all nodes and vertices inside cells contained in $\Delta^{\star}_{g_{\ref{thm:RedrawAndDisperseTM}}(4d,b,k) + h_{\ref{thm:RedrawAndDisperseTM}}(67,4d,k,b)}$.
Once this is done, we have shown that for every member of $(k,d)\text{-}\mathsf{folio}(G,R)$ there exists a rooted minor model avoiding the interior of $\Delta^{\star}_{g_{\ref{thm:RedrawAndDisperseTM}}(4d,b,k) + h_{\ref{thm:RedrawAndDisperseTM}}(67,4d,k,b)}$ which implies in particular that the entire compass of the $r$-central subwall of $W$ is irrelevant.

To reach this point recall the definition of safely $s$-dispersed:
We know that $$(J^{\circ},T^{\circ}) \coloneqq (J^{\star} \cap (G_1 \cap \mathsf{cl}(\Delta^{\star}_{q + h_{\ref{thm:RedrawAndDisperseTM}}(67,4d,k,b)})), T^{\star} \cap (V(G_1) \cap \mathsf{cl}(\Delta^{\star}_{q + h_{\ref{thm:RedrawAndDisperseTM}}(67,4d,k,b)})))$$ is safely $67$-dispersed in $W^{\star}$ which means that
\begin{enumerate}
 \item any two vertices $x,y \in T^{\circ}$ are at face distance of least $135$ in $W^{\star},$ and
 \item for every $x \in T^{\circ}$ with $\mathsf{deg}_{J^{\circ}} = d,$ the graph $J^{\circ}[\mathsf{F}^{67}_{W^{\star}} \cap V(J^{\circ})]$ consists of $d$ paths with $x$ as their unique common endpoint.
\end{enumerate}
It follows that for every $v_c \in Q$ there exists a $134$-subwall $W^{\rightmoon}_c$ of $W^{\star}$ such that $\phi(v_c)$ is at face distance $133$ from the perimeter of $W^{\rightmoon}_c$ and every brick of $W^{\rightmoon}_c$ is also a brick of $W^{\star}$.
Then $W^{\rightmoon}_c$ naturally corresponds to a $134$-subwall $W^{\leftmoon}_c$ of $W$.

Then the paths from $J^{\circ}[\mathsf{F}^{67}_{W^{\star}} \cap V(J^{\circ})]$ form a family of internally vertex-disjoint paths $\mathcal{L}_c$ of order $d_c \coloneqq \mathsf{deg}_{\widetilde{J}}(v_c) = \mathsf{deg}_{J^{\star}}(\phi(v_c)) \leq 3$ in $W^{\rightmoon}_c$ sharing the unique common endpoint $\phi(v_c)$.
Recall that $W$ is $(4d,A)$-folio homogeneous in $G-A$.
Hence, if we denote the innermost brick of $W^{\rightmoon}_c$ by $B^{\star}_c,$ then $B^{\star}_c$ corresponds to a brick $B_c$ of $W$.
Then there exists a cell $c' \in C(\rho)$ contained in the interior of the disc defined by the trace of $B_c$ such that $\mathbf{J}_c \in 4d\text{-}\mathsf{folio}(\mathbf{R}_{c'})$.
Hence, the vertex $v_{c'} \in V(G^{\star})$ exists in the interior of $B^{\star}_c$.
By combining \zcref{lemma:TamePathsComingIntoWall} and \zcref{cor:LinkingSuspendedCell} inside $W^{\rightmoon}_c$ we may replace the paths in $\mathcal{L}_c$ with a family of $d_c$ internally vertex-disjoint paths $\mathcal{L}_{c'}$ such that
\begin{itemize}
 \item each path $L \in \mathcal{L}_{c'}$ has $v_{c'}$ as one endpoint,
 \item the other endpoint of $L \in \mathcal{L}_{c'}$ is an endpoint of some path in $\mathcal{L}_c$ that is not $\phi(v_c),$ and
 \item the paths in $\mathcal{L}_{c'}$ are contained in the disc whose boundary is the perimeter of $W^{\rightmoon}_c$.
\end{itemize}
Now let $\widehat{J}$ be the graph obtained from $J^{\star}$ by replacing the paths $\mathcal{L}_c$ and the vertex $\phi(v_c)$ with the paths $\mathcal{L}_{c'}$ and the vertex $v_{c'}$ for all $v_c \in Q$.
Moreover, if we replace each $\phi(v_c)$ in $T^{\star}$ with $v_{c'}$ for each $v_c \in Q$ resulting in the set $\widehat{T},$ we obtain a tm-pair $(\widehat{J},\widehat{T})$ in $G^{\star}$ such that $\mathsf{Diss}(J^{\star},T^{\star})$ and $\mathsf{Diss}(\widehat{J},\widehat{T})$ are isomorphic.
Indeed, this means that there exists a $Q$-preserving contraction mapping $\psi$ of $(\widetilde{T},\widetilde{J})$ to $(\widehat{J},\widehat{T})$ such that $\psi(v_c) = v_{c'}$ for all $v_c \in Q$.

At this point, we are almost done.
Recall that we had labellings $\lambda_c$ for the nodes of each cell $c \in C(\rho)$ assigned to those nodes in counter clockwise orientation on the boundary of $\mathsf{cl}(c)$.
This means that in our choices for the paths $\mathcal{L}_{c'}$ we may ensure that for every edge $xv_c \in E(\widetilde{J})$ there now is a path $L_{xv_c} \in \mathcal{L}_{c'}$ such that there exists an edge $yv_{c'} \in E(L)$ with $\lambda_c(x) = \lambda_{c'}(y)$.
We stress that \zcref{cor:LinkingSuspendedCell} ensures that this choice is indeed possible.
Let us denote by $\widehat{Q}$ the set
\begin{align*}
 \widehat{Q} \coloneqq \{ \psi(x) ~\colon~ x \in Q \}
\end{align*}
and by $\widehat{Q}^*$ the set
\begin{align*}
 \widehat{Q}^* \coloneqq \big( V(\widehat{J}) \cap \{ v_c \in V(G^{\star}) \text{ and } c\in C(\rho) \} \big) \setminus \widehat{Q}
\end{align*}
Then $\widehat{Q} \cap \widehat{Q}^* = \emptyset$ and $V(\widehat{J}) \setminus (\widehat{Q} \cup \widehat{Q}^*) \subseteq V(G)$.

We now construct a new tm-pair $(J^{\fullmoon},T^{\fullmoon})$ in $G$ as follows.

First, for each $v_c \in \widehat{Q}^*$ let $H_c$ be a minimal connected subgraph containing all vertices from $N(c)$ which are adjacent to $v_c$ in $\widehat{J}$ but none of the other vertices from $N(c)$.
The existence of $H_c$ is guaranteed by the assumption that $W$ is suspended in $G$ as witnessed by $\rho$.
Let $T_c$ be the set of $N(c) \cap V(H_c)$ together with all vertices in $H_c$ of degree at least $3$.
Notice that $\Delta(H_c) \leq 3$ and there is at most one vertex in $H_c$ with degree at least $3$ by minimality.

Second, for each $v_{c'} \in \widehat{Q}$ we know that there exists $v_c \in Q$ with $\psi(v_c) = v_{c'}$ such that $\mathbf{J}_c \in 4d\text{-}\mathsf{folio}(\mathbf{R}_{c'})$.
This means that there exists a tm-pair $(D_{c'},F_{c'})$ in $\mathbf{R}_{c'} = (R_{c'},\mathcal{R}_{c'})$ such that $F_{c'}$ contains all roots $\mathcal{R}_{c'}$ of $\mathbf{R}_{c'}$ and $\mathbf{J}_c$ is a rooted minor of $(D_{c'},\mathcal{R}_{c'})$.

Finally, we let $J^{\fullmoon}$ be the graph obtained from $\widehat{J}$ by deleting, for every $v_c \in \widehat{Q}^*$ the vertex $v_c$ and replacing it by the graph $H_c$.
Similarly, we also delete, for each $v_{c'} \in \widehat{Q},$ the vertex $v_c'$ and replace it by the graph $D_{c'}$.
Moreover, we remove all cell-vertices from $T^{\fullmoon}$ and instead add all members of the sets $F_{c'}$ and $T_c$ has above.
We may now observe three facts:
\begin{enumerate}
 \item $J^{\fullmoon}$ is vertex-disjoint from the compass of the $(g_{\ref{thm:RedrawAndDisperseTM}}(4d,b,k) + h_{\ref{thm:RedrawAndDisperseTM}}(67,4d,k,b) +1)$th layer of $W,$ 
 \item for every vertex $r$ that appears in $\mathcal{T}$ we still have that $r\in V(J^{\fullmoon}),$ and
 \item $(H,\mathcal{R}_H)$ is a rooted minor of $(J^{\fullmoon},\mathcal{T})$.
\end{enumerate}
In particular, this means that, if $X$ is the vertex set of the compass of the $r$-central subwall of $W,$ then $(H,\mathcal{R}_H) \in (k,d)\text{-}\mathsf{folio}(G - X,R)$.
Since our choices of $\mathcal{R}$ and $(H,\mathcal{R}_H)$ were arbitrary based on $R$ this implies that $v$ is strongly irrelevant as desired.
\end{proof}

\subsection{Quickly finding a good flat wall}
\label{subsec:FindingFlatWall}

An important assumption on our flat wall $W$ in \zcref{thm:irrelevantVertex} is that $W$ must be suspended in $G-A$.
The other core assumption is that $W$ is $(k,d)$-folio homogeneous.
In order to guarantee the second condition we may call upon \zcref{thm:FolioHomoWall}.
Indeed, if the input wall for \zcref{thm:FolioHomoWall} was already suspended, then so is the $(k,d)$-folio homogeneous wall that comes out on the other side.
However, in order to apply \zcref{thm:FolioHomoWall} we first need a suspended flat wall whose compass has bounded treewidth.

Finally, there is a third assumption on $W$ in \zcref{thm:irrelevantVertex}, that is, its compass should be disjoint from $R$.
Among all three properties, this is the easiest to guarantee as here, we may use a simple pigeonhole argument based on our bidimensionality bound.
Indeed, here we are able to mimic a result from the work of Gorsky, Protopapas, and Wiederrecht \cite{GorskyPW2026Quickly}.
We slightly change their proof to obtain more clear and better dependencies.

\begin{lemma}\label{lemma:BlankWall}
For all integers $r \geq 2$ and $b \geq 1$ and every annotated graph $(G,R),$ if $W$ is a flat $w$-wall with $w \geq (3b + 1)r + 1$ wall in $G$ as witnessed by the sphere-rendition $\rho,$ then one of the following holds:
\begin{enumerate}
 \item $\mathsf{bidim}(G,R) \geq b,$ or
 \item there exists an $r$-subwall $W'$ of $W$ such that the compass of $W'$ does not contain a vertex from $R$.
\end{enumerate}
Moreover, if $\rho$ witnesses that $W$ is suspended in $G,$ then $\rho$ also witnesses that $W'$ is suspended in $G$.
There also exists an algorithm that, given $G,$ $\rho,$ and $W$ as input finds either a red $((b+1) \times (b+1))$-red-minor in $G$ or the wall $W'$ as above in time $\mathbf{poly}(b+r)|\!|G|\!|$.
\end{lemma}

\begin{proof}
We treat $W$ as a mesh and number its vertical paths as $P_1,\dots,P_w$ and its horizontal paths as $Q_1,\dots,Q_w$ as in the definition of meshes.

Let $\widetilde{W}$ be the wall obtained by taking the horizontal and vertical paths $P_{i \cdot (r+1) + 1}$ and $Q_{i \cdot (r+1) + 1}$ for all $i\in[0,3b+1]$.
Then $\widetilde{W}$ is a $(3b+2)$-subwall of $W$.
Moreover, the compass of every brick $B$ of $\widetilde{W}$ contains an $r$-subwall $W_B$ of $W$ which is disjoint from $B$ itself.

Without loss of generality we may assume that every vertex in a cell $c$ in the compass of $W$ is connected to one of the nodes of $c$.

Now there are two cases.
\medskip

\textbf{Case 1:}
The compass of every brick of $\widetilde{W}$ contains a vertex of $R$ in its interior.
\smallskip

Then, for each of the $(3b+1)^2$ bricks $B$ of $\widetilde{W}$ there exists a vertex of $R$ which can be connected to $B$ by a path that is otherwise disjoint from $\widetilde{W}$.
In this case it follows immediately that $(G,R)$ contains a red $(b \times b)$-grid as a red-minor and thus $\mathsf{bidim}(G,R) \geq b$.
\medskip

\textbf{Case 2:}
There exists a brick $B$ of $\widetilde{W}$ such that the interior of its compass does not contain a vertex from $R$.
\smallskip

In this case it follows that the compass of $W_B$ does not contain a vertex from $R$ and our proof is complete.
\end{proof}

From now on we only need to be concerned with finding a large suspended wall whose compass has small treewidth.
This is because in annotated graphs of bounded bidimensionality, \zcref{lemma:BlankWall} can then easily guarantee the third requirement of \zcref{thm:irrelevantVertex}.

Before we can go there, we need to sit down and have a discussion about the history of the flat wall theorem and its application in algorithm design.
In many papers that utilise the Irrelevant Vertex Technique and the Flat Wall Theorem, the flatness of a wall $W$ under an apex set $A$ is encoded as a so-called ``flatness pair''.
Below we provide a high-level discussion of flatness pairs and some of their properties.
Note that we do not state the definitions involved as we are not working with these definitions and only require very specific implications thereof.
A flatness pair is essentially just a pair $(W,\mathfrak{R})$ where the object $\mathfrak{R}$ holds the information about the graph $G,$ the apex set $A$ and the rendition $\rho$ that witnesses the flatness of $W$ in $G-A$.
Sau, Stamoulis, and Thilikos \cite{SauST2024More}, the authors investigate a set of fundamental properties of the cells inside the compass of a flat wall and on its boundary.
They classify and further refine those cells and show that they are able to always guarantee that their flatness pairs have a property called ``regularity''.
The only core property of regular flatness pairs we are interested in is encapsulated in the following definition:

\emph{Regular flat walls.}
Let $r \geq 4$ be an integer.
Let $G$ be a graph and $W \subseteq G$ be an $r$-wall in $G$.
We say that $W$ is a \emph{regular} flat $r$-wall~--~or simply a \emph{regular flat wall}~--~in $G$ if there exists a sphere-rendition $\rho$ of $G-A$ such that $W$ is flat in $\rho$ and for every cell $c\in C(W)$ and for every non-empty set $S \subseteq N(c)$ there exists a connected subgraph of $\sigma(c)$ containing all vertices from $S$ and avoiding all vertices from $N(c)\setminus S$.
We say that $\rho$ \emph{witnesses} that $W$ is a regular and flat wall of $G-A$.
\medskip

Next, we need the following corollary of Lemma 5.8 of the work of Baste, Sau, and Thilikos on the minor-deletion problem \cite{BasteST2023Hitting}.
It follows directly from regularity and the fact that for any rendition $\rho$ of some graph $G$ and any non-vortex cell $c \in C(\rho),$ if $G$ contains two vertex-disjoint paths $P_1$ and $P_2$ such that none of the endpoints of $P_1$ and $P_2$ belongs to $\sigma(c) - N(c),$ then at most one of the two paths can contain edges or non-boundary vertices of $\sigma(c)$.

\begin{proposition}[Baste, Sau, Thilikos \cite{BasteST2023Hitting}]\label{prop:SuspendingAWall}
Let $r \geq 4$ be an integer.
Moreover, let $G$ be a graph and $W \subseteq G$ be a regular flat $(r+2)$-wall in $G$ as witnessed by the sphere-rendition $\rho$.
Then the central $r$-subwall $W'$ of $G$ is suspended in $G$ and the suspension of $W$ in $G$ is witnessed by $\rho$.
\end{proposition}

In light of \zcref{prop:SuspendingAWall}, in order to find a large suspended wall whose compass has bounded treewidth, all that is required is to find a slightly larger regular flat wall whose compass has bounded treewidth.
This brings is to the core proposition we import from the work of Sau, Stamoulis, and Thilikos \cite{SauST2024More}.

\begin{proposition}[Sau, Stamoulis, Thilikos \cite{SauST2024More}]\label{prop:FindGoodFlatWall}
There exists a functions $f_{\ref{prop:FindGoodFlatWall}},g_{\ref{prop:FindGoodFlatWall}} \colon \mathbb{N} \to \mathbb{N}$ such that for every integer $t \geq 1,$ all odd integers $r \geq 5,$ and every graph $G$ one of the following holds:
\begin{enumerate}
 \item $G$ contains a $K_t$ minor,
 \item $G$ has a tree-decomposition of width at most at most $f_{\ref{prop:FindGoodFlatWall}}(t) \cdot r,$ or
 \item there is a set $A \subseteq V(G),$ where $|A| \leq g_{\ref{prop:FindGoodFlatWall}(t)}$ together with an $r$-wall $W$ that is regular flat in $G-A$ as witnessed by the sphere-rendition $\rho$ and the treewidth of $W$ in $G$ is at most $f_{\ref{prop:FindGoodFlatWall}}(t) \cdot r$.
\end{enumerate}
Moreover, $f_{\ref{prop:FindGoodFlatWall}}(t) \in 2^{\mathbf{O}(t^2 \log t)},$ $g_{\ref{prop:FindGoodFlatWall}}(t) \in \mathbf{O}(t^{24}),$ and there exists an algorithm that finds one of the three outcomes above in time $2^{2^{\mathbf{O}(t^2 \log t)} \cdot r \log r + \mathbf{O}(r^2)} \cdot |G| + 2^{ 2^{\mathbf{O}(t^2 \log t)} \cdot r^3 \log r }$.
\end{proposition}

Notice that the proof of \zcref{prop:FindGoodFlatWall} uses two theorems from the literature.
The first is the flat wall theorem of Kawarabayashi, Thomas and Wollan \cite{KawarabayashiTW2018New} which is the source of the dependency $g_{\ref{prop:FindGoodFlatWall}}(t) \in \mathbf{O}(t^{24})$.
It is possible to replace the use of this theorem with the Flat Wall Theorem as proven by Gorsky, Seweryn, and Wiederrecht \cite{GorskySW2025Polynomial} to obtain the bound $g_{\ref{prop:FindGoodFlatWall}}(t) \in \mathbf{O}(t^{3})$.
Similarly, the proof uses a theorem by Kawarabayashi and Kobayashi \cite{KawarabayashiK2020Linear} stating that every $K_t$-minor-free graph with treewidth in $2^{\Omega(t^2 \log t)} \cdot r$ contains an $r$-wall.
In light of \zcref{prop:ratcatcher}, also this dependency can be tamed to show that $f_{\ref{prop:FindGoodFlatWall}}(t) \in \mathbf{poly}(t)$.
However, these improvements do not impact the final running times and dependencies of our results which is the reason why we do not spend further efforts on implementing these stronger bounds.

\subsection{Proof of \zcref{thm:Main2Folio}}
\label{subsec:MainThmProof}
We are now ready for the proof of \zcref{thm:Main2Folio} which is the core engine of our algorithm for the \textsc{Folio} problem.
This is mostly done by combining our previous findings in this section.
Let us start by preparing the function $f_{\ref{thm:Main2Folio}}$.

For this, we first need to fix the constant for dealing with clique-minors.
\begin{align*}
 t(b,d) \coloneqq \left\lfloor \frac{5}{2} \right\rfloor b^2 + 3d^2 +1
\end{align*}

Now we may proceed with the main dependency of \zcref{thm:Main2Folio}.
\begin{align*}
 f_{\ref{thm:Main2Folio}}(k,b,d) \coloneqq ~ & f_{\ref{prop:FindGoodFlatWall}}( t(b,d) ) \cdot \big( (3b + 4)(f_{\ref{thm:FolioHomoWall}}( g_{\ref{prop:FindGoodFlatWall}}( t(b,d) ), 4d, f_{\ref{thm:irrelevantVertex}}( b, d, k, 3 ) ) +1 ) + 4 \big)
\end{align*}
The claimed bounds now follow directly from the bounds in \zcref{thm:Main2Folio}, \zcref{prop:FindGoodFlatWall}, \zcref{thm:FolioHomoWall}, and \zcref{thm:irrelevantVertex}.
Therefore, from here, all we have to do is explain the algorithm and prove its correctness.

\textbf{Input:} An annotated graph $(G,R)$ with $\mathsf{bidim}(G,R) \leq b$.
\begin{description}
 \item[Step 1:] As long as $|E(G)| \in 2^{\Omega( t(b,d)^2 \log t(b,d) )}|G|$ run the algorithm from \zcref{prop:quicklyfindingaclique} to find a $K_{t(b,d)}$-minor in $G$. Then \textbf{forward} the corresponding minor model to \textbf{Step 7}.

 \item[Step 2:] Run the algorithm from \zcref{prop:FindGoodFlatWall} looking for a regular flat $ \big((3b + 4)(f_{\ref{thm:FolioHomoWall}}( g_{\ref{prop:FindGoodFlatWall}}( t(b,d) ), 4d, f_{\ref{thm:irrelevantVertex}}( b, d, k, 3 ) ) +1 ) + 4\big)$-wall $W_0$. This can have three outcomes.
 \begin{description}
 \item[Case 2.1:] \textsl{The algorithm returns a model of a $K_{ t(b,d) }$-minor}. Then \textbf{forward} the clique-minor model to \textbf{Step 7}.

 \item[Case 2.2:] \textsl{The algorithm returns a tree-decomposition $(T,\beta)$ of width at most $f_{\ref{thm:Main2Folio}}(k,b,d)$.} Then \textbf{return} $(T,\beta)$ and \textbf{terminate}.

 \item[Case 2.3:] \textsl{The algorithm returns an apex set $A \subseteq V(G)$ with $|A| \leq g_{\ref{prop:FindGoodFlatWall}}( t(b,d) )$ together with a regular flat $\big( (3b + 4)(f_{\ref{thm:FolioHomoWall}}( g_{\ref{prop:FindGoodFlatWall}}( t(b,d) ), 4d, f_{\ref{thm:irrelevantVertex}}( b, d, k, 3 ) ) +1 ) + 4 \big)$-wall $W_0$ in $G-A$ and a sphere-rendition $\rho$ of $G-A$ witnessing the regular flatness of $W_0$. Moreover, the algorithm guarantees that the treewidth of $W_0$ in $G$ is at most $f_{\ref{thm:Main2Folio}}(k,b,d)$.}
 Then \textbf{forward} $A,$ $W_0,$ and $\rho$ to \textbf{Step 3}.
 \end{description}

 \item[Step 3:] \textbf{Input:} A set $A \subseteq V(G)$ of size at most $g_{\ref{prop:FindGoodFlatWall}}( t(b,d) ),$ together with a $\big( (3b + 4)(f_{\ref{thm:FolioHomoWall}}( g_{\ref{prop:FindGoodFlatWall}}( t(b,d) ), 4d, f_{\ref{thm:irrelevantVertex}}( b, d, k, 3 ) ) +1 ) + 4 \big)$-wall $W_0$ in $G-A$ and a sphere-rendition $\rho$ of $G-A$ witnessing the regular flatness of $W_0$ such that $W_0$ has treewidth at most $f_{\ref{thm:Main2Folio}}(k,b,d)$ in $G$.
 
 \textbf{Then} let $W_1$ be the $(3b + 4)(f_{\ref{thm:FolioHomoWall}}( g_{\ref{prop:FindGoodFlatWall}}( t(b,d) ), 4d, f_{\ref{thm:irrelevantVertex}}( b, d, k, 3 ) ) +1 )$-central subwall of $W_0$ which is now suspended in $G-A$ as witnessed by $\rho$ due to \zcref{prop:SuspendingAWall}. \textbf{Forward} $A,$ $W_1,$ and $\rho$ to \textbf{Step 4}.

 \item[Step 4:] \textbf{Input:} A set $A \subseteq V(G)$ of size at most $g_{\ref{prop:FindGoodFlatWall}}( t(b,d) ),$ together with a $(3b + 4)(f_{\ref{thm:FolioHomoWall}}( g_{\ref{prop:FindGoodFlatWall}}( t(b,d) ), 4d, f_{\ref{thm:irrelevantVertex}}( b, d, k, 3 ) ) +1 )$-wall $W_1$ in $G-A,$ and a sphere-rendition $\rho$ of $G-A$ witnessing that $W_1$ is suspended in $G-A$ such that $W_1$ has treewidth at most $f_{\ref{thm:Main2Folio}}(k,b,d)$ in $G$.

 \textbf{Then} run the algorithm from \zcref{lemma:BlankWall}.
 Since $\mathsf{bidim}(G,R) \leq b$ it follows that the algorithm returns a suspended $f_{\ref{thm:FolioHomoWall}}( g_{\ref{prop:FindGoodFlatWall}}( t(b,d) ), 4d, f_{\ref{thm:irrelevantVertex}}( b, d, k, 3 )$-subwall $W_2$ of $W_1$ such that $R$ is disjoint from the vertices in the compass of $W_2$.
 Moreover, the treewidth of $W_2$ is still at most $f_{\ref{thm:Main2Folio}}(k,b,d)$ in $G$ and $\rho$ witnesses that $W_2$ is suspended in $G-A$.
 \textbf{Forward} $A,$ $\rho,$ and $W_2$ to \textbf{Step 5}.

 \item[Step 5:] \textbf{Input:} A set $A \subseteq V(G)$ of size at most $g_{\ref{prop:FindGoodFlatWall}}( t(b,d) ),$ together with a $f_{\ref{thm:FolioHomoWall}}( g_{\ref{prop:FindGoodFlatWall}}( t(b,d) ), 4d, f_{\ref{thm:irrelevantVertex}}( b, d, k, 3 ))$-wall $W_2$ in $G-A,$ and a sphere-rendition $\rho$ of $G-A$ witnessing that $W_2$ is suspended in $G-A$ such that $W_2$ has treewidth at most $f_{\ref{thm:Main2Folio}}(k,b,d)$ in $G$ and its compass is disjoint from $R$.

 \textbf{Then} run the algorithm from \zcref{thm:FolioHomoWall} in $W_2$ and $A$ to find a $f_{\ref{thm:irrelevantVertex}}( b, d, k, 3 )$-subwall $W_3$ of $W_2$ that is suspended in $G-A,$ has treewidth at most $f_{\ref{thm:Main2Folio}}(k,b,d)$ in $G,$ whose compass is disjoint from $R,$ and that is $(4d,A)$-folio homogeneous.
 \textbf{Forward} $A,$ $\rho,$ and $W_3$ to \textbf{Step 6}.

 \item[Step 6:] \textbf{Input:} A set $A \subseteq V(G)$ of size at most $g_{\ref{prop:FindGoodFlatWall}}( t(b,d) ),$ together with an $f_{\ref{thm:irrelevantVertex}}( b, d, k, 3 )$-wall $W_3$ in $G-A,$ and a sphere-rendition $\rho$ of $G-A$ witnessing that $W_3$ is suspended in $G-A$ such that $W_3$ has treewidth at most $f_{\ref{thm:Main2Folio}}(k,b,d)$ in $G,$ its compass is disjoint from $R,$ and $W_3$ is $(4d,A)$-folio homogeneous.

 \textbf{Then} let $W_4$ be the $3$-central subwall of $W_3$ and let $v \in V(W_4)$ be any vertex of degree at least $3$ in $W_4$ not on the perimeter of $W_3$.
 \textbf{Return} $v$ and \textbf{terminate}.

 \item[Step 7:] \textbf{Input:} A minor model of $K_{t(b,d)}$ in $G$.
 
 \textbf{Then} run the algorithm from \zcref{lem:bigcliqueirrelevantvertexforfolio} to find an irrelevant vertex $v \in V(G),$ \textbf{return} the irrelevant vertex $v,$ and \textbf{terminate}.
\end{description}

\paragraph{Correctness of the algorithm.}
We begin by discussing the correctness of the algorithm above.

With \textbf{Step 1} we only call upon \zcref{prop:quicklyfindingaclique} and this correctly determines either that the number of edges in $G$ is below the threshold or finds a $K_{t(b,d)}$-minor model in $G$.
If the minor model is found, \zcref{lem:bigcliqueirrelevantvertexforfolio} correctly finds an irrelevant vertex $v$.
Indeed, this vertex $v$ must now be strongly irrelevant due to the way it arises from the clique-minor and thus, \textbf{Step 1} is correct.

Similarly, in \textbf{Step 2}, \textbf{Case 2.1} also refers to \textbf{Step 7} and is therefore correct.
\textbf{Case 2.2} is the only case where the algorithm might determine that the treewidth of $G$ is small which is one of the two valid outcomes of \zcref{thm:irrelevantVertex} and therefore safe.
Finally, \textbf{Case 2.3} finds the initial large wall $W_0$ which is flat in $G-A$ for the additional set $A \subseteq V(G)$ of bounded size as guaranteed by \zcref{prop:FindGoodFlatWall}.
From here, the wall, its rendition, and the apex set are passed through \zcref{prop:SuspendingAWall}, \zcref{lemma:BlankWall}, and \zcref{thm:FolioHomoWall} in order.
The result of the application \zcref{thm:FolioHomoWall} in \textbf{Step 5} is the guarantee that all vertices in the compass of the $3$-central subwall of $W_3$ are irrelevant.
That is, if $x \in V(G)$ is a vertex from the compass of $W_3,$ then $(k,d)\text{-}\mathsf{folio}(G - x,R) = (k,d)\text{-}\mathsf{folio}(G,R)$.
Continuing to \textbf{Step 6} the algorithm now chooses a vertex from $W_4$ that has degree $3$ in $W_4$ but does not belong to the perimeter of $W_4$.
Since $W_4$ is an $r$-wall with $r \geq 3$ such a vertex always exists and by the arguments above, $v$ is indeed irrelevant for the $(k,d)$-folio of $(G,R)$.
This concludes the discussion on the correctness of the algorithm.

\paragraph{Running time analysis.}
What is left is to discuss the running time.
There are essentially only two crucial steps requiring discussion.

The first is an application of \zcref{prop:quicklyfindingaclique} running in time $2^{\mathbf{O}( t(b,d)^2 \log t(b,d) )} | V(G) |$.
If this step is successful it immediately jumps to \textbf{Step 7}, returning an irrelevant vertex in linear time.
Otherwise, the algorithm proceeds to \textbf{Step 2}.
Since we did not immediately jump from \textbf{Step 1} to \textbf{Step 7}, from here we may assume that $|E(G)| \in 2^{\mathbf{O}( t(b,d)^2 \log t(b,d) )} | V(G) |$ which allows us to replace all dependencies on the number of edges in all of the following algorithms with the number of vertices of $G$.
This guarantees that the entire algorithm runs in time linear in $|G|$.

Since $f_{\ref{thm:Main2Folio}}(k,b,d) \in 2^{\mathbf{poly}( b + d )} \cdot \mathbf{poly}( k ),$ the only other step where the running time cannot be deduced immediately from the running times of the subroutines is the combination of \textbf{Step 2} and \textbf{Step 5}.
The algorithm from \zcref{prop:FindGoodFlatWall} runs in time 
\begin{align*}
2^{2^{\mathbf{O}(t(b,d)^2 \log t(b,d))} \cdot f_{\ref{thm:Main2Folio}}(k,b,d) \log f_{\ref{thm:Main2Folio}}(k,b,d) + \mathbf{O}(f_{\ref{thm:Main2Folio}}(k,b,d)^2)} \cdot |G| + 2^{ 2^{\mathbf{O}(t(b,d)^2 \log t(b,d))} \cdot f_{\ref{thm:Main2Folio}}(k,b,d)^3 \log f_{\ref{thm:Main2Folio}}(k,b,d) }.
\end{align*}
With the bound on $f_{\ref{thm:Main2Folio}}(k,b,d)$ above this means that this step of the algorithm runs in time
\begin{align*}
 2^{2^{\mathbf{poly} ( b + d )} \cdot \mathbf{poly} ( k )} \cdot |G|
\end{align*}
which still fits into our predicted bounds.
Notice that the treewidth of $W_i$ in $G$ is at bounded by $f_{\ref{thm:Main2Folio}}(k,b,d)$ for all $i\in[4]$.
The running time of the algorithm from \zcref{thm:FolioHomoWall} is
\begin{align*}
 2^{\mathbf{poly} ( w + a +d)} \cdot \mathbf{poly}( r ) \cdot (|G| + |\!|G|\!|)
\end{align*}
which, due to our analysis of the consequence of \textbf{Step 1}, simplifies to
\begin{align*}
 2^{\mathbf{poly} ( w + a +d)} \cdot \mathbf{poly}( r ) \cdot 2^{\mathbf{poly}( b + d)} \cdot |G|
\end{align*}
in our case.
The value $a$ in the expression above is the size of $A$ which is in $\mathbf{O}( t(b,d)^{24} )$ by \zcref{prop:FindGoodFlatWall}.
The value $w$ on the other hand is the treewidth of $W_3$ which satisfies $f_{\ref{thm:Main2Folio}}(k,b,d) \in 2^{\mathbf{poly}( b + d)} \cdot \mathbf{poly}( k )$.
It follows that \textbf{Step 5} runs in time
\begin{align*}
 2^{2^{ \mathbf{poly}( b + d ) } \cdot \mathbf{poly}( k )}|G|
\end{align*}
as required.
This completes the running time analysis of our algorithm and the proof of \zcref{thm:Main2Folio}.

\section{Two lower bounds for the irrelevant vertex technique}\label{sec:lowerbounds}

In this section we show that the irrelevant vertex technique cannot be expected to yield substantially better bounds than those in our main results. The two lower bounds rely on the same general principle: we construct graphs of very large treewidth in which every vertex is essential for preserving a prescribed minor-folio behaviour. The first construction exploits the graph $\hat{\Gamma}_k$ together with its vital linkage structure and yields an exponential lower bound in terms of bidimensionality. The second construction enriches $\hat{\Gamma}_k$ with many carefully chosen non-isomorphic attachments, forcing any potential irrelevant vertex to preserve a large amount of combinatorial information, which in turn yields a lower bound in terms of the parameter $d$.
\medskip

Our starting point is the family of graphs introduced in \cite{AdlerKKLST2011Tight,AdlerK2019Lower}, which provides canonical obstructions for irrelevant vertex arguments based on linkage uniqueness. Intuitively, the graph $\hat{\Gamma}_k$ contains a highly rigid system of $k$ terminal pairs: there is a linkage connecting these pairs, and this linkage is so constrained that deleting any vertex may destroy the possibility of realizing the same pattern.

In \cite{AdlerKKLST2011Tight,AdlerK2019Lower} a parametric graph $\hat{\Gamma}_{k}, k\in\Nbbb_{\geq 2}$
was defined as follows: 
Let $m_k=2^{k}-1$.
We start with a $(m_k\times m_k)$-grid
whose ``leftmost'' (resp. ``rightmost'') vertices are $v_{1},\ldots,v_{m_k}$ (resp. $u_{1},\ldots,u_{m_k}$), ordered as they appear from above to below.
We first add all edges in the set $\{u_{i}u_{m_k-i+1}\mid i\in[m_{k-1}]\}$
and for $i\in[2,k]$ we add the edges in 
$$E_{i}\coloneqq \{v_{2^{i}+j}v_{2^{i+1}-j}\mid j\in [2^{i}-1]\}.$$

\begin{figure}
\begin{center}
\scalebox{.3}{\includegraphics{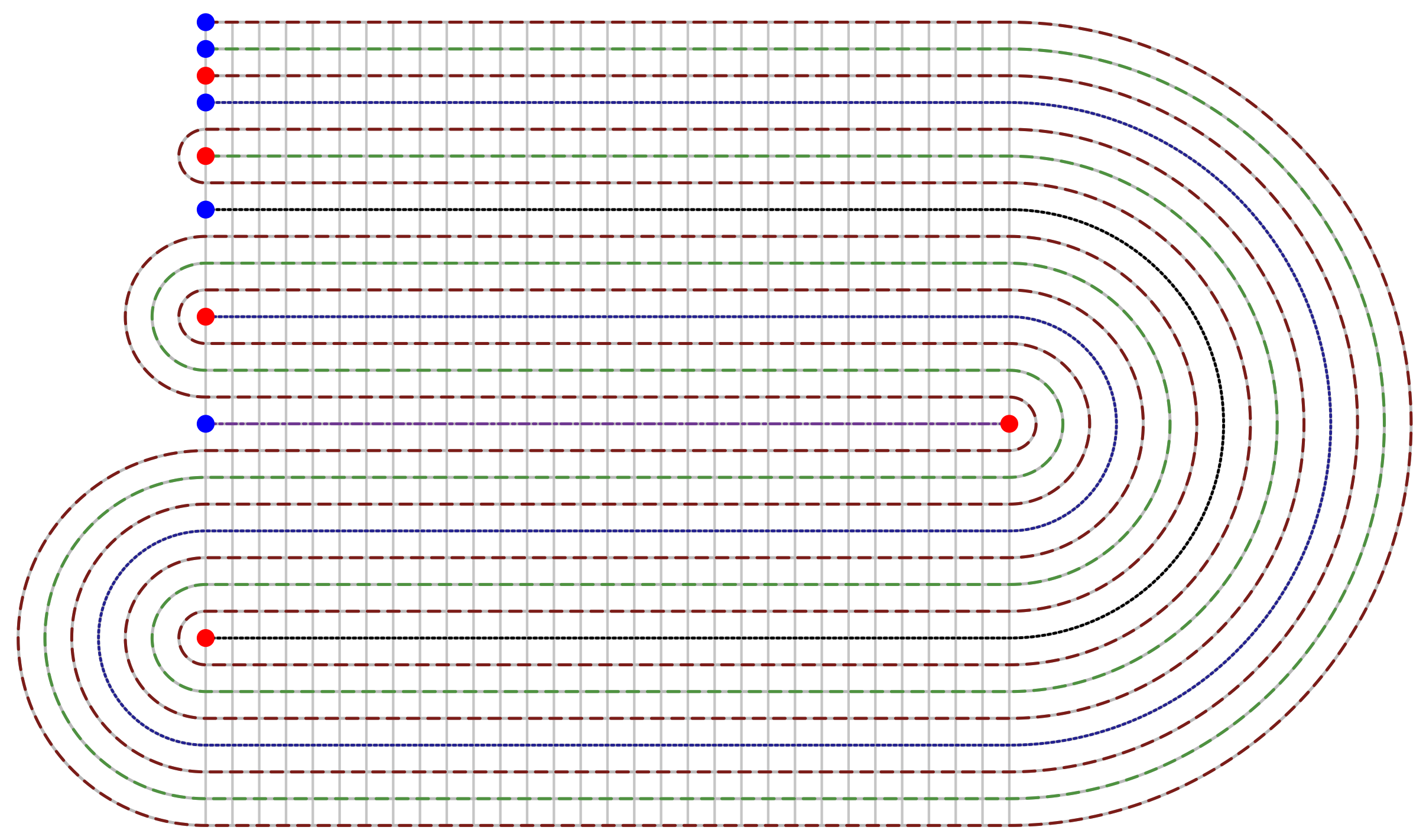}}
\end{center}
\caption{The graph $\hat{\Gamma}_{5}$ with a vital linkage $L$ of order $5$ in it.}
\label{fig_g_vit_lin}
\end{figure}

 We define $\tau_{k}\coloneqq \{(s_{i},t_{i})\mid i\in[k]\}$ as the set of pairs of vertices (terminals) in $\hat{\Gamma}_{k},$ where $s_{i}=v_{2^{i-1}}, i \in[k],$ $t_{i}=v_{3\cdot 2^{i-1}}, i\in[k-1],$ and $t_{k}=u_{2^{i-1}}$.
We also denote by $T_{k}$ the vertices in the pairs of $\tau_{k}$.
For a drawing of $\hat{\Gamma}_{5}$ 
and vital linkage in it with $5$ pairs of terminals, see \cref{fig_g_vit_lin}.

\begin{proposition}
\label{prop_etal_al_all}
For every $k\in \Nbbb_{\geq 2}$ the 
graph $\hat{\Gamma}_{k}$ contains a vital linkage $L$ of order $k$ whose pattern in $\tau_{k}$.
\end{proposition}

The relevance of \cref{prop_etal_al_all} is that it combines two features that are central for lower bounds: on the one hand, the graph $\hat{\Gamma}_k$ has exponentially large treewidth in $k,$ and on the other hand, the linkage pattern prescribed by $\tau_k$ is rigid. This tension is exactly what makes $\hat{\Gamma}_k$ a natural witness against overly optimistic irrelevant vertex bounds.

As observed in \cite{AdlerKKLST2011Tight,AdlerK2019Lower}, 
$\hat{\Gamma}_{k}$ contains as a subgraph a $(m\times m)$-grid, therefore,
$\tw(\hat{\Gamma}_{k})\geq 2^{k}-1$. Also, it is known that every $(m\times m)$-grid
has a path decomposition $\{X_{1},\ldots,X_{m(m-1)}\},$ where $\{v_{1},\ldots,v_{m}\}\subseteq X_{1}$ and $\{u_{1},\ldots,u_{m}\}\subseteq X_{m(m-1)}$.
This readily implies that $\tw(\hat{\Gamma}_{k}) = 2^{k}-1$ (actually the same we can say for the pathwidth of $\hat{\Gamma}_{k}$).

\subsection{Optimality of \cref{th_our_result}}

We first revisit the construction $\hat{\Gamma}_k$ from the perspective of the parameter $\bidim(\hat{\Gamma}_k,T_k)$. The point is that the lower bound from \cite{AdlerKKLST2011Tight,AdlerK2019Lower} is not merely an artifact of the number of terminals: it remains significant even when expressed through bidimensionality. Thus, to understand the limitations of \cref{th_our_result}, it is enough to estimate how the bidimensionality of $(\hat{\Gamma}_k,T_k)$ grows with $k$.

For every $k\in\Nbbb_{\geq 2}$ we set $b_{k}=\bidim(\hat{\Gamma}_{k},T_{k})$.

\begin{lemma}
\label{lem_b_bidim}
There is a constant $c$ such that, for every $k\in \Nbbb_{\geq 2},$ $\nicefrac{1}{2}\sqrt{b_{k}}\leq k\leq c\cdot b_{k}^2$. 
\end{lemma}

\begin{figure}
\begin{center}
\scalebox{.289}{\includegraphics{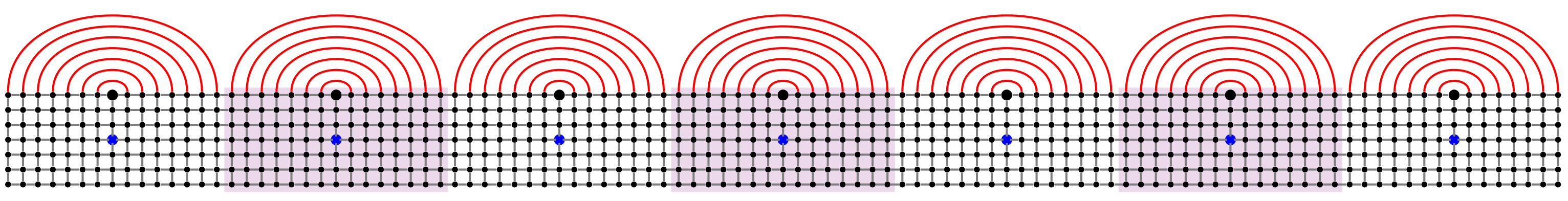}}
\end{center}\caption{The annotated graph $(Z_{7},A_{7})$.}
\label{fig_s_z_5}
\end{figure}

\begin{proof}
Notice that 
for every annotated graph $(G,S),$ $|S|\geq (\bidim(G,S))^2$.
This fact applied for $(\hat{\Gamma}_{k},T_{k})$ and given that $|T_{k}|=2k,$ implies $\nicefrac{1}{2}\sqrt{b_{k}}\leq k$.

We define the graph $Z_{s}$ by considering the $(s\times (s(2s+1))$-grid, where we denote 
the top vertices ordered from the left to the right by $x_{1},\ldots,x_{2s(s+1)},$
and where, for every $i\in[s],$ we add all edges in $\{x_{(2s+1)(i-1)+j}x_{(2s+1)i-j+1}\mid j\in[s]\}$. We also set $A_{s}=\{x_{(2s+1)(i-1)+s+1}\mid i\in[s]\}$ (see \cref{fig_s_z_5} for a drawing of $(Z_{5},A_{5})$ where the added edges are depicted in \red{red}).

We claim that there is a constant $c'\in(0,1)$ 
such that $\bidim(Z_{t^2},A_{t^2})\geq c't$.
To see this, observe first 
that $(Z_{2s},A_{2s}),$ for odd $s,$
contains the $(s\times 2s(s+1))$-grid as a minor where the $s$ vertices in $B_{s}\coloneqq \{((2s+1)(i-1)+s+1,(s+1)/2)\mid i\in[s]\}$ correspond to $s$ connected sets of $(Z_{2s},A_{2s}),$ each containing at least one vertex from $A_{2s}$. 
To visualise $B_{s}$ consider 
the $(s\times 2s(s+1))$-grid drawn in \cref{fig_s_z_5} 
for $s=7$ by omitting the \red{red} edges; the set $B_{s}$ consists of the \blue{blue} vertices.
Also to visualise the minor relations, each shaded square of $Z_{2s}$ is contracted to a single ``vertical'' path containing one of the vertices in $B_{s}$ in its centre. 
Next observe that 
the $(s\times 2s(s+1))$-grid, when 
$s=2\ell^2$ contains the $(\ell\times \ell)$-grid as a topological minor where the branch vertices are chosen to be $\ell^2$ vertices from $B_{s}$.
Combining these two transformations together, we have that there is some constant $c',$ such that
$\bidim(Z_{t^2},A_{t^2})\geq c't$.

By applying the claim above for $\lfloor \sqrt{t}\rfloor$ instead of $t,$ we have that 
$\bidim(Z_{t},A_{t})\geq c'\cdot \lfloor \sqrt{t} \rfloor$.

Next observe that $\hat{\Gamma}_{k}$ contains a subdivision of $Z_{\lfloor k/2\rfloor},$ where the vertices of $A_{\lfloor k/2\rfloor}$ correspond to the terminals $t_{\lfloor k/2\rfloor},\ldots,t_{k}$. Combining these last two facts together, we obtain that 
$$b_{k}=\bidim(\hat{\Gamma}_{k},T_{k})\geq \bidim(\hat{\Gamma}_{k},\{t_{\lfloor k/2\rfloor},\ldots,t_{k}\})\geq \bidim(Z_{\lfloor k/2\rfloor},A_{\lfloor k/2\rfloor})\geq c'\cdot \lfloor\sqrt{k}/2\rfloor.$$ Therefore there is a positive integer $c$ such that $k\leq c\cdot b_{k}^{2}$.
\end{proof}

The estimate in \cref{lem_b_bidim} shows that the bidimensionality of $(\hat{\Gamma}_k,T_k)$ is polynomially related to the order of the vital linkage. In particular, replacing the number of terminals by bidimensionality does not fundamentally weaken the obstruction: the exponential behaviour in $k$ still translates into an exponential lower bound in a polynomial function of $b_k$.

We now observe that the exponential lower bound from \cite{AdlerKKLST2011Tight,AdlerK2019Lower} applies not only in terms of the number of terminals, but also in terms of bidimensionality.

\begin{corollary}
\label{th_l_tight_link}
If \cref{th_our_result} holds for some function $\beta:\mathbb{N}^2\to\mathbb{N},$ then $\beta(0,b)=2^{\Omega(\sqrt{b})}$.
\end{corollary}

\begin{proof}
Assume that \cref{thm:intro_folio_irrelevant_informal} holds for some function
$\beta:\mathbb{N}^2\to\mathbb{N}$. Then, for every $k\geq 2,$ we may apply it to the annotated 
graph $(\hat{\Gamma}_{k},T_{k})$ and deduce by \cref{prop_etal_al_all} that 
$\tw(\hat{\Gamma}_{k})<\beta(k,b)$.
Therefore, as $\tw(\hat{\Gamma}_{k})=2^{k}-1,$ by \cref{lem_b_bidim}, 
$\beta(k,b)=2^{\Omega(\sqrt{b})}$.
\end{proof}

Notice that in \cref{th_l_tight_link} we only used the ``easy'' upper bound to $b_{k}$ in \cref{lem_b_bidim}.
We chose to also argue on how to lower bound $b_k$ in terms of $k,$ in order to illustrate that the gap between the size of the terminal set $T_k$ and its bidimensionality in $\hat{\Gamma}_k$ is polynomial.

\subsection{Optimality of \cref{thm:intro_folio_irrelevant_informal}}

We now turn to lower bounds for the folio version of the irrelevant vertex technique. In contrast with the previous subsection, where the obstruction came from the rigidity of a single linkage instance, here we build a family of graphs that encodes many pairwise non-isomorphic gadgets of small size. The role of these gadgets is to make every vertex indispensable: deleting a vertex necessarily destroys the ability to realise a certain prescribed minor.

\paragraph{Asymptotic enumeration of 5-regular graphs.}
The first ingredient is a supply of many small pairwise non-isomorphic connected graphs of bounded degree. We obtain these from asymptotic enumeration of regular graphs as follows.

\begin{lemma}
\label{lem_ca_antich_5}
For all sufficiently large even integers $n,$ there exist at least $n^{n}$ pairwise non-isomorphic connected $5$-regular graphs on $n$ vertices.
\end{lemma}

\begin{proof}
For fixed degree $r,$ the number $\zeta_n^{(r)}$ of labelled (simple) $r$-regular graphs on $n$ vertices satisfies
$$\zeta_n^{(r)}\sim e^{-(r^2-1)/4}\,\frac{(rn)!}{(rn/2)!\,2^{rn/2}\,(r!)^n},$$
by the asymptotic formula of Bender and Canfield \cite{BenderCanfield1978} (see also \cite{Bollobas1982}). Specializing to $r=5$ gives
$\zeta_n^{(5)}\sim e^{-6}\,\frac{(5n)!}{(5n/2)!\,2^{5n/2}\,120^n},$
for even $n$.
Moreover, Wormald \cite{Wormald1981} proved that for every fixed $r\ge 3,$ almost all labelled $r$-regular graphs are connected; hence almost all labelled $5$-regular graphs are connected, and therefore the labelled connected count $\zeta_n^{(5,\mathrm{conn})}$ satisfies
$\zeta_n^{(5,\mathrm{conn})}\sim \zeta_n^{(5)}$
as $n\to\infty$ through even integers.
Finally, Bollob\'as \cite{Bollobas1982} proved that for fixed $r\ge 3,$ almost all labelled $r$-regular graphs have trivial automorphism group. Hence, if $\tilde{\zeta}_n^{(5,\mathrm{conn})}$ denotes the number of pairwise non-isomorphic connected $5$-regular graphs on $n$ vertices, then
$\tilde{\zeta}_n^{(5,\mathrm{conn})}\sim \frac{\zeta_n^{(5,\mathrm{conn})}}{n!}
\sim e^{-6}\,\frac{(5n)!}{n!\,(5n/2)!\,2^{5n/2}\,120^n},$
again for even $n$.

On the other hand, by Stirling's formula,
$$\frac{1}{n^{n}}\cdot e^{-6}\,\frac{(5n)!}{n!\,(5n/2)!\,2^{5n/2}\,120^n}
\sim
\frac{e^{-6}}{\sqrt{\pi n}}\left(\frac{5\sqrt{5}}{24e^{3/2}}\sqrt{n}\right)^n,$$
and the right-hand side tends to $\infty$ as $n\to\infty$ through even integers. Therefore
$\frac{\tilde{\zeta}_n^{(5,\mathrm{conn})}}{n^{n}}\to\infty$.
In particular, for all sufficiently large even $n,$ one has
$\tilde{\zeta}_n^{(5,\mathrm{conn})}\ge n^{n}$.
\end{proof}

As an immediate consequence of \cref{lem_ca_antich_5}, for every sufficiently large integer $m,$ there exist $m$ pairwise non-isomorphic connected $5$-regular graphs, each on at most
$\min\{\ell\in 2\mathbb{N}:\ \ell^{\ell}\ge m\}$
vertices. In particular, all these graphs may be chosen to have at most
$\left(1+o(1)\right)\frac{\log m}{\log\log m}$
vertices as $m\to\infty$.
We record this in the following form.

\begin{lemma}
\label{lem_ca_antich_5_small}
There exists a function $f_{\ref{lem_ca_antich_5_small}}:\Nbbb\to\Nbbb$ such that, for every sufficiently large integer $k,$ there exist $k$ pairwise non-isomorphic connected $5$-regular graphs, each on $f_{\ref{lem_ca_antich_5_small}}(k)$ vertices. Moreover,
$f_{\ref{lem_ca_antich_5_small}}(k)=O\!\left(\frac{\log k}{\log\log k}\right).$
\end{lemma}

Thus we may think of $f_{\ref{lem_ca_antich_5_small}}(k)$ as the cost of encoding $k$ distinct isomorphism types by bounded-degree connected graphs.

\paragraph{The graph $\hat{\Gamma}^{\Acal^{k}}_{k}$.} We now combine the rigid linkage structure of $\hat{\Gamma}_k$ with the family of non-isomorphic $5$-regular gadgets obtained above. Roughly speaking, for each terminal pair $(s_i,t_i)$ we attach two copies of the same gadget $A_i,$ one on the $s_i$-side and one on the $t_i$-side. Since the gadgets corresponding to different indices are non-isomorphic, any minor model that preserves all of them is forced to match the two copies of $A_i$ to each other. This is what ultimately recovers the terminal pairing pattern $\tau_k$ inside the central graph $\hat{\Gamma}_k$.

Let $\Acal^{k}=\{A_1,\ldots,A_k\}$ be a family of pairwise non-isomorphic connected $5$-regular graphs, each on $f_{\ref{lem_ca_antich_5_small}}(k)$ vertices.

We define the graph $\hat{\Gamma}^{\Acal^{k}}_{k}$ as follows. Start with the disjoint union of $\hat{\Gamma}_{k}$ and, for each $i\in [k],$ two copies $A_i^s$ and $A_i^t$ of $A_i$. Then, for each $i\in [k],$ add two new vertices $s_i'$ and $t_i',$ add all edges between $\{s_i,s_i'\}$ and $V(A_i^s),$ and add all edges between $\{t_i,t_i'\}$ and $V(A_i^t)$. Also add the edge $s_it_i$ for all $i\in[k]$.

\medskip
A \emph{block} of a graph 
is a maximal $2$-connected subgraph of it, i.\@e.\@, every two vertices
are in a cycle.
We say that a subgraph of $G$ is a \emph{bridge}\footnote{Not to be confused by the concept of $H$-bridge, defined in \cref{subsec:HomoMuralis}.} of $G$ 
if it is isomorphic to $K_{2}$ and the removal of its unique edge increases the number of connected components. In what follows, we use the term \emph{bridge} for 
both the subgraph and the edge contained in it.

\begin{observation}
\label{obs_minor_bl}
If $H$ is a minor of $G,$ then every block of $H$ is a minor of some block of $G$.
\end{observation}

Observe that $\hat{\Gamma}^{\Acal^{k}}_{k}$ has exactly $2k+1$ blocks: the block $\hat{\Gamma}_{k},$ together with the $2k$ blocks
$G_1^s,\ldots,G_k^s,G_1^t,\ldots,G_k^t,$
where, for each $i\in [k],$ the graph $G_i^s$ is the subgraph of $\hat{\Gamma}^{\Acal^{k}}_{k}$ induced by $V(A_i^s)\cup\{s_i,s_i'\},$ and similarly $G_i^t$ is the subgraph induced by $V(A_i^t)\cup\{t_i,t_i'\}$.

Let
$\Gcal\coloneqq \{G_1^s,\ldots,G_k^s,G_1^t,\ldots,G_k^t\}$.
Since each graph $A_i$ is $5$-regular, every graph in $\Gcal$ is $3$-connected and has minimum degree at least $7$. In particular, we obtain the following.

\begin{observation}
\label{obs_remain_non_pl_2con}
For every $G\in \Gcal$ and every vertex $v\in V(G),$ the graph $G-v$ is non-planar and $2$-connected.
\end{observation}
 
 Since the treewidth of a graph is the maximum of the treewidths of its blocks, and since every block in
$\{G_1^s,\ldots,G_k^s,G_1^t,\ldots,G_k^t\}$
has order $2+f_{\ref{lem_ca_antich_5_small}}(k),$ each such block has treewidth at most
$1+f_{\ref{lem_ca_antich_5_small}}(k)$.
As
$f_{\ref{lem_ca_antich_5_small}}(k)=O\!\left(\frac{\log k}{\log\log k}\right),$
it follows that
$1+f_{\ref{lem_ca_antich_5_small}}(k)<2^k-1$
for all sufficiently large $k$.
Therefore the block $\hat{\Gamma}_k$ determines the treewidth of $\hat{\Gamma}^{\Acal^k}_k,$ and we obtain the following.

\begin{observation}
\label{obs_ad_conn_tw}
For every sufficiently large $k,$
$\tw(\hat{\Gamma}^{\Acal^{k}}_{k})=2^{k}-1$.
\end{observation}

\paragraph{The graph $H_{k}$.}
We define $H_k$ to be the graph obtained from the disjoint union of
$G_1^s,\ldots,G_k^s,G_1^t,\ldots,G_k^t$
by adding, for each $i\in [k],$ the edge $s_it_i$.
Then $H_k$ has exactly $k$ connected components, and for each $i\in [k],$ exactly one of these components consists of the two blocks
$G_i^s$ and $G_i^t$ joined by the bridge $s_it_i$.

The graph $H_k$ should be viewed as the target minor that records exactly the information carried by the attached gadgets and by the pairing of the terminals.

 \begin{observation}
 \label{obs_new_size_H_k}
For every $k,$ $|V(H_k)|=2k(2+f_{\ref{lem_ca_antich_5_small}}(k)).$
\end{observation}
 
 The next lemma is the heart of the construction. It shows that $H_k$ is present as a minor in $\hat{\Gamma}^{\Acal^k}_k,$ but becomes absent after the deletion of any vertex. In other words, no vertex of $\hat{\Gamma}^{\Acal^k}_k$ is irrelevant with respect to the minor $H_k$.

\begin{lemma}
\label{lem_other_b}
The graph $H_k$ is a minor of $\hat{\Gamma}^{\Acal^{k}}_{k},$ and for every $v\in V(\hat{\Gamma}^{\Acal^{k}}_{k}),$ the graph $H_k$ is not a minor of $\hat{\Gamma}^{\Acal^{k}}_{k}-v$.
\end{lemma}
\begin{proof}
Let $L$ be the linkage in $\hat{\Gamma}_k;$ we also regard it as a linkage in
$\hat{\Gamma}^{\Acal^k}_k$.
Then $H_k$ is a minor of $\hat{\Gamma}^{\Acal^k}_k$: indeed, delete all edges of the planar block $\hat{\Gamma}_k$ that do not belong to $L,$ and then dissolve all vertices of degree $2$ in that block.

Now fix a vertex $v\in V(\hat{\Gamma}^{\Acal^k}_k),$ and suppose for a contradiction that
$\hat{\Gamma}^{\Acal^k}_k-v$ contains a minor model of $H_k$.
Since the blocks of $H_k$ are precisely
$G_1^s,\ldots,G_k^s,G_1^t,\ldots,G_k^t,$
every block of $H_k$ is non-planar.
Hence, by \cref{obs_minor_bl}, every block of $H_k$ must be a minor of a non-planar block of
$\hat{\Gamma}^{\Acal^k}_k-v$.
The only planar block of $\hat{\Gamma}^{\Acal^k}_k-v$ is the block inherited from
$\hat{\Gamma}_k-v,$ while by \cref{obs_remain_non_pl_2con} every other block of
$\hat{\Gamma}^{\Acal^k}_k-v$ is non-planar.

Moreover, every block of $H_k$ has exactly
$2+f_{\ref{lem_ca_antich_5_small}}(k)$ vertices, and the same is true for each graph in
$\{G_1^s,\ldots,G_k^s,G_1^t,\ldots,G_k^t\}$.
Therefore each block of $H_k$ must be a minor of one of these blocks.
Since a connected graph cannot be a proper minor of a connected graph on the same number of vertices, each such minor relation is in fact an isomorphism.

It follows that $v$ does not belong to any graph in
$\{G_1^s,\ldots,G_k^s,G_1^t,\ldots,G_k^t\}$:
indeed, if $v$ belonged to one of them, then that block in
$\hat{\Gamma}^{\Acal^k}_k-v$ would have fewer than
$2+f_{\ref{lem_ca_antich_5_small}}(k)$ vertices, and hence could not contain any block of $H_k$ as a minor.
Thus $v$ is a non-terminal vertex of the planar block $\hat{\Gamma}_k$.

Consequently, the minor model yields a bijection between the blocks of $H_k$ and the blocks
$G_1^s,\ldots,G_k^s,G_1^t,\ldots,G_k^t$
of $\hat{\Gamma}^{\Acal^k}_k-v$.
By construction, for distinct indices the corresponding graphs are non-isomorphic, whereas
$G_i^s\cong G_i^t$ for every $i\in [k]$.
Hence, for each $i\in [k],$ the bridge $s_it_i$ of $H_k$ must be represented by a path in
$\hat{\Gamma}_k-v$ joining the attachment vertices $s_i,t_i$ of the two isomorphic blocks
$G_i^s$ and $G_i^t$.
These paths are pairwise internally vertex-disjoint, so they form a linkage in
$\hat{\Gamma}_k-v$ with pattern $\tau_k$.

This contradicts \cref{prop_etal_al_all}, since $L$ is a vital linkage of
$\hat{\Gamma}_k$.
Therefore $H_k$ is not a minor of $\hat{\Gamma}^{\Acal^k}_k-v$ for any
$v\in V(\hat{\Gamma}^{\Acal^k}_k)$.
\end{proof}

We are now ready to extract the quantitative consequence. The graph $\hat{\Gamma}^{\Acal^k}_k$ still has treewidth $2^k-1,$ while the excluded minor $H_k$ has only
$2k(2+f_{\ref{lem_ca_antich_5_small}}(k))$
vertices. Since
$f_{\ref{lem_ca_antich_5_small}}(k)=O\!\left(\frac{\log k}{\log\log k}\right),$
this converts the exponential treewidth scale in $k$ into the lower bound
$2^{\Omega(d\log\log d/\log d)}$ in terms of the folio parameter $d$.

\begin{theorem}
If \cref{thm:intro_folio_irrelevant_informal} holds for some function $\beta:\mathbb{N}^3\to\mathbb{N},$ then $\beta(0,0,d)=2^{\Omega (\frac{d\log\log d}{\log d})}$.
\end{theorem}

\begin{proof}
Assume that \cref{thm:intro_folio_irrelevant_informal} holds for some function
$\beta:\mathbb{N}^3\to\mathbb{N}$.
For each sufficiently large integer $k,$ let
$d\coloneqq 2k\bigl(2+f_{\ref{lem_ca_antich_5_small}}(k)\bigr).$
Applying \cref{thm:intro_folio_irrelevant_informal} with the first two parameters equal to $0,$ we obtain that for every graph $G$ with
$\tw(G)\ge \beta(0,0,d),$
there exists a vertex $v\in V(G)$ such that
$(0,d)\text{-}\mathsf{folio}(G,\emptyset)=(0,d)\text{-}\mathsf{folio}(G-v,\emptyset).$

By \cref{obs_new_size_H_k}, the graph $H_k$ has exactly $d$ vertices. Hence, if $G$ contains $H_k$ as a minor, then $H_k\in (0,d)\text{-}\mathsf{folio}(G,\emptyset),$ and therefore also
$H_k\in (0,d)\text{-}\mathsf{folio}(G-v,\emptyset)$.
Equivalently, $G-v$ contains $H_k$ as a minor.

Now apply this with
$G\coloneqq \hat{\Gamma}^{\Acal^{k}}_{k}.$
Since $H_k$ is a minor of $\hat{\Gamma}^{\Acal^{k}}_{k}$ but is not a minor of
$\hat{\Gamma}^{\Acal^{k}}_{k}-v$ for any
$v\in V(\hat{\Gamma}^{\Acal^{k}}_{k})$ by \cref{lem_other_b}, it follows that
$\tw(\hat{\Gamma}^{\Acal^{k}}_{k})<\beta(0,0,d).$
Combining this with \cref{obs_ad_conn_tw}, we obtain
$\beta(0,0,d)>2^k-1$
for all sufficiently large $k$.

Finally, since
$f_{\ref{lem_ca_antich_5_small}}(k)=O\!\left(\frac{\log k}{\log\log k}\right),$
we have
$d=2k\bigl(2+f_{\ref{lem_ca_antich_5_small}}(k)\bigr)
=O\!\left(k\frac{\log k}{\log\log k}\right).$
By inversion, this yields
$k=\Omega\!\left(\frac{d\log\log d}{\log d}\right).$
Therefore
$\beta(0,0,d)>2^k-1=2^{\Omega\left(\frac{d\log\log d}{\log d}\right)},$
and hence
$\beta(0,0,d)=2^{\Omega\left(\frac{d\log\log d}{\log d}\right)}.$
\end{proof}

\section{Conclusion}\label{sec:conclusion}

This paper provides strong evidence that the algorithmic side of the Graph Minors Series is much closer to the ``solar system'' of mainstream parameterized algorithms than its longstanding ``galactic'' reputation may suggest.

\subsection{Résumé}

Our main contribution is a quantitatively controlled version of the irrelevant-vertex machinery. More precisely, we prove that the vital-linkage threshold and, more generally, the irrelevant-vertex threshold for folios admit bounds of the form
$\beta(k,b,d)\in 2^{\poly(b+d)}\cdot \poly(k),$
where the exponential dependence is governed by the bidimensionality parameter and the detail of the folio, rather than the number of terminals.
Our lower bounds further show that the combinatorial source of this exponential behaviour can be identified precisely, isolated, and bounded in an essentially optimal way.

Beyond the technical core of the proof, we see our results as part of a broader program whose goal is to make the structural and algorithmic content of the graph minors framework quantitatively explicit.
Our arguments build on a sequence of recent works \cite{Chuzhoy2015Excluded,GorskySW2025Polynomial,GorskySTW2025Catching,GorskySW2026Price,GorskyPW2026Quickly} that replaced enormous or implicit dependencies in several key ingredients of the series by polynomial ones (see also \cite{KawarabayashiTW2021Quickly,SauST21amor,ThilikosW2023Graph,PaulPTW2025Localstructure,ProtopapasTW2025Colorful,ThilikosW2026Excludingsurfaces} for precursors).
In this sense, the present paper may be viewed as the culmination of this line of research: once these ingredients are brought under polynomial control, the algorithmic heart of the graph minors framework~--~namely, the irrelevant-vertex technique for routing and rooted pattern detection~--~can finally be restated with bounds that are explicit, structurally meaningful, and essentially optimal.

\subsection{Bidimensionality, annotation irrelevancy, and more}
\label{subsec_comp_folios_un_an}

In the introduction we mentioned, as a direct consequence of \cref{thm:intro_folio_algorithm_informal}, the following.

\begin{lemma}
\label{lem_with_r_in_it}
There exists an algorithm that, given an annotated graph $(G,R)$ with $\textsf{bidim}(G,R)\leq~b$ and $|R|\leq k,$ computes $(d,k)\text{-}\mathsf{folio}(G,R)$ in
$2^{\mathbf{poly}(k) \cdot 2^{\mathbf{poly}(b + d)}}\cdot |G|^2$ time.
\end{lemma}

Since $\bidim(G,R)\leq \sqrt{|R|},$ this yields the following corollary.

\begin{corollary}
\label{cor_wsithr_mor_fin_it}
There exists an algorithm that, given an annotated graph $(G,R)$ with $|R|\leq k,$ computes $(d,k)\text{-}\mathsf{folio}(G,R)$ in
$2^{2^{\mathbf{poly}(k + d)}}\cdot |G|^2$ time.
\end{corollary}

The algorithm in \cref{lem_with_r_in_it} is based on removing a set of vertices that produces an equivalent instance $(G’,R),$ where $G’$ is a subgraph of $G,$ $t\coloneqq \tw(G’)\in 2^{\mathbf{poly}(b + d)}\mathbf{poly}(k),$ and
\begin{eqnarray}
(k,d)\textsf{-folio}(G’,R) & = & (k,d)\textsf{-folio}(G,R). \label{eq_olk_last_equiv}
\end{eqnarray}
Since each removed vertex is \textsl{strongly irrelevant}, as certified by \cref{thm:intro_folio_irrelevant_informal}, the equality in \eqref{eq_olk_last_equiv} does not require any bound on the size of $R$.

By standard dynamic programming, computing $(k,d)\text{-}\textsf{folio}(G’,R)$ can be done in $2^{\poly(t,k,d)}\cdot |G’|,$ as the tables of the dynamic program encode $(t+k+d)$-rooted graphs of detail $d,$ of which there are $2^{\poly(t,k,d)}$. This implies the following.

\begin{theorem}
\label{lem:more_out_intro}
There exists an algorithm that, given an annotated graph $(G,R)$ with $\textsf{bidim}(G,{R})\leq~b,$ computes $(d,k)\text{-}\mathsf{folio}(G,{R})$ in
$2^{\mathbf{poly}(k) \cdot 2^{\mathbf{poly}(b + d)}}\cdot |G|^2$ time.
\end{theorem}

A natural challenge beyond \cref{lem:more_out_intro} would be to eliminate the dependence on $b$ from the complexity of the algorithm, again without assuming any bound on $R$. We believe that there exists an algorithm computing the $(k,d)$-folio of an annotated graph $(G,R)$ in
\begin{eqnarray}
2^{2^{2^{\mathbf{poly}(k + d)}}}\cdot |G|^2 \label{eq_ola_dp_flolio}
\end{eqnarray}
time. In other words, we conjecture that \cref{cor_wsithr_mor_fin_it} holds without any restriction on $R$.
An $f(k,d)\cdot |G|^{2}$-time algorithm~-- albeit with a huge and unspecified $f$ --~already follows from the meta-algorithmic results of \cite{GolovachST2023ModelChecking,SauST2025Parameterizing}.

We believe that such an algorithm can be achieved by employing the ``annotation irrelevant’’ argument that is implicit in \cite{SauST2025Parameterizing} (see also \cite{GolovachST2023ModelChecking,GolovachSG2023Hitting}), where bidimensionality plays a pivotal role. (For the combinatorial and algorithmic properties of bidimensionality, see \cite{PaulPTW2024Obstructions,PaulPTW2025Localstructure,ProtopapasTW2025Colorful}.)

Roughly speaking, the approach is as follows: if $R$ has sufficiently large bidimensionality, then it is possible to find a vertex in $R$ that is \textsl{annotation irrelevant}; removing it from $R$~-- while keeping it in the graph --~yields an equivalent instance. By iteratively removing such vertices from $R,$ we can eventually bound the bidimensionality of the annotated set $R,$ at which point \cref{lem:more_out_intro} becomes applicable.
Working out the details of this procedure is technically involved, as one must carefully control dependencies to remain within polynomial bounds. Addressing these issues is beyond the scope of this paper.

Nevertheless, we emphasise that the idea of removing ``annotation irrelevant’’ vertices from the annotation set is not specific to folio computation, but rather indicative of a broader principle. It underlies the meta-algorithmic framework of \cite{SauST2025Parameterizing}, which is driven by low bidimensionality and whose parameter dependencies are strongly influenced by the linkage function, which we now optimise. We expect that many of the algorithmic consequences of \cite{SauST2025Parameterizing} can now be translated into more streamlined algorithms using our optimal bounds on the linkage function.

More generally, for annotated graphs (or even colourful graphs \cite{ProtopapasTW2025Colorful}), bidimensionality plays a role analogous to that of treewidth in many graph problems: it provides a threshold parameter above which annotation-irrelevant vertices can be found and below which irrelevant vertices can be identified, until eventually a bounded-treewidth equivalent instance is obtained.

Based on this general algorithmic scheme, we expect that our bounds and algorithms will have sufficient gravitational pull to liberate more algorithms from their “galactic” reputation. This, however, should be the subject of a broader future project.

\subsection{Other research directions}
\label{subsec_ta_more_oro}

\paragraph{Reducing degrees of polynomials.}
A natural question is how far the quantitative bounds obtained here can be sharpened further. The near-optimality of our results leaves limited room for improving the general asymptotic form of the Linkage Function from \cref{th_our_result}. Still, even within a bound of the form
$\beta(k,b)\in 2^{\poly(b)}\cdot \poly(k),$
it would be very interesting to reduce the degrees of the polynomials appearing both in the exponent and in the multiplicative factor. Our current proof still produces rather large constants in these degrees, and reducing them would require a careful quantitative refinement of several parts of the argument and the results which feed into our proof ~--~ particularly \zcref{thm:localstructure}. Such improvements would be valuable not only from a technical point of view, but also for applications in more specialised settings. In particular, when one studies the Linkage Function under the exclusion of specific parametric graph classes, additional structure may become available and could lead to substantially better estimates.

\paragraph{Other linkage-like constellations.}
Moving beyond graph minors, it would be very interesting to establish analogous linkage theorems for other rooted containment relations, such as labelled minors, immersions, or containment relations in directed graphs.
Each of these settings appears to require genuinely new ideas, but the present work suggests that one should again search for the structural parameter that governs the exponential part of the bound. It is also natural to ask whether our proof techniques can be used to improve the parameter dependencies in several recent results in structural algorithmic graph theory, such as 
the edge-disjoint linkage theorem for the strong immersion relation in Eulerian digraphs from \cite{CavallaroKK2025wellquasiordering}, and the group-labelled linkages studied in \cite{LiuY2025Disjoint}.

\paragraph{Meta-folios.}
We define the \emph{$(k,d)$-meta-folio} of an annotated graph $(G,R)$ as follows
\begin{eqnarray*}
 (k,d)\text{-meta-folio}(G,R)\coloneqq \big\{(k,d)\text{-folio}(G,R')\mid R'\subseteq R\big\} .
\end{eqnarray*}
Notice that the meta-folio of an annotated graph can be seen as a ``2nd order object'', as it contains a collection of collections of $k$-rooted graphs. 
By the strong irrelevance property in \cref{thm:intro_folio_irrelevant_informal}, the repeated application of the irrelevant vertex argument produces, in $$2^{2^{\mathbf{poly}(b + d)} \cdot \mathbf{poly}(k)} \cdot |G|^2$$ time, an annotated graph $(G',R)$ such that $G'$ is an induced subgraph of $G,$ $\tw(G')\in 2^{\mathbf{poly}(b + d)}\mathbf{poly}(k),$ and, moreover, the $(k,d)$-meta-folio of $(G,R)$ coincides with that of $(G',R)$.
Since standard dynamic programming can compute the $(k,d)$-meta-folio of $(G',R)$ in $$2^{2^{\mathbf{poly}(\tw(G'))}}\cdot |G'|$$ time, this yields an algorithm for computing the $(k,d)$-meta-folio of $(G,R)$ in
\begin{eqnarray}
 2^{2^{\mathbf{poly}(b + d)} \cdot \mathbf{poly}(k)} \cdot |G|^2+2^{2^{2^{\mathbf{poly}(b + d)} \cdot \mathbf{poly}(k)}} \cdot |G|\label{eq_many_dp_bp_dr}
\end{eqnarray}
time. 
This establishes an interesting trade-off between the parametric dependence and the polynomial dependence for the computation of the meta-folio, seen as a 2nd order object.
Can the dependencies in (\ref{eq_many_dp_bp_dr}) be improved? More broadly, do such ``higher-order'' extensions of folios inherently involve higher towers of exponentiations and under what trade-offs?

\paragraph{Reducing levels of exponentiality.}
On the algorithmic side, the most immediate challenge is to improve the running time of the folio-computation algorithm in \cref{thm:intro_folio_algorithm_informal} to
$$2^{\poly(d,k)}\cdot |G|^{\mathbf{O}(1)},$$
or to obtain algorithms for $k$-\textsc{Disjoint Paths} with running time
$$2^{\poly(k)}\cdot |G|^{\mathbf{O}(1)}.$$
For $k$-\textsc{Disjoint Paths}, some progress has already been made for planar graphs in \cite{LokshtanovMPSZ2020Exponential} and \cite{ChoOO2023Parameterized}, where algorithms with running times $2^{\mathbf{O}(k^2)}|G|^{6}$ and $2^{\mathbf{O}(k^2)}|G|,$ respectively, were obtained. We believe that ideas from the proof of this paper, combined with those of \cite{ChoOO2023Parameterized}, may help extend these results to more general minor-closed graph classes or even in general and for folios as well.
In fact we are curious whether the (heavy) dynamic programming part in \eqref{eq_many_dp_bp_dr} 
or the running time in \eqref{eq_ola_dp_flolio} for computing $(k,d)$\textsf{-folio}$(G,R)$ 
can be reduced by one level of exponentiality.

\paragraph{Reducing polynomial dependence.}
As a final challenge, we would like to highlight the recent result of Korhonen, Mi.\ Pilipczuk, and Stamoulis \cite{KorhonenPK2024Minorcontainment}. They gave an algorithm that, given an annotated graph $(G,R)$ with $|R|\leq k,$ computes $(k,d)\mbox{-}\mathsf{folio}(G,R)$ in time
$$f(k,d)\cdot m^{1+o(1)},$$
where $m$ is the total number of vertices and edges of $G$. Moreover, it follows that this running time can be improved to
$\mathbf{O}(f(k,d)\cdot |\!|G|\!|\cdot \log |G|),$ 
for some (galactic) function $f:\Nbbb^2\to\Nbbb$\cite{KorhonenPG2026Personal}. 
For surface-embeddable graphs, the dependence on the input size can even be made linear, as proved by Golovach, Kolliopoulos, Stamoulis, and Thilikos in \cite{GolovachKST25Finding}.
The parameter dependence hidden in the function $f,$ however, remains galactic. This naturally raises the following question: is it possible to obtain a true ``best of both worlds'' theorem, combining the near-linear dependence on the input size achieved in \cite{KorhonenPK2024Minorcontainment} (or \cite{GolovachKST25Finding}) with the controlled parameter dependence of \cref{thm:intro_folio_algorithm_informal}?

\newpage
\phantomsection
\addcontentsline{toc}{section}{References}
\bibliographystyle{alphaurl}
\bibliography{optlink_arXiv}

\end{document}